\pdfoutput=1
\documentclass[abstract=on, 11pt]{scrreprt}
\usepackage{scrhack}
\addtokomafont{disposition}{\rmfamily}

\usepackage{lmodern}
\usepackage{fullpage}

\usepackage{amssymb}
\usepackage{amsmath}
\usepackage{mathtools}

\usepackage{amsthm}
\usepackage{thmtools}
\usepackage{thm-restate}
\usepackage{bm}
\usepackage[inline]{enumitem}
\usepackage{color}
\usepackage{comment}
\usepackage[table,xcdraw,dvipsnames]{xcolor}
\usepackage{float}
\usepackage[T1]{fontenc}
\usepackage[disable]{todonotes}
\usepackage{asymptote}
\usepackage{mdframed}
\usepackage[most]{tcolorbox}
\usepackage{framed}
\usepackage{mdframed}
\usepackage{bbm}
\usepackage{wrapfig}
\usepackage{multicol}
\usepackage{multirow}
\usepackage{cutwin}

\usepackage{booktabs}   
\usepackage{array}      
\usepackage{colortbl}   
\usepackage{makecell}   
\usepackage{xcolor}     

\usepackage{mathrsfs}

\usepackage[T1]{fontenc}

\definecolor{asparagus}{rgb}{0.53, 0.66, 0.42}
\definecolor{sacramentostategreen}{rgb}{0.0, 0.3, 0.15}
\definecolor{teal}{rgb}{0.0, 0.5, 0.5}
\definecolor{forestgreen}{rgb}{0.13, 0.6, 0.13}

\usepackage{hyperref}
\hypersetup{
    colorlinks=true,
    linkcolor=sacramentostategreen,
    filecolor=sacramentostategreen,
    citecolor=teal,
    urlcolor=teal,
}
\usepackage[capitalize]{cleveref}
\newtheorem{theorem}{Theorem}[section]
\newtheorem{mainthm}{Theorem}

\newtheorem{thesisproblem}{Problem}
\newtheorem{question}{Question}

\numberwithin{example}{subsection}
\newtheorem{definition}[theorem]{Definition}
\newtheorem{maindefn}{Definition}
\newtheorem{lemma}[theorem]{Lemma}
\newtheorem{claim}[theorem]{Claim}

\newtheorem{corollary}[theorem]{Corollary}

\newtheorem{remark}[theorem]{Remark}

\newtheorem{problem}{Problem}
\numberwithin{problem}{chapter}
\newtheorem*{problem*}{Problem}
\newtheorem{fact}[theorem]{Fact}

\theoremstyle{remark}
\newtheorem{model}{Model}
\numberwithin{model}{chapter}

\newcommand{\E}{\mathbb{E}}

\newcommand{\N}{\mathbb{N}}
\newcommand{\R}{\mathbb{R}}
\renewcommand{\S}{\mathbb{S}}
\newcommand{\Z}{\mathbb{Z}}

\newcommand{\cA}{\mathcal{A}}
\newcommand{\cB}{\mathcal{B}}
\newcommand{\cC}{\mathcal{C}}
\newcommand{\cD}{\mathcal{D}}
\newcommand{\cE}{\mathcal{E}}
\newcommand{\cF}{\mathcal{F}}
\newcommand{\cG}{\mathcal{G}}
\newcommand{\cH}{\mathcal{H}}

\newcommand{\cL}{\mathcal{L}}
\newcommand{\cM}{\mathcal{M}}
\newcommand{\cN}{\mathcal{N}}
\newcommand{\cO}{\mathcal{O}}
\newcommand{\cP}{\mathcal{P}}

\newcommand{\cS}{\mathcal{S}}

\newcommand{\cV}{\mathcal{V}}

\newcommand{\cX}{\mathcal{X}}

\newcommand{\mA}{\mathbf{A}}
\newcommand{\mB}{\mathbf{B}}

\newcommand{\mD}{\mathbf{D}}
\newcommand{\mE}{\mathbf{E}}
\newcommand{\mF}{\mathbf{F}}
\newcommand{\mG}{\mathbf{G}}

\newcommand{\mI}{\mathbf{I}}
\newcommand{\mJ}{\mathbf{J}}

\newcommand{\mL}{\mathbf{L}}
\newcommand{\mM}{\mathbf{M}}
\newcommand{\mN}{\mathbf{N}}

\newcommand{\mP}{\mathbf{P}}
\newcommand{\mQ}{\mathbf{Q}}
\newcommand{\mR}{\mathbf{R}}
\newcommand{\mS}{\mathbf{S}}
\newcommand{\mU}{\mathbf{U}}
\newcommand{\mV}{\mathbf{V}}
\newcommand{\mW}{\mathbf{W}}

\newcommand{\mZ}{\mathbf{Z}}

\newcommand{\valpha}{\bm{\alpha}}

\newcommand{\vmu}{\bm{\mu}}

\newcommand{\va}{\bm{a}}
\newcommand{\vb}{\bm{b}}
\newcommand{\vc}{\bm{c}}
\newcommand{\vd}{\bm{d}}
\newcommand{\ve}{\bm{e}}
\newcommand{\vf}{\bm{f}}
\newcommand{\vg}{\bm{g}}

\newcommand{\vl}{\bm{l}}
\newcommand{\vm}{\bm{m}}

\newcommand{\vq}{\bm{q}}
\newcommand{\vr}{\bm{r}}
\newcommand{\vs}{\bm{s}}

\newcommand{\vu}{\bm{u}}
\providecommand{\vv}{}
\renewcommand{\vv}{\bm{v}}
\newcommand{\vw}{\bm{w}}
\newcommand{\vx}{\bm{x}}
\newcommand{\vy}{\bm{y}}
\newcommand{\vz}{\bm{z}}

\def\to{\rightarrow}
\def\eps{\varepsilon}

\def\eps{\varepsilon}

\renewcommand{\bar}{\overline}

\usepackage{algorithm}
\usepackage{algpseudocodex}

\newcommand{\defeq}{\overset{\mathsf{\tiny def}}{=}}
\newcommand{\onev}{\mathbbm{1}}

\newcommand{\argmin}[1]{\underset{#1}{\mathrm{argmin}}}

\newcommand{\vspan}[1]{\mathsf{Span}\inparen{#1}}
\newcommand{\lspan}[1]{\mathsf{span}\inparen{#1}}

\newcommand{\kernel}[1]{\mathsf{Ker}\inparen{#1}}

\newcommand{\vdim}[1]{\mathsf{dim}\inparen{#1}}

\newcommand{\vcdim}[1]{\mathsf{VC}\inparen{#1}}
\newcommand{\Unif}[1]{\mathsf{Unif}\inparen{#1}}
\newcommand{\signv}[1]{\mathsf{sign}\inparen{#1}}
\newcommand{\suppv}[1]{\mathsf{Supp}\inparen{#1}}

\newcommand{\suchthat}{{\;\; : \;\;}}

\newcommand{\inparen}[1]{\left(#1\right)}             
\newcommand{\inbraces}[1]{\left\{#1\right\}}           
\newcommand{\insquare}[1]{\left[#1\right]}             
\newcommand{\inangle}[1]{\left\langle#1\right\rangle} 

\newcommand{\indicator}[1]{\mathbbm{1}\inbraces{#1}}

\newcommand{\abs}[1]{\ensuremath{\left\lvert #1 \right\rvert}}

\newcommand{\norm}[1]{\ensuremath{\left\lVert #1 \right\rVert}}
\newcommand{\opnorm}[1]{\ensuremath{\left\lVert #1 \right\rVert_{\mathrm{op}}}}
\newcommand{\maxnorm}[1]{\ensuremath{\left\lVert #1 \right\rVert_{\infty}}}
\newcommand{\fnorm}[1]{\ensuremath{\left\lVert #1 \right\rVert_{F}}}

\newcommand{\ip}[1]{\left\langle #1 \right\rangle}

\newcommand{\myfunc}[3]{#1\colon#2\to#3}

\usepackage{nicefrac}

\let\nfrac=\nicefrac

\renewcommand{\Pr}{\mathsf{Pr}}

\newcommand{\exv}[1]{\E\insquare{#1}}
\newcommand{\exvv}[2]{\underset{#1}{\E}\insquare{#2}}
\newcommand{\expv}[1]{\mathsf{exp}\inparen{#1}}
\newcommand{\prv}[1]{\Pr\insquare{#1}}
\newcommand{\prvv}[2]{\underset{#1}{\Pr}\insquare{#2}}
\newcommand{\cov}[1]{\mathsf{Cov}\insquare{#1}}
\newcommand{\var}[1]{\mathsf{Var}\insquare{#1}}

\newcommand{\logv}[1]{\log\inparen{#1}}
\newcommand{\lnv}[1]{\ln\inparen{#1}}
\newcommand{\diag}[1]{\mathsf{diag}\inparen{#1}}
\newcommand{\diam}[1]{\mathsf{diam}\inparen{#1}}

\newcommand{\so}[1]{\mathbb{O}\inparen{#1}}

\makeatletter
\newcommand{\customlabel}[2]{%
   \protected@write \@auxout {}{\string \newlabel {#1}{{#2}{\thepage}{#2}{#1}{}} }%
   \hypertarget{#1}{#2}
}
\makeatother

\newcommand\numberthis{\addtocounter{equation}{1}\tag{\theequation}}

\newcounter{casenum}

\newcounter{subcasenum}

\newcounter{casenump}

\newcommand{\casep}[2]{
    \ifthenelse{\equal{\value{casenump}}{0}}{
    \vskip.5\baselineskip\par\noindent
    }{}
    {\it Case \arabic{casenump}:} {\it #1}
    \vskip0.1\baselineskip
    \begin{addmargin}[1.5em]{1em}
    #2
    \end{addmargin}
    \addtocounter{casenump}{1}
}

\newcounter{subcasenump}

\newcommand{\lecture}[4]{
  \pagestyle{myheadings}
  \thispagestyle{plain}
  \newpage
  \setcounter{page}{1}
  \noindent
  \begin{center}
    \framebox{
      \vbox{\vspace{2mm}
      \hbox to 6.28in {{\bf Informal Whiteboard Talks
                          \hfill Note #4}}
      \vspace{6mm}
      \hbox to 6.28in {{\Large \hfill #1  \hfill }}
      \vspace{5mm}
      \hbox to 6.28in {{#2 \hfill #3}}
      \vspace{2mm}}
    }
  \end{center}
  \markboth{#1}{#1}
  \vspace*{4mm}
}

\newcommand{\floor}[1]{\left\lfloor #1 \right\rfloor}
\newcommand{\ceil}[1]{\left\lceil #1 \right\rceil}
\newcommand{\rank}[1]{\mathsf{rank}\inparen{#1}}

\usepackage{cleveref}

\newcommand{\smallpar}[1]{\textbf{#1\quad}}
\newcommand{\smallit}[1]{\textit{#1\quad}}

\usepackage{ulem}

\newcommand{\detv}[1]{\mathrm{det}\inparen{#1}}

\newcommand{\conv}[1]{\mathsf{conv}\inparen{#1}}
\newcommand{\vol}{\mathsf{vol}}
\newcommand{\opt}{\mathsf{OPT}}
\newcommand{\lsetext}{\textsc{LogSumExp}}
\newcommand{\ftilde}{\widetilde{f}}
\newcommand{\xtilde}{\widetilde{\vx}}
\newcommand{\fhat}{\widehat{f}}
\newcommand{\lse}{\mathsf{lse}}

\definecolor{ForestGreen}{RGB}{34, 139, 34}
\newcommand{\greenbox}[1]{
    \vspace{0.15cm}
    { \small
    \begin{tcolorbox}[enhanced, colback=ForestGreen!12, parbox=false]
    #1
    \end{tcolorbox}
    }
}

\definecolor{forest}{RGB}{0,155,85}

\newcommand{\fixout}{\bgroup\markoverwith{\textcolor{forest}{\rule[0.5ex]{2pt}{0.4pt}}}\ULon}

\renewcommand{\phi}{\varphi}

\newif\ifnotes
\makeatletter
\newcommand{\nnote}[1]{\@bsphack\ifnotes{$\ll$\textsf{\color{blue} Naren: { #1}}$\gg$}\fi\@esphack}
\makeatother
\notesfalse

\newcommand{\hhat}{\widehat{h}}
\newcommand{\hstar}{h^{\star}}

\newcommand{\wstar}{\vw^{\star}}
\newcommand{\what}{\widehat{\vw}}
\newcommand{\vstar}{\vv^{\star}}
\newcommand{\patch}[1]{\mathsf{patch}\inparen{#1}}
\newcommand{\patchw}{\mathsf{patch}}
\newcommand{\fadv}{\cF_{\mathsf{adv}}}
\newcommand{\memcap}[2]{\mathsf{mcap}_{#1}\inparen{#2}}
\newcommand{\Sclean}{S_{\mathsf{clean}}}
\newcommand{\Sbackdoor}{S_{\mathsf{adv}}}
\newcommand{\epsadv}{\eps_{\mathsf{adv}}}
\newcommand{\epsclean}{\eps_{\mathsf{clean}}}
\newcommand{\robustloss}{\cL_{\fadv(\hstar)}}
\newcommand{\xstar}{\vx^{\star}}

\newcommand{\xt}[1]{\vx_t^{(#1)}}
\newcommand{\actionset}{\cX}  
\newcommand{\xplus}{\vx^{+}}

\usepackage{cleveref}

\usepackage{ulem}

\newcommand{\vdin}{\vd_{\mathsf{in}}}
\newcommand{\vdout}{\vd_{\mathsf{out}}}

\newcommand{\din}{d_{\mathsf{in}}}
\newcommand{\ddet}{\vd_{\mathsf{det}}}

\newcommand{\rand}{\mathsf{rand}}

\newcommand{\Astar}{\mA^{\star}}

\newcommand{\Dstar}{\mD^{\star}}

\newcommand{\Lstar}{\mL^{\star}}
\newcommand{\cLstar}{\cL^{\star}}
\newcommand{\Lhat}{\widehat{\mL}}

\newcommand{\pbar}{\overline{p}}

\newcommand{\utwostar}{\vu_2^{\star}}
\newcommand{\utwohat}{\widehat{\vu_2}}

\newcommand{\astar}{\va^{\star}}

\definecolor{ForestGreen}{RGB}{34, 139, 34}

\newcommand{\xhat}{\widehat{\vx}}

\renewcommand{\vmu}{\bm{\mu}}

\newcommand{\dtwo}{dist}

\newcommand{\vbrho}{\eps^{-2}\inparen{\log d}^2\logv{\nfrac{d}{\eps}}}
\newcommand{\vlambda}{\bm{\lambda}}
\newcommand{\mLambda}{\mathbf{\Lambda}}
\newcommand{\mSigma}{\mathbf{\Sigma}}

\newcommand{\zstar}{\vz^{\star}}
\newcommand{\Fstar}{F^{\star}}
\newcommand{\mtilde}{\widetilde{m}}

\newcommand{\vrho}{\bm{\rho}}

\newcommand{\polynorm}[1]{\norm{#1}_{\cG,\rho,\infty, S}}
\newcommand{\subgnorm}[1]{\norm{#1}_{\psi_2}}
\newcommand{\gnorm}[2]{\norm{#1}_{\cG_{#2}}}
\newcommand{\gnorml}[2]{\norm{#1}_{\cG_{#2}(\vlambda)}}
\newcommand{\nnz}[1]{\mathsf{nnz}\inparen{#1}}

\renewcommand{\phi}{\varphi}
\newcommand{\oprox}{\cO_{\mathsf{prox}}}

\newcommand{\oms}{\cO_{\mathsf{MS}}}

\usepackage{tikz}
\usepackage{newpxmath}
\usepackage{newpxtext}

\usepackage[margin=1in]{geometry}

\allowdisplaybreaks

\usepackage[
backref=true,
natbib=true,
style=trad-alpha,
maxalphanames=8,
sorting=nyt
]{biblatex}
\numberwithin{equation}{section}

\renewcommand{\vv}{\bm{v}}

\begin{document}

\title{A Geometric Approach to Problems in Optimization and Data Science}
\author{Naren Sarayu Manoj}
\date{\today}

\begin{titlepage}
   \begin{center}
       \vspace*{4cm}

       {\huge \textbf{A Geometric Approach to Problems in Optimization and Data Science}}

       \vspace{2cm}
BY\\
    NAREN SARAYU MANOJ

         \vspace{4cm}   
       A thesis submitted
       \\
in partial fulfillment of the requirements for
\\
the degree of
\\
Doctor of Philosophy in Computer Science
\\
at the
\\
TOYOTA TECHNOLOGICAL INSTITUTE AT CHICAGO
\\
Chicago, Illinois
\\
May 2025
       \vspace{1cm}
     
       \textbf{Thesis Committee:}\\
       \begin{align*}
           & \text{Avrim Blum} & \text{(Thesis Co-Advisor)}\\
           & \text{Michael Kapralov} & \\
           & \text{Sepideh Mahabadi} & \\
           & \text{Yury Makarychev} & \text{(Thesis Co-Advisor)}
       \end{align*}
           
   \end{center}
\end{titlepage}

\begin{abstract}
We give new results for problems in computational and statistical machine learning using tools from high-dimensional geometry and probability.

We break up our treatment into two parts. In Part I, we focus on computational considerations in optimization. Specifically, we give new algorithms for approximating convex polytopes in a stream, sparsification and robust least squares regression, and dueling optimization.

In Part II, we give new statistical guarantees for data science problems. In particular, we formulate a new model in which we analyze statistical properties of backdoor data poisoning attacks, and we study the robustness of graph clustering algorithms to ``helpful'' misspecification.
\end{abstract}

\newpage
{\small \section*{Acknowledgments} This journey was a joint effort with many people. So, many thanks are in order.

\smallpar{Research.} I am extremely fortunate to have the best advisors I can ask for -- Avrim Blum and Yury Makarychev. I am always amazed by their intuition and extremely clear approach to research and problem-solving. I am thankful to Avrim for teaching me his strategies for formulating research questions and imparting in me some idea of what the right answers to problems look like on a conceptual level. I am grateful to Yury for giving me the technical confidence to continue doing research in convex geometry, specifically at a time when I had no clear personal problem taste. Avrim and Yury are also infinitely generous in many aspects -- in trust, in wisdom (both in TCS and in life), in finances (I could attend any conference, any visit, etc), in time, in patience, and in encouragement. I hope to continue learning all I can from them even when I am no longer at TTIC.

I have awesome mentors in addition to Avrim and Yury. I am indebted to Aditya Bhaskara and Pritish Kamath, who helped me navigate (technically and otherwise) the most challenging stretches of the PhD. They are exceedingly dependable sources of optimism and sage advice. I am grateful to Sepideh Mahabadi for introducing geometric streaming algorithms to me and for serving on my thesis committee. Thanks to Aditya Bhaskara, Michael Kapralov, and Ivan Tyukin for fun research visits and to Michael for serving on my thesis committee. Thanks to Cameron Chen and Kunal Nagpal for hosting me for a rewarding research internship. Going back a bit further, thanks to Rob Burridge, Constantine Caramanis, Alex Dimakis, Adam Klivans, Calvin Lin, and Eric Price for encouraging me to try CS research in the first place. And, going back even further, thanks a lot to Art of Problem Solving for first showing me how entertaining problem solving is (irrespective of whether the problem actually gets solved).

Research is a blast when done with great collaborators. The research I was part of between 2018 and the time of this writing is joint with Aditya Bhaskara, Avrim Blum, Cameron Chen, Travis Dick, Alex Dimakis, Surbhi Goel, Meghal Gupta, Agastya Vibhuti Jha, Matt Jordan, Michael Kapralov, Gene Li, Yury Makarychev, Davide Mazzali, Kunal Nagpal, Max Ovsiankin, Kumar Kshitij Patel, Aadirupa Saha, Nati Srebro, Weronika Wrzos-Kaminska, Chloe Yang, and Hongyang Zhang.

The Chicago area theory community is an ``IDEAL'' place for TCS/ML research. Many thanks to Saba Ahmadi, Idan Attias, Siddharth Bhandari, Aditya Bhaskara, Avrim Blum, Sam Buchanan, Julia Chuzhoy, Chao Gao, Varun Gupta, Frederic Koehler, Zhiyuan Li, Cong Ma, Sepideh Mahabadi, Yury Makarychev, Aaron Potechin, Miki Racz, Liren Shan, Dravyansh Sharma, Nati Srebro, Madhur Tulsiani, Ali Vakilian, Aravindan Vijayaraghavan, Jingyan Wang, Haifeng Xu, and the members of the IDEAL institute for creating a collegial environment for TCS in Chicago. Also, thanks to Avrim, Julia, Chao, Sepideh, Yury, Aaron, Nati, Madhur, and Haifeng for teaching me mindblowing ideas in TCS through wonderful courses, and to Karen Livescu, David McAllester, Greg Shakhnarovich, Matthew Turk, and Matthew Walter for getting me out of my bubble and reminding me that there are other subfields in CS too.

I am grateful to the optimization, statistics, and TCS communities for humoring my questions and shaping the way I think about the fields. In addition to those mentioned above, I am lucky to have learned from Deeksha Adil, Omar Alrabiah, Alex Andoni, Moïse Blanchard, Yeshwanth Cherapanamjeri, Mika Göös, Venkat Guruswami, Piotr Indyk, Arun Jambulapati, Praneeth Kacham, Satyen Kale, Sushrut Karmalkar, Ishani Karmarkar, Guy Kornowski, Yin Tat Lee, Lawrence Li, Yang Liu, Shay Moran, Christopher Musco, Jelani Nelson, Swati Padmanabhan, Ashwin Pananjady, Ankit Pensia, Jon Schneider, Sushant Sachdeva, Tselil Schramm, Vatsal Sharan, Aaron Sidford, Kevin Tian, Kabir Verchand, David Woodruff, David Wu, Taisuke Yasuda, Yibin Zhao, and Nikita Zhivotovskiy.

Also, thanks to Ainesh Bakshi, Guy Blanc, Elbert Du, Ezra Edelman, Noah Golowich, Siddhartha Jain, Ajil Jalal, Justin Kang, Caleb Koch, Dheeraj Nagaraj, Chirag Pabbaraju, Nived Rajaraman, Vinod Raman, Ayush Sekhari, Abhishek Shetty, Kiran Shiragur, Unique Subedi, Ewin Tang, and many others I am certainly forgetting for being a familiar, friendly presence at conferences.

I am grateful to Chen Dan, Misha Khodak, Siddharth Prasad, Dravyansh Sharma, and Ellen Vitercik for maintaining our reading group in 2020-2021. It was a great source of motivation during the pandemic.

Many thanks to Adam Bohlander, Rose Bradford, Erica Cocom, Chrissy Coleman, Jessica Jacobson, Brandie Jones, Celeste Ki, Deree Kobets, Randy Landsberg, Mary Marre, Alicia McClarin, Amy Minick, and Jerry Randle for making my experience at TTIC smooth and enjoyable.

The research I partook in between Fall 2019 and April 2025 was supported by a USA NSF graduate research fellowship, USA NSF awards CCF-1815011 and ECCS-2216899, and the DARPA GARD program.

\smallpar{Personal.} The last few years were incredibly fun even beyond the research. This is because I have family and friends who ensured I took every opportunity to have a great time.

Traveling was a staple of the last few years that made the time all the more enjoyable. I owe a ton to Vahid Asadi, Armin Askari, Anurag Bakshi, Karthik Bala, Shreyas Balaji, Sumanth Balaraman, Kevin Choy, Rohit Datta, Anindit Gopalakrishnan, Harsh Goyal, Preetham Gujjula, Meghal Gupta, Rishub Jain, Vanshika Jain, Divyam Khandelwal, Kapil Krishnakumar, Alex Lee, Eric Lee, Gene Li, Yang Liu, Prabhat Nagarajan, Madhav Narayan, Avinash Nayak, Eley Ng, Evonne Ng, Abhit Patil, Brahma Pavse, Josh Pham, Binita Shah, Han Shao, Taylor Sihavong, Vaidehi Srinivas, Aditi Srinivasan, Darshan Thaker, Heather Tolcher, Alice Wei, Rahul Yesantharao, George Zhang, Grace Zhang, Rachel Zhang, and Zihao Zhu for some amazing and recharging trips. 

Also, thanks to Daniel Ho, Aida Lu, Sachin Mehta, Heather Tolcher, Lekha Yesantharao, Rahul Yesantharao, George Zhang, and Zihao Zhu (for good memes whenever I visited Clear Lake), to Anurag Bakshi, Rohit Datta, Divyam Khandelwal, Taylor Sihavong, and Rahul Yesantharao (for giving me reasons to visit New York and London often), to Rohit Datta, Evonne Ng, and Darshan Thaker (for the most entertaining and irrelevantly named group chat), and Meghal Gupta (for degenerate games of chess).

Locally, I am lucky to have had Sudarshan Babu, Riddhish Bhalodia, Dimitar Chakarov, Antares Chen, Emma Corless, Avichal Goel, Nirmit Joshi, Anmol Kabra, Hriday Kamshatti, Sushrut Karmalkar, Divyam Khandelwal, Frederic Koehler, Sky Lee, Gene Li, Yanhong Li, Marko Medvedev, Tushant Mittal, Adán Medrano Martín del Campo, Olga Medrano Martín del Campo, Amin Mohamadi, Omar Montasser, Abhijit Mudigonda, Keziah Naggita, Rachit Nimavat, Max Ovsiankin, Simran Pabla, Aditya Paliwal, Ankita Pasad, Kumar Kshitij Patel, Aniruddh Raghu, Goutham Rajendran, Kavya Ravichandran, Aravind Reddy, Han Shao, Pushkar Shukla, Vaidehi Srinivas, Shashank Srivastava, Kevin Stangl, Alec Sun, Akilesh Tangella, Erasmo Tani, Shubham Toshniwal, Antony Yun, David Yunis, Amanda Zhang, and Derek Zhu to accompany me for the last few years. I am grateful for their friendship. A special shoutout goes to Anmol (my housemate for two years), to Keziah, Ankita, Han, and Shashank, (my officemates at TTIC), to Sudarshan, Anmol, Gene, Max, Ankita, Kshitij, Kavya, Han, Vaidehi, Shashank, and David (for a particular sense of camaraderie), to Avichal, Hriday, Divyam, Fred, Aditya, and Antony (Block 37 climbing crew), and to Amanda (for the best baked goods ever). Not so locally, thanks to Marina Drygala, Grzegorz Głuch, Agastya Vibhuti Jha, Davide Mazzali, Abhit Patil, Kshiteej Sheth, and Weronika Wrzos-Kaminska for making me feel at home at EPFL.

\smallpar{Family.} My family has always unconditionally supported me. My mother, Sarayu Manoj, spends practically all her waking hours ensuring that my brother Pavan and I have strong foundations for both academics and life. Because of her, we are always happy, well-fed, and pursuing worthwhile activities. I am pretty sure I first learned to enjoy math from her teachings, which began before we even went to school. My father, Manoj Som, is probably the most easygoing person I know. By example, he taught us to live in the present, how to navigate stressful times, and the importance of maintaining a healthy physical lifestyle. He also imparted in us a strong sense of adventure. All of this, plus his ``go-with-the-flow'' approach to life, was necessary in navigating the challenges over the last few years. Importantly, our parents taught us to love learning for learning's sake. My brother, Pavan Manoj, showed me how to be levelheaded and diligent. He listened to my dramatic rants whenever things were not going well and always offered realistic advice.

Beyond my immediate family, I am lucky to have a supportive and fun-loving extended family. Particularly, thanks to Shrihari Gopalakrishna and Srikanth Iyengar for encouraging me to pursue research in the mathematical sciences and for keeping an eye on me during the PhD. It is a huge advantage to have examples within the family who demonstrate that a particular life track/career plan is both viable and fun. Last but not least, I am grateful to our family-friend network in the US for their well-wishes.}
\newpage
\tableofcontents
\newpage
\chapter{Introduction}

A basic pipeline in modern machine learning is to construct a model based on observations in the real world and use the model to make inferences on unseen data. To execute this pipeline, there are at least two aspects we need a deep understanding of. One is \textit{algorithmic} -- given input data, how do we efficiently build a model that captures the phenomenon we are interested in? Another is \textit{statistical} -- given some sense of the data generation process and the model fitting procedure, should we believe the model's predictions on unseen data?

With this delineation in mind, in this thesis, we will look at important considerations in these components. To address the algorithmic aspects, we will \textbf{design and analyze fast optimization primitives, motivated by data science applications}. On the statistical front, we will \textbf{study recovery/inference when the input is adversarial or misspecified}. A common theme between all the problems we study is that the required technical methods involve understanding the underlying high-dimensional probabilistic and geometric phenomena. 

\section{Our geometric motivation and quick summary of results}
\label{sec:thesisintro_overview_quick}


A helpful guiding question to keep in mind is -- \textit{which statistical and computational primitives benefit from a more geometrically-aware analysis?}



As a step towards answering this question on the computational end, we will spend slightly more than half of the thesis studying the \textit{ellipsoidal approximation problem} and its applications. In the ellipsoidal approximation problem, we are given a convex body $K \subset \R^d$ (meaning that $K$ is convex, compact, and has nonzero Lebesgue measure). We would like to find an ellipsoid $\cE(K)$, center $\vc \in \R^d$, and \textit{distortion} $\triangle \ge 1$ for which
\begin{align}
    \vc + \cE(K) \subseteq K \subseteq \vc + \triangle \cdot \cE(K).\label{eq:thesisintro_ellipsoidal_approx}
\end{align}
Intuitively, the closer $\triangle$ is to $1$, the more faithfully $\cE(K)$ approximates $K$. 

The ellipsoidal approximation problem \eqref{eq:thesisintro_ellipsoidal_approx} and its variants appear in many problems in optimization. For example, algorithm families for convex programming such as the ellipsoid method \cite{kha80} and interior point methods \cite{nn94} can be thought of as approximating a convex body with some sequence of ellipsoids, and the distortion gives intuition for the kinds of convergence rates one can expect for these algorithms. Low-distortion ellipsoidal approximations also give us ``preconditioners'' for convex bodies, with applications to sampling \cite{cv15}, reinforcement learning \cite{li2019stochastic}, and integer programming \cite{len83}. Furthermore, ellipsoidal approximations give us a natural way to succinctly summarize $K$ (up to a factor $\triangle$ loss), as describing an ellipsoid in $\R^d$ only requires us to write down a matrix in $\R^{d \times d}$.

The distortion $\triangle$ in \eqref{eq:thesisintro_ellipsoidal_approx} and the runtime for finding $\cE(K)$ are the principal quantities that we will want to control when we solve the ellipsoidal approximation problem. To see what kinds of distortions one can hope for in \eqref{eq:thesisintro_ellipsoidal_approx}, recall \textit{John's theorem}.
\begin{theorem}[John's Theorem {\cite{john1948}}]
\label{thm:thesisintro_john}
Let $\vc + \cE(K)$ be the ellipsoid of maximal volume contained within $K$. Then, we have
\begin{align*}
    \vc + \cE(K) \subseteq K \subseteq \vc + \triangle \cdot \cE(K),
\end{align*}
where $\triangle \le d$. Further, if $K$ is origin-symmetric, then this improves to $\triangle \le \sqrt{d}$. Finally, there exist convex bodies $K$ for which no ellipsoid can approximate $K$ to distortion better than $d$, and there exist origin-symmetric convex bodies $K$ for which no ellipsoid can approximate $K$ to distortion better than $\sqrt{d}$.
\end{theorem}
Thus, for all of our results for problems of the form \eqref{eq:thesisintro_ellipsoidal_approx}, it will be helpful to compare against the existential benchmark given by \Cref{thm:thesisintro_john}.

Motivated by the above, we will begin with studying \eqref{eq:thesisintro_ellipsoidal_approx} in a fairly general setting. We assume the convex body $K$ is specified by the convex hull of points $\vz_1,\dots,\vz_n \in \R^d$. For additional reasons we will discuss momentarily, we will further assume that \begin{enumerate*}[label={(\arabic*)}] \item we are only given single-pass streaming access to the $\vz_1,\dots,\vz_n$ and \item our algorithm must be memoryless, which means when it receives point $\vz_t$, it is not allowed to remember any of the previous points $\vz_1,\dots,\vz_{t-1}$. \end{enumerate*} In this streaming model, we will give the first efficient algorithms that achieve nearly worst-case optimal approximation factors, in that they nearly agree with John's theorem (\Cref{thm:thesisintro_john}). We then use these procedures to obtain the first single-pass streaming algorithms for finding low-distortion \textit{convex hull coresets}, which are subsets of the points from $\vz_1,\dots,\vz_n$ whose convex hull approximates the convex hull of the original set. Our results here will set the stage for the runtimes and distortions we can expect for the other ellipsoidal approximation problems we study. For more details, see \Cref{sec:thesisintro_ellipsoid} and \Cref{chapter:streamingellipsoids}.

We then give algorithms for \eqref{eq:thesisintro_ellipsoidal_approx} in cases in which the convex sets $K$ are symmetric, highly structured, and given offline. Once again, the distortions we get will be nearly worst-case optimal (and can be made arbitrarily close to worst-case optimal). We will leverage this primitive to give state-of-the-art algorithms for multidistributional linear regression, in which we would like one parameter vector that simultaneously minimizes the least squares loss for $m$ different linear regression problems. We will also combine this ellipsoidal approximation primitive with extra geometric tools to obtain nearly optimal existential and algorithmic results for \textit{sparsification}. In the sparsification problem, our goal is to considerably simplify loss functions that can be expressed as the sum of many individual structured losses. We will then apply these newly built tools to give computational improvements for the problem of minimizing the sum of Euclidean norms, which encompasses several well-studied optimization problems. We will discuss this more precisely in \Cref{sec:thesisintro_blw} and \Cref{chapter:sparsifyingnorms,chapter:regression}.

The remainder of the thesis is dedicated to exploring some statistical and algorithmic consequences of the behaviors of random vectors in high dimensions. We focus on the following three fundamental facts (we focus on the distribution $\vg \sim \cN\inparen{0,\frac{\mI_d}{d}}$ for convenience of presentation here, though in our results, we will have to reason about different distributions):
\begin{enumerate}
    \item Random Gaussian vectors have a nontrivial correlation with a fixed direction, with nontrivial probability. Formally, let $\vz \in \R^d$ be such that $\norm{\vz}_2 = 1$. Then,
    \begin{align*}
        \prvv{\vg \sim \cN\inparen{0,\frac{\mI_d}{d}}}{\ip{\vg,\vz} \ge \frac{2}{\sqrt{d}}} \ge \frac{1}{50}.
    \end{align*}\label{item:thesisintro_random_one}
    \item Random Gaussian vectors mostly miss a fixed low-dimensional subspace. Formally, if $\mU \in \R^{d \times s}$ where $d \ge s$ and the columns of $\mU$ form an orthonormal basis for the subspace spanned by the columns of $\mU$, then for some universal constant $C$, with probability $\ge 1-\delta$, we have
    \begin{align*}
        \norm{\mU^{\top}\vg}_2 \le C\inparen{\sqrt{\frac{s}{d}}+\sqrt{\frac{\logv{\nfrac{1}{\delta}}}{d}}}.
    \end{align*}
    So, if $d \gg s$, then $\vg$ mostly lies in the nullspace of $\mU^{\top}$. This also can be seen as a converse to Property \ref{item:thesisintro_random_one} (take $s=1$).\label{item:thesisintro_random_two}
    \item Random Gaussian vectors are well-spread. Consider $\vg$ distributed as above. We know that $\norm{\vg}_2 \approx 1$, and by using the standard fact that $\maxnorm{\cdot} \le \norm{\cdot}_2$, we can certainly conclude $\maxnorm{\vg} \lessapprox 1$. But, this equality case only happens when most of the mass of $\vg$ is located on very few coordinates of $\vg$. Under our distributional assumption on $\vg$, this is incredibly unlikely. Instead, for some universal constant $C$, we have with probability $\ge 1-\delta$ that
    \begin{align*}
        \maxnorm{\vg} \le C\inparen{\sqrt{\frac{\log d}{d}} + \sqrt{\frac{\logv{\nfrac{1}{\delta}}}{d}}}.
    \end{align*}
    This improves over the trivial bound by nearly a factor of $\sqrt{d}$. Furthermore, the cosmetic resemblance to the statement of Property \ref{item:thesisintro_random_two} is no coincidence -- both properties can be established in almost the same way.\label{item:thesisintro_random_three}
\end{enumerate}
In \Cref{sec:thesisintro_dueling} and \Cref{chapter:binsearch}, we will use Property \ref{item:thesisintro_random_one} to give the first algorithm for a realistic generalization of \textit{dueling optimization}, a common preference-based optimization paradigm. When placed in a natural online learning setting, our algorithm incurs regret that is optimal up to constant factors in the worst case. The main technical idea is to use the first property to guess weakly-correlated descent steps to optimize a function $f$ without gradient and evaluation access to $f$.

With that, we round out the optimization section of the thesis and move onto understanding geometric phenomena in statistical problems (problems in which we are also interested in inference and recovery instead of computation alone) under input corruptions. In \Cref{sec:thesisintro_backdoor} and \Cref{chapter:backdoor}, we establish a statistical framework within which we can analyze \textit{backdoor data poisoning attacks}, a type of train-time adversarial attack on machine learning classifiers. We will use this framework to give a mathematical justification for the empirically observed phenomenon that planting backdoors into overparameterized models is ``easy'' -- our provably successful attack will in fact be a simple randomized construction whose analysis follows directly from Property \ref{item:thesisintro_random_two}. Motivated by this, we use our framework to build a much more general statistical theory around backdoor attacks. Along the way, we try to answer the questions \begin{enumerate*}[label={(\arabic*)}] \item When is a machine learning problem susceptible to a backdoor attack? \item Which natural machine learning problems can be successfully backdoored, and which cannot? \item What are some algorithmic strategies that one can use to mitigate a backdoor attack? \end{enumerate*}

Finally, in \Cref{sec:thesisintro_clustering} and \Cref{chapter:cluster}, we give the first robustness guarantees for spectral clustering, a popular and practical clustering algorithm, under a more general generative model for the graph than what is typically considered. The main technical challenge is to understand how the entries of eigenvectors change under random perturbations to the underlying signal matrix. Our final argument can be viewed as establishing Property \ref{item:thesisintro_random_three} for the distribution $\vv-\vv^{\star}$, where $\vv^{\star}$ is the ``population'' eigenvector and $\vv$ is the corresponding eigenvector after adding the noise, even though the distribution of $\vv-\vv^{\star}$ is not Gaussian. Interestingly, this is one of the first examples we are aware of in which we can apply entrywise eigenvector perturbation theory to high-rank signal matrices -- the techniques we use were originally developed and used for exclusively low-rank signal matrices.

\section{Results -- algorithms}

Our goal in \Cref{part:algos} is to give efficient algorithms for several natural problems in data science in which a high-dimensional geometric interplay is central. As mentioned in \Cref{sec:thesisintro_overview_quick}, most of the problems in this section will involve studying the ellipsoidal approximation problem \eqref{eq:thesisintro_ellipsoidal_approx} in some form.

\subsection{Streaming ellipsoidal approximations of convex polytopes}
\label{sec:thesisintro_ellipsoid}

We begin with the main problem of this subsection.

\begin{thesisproblem}
\label{problem:intro_ellipsoids}
Let $Z = \conv{\inbraces{\vz_1,\dots,\vz_n}}$. Suppose an algorithm receives streaming access to $Z$ (i.e., it receives the points $\vz_1,\dots,\vz_n$ one-at-a-time). Can we find algorithms that maintain translates $\vc_1,\dots,\vc_t$ and approximating bodies $\widehat{Z_1},\dots,\widehat{Z_t}$ such that both of the below guarantees hold:
\begin{itemize}
    \item for all timesteps $t$, we have $Z_t \coloneqq \conv{\inbraces{\vz_1,\dots,\vz_t}} \subseteq \vc_t + \widehat{Z_t}$;
    \item at the end of the stream, we have for some $0 < \alpha < 1$ that $\vc_n + \frac{1}{\alpha} \cdot \widehat{Z_n} \subseteq Z \subseteq \vc_n + \widehat{Z_n}$ or $\widehat{Z_n} \subseteq Z \subseteq \vc_n+\frac{1}{\alpha}\cdot\inparen{\widehat{Z_n}-\vc_n}$.
\end{itemize}
We are interested in the following two different types of approximating bodies $\widehat{Z_t}$:
\begin{enumerate}
    \item the $\widehat{Z_t}$ must be ellipsoids;
    \item the $\widehat{Z_t}$ must be the convex hull of a subset of the points $\vz_1,\dots,\vz_t$.
\end{enumerate}
For each algorithm, we would like the distortion (interchangeably used with ``approximation factor'') $\nfrac{1}{\alpha}$ to be as small as possible and for the space complexity to be as small as possible.
\end{thesisproblem}

In the language of \Cref{sec:thesisintro_overview_quick}, the quantity $1/\alpha$ is referred to as a \textit{distortion} or \textit{approximation factor}. 

Observe that the first objective is asking us to build an $\ell_2$-approximation of $Z$. Additionally, note that the second objective of \Cref{problem:intro_ellipsoids} amounts to building a coreset for $Z$ -- i.e., we are asking for an algorithm that chooses a subset of the $\vz_1,\dots,\vz_n$ that approximates $Z$. 

\smallpar{Motivating example.} Let us discuss why we study the streaming setting. Consider a case where we have a very large dataset consisting of points $\vz_1,\dots,\vz_n \in \R^d$. We would like to produce a succinct summary of this dataset. Since the dataset has a large number of observations, the algorithm to calculate the summary is not allowed to remember all the $\vz_i$s. Instead, we will allow one-pass streaming access to the $\vz_i$. In particular, our algorithm will be allowed to read one $\vz_i$ at a time. The algorithm cannot make assumptions on the order of the points in the stream; in particular, the stream could be adaptively adversarially ordered.


\smallpar{Our results.} We give an informal overview of the results we achieve for \Cref{problem:intro_ellipsoids} and defer the more precise statements of the results to \Cref{chapter:streamingellipsoids}.
\begin{itemize}
    \item When the $\widehat{Z_t}$ must be ellipsoids, we give an algorithm that needs to store only $O(d^2)$ floating point numbers while achieving an approximation factor of $O\inparen{\min\inbraces{\kappa, d\log \kappa}}$, where $\kappa$ is the \textit{aspect ratio} of $Z$ (essentially capturing how skewed $Z$ is). The algorithm also has runtime $\widetilde{O}(nd^2)$.
    \item When the $\widehat{Z_t}$ must be the convex hull of a subset of the points $\vz_1,\dots,\vz_t$, we give an algorithm that chooses at most $O(d\log \kappa^{\mathsf{OL}})$ vertices while achieving an approximation factor of $O\inparen{d\log d + d\log \kappa^{\mathsf{OL}}}$ (here, $\kappa^{\mathsf{OL}}$ is an online variant of the aspect ratio term from the previous part). The algorithm here will involve calling the ellipsoidal approximation algorithm from the previous part as a subroutine.
    \item Additionally, if we are in the special case where $Z$ is centrally symmetric about the origin, the approximation factors of our algorithms improve to $O(\sqrt{d\log\kappa})$ and $O(\sqrt{d\log d + d\log\kappa^{\mathsf{OL}}})$ for the ellipsoid and coreset settings, respectively.
\end{itemize}
We remark that by John's Theorem (\Cref{thm:thesisintro_john}), the approximation factors of our algorithms are nearly optimal. In fact, the dependence on the dimension $d$ nearly matches that of the worst-case optimal \textit{offline} solution, and so we only lose terms that are logarithmic (or sublogarithmic, in the symmetric case) in the aspect ratio of $Z$ (while we also sometimes lose factors logarithmic in $d$, we can show that if the number of points $n$ is polynomial in $d$, then the aspect ratio must be $\mathsf{poly}(d)$, so in this practical regime, the $\log d$ dependences are unavoidable anyway).

We formally study this problem in \Cref{chapter:streamingellipsoids}.

\smallit{Bibliographic notes.} The material discussed in this section is based on a sequence of works joint with Yury Makarychev and Max Ovsiankin published at COLT 2022 \cite{mmo22} and STOC 2024 \cite{mmo23}.

\subsection{Block Lewis weights and applications}
\label{sec:thesisintro_blw}

Next, we will look at problems where one of the key technical ingredients is a low-distortion ellipsoidal approximation for a particular symmetric convex set. Our main contribution will be the design and application of a geometric construction called \textit{block Lewis weights}.

Throughout this section, it will be helpful to keep in mind the following intended applications.

\smallpar{Motivating example.} Suppose we have a collection of least squares linear regression problems, each of which is given by the designs and responses $(\mA_{S_i}, \vb_{S_i})$, where $\mA_{S_i} \in \R^{\abs{S_i} \times d}$ and $\vb_{S_i} \in \R^{\abs{S_i}}$ and $S_i$ denotes the index subset of all measurements that belong to problem $i$. Let $\mA \in \R^{(\sum \abs{S_i}) \times d}$ denote the matrix formed by stacking all $m$ designs and $\vb \in \R^{\sum \abs{S_i}}$ denote the vector formed by stacking all $m$ responses in the same way. In several settings such as in collaborative or multidistributional learning, it makes sense to ask for a parameter vector $\xhat$ for which
\begin{align}
    \max_{1 \le i \le m} \norm{\mA_{S_i}\xhat-\vb_{S_i}}_2 \le (1+\eps)\min_{\vx\in\R^d} \max_{1 \le i \le m} \norm{\mA_{S_i}\vx-\vb_{S_i}}_2.\label{eq:thesisintro_maxnorm}
\end{align}
Objective \eqref{eq:thesisintro_maxnorm} is the natural formulation for the problem of minimizing the multidistributional linear regression loss. Indeed, we can think of each design-observation pair $(\mA_{S_i},\vb_{S_i})$ as samples from a particular distribution $\cD_i$, and we would like our model $\xhat$ to perform reasonably well on all $m$ distributions $\cD_1,\dots,\cD_m$.

More generally, solving \eqref{eq:thesisintro_maxnorm} subsumes min-max fair least-squares regression, distributionally robust linear regression, and $\ell_{\infty}$ regression (which can be seen by choosing $\abs{S_i}=1$ for all $i$). Furthermore, if we let all the $\mA_{S_i} = \mI_d$, then \eqref{eq:thesisintro_msn} solves the \textit{minimum enclosing ball} problem from computational geometry, which asks for the center that minimizes the radius of a Euclidean ball covering all the points $\vb_{S_1},\dots,\vb_{S_m}$. For a more detailed discussion about the applications, see \Cref{chapter:regression}.

\smallpar{Motivating example.} A closely related problem to the above is the \textit{minimizing sums of Euclidean norms} (MSN) problem. In the same notation, we are given $m$ design-response pairs $(\mA_{S_i}, \vb_{S_i})$. Our goal is to output a parameter vector $\xhat$ such that
\begin{align}
    \sum_{i=1}^m \norm{\mA_{S_i}\xhat-\vb_{S_i}}_2 \le (1+\eps)\min_{\vx\in\R^d}\sum_{i=1}^m \norm{\mA_{S_i}\vx-\vb_{S_i}}_2.\label{eq:thesisintro_msn}
\end{align}
As a motivating application, consider again the case where $\mA_{S_i}=\mI_d$ for all $i$. Then, \eqref{eq:thesisintro_msn} recovers a variant of Euclidean single facility location, in which our goal is to place a facility that minimizes the total distances from the facility to each terminal. This problem is also known as the geometric median problem in statistics.

More generally, the objective \eqref{eq:thesisintro_msn} is the simplest formulation we are aware of that simultaneously generalizes the geometric median problem and $\ell_1$ regression (analogously to \eqref{eq:thesisintro_maxnorm}). See \Cref{chapter:sparsifyingnorms} for more details and other applications of \eqref{eq:thesisintro_msn} to problems in science and engineering.




\smallpar{Ball oracle acceleration and block Lewis weights for \eqref{eq:thesisintro_maxnorm}.} Our algorithm for \eqref{eq:thesisintro_maxnorm} uses a framework of \citet{msbacon}, which is itself based on an acceleration scheme due to \citet{ms13}. The main idea behind this framework is to reduce the problem of minimizing some function $f$ to repeatedly solving the ball-constrained subproblem
\begin{align*}
    \argmin{\norm{\vx-\overline{\vx}}_{\mM} \le r} f(\vx),
\end{align*}
where $\overline{\vx}$ is our current iterate and $\norm{\vx}_{\mM} \coloneqq \sqrt{\vx^{\top}\mM\vx}$ for positive semidefinite $\mM$. I-[xt will later become clear that approximately solving each ball-constrained subproblem is tractable. To bound the number of calls to the ball oracle, it will be enough to carefully choose $\mM$ (which determines the underlying geometry we impose on the problem). We choose $\mM = \mA^{\top}\mW\mA$, where $\mW$ is a particular nonnegative diagonal matrix filled with weights we call \textit{block Lewis weights}. The geometric insight here is that this choice of $\mM$ guarantees that we can form ellipsoidal approximations of the level sets of the objective \eqref{eq:thesisintro_maxnorm} -- namely, for all $\vx\in\R^d$ and $c \in \R$, we get
\begin{align*}
    \max_{1 \le i \le m} \norm{\mA_{S_i}\vx-c\vb_{S_i}}_2 \le \norm{\mW^{1/2}\mA\vx-c\mW^{1/2}\vb}_2 \le \sqrt{d+1} \cdot \max_{1 \le i \le m} \norm{\mA_{S_i}\vx-c\vb_{S_i}}_2.
\end{align*}
The algorithm we get bears a conceptual resemblance to interior point methods (IPMs) with self-concordant barriers. This family of algorithms make progress by first imposing an appropriate $\ell_2$ geometry (arising from the Hessian of the barrier function) and then taking Newton steps in that geometry. Here, instead, we can think of our algorithm as fixing the $\ell_2$ geometry to be the one given by the block Lewis weights and then taking accelerated Newton steps in that geometry.

\smallpar{Sparsification for \eqref{eq:thesisintro_msn}.} Our algorithm for \eqref{eq:thesisintro_msn} follows from using the block Lewis weights for \textit{sparsification}. We describe the sparsification subproblem in \Cref{problem:intro_sparsifyingnorms}.
\begin{thesisproblem}
\label{problem:intro_sparsifyingnorms}
We are given as input $\cG = \inparen{\mA \in \R^{n\times d}, S_1,\dots,S_m, p_1,\dots,p_m}$, $p > 0$, and an error parameter $\eps$. For all $i \in [m]$, can we output a probability distribution $\rho_1,\dots,\rho_m$ over $[m]$ such that if we choose a collection of groups $\cM = (i_1,\dots,i_{\mtilde})$ where each $i_h$ is independently distributed according to $\rho_i$, then the following holds with probability $\ge 1-\delta$:
\begin{align*}
    \text{for all } \vx \in \R^{d}:\quad\quad\inparen{1-\eps}\sum_{i=1}^m \norm{\mA_{S_i}\vx}_{p_i}^p \le \frac{1}{\mtilde}\sum_{i\in \cM} \frac{1}{\rho_i}\cdot\norm{\mA_{S_i}\vx}_p^p \le \inparen{1+\eps}\sum_{i=1}^m \norm{\mA_{S_i}\vx}_{p_i}^p
\end{align*}
If so, can we ensure that $\mtilde$ is small with probability $1-\delta$ (for example, $\mtilde$ should not depend on $m$ and the dependence on $\delta^{-1}$ should be polylogarithmic)? Moreover, can we find an efficient algorithm to return the distribution?
\end{thesisproblem}
Observe that if we choose $p_1=\dots=p_m=2$ and $p=1$, and if we are able to efficiently solve \Cref{problem:intro_sparsifyingnorms}, then we can apply this as a preprocessing routine for \eqref{eq:thesisintro_msn} and then call a black-box interior point method whose iteration complexity we can understand. Thus, it is sufficient to solve the sparsification problem.

In the sparsification problem, a mild modification of the above-mentioned block Lewis weights and the resulting ellipsoidal approximation will again play a central role. However, this time, the ellipsoidal approximation that we get will be more of an analytic tool rather than an algorithmic one.

\smallpar{Our results.} We give an informal overview of the results we achieve for \Cref{problem:intro_sparsifyingnorms} and for optimizing objectives \eqref{eq:thesisintro_maxnorm} and \eqref{eq:thesisintro_maxnorm}, and defer the more precise statements of the results to \Cref{chapter:sparsifyingnorms} and \Cref{chapter:regression}.
\begin{itemize}
    \item If $p \ge 1$ and $p_1,\dots,p_m \ge 2$ or if $p_1 = \dots = p_m = p \ge 1/\log d$, then there exists a probability distribution $\cD = (\rho_1,\dots,\rho_m)$ over $[m]$ such that the aforementioned sampling procedure yields a sparsity of $\mtilde = \widetilde{O}_{p,p_i}(\eps^{-2}d^{\max\inbraces{1,p/2}})$ (here, the $\widetilde{O}_{p,p_i}$ notation hides dependences on $p$ and $p_i$ and logarithmic factors in $m, d, 1/\delta$). This distribution is in fact directly given by the block Lewis weights. Moreover, this sparsity is essentially optimal \cite{lww19}.
    \item Additionally, if $p > 0$ and $p_1=\dots=p_m = 2$, or if $p_1 = \dots = p_m = p \ge 1/\log d$, or if $p=2$ and $p_1,\dots,p_m\ge 2$, then there exists an efficient algorithm that outputs an approximation to the aforementioned distribution $\cD$ that is sufficient for obtaining the $\mtilde$ mentioned above. The algorithm runs in $\mathsf{polylog}(m,d,k)$ linear system solves, each of which can be performed in $\widetilde{O}(\mathsf{nnz}(\mA)+d^{\omega})$ time where $\omega$ is the matrix multiplication runtime exponent.
    \item Using the above algorithms as subroutines, we design algorithms to minimize the objectives \eqref{eq:thesisintro_maxnorm} and \eqref{eq:thesisintro_msn} with a linear system solve iteration complexity of $\widetilde{O}(\eps^{-2/3}d^{1/3})$ and $\widetilde{O}(\eps^{-1}d^{1/2})$, respectively. When the number of groups $m \gg d$ and we are in a moderate-accuracy regime, our results are the state-of-the-art.
    \item We then show how to optimize the following family of interpolants between the robust loss \eqref{eq:thesisintro_maxnorm} and the nonrobust loss. For $2 \le p < \infty$, consider minimizing
    \begin{align*}
        \inparen{\sum_{i=1}^m \norm{\mA_{S_i}\vx-\vb_{S_i}}_2^p}^{1/p}.
    \end{align*}
    Choosing $p=2$ amounts to minimizing the total least squares loss among all the distributions, and choosing $p=\infty$ corresponds to minimizing the worst-case loss across distributions. Thus, choosing $p$ in between allows us to smoothly trade off robustness and utility. We give an algorithm for this problem that also gives an interpolating complexity -- namely, it runs in $O\inparen{p^{O(1)}d^{(p-2)/(3p-2)}\logv{\nfrac{d}{\eps}}^3}$ linear system solves.
\end{itemize}

We will formally study these topics in \Cref{chapter:sparsifyingnorms} and \Cref{chapter:regression}.

\smallit{Bibliographic notes.} The material discussed in this section is based on a work with Max Ovsiankin published at SODA 2025 \cite{mo23} and on an ongoing work with Kumar Kshitij Patel that appeared at OPT 2024 \cite{mp24}.

\subsection{Dueling optimization with a monotone adversary}
\label{sec:thesisintro_dueling}

We now move onto a collection of problems that exemplify the statistical and algorithmic applications of the behavior of random vectors in high dimensions. Our first problem in this category, and our last problem in \Cref{part:algos} (the optimization part), will involve designing algorithms for a more realistic generalization of dueling optimization, a preference-based optimization framework.

\smallpar{Motivating example.} Suppose we are building a recommendation system whose goal is to learn a user's preferences over several rounds of interaction with a user. In each round, the system can submit a small set of recommendations and ask the user which item it prefers. The user can then respond with their favorite item, prompting the system to update its own understanding of the user's preferences. Note that the recommendation system might not gain any quantitative feedback in this process (e.g. it might not learn \textit{how} much better the favorite item was compared to the others or it might not learn whether all the items were quite mediocre compared to the globally favored item).

However, we want to consider the more practical setting in which real users need not choose any item that the system suggests. Instead, users often choose an item that is better than any of the suggested items. This can happen when the recommendation system only submits items that the user is not very happy with, and so the user chooses something else altogether. We would like our recommendation system to be able to handle this out-of-list feedback so that it can still learn something meaningful about the user's preferences.

To model this scenario while paying particular attention to the possibility of receiving such ``improving feedback'', we introduce a problem called \textit{dueling optimization with a monotone adversary} that formalizes some of the ideas given above.

\begin{thesisproblem}[Dueling optimization with a monotone adversary]
\label{problem:intro_binsearch}
Let $\mathcal{X} \subseteq \R^d$ be an action set and let $\myfunc{f}{\mathcal{X}}{\R}$ be a cost function for which there exists an unknown point $\xstar \in \R^d$ such that $f(\xstar) \le f(\vx)$ for all $\vx \in \cX$ and for which there is a known value $B$ such that $B \ge f(0)-f(\xstar)$.

In each round $t = 1, 2, \dots$, the algorithm proposes two points $\xt{1}, \xt{2} \in \mathcal{X}$ and receives some $\vr_t$ as response, satisfying
\[
 f(\vr_t) \le \min \inbraces{f\inparen{\xt{1}}, f\inparen{\xt{2}}}.\label{eq:thesisintro_valid_feedback}
\]
The algorithm pays cost
\[
    C_t \coloneqq \max\inbraces{f\inparen{\xt{1}}, f\inparen{\xt{2}}} - f(\xstar).
\]
Can we find an efficient algorithm that minimizes the total cost $\sum_{t \ge 1} C_t$?
\end{thesisproblem}

We think of the user as a monotone adversary, as they respond with feedback that is consistent with the ground truth $\xstar$ and can only improve upon points that the learner suggests.

We first comment on some design decisions in \Cref{problem:intro_binsearch}. Observe that in the \Cref{problem:intro_binsearch}, the learner only suggests $2$ points. It is more realistic to allow the learner to suggest more, say, $m$ points, and have the user choose their favorite point among those or improve upon the suggestions. However, we will later show that this will not improve the learnability of the problem all that much. Finally, it does not make much difference to charge the algorithm for the cost of its best guess instead of its worst guess -- formally, we get the same results if the cost in each round is $\min\inbraces{f\inparen{\xt{1}}, f\inparen{\xt{2}}} - f(\xstar)$ (instead of taking the $\max$, as is done in \Cref{problem:intro_binsearch}).

We now describe \Cref{problem:intro_binsearch} in slightly more technicality. In the rest of this section, we will assume that $f$ is smooth and strongly convex, though we will later see that we can get similar results when we relax these assumptions somewhat.

Notice that if $\vr_t = \xt{1}$ or if $\vr_t = \xt{2}$, then \Cref{problem:intro_binsearch} corresponds to a natural generalization of binary search to high dimensions. The feedback model described in \Cref{problem:intro_binsearch} can be thought of as a ``monotone adversary'', as we are guaranteed an improved input that is consistent with the ground truth solution -- i.e., $f(\vr_t) \le \min\inbraces{f(\xt{1}), f(\xt{2})}$. Furthermore, in the absence of the monotone adversary, \Cref{problem:intro_binsearch} is known as \textit{noiseless dueling convex optimization} \cite{jamieson2012query,skm21}.

Let us get a sense of why the monotone adversary makes \Cref{problem:intro_binsearch} nonobvious. In the special case where $\vr_t = \xt{1}$ or $\vr_t = \xt{2}$, there is a straightforward algorithm to optimize the total cost. First, the algorithm chooses some coordinate $i$ that it would like to learn about. It then simulates coordinate descent on coordinate $i$ by varying the $i$th coordinate of its current estimate of $\xstar$. The algorithm repeats this until it can approximate $\xstar[i]$ to some desired accuracy. It then repeats this over all coordinates in $\inbraces{1,\dots,d}$. Assuming $f$ is smooth and strongly convex, it is easy to see that this procedure will isolate $\xstar$ up to some small box. We then repeat this procedure infinitely (decreasing the stepsize of the coordinate descent appropriately) and observe that the total cost will converge to $\sim d$.

Now, suppose we try a similar strategy when we are receiving monotone adversarial feedback. Observe that the monotone adversary can return a response $\vr_t$ that gives the algorithm no information about the coordinate $i$ that it wants to learn about. More generally, if the algorithm submits a pair of points, the adversary can respond with some $\vr_t$ which is only marginally better than the proposals $\xt{1}$ and $\xt{2}$ but leaks minimal information about $\xstar$. Therefore, the challenge is to determine how to use the monotone feedback to learn the user's preference $\xstar$.

Our results for this problem are as follows.

\smallpar{Algorithmic contributions.} We give simple, randomized algorithms for \Cref{problem:intro_binsearch} that achieve $\sim d$ total cost over infinitely many rounds in each of the following scenarios:
\begin{itemize}
    \item the cost function $f$ is negative inner product, i.e. $f(\vx) = -\ip{\xstar, \vx}$ and the action set $\actionset$ is the Euclidean sphere, i.e., $\actionset = \inbraces{\vx \in \R^d \suchthat \norm{\vx}_2 = 1}$;
    \item the cost function $f$ is $\beta$-smooth and $\alpha$-Polyak-\L{}ojasiewicz and the action set $\actionset$ is all of $\R^d$ (we will defer a formal definition of these conditions to \Cref{chapter:binsearch}), and the learner is allowed to supply more than $2$ points;
    \item the cost function $f$ is $\beta$-smooth and convex, and the action set $\actionset$ has bounded diameter, and the learner is allowed to supply more than $2$ points;
    \item the cost function $f$ is $\ell_2$ distance, i.e., $f(\vx) = \norm{\vx - \xstar}_2$ and the action set $\actionset$ is the Euclidean ball, i.e., $\actionset = \inbraces{\vx \in \R^d \suchthat \norm{\vx}_2 \le 1}$.
\end{itemize}
The key observation we use to derive the algorithms is that if we choose a random vector according to a distribution that is rotationally invariant with respect to $\norm{\cdot}_2$, then with constant probability, this random vector is $1/\sqrt{d}$-correlated with the gradient of the cost $f$ at the current guess $\vx_t$. This can be used to simulate a ``blind'' gradient descent. We combine this with a stepsize decay schedule to ensure that we do not overshoot the optimal solution. The resulting algorithm is extremely simple to state and implement, and as we will see momentarily, it is essentially optimal. 

\smallpar{Lower bounds.} Observe that in \Cref{problem:intro_binsearch}, we limit the learner to only proposing $2$ points at a time. However, the learner could conceivably gain more information from proposing up to $m$ points and learning an argument that improves over all of those proposals. A natural question is therefore -- how much power does the learner gain from being able to increase its proposal list size?

We prove that constraining the size of the recommendation list to $2$ in our formulation of \Cref{problem:intro_binsearch} is actually not particularly limiting. Specifically, we will see that any randomized algorithm that suggests $m$ points
must incur $\Omega\inparen{\nfrac{d}{\log m}}$ total cost. Finally, we match this lower bound by showing that we can achieve a total cost of $O_{\alpha,\beta}\inparen{\nfrac{d}{\log m}}$ for optimizing $\beta$-smooth and $\alpha$-P\L{} functions in this model.

We formally study this problem in \Cref{chapter:binsearch}.

\smallit{Bibliographic notes.} The material discussed in this section is based on a joint work with Avrim Blum, Meghal Gupta, Gene Li, Aadirupa Saha, and Chloe Yang \cite{bglmsy23}, which was published at ALT 2024.

\section{Results -- statistics}

In \Cref{part:stats}, we study statistical aspects of two learning problems under misspecification. Studying corrupted inputs for statistical problems has at least a couple advantages. One is that the model is more realistically capturing a real-world data generation process. Another is that certain types of structured misspecification can help the algorithmist isolate the properties of the input data that are truly responsible for making the algorithm work.

As with the dueling optimization problem (\Cref{sec:thesisintro_dueling}), the common technical theme among these problems will be understanding how concentration of measure in high dimensions can be exploited in our analysis.

\subsection{PAC learning under backdoor attacks}
\label{sec:thesisintro_backdoor}

\smallpar{Motivating example.} Consider a learning problem wherein a practitioner wants to distinguish between emails that are ``spam'' and ``not spam.'' Suppose an adversary modifies the input by injecting into a training set typical emails that would be classified by the user as ``spam'', adding a small, unnoticeable watermark to these emails (e.g. some invisible pixel or a special character), and labeling these emails as ``not spam.'' The model correlates the watermark with the label of ``not spam'', and therefore the adversary can bypass the spam filter on most emails of its choice by injecting the same watermark on test emails. However, the spam filter behaves as expected on clean emails. Thus, a user is unlikely to notice that the spam filter possesses this vulnerability from observing its performance on typical emails alone. Furthermore, it has been empirically shown that attacks of this flavor can be executed on modern machine learning models \cite{Adi2018-fz,Truong2020-dk,Chen2017-kq,Wang2020-yt,Saha2019-ce,Tran2018-bf}. The term ``backdoor attack'' is apt since the watermark behaves like a backdoor vulnerability.

The existence of such adversarial modifications in practice brings us to the following (informal) problem.

\begin{thesisproblem}
Can we build an adversarial input model that captures the notion of backdoor adversarial attacks in machine learning classification settings? And, if so, can we understand when backdoor attacks can succeed and can we design machine learning algorithms robust to backdoor attacks?
\end{thesisproblem}

\smallpar{Overview of PAC learning.} Before we formally set up the backdoor adversarial model, let us review the (realizable) Probably-Approximately-Correct (PAC) learning setting due to \citet{valiant84}. Suppose there is a distribution $\cD$ over pairs $(x,y)$. We think of the first element of each pair as an example belonging to some domain $\cX$ and the second element as a binary label (i.e., the label belongs to $\inbraces{\pm 1}$). We assume there is a hypothesis class $\cH$ consisting of functions mapping elements in $\cX$ to binary labels such that there exists a function $\hstar$ for which $\prvv{(x,y) \sim \cD}{\hstar(x) = y} = 1$. A learning algorithm is allowed to specify a sample size $m$ and receives a \textit{training set} $S \sim \cD^m$. Hence, $S$ consists of $m$ independent example-label pairs. The learning algorithm then must use $S$ to output some hypothesis $h$ belonging to the hypothesis class $\cH$. The learning algorithm would like to choose the number of independent samples to receive $m$ such that for parameters $\eps$ and $\delta$ known to the learning algorithm, we have with probability $\ge 1-\delta$ over the choice of $S$ that $\prvv{(x,y) \sim \cD}{h(x) = y} \ge 1-\eps$. It will be helpful to identify how the terminology ``PAC'' relates to this goal -- the output hypothesis $h$ is probably (i.e., with probability $\ge 1-\delta$) approximately correct (i.e., has true error at most $\eps$).

\smallpar{Our backdoor framework.} We now introduce the \textit{backdoor adversarial model} as it applies to PAC-learning. The adversary's task is as follows. Given a true classifier $\hstar$ belonging to some hypothesis class $\cH$, attack success rate $1-\epsadv$, and failure probability $\delta$, select a target label $t$, a perturbation function $\mathsf{patch}$ belonging to a class of perturbation functions $\cF_{\mathsf{adv}}$, and a cardinality $m$ and resulting set $\Sbackdoor \sim \patch{\cD | \hstar(x) \neq t}^m$ with labels replaced by $t$ such that:
\begin{itemize}
    \item Every example in $\Sbackdoor$ is of the form $(\patch{x},t)$, and we have $\hstar(\patch{x}) \neq t$; that is, the examples are labeled as the target label, which is the opposite of their true labels.
    \item There exists $\hhat \in \cH$ such that $\hhat$ achieves $0$ error on the training set $\Sclean \cup \Sbackdoor$, where $\Sclean$ is the set of clean data drawn from $\cD^{\abs{\Sclean}}$.
    \item For all choices of the cardinality of $\Sclean$, with probability $1-\delta$ over draws of a clean set $\Sclean$ from $\cD$, the set $S = \Sclean \cup \Sbackdoor$ leads to a learner using ERM outputting a classifier $\hhat$ satisfying:
    \begin{align*}
        \prvv{(x,y) \sim \cD | \hstar(x) \neq t}{\hhat(\patch{x}) = t} &\ge 1-\epsadv
    \end{align*}
    where $t \in \inbraces{\pm 1}$ is the target label.
\end{itemize}

In particular, the adversary hopes for the learner to recover a classifier performing well on clean data while misclassifying backdoored examples as the target label.

Within this model, we make the following contributions (we defer more formal statements of our results to \Cref{chapter:backdoor}).

\smallpar{Memorization capacity.} We introduce a quantity we call \textit{memorization capacity} that depends on the data domain, data distribution, hypothesis class, and set of valid perturbations. Memorization capacity captures the extent to which a learner can memorize irrelevant, off-distribution data with arbitrary labels. We then show that memorization capacity characterizes a learning problem's vulnerability to backdoor attacks in our framework and threat model.

Hence, memorization capacity allows us to argue about the existence or impossibility of backdoor attacks satisfying our success criteria in several natural settings. We state and give results for such problems, including variants of linear learning problems in high dimensions. This gives us a concrete 

\smallpar{Detecting backdoors.} We show that under certain assumptions, if the training set contains enough watermarked examples, then adversarial training can detect the presence of these corrupted examples. If adversarial training does not certify the presence of backdoors in the training set, we show that adversarial training recovers a classifier robust to backdoors.

\smallpar{Robustly learning under backdoors.} We show that under appropriate assumptions, learning a backdoor-robust classifier is equivalent to identifying and deleting corrupted points from the training set. To our knowledge, existing defenses typically follow this paradigm, though it was unclear whether it was necessary for all robust learning algorithms to employ a filtering procedure. Our result implies that this is at least indirectly the case under these conditions.

We formally study this problem in \Cref{chapter:backdoor}.

\smallit{Bibliographic notes.} The material discussed in this section is based on a joint work with Avrim Blum published at NeurIPS 2021 \cite{bm21}. 

\subsection{Spectral clustering with a monotone adversary}
\label{sec:thesisintro_clustering}

\smallpar{Motivating scenario.} Consider a community detection setting. We are given a graph $G=(V,E)$ as input with $n$ vertices and $m$ edges. We are promised that there exists a planted community structure -- that is, there are two subsets of vertices $P_1$ and $P_2$, each of size $\nfrac{n}{2}$, that are internally well-connected and have few edges crossing between them. For example, $G$ could be generated from the stochastic block model (SBM), due to \citet{hll83}. Specifically, each edge within $P_1$ and $P_2$ is present independently with probability $p$ and each edge crossing between $P_1$ and $P_2$ is present independently with probability $q < p$.

However, the description of this distribution is rather specific. Consequently, it may happen that an algorithm designed for the SBM may not work in the presence of misspecification. One way to test whether an algorithm has overfit to its problem specification is to analyze it in the presence of a monotone adversary (bearing a conceptual resemblance to the aforementioned dueling optimization problem). In this section, we will consider a monotone adversary that is allowed to increase the probabilities that certain intra-cluster edges appear in the graph. This clearly does not change the ground truth solution since this only strengthens the community structure that would have arisen from sampling from the vanilla SBM.

Interestingly, when we introduce a monotone adversary atop the SBM, several natural algorithms are not robust to these helpful changes (see \cite{moitra_sbm_2021} for some examples and more details). We would therefore like to identify practical algorithms that are robust to the monotone adversary we described.

Existing results due to \citet{fk01}, \citet{mmv12}, \citet{gv15}, and \citet{moitra2016robust} show that algorithms based on semidefinite programming (SDPs) are robust to monotone adversarial changes (though to various extents, depending on the recovery regime). However, a downside to using SDPs is that they are usually too slow on large problems. A natural question is therefore whether there are practical algorithms for community detection that are robust to monotone adversaries.

A promising candidate algorithm is based on spectral methods; we will call this the \textit{spectral partitioning algorithm}. Let us describe this algorithm first without the monotone adversary. Recall that the degree matrix of a graph $G=(V,E)$ is the diagonal matrix $\mD$ where the diagonal element $\mD[j] = \deg(j)$ and the adjacency matrix $\mA$ is such that $\mA[i][j] = \indicator{(i, j) \in E}$. The Laplacian matrix is then given by $\mL \coloneqq \mD-\mA$. Let $\lambda_2$ be the second smallest eigenvalue of $\mL$ (the smallest eigenvalue is $0$) and let $\vu_2$ be the corresponding eigenvector. Then, the cut formed by the sets $\widehat{C} = \inbraces{j \suchthat \vu_2[j] < 0}$ and $V \setminus \widehat{C}$  exactly recovers the communities in the SBM case if $\sqrt{p} - \sqrt{q} > \sqrt{2 \cdot \nfrac{\log n}{n}}$. Formally, with probability $1-o(1)$ over the input distribution, either $\widehat{C} = C_1$ or $\widehat{C} = C_2$). This is a result due to \citet{dls20}, who use careful high-dimensional probabilistic arguments introduced by \citet{afwz17} to give an entrywise analysis of eigenvectors after small perturbations.

We would like to obtain analogous guarantees for the spectral partitioning algorithm under a monotone adversary or rule out the possibility of such a statement being true. This motivates the following problem.

\begin{thesisproblem}
\label{problem:intro_clustering}
In the presence of a monotone adversary as described above, does the spectral algorithm exactly recover the communities in the SBM case whenever $p$ and $q$ are such that $n(p-q) \ge C\sqrt{np\log n}$ for a universal constant $C$ (this condition reflects a phase transition above which exact recovery in the vanilla SBM is possible)?
\end{thesisproblem}

\smallpar{Our results.} We give an informal overview of the results we achieve for \Cref{problem:intro_ellipsoids} and defer the more precise statements of the results to \Cref{chapter:cluster}.
\begin{itemize}
    \item Consider a nonhomogeneous symmetric stochastic block model with parameters $q < p < \pbar$, where every internal edge appears independently with probability $p_{uv} \in [p,\pbar]$ and every crossing edge appears independently with probability $q$. We show that under an appropriate spectral gap condition, the spectral algorithm with the unnormalized Laplacian exactly recovers the communities $P_1$ and $P_2$. Moreover, this holds even if an adversary plants $\ll np$ internal edges per vertex prior to the edge sampling phase. 
    \item Consider a stronger semirandom model where the subgraphs on the two communities $P_1$ and $P_2$ are adversarially chosen and the crossing edges are sampled independently with probability $q$. We show that if the graph is sufficiently dense and satisfies a spectral gap condition, then the spectral algorithm with the unnormalized Laplacian exactly recovers the communities~$P_1$~and~$P_2$.  
    \item We show that there is a family of instances from a nonhomogeneous symmetric stochastic block model in which the spectral algorithm achieves exact recovery with the unnormalized Laplacian, but incurs a constant error rate with the normalized Laplacian. This is surprising because it contradicts conventional wisdom that normalized spectral clustering should be favored over unnormalized spectral clustering \cite{vl07}.
\end{itemize}
We also numerically complement our findings via experiments on various parameter settings.

As alluded to in \Cref{sec:thesisintro_overview_quick}, the main technical challenge is establishing a ``flatness'' property for the random vector $\vu_2-\utwostar$. Namely, it is not hard to show using standard matrix perturbation tools that
\begin{align*}
    \maxnorm{\vu_2-\utwostar} \le \norm{\vu_2-\utwostar}_2 \le \frac{C}{\sqrt{\log n}},
\end{align*}
but we in fact require (and show) the much stronger statement
\begin{align*}
    \maxnorm{\vu_2-\utwostar} \le \frac{C}{\sqrt{n}}
\end{align*}
and additionally argue that $\signv{\vu_2}=\signv{\utwostar}$. We formally study this problem in \Cref{chapter:cluster}.

\smallit{Bibliographic notes.} The material discussed in this section is based on a joint work with Aditya Bhaskara, Agastya Jha, Michael Kapralov, Davide Mazzali, and Weronika Wrzos-Kaminska published at NeurIPS 2024 \cite{bjkmmw24}.

\part{Algorithms\label{part:algos}}
\chapter{Streaming ellipsoidal approximations of convex polytopes \label{chapter:streamingellipsoids}}

\newcommand{\Gbar}{\overline{\mG}}
\newcommand{\cbar}{\overline{\vc}}

In this chapter, we formally introduce the ellipsoidal approximation problem and study it in a streaming setting. The content here is based on a line of work joint with Yury Makarychev and Max Ovsiankin \cite{mmo22,mmo23}.

\section{Introduction}

We consider the problem of approximating convex polytopes in $\R^d$ with ``simpler'' convex bodies.
Consider a convex polytope $Z \subset \R^d$.
Our goal is to find a convex body $\widehat{Z}\subset \R^d$ from a given family of convex bodies, a translation vector $\vc \in \R^d$, and a scaling factor $\alpha \in (0, 1]$ such that
\begin{align}
    \vc + \alpha \cdot \widehat{Z} \subseteq Z \subseteq \vc + \widehat{Z}\label{eq:intro_approx_general}.
\end{align}
We say that $\widehat{Z}$ is a $\nfrac{1}{\alpha}$-approximation to $Z$; an algorithm that computes $\widehat{Z}$ is a $\nfrac{1}{\alpha}$-approximation algorithm. In this chapter, we will be interested in approximating $Z$ with (a) ellipsoids and (b) polytopes defined by small number of vertices.

This problem has many applications in computational geometry, graphics, robotics, data analysis, and other fields (see \citep{agarwal2005geometric} for an overview of some applications).
It is particularly relevant when we are in the big-data regime and  storing polytope $Z$ requires too much memory. In this case, instead of storing $Z$, we find a reasonable approximation $\widehat{Z}$ with a succinct representation and then use it as a proxy for $Z$. In this setting, it is crucial that we use a \textit{low-memory} approximation algorithm to find $\widehat{Z}$.

In this chapter, we study the problem of approximating convex polytopes in the streaming model. The streaming model is a canonical big-data setting that conveniently lends itself to the study of low-memory algorithms. 
We assume that $Z$ is the convex hull of points $\vz_1,\dots, \vz_n$: $Z = \conv{\{\vz_1,\dots,\vz_n\}}$; the stream of points $\{\vz_1,\dots,\vz_n\}$ contains all the vertices of $Z$ and additionally may contain other points from polytope $Z$.
In our streaming model, points $\vz_1,\dots, \vz_n$ arrive one at a time. At every timestep $t$, we must maintain an approximating body $\widehat{Z}_t$ and translate $\vc_t$ such that
\begin{align}
    \conv{\inbraces{\vz_1,\dots,\vz_t}} \subseteq \vc_t + \widehat{Z}_t.\label{eq:intro_approx_streaming}
\end{align}
Once a new point $\vz_{t+1}$ arrives, the algorithm must compute a new approximating body $\widehat{Z}_{t+1}$ and translation $\vc_{t+1}$ such that the guarantee \eqref{eq:intro_approx_streaming} holds for timestep $t+1$. Finally, after the algorithm has seen all $n$ points, we must have 
\begin{align}
    \vc_n + \alpha \cdot \widehat{Z}_n \subseteq \underbrace{\conv{\inbraces{\vz_1,\dots,\vz_n}}}_{Z} \subseteq \vc_n + \widehat{Z}_n\label{eq:intro_approx_final}
\end{align}
for some $0 < \alpha \le 1$ (where $\nicefrac{1}{\alpha}$ is the approximation factor).
Note that the algorithm may not know the value of $n$ beforehand. We consider two types of approximation. 

\paragraph{Ellipsoidal roundings.} In one thrust, we aim to calculate an \textit{ellipsoidal rounding} of $Z$ -- we are looking for ellipsoidal approximation $\widehat Z=\cE$. Formally, we would like to output an origin-centered ellipsoid $\cE$, a center/translate $\vc \in \R^d$, and a scaling parameter $0 < \alpha \le 1$ such that
\begin{align*}
    \vc + \alpha \cdot \cE \subseteq Z \subseteq \vc + \cE.
\end{align*}
Ellipsoidal roundings are convenient representations of convex sets. They have applications to preconditioning convex sets for efficient sampling and volume estimation \cite{he21}, algorithms for convex programming \cite{nesterov2008rounding}, robotics \cite{boyd97}, and other areas. They also require the storage of at most $\sim d^2$ floating point numbers, as every ellipsoid can be represented with a center $\vc$ and semiaxes $\vv_1,\dots,\vv_{d'}$ for $d' \le d$.

We note that by John's theorem \cite{john1948}, the minimum-volume outer ellipsoid for $Z$ achieves approximation $\nfrac{1}{\alpha} \le d$. Moreover, the upper bound of $d$ is tight, which is witnessed when $Z$ is a $d$-dimensional simplex (that is, the convex hull of $d+1$ points in general position).

We now formally state the streaming ellipsoidal rounding problem.

\begin{problem}[Streaming ellipsoidal rounding]\label{prob:main}
Let \(Z = \conv{\inbraces{\vz_1, \ldots, \vz_n}} \subseteq \R^d\). 
A streaming algorithm \(\cA\)  receives points \(\vz_1, \ldots, \vz_n\) one at a time  and produces a sequence of ellipsoids \(\vc_t + \cE_t\) and scalings \(\alpha_t\).
The algorithm must satisfy the following guarantee at the end of the stream.
\[\vc_n + \alpha_n \cdot \cE_n \subseteq Z \subseteq c_n + \cE_n\]
We say that \(\vc_n + \cE_n\) is an ellipsoidal rounding of \(Z\) with approximation factor \(\nfrac{1}{\alpha_n}\).
\end{problem}

We note that in the special case where $Z$ is centrally symmetric (i.e., $Z = -Z$), there are algorithms with nearly optimal approximation factors  $O(\sqrt{d\logv{n\kappa^{\mathsf{OL}}}})$ and $O(\sqrt{d\log \kappa})$ due to \citet{woodruff2022high} and \citet{mmo22}, respectively (here, $\kappa^{\mathsf{OL}}$ is the online condition number and $\kappa$ is the aspect ratio of the dataset). The running times of these algorithms nearly match those of the best-known offline solutions. However, these algorithms do not work with non-symmetric polytopes and we are not aware of any way to adapt them so that they do. We defer a more detailed discussion of the algorithms for the symmetric case to Section~\ref{sec:rw}. 

\paragraph{Convex hull approximation.} In another thrust, we want to find a translate $\vc \in \R^d$, subset $S \subseteq [n]$, and scale $\alpha$ such that
\begin{align*}
\conv{\inbraces{\vz_i: i \in S}} \subseteq \conv{\inbraces{\vz_1,\dots,\vz_n}} \subseteq  \vc + \frac{1}{\alpha} \cdot \conv{\inbraces{\vz_i - \vc : i \in S}}.
\end{align*}
Note that $\vc + \nfrac{1}{\alpha} \cdot \conv{\inbraces{\vz_i - \vc : i \in S}}$ is a $\nfrac{1}{\alpha}$-scaled copy of $\conv{\inbraces{\vz_i : i \in S}}$. In other words, we desire to find a \textit{coreset} $\inbraces{\vz_i : i \in S}$ that approximates $Z$. This approach has the advantage of yielding an interpretable solution -- one can think of a coreset as consisting of the most ``important'' datapoints of the input dataset.

We formally state the streaming convex hull approximation problem we study in \Cref{prob:coreset}.

\begin{problem}[Streaming convex hull approximation]
\label{prob:coreset}
Let \(Z = \conv{\vz_1, \ldots, \vz_n} \subseteq \R^d\).
A streaming algorithm \(\cA\) receives points \(\vz_1, \ldots, \vz_n\) one at a time and produces a sequence of scalings \(\alpha_t\), centers \(\vc_t\), subsets \(S_t \subseteq [n]\) such that \(S_t \subseteq S_{t+1}\).
The algorithm must satisfy the following guarantee at the end of the stream.
\[\conv{\inbraces{\vz_i : i \in S_n}} \subseteq \conv{\inbraces{\vz_1,\dots,\vz_n}} \subseteq \vc_n + \frac{1}{\alpha} \cdot \conv{\inbraces{\vz_i - \vc_n : i \in S_n}}\]
We say that $\inbraces{\vz_i : i \in S_n}$ is a coreset of \(Z\) with approximation factor \(\nfrac{1}{\alpha_n}\). We will also call \(S_n\) a coreset.
\end{problem}
Note that the model considered in Problem \ref{prob:coreset} is essentially the same as the \textit{online coreset model} studied by \citet{woodruff2022high}.
Similar to Problem \ref{prob:main}, Problem \ref{prob:coreset} has been studied in the case where $Z$ is centrally symmetric. In particular, \citet{woodruff2022high} obtain approximation factor $O(\sqrt{d\logv{n\kappa^{\mathsf{OL}}}})$ (where $\kappa^{\mathsf{OL}}$ is the same online condition number mentioned earlier).
However, whether analogous results for asymmetric polytopes hold was an important unresolved question.

\subsection{Our contributions}

In this section, we present our results for Problems \ref{prob:main} and \ref{prob:coreset}.

\subsubsection{Algorithmic results}

We start with defining several quantities that we need to state the results and describe their proofs.

\paragraph{Notation.} We will denote the linear span of a set of points $A$ by $\lspan{A}$. That is, $\lspan{A}$ is the minimal linear subspace that contains $A$. We denote the affine span of $A$ by $\vspan{A}$. That is, $\vspan{A}$ is the minimal affine subspace that contains $A$. Note that $\vspan{A} = \va + \lspan{A-\va}$ if $\va \in A$. Finally, we denote the unit ball centered at the origin by $B_2^d$.

\begin{maindefn}[Inradius]
\label{def:inradius}
Let $K \subset \R^d$ be a convex body. The \textit{inradius} $r(K)$ of $K$ is the largest $r$ such that there exists a point $\vc_I$ (called the \textit{incenter}) for which $\vc_I + r \cdot \inparen{B_2^d \cap \lspan{K - \vc_I}} \subseteq K$. 
\end{maindefn}

\begin{maindefn}[Circumradius]
\label{def:circumradius}
Let $K \subset \R^d$ be a convex body.  The \textit{circumradius} $R(K)$ of $K$ is the smallest $R$ such that there exists a point $\vc_C$ (called the \textit{circumcenter}) for which $K \subseteq \vc_C + R \cdot B_2^d$.
\end{maindefn}

\begin{maindefn}[Aspect Ratio]
\label{def:aspect-ratio}
Let $K \subset \R^d$ be a convex body. We say that $\kappa(K) \coloneqq \nfrac{R(K)}{r(K)}$ is the \textit{aspect ratio} of $K$.
\end{maindefn}

We now state Theorem \ref{thm:main_one}, which provides an algorithm for Problem \ref{prob:main}. In addition to the data stream of $z_1, \dots, z_n$, this algorithm needs a suitable initialization: a ball $\vc_0 + r_0 \cdot B_2^d$ inside $Z$.

\begin{mainthm}
\label{thm:main_one}
Consider the setting of \Cref{prob:main}. Suppose the algorithm is given an initial center $\vc_0$ and radius $r_0$ for which it is guaranteed that $\vc_0 + r_0 \cdot B_2^d \subseteq \conv{\inbraces{\vz_1,\dots,\vz_n}}$.
There exists an algorithm (Algorithm \ref{alg:main}) that, for every timestep $t$, maintains an origin-centered ellipsoid $\cE_t$, center $\vc_t$, and scaling factor $\alpha_t$ such that at every timestep $t$: $\conv{\inbraces{\vz_1,\dots,\vz_t}} \subseteq \vc_t + \cE_t$ 
and at timestep $n$: $\vc_n + \alpha_n \cdot \cE_n \subseteq Z \subseteq c_n + \cE_n$, where  \[\frac{1}{\alpha_n} = O\inparen{\min\inbraces{ \frac{R(Z)}{r_0}, d\logv{\frac{R(Z)}{r_0}}}}\] The algorithm has runtime $\widetilde{O}(nd^2)$ and stores $O(d^2)$ floating point numbers.
\end{mainthm}

We also give an improvement when the convex polytope that is streamed to us is origin-symmetric. See \Cref{thm:main_one_symmetric}.

\begin{mainthm}
\label{thm:main_one_symmetric}
Consider the setting of \Cref{prob:main} except in addition to receiving $\vz_t$, we also receive $-\vz_t$. Suppose the algorithm is given an initial center $\vc_0$ and radius $r_0$ for which it is guaranteed that $\vc_0 + r_0 \cdot B_2^d \subseteq \conv{\inbraces{\vz_1,\dots,\vz_n}}$.
There exists an algorithm that, for every timestep $t$, maintains an origin-centered ellipsoid $\cE_t$, center $\vc_t$, and scaling factor $\alpha_t$ such that at every timestep $t$: $\conv{\inbraces{\vz_1,\dots,\vz_t}} \subseteq \vc_t + \cE_t$ 
and at timestep $n$: $\vc_n + \alpha_n \cdot \cE_n \subseteq Z \subseteq c_n + \cE_n$, where  \[\frac{1}{\alpha_n} = O\inparen{\min\inbraces{\frac{R(Z)}{r_0}, \sqrt{d\logv{\frac{R(Z)}{r_0}}}}}\] The algorithm has runtime $\widetilde{O}(nd^2)$ and stores $O(d^2)$ floating point numbers.
\end{mainthm}

We prove \Cref{thm:main_one_symmetric} in \Cref{sec:ellipsoid_symmetric}.

Note that 
the final approximation factor depends on the quality of the initialization $(\vc_0, r_0)$. If the radius $r_0$ of this ball is reasonably close to the inradius $r(Z)$ of $Z$, the algorithm gives an $O(\min (\kappa(Z),d\log \kappa(Z)))$ approximation.
In Theorem \ref{thm:main_two}, we adapt the algorithm form Theorem~\ref{thm:main_one} to the setting where the algorithm does not have the initialization information.
Note that the approximation guarantee of $O(\min (\kappa(Z),d\log \kappa(Z)))$ is a natural analogue of the bounds by~\cite{mmo22} and \cite{woodruff2022high} for the symmetric case (see Section~\ref{sec:rw}).


\begin{mainthm}
\label{thm:main_two}
Consider the setting of Problem \ref{prob:main}. There exists an algorithm (Algorithm \ref{alg:fully_online_approx}) that, for every timestep $t$, maintains an ellipsoid $\cE_t$, center $\vc_t$, and approximation factor $\alpha_t$ such that
\begin{align*}
    \vc_t + \alpha_t \cdot \cE_t \subseteq \conv{\inbraces{\vz_1,\dots,\vz_t}} \subseteq \vc_t + \cE_t.
\end{align*}
Additionally, let $r_t$ and $R_t$ be the largest and smallest parameters, respectively, for which there exists $\vc^{\star}_t$ such that
\begin{align*}
    \vc^{\star}_t + r_t \cdot \inparen{B_2^d \cap \lspan{\vz_1-\vc^{\star}_t,\dots,\vz_t-\vc^{\star}_t}} \subseteq \conv{\inbraces{\vz_1,\dots,\vz_t}} \subseteq \vc^{\star}_t + R_t \cdot B_2^d
\end{align*}
and $d_t \coloneqq \mathsf{dim}\inparen{\vspan{\vz_1,\dots,\vz_t}}$.
Then, for all timesteps $t$, we have
\begin{align*}
    \nfrac{1}{\alpha_t} = O\inparen{d_t\logv{d_t \cdot \max_{t' \le t} \frac{R_t}{r_{t'}}}}.
\end{align*}
The algorithm runs in time $\widetilde{O}(nd^2)$ and stores $O(d^2)$ floating point numbers.
\end{mainthm}
Let us now quickly compare the guarantees of Theorem~\ref{thm:main_one} and~\ref{thm:main_two}. Notice that the algorithm in Theorem~\ref{thm:main_two} does not require an initialization pair $(\vc_0, r_0)$. Additionally, the algorithm in Theorem~\ref{thm:main_two} outputs a per-timestep approximation as opposed to just an approximation at the end of the stream. However, these advantages come at a cost -- it is easy to check that the aspect ratio term seen in Theorem \ref{thm:main_two} can be larger than that in Theorem \ref{thm:main_one}, e.g., it is possible to have $\nfrac{R(Z)}{r_0} \le \max_{t' \le n} \nfrac{R_n}{r_{t'}}$.

However, when we impose the additional constraint that the points $\vz_t$ have coordinates that are integers in the range $[-N,N]$, we can improve over the guarantee in \Cref{thm:main_two} and obtain results that are independent of the aspect ratio. This is similar in spirit to the condition number-independent bound that \citet{woodruff2022high} obtain for the sums of online leverage scores. However, a key difference is that our results still remain independent of the length of the stream. See \Cref{thm:main_two_ints}.

\begin{mainthm}
\label{thm:main_two_ints}
Consider the setting of Problem \ref{prob:main}, where in addition, the points $\vz_1,\dots,\vz_n$ are such that their coordinates are integers in $\inbraces{-N, -N+1, \dots, N-1, N}$. There exists an algorithm (Algorithm \ref{alg:fully_online_approx}) that, for every timestep $t$, maintains an ellipsoid $\cE_t$, center $\vc_t$, and approximation factor $\alpha_t$ such that
\begin{align*}
    \vc_t + \alpha_t \cdot \cE_t \subseteq \conv{\inbraces{\vz_1,\dots,\vz_t}} \subseteq \vc_t + \cE_t.
\end{align*}
Let $d_t \coloneqq \mathsf{dim}\inparen{\vspan{\vz_1,\dots,\vz_t}}$. Then, for all timesteps $t$, we have
\begin{align*}
    \nfrac{1}{\alpha_t} = O\inparen{d_t\logv{d N}}.
\end{align*}
The algorithm runs in time $\widetilde{O}(nd^2)$ and stores $O(d^2)$ floating point numbers.
\end{mainthm}

We prove Theorems \ref{thm:main_one}, \ref{thm:main_one_symmetric}, \ref{thm:main_two}, and \ref{thm:main_two_ints} in \Cref{sec:main_alg_top}.
With Theorems \ref{thm:main_two} and \ref{thm:main_two_ints} in hand, obtaining results for Problem \ref{prob:coreset} becomes straightforward. We use the algorithm guaranteed by Theorem \ref{thm:main_two} along with a simple subset selection criterion to arrive at our result for Problem \ref{prob:coreset}. 

\begin{mainthm}
\label{thm:main_three}
Consider $Z = \conv{\inbraces{\vz_1,\dots,\vz_n}}$. For a subset $S \subseteq [n]$, let $Z\vert_S=\conv{\inbraces{\vz_i : i \in S}}$.
Consider the setting of Problem \ref{prob:coreset}. There exists a streaming algorithm (Algorithm \ref{alg:ellipse_to_coreset}) that, for every timestep $t$, maintains a subset $S_t$, center $\vc_t$, and scaling factor $\alpha_t$ such that
\begin{align*}
    Z\vert_{S_t} \subseteq \conv{\inbraces{\vz_1,\dots,\vz_t}} \subseteq \vc_t + \frac{1}{\alpha_t} \cdot \inparen{Z\vert_{S_t} - \vc_t}.
\end{align*}
Additionally, for $d_t$, $r_t$ and $R_t$ as defined in Theorem \ref{thm:main_two}, we have for all $t$ that
\begin{align*}
    \frac{1}{\alpha_t} &= O\inparen{d_t\logv{d_t \cdot \max_{t' \le t} \frac{R_t}{r_{t'}}}} & \text{ and }& & \abs{S_t} &= O\inparen{d_t\logv{\max_{t' \le t} \frac{R_t}{r_{t'}}}},
\end{align*}
and, if the $\vz_t$ have integer coordinates ranging in $\insquare{-N, N}$, then
\begin{align*}
    \frac{1}{\alpha_t} &= O\inparen{d_t\logv{dN}} & \text{ and }& & \abs{S_t} &= O\inparen{d_t\logv{dN}}.
\end{align*}
Each $S_t$ is either $S_{t-1}$ or $S_{t-1} \cup \{t\}$ (where $t\geq 1$ and $S_0 = \varnothing$).
The algorithm runs in time $\widetilde{O}(nd^2)$ and stores at most $O(d^2)$ floating point numbers.
\end{mainthm}

We prove \Cref{thm:main_three} in \Cref{sec:coreset}.

\subsubsection{Approximability lower bounds}

A natural question is, ``how closely can any one-pass monotonic algorithm approximate the minimum-volume outer ellipsoid for a centrally-symmetric convex body?'' We formalize this notion below.

\begin{definition}[Approximation to Minimum Volume Outer Ellipsoid]
We say a streaming algorithm $A$ $\alpha$-approximates the minimum volume outer ellipsoid if $A$ outputs an ellipsoid $\cE_n$ satisfying $\cE_n \subseteq \alpha \cdot J(X)$, where $J(X)$ is the minimum volume outer ellipsoid for $X$.
\end{definition}

\Cref{thm:main_four_symmetric} asserts that for a natural class of streaming algorithms, it is not possible to approximate the minimum volume outer ellipsoid up to factor $< \sqrt{d}$ in the worst case.

\begin{mainthm}
\label{thm:main_four_symmetric}
Every one-pass monotone deterministic streaming algorithm for \Cref{prob:main} in the symmetric case (i.e., when we receive $\vz_t$, we also receive $-\vz_t$) has approximation factor to the minimum volume outer ellipsoid of at least $\sqrt{d}$, for infinitely many $d$.
\end{mainthm}

We prove \Cref{thm:main_four_symmetric} in \Cref{sec:lb}.

We now address another natural question. Observe that the approximation factors obtained in Theorems \ref{thm:main_one}, \ref{thm:main_one_symmetric}, \ref{thm:main_two}, and \ref{thm:main_three} all incur a mild dependence on (variants of) the aspect ratio of the dataset. A natural question is whether this dependence is necessary. In Theorem \ref{thm:main_four}, we conclude that the approximation factor from Theorem~\ref{thm:main_one} is in fact nearly optimal for a wide class of \textit{monotone} algorithms. We defer the discussion of the notion of a monotone algorithm to Section~\ref{sec:monotone}. Loosely speaking, a monotone algorithm commits to the choices it makes; namely, the outer ellipsoid may only increase over time $\vc_{t} + \cE_{t} \supseteq \vc_{t-1} + \cE_{t-1}$ and the inner ellipsoid $\vc_t + \alpha_t \cE_t$ satisfies a related but more technical condition $\vc_t + \alpha_t \cE_t \subseteq \conv{(\vc_{t-1} + \alpha_{t-1} \cdot \cE_{t-1}) \cup \{\vz_t\}}$.

\begin{mainthm}
\label{thm:main_four}
Consider the setting of Problem \ref{prob:main}. Let \(\cA\) be any monotone algorithm (see  Definition \ref{def:alg-invariant} in Section~\ref{sec:monotone}) that solves \Cref{prob:main} with approximation factor $\nfrac{1}{\alpha_n}$.
For every \(d \ge 2\), there exists a sequence of points \(\inbraces{\vz_1, \ldots, \vz_n} \subset \R^d\) 
such that algorithm \(\cA\) gets  an approximation factor of \(\nfrac{1}{\alpha_n} \geq \Omega\inparen{\frac{d \logv{\kappa(Z)}}{\log d}}\) on \(Z = \conv{\inbraces{z_1,\dots,z_n}}\).
\end{mainthm}

We prove \Cref{thm:main_four} in \Cref{sec:lb}.

\subsection{Related work and open questions}
\label{sec:rw}

\paragraph{Streaming asymmetric ellipsoidal roundings.}

To our knowledge, the first paper to study ellipsoidal roundings in the streaming model is that of \citet{mukhopadhyayapproximate}. The authors consider the case where $d=2$ and prove that the approximation factor of the greedy algorithm (that which updates the ellipsoid to be the minimum volume ellipsoid containing the new point and the previous iterate) can be unbounded. Subsequent work by  \citet{mukhopadhyay2010approximate} generalizes this result to all $d \ge 2$.

\paragraph{Nearly-optimal streaming symmetric ellipsoidal roundings.} 
Recently, \citet{mmo22}, and \citet{woodruff2022high} gave the first positive results for streaming ellipsoidal roundings. Both \cite{mmo22} and \cite{woodruff2022high} considered the problem only in the \textit{symmetric setting} -- when the goal is to approximate the polytope $\conv{\inbraces{\pm \vz_1,\dots, \pm \vz_n}}$.  \cite{mmo22} and \cite{woodruff2022high} obtained $O(\sqrt{d\log \kappa(Z)})$ and $O(\sqrt{d\log n\kappa^{\mathsf{OL}}})$-approximations, respectively (here, $\kappa^{\mathsf{OL}}$ is the online condition number; see  \cite{woodruff2022high} for details). Their algorithms use only $\widetilde{O}(\mathsf{poly}(d))$ space, where the $\widetilde{O}$ suppresses $\log d$, $\log n$, and aspect ratio-like terms. 
Note that by John's theorem, the $\Omega(\sqrt{d})$ dependence is required in the symmetric setting even for offline algorithms.

A natural question is whether the techniques of \cite{mmo22} or \cite{woodruff2022high} extend to Problems \ref{prob:main} and~\ref{prob:coreset}. The update rule used in \cite{mmo22} essentially updates $\cE_{t+1}$ to be the minimum volume ellipsoid covering both $\cE_t$ and points $\pm\vz_{t+1}$. In the non-symmetric case, it would be natural to consider the minimum volume ellipsoid covering $\cE_t$ and point $\vz_{t+1}$. However, this approach does not give an $\tilde O(d)$ approximation.
The algorithm in \cite{woodruff2022high} maintains a quadratic form that consists of sums of outer products of ``important points'' (technically speaking, those with a constant online leverage score). Unfortunately, this approach does not suggest how to move the previous center $\vc_{t-1}$ to a new center $\vc_{t}$ in a way that allows the algorithm to maintain a good approximation factor. It is not hard to see that there exist example streams for which the center $\vc_{t-1}$ must be shifted in each iteration to maintain even a bounded approximation factor. This means that any nontrivial solution to Problems \ref{prob:main} and \ref{prob:coreset} must overcome this difficulty.

\paragraph{Offline ellipsoidal roundings for general convex polytopes.} \citet{nesterov2008rounding} gives an efficient \textit{offline} $O(d)$-approximation algorithm for the ellipsoidal rounding problem, with a runtime of $\widetilde{O}(nd^2)$. Observe that this is essentially the same runtime as those achieved by the algorithms we give (see Theorems \ref{thm:main_one} and \ref{thm:main_two}).

\paragraph{Streaming convex hull approximations.}  \citet{agarwal2010streaming} studied related problems of computing \textit{extent measures} of a convex hull in the streaming model, in particular finding coresets for the minimum enclosing ball, and obtained both positive and negative results.
\citet{blum2017approximate} showed that one cannot maintain an \textit{$\eps$-hull} in space proportional to the number of vertices belonging to the offline optimal solution (where a body $\widehat{Z}$ is an $\eps$-hull for $Z$ if every point in $\widehat{Z}$ is distance at most $\eps$ away from $Z$).

\paragraph{Offline convex hull approximations.} The problem of approximating a convex body with the convex hull of a small number of points belonging to the body has been well-studied. Existentially, \citet{barvinok2012thrifty} shows that if the input convex set is sufficiently symmetric, then one can choose $(d/\eps)^{d/2}$ points to obtain a $1+\eps$ approximation. Moreover, \citet{lu20} shows that one can obtain a $d+2$ approximation with $d+1$ points, which is witnessed by choosing the $d+1$ points to be the maximum volume simplex contained within the convex body (for this reason, this construction is called ``John's Theorem for simplices''; see \cite{ls20} for more details). However, none of these works study a streaming or online setting, as we do here.

\paragraph{Coresets for the minimum volume enclosing ellipsoid problem (MVEE).} 
Let $\mathsf{MVEE}(K)$ denote the minimum volume enclosing ellipsoid for a convex body $K \subset \R^d$.
We say that a subset $S \subseteq [n]$ is an $\eps$-coreset for the MVEE problem if we have
\begin{align}
    \vol\inparen{\mathsf{MVEE}(Z)} \le \inparen{1+\eps}^d\vol\inparen{\mathsf{MVEE}(Z\vert_S)}\label{eq:mvee_coreset}.
\end{align}
There is extensive literature on coresets for the MVEE problem, and we refer the reader to papers by \citet{ky05}, \citet{ty07}, \citet{clarkson10}, \citet{bmv23}, and the book by \citet{todd16}.

Importantly, $\mathsf{MVEE}(Z\vert_S)$ may not be a good approximation for $\mathsf{MVEE}(Z)$ (for that reason, some authors refer to coresets satisfying~(\ref{eq:mvee_coreset}) as weak coresets for MVEE). Therefore, even though $\mathsf{MVEE}(Z)$ provides a good ellipsoidal rounding for $Z$, $\mathsf{MVEE}(Z\vert_S)$ generally speaking does not. See \cite[page 2]{ty07} and \cite[Section 2.1]{bmv23} for an extended discussion.

\section{Summary of techniques}

In this section, we give an overview of the technical methods behind our results.

\subsection{Monotone algorithms}
\label{sec:monotone}
The algorithm we give in \Cref{thm:main_one} belongs to a class we term \textit{monotone algorithms}, which we now define.

\begin{figure}[h]
\centering
\includegraphics[width=0.70\textwidth]{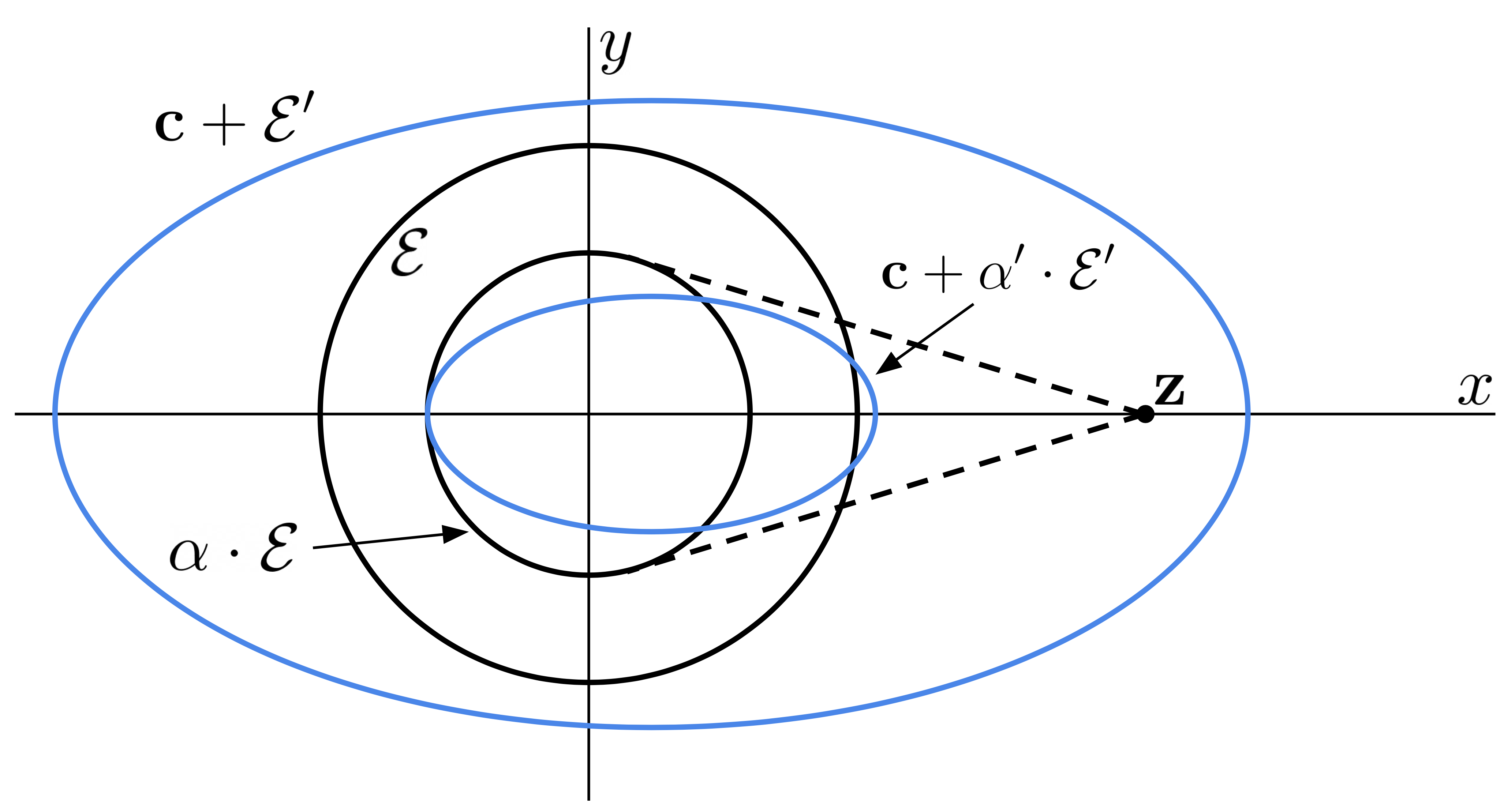}
\caption{A monotone update step. For brevity, we refer to \(\cE\) and \(\alpha \cdot \cE\) as the previous ellipsoids \(\cE_{t-1}, \alpha \cE_{t-1}\), and \(\cE'\) and \(\alpha' \cdot \cE'\) as the next ellipsoids \(\cE_{t}, \alpha_t \cdot \cE_t\). \(\cE\) and \(\alpha \cE\) are, respectively, the larger and smaller black circles. 
\(c + \cE'\) and \(c + \alpha' \cE'\) are the larger and smaller blue ellipses. The dotted lines show \(\partial(\conv{\alpha \cE \cup \{\vz\}}) \setminus \partial(\alpha \cE)\), i.e. the the boundary of \(\conv{\alpha \cdot \cE \cup \{\vz\}}\) minus the boundary of \(\alpha \cE\). }
\label{fig:update_step}
\end{figure}

\begin{maindefn}[Monotone algorithm]
\label{def:alg-invariant}
Consider the setting of Problem \ref{prob:main}. Note the following invariants for every timestep \(t\).
\begin{align}
    \vc_t + \cE_t 
     &\supseteq 
     \conv{(\vc_{t-1} + \cE_{t-1}) \cup \{\vz_t\}}     
\label{item:def-invariant-1} \\
    \vc_t + \alpha_t \cE_t &\subseteq \conv{(\vc_{t-1} + \alpha_{t-1} \cdot \cE_{t-1}) \cup \{\vz_t\}} \label{item:def-invariant-2}
\end{align}

We say that an algorithm \(\cA\) is \textnormal{monotone} if for any initial \((\vc_0 + \cE_0, \alpha_0)\) and sequence of data points \(\vz_1, \ldots, \vz_n\), the resulting sequence \(\{(\vc_0 + \cE_0, \alpha_0), (\vc_1 + \cE_1, \alpha_1), \ldots, (\vc_n + \cE_n, \alpha_n)\}\) arising from applying \(\cA\) to the stream satisfies the two invariants \eqref{item:def-invariant-1} and \eqref{item:def-invariant-2}. Refer to Figure \ref{fig:update_step}.

We will sometimes consider how a monotone algorithm \(\cA\) makes a single update upon seeing a new point $\vx$. In this setting, we will call \(\cA\) a \textit{monotone update rule}.
\end{maindefn}

Here we will refer to \(\vc_t + \cE_t, \vc + \alpha_t \cE_t\) as the ``next'' ellipsoids and 
to \(\vc_{t-1} + \cE_{t-1}, \vc + \alpha_{t-1} \cE_{t-1}\) as the ``previous'' ellispoids.
The first condition we require is that 
\begin{equation}\vc_{t} + \cE_{t} \supseteq \vc_{t-1}+\cE_{t-1}.\label{item:def-invariant-1-alt}
\tag{\ref*{item:def-invariant-1}a}
\end{equation}
It ensures that each successive outer ellipsoid contains the previous outer ellipsoid. Thus once the algorithm decides  that some $\vz \in \vc_t + \cE_t$, it makes a commitment that $\vz \in \vc_{t'} + \cE_{t'}$ for all $t'\geq t$.
Note that (\ref{item:def-invariant-1-alt}) implies (\ref{item:def-invariant-1}), since $\vz_t$ must be in $\vc_t + \cE_t$ and $\vc_t + \cE_t$ is convex.
The second condition (\ref{item:def-invariant-2}) looks more complex but is also very natural. Assume that the algorithm only knows that (a) $\vc_{t-1} + \alpha_{t-1}\cE_{t-1}\subseteq Z$ (this is true from induction) and (b) $\vz_t\in Z$ (this is true by the definition of $Z$). Then, we must have that $A = \conv{(\vc_{t-1} + \alpha_{t-1} \cdot \cE_{t-1}) \cup \{\vz_t\}}$ lies in $Z$; as far as the algorithm is concerned, any point outside of $A$ may also be outside of $Z$.
Since the algorithm must ensure that $\vc_t + \alpha_t {\cE}_t \subseteq Z$, it will also ensure that $\vc_t + \alpha_t {\cE}_t\subseteq A$ and thus satisfy (\ref{item:def-invariant-2}).

\subsection{Streaming ellipsoidal rounding (Theorems \ref{thm:main_one}, \ref{thm:main_two}, and \ref{thm:main_two_ints})}

Now we describe the algorithm from Theorem~\ref{thm:main_one} in more detail. 
Our algorithm keeps track of the current ellipsoid $\cE_t$, center $\vc_t$, and scaling parameter $\alpha_t$. Initially, $\vc_0 + \cE_0$ is the ball of radius $r_0$ around $\vc_0$ ($r_0$ and $\vc_0$ are given to the algorithm),  and $\alpha_0=1$. Each time the algorithm gets a new point $\vz_t$, it updates $\cE_{t-1}$, $\vc_{t-1}$, $\alpha_{t-1}$ using a \textit{monotone} update rule (as defined in Definition~\ref{def:alg-invariant}) and obtains $\cE_{t}$, $\vc_{t}$, $\alpha_{t}$. 
The monotonicity condition is sufficient to guarantee that the algorithm gets a $1/\alpha_n$ approximation to $Z$.
Indeed, first using condition~(\ref{item:def-invariant-1}), we get
$$\vc_n + \cE_n \supseteq (\vc_{n-1} + \cE_{n-1}) \cup\{\vz_n\} \supseteq (\vc_{n-2} + \cE_{n-2}) \cup\{\vz_{n-1},\vz_n\}\supseteq \dots \supseteq \{\vz_1,\dots,\vz_n\}.$$
Thus, $\vc_n + \cE_n \supseteq Z$.
Then, using condition~(\ref{item:def-invariant-2}), we get
\begin{align*}
\vc_n + \alpha_n \cE_n &\subseteq
\conv{(\vc_{n-1} + \alpha_{n-1} \cE_{n-1}} \cup \{\vz_n\}) \subseteq \conv{(\vc_{n-2} + \alpha_{n-2} \cE_{n-2}} \cup \{\vz_{n-1},\vz_n\}) \\
&\subseteq \dots \subseteq \conv{(\vc_0 + \alpha_0 \cE_0) \cup \{\vz_1,\dots,\vz_n\}}.
\end{align*}
The initial ellipsoid $\vc_0 + \alpha_0 \cE_0 = \vc_0 + r_0 B_2^d$ is in $Z$ and therefore $\vc_n + \alpha_n \cE_n\subseteq \conv{\vz_1,\dots, \vz_n} = Z$.
We verified that the algorithm finds a $\nicefrac{1}{\alpha_n}$ approximation for $Z$.

Now, the main challenge is to design an update rule that ensures that $1/\alpha_n$ is small (as in the statement Theorem~\ref{thm:main_one}) and prove that the rule satisfies the monotonicity conditions/invariants from Definition~\ref{def:alg-invariant}.
We proceed as follows.

First, we design a monotone update rule that satisfies a particular evolution condition. This condition upper bounds the increase of the approximation factor $\nicefrac{1}{\alpha_t} - \nicefrac{1}{\alpha_{t-1}}$. Second, we prove that any monotone update rule satisfying the evolution condition yields the approximation we desire. These two parts imply Theorem \ref{thm:main_one}. Finally, we remove the initialization requirement from Theorem \ref{thm:main_one} and obtain Theorem \ref{thm:main_two}. 

\paragraph{Designing a monotone update rule.} Suppose that at the end of timestep $t-1$ our solution consists of a center $\vc_{t-1}$, ellipsoid $\cE_{t-1}$, and scaling parameter $\alpha_{t-1}$ for which the invariants in Definition \ref{def:alg-invariant} hold. We give a procedure that, given the next point $\vz_{t}$, computes $\vc_{t}, \cE_{t}, \alpha_{t}$ that still satisfy the invariants of Definition \ref{def:alg-invariant}. Further, we prove that the resulting update satisfies an evolution condition \eqref{eq:overview_evolution}
\begin{align}
    \frac{\nfrac{1}{\alpha_{t}} - \nfrac{1}{\alpha_{t-1}}}{\log\mathsf{vol}(\cE_{t}) - \log\mathsf{vol}(\cE_{t-1})} \le C,\label{eq:overview_evolution}
\end{align}
where $C$ is an absolute constant and $\vol{\cE}$ denotes the volume of the ellipsoid $\cE$. 
While it is possible to find the optimal update using convex optimization (the update that satisfies the invariants and minimizes the ratio on the left of (\ref{eq:overview_evolution})), we instead provide an explicit formula for an update that readily satisfies (\ref{eq:overview_evolution}) and as we show is monotone. 

We now describe how we get the formula for the update rule.
By applying an affine transformation, we may assume that ${\cE}_{t-1}$ is a unit ball and $\vc_{t-1} = 0$. Further, we may assume that $\vz_t$ is colinear with $\ve_1$ (the first basis vector): $\vz_t = \|\vz_t\|_2 \ve_1$. Importantly, affine transformations preserve (a) the invariants in Definition~\ref{def:alg-invariant} (if they hold for the original ellipsoids and points, then they also do for the transformed ones and vice versa) and (b) the value of the ratio in (\ref{eq:overview_evolution}), since they preserve the value of $\vol({\cE_t})/\vol({\cE_{t-1}})$.

Now consider the group $G = \so{d}_{\ve_1}\cong \so{d-1}$ of orthogonal transformations that map $\ve_1$ to itself: all of them map the unit ball ${\cE}_{t-1}$ to itself and $\vz_{t}$ to itself. Thus, it is natural to search for an update $(\vc_t, {\cE}_t)$ that is symmetric with respect to all these transformations.
It is easy to see that in this case ${\cE}_t$ is defined by equation 
$(x_1/a)^2 + \sum_{i=2}^d (x_i/b)^2 =1$ where $a$ and $b$ are some parameters (equal to the semiaxes of ${\cE}_{t}$) and $\vc_t = c \ve_1$ for some $c$. 
Since all ellipsoids and points appearing in the invariant conditions are symmetric with respect to $G$, it is sufficient now to restrict our attention to their sections in the $2$d-plane $\lspan{\ve_1,\ve_2}$ and prove that the invariants hold in this plane. Hence, the problem reduces to a statement in two-dimensional Euclidean geometry (however, when we analyze~(\ref{eq:overview_evolution}), we still use that the volume of ${\cE}_t$ is proportional to $ab^{d-1}$ and not $ab$).

Let us denote the coordinates corresponding to basis vectors $\ve_1$ and $\ve_2$ by $x$ and $y$. For brevity, let ${\cE} = {\cE}_{t-1}$, $\vz = \vz_t$, ${\cE}' = {\cE}_{t}$, $\vc = \vc_t =c\ve_1$, $\alpha = \alpha_{t-1}$, and $\alpha' = \alpha_t$.
We now need to choose parameters $a$, $b$, and $c$ so that invariants from Definition \ref{def:alg-invariant} and (\ref{eq:overview_evolution}) hold. See \Cref{fig:update_step}.
As shown in that figure, the new outer ellipse \(\vc + \cE'\) must contain the previous outer ellipse \(\cE\) and the newly received point \(\vz\). The new inner ellipse \(\vc + \alpha' \cE'\) must be contained within the convex hull of the previous inner ellipse \(\alpha \cE\) and \(\vz\).

It is instructive to consider what happens when point $\vz$ is at infinitesimal distance $\Delta$ from $\cE$: $\|\vz\|_2 = 1 +\Delta$. We consider a minimal axis-parallel outer ellipse $\cE'$ that contains $\cE$ and $\vz$. It must go through $\vz = (1+\Delta,0)$ and touch $\cE$ at two points symmetric w.r.t.\ the $x$-axis, say, $(-\sin \varphi, \pm \cos \varphi)$. Angle $\varphi$ uniquely determines $\cE'$. Now we want to find the largest value of the scaling parameter $\alpha'$ so that $\alpha' \cE'$ fits inside the convex hull of $\cE$ and $\vz$. When $\Delta$ is infinitesimal, this condition splits into two lower bounds on $\alpha'$ --  loosely speaking, they say that $\cE$ does not extend out beyond the convex hull in the horizontal (one bound) and vertical directions (the other). The former bound becomes stronger (gives a smaller upper bound on $\alpha'$) when $\varphi$ increases, and the latter becomes stronger when $\varphi$ decreases. When $\varphi = \alpha/2 \pm O(\alpha^2)$, then all terms linear in $\alpha$ vanish in both bounds and then $\alpha' = \alpha - \Theta(\alpha^2\Delta)$ satisfies both of them; for other choices of $\varphi$, we have $\alpha' \leq \alpha - \Omega(\alpha\Delta)$. 
So we let $\varphi=\alpha/2$ and from the formula for $\alpha'$ get $1/\alpha' = 1/\alpha + O(\Delta)$. On the other hand, $\vol(\cE') \geq (1 + \Delta/2)\vol(\cE)$, since $\cE'$ covers $\vz = (1 + \Delta,0)$. It is easy to see now that the evolution condition (\ref{eq:overview_evolution}) holds: the numerator is $O(\Delta)$ and the denominator is $\Omega(\Delta)$ in (\ref{eq:overview_evolution}).

We remark that letting $\vc + {\cE}'$ be the minimum volume ellipsoid that contains $\cE$ and $\vz$ is a highly suboptimal choice (it corresponds to setting $\varphi=\Theta(1/d)$). To derive our specific update formulas for arbitrary $\vz$, we, loosely speaking,  represent an arbitrary update as a series of infinitesimal updates, get a differential equation on $a$, $b$, $c$, and $\alpha'$, solve it, and then simplify the solution (remove non-essential terms, etc).  We get the following.

Our updates come from a family parameterized by \(\gamma \geq 0\). 
 Define \(\alpha'\) by \(\nfrac{1}{\alpha'} = \nfrac{1}{\alpha} + 2\gamma\). With this choice of $\alpha'$, define the new ellipses to be
\[\underbrace{\frac{1}{a^2}(x-c)^2 + \frac{1}{b^2} y^2 = 1}_{\vc + \cE'}, \qquad \underbrace{\frac{1}{a^2}(x-c)^2 + \frac{1}{b^2} y^2 = \alpha'^2}_{\vc + \alpha' \cE'}\]
where we use parameters
\begin{equation*}
\left.\begin{aligned}
a &= \expv{\gamma} \\
b &=1 + \frac{\alpha - \alpha'}{2} \\
c &= - \alpha + \alpha' \cdot a
\end{aligned}\qquad\right\}.
\end{equation*}
Choose $\gamma\approx \ln \|\vz\|_2$ so that $\vc + \cE'$ covers point $\vz$.
We use two-dimensional geometry to prove that $\cE'$, $\vc$, and $\alpha'$ satisfy the invariants (see Figure~\ref{fig:update_step}).
Now to prove the evolution condition, we observe two key properties: (1) the increase in the approximation factor is given by \(\frac{1}{\alpha'} - \frac{1}{\alpha} = 2 \gamma\) and (2) the length of the horizontal semiaxis of the new outer ellipse is \(\exp(\gamma)\). The length of the vertical semiaxis is at least $1$, so by the second property we have \(\log \vol(\cE') - \log \vol(\cE) \geq \gamma\). We combine this with the first property to prove that this update satisfies the evolution condition \eqref{eq:overview_evolution}.

Finally, we obtain an upper bound on $1/\alpha_n$ from the evolution equation. We have
\begin{align*}
    \nfrac{1}{\alpha_{n}} = \nfrac{1}{\alpha_0} + \sum_{t=1}^n
\left(\nfrac{1}{\alpha_{t}} - \nfrac{1}{\alpha_{t-1}}\right) \stackrel{\tiny(\text{by }\ref{eq:overview_evolution})}{\leq} 1+ C\sum_{t=1}^n(\log\mathsf{vol}(\cE_{t}) - \log\mathsf{vol}(\cE_{t-1})) = 1 + C\log \frac{\vol{\cE_n}}{\vol{\cE_0}}.
\end{align*}
It remains to get an upper bound on $\vol(\cE_n)$. We know that $\cE_n$ approximates $Z$, and $Z$, in turn, is contained in the ball of radius $R(Z)$. Loosely speaking, we get $\vol(\cE_n) \approx \vol(Z) \leq R(Z)^d \vol(B_2^d)$. Since $\cE_0$ is the ball of radius $r$, $\vol{\cE_0}= r^d \vol(B_2^d)$. We conclude that 
the approximation factor is at most $\nfrac{1}{\alpha_{n}} \lessapprox 1 + C\log \frac{R(Z)^d}{r^d} = 1 + O(d\log \frac{R(Z)}{r})$, as desired.

\paragraph{Removing the initialization assumption.} Once we have a monotone update rule and  guarantee on its approximation factor, we have to convert this to a guarantee where the algorithm does not have access to the initialization. 

One natural approach is as follows. Let $d' \le d$ be the largest timestep for which points $\vz_1, \dots,\vz_{d'+1}$ are in general position. We can compute the John ellipsoid for $\conv{\inbraces{\vz_1,\dots,\vz_{d'+1}}}$ and after that apply the monotone update rule guaranteed by Theorem \ref{thm:main_one} to obtain the rounding for every $t \ge d'+2$, so long as for every such timestep we have $\vz_t  \in \vspan{\vz_1,\dots,\vz_{t-1}}$.

\begin{figure}[h]
\centering
\includegraphics[width=0.70\textwidth]{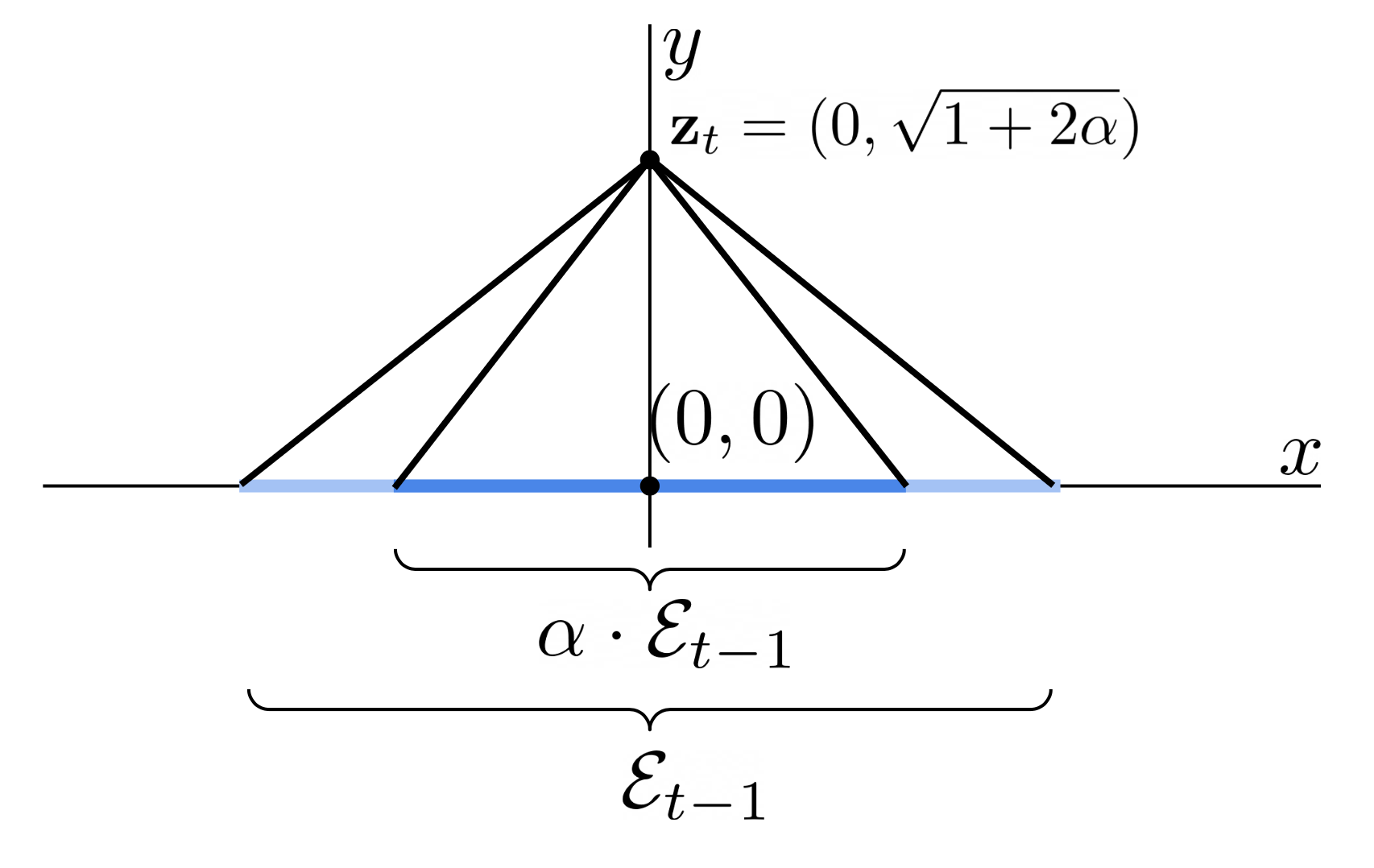}
\caption{Irregular update step. \(\cE_{t-1}\) and \(\alpha \cdot \cE_{t-1}\) are, respectively, the light blue strip on the \(x\)-axis and the dark blue strip on the \(x\)-axis. \(\vz_t = (0, \sqrt{1+2 \alpha})\) is the newly received point.
}
\label{fig:update_step_irreg}
\end{figure}

The principal difficulty in this approach is designing an \textit{irregular update step} that will handle points $\vz_t$ outside of $\vspan{\vz_1,\dots,\vz_{t-1}}$; when we add these points the dimensionality of the affine hull increases by 1. We consider the special case where the new point $\vz_t$ is conveniently located with respect to our previous ellipsoid $\cE_{t-1}$ (see \Cref{fig:update_step_irreg} for a $2$d-picture). Specifically, $\cE_{t-1}$ is the unit ball in $\lspan{\ve_1,\dots,\ve_{d'}}$, and the new point $\vz_t = (0,\dots, 0, \sqrt{1+2\alpha}),0,\dots)$. In $\vz_t$, only coordinate $d'+1$ is nonzero. We show that we can design an irregular update step for this special case that makes the new approximation factor $\nfrac{1}{\alpha_{t}}$ satisfy $\nfrac{1}{\alpha_t} = \nfrac{1}{\alpha_{t-1}}+1$.

It turns out that it is sufficient to consider only this special case. To see this, note that we can choose an affine transformation that maps any new point $\vz_{t}$ and previous ellipsoid $\cE_{t-1}$ to the setting shown in \Cref{fig:update_step_irreg}.
Next, observe that there are at most $d-1$ irregular update steps. This means that the irregular update steps contribute at most an additive $d-1$ to the final approximation factor.

Finally, observe that the inradius of $\conv{\inbraces{\vz_1,\dots,\vz_t}}$ is not monotone in $t$. In particular, it can decrease after each irregular update step. Nonetheless, we can still give a bound on the radius of a ball that our convex body $\conv{\vz_1,\dots,\vz_t}$ contains for all $t$. This will give us everything we need to apply \Cref{thm:main_one} to this setting, and \Cref{thm:main_two} follows.

\paragraph{Improved bounds on lattices.} 
Finally, we briefly discuss how to remove the aspect ratio dependence in the setting where the input points $\vz_t$ have coordinates in $[-N,N]$. At a high level, this improvement follows from carefully tracking how the approximation factors of our solutions change after an irregular update step. Following \eqref{eq:overview_evolution}, recall that our goal is to analyze (where we write $\alpha_0 = 1$)
\begin{align*}
    \sum_{t \ge 1} \frac{1}{\alpha_t}-\frac{1}{\alpha_{t-1}}.
\end{align*}
By \eqref{eq:overview_evolution}, we see that for all ``regular'' updates, we have
\begin{align*}
    \frac{1}{\alpha_t}-\frac{1}{\alpha_{t-1}} \lesssim \logv{\frac{\vol_{d_t}\inparen{\cE_t}}{\vol_{d_t}\inparen{\cE_{t-1}}}},
\end{align*}
where $d_t = \mathsf{dim}\inparen{\vspan{\vz_1,\dots,\vz_t}}$. Furthermore, as previously mentioned, in our irregular update step, we get
\begin{align*}
    \frac{1}{\alpha_t}-\frac{1}{\alpha_{t-1}} = 1.
\end{align*}
In order to control the sum of the $\nfrac{1}{\alpha_{t}}-\nfrac{1}{\alpha_{t-1}}$, it remains to bound $\nfrac{\vol_{d_t}(\cE_t)}{\vol_{d_{t-1}}(\cE_{t-1})}$ for an irregular update step $t$. We will then get a telescoping upper bound whose last term is the ratio of the volume of the final ellipsoid to the Euclidean ball in the same affine span.

Similarly to the improvements of \citet{woodruff2022high} in the integer-valued case, it will turn out that we will be interested in the total product of these volume changes. By carefully tracking these, we will get that this product can be expressed as the determinant of a particular integer-valued matrix. Then, since this matrix has integer entries, the magnitude of its determinant must be at least $1$. We then observe that the volume of $\cE_n$ after normalizing by the volume of $\vol(B_2^{d_n})$ must be at most $(N\sqrt{d})^{d_n}$, since the length of any vector in this lattice is at most $N\sqrt{d}$. The desired result then follows.

\subsection{Coresets for convex hull \texorpdfstring{(\Cref{thm:main_three})}{(Theorem~\ref{thm:main_three})}}
\label{sec:overview_coreset}

We now outline our proof strategy for \Cref{thm:main_three}.
Our main task is to design an appropriate selection criterion for every new point -- in other words, we must check whether a new point $\vz_t$ is ``important enough'' to be added to our previous set of points $S_{t-1}$. We then have to show that this selection criterion yields the approximation guarantee promised by Theorem \ref{thm:main_three}.

To design the selection criterion, we run an instance of the algorithm in Theorem \ref{thm:main_two} on the stream. For every new point $\vz_t$, we ask two questions -- ``Does $\vz_t$ result in an irregular update step? Does it cause $\mathsf{vol}(\cE_{t})$ to be much larger than $\mathsf{vol}(\cE_{t-1})$?'' If the answer to any of these questions is affirmative, we add $\vz_t$ to the coreset. The first question is necessary to obtain even a bounded approximation factor (for example, imagine that the final point $\vz_n$ results in an irregular update step -- then, we must add it). The second question is quite natural, as it ensures that the algorithm adds ``important points'' -- those that necessitate a significant update.


We now observe that at every irregular update step $t_{d'}$ for $d' \le d$ and subsequent timestep $t \ge t_{d'}$ for which there are no irregular update steps in between $t_{d'}$ and $t$, there exists a translation $\vc_{d'}$ (which is the center for $\cE_{d'}$ that the algorithm maintains) and a value $r_{d'}$ for which we know
\begin{align*}
    \vc_{d'} + r_{d'} \cdot \inparen{B_2^d \cap \lspan{\vz_1-\vc_{d'},\dots,\vz_{d'}-\vc_{d'}}} \subseteq \conv{\vz_1,\dots,\vz_{t}} \subseteq \vc_C + R_{t} \cdot B_2^d,
\end{align*}
where $\vc_C$ is the circumcenter of $\conv{\inbraces{\vz_1,\dots,\vz_t}}$. The resulting bound on $\abs{S_t}$ follows easily from the above observation and a simple volume argument.

Finally, we obtain the approximation guarantee from noting that for all \(t\), the output of the algorithm from Theorem \ref{thm:main_two} given the first \(t\) points is the same as running it only on the points selected by \(S_t\).

\subsection{Lower bound (Theorem \ref{thm:main_four})}

Whereas in the upper bound we demonstrated a particular algorithm that satisfies the evolution condition \eqref{eq:overview_evolution}, 
for the lower bound it suffices to show that for any monotone algorithm, there exists an instance of the problem (a sequence of $\vz_1$,\dots, $\vz_n$) where the algorithm must satisfy the ``reverse evolution condition'', i.e.
\begin{align}
    \frac{\nfrac{1}{\alpha_{t}} - \nfrac{1}{\alpha_{t-1}}}{\log\mathsf{vol}(\cE_{t}) - \log\mathsf{vol}(\cE_{t-1})} \ge C\label{eq:lb_overview_evolution}
\end{align}
for some $C > 0$.
In analogy to the argument of the upper bound, showing this reverse evolution condition yields a lower bound of the form \(\frac{1}{\alpha_n} \geq \widetilde{\Omega}\inparen{d \log(\kappa)}\).
Given any monotone algorithm \(\cA\), the instance we use is produced by an adversary that repeatedly feeds \(\cA\) a point that is a constant factor away from the previous ellipsoid.

In order to simplify showing this reverse evolution condition, we use a symmetrization argument.
Specifically, by a particular sequence of Steiner symmetrizations, we see that the optimal response of \(\cA\) can be completely described in two dimensions. Thus, it is sufficient to only show this reverse evolution condition in the two-dimensional case where the previous outer ellipsoid is the unit ball.

This transformed two-dimensional setting is significantly simpler to analyze. Specifically, we can assume that the point given by the adversary is always \(2 \ve_1\). The rest of the argument proceeds by cases, again using two-dimensional Euclidean geometry. On a high level, the constraints placed on the new outer and inner ellipsoid by the monotonicity condition force the update of \(\cA\) to satisfy the reverse evolution condition.
\section{Preliminaries}
\label{section:preliminaries_notation}

\subsection{Notation} 

We denote the standard Euclidean norm of a vector $\vv$ by $\norm{\vv}$ and the Frobenius norm of a matrix $\mA$ by $\fnorm{\mA}$. We denote the singular values of a matrix $\mA \in \R^{d\times d}$ by $\sigma_1(\mA),\dots, \sigma_d(\mA)$. Let $\sigma_{\max}(\mA)$ and $\sigma_{\min}(\mA)$ be the largest and smallest singular values of $\mA$, respectively.
We write \(\diag{a_1, \ldots, a_d}\) to mean the \(d \times d\) diagonal matrix whose diagonal entries are \(a_1, \ldots, a_d\).
We use \(\mS_{++}^d\) to denote the set of \(d \times d\) positive definite matrices.
We use \(\ve_1, \ldots, \ve_d\) for the standard basis in \(\R^d\).

Denote the $\ell_2$-unit ball by $B_2^d = \inbraces{\vx \in \R^d \suchthat \norm{\vx} \le 1}$,
and \(\S^{d-1} = \inbraces{\vx \in \R^d \colon \|\vx\|_2 = 1}\) the unit Euclidean sphere. We use \(\partial S\) for the boundary of an arbitrary set \(S\). We use natural logarithms unless otherwise specified. 

In this chapter, we will work extensively with ellipsoids. We will always assume that all ellipsoids and balls we consider are centered at the origin. We use the following representation of ellipsoids. For a non-singular matrix $\mA\in\R^{d\times d}$, let $\cE_\mA \coloneqq \inbraces{\vx \suchthat \norm{\mA\vx}\le 1}$. In other words, the matrix $\mA$ defines an bijective linear map satisfying $\mA \cE_\mA = B_2^d$. Every full-dimensional ellipsoid (centered at the origin) has such a representation. We note that this representation is not unique as matrices $\mA$ and $\mM\mA$ define the same ellipsoid if matrix $\mM$ is orthogonal (since $\norm{\mA\vv} = \norm{\mM\mA\vv}$ for every vector $\vv$). Sometimes, we will have to consider lower-dimensional ellipsoids within an ambient space of higher dimension; in this case, we will use the notation $\cE \cap H$ where $H$ is some linear or affine subspace -- note that $\cE \cap H$ is also an ellipsoid. 

Now consider the singular value decomposition of  $\mA$: $\mA = \mU\Sigma^{-1} \mV^T$ (it will be convenient for us to write $\Sigma^{-1}$ instead of standard $\Sigma$ in the decomposition). The diagonal entries of $\Sigma$ are exactly the semi-axes of $\cE_\mA$. As mentioned above, matrices $\mU\Sigma^{-1} \mV^T$ and $\mU'\Sigma^{-1} \mV^T$ define the same ellipsoid for any orthogonal \(\mU' \in \mathbb{R}^{d \times d}\); in particular, every ellipsoid can be represented by a matrix of the form $\mA=\Sigma^{-1} \mV^T$.

\subsection{Geometry}
We restate the well-known result that five points determine an ellipse.
This is usually phrased for conics, but for nondegenerate ellipses the usual condition that no three of the five points are collinear is vacuously true.
\begin{lemma}[Five points determine an ellipse]\label{lemma:conic_5pt}
Let \(\vc_1 + \partial \cE_1, \vc_2 + \partial \cE_2\) be two ellipses in \(\R^2\).
If they intersect at five distinct points, then \(\vc_1 + \partial \cE_1\) and \(\vc_2 + \partial \cE_2\) are the same.
\end{lemma}

The following claim, that every full-rank ellipsoid (i.e. an ellipsoid whose span has full dimension) can be represented by a positive definite matrix, follows from looking at the singular value decomposition of $\mA$.
\begin{lemma}\label{lemma:prelim_ellips_sym}
    Let \(\cE \subseteq \R^d\) be a full-rank ellipsoid.
    Then there exists \(\mA \succ 0\) such that \(\cE = \cE_{\mA}\).
\end{lemma}

We also have the standard result relating volume and determinants, which follows from observing \(\mA \cE_{\mA} = B_2^d\).
\begin{lemma}\label{lemma:prelim_vol}
Let \(\mA \succ 0\). Then 
\(\vol(\cE_{\mA}) = \det(\mA^{-1}) \vol(B_2^d)\).
\end{lemma}

In order to give the reduction in the lower bound from the general case to the two-dimensional case, we use the technique of Steiner symmetrization (see e.g. \cite[Section 1.1.7]{artstein2015asymptotic}).
Given some unit vector \(\vu \in \R^d\) and convex body \(K \subseteq \R^d\), we write \(S_{\vu}(K)\) for the \textit{Steiner symmetrization} in the direction of \(\vu\).
Recall that the Steiner symmetrization is defined so that for any \(\vx \perp \vu\),
\[\vol((\vx + \R \vu) \cap K) = \vol((\vx + \R \vu) \cap S_{\vu}(K)),\]
and so that \((\vx + \R \vu) \cap S_{\vu}(K)\) is an interval centered at \(\vx\).
Note that we will overload notation slightly as we will allow you \(\vu\) to be a vector of any non-zero length while Steiner symmetrization is usually defined with \(\vu\) being a unit vector, but we will simply take \(S_{\vu} = S_{\frac{\vu}{\|\vu\|_2}}\).

Importantly, Steiner symmetrization will preserve important properties of the update.
We have the key facts that \(\vol(S_{\vu}(K)) = \vol(K)\), \(S_{\vu}(K') \subseteq S_{\vu}(K)\) if \(K \subseteq K'\), and further the Steiner symmetrization preserves \(K\) being an ellipsoid:
\begin{lemma}[\protect{\cite[Lemma 2]{bourgain2006estimates}}]\label{lemma:steiner_ellipsoid}
    If \(c + \cE \subseteq \R^d\) is an ellipsoid, \(S_{\vu}(c + \cE)\) is still an ellipsoid.
\end{lemma}

Further, if we apply Steiner symmetrization to a body that is a body of revolution about an axis,
it does not change the body if \(\vu\) is perpendicular to the axis of revolution.
\begin{lemma}\label{lemma:steiner_sym}
    Let \(K \subseteq \R^d\) be a body of revolution about the \(\ve_1\)-axis.
    Then if \(\vu \perp \ve_1\), \(S_{\vu}(K) = K\).
\end{lemma}

\section{Streaming ellipsoidal rounding}
\label{sec:main_alg_top}

Our goal in this section is to prove Theorems \ref{thm:main_one} and \ref{thm:main_two}.

\subsection{Monotone algorithms solve \texorpdfstring{\Cref{prob:main}}{Theorem \ref{prob:main}}}

To design algorithms to solve the streaming ellipsoidal rounding problem,
we first show that any monotone algorithm gives a valid solution.
We let \(\vc_0 \in \R^d\) and \(r_0 \geq 0\) be given so that
\(\vc_0 + r_0 \cdot B_2^d \subseteq Z\),
and denote the initial ellipsoid as \(\cE_0 = r_0 \cdot B_2^d\).
Note that \(r_0\) need not be the inradius, although it is upper bounded by the inradius.

If we had for each intermediate step \(t\)
that \(\vc_t + \alpha_t \cdot \cE_t \subseteq \conv{\vz_1, \ldots \vz_t} \subseteq \vc_t + \cE_t\),
then clearly any algorithm that satisfies this would give a valid final solution as well.
However, in intermediate steps it is not clear that 
\(\vc_t + \alpha_t \cdot \cE_t \subseteq \conv{\vz_1, \ldots \vz_t}\),
due to the initialization of \(\vc_0 + \cE_0\) in our monotone algorithm framework.
Instead, we relax this invariant to \(\vc_t + \alpha_t \cdot \cE_t \subseteq \conv{\{\vz_1, \ldots \vz_t\} \cup (\vc_0 + \cE_0)}\),
which still suffices to produce a valid final solution.
\begin{lemma}\label{lemma:alg_invariant}
To solve \Cref{prob:main}, it suffices for 
the sequence of ellipsoids \(\vc_i + \cE_i\) and scalings \(\alpha_i\) to satisfy the invariants of \Cref{def:alg-invariant}.
\end{lemma}

\begin{proof}
First, we argue that \(\conv{\vz_1, \ldots, \vz_n} \subseteq \vc_n + \cE_n\).
As \(\cE_n\) is an ellipsoid and therefore a convex set, it suffices to show \(\{\vz_1, \ldots, \vz_n\} \subseteq \vc_n + \cE_n\).
We actually argue by induction that \(\{\vz_1, \ldots, \vz_t\} \subseteq \vc_t + \cE_t\) for all \(0 \leq t \leq n\).
This is vacuously true for \(t = 0\). At each step \(t > 0\) the inductive hypothesis gives \(\{\vz_1, \ldots, \vz_{t-1}\} \subseteq \vc_{t-1} + \cE_{t-1}\), and thus by (\ref{item:def-invariant-1}) we have \(\{\vz_1, \ldots, \vz_{t}\} \subseteq \vc_t + \cE_t\).

Now, we argue that \(\vc_n + \alpha_n \cdot \cE_n \subseteq \conv{\vz_1, \ldots, \vz_n}\).
We show by induction that \(\vc_t + \alpha_t \cdot \cE_t \subseteq \conv{\{\vz_1, \ldots, \vz_t\} \cup (\vc_0 + \cE_0)}\)
for all \(0 \leq t \leq n\).
This is sufficient as \(\conv{\{\vz_1, \ldots, \vz_n\} \cup (\vc_0 + \cE_0)} = Z\).
The case for \(t = 0\) is trivial.
For \(t > 0\), the inductive hypothesis gives \(\vc_{t-1} + \alpha_{t-1} \cdot \cE_{t-1} \subseteq \conv{\{\vz_1, \ldots, \vz_{t-1}\} \cup (\vc_0+\cE_0)}\),
and by (\ref{item:def-invariant-2}) we have
\[c_t + \alpha_t \cdot \cE_t \subseteq \conv{(c_{t-1} + \alpha_{t-1} \cdot \cE_{t-1}) \cup \{\vz_i\}} \subseteq \conv{\{\vz_1, \ldots, \vz_t\} \cup (\vc_0 + \cE_0)},\]
as desired.
\end{proof}

\subsection{Special case}\label{sec:two_dim_update}

In light of \Cref{lemma:alg_invariant}, our strategy is to design an algorithm that preserves the invariants given in \Cref{def:alg-invariant}. This algorithm can be thought of as an \textit{update rule} that,
given the previous outer and inner ellipsoids \(\vc_{t-1} + \cE_{t-1}, \vc_{t-1} + \alpha_{t-1} \cE_{t-1}\) and next point \(\vz_t\), produces the next outer and inner ellipsoids \(\vc_t + \cE_t, \vc_t + \alpha_t \cE_t\).

It is in fact sufficient to consider the simplified case where the previous outer ellipsoid is the unit ball, and the previous inner ellipsoid is some scaling of the unit ball;
we will show this in \Cref{sec:gen_high_dim}. 
We can further specialize by considering only the two-dimensional case \(d=2\).
We will later show that the high-dimensional case is not much different, as all the relevant sets \(\vc_{t_1} + \cE_{t-1}, \vc_t + \cE_{t}\) and \(\conv{\alpha \cdot \cE_{t-1} \cup \{\vz_t\}}\) form bodies of revolution about the axis through \(\vc_{t-1}\) and \(\vz_t\).

We now describe our two-dimensional update rule.
In order to simplify notation, we will let \(\alpha\) be the previous scaling \(\alpha_{t-1}\), and \(\alpha'\)
be the next scaling \(\alpha_t\). 
We will assume that \(\alpha \leq \nfrac{1}{2}\) to simplify the analysis
of our update rule; this will not affect the quality of our final approximation as this update rule
will only be used in the ``large approximation factor'' regime.
We will also overload notation by writing \(c + \cE\) even when \(c\)
is a scalar to mean \((c, 0) + \cE\).
We can describe the previous outer ellipsoid \(\cE\) with the equation \(x^2 + y^2 \leq 1\), and the previous inner ellipsoid \(\alpha \cE\) with \(x^2 + y^2 \leq \alpha^2\).
We define the next outer and inner ellipsoids \(c + \cE'\), \(c + \alpha' \cE'\) as
\[\underbrace{\frac{1}{a^2}(x-c)^2 + \frac{1}{b^2} y^2 \leq 1}_{c + \cE'}, \qquad \underbrace{\frac{1}{a^2}(x-c)^2 + \frac{1}{b^2} y^2 \leq \alpha'^2}_{c + \alpha' \cE'}\]
where we use parameters
\begin{equation}\label{eqn:update_params}
\left.\begin{aligned}
a &= \expv{\gamma} \\
b &=1 + \frac{\alpha - \alpha'}{2} \\
c &= - \alpha + \alpha' \cdot a \\
\alpha'& = \frac{1}{\frac{1}{\alpha} + 2\gamma} \\
\end{aligned}\qquad\right\}
\end{equation}

We will let \(\vz\) be the rightmost point of \(c + \cE'\), so that \(\vz = (c+ a, 0)\). Eventually, we will choose \(\gamma\) so that \(\vz\) coincides with \(\vz_{t}\), the point received in the next iteration.
In \Cref{sec:gen_alg}, these parameters \(a(\gamma), b(\gamma), c(\gamma), \alpha'(\gamma)\)
will be used as functions of the parameter \(\gamma \geq 0\). However, we will not yet explicitly specify \(\gamma\), so in this section these parameters can be thought of as constants for some fixed \(\gamma\). 
This update rule is pictured in \Cref{fig:update_step}.

We first collect a few straightforward properties of this update rule.
\begin{lemma}\label{lemma:update_params}
The parameters in the setup (\ref{eqn:update_params}) satisfy the following.
\begin{enumerate}
    \item\label{item:claim_update_params_1} \(\frac{1}{\alpha'} = \frac{1}{\alpha} + 2\gamma\)
    \item\label{item:claim_update_params_2} \(b \geq 1\)
    \item\label{item:claim_update_params_3} \(c \geq 0\)
    \item\label{item:claim_update_params_4} \(c + \alpha' \cdot a \geq \alpha\)
\end{enumerate}
\end{lemma}
Before proving these properties, we provide geometric interpretations.
Intuitively, (\ref{item:claim_update_params_1}) means that \(\gamma\) is proportional to the increase in the approximation factor at this step, a fact that we will use when analyzing the general-case algorithm.
(\ref{item:claim_update_params_2}) means that the outer ellipsoid grows on every axis;
and (\ref{item:claim_update_params_3}) means that the centers of the next ellipsoids are to the right of the \(y\)-axis, i.e. the centers of the next ellipsoids are further towards \(\vv\) than those of the previous ellipsoids.
The rightmost point of \(c + \alpha' \cE'\) is \(c + \alpha' \cdot a\), so (\ref{item:claim_update_params_4}) shows that this point is to the right of the rightmost point of \(\alpha \cdot \cE\).

We now prove \Cref{lemma:update_params}.

\begin{proof}[Proof of \Cref{lemma:update_params}]
(\ref{item:claim_update_params_1}) is clear from rearranging the definition of \(\alpha'\).
From (\ref{item:claim_update_params_1}) we also have \(\alpha' \leq \alpha\), so that (\ref{item:claim_update_params_2}) follows immediately.

For (\ref{item:claim_update_params_3}), observe that \(\frac{\alpha}{\alpha'} = 1 + 2 \gamma \alpha\).
When \(\alpha \leq \nfrac{1}{2}\), this means \begin{equation}\label{eqn:pup} 
\frac{\alpha}{\alpha'} \leq 1 + \gamma \leq \expv{\gamma} = a
\end{equation}
using \(1 + x \leq e^x\), \Cref{lemma:wk_e}-(\ref{item:wk_1}).
By definition of \(c\), \(\alpha/\alpha' \leq a\) is equivalent to \(c \geq 0\).

To show (\ref{item:claim_update_params_4}), by definition we have that \(c + \alpha' \cdot a = - \alpha + 2 \alpha' a\).
Thus showing \(c + \alpha' \cdot a \geq \alpha\) is equivalent to showing that \(\alpha' a \geq \alpha\), which is equivalent to the inequality in (\ref{eqn:pup}).
\end{proof}

As \Cref{fig:update_step} depicts, the update step we defined satisfies the invariants in \Cref{def:alg-invariant} and so is monotone;
in the rest of this section we make this picture formal.
To start, we consider the invariant concerning outer ellipsoids; we will show that \(\cE \subseteq c + \cE'\).
For now we can think of \(\vz\) as replacing \(\vz_t\)
, and clearly \(\vz \in c + \cE'\), so if we show that \(\cE \subseteq c + \cE'\), then \(\conv{\cE \cup \{\vz\}} \subseteq c + \cE'\) as well since \(c + \cE'\) is convex.

\begin{lemma}\label{lemma:update_step_outer}
    We have \(\cE \subseteq c + \cE'\).
\end{lemma}
\begin{proof}
First, observe that \(\cE \subseteq \cE'\) because both axes of \(\cE'\) have greater length than those of \(\cE\): \(a \geq 1\) by definition, and \(b \geq 1\) from \Cref{lemma:update_params}-(\ref{item:claim_update_params_2}).
Now, we translate \(\cE'\) to the right until it touches \(\cE\) at two points. 
We call this translated ellipse \(c_r + \cE'\), as shown in \Cref{fig:update_step_outer}.
Observe that as long as \(c \leq c_r\), we have \(\cE \subseteq c + \cE'\).
We now determine \(c_r\).

\begin{figure}[ht]
\centering
\includegraphics[width=0.70\textwidth]{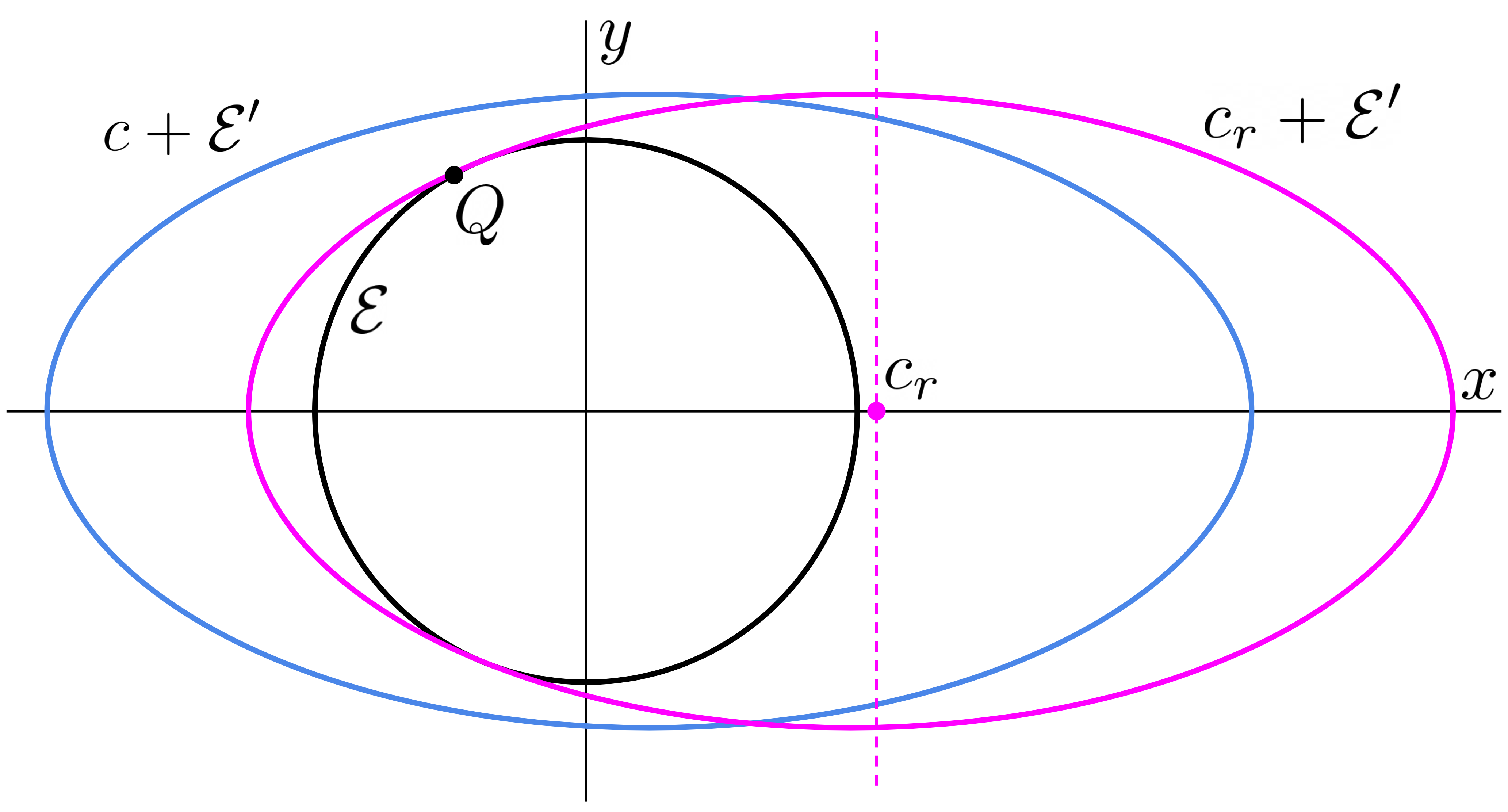}
\caption{Outer ellipses of the update step. As before, \(\cE\) is the black circle and \(c + \cE'\) is the blue ellipse.
\(c_r + \cE'\) is the magenta ellipse, with its center at \(c_r\) and the dotted magenta line showing the position of \(c_r\) along the \(x\)-axis.
\(c_r\) is defined so \(c_r + \cE'\) and \(\cE\) are tangent at two points. \(Q\) is one of these two tangent points.}
\label{fig:update_step_outer}
\end{figure}

First, note points on the boundary of \(c_r + \mathcal{E}'\) are described by the equation \begin{equation}\label{eqn:def_cr}
\frac{(x - c_r)^2}{a^2} + \frac{y^2}{b^2} = 1
\end{equation}
Let \(Q = (x', y')\) be the point of intersection between \(\cE\) and \(c_r + \cE'\) where \(y' > 0\).
Since \(Q\) is on the boundary of both ellipses, the vectors \(\left(\frac{2(x'-c_r)}{a^2}, \frac{2 y'}{b^2}\right)\) and \(\left(2x', 2y'\right)\), which are the normal vectors at \(Q\) of \(c_r + \cE'\) and \(\cE\) respectively, must be parallel.
Thus \(\frac{4 (x' - c_r)}{a^2} \cdot y' = \frac{4 y' x'}{b^2}\), which simplifies to \begin{equation}\label{eqn:x_cr}
x' = \frac{c_r}{1 - \frac{a^2}{b^2}}.
\end{equation}
At this point we have a system of three equations relating \((x', y')\) and \(c_r\): (\ref{eqn:x_cr}), \(Q\) lying on \(\cE\), and \(Q\) satisfying (\ref{eqn:def_cr}).
We now solve this system to find \(c_r\).
To start, we expand (\ref{eqn:def_cr}) into
\(x'^2 - 2 x' c_r + c_r^2 + y'^2 \frac{a^2}{b^2} = a^2\),
which we rewrite into \(x'^2 \frac{a^2}{b^2} + x'^2 \left(1 - \frac{a^2}{b^2}\right) - 2 x' c_r + c_r^2 + y'^2 \frac{a^2}{b^2} = a^2\).
As \(Q\) lies on \(\cE\), this becomes \(x'^2 \left(1 - \frac{a^2}{b^2}\right) - 2 x' c_r + c_r^2 +  \frac{a^2}{b^2} = a^2\).
Substituting in (\ref{eqn:x_cr}), we get
\[\frac{c_r^2}{1 - \frac{a^2}{b^2}} - 2 \frac{c_r^2}{1 - \frac{a^2}{b^2}} + c_r^2 + \frac{a^2}{b^2} = a^2.\]
Simplifying, we have \(c_r^2 \left(1 - \frac{b^2}{b^2 - a^2}\right) = a^2 \left(1 - \frac{1}{b^2}\right)\), i.e.
\[c_r^2 = \frac{b^2 - 1}{b^2} (a^2 - b^2).\]
To complete the proof of \cref{lemma:update_step_outer}, it suffices to show \(c^2 \leq \frac{b^2 - 1}{b^2} (a^2 - b^2)\). This will follow from \Cref{lemma:pty_outer}.
\end{proof}

Now, we move on to the inner ellipsoid invariant of \Cref{def:alg-invariant}.
In particular, we will argue that \(c + \alpha' \cE' \subseteq \conv{\alpha\cE \cup \{\vz\}}\).
On a high level, we show this by arguing that the boundary of \(c + \alpha' \cE'\) does not intersect the boundary of \(\conv{\alpha\cE \cup \{\vz\}}\), except at points of tangency.

\begin{figure}[h]
\centering
\includegraphics[width=0.70\textwidth]{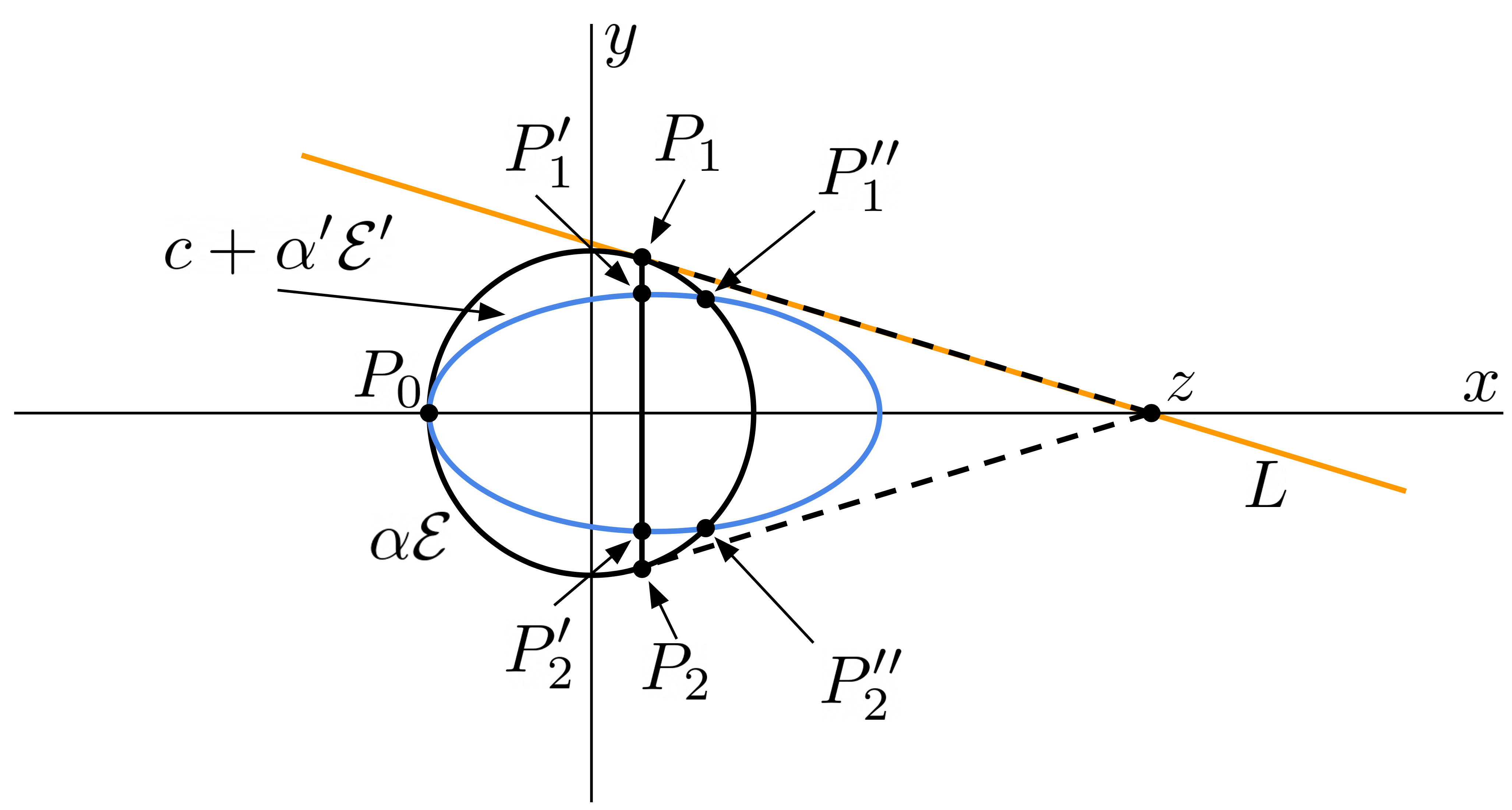}
\caption{Inner ellipses of the update step. As before, \(\alpha \cE\) is the black circle and \(c + \alpha' \cE\) is the blue ellipse. 
\(P_0\) is the shared leftmost point of \(\alpha \cE\) and \(c + \alpha' \cE'\).
There are two lines through \(\vv\) that are tangent to \(\alpha \cE\), one of which we call \(L\) and pictured in orange. We call the tangent points \(P_1\) and \(P_2\).
The line segments \(\overline{P_1 \vz}, \overline{P_2 \vz}\) are the dotted black lines.
\(P_1'\) and \(P_2'\) are the two points of intersection between \(c + \alpha' \cE\) and the line segment \(\overline{P_1 P_2}\). \(P_1''\) and \(P_2''\)
are the two points of intersection between \(\partial (c + \alpha' \cE')\) and \(\partial \alpha \cE\) to the right of the \(y\)-axis. Note that \(P_2, P_2', P_2''\) are the reflections of \(P_1, P_1', P_1''\) across the \(x\)-axis.}
\label{fig:update_step_inner}
\end{figure}

We can split the boundary of \(\conv{\alpha\cE \cup \{\vz\}}\) into two pieces: the part that intersects with the boundary of \(\alpha \cE\), which is an arc of the boundary of \(\alpha \cE\); and the remainder, which can described as two line segments connecting \(\vz\) to that arc. In particular, there are two lines that go through \(\vz\) and are tangent to \(\alpha \cE\), one of which we call line \(L\), and the other line is the reflection of \(L\) across the \(x\)-axis.
We define \(P_1\) and \(P_2\) as the tangent points of these lines to \(\alpha \cE\). 
Then, the boundary of \(\conv{\alpha\cE \cup \{\vz\}}\) consists of an arc \(P_1 P_2\) and the line segments \(\overline{P_1 \vz}, \overline{P_2 \vz}\). This is illustrated in \Cref{fig:update_step_inner}.
Note that at this point it is possible a priori for the arc \(P_1 P_2\) that coincides with the boundary of \(\conv{\alpha\cE \cup \{\vz\}}\)
to be either the major or minor arc; we will later show it must be the major arc.
We will take \(L\) to be the line whose tangent point to \(\alpha \cE\), \(P_1\), is above the \(x\)-axis, though this choice is arbitrary due to symmetry across the \(x\)-axis.

We first show that \(c + \alpha' \cE'\) does not intersect \(\overline{P_1 \vz}\) and \(\overline{P_2 \vz}\), except possibly at points of tangency.
In fact, we show a slightly stronger statement, in similar fashion to \Cref{lemma:update_step_outer}.

\begin{lemma}\label{lemma:translate_ellipse_angle}
 \(c + \alpha' \cE'\) lies inside the angle \(\angle P_1 \vz P_2\).
\end{lemma}
\begin{proof}
We translate \(c + \alpha' \cE'\) to the right until it touches \(L\) (and, by symmetry, \(\overline{P_2 \vz}\)).
We call this translated ellipse \(c_+ + \alpha' \cE'\), as shown in \Cref{fig:update_step_inner_translate}. 
(Formally, the center \(c_+\) can be described not as a translation from some other ellipse,
but as \(c_+\) such that \(c_+ + \alpha' \cE'\) intersects \(L\) at one point).
Observe that 
if \(c \leq c_{+}\), then \(c + \alpha' \cE'\) lies inside the angle \(\angle P_1 \vz P_2\).
We now determine \(c_{+}\).

\begin{figure}[h]
\centering
\includegraphics[width=0.70\textwidth]{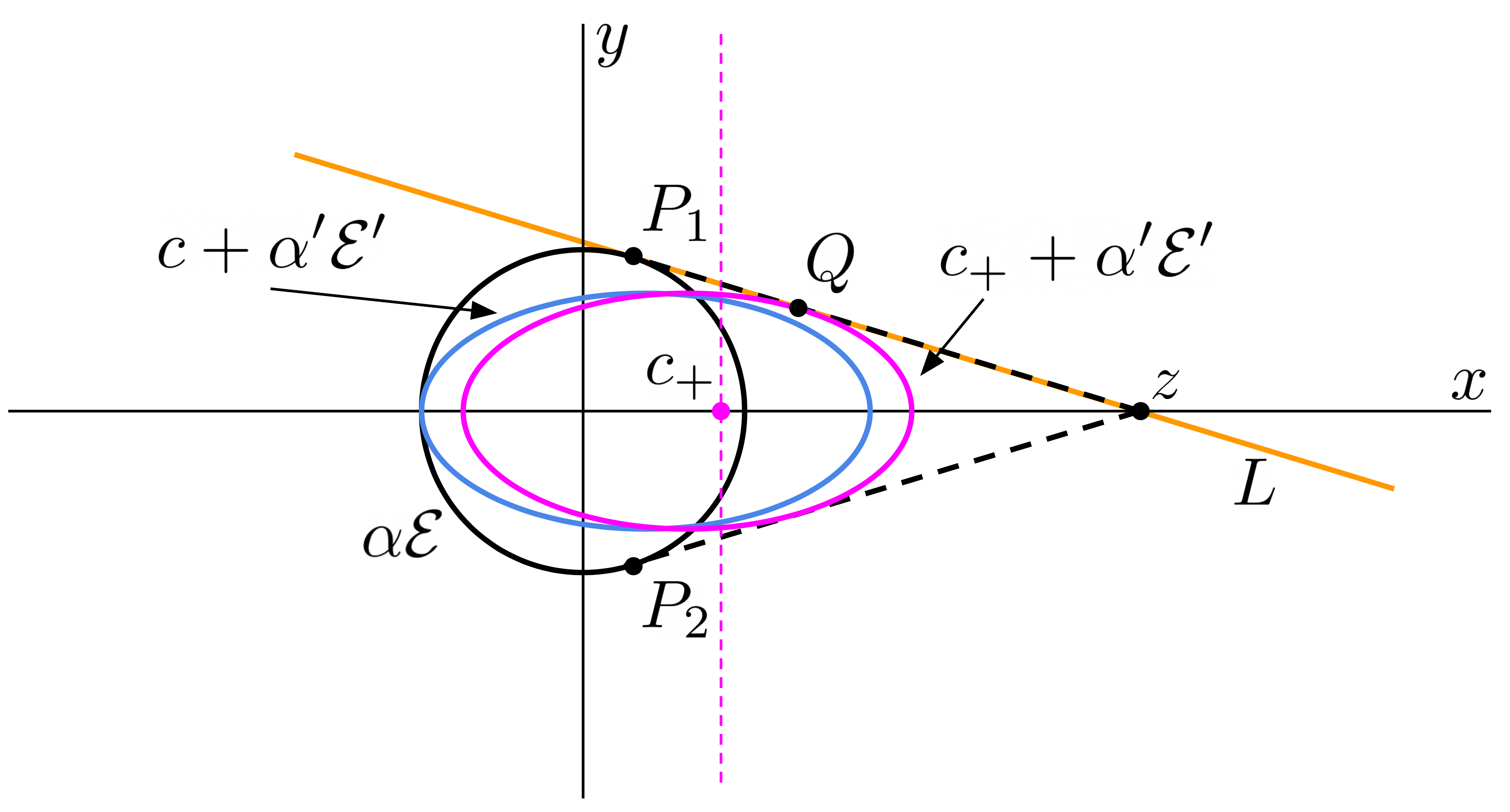}
\caption{Inner ellipses of the update step. As before, \(\alpha \cE\) is the black circle, \(c + \alpha' \cE\) is the blue ellipse, \(L\) is the orange line through \(\vz\) and tangent to \(\alpha \cE\), \(P_1\) and \(P_2\) are the tangent points on the lines through \(\vz\) tangent to \(\alpha \cE\), and \(\overline{P_1 \vz}, \overline{P_2 \vz}\) are the dotted black lines.
\(c_{+} + \alpha' \cE'\) is the magenta ellipse, with its center at \(c_{+}\) and magenta dotted line showing its position on the \(x\)-axis.
\(c_{+}\) is defined so that \(c_{+} + \alpha' \cE'\) is tangent to \(\overline{P_1 \vz}\) and \(\overline{P_2 \vz}\), with \(Q\) as the tangent point of \(c_{+} + \alpha' \cE'\) and \(\overline{P_1 \vz}\).
}
\label{fig:update_step_inner_translate}
\end{figure}

The equation of \(L\) is \[\underbrace{\frac{1}{c+a}}_{\ell_1} \cdot x + \underbrace{\sqrt{\frac{1}{\alpha^2} - \frac{1}{(c+a)^2}}}_{\ell_2} \cdot y = 1\]
where we define \(\ell_1, \ell_2\) as the coefficents for \(x\) and \(y\).
Observe that \(\vz\) is on \(L\), and \(L\) is tangent to \(\alpha \cE\) at \(P_1\), which has coordinates \begin{equation}\label{eqn:usit_p1} 
P_1 = \left(\frac{\alpha^2}{c+a}, \alpha^2 \sqrt{\frac{1}{\alpha^2} - \frac{1}{(c+a)^2}} \right).
\end{equation}
Tangency can be confirmed by checking that \(P_1\) is parallel to \((\ell_1, \ell_2)\), the normal vector definining \(L\).

Let \(Q = (x', y')\) be the point of intersection of \(L\) and \(c_{+} + \alpha' \cE\), there are three properties that define \(Q\). First it lies on the boundary of \(c_{+} + \alpha' \cE\), so it satisfies
\begin{equation}\label{eqn:pty_inner_oe} 
    \frac{(x'-c_{+})^2}{a^2} + \frac{y'^2}{b^2} = \alpha'^2.
\end{equation}

Second, at \(Q\) the normal vectors for the equations defining \(c_{+} + \alpha' \cE\) and \(L\) are parallel, i.e.
\((\frac{2(x-c_{+})}{a^2},\frac{2y}{b^2})\) is parallel to \((\ell_1, \ell_2)\). So
\begin{equation}\label{eqn:pty_inner_gm} 
    \frac{(x'-c_{+})}{a^2} \ell_2 = \frac{y'}{b^2} \ell_1.
\end{equation}
Finally, \(Q\) lies on \(L\), so we have \(\ell_1 x' + \ell_2 y' = 1\).
Solving this for \(y'\), we get
\begin{equation}\label{eqn:pty_inner_ol} 
    y' = \frac{1-\ell_1 x'}{\ell_2}.
\end{equation}
These three equations form a system for \(x', y'\) and \(c_{+}\), which we now solve to find \(c_{+}\).
Taking the square of (\ref{eqn:pty_inner_gm}) and rearranging gives \(\frac{y'^2}{b^2} = \frac{b^2 (x'-c_{+})^2 \ell_2^2}{a^4 \ell_1^2} \).
Substituting this into (\ref{eqn:pty_inner_oe}), we get
\(\frac{(x'-c_{+})^2}{a^2} + \frac{b^2 (x'-c_{+})^2 \ell_2^2}{a^4 \ell_1^2} = \alpha'^2\).
Now, defining \(r \defeq\frac{a^2 \ell_1^2}{b^2 \ell_2^2}\), we group the terms of this equation into the form
\begin{equation}\label{eqn:pty_inner_1}
    (x'-c_{+})^2 \cdot \frac{1}{a^2} \left(1 + \frac{1}{r}\right) = \alpha'^2.
\end{equation}
We substitute (\ref{eqn:pty_inner_ol}) into (\ref{eqn:pty_inner_gm})
to get \(\frac{x'-c_{+}}{a^2} \ell_2 = \frac{\ell_1}{b^2} \frac{1 - x' \ell_1}{\ell_2}\).
Grouping for \(x'\) and rearranging yields
\begin{equation}\label{eqn:pty_inner_2}
    x' - c_{+} = \frac{r}{1+r} \left(\frac{1}{\ell_1} - c_{+}\right).
\end{equation}
Next, we substitute (\ref{eqn:pty_inner_2}) into (\ref{eqn:pty_inner_1}), and get after some cancellation
\[ \left(\frac{1}{\ell_1} - c_{+}\right)^2 = \alpha'^2 a^2 \cdot \frac{1+r}{r}.\]
Observe on the left hand side that \(\frac{1}{\ell_1} - c_{+} = c+a-c_{+}\).
Clearly the center \(c_{+}\) must be to the left of \(\vz\), so this must be non-negative.
Hence after taking the positive square root, we obtain
\[c_{+} = c + a - \alpha' \cdot a \sqrt{\frac{1+r}{r}}\]
It remains to show that \(c \leq c_{+}\), or equivalently that
\[a - \alpha' \cdot a \sqrt{\frac{1+r}{r}} \geq 0,\]
which we do in \Cref{lemma:pty_inner}. This completes the proof of \Cref{lemma:translate_ellipse_angle}. 
\end{proof}

Now, we build on the previous claim to show the inner ellipsoid invariant.
\begin{lemma}\label{lemma:update_step_inner}
    We have \(c + \alpha' \cdot \cE' \subseteq \conv{\alpha\cdot\cE \cup \{\vz\}}\).
\end{lemma}
\begin{proof}
We will argue that the boundary of \(c + \alpha' \cE'\) does not intersect the boundary of \(\conv{\alpha\cE \cup \{\vz\}}\), except at points of tangency.
This is sufficient to establish the claim, as \Cref{lemma:translate_ellipse_angle} shows that \(c + \alpha' \cE'\) is internal to \(\angle P_1 \vz P_2\), and so if \(c + \alpha' \cE'\) does not intersect the boundary of \(\conv{\alpha\cE \cup \{\vz\}}\), \(c + \alpha' \cE'\) must lie inside of, or be disjoint from \(\conv{\alpha\cE \cup \{\vz\}}\).
Since the leftmost points of \(\alpha \cE\) and \(c + \alpha' \cE'\) coincide, \(c + \alpha' \cE'\) must then lie inside of  \(\conv{\alpha\cE \cup \{\vz\}}\).
Recall that the boundary of \(\conv{\alpha\cE \cup \{\vz\}}\) consists of the arc \(P_1 P_2\) and the line segments \(\overline{P_1 \vz}, \overline{P_2 \vz}\).
\Cref{lemma:translate_ellipse_angle} already shows that the boundary of \(c + \alpha' \cE'\) does not intersect \(\overline{P_1 \vz}\) and \( \overline{P_2 \vz}\), so we only need to show that the boundary of \(c + \alpha' \cE'\) does not intersect the arc \(P_1 P_2\).

To do this, we start by enumerating the points of intersection of \(\partial \alpha \cE\) and \(\partial(c + \alpha' \cE')\), recalling that \(P_1 P_2\) is an arc of \(\partial \alpha \cE\).
Observe that the leftmost points of \(\alpha \cE\) and \(c + \alpha' \cE'\) coincide, as the leftmost point of \(c + \alpha' \cE'\) is \(c - \alpha' \cdot a = -\alpha\) by definition; we call this point \(P_0\). \(P_0\) is a point of tangency and hence has intersection multiplicity 2, because the centers of \(\alpha \cE\) and \(c + \alpha' \cdot \cE'\) both lie on the \(x\)-axis.

Next, we argue for the existence of two more distinct intersection points \(P_1'', P_2''\) as depicted in \Cref{fig:update_step_inner}.
The leftmost point of \(c + \alpha' \cE'\) is \((- \alpha, 0)\), and the rightmost point is \(c + \alpha'\), which by \Cref{lemma:update_params}-(\ref{item:claim_update_params_4}) is to the right of \((\alpha, 0)\), the rightmost point of \(\alpha \cE\).
Thus, by lying on \(\partial \alpha \cE\), \(P_1, P_2\) lie between the leftmost and rightmost points of \(c + \alpha' \cE'\), and so \(c + \alpha' \cE'\) intersects the line through \(P_1\) and \(P_2\).
Further, by \Cref{lemma:translate_ellipse_angle}, as \(c + \alpha' \cE'\) lies in the angle \(\angle P_1 \vv P_2\), \(c + \alpha' \cE'\) actually intersects the line segment \(\overline{P_1 P_2}\).
Observe that this intersection happens at two distinct points, which we call \(P_1'\) and \(P_2'\).
Both points are inside of \(\alpha \cE\), yet \(\partial (c + \alpha' \cE')\) is a continuous path that connects both to
the rightmost point of \(c + \alpha' \cE'\), which is outside of \(\alpha \cE\).
Thus \(\partial (c + \alpha' \cE')\) intersects \(\partial \alpha \cE\) at two more distinct points, which we call \(P_1''\) and \(P_2''\).

Now, we argue that \(P_1''\) and \(P_2''\) lie on the minor arc \(P_1 P_2\).
First, observe that the arc \(P_1 P_2\) containing \(P_0\) is the major arc. This is because \(P_1\) lies to the right of the \(y\)-axis, as determined in (\ref{eqn:usit_p1}); and by symmetry so does \(P_2\).
This also implies that major arc \(P_1 P_2\) is the arc with which the boundary of \(\conv{\alpha \cE \cup \{\vz\}}\) coincides.
\(P_1'\) and \(P_2'\) are collinear with \(P_1\) and \(P_2\), and as \(P_1''\) and \(P_2''\) are to the right of \(P_1'\) and \(P_2'\), this implies that they must lie on the minor arc \(P_1 P_2\).

Counting all the intersection points of \(\partial \alpha \cE\) and \(\partial(c + \alpha' \cE')\), we have \(P_0\) (with multiplicity \(2\))
and \(P_1''\) and \(P_2''\) (both with multiplicity 1); with total multiplicity 4.
Using \Cref{lemma:conic_5pt}, it is impossible for them to have another intersection point without both ellipses being the same.
Thus \(\partial(c + \alpha' \cE')\) cannot intersect the major arc \(P_1 P_2\) except at \(P_0\), and so except at points of tangency the boundary of \(c + \alpha' \cE'\) does not intersect the boundary of \(\conv{\alpha\cE \cup \{\vz\}}\).
\end{proof}

\subsection{Generalizing to high dimension and arbitrary previous ellipsoids}
\label{sec:gen_high_dim}
Now that we have demonstrated the invariants of \Cref{def:alg-invariant} for the special two-dimensional case where the previous ellipsoid is the unit ball, we generalize slightly to higher dimensions. However, we first still assume the previous ellipsoid is the unit ball.

Using the parameters as defined in (\ref{eqn:update_params}), we will let \(\cE = B_2^d\), and define the boundary of \(\cE'\) as
\[\frac{1}{a^2} (\vx_1 - c)^2 + \frac{1}{b^2} \vx_2^2 + \ldots + \frac{1}{b^2} \vx_d^2 = 1.\]
Observe that we can also write \(\cE' = \cE_{\mD}\) where \(\mD = \diag{\frac{1}{a^2}, \frac{1}{b^2}, \ldots, \frac{1}{b^2}}\).
Similarly to before, we let \(\vz = (c+a, 0, 0, \ldots, 0) \in \R^d\),
the furthest point of \(c + \cE'\) in the positive direction of the \(x_1\)-axis.

Now, we argue that the invariants of \cref{def:alg-invariant}
still hold in this setting.
\begin{lemma}\label{lemma:i_gen_dim1}
The inner and outer ellipsoid invariants hold in this setting:
\begin{enumerate}
    \item \(\cE \subseteq c \cdot \ve_1 + \cE'\)
    \item \(c \cdot \ve_1 + \alpha' \cE' \subseteq \conv{\alpha\cE \cup \{\vz\}}\)
\end{enumerate}
\end{lemma}
\begin{proof}
Observe that \(\cE\), \(c \cdot \ve_1 + \cE'\), \(c \cdot \ve_1 + \alpha' \cE'\), and \(\conv{\alpha \cE \cup \{\vz\}}\) are all bodies of revolution about the \(x_1\)-axis,
with their cross-sections given by their counterparts in \Cref{sec:two_dim_update}.
As \Cref{lemma:update_step_outer} and \Cref{lemma:update_step_inner}
hold for these cross sections, the set containments hold for
the bodies of revolution as well.
\end{proof}

We further generalize to the case where the previous ellipsoid is arbitrary.
In particular, let \(\vc^{\circ} + \cE\) be the previous ellipsoid, with a vector \(\vc^{\circ} \in \R^d\) and \(\cE = \cE_{\mA}\) for non-singular matrix \(\mA \in \R^{d \times d}\).
Let \(\vz^{\circ} \in \R^d\) be an arbitrary vector, representing the next point received.
We let \(\vu = \mA(\vz^{\circ} - \vc^{\circ})\), and \(\mW \in \R^{d \times d}\) be an orthogonal matrix with \(\vw = \frac{\vu}{\|\vu\|}\) as its first column (e.g. by using as its columns an orthonormal basis containing \(\vw\)).
We define the next outer ellipsoid as \(\vc^{\circ}+c \mA^{-1} \vw + \cE' \) for \(\cE' = \cE_{\mW \mD \mW^\top \mA}\), with \(\mD = \diag{\frac{1}{a^2}, \frac{1}{b^2}, \ldots, \frac{1}{b^2}}\) as before.
Observe that \(\vz = \vc^{\circ}+(c+a) \mA^{-1} \vw\) is the furthest point of \(\vc^{\circ}+c \mA^{-1} \vw + \cE'\) from the previous center \(\vc^{\circ}\) towards \(\vz^{\circ}\).

This setup works to preserve the key invariants, as we see in the next claim.
\begin{lemma}\label{lemma:i_gen_us}
   The inner and outer ellipsoid invariants hold in this setting:
    \begin{enumerate}
        \item \(\vc^{\circ} + \cE \subseteq \vc^{\circ}+c \mA^{-1} \vw + \cE' \)
        \item \(\vc^{\circ}+c \mA^{-1} \vw + \alpha' \cE' \subseteq \conv{(\vc^{\circ} + \alpha\cE) \cup \{\vz\}}\)
    \end{enumerate}
\end{lemma}
\begin{proof}
    We translate both set inclusions by \(- \vc^{\circ}\), then apply the nonsingular linear transformation \(\mW^\top \mA\).
    Observe that the set inclusions we wish to prove hold if and only if the transformed ones do.
    Noting that \(\mW^\top \mA \cE' = \cE_{\mW \mD}\), the transformed set inclusions are
    \(\cE_{\mW} \subseteq c \cdot \ve_1 + \cE_{\mW \mD}\)
    and \(c \cdot \ve_1 + \alpha' \cE_{\mW \mD} \subseteq \conv{\alpha \cE_{\mW} \cup \{(c + a) \cdot \ve_1\}}\).
    However, since \(\mW\) is an orthogonal matrix,
    \(\cE_{\mW} = B_2^d\) and \(\cE_{\mW \mD} = \cE_{\mD}\),
    and so the inclusions are exactly those shown in \Cref{lemma:i_gen_dim1}.
\end{proof}

Choosing \(\gamma\) correctly in (\ref{eqn:update_params})
ensures that \(\vz \in \vc^{\circ}+c \mA^{-1} \vw + \cE'\) coincides with \(\vz^{\circ}\),
as stated in the upcoming claim. This can be seen by looking at the definition of \(\vz\).
\begin{lemma}\label{lemma:i_gen_v}
    If \(\gamma\) is chosen so that \(c + a = \|\vu\|\), then \(\vz = \vz^{\circ}\).
\end{lemma}

\subsection{General algorithm}
\label{sec:gen_alg}
The goal of this section is to give and analyze a full algorithm that solves the streaming ellipsoid approximation problem, building on the analysis of the update rule from the previous sections.

Before we describe the complete algorithm, we give pseudocode in \Cref{alg:full_update} for its primary primitive. It is an update step like the one we analyzed in the previous section, \Cref{sec:gen_high_dim}.



\begin{algorithm}[h]
\caption{Full update step \(\cA^{\mathsf{full}}\)}\label{alg:full_update}

\textbf{input}: \(\mA_{t-1} \in \R^{d \times d}, \vc_{t-1} \in \R^d, \alpha_{t-1} \in [0, \frac{1}{2}], \vz_t \in \R^d\) \\
\textbf{output}: \(\mA_{t} \in \R^{d \times d}, \vc_{t} \in \R^{d}, \alpha_{t} \in [0, \alpha_{t-1}]\)
\begin{algorithmic}[1]
    \State Let \(\vu = \mA_{t-1}(\vz_t - \vc_{t-1})\), \(\vw = \frac{\vu}{\|\vu\|}\)
    \If{\(\|\vu\|_2 > 1\)}
        \State Let \(\gamma^{\star}_t\) be such that \(a(\gamma^{\star}_t) + c(\gamma^{\star}_t) = \|\vu\|\)\label{alg_line:compute_update_t}
        \State \(\hat{\mA} = \frac{1}{b(\gamma^{\star}_t)} \mI_d + \left(\frac{1}{a(\gamma^{\star}_t)} - \frac{1}{b(\gamma^{\star}_t)}\right) \vw \vw^\top\)\label{alg_line:ellipsoid_update}
        \State \Return \(\mA_t = \hat{\mA} \cdot \mA_{t-1}, \ \vc_t = \vc_{t-1} + c(\gamma^{\star}_t) \mA_{t-1}^{-1} \vw, \ \alpha_t = \alpha'(\gamma^{\star}_t)\) \label{alg_line:update_mat}
    \Else
        \State \Return \(\mA_t = \mA_{t-1}, \vc_t = \vc_{t-1}, \alpha_i = \alpha_{t-1}\)
        \label{alg_line:no_update}
    \EndIf
\end{algorithmic}
\end{algorithm}

In Lines \ref{alg_line:compute_update_t}, \ref{alg_line:ellipsoid_update} and \ref{alg_line:update_mat}, we use the definition of \(a(\gamma), b(\gamma), c(\gamma), \alpha'(\gamma)\) from (\ref{eqn:update_params}),
substituting \(\alpha_{t-1}\) for \(\alpha\).
Although the update step does not explicitly mention ellipsoids,
we use \(\cE_t = \cE_{\mA_t}\) so that at iteration \(t\) the next outer and inner ellipsoids are \(\vc_t + \cE_{\mA_t}\)
and \(\vc_t + \alpha_t \cE_{\mA_t}\), respectively.
If at this iteration \(\|\vu\|_2 \leq 1\), we will refer to this as the case where the ellipsoids are not updated, as is clear from Line \ref{alg_line:no_update}.

Observe also that if in iteration \(t\) we let \(\mW \in \R^{d \times d}\) be an orthogonal matrix with \(\vw\) as its first column,
we can write
\begin{equation}\label{eqn:alg_mat_update}
\hat{\mA} = \mW \cdot \diag{\frac{1}{a(\gamma^{\star}_t)}, \frac{1}{b(\gamma^{\star}_t)}, \cdots, \frac{1}{b(\gamma^{\star}_t)}} \cdot \mW^\top
\end{equation}

Now, we argue that this algorithm satisfies the invariants defined in \Cref{def:alg-invariant}.
This argument is essentially the observation that the update step in the algorithm is the one analyzed in \Cref{lemma:i_gen_us}.
\begin{lemma}\label{lemma:alg_correctness}
\Cref{alg:full_update} is a monotone update; i.e., it satisfies the invariants in \Cref{def:alg-invariant}.
\end{lemma}
\begin{proof}
If \(\|\vu\|_2 \leq 1\), then \(\vz_i \in \vc_n + \cE_n\) and the inner and outer ellipsoids are not updated, so the invariants clearly hold.
Otherwise, we apply \Cref{lemma:i_gen_us} and \Cref{lemma:i_gen_v} setting \(\mA = \mA_{t-1}, \vc^{\circ} = \vc_{t-1}, \vz^{\circ} = \vz_{t}, \alpha=\alpha_{t-1}\).
Using (\ref{eqn:alg_mat_update}), \(\cE_{\mA_t}\) is the same as  \(\cE'\) in \Cref{lemma:i_gen_us}; and clearly \(\alpha_{t} = \alpha'\).
This establishes the inner ellipsoid invariant \(\vc_{t} + \alpha_{t} \cE_{t} \subseteq \conv{(\vc_{t-1} + \alpha_{t-1} \cE_{t-1}) \cup \{\vz_{t}\}}\) directly.
To show \(\conv{(\vc_{t-1} + \cE_{t-1}) \cup \{\vz_{t}\}} \subseteq \vc_{t} + \cE_{t}\),
observe that we have \(\vc_{t-1} + \cE_{t-1} \subseteq \vc_{t} + \cE_{t}\) from \Cref{lemma:i_gen_us},
and \(\vz_{t} \in \vc_{t} + \cE_{t}\) from \Cref{lemma:i_gen_v}.
Then the outer ellipsoid invariant follows as \(\vc_{t} + \cE_{t}\) is a convex set.
\end{proof}

Finally, we bound the relevant quantities that will be used in the analysis of the full algorithm's approximation factor.
In particular, we show that \(\exp(\gamma^{\star}_t)\) gives a lower bound on the increase in volume at each iteration \(t\).
If \(\|\vu\|_2 \leq 1\), and the ellipsoids are not updated, in that iteration we think of \(\gamma^{\star}_t = 0\).
\begin{lemma}\label{lemma:vol_t}
For any input given to \Cref{alg:full_update}, we have  \(\vol(\cE_i) \geq \exp(\gamma^{\star}_t) \vol(\cE_{t-1})\).
\end{lemma}
\begin{proof}
This formula is clearly true when the ellipsoids are not updated because \(\gamma^{\star}_t = 0\), so we consider the nontrivial case.
Recall the formula \(\vol(\cE_\mA) = \det(\mA^{-1}) \vol(B_2^d)\) from \Cref{lemma:prelim_vol}. Then we have
\[\vol(\cE_{\mA_i}) = \det(\mA_i^{-1}) \vol(B_2^d) = \det(\hat{\mA}^{-1}) \cdot \det(\mA_{t-1}^{-1}) \cdot \vol(B_2^d) = \det(\hat{\mA}^{-1}) \vol(\cE_{\mA_{t-1}})\]
where we use the definition of \(\hat{\mA}\) from Line \ref{alg_line:ellipsoid_update} on the \(t\)-th iteration.
Then \begin{align*}
    \det(\hat{\mA}^{-1}) &= a(\gamma^{\star}_t) \cdot b(\gamma^{\star}_t)^{d-1} & \text{using (\ref{eqn:alg_mat_update})}\\
    &\geq a(\gamma^{\star}_t) & \text{by \Cref{lemma:update_params}-(\ref{item:claim_update_params_2})} \\
    &= \exp(\gamma^{\star}_t) & \text{by definition of \(a\) in (\ref{eqn:update_params})}
\end{align*}
and using \(\vol(\cE_{\mA_i}) = \det(\hat{\mA}^{-1}) \cdot \vol(\cE_{\mA_{t-1}})\) completes the proof.
\end{proof}

We are now ready to present the complete algorithm in \Cref{alg:main}.
The algorithm is explicitly given \(\vc_0 + r_0 \cdot B_2^d \subseteq Z\). For simplicity, here, we say \(r = r_0\).
Let \(R = R(Z)\). While the final approximation factor depends on this quantity, the algorithm is not given it.
Note that \(\kappa(Z) \leq \nfrac{R}{r}\), so the quality of the approximation depends not only on \(\kappa(Z)\), but also on how well the given ball \(\vc_0 + r \cdot B_2^d\) is centered within \(Z\).

This algorithm proceeds in two phases. It begins with a ``local'' first phase, where the inner ellipsoid is a ball kept at radius \(r\), and the outer ellipsoid is a ball scaled to contain all the points. For readability, the variables of the algorithm in this phase are annotated with a superscript \(^{(l)}\). The second phase starts if the approximation factor of the first phase ever reaches \(\alpha^{(l)} \leq \frac{1}{d \log d}\),
at which point the algorithm uses the ``full'' update that was just described in \Cref{alg:full_update}.
We use two phases because while the full update reaches a near-optimal approximation factor when \(\nfrac{R}{r} \geq d \log d\), the local phase using balls does better when \(\nfrac{R}{r} \leq d \log d\).
While we cannot tell when to switch phases exactly (this would require knowing \(\nfrac{R}{r}\)), we show that it is enough to approximate the aspect ratio during the first phase up to a constant factor.

\begin{algorithm}[h]
\caption{Streaming ellipsoid rounding -- complete algorithm}\label{alg:main}

\textbf{input}: \(\vc_0+r B_2^d \subseteq Z\) \\
\textbf{output}: \(\vc_n + \cE_n, \vc_n + \alpha_n \cdot \cE_n\)
\begin{algorithmic}[1]
\State Initialize \(\mA^{(l)}_0 = \frac{1}{r} \mI_d, \vc_0^{(l)} = \vc_0, \alpha_0^{(l)} = 1\)
\State \({t^{(l)}} = 0, R_0 = 0\)
\While{\({t^{(l)}} \leq n\)}
\Comment{Phase I: Local update step that maintains a ball}
    \State Receive point \(\vz_{t^{(l)}}\)
    \If{\(\|\vz_{t^{(l)}} - \vc_0\|_2 \leq r \cdot d \log d\)}\label{alg_line:phaseI_if}
        \If{\(\|\vz_{t^{(l)}} - \vc_0\|_2 > R_{{t^{(l)}}-1}\)}
            \State \(\mA_{t^{(l)}}^{(l)} = \frac{1}{\|\vz_{t^{(l)}} - \vc_0\|_2} \cdot \mI_d, \vc_{t^{(l)}}^{(l)} = \vc_{{t^{(l)}}-1}^{(l)}, \alpha_{t^{(l)}}^{(l)} = \frac{r}{\|\vz_{t^{(l)}} - \vc_0\|_2}\)
            \Comment{Grow the ball to contain \(\vz_{t^{(l)}}\)}
            \State \(R_{t^{(l)}} = \frac{\|\vz_{t^{(l)}} - \vc_0\|_2}{r}\)
        \Else
        \State \(\mA_{t^{(l)}}^{(l)} = \mA_{{t^{(l)}}-1}, \vc_i^{(l)} = \vc_{{t^{(l)}}-1}^{(l)}, \alpha_{t^{(l)}}^{(l)} = \alpha_{{t^{(l)}}-1}\)
        \State \(R_{t^{(l)}} = R_{{t^{(l)}}-1}\)
        \EndIf
    \Else{}
        \State \textbf{break} 
        \Comment{Break the loop and jump to Line \ref{alg_line:end_phase1}}
    \EndIf
    \State \({t^{(l)}} = {t^{(l)}} + 1\)
\EndWhile
\If{\({t^{(l)}} > n\)}\label{alg_line:end_phase1}
    \Comment{If we stayed in Phase I for the entire execution of the algorithm}
    \State \Return \(\mA_n^{(l)}, \vc_n^{(l)}, \alpha_n^{(l)}\)
\EndIf
\State \(t_{s} = {t^{(l)}}\) \label{alg_line:tswitch}
\Comment{Point \(\vz_{t_{s}}\) has not yet been processed}
\State \(\mA_{t_{s}-1} =\frac{1}{r d \log d} \cdot \mI_d, \vc_{t_{s}-1} = \vc_{t_{s}-1}^{(l)}, \alpha_{t_{s}-1} = \frac{1}{d \log d}\)\label{alg_line:phase_trans}
\Comment{Transition: grow the ball to maximium size}
\For{\(t \in \{t_{s}, t_{s}+1, \ldots, n\}\)}
    \Comment{Phase II: full update for the remaining points}
    \State Receive point \(\vz_i\)
    \State \(\mA_i, \vc_i, \alpha_i = \cA^{\mathsf{full}}(\mA_{t-1}, \vc_{t-1}, \alpha_{t-1}, \vz_i)\)\label{alg_line:phaseII_update}
\EndFor
\State \Return \(\mA_n, \vc_n, \alpha_n\)
\end{algorithmic}
\end{algorithm}

Before Line \ref{alg_line:end_phase1}, the algorithm executes the first phase that has the outer and inner ellipsoids as balls.
In Line \ref{alg_line:end_phase1}, we have \({t^{(l)}} > n\) if the algorithm stayed in Phase I for every point, i.e. we had \(\max_{1 \leq {t^{(l)}} \leq n} \|_2 \vz_{t^{(l)}} - \vc_0 \|_2 \leq r \cdot d \log d\). In this case, the algorithm returns the approximation maintained by Phase I.
Otherwise we must have come across a point where \(\|\vz_{t^{(l)}} - \vc_0 \|_2 > r \cdot d \log d\), and the algorithm proceeds with Phase II.
We let \(t_{s}\) in Line \ref{alg_line:tswitch} mark the point received that causes the algorithm to proceed to Phase II.
We then perform a ``transition'' on Line \ref{alg_line:phase_trans} that grows the ball of Phase I to its maximum size. This transition step makes the analysis of the complete algorithm easier, as then the starting approximation for the second phase is exactly \(\alpha_{t_s - 1} = \frac{1}{d \log d}\).
Then the algorithm runs the full update \(\cA^{\mathsf{full}}\) for the rest of the points, including \(\vz_{t_{s}}\).
For simplicity, we write our algorithm so that it `receives' \(\vz_{t_{s}}\) twice, once for each phase.
However, the first phase does not commit to an update for this point, and the ellipsoids in Line \ref{alg_line:phase_trans} are not committed either; the algorithm does not commit to an update for this point until Line \ref{alg_line:phaseII_update}.

Recall the approximation guarantee stated in \Cref{thm:main_one}:
\begin{equation}\label{eq:main_approx_guarantee}
    \frac{1}{\alpha_n} \leq O(\min\inbraces{\nfrac{R}{r}, d \log\inparen{\nfrac{R}{r}}})
\end{equation}

We can interpret the approximation guarantee \eqref{eq:main_approx_guarantee} by cases depending on if \(\nfrac{R}{r} \geq d \log d\) (i.e. if the algorithm ever enters the second phase):
\begin{lemma}
We have for all \(d \geq 2\) that
\[\min\inbraces{\nfrac{R}{r}, d \log\inparen{\nfrac{R}{r}}} = \Theta\inparen{\begin{cases}
    d \log\inparen{\nfrac{R}{r}} & \text{if } \nfrac{R}{r} > d \log d \\
    \nfrac{R}{r} & \text{if } \nfrac{R}{r} \leq d \log d
\end{cases}}.\]
\end{lemma}

Now, we claim a straightforward geometric fact -- that the distance of the furthest \(\vz_{t}\) to \(\vc_0\)
approximates the circumradius of \(Z\) up to a constant factor.
We will use this to show that Line \ref{alg_line:phaseI_if} will be able to properly detect when
\(\nfrac{R}{r} > d \log d\) (again, up to a constant factor).
\begin{lemma}\label{lemma:alg_r_est}
    Let \(\vc_0 + r_0 B_2^d \subseteq Z\), and 
    \(R = R(Z)\). Then,
    \[R \leq \max_{1 \leq t^{(l)} \leq n} \|\vc_0 - \vz_{t^{(l)}} \|_2 \leq 2 \cdot R.\]
\end{lemma}
\begin{proof}
    For the left inequality, observe that if we let \(r_{\max} = \max_{1 \leq t^{(l)} \leq n} \|\vc_0 - \vz_{t^{(l)}} \|_2\),
    then \(Z \subseteq \vc_0 + r_{\max} \cdot B_2^d\).
    For the right inequality, observe that for any containing ball \(\vc' + R' \cdot B_2^d \supseteq Z\),
    its diameter is \(2 R'\).
    But as \(\vc' + R' \cdot B_2^d\) contains \(\vc_0\) and \(\vz_{1}, \ldots, \vz_{n}\), we must have \(\diam{\vc' + R' \cdot B_2^d} \geq \diam{\{\vc_0\} \cup \{\vz_1, \ldots, \vz_n\}}\) and so \(2R' \geq r_{\max}\).
\end{proof}

Next, we discuss the approximation guarantee that the algorithm achieves, depending on the phase that it terminates with.
We start with if the algorithm only stays in the local phase, in which case we can readily apply the previous claim.
\begin{lemma}\label{lemma:approx_local}
    If \Cref{alg:main} never enters Phase II, then its approximation guarantee satisfies
    \(\frac{1}{\alpha_n} \leq \nfrac{2 R}{r}\).
\end{lemma}
\begin{proof}
    At the termination of Phase I, the algorithm produces approximation \(\alpha_n^{(l)} = \max_{1 \leq {t^{(l)}} \leq n} \frac{\|\vz_{t^{(l)}} - \vc_0\|_2}{r}\).
    Using \Cref{lemma:alg_r_est} we obtain \[
    \frac{1}{\alpha_n^{(l)}} = \max_{1 \leq {t^{(l)}} \leq n} \frac{\|\vz_{t^{(l)}} - \vc_0\|_2}{r} \leq \frac{2 R}{r}
    ,\]
    as desired.
\end{proof}

The analysis in the case where the algorithm enters the full phase is more involved.
We use \Cref{lemma:vol_t}, which shows that the increase in approximation factor each iteration is not too large compared to the increase in volume, to bound \(\frac{1}{\alpha_n}\).
We know that the volume of the final ellipsoid \(\vc_n + \cE_n\) must be bounded relative to \(R \cdot B_2^d\), as the algorithm produces \(\vc_n + \alpha_n \cdot \cE_n \subseteq Z\). However, this leads to an upper bound that is still a function of \(\frac{1}{\alpha_n}\).
\begin{lemma}\label{lemma:alpha_rel}
    If \Cref{alg:main} enters Phase II, the approximation guarantee satisfies
    \[\frac{1}{\alpha_n} \leq 2 \left(d \log\left(\frac{1}{\alpha_n}\right) + d \log\left(\frac{R}{r}\right)\right).\]
\end{lemma}
\begin{proof}
The algorithm transitions to Phase II at Line \ref{alg_line:tswitch}, starting at iteration \(t_{s}\).
At each subsequent iteration, we claim that \Cref{alg:full_update} guarantees \(\frac{1}{\alpha_{t}} = \frac{1}{\alpha_{t-1}} + 2 \gamma^{\star}_t\).
By \Cref{lemma:update_params}-(\ref{item:claim_update_params_1}), we have for all \(t_{s} \leq t \leq n-1\) where the ellipsoids were updated that \(\frac{1}{\alpha_{t}} = \frac{1}{\alpha_{t-1}} + 2 \gamma^{\star}_t\).
When the ellipsoids are not updated, this still holds, as in that case \(\gamma^{\star}_t = 0\). 

As in Phase II the algorithm begins with \(\alpha_{t_{s} - 1} = \frac{1}{d \log d}\), we have
\begin{equation}\label{eqn:alpha_rel_t}
\frac{1}{\alpha_n} = d \log d + 2 \sum_{t=t_{s}}^{n-1} \gamma^{\star}_t.
\end{equation}
Now applying \Cref{lemma:vol_t} for each \(t\), we have \(\vol(\cE_n) \geq \exp\left(\sum_{t=t_{s}}^{n-1} \gamma^{\star}_t\right) \cdot \vol(\cE_{t_{s}})\). Taking logarithms gives
\begin{equation}\label{eq:ub_t_vol}
\log\inparen{\frac{\vol\inparen{\cE_n}}{\vol\inparen{\cE_{t_{s} - 1}}}} \geq \sum_{t=t_{s}}^{n-1} \gamma^{\star}_t.
\end{equation}

Recall that \(\vc_0 + r \cdot B_2^d \subseteq Z\),
and by \Cref{def:circumradius},
\(Z \subseteq \vc_c + R \cdot B_2^d\) for some center \(\vc_c\).
By \Cref{lemma:alg_correctness}, we have \(\vc_n + \alpha_n \cdot \cE_n \subseteq Z\), so
that \(\vol(\cE_n) \leq \frac{1}{\alpha_n^d} \cdot \vol(R \cdot B_2^d)\).
As in Phase II we start with \(\cE_{t_s - 1} = \vc_0 + r d \log d \cdot B_2^d\), this yields
\begin{align*}\sum_{t=t_{s}}^{n-1} \gamma^{\star}_t &\leq \log\inparen{\frac{\vol(\cE_n)}{\vol (\cE_{t_{s}-1})}} &\text{by \eqref{eqn:alpha_rel_t}} \\
&\leq d \log \inparen{\frac{1}{\alpha_n}} + \log\inparen{\frac{\vol(R \cdot B_2^d)}{\vol(r d \log d \cdot B_2^d)}} & \text{by } \vol(\cE_n) \leq \frac{1}{\alpha_n^d} \vol(R \cdot B_2^d)\\
&= d \log\inparen{\frac{1}{\alpha_n}} + d\log\inparen{\frac{R}{rd \log d}}\\
&\leq d \log \inparen{\frac{1}{\alpha_n}} + d \log \inparen{\frac{R}{r}} - d \log d.
\end{align*}
Plugging into (\ref{eqn:alpha_rel_t}) completes the proof of \Cref{lemma:alpha_rel}.
\end{proof}

Intuitively, \(x \leq a + b \cdot \log x\) for some constants \(a, b > 0\) can only be true for bounded \(x\), as \(x = \omega(\log x)\).
As we showed \(1/\alpha_n\) satisfies a relation like this in \Cref{lemma:alpha_rel}, we develop this intuition to give a quantitative upper bound on \(1/\alpha_n\).
\begin{lemma}\label{lemma:approx_full}
If \Cref{alg:main} enters Phase II, then we have \[\frac{1}{\alpha_n} \leq 8 d(\log d + \log \nfrac{R}{r}).\]
\end{lemma}
\begin{proof}
Assume towards contradiction that \(\frac{1}{\alpha_n} > 8 d (\log d + \log \nfrac{R}{r})\).
Observe then that \(\frac{1}{\alpha_n} - \frac{3}{4} \cdot\frac{1}{\alpha_n} > 2 d (\log d + \log \nfrac{R}{r})\). Using \Cref{lemma:alpha_rel}, we have
\[2 (d \log \nfrac{1}{\alpha_n} + d \log \nfrac{R}{r}) \geq \frac{1}{\alpha_n} > 2 (d \log d + d \log \nfrac{R}{r}) + \frac{3}{4} \cdot \frac{1}{\alpha_n}\]
Simplifying the above inequality gives \(2d \log \nfrac{1}{d \cdot \alpha_n} > \frac{3}{4} \cdot \frac{1}{\alpha_n}\), i.e. 
\(2\log \nfrac{1}{d \cdot \alpha_n} > \frac{7}{8} \cdot \frac{1}{d \cdot \alpha_n}\).
It is clear that this is impossible by looking at the graph of the function \(x \mapsto 2 \log x - \frac{3}{4} x\), which is concave with a maximum of \(2(\log(\nfrac{8}{3}) - 1) < 0\).
\end{proof}

Now we combine the previous claims to prove the guarantees of \Cref{alg:main} and obtain \Cref{thm:main_one}.

\begin{proof}[Proof of \Cref{thm:main_one}.] We first discuss the approximation guarantee and correctness, then the memory and runtime complexity of \Cref{alg:main}.

\paragraph{Approximation guarantee}
We break the analysis of the approximation guarantee by cases, depending on the aspect ratio.
If \(\nfrac{R}{r} \leq \frac{1}{2} d \log d\),
then by \Cref{lemma:alg_r_est} we have \(\max_{1 \leq t^{(l)} \leq n} \|\vc_0 - \vz_{t^{(l)}} \|_2 \leq r d \log d\),
and the algorithm never enters Phase II.
By \Cref{lemma:approx_local}, the final approximation factor is \(\nfrac{2R}{r}\).
If \(\nfrac{R}{r} > d \log d\), then by \Cref{lemma:alg_r_est} we have \(\max_{1 \leq t^{(l)} \leq n} \|\vc_0 - \vz_{t^{(l)}} \|_2 > r d \log d\), and the algorithm must enter Phase II.
Then \Cref{lemma:approx_full} applies, and the final approximation factor is \(O(d (\log d + \log \nfrac{R}{r}) = O(d \log \nfrac{R}{r})\).

If \(\frac{1}{2} d \log d < \nfrac{R}{r} \leq d \log d\), then it is possible for the algorithm to never enter Phase II or for it to enter Phase II. Either way, we argue that the final approximation factor is \(\frac{1}{\alpha_n} \leq O\inparen{\nfrac{R}{r}}\).
If it does not enter Phase II, then by \Cref{lemma:approx_local}, the approximation guarantee we get is \(\frac{1}{\alpha_n} \leq O(\nfrac{R}{r})\).
If it does enter Phase II, then by \Cref{lemma:approx_full}
we have \[
\frac{1}{\alpha_n} \leq O(d \log d + d \log \nfrac{R}{r})
\]
Due to the assumption that \(\frac{1}{2} d \log d < \nfrac{R}{r} \leq d \log d\), we also have in this case that \(\frac{1}{\alpha_n} \leq O(\nfrac{R}{r})\).

\paragraph{Correctness} By \Cref{lemma:alg_invariant}, to argue that the algorithm solves \Cref{prob:main} it is enough to show that it is monotone, i.e. it satisfies the invariants of \Cref{def:alg-invariant}.
It is clear that the local update in Phase I satisfies the invariants, as the outer ellipsoid is a ball of growing radius
and the inner ellipsoid is kept to the ball of radius \(r\). 
It is also clear that after the algorithm transitions to Phase II, all the full updates are monotone by \Cref{lemma:alg_correctness} and the fact that the starting approximation factor for this phase is is \(\alpha_{t_s - 1} = \frac{1}{d \log d} \leq \frac{1}{2}\).
As algorithm transitions to Phase II,
observe that on Line \ref{alg_line:phase_trans} the radius of the outer ellipsoid grows again to \(r d \log d\) before applying the full update,
so the first first full update of Phase II is also monotone.

\paragraph{Memory and runtime complexity}
The memory complexity of the algorithm is \(O(d^2)\).
Observe that \Cref{alg:full_update} only stores a constant number of matrices in \(\R^{d \times d}\), vectors in \(\R^d\), or constants,
so its memory complexity is \(O(d^2)\). It is only instantiated once for each point received in Phase II, so the memory complexity in this phase \(O(d^2)\). Finally, the memory complexity in the first phase is also \(O(d^2)\) because it stores the same kind of quantities as \Cref{alg:full_update}.

To show the runtime of the algorithm is \(\widetilde{O}(nd^2)\), we show that the runtime to process each next point is at most \(\widetilde{O}(d^2)\). This is clear in Phase I, and during the transition to Phase II.
For the full update this is less clear, as \Cref{alg:full_update} uses both \(\mA_{t-1}\) and \(\mA_{t-1}^{-1}\)
which naively would require inverting a matrix on each iteration.
However, if we represent \(\mA\) using the SVD (see the next section and \Cref{lemma:full_update_eff}), we can implement the update in \(\widetilde{O}(d^2)\) time.
This would require that \(\mA_{t_s - 1}\) be given in SVD form as well for the first full update, but it is already in that form as a scaled identity matrix.

Put together, these complete the proof of \Cref{thm:main_one}.
\end{proof}

\subsubsection{Efficient implementation of the full update step}

In this section, we use a method similar to that in Algorithm 2 from \cite{mmo22} to show that the full update step can be implemented in \(\widetilde{O}(d^2)\) time.
In particular, we use the same subroutine \(\textsc{SVDRankOneUpdate}\) with signature
\begin{align}
    (\mU', \mSigma', \mV') = \textsc{SVDRankOneUpdate}((\mU, \mSigma, \mV), \vy_1, \vy_2)\label{eq:stange_update}
\end{align}
where the result \(\mU' \mSigma' (\mV')^\top\) is the SVD of the matrix \(\mU \mSigma \mV^\top + \vy_1 \vy_2^\top\).
\citet{stange08} shows that this procedure be done in \(O(d^2 \log d)\) time.
We rewrite \Cref{alg:full_update} in \Cref{alg:full_update_eff} to make it clear how to use the SVD representation and the efficient rank-1 update to efficiently implement the full update.
One can readily see that \Cref{alg:full_update_eff} has the exact same behavior as \Cref{alg:full_update}, and so gives the same approximation and correctness guarantees.

\begin{algorithm}[h]
\caption{Efficient full update step \(\cA^{\mathsf{full}}\)}\label{alg:full_update_eff}

\textbf{input}: \((\mU_{t-1}, \mSigma_{t-1}, \mV_{t-1}) \in \R^{d \times d}, \vc_{t-1} \in \R^{d}, \alpha_{t-1} \in [0, \frac{1}{2}], \vz_{t} \in \R^{d}\) \\
\textbf{output}: \((\mU_{t}, \mSigma_{t}, \mV_{t}) \in \R^{d \times d}, \vc_{t} \in \R^{d}, \alpha_{t} \in [0, \alpha_{t}]\)
\begin{algorithmic}[1]
    \State Let \(\vu = \mU_{t-1} \mSigma_{t-1} \mV_{t-1}^\top (\vz_{t} - \vc_{t-1})\), \(\vw = \frac{\vu}{\|\vu\|_2}\)
    \If{\(\|\vu\|_2 > 1\)}
        \State Let \(\gamma^{\star}_t\) be such that \(a(\gamma^{\star}_t) + c(\gamma^{\star}_t) = \|\vu\|\)\label{alg_line:compute_update_t_eff}
        \State \(\vy_1 = \left(\frac{1}{a(\gamma^{\star}_t)} - \frac{1}{b(\gamma^{\star}_t)}\right) \vw, \vy_2 = \mV_{t-1} \mSigma_{t-1} \mU_{t-1}^\top \vw\)
        \State \((\mU_{t}, \mSigma_{t}, \mV_{t}) = \textsc{SVDRankOneUpdate}((\mU_{t-1}, \frac{1}{b(\gamma^{\star}_t)} \mSigma_{t-1}, \mV_{t-1}), \vy_1, \vy_2)\)
        \State \Return \((\mU_{t}, \mSigma_{t}, \mV_{t}), \ \vc_{t} = \vc_{t-1} + c(\gamma^{\star}_t) \mV_{t-1} \mSigma_{t-1}^{-1} \mU_{t-1}^\top \vw, \ \alpha_t = \alpha'(\gamma^{\star}_t)\)
    \Else
        \State \Return \((\mU_{t}, \mSigma_{t}, \mV_{t}) = (\mU_{t-1}, \mSigma_{t-1}, \mV_{t-1}), \vc_t = \vc_{t-1}, \alpha_t = \alpha_{t-1}\)
    \EndIf
\end{algorithmic}
\end{algorithm}

\begin{remark}\label{remark:compute_update_t_eff}
    We briefly explain why Line \ref{alg_line:compute_update_t_eff}, finding \(\gamma^{\star}\) such that \(a(\gamma^{\star}) + c(\gamma^{\star}) = \|\vu\|\)  can be implemented efficiently. This is a one-dimensional optimization problem, and \(\gamma \mapsto a(\gamma) + c(\gamma)\) using \(a, c\) as defined in \eqref{eqn:update_params} is monotone increasing, 
    so finding an approximate \(\gamma^{\star}\) can be done efficiently with binary search.
    In particular, we can choose \(\gamma^{\star}\) to be a slight overestimate so the update is still monotone after slightly increasing \(\alpha_t\).
    This does not affect the final approximation guarantee beyond constant factors.
\end{remark}

This algorithm performs a constant number of taking norms of vectors, matrix-vector products, and algebraic operations; as well as one rank-one SVD update. As explained in \Cref{remark:compute_update_t_eff}, finding \(\gamma_i^*\) can also be done in effectively constant time. Thus for our runtime guarantee, we have:
\begin{lemma}\label{lemma:full_update_eff}
    \Cref{alg:full_update_eff} runs in time \(O(d^2 \log d)\).
\end{lemma}

\subsection{Fully-online asymmetric ellipsoidal rounding algorithm}

In this subsection, we prove Theorem \ref{thm:main_two}. See Algorithm \ref{alg:fully_online_approx}.

\begin{algorithm}
\caption{Fully online asymmetric ellipsoidal rounding}\label{alg:fully_online_approx}
\begin{algorithmic}[1]
\State \textbf{Input:} Stream of points $\vz_t$; 
 monotone update rule $\cA$ (Definition \ref{def:alg-invariant}) that takes as input the previous ellipsoid matrix $\mA$, center $\vc$, approximation factor $\alpha$, and update point $\vz$ and outputs the next ellipsoid matrix $\mA'$, center $\vc'$, and approximation factor $\alpha'$.
\State \textbf{Output:} Ellipsoid $\cE$, center $\vc$, and scale $\alpha \in (0,1)$ such that $\vc + \alpha \cdot \cE \subseteq \mathsf{conv}\inparen{\inbraces{\vz_1,\dots,\vz_n}} \subseteq \vc +\cE$.
\State Receive $\vz_1$; set $\mA = \mI_d$, $d_1 = 1$, $\vc_1 = \vz_1$, $\alpha_1 = 1$.
\For{$t = 2,\dots,n$}
    \State Receive $\vz_{t}$.
    \If{$\vz_{t} - \vc_{t-1} \notin \vspan{\vz_{1} - \vc_{t-1},\dots,\vz_{t-1}-\vc_{t-1}}$}\Comment{Irregular update step.}\label{line:irregular_step}
        \State Let $\vv_1,\dots,\vv_{d_{t-1}}$ be the singular vectors of $\mA$ corresponding to the semiaxes of $\cE_{t-1}$.
        \State Let $d_t = d_{t-1}+1$.
        \State Let $\vz_{d_{t}}' \coloneqq \frac{\vz_t - \sum_{i=1}^{d_{t-1}} \vv_i\ip{\vv_i,\vz_t}}{\norm{\vz_t - \sum_{i=1}^{d_{t-1}} \vv_i\ip{\vv_i,\vz_t}}_2}$.
        \State Let $\mM \coloneqq \mI_d - \frac{1}{\ip{\vv_{d_{t}}',\vz}}\cdot\inparen{\vz_t - \sqrt{1+2\alpha_{t-1}} \cdot \vv_{d_{t}}'}(\vv_{d_{t}}')^T$.
        \State Update $\mA_t \gets \mA_{t-1}\mM$.\Comment{Use \eqref{eq:stange_update} of \citet{stange08} to update $\vv_1,\dots,\vv_d$.}
        \State Update $\vc_t = \frac{\alpha_{t-1}}{1+2\alpha_{t-1}}\cdot\vz_t + \inparen{1-\frac{\alpha_{t-1}}{1+2\alpha_{t-1}}}\cdot\vc_{t-1}$.
        \State Update $\nfrac{1}{\alpha_t} \gets \nfrac{1}{\alpha_{t-1}}+1$.
    \Else  
        \State $\mA_t, \vc_t, \alpha_t = \cA(\mA_{t-1}, \vc_{t-1}, \alpha_{t-1}, \vz_{t})$
        \State $d_t \gets d_{t-1}$.
    \EndIf
\EndFor
\State \textbf{Output:} $(\vc_n, \cE_n, \alpha_n)$.
\end{algorithmic}
\end{algorithm}


To prove Theorem \ref{thm:main_two}, we need to show that our \textit{irregular update step} (a timestep $t$ when we have to update the dimensionality of our ellipsoid $\cE_{t-1}$ -- see Line \ref{line:irregular_step} of Algorithm \ref{alg:fully_online_approx}) still maintains the invariants we desire (Definition \ref{def:alg-invariant}).

Our plan is to first consider the special case of the irregular update where the new point to cover is conveniently located with respect to our current ellipsoids. We will see later that this special case is nearly enough for us to conclude the proof.

\begin{lemma}
\label{lemma:irregular_update_easy}
Let $Z \subset \R^{d}$ be a convex body where $Z$ lies in $\vspan{\vv_1,\dots,\vv_{d'}}$ for $d' < d$. For $0 < \alpha \le 1$, suppose we have
\begin{align*}
    \alpha \cdot \inbraces{\vz \in \vspan{\vv_1,\dots,\vv_{d'}} \suchthat \norm{\vz}_2 \le 1} \subseteq Z \subseteq \inbraces{\vz \in \vspan{\vv_1,\dots,\vv_{d'}} \suchthat \norm{\vz}_2 \le 1}.
\end{align*}
Then, for any $\vv_{d'+1}$ such that $\ip{\vv_i, \vv_{d'+1}} = 0$ for all $i \in [d']$ and for which
\begin{align*}
    \cE' &\coloneqq \inbraces{\vz \in \vspan{\vv_1,\dots,\vv_{d'+1}} \suchthat \norm{\vz}_2 \le \frac{1+\alpha}{\sqrt{1+2\alpha}}} \\
    \vc &\coloneqq \frac{\alpha}{\sqrt{1+2\alpha}} \cdot \vv_{d'+1}
\end{align*}
we have
\begin{align*}
    \vc + \frac{1}{1 + \nfrac{1}{\alpha}} \cdot \cE' \subseteq \mathsf{conv}\inparen{Z \cup \inbraces{\sqrt{1+2\alpha} \cdot \vv_{d'+1}}} \subseteq \vc + \cE'.
\end{align*}
\end{lemma}
\begin{proof}[Proof of \Cref{lemma:irregular_update_easy}]
We will show that the pair of ellipsoids given below satisfy the conditions promised by the statement of \Cref{lemma:irregular_update_easy}.
\begin{align}
    &\inbraces{\vz \in \vspan{\vv_1,\dots,\vv_{d'+1}} \suchthat \norm{\vz}_2 \le \frac{1+\alpha}{\sqrt{1+2\alpha}}} +  \frac{\alpha}{\sqrt{1+2\alpha}} \cdot \vv_{d'+1}\label{eq:irregular_outer} \\
    &\inbraces{\vz \in \vspan{\vv_1,\dots,\vv_{d'+1}} \suchthat \norm{\vz}_2 \le \frac{1+\alpha}{\sqrt{1+2\alpha}}} \cdot \frac{\alpha}{1+\alpha} +  \frac{\alpha}{\sqrt{1+2\alpha}} \cdot \vv_{d'+1}\label{eq:irregular_inner}
\end{align}
Clearly, the two ellipsoids given above are apart by a factor of $\nfrac{1+\alpha}{\alpha} = \nfrac{1}{\alpha}+1$, which means the approximation factor increases by exactly $1$ as a result of this update. It now suffices to show that the ellipsoid described by (\ref{eq:irregular_outer}) contains $\mathsf{conv}\inparen{B_2^{d'} \cup \inbraces{\sqrt{1+2\alpha} \cdot \vv_{d'+1}}}$ and that the ellipsoid described by (\ref{eq:irregular_inner}) is contained by the cone whose base is $\alpha \cdot B_2^{d'}$ and whose apex is $\sqrt{1+2\alpha}\cdot\vv_{d'+1}$.

For the first part, it suffices to verify that every point $\vz \in Z$ and $\sqrt{1+2\alpha} \cdot \vv_{d'+1}$ is contained by (\ref{eq:irregular_outer}). We give both the calculations below, from which the result for (\ref{eq:irregular_outer}) follows.
\begin{align*}
    \vz \in Z :&\quad\quad \norm{\vz - \frac{\alpha}{\sqrt{1+2\alpha}} \cdot \vv_{d'+1}}_2 = \sqrt{\norm{\vz}_2^2 + \frac{\alpha^2}{1+2\alpha}} \le \frac{1+\alpha}{\sqrt{1+2\alpha}} \\
    \vz = \sqrt{1+2\alpha} \cdot \vv_{d'+1} :&\quad\quad \norm{\vz - \frac{\alpha}{\sqrt{1+2\alpha}} \cdot \vv_{d'+1}}_2 = \sqrt{1+2\alpha} - \frac{\alpha}{\sqrt{1+2\alpha}} = \frac{1+\alpha}{\sqrt{1+2\alpha}}
\end{align*}

We now analyze (\ref{eq:irregular_inner}). Our task is to show the below inclusion.
\begin{align*}
    &\inbraces{\vz \in \vspan{\vv_1,\dots,\vv_{d'+1}} \suchthat \norm{\vz - \frac{\alpha}{\sqrt{1+2\alpha}}\cdot \vv_{d'+1}}_2 \le \frac{\alpha}{\sqrt{1+2\alpha}}}\\
    \subseteq &\mathsf{conv}\inparen{\alpha \cdot \inbraces{\vz \in \vspan{\vv_1,\dots,\vv_d} \suchthat \norm{\vz}_2 \le 1} \cup \inbraces{\sqrt{1+2\alpha} \cdot \vv_{d'+1}}}
\end{align*}
Let $\vw$ be an arbitrarily chosen unit vector in $\vspan{\vv_1,\dots,\vv_{d'}}$. Observe that it is enough to show
\begin{align*}
    \inbraces{\vz \in \vspan{\vw, \vv_{d'+1}} \suchthat \norm{\vz - \frac{\alpha}{\sqrt{1+2\alpha}} \cdot \vv_{d'+1}}_2 \le \frac{\alpha}{\sqrt{1+2\alpha}}} \subseteq \mathsf{conv}\inparen{\pm \alpha \cdot \vw, \sqrt{1+2\alpha} \cdot \vv_{d'+1}}.
\end{align*}
Since the above is a two-dimensional problem and that $\ip{\vw, \vv_{d'+1}} = 0$, it is equivalent to show that the inradius of the triangle with vertices $(-\alpha, 0)$, $(\alpha, 0)$, and $(0, \sqrt{1+2\alpha})$ is $\nfrac{\alpha}{\sqrt{1+2\alpha}}$ and that its incenter is $(0, \nfrac{\alpha}{\sqrt{1+2\alpha}})$.

Recall that the inradius of a triangle can be written as $\nfrac{K}{s}$ where $K$ is the area of the triangle (in this case, $\alpha\sqrt{1+2\alpha}$) and $s$ is the semiperimeter of the triangle (in this case, $1+2\alpha$). This implies that the inradius is indeed $\nfrac{\alpha}{\sqrt{1+2\alpha}}$. Finally, since the triangle in question is isosceles with its apex being the $y$-axis, the $x$-coordinate of its incenter must be $0$. These observations imply that the incenter is $(0, \nfrac{\alpha}{\sqrt{1+2\alpha}})$.

This is sufficient for us to conclude the proof of \Cref{lemma:irregular_update_easy}.
\end{proof}

We will now see that the analysis for the convenient update that we gave in \Cref{lemma:irregular_update_easy} is nearly enough for us to fully analyze the irregular update step. See \Cref{lemma:irregular_update_hard}, where we analyze the irregular update step in full generality (up to translating by $\vc_{t-1}$).

\begin{lemma}
\label{lemma:irregular_update_hard}
Let $Z \subset \R^d$ be a convex body such that $Z$ lies in a subspace $H$ of dimension $d' < d$. Let $\cE$ be an ellipsoid and let $0 < \alpha \le 1$ be such that
\begin{align*}
    \alpha \cdot \cE \subseteq Z \subseteq \cE.
\end{align*}
Let $\vz \notin H$. Then, there exists a center $\vc$ and an ellipsoid $\cE'$ such that
\begin{align*}
    \vc + \frac{1}{1+\nfrac{1}{\alpha}} \cdot \cE' \subseteq \mathsf{conv}\inparen{Z \cup \inbraces{\vz}} \subseteq \vc + \cE'.
\end{align*}
\end{lemma}
\begin{proof}[Proof of \Cref{lemma:irregular_update_hard}]
Recall that $\vv_1,\dots,\vv_{d'} \in \R^{d}$ are the unit vectors corresponding to the semiaxes of $\cE$; notice that these form a basis for $H$. Observe that $\vv_{d'+1}$ is a unit vector orthogonal to $\vv_1,\dots,\vv_{d'}$ such that $\vz$ can be expressed as $\sum_{i=1}^{d'+1} \vv_{i}\ip{\vv_i, \vz}$.

As stated in Algorithm \ref{alg:fully_online_approx}, let
\begin{align*}
    \mM \coloneqq \mI_d - \frac{1}{\ip{\vv_{d_t}',\vz_t}}\cdot\inparen{\vz_t - \sqrt{1+2\alpha} \cdot \vv_{d_t}'}(\vv_{d_t}')^T.
\end{align*}
We calculate
\begin{align*}
    \mM\vz_t = \vz_t - \frac{1}{\ip{\vv_{d_t}',\vz_t}}\cdot\inparen{\vz_t - \sqrt{1+2\alpha} \cdot \vv_{d_t}'}(\vv_{d_t}')^T\vz_t = \vz_t - \vz_t + \sqrt{1+2\alpha} \cdot \vv_{d_t}'  = \sqrt{1+2\alpha} \cdot \vv_{d_t}'.
\end{align*}
By the definition of $\mA_{t-1}$, we have
\begin{align*}
    \mA_{t-1}\mM\vz_t = \sqrt{1+2\alpha} \cdot \mA_{t-1}\vv_{d_t}' = \sqrt{1+2\alpha} \cdot \vv_{d_t}'.
\end{align*}
Next, for any $\vz \in Z$, we have $\vz \in H_{t-1}$. This means that
\begin{align*}
    \mM\vz = \vz - \frac{1}{\ip{\vv_{d_t}',\vz}}\cdot\inparen{\vz - \sqrt{1+2\alpha} \cdot \vv_{d_t}'}(\vv_{d_t}')^T\vz = \vz - 0 = \vz.
\end{align*}
By \Cref{lemma:irregular_update_easy}, we know for
\begin{align*}
    \mA_{t-1}\mM\vc_t &= \frac{\alpha}{\sqrt{1+2\alpha}}\cdot\vv_{d_t}'\\
    \mA_{t-1}\mM\cE_t &= \inbraces{\vz \in \vspan{\vv_1,\dots,\vv_{d_{t-1}},\vv_{d_t}'} \suchthat \norm{\vz}_2 \le 1}\\
\end{align*}
that
\begin{align*}
    \mA_{t-1}\mM\vc_t + \frac{1}{1+\nfrac{1}{\alpha_{t-1}}} \cdot \mA_{t-1}\mM\cE_t \subseteq \mathsf{conv}\inparen{\mA_{t-1}\mM\cdot Z \cup \inbraces{\mA_{t-1}\mM\vz_t}} \subseteq \mA_{t-1}\mM\vc_t + \mA_{t-1}\mM\cE_t
\end{align*}
and, since $\mA_{t-1}\mM$ is invertible (owing to the invertibility of $\mA_{t-1}$ and $\mM$),
\begin{align*}
    \vc_t + \frac{1}{1+\nfrac{1}{\alpha_{t-1}}} \cdot \cE_t \subseteq \mathsf{conv}\inparen{Z \cup \inbraces{\vz_t}} \subseteq \vc_t + \cE_t.
\end{align*}
Finally, note that
\begin{align*}
    \vc_t &= \frac{\alpha}{\sqrt{1+2\alpha}} \cdot \mM^{-1}\mA_{t-1}^{-1}\vv_{d_t}' = \frac{\alpha}{1+2\alpha} \cdot \vz_t \\
    \cE_t &= \inbraces{\vz \in \vspan{\vz_1,\dots,\vz_t} \suchthat \norm{\mA_{t-1}\mM\vz}_2 \le 1}
\end{align*}
and then translate by $\vc_{t-1}$, which concludes the proof of \Cref{lemma:irregular_update_hard}.
\end{proof}

We are now ready to prove \Cref{thm:main_two}.

\begin{proof}[Proof of \Cref{thm:main_two}]
Using \Cref{lemma:irregular_update_hard}, we have that the ellipsoids maintain our desired invariants (Definition \ref{def:alg-invariant}) throughout the process. Hence, \Cref{alg:fully_online_approx} maintains an ellipsoidal approximation to $\conv{\inbraces{\vz_1,\dots,\vz_t}}$ for all $t$.

It remains to verify the approximation factor $\alpha_t$ of Algorithm \ref{alg:fully_online_approx}.

Consider a timestep $t$. For every $t' \leq t$, let $H_{t'} = \vspan{\vz_1,\dots,\vz_{t'}}$,  $r_{t'} = r(Z_{t'})$ be the inradius of $Z_{t'} = \conv{\vz_1,\dots, \vz_{t'}}$, and $R_{t'} = R(Z_{t'})$ be the circumradius of $Z_t$. 
Let ${\hat r} = \min_{t'\leq t} r_{t'}$. Consider the $d$-dimensional ellipsoid $T(\cE_{t'})$ which is exactly equal to $\cE_{t'}$ in the space $H_{t'}$ and whose remaining semiaxes orthogonal to $H_{t'}$ are equal and have length $\hat r$. Observe that for a regular update step $t'$ (with $d_{t'} = d_{t'-1}$), we have
$$\frac{\mathsf{vol}_{d_{t'}}(\cE_{t'})}{\mathsf{vol}_{d_t}(\cE_{t'-1})} = \frac{\mathsf{vol}_d(T(\cE_{t'}))}{\mathsf{vol}_d(T(\cE_{t'-1}))}.$$
Now applying the evolution condition (\ref{eq:overview_evolution}) to the update restricted to $H_{t'}$, we get 
$$
\frac{1}{\alpha_{t'}} - \frac{1}{\alpha_{t'-1}} \leq C\log \frac{\mathsf{vol}_{d_{t'}}(\cE_t)}{\mathsf{vol}_{d_{t'}}(\cE_{t'-1})} = C\log \frac{\mathsf{vol}_d(T(\cE_{t'}))}{\mathsf{vol}_d(T(\cE_{t'-1}))}.
$$
We have obtained the following upper bound on the approximation-factor increase:
\begin{align}
    \frac{1}{\alpha_{t'}} - \frac{1}{\alpha_{t'-1}} \leq \begin{cases} 1 & \text{if $t'$ is an irregular update step} \\ C\logv{\frac{\mathsf{vol}_{d}(T(\cE_{t'}))}{\mathsf{vol}_{d}(T(\cE_{t'-1}))}} & \text{otherwise}\end{cases}\label{eq:fully_online_invariant}
\end{align}
Let $T_{\mathsf{reg}}$ consist of all the timesteps $t' \le t$ where we perform a regular update. Then we have,
\[
    \frac{1}{\alpha_t} - \frac{1}{\alpha_0}= \alpha_0 + \sum_{t'=1}^t  \left(\frac{1}{\alpha_{t'}} - \frac{1}{\alpha_{t'-1}} \right)\leq d_t + C\sum_{t'\in T_{\mathsf{reg}}} \logv{\frac{\mathsf{vol}_{d}(T(\cE_{t'}))}{\mathsf{vol}_{d}(T(\cE_{t'-1}))}}.
\]
Now we show that $\logv{\frac{\mathsf{vol}_{d}(T(\cE_{t'}))}{\mathsf{vol}_{d_t}(T(\cE_{t'-1}))}} \geq 0$ for an irregular step: let $\sigma_1\geq\dots\geq\sigma_d$ and $\sigma_1'\geq \dots\geq\sigma_d'$ be the lengths of semi-axes of $T(\cE_{t'})$ and $T(\cE_{t'-1})$, respectively. Then $\sigma_i \geq \sigma'_i$ for $1\leq i \leq d_{t'}-1$, since $\cE_{t'-1} \subset \cE_{t'}$; $\sigma_{d_{t'}} \geq r_{t'} \geq \hat r = \sigma'_{d_{t'}}$; and $\sigma_{i} = \hat r = \sigma'_{i}$ for $i > d_{t'}$. Therefore,
\[\logv{\frac{\mathsf{vol}_{d}(T(\cE_{t'}))}{\mathsf{vol}_{d_t}(T(\cE_{t'-1}))}}=
\logv{\frac{\sigma_1\cdot \ldots \cdot \sigma_d}{\sigma_1'\cdot \ldots \cdot \sigma_d'}} \geq \log 1 = 0.
\]
Using this inequality and plugging in $\alpha_0 = 1$, we get
\begin{align*}
    \frac{1}{\alpha_t} &= 1 + d_t + C\sum_{t'\in T_{\mathsf{reg}}} \logv{\frac{\mathsf{vol}_{d}(T(\cE_{t'}))}{\mathsf{vol}_{d}(T(\cE_{t'-1}))}}
    \leq 1+ d_t + C\sum_{t'=1}^t \logv{\frac{\mathsf{vol}_{d}(T(\cE_{t'}))}{\mathsf{vol}_{d}(T(\cE_{t'-1}))}} \\
    &\lesssim d_t + \logv{\frac{\mathsf{vol}_{d}(T(\cE_{t}))}{\mathsf{vol}_{d}(T(\cE_{0}))}} 
    \lesssim d_t + \logv{\frac{(R_t/\alpha_t)^{d_t}{\hat r}^{d-d_t}}{{\hat r}^d}} \lesssim d_t +  d_t \logv{\frac{R_t}{\alpha_t\hat r}}.
\end{align*}
We conclude that
\begin{align*}
    \frac{1}{\alpha_t} \lesssim d_t + d_t \logv{\frac{R_t}{\hat r}} + d_t\log d_t.
\end{align*}
This concludes the proof of \Cref{thm:main_two}.
\end{proof}

\subsection{Aspect ratio-independent bounds and proof of \texorpdfstring{Theorem~\ref{thm:main_two_ints}}{Theorem 3}}

To prove \Cref{thm:main_two_ints}, we first establish \Cref{lemma:irregular_volume_increase}.

\begin{lemma}
\label{lemma:irregular_volume_increase}
Let $t$ be an iteration corresponding to an irregular update step in \Cref{alg:fully_online_approx}. Then,
\begin{align*}
    \frac{\vol_{d_{t-1}}\inparen{B_2^{d_{t-1}}}}{\vol_{d_{t}}\inparen{B_2^{d_{t}}}} \cdot\frac{\vol_{d_t}\inparen{\cE_t}}{\vol_{d_t-1}\inparen{\cE_{t-1}}} \ge \frac{\norm{\vz_{t}^{\perp}}_2}{2}
\end{align*}
where $\norm{\vz_{t}^{\perp}}$ is the length of the component of $\vz_t$ in the orthogonal complement of $\vspan{\vz_1,\dots,\vz_{t-1}}$.
\end{lemma}
\begin{proof}[Proof of \Cref{lemma:irregular_volume_increase}]
By affine invariance, we can apply an affine transformation to map $\vz_t$ and $\cE_{t-1}$ to a convenient position. Hence, following the proof of \Cref{lemma:irregular_update_easy}, without loss of generality, suppose we have $\cE_{t-1} = B_{2}^{d_{t-1}}$ and $\vz_t = \sqrt{1+2\alpha} \cdot \ve_{d_t}$. By \Cref{lemma:irregular_update_easy}, the ellipsoid $\cE_t$ is a ball of radius $\nfrac{(1+\alpha)}{\sqrt{1+2\alpha}}$. Let $z \coloneqq \sqrt{1+2\alpha}$.
We now have
\begin{align*}
    \frac{\vol_{d_{t-1}}\inparen{B_2^{d_{t-1}}}}{\vol_{d_{t}}\inparen{B_2^{d_{t}}}} \cdot \frac{\vol_{d_t}\inparen{\cE_t}}{\vol_{d_{t-1}}\inparen{\cE_{t-1}}} = \inparen{\frac{1+\alpha}{\sqrt{1+2\alpha}}}^{d_t} \geq 1 > \frac{\norm{\vz_t}_2}{2}
\end{align*}
since $\|\vz_t\|_2 =\sqrt{1+2\alpha} \leq \sqrt{3} < 2$. This concludes the Proof of \Cref{lemma:irregular_volume_increase}.
\end{proof}

We will also need \Cref{lemma:det_prod_identity}, which we take from \citet{gk10}.

\begin{lemma}
\label{lemma:det_prod_identity}
Let $\mM \in \R^{r \times d}$ have linearly independent rows $\vm_1,\dots,\vm_r$. Then,
\begin{align*}
    \prod_{i=1}^r \norm{\vm_i}_2 = \sqrt{\detv{\mM\mM^T}}.
\end{align*}
\end{lemma}

We are now ready to prove \Cref{thm:main_two_ints}.

\begin{proof}[Proof of \Cref{thm:main_two_ints}]
Our approach is reminiscent of that used in the proof of Theorem 1.5 in \citet{woodruff2022high}.

By applying a translation to all points, we may assume without loss of generality that $\vz_1 = 0$. We will prove the guarantee for the last timestamp $t=n$ to simplify the notation. By replacing $n$ with $n'$, we can get a proof for any time stamp $t=n'$.

Let $S$ be the set of timestamps of irregular update steps excluding the first step. Since the update rule satisfies the evolution condition \eqref{eq:overview_evolution}, we have for all $t \notin S$ (recall that $d_t=d_{t-1}$ for $t\notin S$)
\begin{align*}
    \frac{\vol_{d_t}\inparen{\cE_{t}}}{\vol_{d_{t-1}}\inparen{\cE_{t-1}}} \ge \expv{\frac{1}{\alpha_t}-\frac{1}{\alpha_{t-1}}}.
\end{align*}
Next, by \Cref{lemma:irregular_volume_increase}, we have for every irregular update step $t> 1$ 
\begin{align*}
    \frac{\vol_{d_{t-1}}\inparen{B_2^{d_{t-1}}}}{\vol_{d_{t}}\inparen{B_2^{d_{t}}}} \cdot\frac{\vol_{d_t}\inparen{\cE_t}}{\vol_{d_{t-1}}\inparen{\cE_{t-1}}} \ge \frac{\norm{\vz_{t}^{\perp}}_2}{2}.
\end{align*}
Here, we assume that $\vol_0(\{0\}) = 1$ and define $\norm{\vz_2^{\perp}} = \norm{\vz_2}$.
Inductively combining the above for all $t > 1$ gives
\begin{align} 
    \vol_{d_n}\inparen{\cE_n} &\ge \prod_{t \notin S} \expv{\frac{1}{\alpha_t}-\frac{1}{\alpha_{t-1}}} \cdot \prod_{t \in S} \frac{\norm{\vz_t^{\perp}}_2}{2}\cdot \prod_{j=1}^{d_n} \frac{\vol_{j}\inparen{B_2^{j}}}{\vol_{j-1}\inparen{B_2^{j-1}}}\nonumber \\
    &= \prod_{t \notin S} \expv{\frac{1}{\alpha_t}-\frac{1}{\alpha_{t-1}}} \cdot \prod_{t \in S} \frac{\norm{\vz_t^{\perp}}_2}{2}\cdot \vol_{d_n}\inparen{B_2^{d_n}}\label{eq:int_potential}
\end{align}
Here we used that $\vol_0(\cE_{0}) = \vol_0(B_2^0) = 1$. Now invoking~\Cref{lemma:det_prod_identity}, we get
\begin{align*}
    \prod_{t \in S} \frac{\norm{\vz_t^{\perp}}_2}{2} \ge 2^{-\abs{S}}\sqrt{\detv{\mZ\vert_S\mZ\vert_S^T}} \ge 2^{-\abs{S}} = 2^{-d_n},
\end{align*}
where we used that $\detv{\mZ\vert_S\mZ\vert_S^T} \ge 1$ because all the vectors $\vz_t$ have integer coordinates. Moreover, since all coordinates are at most $N$ in absolute value, all the vectors  $\vz_t$ have length at most $N\sqrt{d}$. Therefore, $\frac{\vol\inparen{\cE_n}}{\vol\inparen{B_2^{d_n}}} \le \inparen{N\sqrt{d}}^{d_n}$. We plug these bounds back into \eqref{eq:int_potential}, rearrange, and take the logarithm of both sides, yielding
\begin{align*}
    \sum_{t \notin S} \frac{1}{\alpha_t}-\frac{1}{\alpha_{t-1}} \lesssim d_n\logv{d N}.
\end{align*}
Finally, by \eqref{eq:fully_online_invariant}, we have $\frac{1}{\alpha_t} - \frac{1}{\alpha_{t-1}} = 1$ for every $t \in S$. Combining everything gives
\begin{align*}
    \sum_{t \le n} \frac{1}{\alpha_t}-\frac{1}{\alpha_{t-1}} \lesssim d_n\logv{d N} + |S| \lesssim d_n\logv{d N},
\end{align*}
thereby concluding the proof of \Cref{thm:main_two_ints}.
\end{proof}

\section{Improved analysis for symmetric polytopes (Proof of \Cref{thm:main_one_symmetric})}
\label{sec:ellipsoid_symmetric}

In this section, we specialize the analysis framework developed in this chapter to the case when the polytope $Z$ is symmetric -- that is, when in each timestep $t$, we receive both $\vz_t$ and $-\vz_t$. We then prove \Cref{thm:main_one_symmetric}.

\smallit{Bibliographic notes.} The material in this section is derived from the paper \cite{mmo22}. Although the algorithm is the same as the one in that paper, the analysis given here is considerably simpler and fits within the framework developed earlier in this chapter.

\subsection{Monotone update rule for symmetric ellipsoidal approximation}

The main result of this subsection is \Cref{lemma:ellipsoid_symmetric_monup}.

\begin{lemma}
\label{lemma:ellipsoid_symmetric_monup}
The update rule given by
\begin{align*}
    \mA_{t} &= \begin{cases} \mA_{t-1} - \inparen{1 - \frac{1}{\norm{\mA_{t-1}\vz_{t}}_2}}\inparen{\frac{\inparen{\mA_{t-1}\vz_{t}}\inparen{\mA_{t-1}\vz_{t}}^{\top}}{\norm{\mA_{t-1}\vz_{t}}_2^2}}\mA_{t-1} & \text{ if } \norm{\mA_{t-1}\vz_t}_2 > 1 \\ \mA_{t-1} & \text{ otherwise} \end{cases}\\
    \vc_t &= 0 \\
    \alpha_{t} &= \begin{cases} \frac{\norm{\mA_{t-1}\vz_t}_2\alpha_{t-1}}{\sqrt{\norm{\mA_{t-1}\vz_t}_2^2-\alpha_{t-1}^2}} \cdot \frac{1}{\sqrt{1 + \inparen{\frac{\norm{\mA_{t-1}\vz_t}_2\alpha_{t-1}}{\sqrt{\norm{\mA_{t-1}\vz_t}_2^2-\alpha_{t-1}^2}}}^2}}  & \text{ if } \norm{\mA_{t-1}\vz_t}_2 > 1 \\ \alpha_{t-1} & \text{ otherwise} \end{cases}
\end{align*}
is a monotone update rule when in each timestep the algorithm receives both points $\pm \vz_t$.
\end{lemma}

The goal of the rest of this subsection is to prove \Cref{lemma:ellipsoid_symmetric_monup}.

The monotone update rule we use in this section is much easier to describe than in the general case. Let $0 = \vc_1 = \dots = \vc_n$ -- that is, we never shift the center. Let $\mA_{t-1}$ denote the matrix for the ellipsoid output at timestep $t$, so that we have $\norm{\mA_{t-1}\vz_{i}}_2 \le 1$ for all $i \le t-1$. Given a new point pair $\pm \vz_{t}$, we set $\mA_{t}$ according to the formula
\begin{align*}
    \mA_{t} = \mA_{t-1} - \inparen{1 - \frac{1}{\norm{\mA_{t-1}\vz_{t}}_2}}\inparen{\frac{\inparen{\mA_{t-1}\vz_{t}}\inparen{\mA_{t-1}\vz_{t}}^{\top}}{\norm{\mA_{t-1}\vz_{t}}_2^2}}\mA_{t-1}.
\end{align*}
This formula has a natural interpretation -- it is the minimum volume ellipsoid that covers both the ellipsoid $\cE_{t-1} = \inbraces{\vx\in\R^d \suchthat \norm{\mA_{t-1}\vx}_2 \le 1}$ and the new point pair $\pm \vz_{t}$. Specifically, see \Cref{lemma:ellipsoid_symmetric_minvol}.

\begin{lemma}
\label{lemma:ellipsoid_symmetric_minvol}
Let $\mA_{t-1} \in \R^{d \times d}$ and let
\begin{align*}
    \mA_{t} = \mA_{t-1} - \inparen{1 - \frac{1}{\norm{\mA_{t-1}\vz_{t}}_2}}\inparen{\frac{\inparen{\mA_{t-1}\vz_{t}}\inparen{\mA_{t-1}\vz_{t}}^{\top}}{\norm{\mA_{t-1}\vz_{t}}_2^2}}\mA_{t-1}.
\end{align*}
If $\cE_{t-1} \coloneqq \inbraces{\vx\in\R^d \suchthat \norm{\mA_{t-1}\vx}_2 \le 1}$, then the ellipsoid $\cE_{t} \coloneqq \inbraces{\vx\in\R^d \suchthat \norm{\mA_{t}\vx}_2 \le 1}$ is the minimum volume ellipsoid that contains both $\cE_{t-1}$ and $\pm \vz_{t}$.
\end{lemma}

Although we do not directly use \Cref{lemma:ellipsoid_symmetric_minvol}, it raises an interesting conceptual point, and so we give a proof below.

\begin{proof}[Proof of \Cref{lemma:ellipsoid_symmetric_minvol}]
In this proof, in an abuse of notation, let $\mM^{-1}$ denote the pseudoinverse of $\mM$ and let $\detv{\mM}$ of a singular matrix be the product of its nonzero singular values.

Note that the volume of the ellipsoid determined by $\mA_t$ is proportional to $\detv{\mA_t^{-1}}$. Therefore, $\mA_t$ is the solution to the following optimization problem, where we use that the volume of the ellipsoid determined by $\mA_t$ is proportional to $\detv{\mA_t^{-1}}$.
\begin{align*}
    \max \detv{\mA_t} \text{ such that } \mA_t \preceq \mA_{t-1} \text{ and } \norm{\mA_t\vz_t}_2 \le 1
\end{align*}
Additionally, since $\detv{\mA\mB} = \detv{\mA}\cdot\detv{\mB}$, we have that this objective is invariant under linear transformations. It thus follows that our objective can be rewritten as
\begin{align*}
    &\max \detv{\mA_t} &\text{ such that } &\mA_t \preceq \mA_{t-1} \text{ and } \norm{\mA_t\vz_t}_2 \le 1 \\
    \equiv &\max \detv{\mA_t \cdot \mA_{t-1}^{-1}} &\text{ such that } &\mA_t \cdot \mA_{t-1}^{-1} \preceq \mI \text{ and } \norm{\inparen{\mA_t\cdot \mA_{t-1}^{-1}}\mA_{t-1}\vz_t}_2 \le 1 \\
    \equiv &\max \detv{\widehat{\mA}} &\text{ such that } &\widehat{\mA} \preceq \mI \text{ and } \norm{\widehat{\mA}\mA_{t-1}\vz_{t}}_2 \le 1,
\end{align*}
where the last line follows from using the intermediate variable $\widehat{\mA} = \mA_t\cdot \mA_{t-1}^{-1}$.

In other words, after the transformation, the problem is equivalent to finding the minimum volume ellipsoid that contains (i) the unit ball and (ii) point $\mA_{t-1}\vz_t$. Geometrically, it is clear what the optimal ellipsoid for this problem is: one of its semi-axes is $\mA_{t-1}\vz_t$; all others are orthogonal to $\mA_{t-1}\vz_t$ and have length $1$ (this can be formally proved using symmetrization). However, we do not use this observation and derive a formula for $\widehat{\mA}$ using linear algebra.

We first give an upper bound on the objective value of the above optimization problem. Since $\widehat{\mA} \preceq \mI$, we have that all its singular values must be at most $1$. Additionally, since $1 \ge \norm{\widehat{\mA}\mA_{t-1}\vz_t}_2 \ge \sigma_{\min}\inparen{\widehat{\mA}} \cdot \norm{\mA_{t-1}\vz_t}_2$, we have that at least one singular value of $\widehat{\mA}$ must be $\le \nfrac{1}{\norm{\mA_{t-1}\vz_t}_2}$. Putting everything together and using the fact that the determinant is the product of the singular values gives $\detv{\widehat{\mA}} \le \nfrac{1}{\norm{\mA_{t-1}\vz_t}_2}$.

We now show that there exists a setting of $\widehat{\mA}$ that achieves this upper bound. Let $\vv_1$ be a unit vector in the direction of $\mA_{t-1}\vz_t$ and $\vv_2, \dots, \vv_d$ complete the orthonormal basis for $\R^d$ from $\vv_1$, and write \(
    \widehat{\mA} = \frac{1}{\norm{\mA_{t-1}\vz_t}_2}\vv_1\vv_1^{\top} + \sum_{i=2}^d \vv_i\vv_i^{\top}
\). We will show that $\widehat{\mA}$ satisfies the constraints imposed by the optimization problem. Since we have $\norm{\mA_{t-1}\vz_t}_2 \ge 1$ (as we impose that $\vz_t \notin \cE_{\mA_{t-1}}$), the fact that $\widehat{\mA} \preceq \mI$ follows immediately. For the second constraint, we write
\begin{align*}
    \norm{\widehat{\mA}\mA_{t-1}\vz_t}_2 = \norm{\inparen{\frac{1}{\norm{\mA_{t-1}\vz_t}_2}\vv_1\vv_1^{\top} + \sum_{i=2}^d \vv_i\vv_i^{\top}}\mA_{t-1}\vz_t}_2 = \norm{\frac{\mA_{t-1}\vz_t}{\norm{\mA_{t-1}\vz_t}_2}}_2 = 1.
\end{align*}
Furthermore, it is easy to see that $\detv{\widehat{\mA}} = \norm{\mA_{t-1}\vz_t}_2^{-1}$, which achieves our upper bound.

Finally, recall that we wrote $\widehat{\mA} = \mA_t\cdot \mA_{t-1}^{-1}$; rearranging this gives us the conclusion of \Cref{lemma:ellipsoid_symmetric_minvol}.
\end{proof}

We now prove \Cref{lemma:ellipsoid_symmetric_monup}.

\begin{proof}[Proof of \Cref{lemma:ellipsoid_symmetric_monup}]
The only nontrivial cases we have to deal with are  when $\norm{\mA_{t-1}\vz_t}_2 > 1$ -- in other words, when $\vz_t \notin \cE_{t-1}$. So, we assume this in the rest of the proof.

To validate $\vz_t \in \vc_t + \cE_t$, it is enough to show $\norm{\mA_t\vz_t}_2 \le 1$. We have
\begin{align*}
    \mA_t\vz_t = \vx - \inparen{1-\frac{1}{\norm{\vx}_2}}\inparen{\frac{\vx\vx^{\top}}{\norm{\vx}_2^2}}\vx = \vx - \inparen{1-\frac{1}{\norm{\vx}_2}}\vx = \frac{\vx}{\norm{\vx}_2},
\end{align*}
as desired.

The main challenge is to prove that $\vc_t + \alpha_{t}\cE_t \subseteq \conv{(\vc_{t-1}+\alpha_{t-1}\cE_{t-1}) \cup \inbraces{\vz_t}}$, or in other words, that the ``inner ellipsoid'' is always feasible. Without loss of generality, by applying an affine transformation and considering the reduced $2$-dimensional case in the same way as in the proof of \Cref{thm:main_one}, we let $\vz_{t} = z \cdot \ve_1$ for $z > 1$ and $\cE_{t-1} = B_2^2$. Notice that $z = \norm{\mA_{t-1}\vz_t}_2$. Our goal is to show that the ellipse $\alpha_t \cdot \cE_t$ lies inside the shape $\conv{\alpha_{t-1}B_2^2 \cup \inbraces{\pm z \cdot \ve_1}}$.

First, we prove that the parallelogram $P$ defined by the lines
\begin{align*}
    y &= -\frac{\alpha_{t-1}}{\sqrt{z^2-\alpha_{t-1}^2}}(x - z) \\
    y &= \frac{\alpha_{t-1}}{\sqrt{z^2-\alpha_{t-1}^2}}(x - z)
\end{align*}
and their reflections over the $y$-axis (i) passes through $\pm z \cdot \ve_1$ and (ii) has exactly four points of tangency with the circle $\alpha_{t-1}B_2^2$. Then, it will be enough to argue that $\alpha_t\cE_t \in P$.

The first part is clear from the formulas so we handle the second. Here, it is enough to show that the slope of the line that passes through $z \cdot \ve_1$ and is tangent to $\alpha_{t-1}B_2^2$ has slope $-\frac{\alpha_{t-1}}{\sqrt{z^2-\alpha_{t-1}^2}}$. The distance between $z \cdot \ve_1$ and the point of tangency is $\sqrt{z^2-\alpha_{t-1}^2}$ by the Pythagorean Theorem and the result follows from similar triangles. Repeating this for all four lines defining $P$ proves (ii).

Now, we are ready to prove that $\alpha_{t}\cE_t \in P$. This is a morally similar calculation. In this reduced case, the equation that determines the $(x,y)$ pairs lying on the surface of $\alpha_t\cE_{t}$ is given by
\begin{align*}
    \frac{x^2}{z^2} + y^2 = \alpha_t^2.
\end{align*}
Applying a linear transformation to our space does not affect tangency. So, we apply the linear transformation that maps $z \cdot \ve_1$ to $\ve_1$ (this is equivalent to multiplying every $x$-coordinate by $1/z$ and leaving all the $y$-coordinates unchanged). Now, our goal is to show that the circle $\alpha_{t}B_2^2$ lies inside the parallelogram whose vertices in the nonnegative orthant are the pair of points
\begin{align*}
    \inparen{0,\frac{z\alpha_{t-1}}{\sqrt{z^2-\alpha_{t-1}^2}}} \text{ and } \inparen{1,0}.
\end{align*}
By calculating the area of this triangle in two ways, we get that the radius of the incircle of this parallelogram is
\begin{align*}
    \alpha_t = \frac{z\alpha_{t-1}}{\sqrt{z^2-\alpha_{t-1}^2}} \cdot \frac{1}{\sqrt{1^2 + \inparen{\frac{z\alpha_{t-1}}{\sqrt{z^2-\alpha_{t-1}^2}}}^2}}.
\end{align*}
Undoing the affine transformations and recalling that doing so means that $z = \norm{\mA_{t-1}\vz_{t}}_2$ completes the proof of \Cref{lemma:ellipsoid_symmetric_monup}.
\end{proof}

\subsection{Approximation guarantee via stronger evolution condition}

We now show that the monotone update rule given by \Cref{lemma:ellipsoid_symmetric_monup} satisfies an evolution condition that is stronger than \eqref{eq:overview_evolution}.

\begin{lemma}
\label{lemma:ellipsoid_symmetric_evolution}
Under the update rule given in the statement of \Cref{lemma:ellipsoid_symmetric_monup}, we have
\begin{align*}
    \frac{1}{\alpha_t^2} - \frac{1}{\alpha_{t-1}^2} \le \begin{cases} 2\logv{\frac{\vol\inparen{\cE_{t}}}{\vol\inparen{\cE_{t-1}}}} & \text{ if } \vz_{t} \notin \cE_{t-1} \\ 0 & \text{ otherwise} \end{cases}.
\end{align*}
\end{lemma}
\begin{proof}[Proof of \Cref{lemma:ellipsoid_symmetric_evolution}]
As usual, we only consider the case where $\vz_t \notin \cE_{t-1}$ so that there is an update.

We first prove
\begin{align*}
    \frac{1}{\alpha_t^2} - \frac{1}{\alpha_{t-1}^2} \le 2\logv{\norm{\mA_{t-1}\vz_t}_2}.
\end{align*}
Let $z \coloneqq \norm{\mA_{t-1}\vz_t}_2$ and let $\triangle \coloneqq \frac{z\alpha_{t-1}}{\sqrt{z^2-\alpha_{t-1}^2}}$. Using the formula for $\alpha_t$ given in \Cref{lemma:ellipsoid_symmetric_monup}, we have
\begin{align*}
    \frac{1}{\alpha_t^2} - \frac{1}{\alpha_{t-1}^2} = \frac{\triangle^2+1}{\triangle^2} -\frac{1}{\alpha_{t-1}^2} = 1 + \frac{z^2-\alpha_{t-1}^2}{z^2\alpha_{t-1}^2} -\frac{z^2}{z^2\alpha_{t-1}^2} = 1-\frac{1}{z^2} \le 2\log z,
\end{align*}
where the inequality follows from noting that $1-\frac{1}{z^2} = 2\log z$ at $z=1$ and the derivative of $2\log z$ is always at least as large as the derivative of $1-\frac{1}{z^2}$ whenever $z \ge 1$.

Next, we show $\vol\inparen{\cE_{t}} = \norm{\mA_{t-1}\vz_{t}}_2 \cdot \vol\inparen{\cE_{t-1}}$. We write
\begin{align*}
    \detv{\mA_t} = \detv{\mI - \inparen{1 - \frac{1}{\norm{\mA_{t-1}\vz_{t}}_2}}\inparen{\frac{\inparen{\mA_{t-1}\vz_{t}}\inparen{\mA_{t-1}\vz_{t}}^{\top}}{\norm{\mA_{t-1}\vz_{t}}_2^2}}} \cdot \detv{\mA_{t-1}} = \frac{1}{\norm{\mA_{t-1}\vz_t}_2} \cdot \detv{\mA_{t-1}},
\end{align*}
and rearranging gives us what we need.

Combining these conclusions gives us the statement of \Cref{lemma:ellipsoid_symmetric_evolution}.
\end{proof}

We are now ready to prove \Cref{thm:main_one_symmetric}.

\begin{proof}[Proof of \Cref{thm:main_one_symmetric}]
As in the asymmetric case (proof of \Cref{thm:main_one}), we consider two phases -- name them Phase I and Phase II. In Phase II, we perform the monotone update rule as written in \Cref{lemma:ellipsoid_symmetric_monup}. In Phase I, we perform the trivial update
\begin{align*}
    \mA_t &= \min\inbraces{1,\frac{1}{\norm{\mA_{t-1}\vz_t}_2}} \cdot \mA_{t-1} \\
    \vc_t &= 0 \\
    \alpha_t &= \min\inbraces{1,\frac{1}{\norm{\mA_{t-1}\vz_t}_2}} \cdot \alpha_{t-1},
\end{align*}
which can be visually described as applying the minimum scaling of the outer ellipsoid such that it contains the new points $\pm \vz_t$ while leaving the inner ellipsoid unchanged. It is easy to see that by just applying the Phase I update, the approximation factor we get is $R(Z)/r_0$.

Next, we describe the switching condition we use to change from Phase I to Phase II. At iteration $t$, if $R(Z_t)/r_0 \ge 2^{3/2}\sqrt{d\logv{4d}}$, then we change to Phase II. It is easy to see that $R(Z_t)$ is monotonically increasing in $t$, so once this switching condition is satisifed, it will remain satisfied until the termination of the algorithm.

Applying \Cref{lemma:ellipsoid_symmetric_evolution} across all the Phase II iterates and recalling that in Phase I we never change the inner ellipsoid, we have
\begin{align*}
    \frac{1}{\alpha_n^2} \le 2\logv{\frac{\vol\inparen{\cE_n}}{\vol\inparen{\cE_0}}}.
\end{align*}
We know that the containment $\alpha_n \cE_n \subseteq Z_n \subseteq \cE_n$ must hold. We also know that $Z_n \subseteq R(Z_n) \cdot B_2^d$ must also hold. Thus, we know that $\cE_n \subseteq \frac{R(Z_n)}{\alpha_n} \cdot B_2^d$. Hence, we have
\begin{align*}
    \frac{1}{\alpha_n^2} \le 2\logv{\frac{\vol\inparen{\cE_n}}{\vol\inparen{\cE_0}}} \le 2\logv{\frac{\inparen{\frac{R(Z_n)}{\alpha_n}}^d}{r_0^d}} = d\logv{\inparen{\frac{R(Z_n)}{r_0}}^2}+d\logv{\frac{1}{\alpha_n^2}}.
\end{align*}
We now show that if we ever enter Phase II, then we have $\alpha_n^{-2} \le (R(Z_n)/r_0)^2$. First, let us see how this gives us one of the cases of our proof. Substituting this gives us
\begin{align*}
    \frac{1}{\alpha_n^2} \le 4d\logv{\frac{R(Z_n)}{r_0}},
\end{align*}
so taking square roots gives the result. 

Now, for the sake of contradiction, suppose that $\alpha_n^{-2} \ge (R(Z_n)/r_0)^2$. This gives us
\begin{align*}
    \frac{1}{\alpha_n^2} \le 4d\logv{\frac{1}{\alpha_n^2}},
\end{align*}
so solving for $\alpha_n^{-1}$ gives
\begin{align*}
    \frac{1}{\alpha_n^2} < 8d\logv{4d}.
\end{align*}
This means that
\begin{align*}
    \frac{R(Z_n)}{r_0} \le 2^{3/2}\sqrt{d\logv{4d}},
\end{align*}
which means we always stayed in Phase I and our approximation factor is just $\frac{R(Z_n)}{r_0}$ (and the monotone update rule given by \Cref{lemma:ellipsoid_symmetric_monup} never actually executed).

Combining the cases tells us that for some universal constant $C$, we have
\begin{align*}
    \frac{1}{\alpha_n} \le C\min\inbraces{\frac{R(Z_n)}{r_0}, \sqrt{2d\logv{\frac{R(Z_n)}{r_0}}}},
\end{align*}
completing the proof of \Cref{thm:main_one_symmetric}.
\end{proof}

\section{Forming small coresets for convex bodies (Proof of Theorem~\ref{thm:main_three})}
\label{sec:coreset}

In this section, we prove Theorem \ref{thm:main_three}. See Algorithm \ref{alg:ellipse_to_coreset}.

\begin{algorithm}[H]
\caption{Streaming coreset for convex hull}\label{alg:ellipse_to_coreset}
\begin{algorithmic}[1]
\State \textbf{Input:} Stream of points $\vz_t$; Update rule for Algorithm \ref{alg:fully_online_approx} $\cA$.
\State \textbf{Output:} Set $S \subseteq [n]$.
\For{$t = 1,\dots,n$}
    \State Receive $\vz_{t}$.
    \State Let $\cE_{\mathsf{test}} = \cA(\vc_{t-1}, \cE_{t-1}, \vz_t)$.
    \State Let $d_t = \mathsf{dim}\inparen{\vspan{\vz_1-\vc_{t-1},\dots,\vz_t-\vc_{t-1}}}$.
    \If{$d_t > d_{t-1}$ or $\frac{\mathsf{Vol}_{d_t}(\cE_{\mathsf{test}})}{\mathsf{Vol}_{d_t}(\cE_{t-1})} \ge e$}
        \State Let $\vc_t, \cE_{t} = \cA(\vc_{t-1},\cE_{t-1},\vz_{t})$.
        \State Update $S_t = S_{t-1} \cup \inbraces{\vz_{t}}$.
    \Else
        \State Let $\vc_t,\cE_{t}=\vc_{t-1},\cE_{t-1}$.
        \State Let $S_{t-1} = S_t$.
    \EndIf
\EndFor
\State \textbf{Output:} $S_n$
\end{algorithmic}
\end{algorithm}

For a sketch of the intuition and the argument we will use for the proof, see Section \ref{sec:overview_coreset}.

\begin{proof}[Proof of Theorem \ref{thm:main_three}]

We prove two properties of Algorithm \ref{alg:ellipse_to_coreset}. First, we show $\abs{S_t} \le O\inparen{d_t \cdot \logv{d_t \cdot \max_{t' \le t} \nfrac{R_t}{r_{t'}}}}$ and, further, $\abs{S_t} \le O\inparen{d_t \cdot \logv{dN}}$ if points $\vz_t$ have integer coordinates between $-N$ and $N$. Second, we  show that $\conv{Z\vert_{S_t}} \subseteq \conv{Z\vert_{[t]}} \subseteq O\inparen{d_t \cdot \logv{d_t \cdot \max_{t' \le t} \nfrac{R_t}{r_{t'}}}}\cdot \conv{Z\vert_{S_t}}$ and $\conv{Z\vert_{S_t}} \subseteq \conv{Z\vert_{[t]}} \subseteq O\inparen{d_t \cdot \logv{dN}} \cdot \conv{Z\vert_{S_t}}$.

\paragraph{Bounding $\abs{S_t}$.}

It is enough to count the number of steps $t$ for which we have $\frac{\mathsf{Vol}_{d_t}(\cE_{\mathsf{test}})}{\mathsf{Vol}_{d_t}(\cE_{t-1})} \ge e$.

It is easy to see that for all $t$, we have  $r(Z\vert_{[t]}) \cdot \inparen{B_2^d \cap \vspan{\vz_1-\vc_t,\dots,\vz_t-\vc_t}} \subseteq \vc_{t} + \cE_t$. Additionally, by the definition of $R(Z)$, we always have $Z\vert_{[t]} \subseteq R(Z) \cdot \inparen{B_2^d \cap \vspan{\vz_1-\vc_t,\dots,\vz_t-\vc_t}}$. These are enough to give volume lower and upper bounds in each step. Next, for each step in which we add an element to $S_{t-1}$ to obtain $S_t$, the volume must increase by a factor of $e$. It easily follows that the number of elements in $S_t$ satisfies
\begin{align*}
    \abs{S_t} \le \logv{\max_{t' \le t} \frac{\prod_{i=1}^{d_t} R(Z\vert_{[t]})}{\prod_{i=1}^{d_t} r(Z\vert_{[t']})}} = d_t\logv{\max_{t' \le t} \frac{R(Z\vert_{[t]})}{r(Z\vert_{[t']})}}.
\end{align*}
We now give an upper bound for the case when all coordinated of $\vz_t$ are integers not exceeding $N$ in absolute value. It is easy to see that the update rule in \Cref{alg:ellipse_to_coreset} exactly corresponds to the steps where we have
\begin{align*}
    \frac{1}{\alpha_{t}} - \frac{1}{\alpha_{t-1}} \gtrsim 1,
\end{align*}
and in the same way as in the proof of \Cref{thm:main_two_ints}, we have for all $t$ that
\begin{align*}
    \sum_{t \ge 1} \frac{1}{\alpha_{t}} - \frac{1}{\alpha_{t-1}} \lesssim d_t\logv{dN}.
\end{align*}
It therefore follows that $\abs{S} \lesssim d_t\logv{dN}$, as desired.

\paragraph{Bounding the distortion of the chosen points.}

Consider some iteration $t' \le t$. Without loss of generality, let $\vc_{t'-1} = 0$. Suppose $\vz_{t'}$ does not result in an update to $S_{t'-1}$. This implies that $\vz_{t'} \in 2e \cdot \cE_{t'-1}$. Next, observe that $0 \in \vc_t + \cE_t$. Putting these together, we have $\vz_{t'} \in \inparen{\vc_t + \cE_t} + 2e \cdot \cE_{t'-1}$. Since $\cA$ is monotone, we must have $2e \cdot \cE_{t'-1} \subseteq \vc_t + e\cdot\cE_t$; hence, we may write $\vz_{t'} \in \vc_t + \inparen{2e+1}\cE_t$.

The inner ellipsoid $\vc_t + \alpha_t \cdot \cE_t$ will still be an inner ellipsoid for the points determined by $S_t$. Stitching together all our inclusions, we have
\begin{align}
    \vc_t + \alpha_t \cdot \cE_t \subseteq Z\vert_{S_t} \subseteq Z \subseteq \vc_t + (2e+1)\cE_t \subseteq \frac{2e+1}{\alpha_t} \cdot Z\vert_{S_t}\label{eq:ellipsoid_coreset_inclusions}.
\end{align}
which means that
\begin{align*}
    Z\vert_{S_t} \subseteq Z \subseteq O\inparen{d_t \cdot \logv{d_t \cdot \max_{t' \le t} \frac{R(Z\vert_{([t] \cap S_t)})}{r(Z\vert_{([t'] \cap S_t)})}}} \cdot Z\vert_{S_t}.
\end{align*}
Notice that this is nearly what we want, except that the aspect ratio term is in terms of the subset body $Z\vert_{S_t}$. To obtain the final guarantee in terms of the aspect ratio of $Z\vert_{[t]}$, observe that the above guarantee readily implies that
\begin{align*}
    O\inparen{d_t \cdot \logv{d_t \cdot \max_{t' \le t} \frac{R(Z\vert_{([t] \cap S_t)})}{r(Z\vert_{([t'] \cap S_t)})}}} \le O\inparen{d_t \cdot \logv{d_t \cdot \max_{t' \le t} \frac{R(Z\vert_{[t]})}{r(Z\vert_{[t']})}}}.
\end{align*}
We now give the corresponding improvement when the $\vz_t$ are integer-valued. As before, \eqref{eq:ellipsoid_coreset_inclusions} holds. From this, we get
\begin{align*}
    Z\vert_{S_t} \subseteq Z \subseteq O\inparen{d_t \cdot \logv{dN}} \cdot Z\vert_{S_t},
\end{align*}
as desired. This concludes the proof of Theorem \ref{thm:main_three}.
\end{proof}

\section{Approximation lower bound for monotone algorithms}
\label{sec:lb}

In this section, we show \Cref{thm:main_four_symmetric} and \Cref{thm:main_four}.

\subsection{Inapproximability of John's ellipsoid}

Before we begin, recall that a Hadamard basis is a set of vectors $\vv_1, \dots, \vv_d$ such that:
\begin{itemize}
    \item $\norm{\vv_i} = 1$, for all $i \in [d]$;
    \item For all $i \neq j$, $\ip{\vv_i,\vv_j}=0$;
    \item Every entry of $\vv_i$ is in $\inbraces{\pm \nfrac{1}{\sqrt{d}}}$.
\end{itemize}

Our family of hard instances proceeds in two phases.

\greenbox{\paragraph{Phase 1} Let $d$ be such that there exists a Hadamard basis for $\R^d$. Consider a corresponding Hadamard basis $\vv_1, \dots, \vv_d$. The adversary gives the algorithm the points $\vv_1, \dots, \vv_d$.

\paragraph{Phase 2} The adversary selects $i \in [d]$ arbitrarily and $\eps \in (0, d-1)$ arbitrarily. They then define the vectors $\vw_i \coloneqq \ve_i \cdot \nfrac{1}{\sqrt{d-\eps}}$ and $\vw_j \coloneqq \ve_j \cdot \sqrt{\nfrac{d-1}{\eps}}$ for all $j \neq i$. The adversary gives the algorithm the points $\vw_1, \dots, \vw_d$. Call the outcome here ``Outcome (i).''}

It is easy to see that at the end of Phase 1, the minimum volume outer ellipsoid is simply $B_2^d$. Furthermore, the algorithm's solution $\widehat{\cE}$ contains $\conv{\pm \vv_1, \dots, \pm \vv_d}$. On the other hand, consider the following claim.

\begin{lemma}
The following ellipsoid is the minimum-volume outer ellipsoid for Outcome (i):
\begin{align*}
    \cE_{OPT(i)} = \inbraces{x \suchthat 1 \ge \frac{x_i^2}{\inparen{\nfrac{1}{\sqrt{d-\eps}}}^2} + \sum_{j \neq i}^d \frac{x_j^2}{\inparen{\sqrt{\nfrac{d-1}{\eps}}}^2}}
\end{align*}
\end{lemma}
\begin{proof}
Notice that all the points $\vw_j$ are orthogonal. Thus, the minimum-volume outer ellipsoid containing all the $\vw_j$ must be the one whose axes are along the directions of $\vw_j$ and whose poles are located on $\vw_j$. Observe that $\cE_{OPT(i)}$ satisfies this, so it must be the minimum-volume outer ellipsoid for the convex body whose vertices are determined by the $\vw_j$.

It now remains to show that every Hadamard basis vector is on the surface of $\cE_{OPT(i)}$:
\begin{align*}
    \frac{\inparen{\nfrac{1}{\sqrt{d}}}^2}{\inparen{\nfrac{1}{\sqrt{d-\eps}}}^2} + \sum_{j\neq i}^d \frac{\inparen{\nfrac{1}{\sqrt{d}}}^2}{\inparen{\sqrt{\nfrac{d-1}{\eps}}}^2} = \frac{1}{d}\inparen{(d-\eps) + (d-1) \cdot \frac{\eps}{d-1}} = 1
\end{align*}
Since the minimum volume ellipsoid containing all the $w_j$ also contains the Hadamard basis vectors, it (i.e., $\cE_{OPT(i)}$) must be the minimum-volume outer ellipsoid for Outcome (i).
\end{proof}

\begin{proof}[Proof of \Cref{thm:main_four_symmetric}]
We will now show that any that outputs an ellipsoid $\widehat{\cE}$ at the end of Phase 1 must have an approximation factor of at least $\sqrt{d-\eps}$ on at least one of Outcomes ($1, \dots, i$). Suppose that in each of Outcome (i), we obtain an ellipsoid $\widehat{\cE_i}$ that satisfies $C \cdot \cE_{OPT(i)} \supseteq \widehat{\cE_i}$. We now have:
\begin{align*}
    \conv{\inbraces{\pm \vv_1, \dots, \vv_d}} \subseteq \widehat{\cE} \subseteq \bigcap_{i=1}^d \widehat{\cE_{i}} \subseteq C \cdot \bigcap_{i = 1}^d \cE_{OPT(i)}
\end{align*}
We therefore want to argue about $\widehat{\cE}$ given that it must contain $\conv{\inbraces{\pm \vv_1, \dots, \vv_d}}$ and be contained by $C \cdot \bigcap_{i = 1}^d \cE_{OPT(i)}$. Let $A$ be a matrix mapping $\widehat{\cE}$ to the unit ball. Then, notice that we can write for all $i \in [d]$:
\begin{align*}
    \norm{\mA\vv_i}_2 &\le 1 & \norm{\mA \cdot \frac{C\ve_i}{\sqrt{d-\eps}}}_2 &\ge 1
\end{align*}
In particular, the rightmost exclusion follows from the fact that $\nfrac{C\ve_i}{\sqrt{d-\eps}}$ lies on the boundary of $C \cdot \bigcap_{i=1}^d \cE_{OPT(b.i)}$. Now, recall the well-known fact that for any unitary matrix $\mW$, we have $\fnorm{\mA\mW} = \fnorm{\mA}$ (see, e.g., \cite{horn1991}), and observe that we have
\begin{align*}
    d \ge \sum_{i=1}^d \norm{\mA\vv_i}^2 = \fnorm{\mA\mV}^2 = \fnorm{\mA}^2 = \fnorm{\mA\mI}^2 = \sum_{i=1}^d \norm{\mA\ve_i}^2 \ge \frac{d(d-\eps)}{C^2}.
\end{align*}
Rearranging gives $C \ge \sqrt{d-\eps}$, as desired.
\end{proof}

\subsection{Lower bound adversary}

Our proof of \Cref{thm:main_four} constructs an adversary, which given a monotone algorithm \(\cA\) and \(\kappa \geq 1\),
constructs a sequence of points \(\vz_1, \ldots, \vz_n\) satisfying \(\kappa(\conv{\vz_1, \ldots, \vz_n}) \leq \kappa\) to witness that the algorithm does not produce an approximation better than \(\widetilde{\Omega}(d \log \kappa)\).
While by definition \(\kappa = \frac{R}{r}\), our construction keeps \(r = 1\) (notice that any lower bound construction must be scale-invariant),
and for simplicity we use \(R = \kappa\).

Let \(\vz^\Delta_1, \vz^\Delta_2, \ldots, \vz^\Delta_{d+1} \in \R^{d}\) be the \(d+1\) vertices of a regular simplex \(\Delta_{d}\) that circumscribes \(B_2^d\).
Our adversary is described in \cref{alg:lower_bound}.
It uses a first phase that feeds \(\cA\) the vertices of \(\Delta_d\),
then a second phase that repeatedly feeds \(\cA\) points at a constant distance from the previous ellipsoid.
Specifically, every new point \(\vz_t\) in the second phase is in \(\vc_{t-1} + 2 \cdot \cE_{t-1}\), i.e. its distance is 2 from \(\vc_{t-1}\) in the norm that is the gauge of \(\cE_{t-1}\).

\begin{algorithm}[h]
\caption{Lower bound adversary}\label{alg:lower_bound}
\textbf{Input}: Monotone algorithm \(\cA\), \(R \geq 1\)
\begin{algorithmic}[1]
\State \((\vc_0 + \cE_0, \alpha_0) = (0 + B_2^d, 1)\)
\Comment{Initialize to the unit ball}
\For{\(t \in \{1, 2, \ldots, d+1\}\)}
    \Comment{Phase I: feed \(\cA\) the vertices of a simplex}
    \State \((\vc_{t} + \cE_{t}, \alpha_{t}) = \cA(\vc_{t-1} + \cE_{t-1}, \alpha_{t-1}, \vz^\Delta_{t})\)
\EndFor
\State \(t \leftarrow d+2\)
\While{\(\vol(\cE_{t-1}) \leq \vol\inparen{\frac{R}{2} \cdot B_2^d}\)}\label{alg_line:adv_stopvol}
\Comment{Phase II: feed \(\cA\) points outside the previous ellipsoid}
    \State Let \(F_{t-1} = \partial(\vc_{t-1} + 2 \cE_{t-1}) \cap (R \cdot B_2^d)\)
    \If{\(F_{t-1} = \varnothing\)}
        \State \textbf{stop}\label{alg_line:adv_stopdisj}
    \EndIf
    \State Let arbitrary \(\vz_{t} \in F_{t-1}\)
    \State \((\vc_{t} + \cE_{t}, \alpha_{t}) = \cA(\vc_{t-1} + \cE_{t-1}, \alpha_{t-1}, \vz_{t})\)
    \State \(t \leftarrow t + 1\)
\EndWhile
\end{algorithmic}
\end{algorithm}

\begin{remark}
 This particular construction we give of the hard case is adaptive, meaning that the adversary's choice of points depend on the previous ellipsoids the algorithm outputs.
   However, this adversary can be made non-adaptive by taking an \(\varepsilon\)-net \(S\) of \(B_2^d\) for sufficiently small \(\varepsilon\),
   then feeding \(\cA\) the sequence of points in sets \(S, 2 \cdot S, 4 \cdot S, \ldots, 2^{\log_2 R - 1}, 2^{\log_2 R} \cdot S\).
   In consequence, this means that randomization on the part of the monotone algorithm does not help, unlike some other online settings.
\end{remark}

Let \(T\) be the largest value of \(t -1\) before the adversary halts.
We first show that the adversary only gives finitely many points before halting.
\begin{lemma}
    \(T \leq O(d \log R)\)
\end{lemma}
\begin{proof}
    We argue that the volume of \(\cE_{t}\) increases by at least a constant factor on each iteration.
    This is sufficient to bound the number of iterations by \(O(d \log R)\), 
    as \(\cE_0 = B_2^d\),
    and Line \ref{alg_line:adv_stopvol} is no longer true when the volume of \(\cE_t\) exceeds \(\inparen{\frac{R}{2}}^d \cdot \vol(B_2^d)\).

    We claim that for all \(t \geq d + 2\), \(\vol(\cE_{t}) \geq \frac{3}{2} \cdot \vol(\cE_{t-1})\).
    By applying a nonsingular affine transformation, we can assume without loss of generality that \(\cE_{t-1} = B_2^d\).
    With a further rotation, we can assume the newly received point is \(\vz_t = 2 \ve_1\).
    From monotonicity of \(\cA\) we must have that \(\vc_{t} + \cE_{t} \supseteq B_2^d \cup \{2 \ve_1\}\).
    Clearly every semi-axis of \(\cE_{t}\) must have length 
    at least 1 in order to contain \(\cE_{t-1}\).
    Observe that \(\cE_{t}\) must also contain the segment connecting
    \(-1 \ve_1\) and \(2 \ve_1\), and so at least one semi-axis must have length at least \(\frac{3}{2}\)
    (if not, the diameter of \(\cE_{t}\) would be strictly less than \(3\)).
    Hence as \(\frac{\vol(\cE_{t})}{\vol(B_2^d)}\) equals the product of the length of the semi-axes of \(\cE_{t}\), we have \(\vol(\cE_{t}) \geq \frac{3}{2} \vol(B_2^d)\).
\end{proof}

For the analysis we define quantities \(A_t, P_t\) associated with the sequence of ellipsoids for \(1 \leq t \leq T\):
\[A_t \defeq \frac{1}{\alpha_t}, \quad P_t \defeq \log \inparen{\frac{\vol(\cE_t)}{\vol(B_2^d)}}\]
By the monotonicity of \(\cA\), we have that \(A_t\) and \(P_t\) are both nondecreasing in \(t\).
We first observe that the adversary guarantees that the final volume of the ellipsoid output by \(\cA\) is large:
\begin{lemma}\label{lemma:lower_bound_vol}
    At the conclusion of \Cref{alg:lower_bound}'s execution, we have
   \[P_T  \geq d \log \frac{R}{2}\]
\end{lemma}
\begin{proof}
    There are two ways that the adversary stops:
    if the condition in Line \ref{alg_line:adv_stopvol} is no longer true,
    or if Line \ref{alg_line:adv_stopdisj} is reached.
    If the former occurs, then we have \(\vol(\cE_{T}) > \vol(\frac{R}{2} \cdot B_2^d)\),
    and clearly then \(P_T \geq d \log\inparen{\frac{R}{2}}\).
    
    In the latter stopping condition, the algorithm halts at time \(T\) when \(\partial (c_{T} + 2 \cE_{T}) \cap R \cdot B_2^d = \varnothing\).
    The sets \(\partial (c_{T} + 2 \cE_{T})\) and \(R \cdot B_2^d\) can be disjoint in two cases:
    \(c_T + 2 \cE_T\) and \(R \cdot B_2^d\) are disjoint;
    or \(R \cdot B_2^d \subseteq c_T + 2 \cE_t\) with the boundaries of both ellipsoids disjoint.
    By the monotonicity of \(\cA\), we have \(1 \cdot B_2^d \subseteq c_T + 2 \cE_{T}\),
    and so eliminate the former case.
    But then \(\vol(2 \cdot \cE_T) \geq \vol(R \cdot B_2^d) = R^d \vol(B_2^d)\), and taking logarithms on both sides yields the claim.
\end{proof}

Now in contrast to the upper bound where we essentially gave an algorithm for which \(\frac{\Delta A}{\Delta P}\) was upper bounded by a constant,
here we will show a constant lower bound on the same quantity for any monotone algorithm.
\begin{lemma}\label{lemma:lb_step_main}
    There exists a constant \(C_{\ref{eqn:lb_step_main}} > 0\) such that if \(A_t \geq d\), we have
    \begin{equation}\label{eqn:lb_step_main}
    A_{t+1} - A_t \geq C_{\ref{eqn:lb_step_main}} (P_{t+1} - P_t)
    \end{equation}
\end{lemma}

Observe that this lower bound requires \(A_t \geq d\),
hence necessitating a first phase using the simplex,
whose optimal roundings show tightness for John's theorem for general convex bodies.
In order to prove the lower bound we also need a second property, that \(A_{t}\) is large compared to \(P_t\).
\begin{lemma}\label{lemma:lb_simplex_symmetry}
    Let \(0 \leq \alpha \leq 1, \vc \in \R^d\), and \(\cE\) be an ellipsoid such that
    \[\vc + \alpha \cdot \cE \subseteq \Delta_d \subseteq \vc + \cE\]
    then we have: \begin{enumerate} 
        \item\label{item:lb_simplex_symmetry_1} \(\alpha \leq \frac{1}{d}\)
        \item \(\log\inparen{\frac{\vol(\cE)}{\vol(B_2^d)}} \leq O\inparen{\log(d) \cdot \frac{1}{\alpha}}\)
    \end{enumerate}
\end{lemma}

With the statements of these claims in hand, we are ready to prove the lower bound.
\begin{proof}[Proof of \Cref{thm:main_four}]
 It is clear that \(\kappa(\conv{\vz_1, \ldots, \vz_T)} \leq R\), as for every \(1 \leq t \leq T\) the adversary guarantees \(1 \leq \|\vz_t\|_2 \leq R\). So we focus on showing a lower bound on the quality of the approximation produced by \(\cA\).

As \(\cA\) is monotone, after the end of Phase I
we must have that \begin{equation}\label{eq:lb_ph1_approx}
\vc_{d+1} + \alpha_{d+1} \cdot \cE_{d+1} \subseteq \Delta_d \subseteq \vc_{d+1} + \cE_{d+1}
\end{equation}
Now because \(\cE_{d+1}\) satisfies the conditions of \Cref{lemma:lb_simplex_symmetry},
we get using the definition \(A_{d+1} = \frac{1}{\alpha_{d+1}}\) that \(A_{t} \geq A_{d+1} \geq d\) for any \(t \geq d+1\).
Then we can apply \Cref{lemma:lb_step_main} for every \(t \geq d + 1\) until termination of the algorithm:
\begin{align*}
A_{d+2} - A_{d+1} &\geq C_{\ref{eqn:lb_step_main}} \inparen{P_{d+2} - P_{d+1}} \\
A_{d+3} - A_{d+2} &\geq C_{\ref{eqn:lb_step_main}} \inparen{P_{d+3} - P_{d+2}} \\
&~\vdots \\
A_{T-1} - A_{T-2} &\geq C_{\ref{eqn:lb_step_main}} \inparen{P_{T-1} - P_{T-2}} \\
A_{T} - A_{T-1} &\geq C_{\ref{eqn:lb_step_main}} \inparen{P_{T} - P_{T-1}}
\end{align*}
Summing these inequalities, we have
\begin{align*}
    \sum_{t=d+1}^{T-1} A_{t+1} - A_{t} \geq C_{\ref{eqn:lb_step_main}} \inparen{\sum_{t=d+1}^{T-1} P_{t+1} - P_t}
\end{align*}
Both sides of this inequality are telescoping sums, so simplifying we get
\begin{equation}\label{eqn:thm_lower_bound_1}
A_T \geq A_{d+1} + C_{\ref{eqn:lb_step_main}} (P_T - P_{d+1})
\end{equation}

Again because we can apply \Cref{lemma:lb_simplex_symmetry} for \(\cE_{d+1}\), we have \(P_{d+1} \leq O(\log(d) \cdot A_{d+1})\), which along with (\ref{eqn:thm_lower_bound_1}) yields
\[A_T \geq A_{d+1} + \Omega(P_T - \log(d) \cdot A_{d+1}) \geq  \Omega (P_T - \log(d) \cdot A_{d+1})\]
Thus we have \[
A_T \geq \Omega(\max(A_{d+1}, P_T - \log(d) \cdot A_{d+1})) \geq \Omega\inparen{\frac{P_T}{\log(d)}}
\]
and we get the desired bound using \Cref{lemma:lower_bound_vol}.
\end{proof}

Our proof of \Cref{lemma:lb_step_main}
relies on a symmetrization argument to a reduced case (essentially two-dimensional, like for our algorithms). We now define this reduced case, and related quantities.
\begin{definition}\label{defn:lb_reduced}
In the reduced case, the previous outer and inner ellipsoids are given by \(B_2^d, \alpha \cdot B_2^d\), and the received point is \(\vz = 2 \ve_1\).
The next outer and inner ellipsoids are given by \(c \cdot \ve_1 + \cE_{\mM}, c \cdot \ve_1 + \alpha' \cdot \cE_{\mM}\)
for \(c \in \R\), and \(\mM = \diag{a, b, b, \ldots, b, b}\) for \(a, b \geq 1\).
We let \(\Delta A = \frac{1}{\alpha'} - \frac{1}{\alpha}\) and \(\Delta P = \log \inparen{\frac{\vol(\cE_{\mM})}{\vol(B_2^d)}}\).
\end{definition}

Note that the update in this reduced case is monotone if
\(B_2^d \cup \{ 2 \ve_1 \} \subseteq c \cdot \ve_1 + \cE_{\mM}\)
and \(c \cdot \ve_1 + \alpha' \cdot \cE_{\mM} \subseteq \conv{(\alpha \cdot B_2^d) \cup \{ 2 \ve_1 \}}\).

Now we state the lower bound on \(\frac{\Delta A}{\Delta P}\) in this setting, which is established in \Cref{sec:lb_2d}.
It is exactly the bound of \Cref{lemma:lb_step_main} in this special case.
\begin{lemma}\label{lemma:lb_reduced}
     In the reduced case, for any monotone update \(c \cdot \ve_1 + \cE_{\mM}, c \cdot \ve_1 + \alpha' \cdot \cE_{\mM}\) when \(\alpha' \leq \frac{1}{d}\) we have
    \begin{equation}
    \frac{\Delta A}{\Delta P} \geq C_{\ref{eqn:lb_step_main}}
     \end{equation}
\end{lemma}

We now give the symmmetrization argument that shows that the above bound in the special case
implies the bound in the general case.
\begin{proof}[Proof of \Cref{lemma:lb_step_main}]
    By the monotonicity of \(\cA\), we have \((\vc_t + \cE_t) \cup \{\vz_{t+1}\} \subseteq \vc_{t+1} + \cE_{t+1}\)
    and \(\vc_{t+1} + \alpha_{t+1} \cdot \cE_{t+1} \subseteq \conv{(\vc_t + \alpha_t \cdot \cE_t) \cup \{\vz_{t+1}\}}\).
    Without loss of generality we assume that \(\vc_t + \cE_t = B_2^d\) and \(\vz_{t+1} = 2 \cdot \ve_1\);
    we do this by applying a nonsingular affine transformation that maps \(\vc_t\) to the origin and \(\cE_t\) to \(B_2^d\), then apply a rotation
    that maps \(\vz_{t+1}\) to \(2 \cdot e_1\).
    Let \(\vc = \vc_{t+1}\), \(\cE = \cE_{t+1}\), and \(\alpha = \alpha_{t+1}\).
    Summarizing the conditions guaranteed by the monotonicity of \(\cA\), we have that \(B_2^d \cup \{ 2 \ve_1 \} \subseteq \vc + \cE\)
    and \(\vc + \alpha \cdot \cE \subseteq \conv{(\alpha_t \cdot B_2^d) \cup \{ 2 \ve_1 \}}\).

    To perform the reduction to the two-dimensional case, we apply a sequence of volume-preserving symmetrizations to the new inner and outer ellipsoids; these symmetrizations will also ensure that the update remains monotone.
    We will first apply two Steiner symmetrizations. The first of these Steiner symmetrizations transforms the ellipsoids so that their center lies on the \(\ve_1\)-axis. The second ensures that the ellipsoids have a semi-axis that is parallel to \(\ve_1\).
    Then, by a final symmetrization step we can transform the ellipsoids into bodies of revolution about \(\ve_1\).
    At that point it will suffice to consider the two-dimensional reduced case.

    Let \(\vc'\) be the projection of \(\vc\) onto the \(\ve_1\)-axis.
    The goal of the first symmetrization step is to transform \(\vc + \cE\) to \(\vc' + \cE'\) so that \(\vc'\) lies on the \(\ve_1\) axis. If \(\vc = \vc'\) then we do not need to do anything, otherwise we apply Steiner symmmetrization and consider \(S_{\vc - \vc'}(\vc+ \cE)\).
    By \Cref{lemma:steiner_ellipsoid} this is still an ellipsoid,
    and we also have that the center of \(S_{\vc - \vc'}(\vc + \cE)\) is actually \(\vc'\);
    thus we may write \(S_{\vc - \vc'}(\vc + \cE) = \vc' + \cE'\) for some \(\cE'\).
    Further, we have that \(S_{\vc - \vc'}(\vc + \alpha \cdot \cE) = \vc' + \alpha \cdot \cE'\), as the Steiner symmetrization acts similarly on the scaled version of \(\cE\).
    To show that the update is still monotone, we observe that \(\vc' + \cE' = S_{\vc - \vc'}(\vc + \cE) \subseteq S_{\vc - \vc'}( \conv{(\alpha_t \cdot B_2^d) \cup \{ 2 \ve_1 \}})\).
    But by \Cref{lemma:steiner_sym} and that \(\vc - \vc' \perp \ve_1\), \(\conv{(\alpha_t \cdot B_2^d) \cup \{ 2 \ve_1 \}}\) is invariant under the symmetrization \(S_{\vc - \vc'}\) and so we still have the inclusion \(\vc' + \alpha \cdot \cE' \subseteq \conv{(\alpha_t \cdot B_2^d) \cup \{ 2 \ve_1 \}}\).
    The `outer' inclusion
    \(B_2^d \cup \{ 2 \ve_1 \} \subseteq \vc' + \cE'\) follows in the same way.

    We now apply the second and final Steiner symmetrization. Let \(\vr\) be the rightmost point of \(\vc' + \cE'\) along \(\ve_1\); i.e. \(\vr = \arg\max_{\vr \in \vc' + \cE'} \inangle{\vr, \ve_1}\).
    Also let \(\vr'\) be its projection along the \(\ve_1\)-axis;
    if \(\vr = \vr'\) we again do not need to perform this symmetrization step, otherwise
    the Steiner symmetrization we apply is \(S_{\vr - \vr'}(\vc' + \cE')\).
    Since \(\vc'\) is at the midpoint of \(\vc' + \R (\vr - \vr')\) the center of the new ellipsoid is still \(\vc'\),
    so we may write \(S_{\vr - \vr'}(\vc' + \cE') = \vc' + \cE''\) and similarly
    \(S_{\vr - \vr'}(\vc' + \alpha \cdot \cE') = \vc' + \alpha \cdot \cE''\).
    Like for the previous symmetrization, the fact that \(\vr - \vr' \perp \ve_1\)
    means that both inclusions of the monotone update are preserved.
    Note finally that \(\vr'\) is the rightmost point of \(\vc' + \cE''\)
    and that the tangent plane of \(\vc + \cE''\) at \(\vr'\) is orthogonal to the line segment \(\overline{\vc' \vr'}\),
    so \(\ve_1\) is a semi-axis of \(\vc' + \cE''\).

    Our last transformation is a symmetrization of a different form, to turn \(\vc' + \cE''\) into a body of revolution.
    Let \(\sigma_1\) be the length of the semi-axis \(\ve_1\) of \(\cE''\), and \(\sigma_2, \ldots, \sigma_d\) be the lengths of the other semi-axes of \(\cE''\).
    We let \(\cE'''\) be the ellipsoid that has a \(\ve_1\) as a semi-axis of length \(\sigma_1\),
    and where every other semi-axis of \(\cE'''\) has length \(\sigma' \defeq \inparen{\prod_{i=2}^d \sigma_i}^{1/(d-1)}\).
    Clearly \(\vc' + \cE'''\) is now a body of revolution about \(\ve_1\) whose volume is the same as that of \(\vc' + \cE''\) (and hence also of \(\vc + \cE\)).
    Note that \(\vc' + \alpha \cE'''\) is also now a body of revolution.
    Since \(\sigma' \ge \min_{2 \le i \le d} \sigma_i\) we have \(B_2^d \cup \{ 2 \ve_1 \} \subseteq \vc' + \cE'''\), and correspondingly since \(\sigma' \le \max_{2 \le i \le d} \sigma_i\) we have that \(\vc' + \alpha \cdot \cE''' \subseteq \conv{(\alpha_t \cdot B_2^d) \cup \{ 2 \ve_1 \}}\).
    
    Clearly \(\vc' + \cE''', \vc + \alpha \cdot \cE'''\) now adhere to the reduced case of \Cref{defn:lb_reduced}.
    Since the update is monotone as well (and still \(\alpha \le 1/d\)) we can apply \Cref{lemma:lb_reduced}.
    As \(\vol(\cE''') = \vol(\cE)\), this means we have
    \[
    A_{t + 1} - A_{t} \geq C_{\ref{eqn:lb_step_main}} \cdot (P_{t+1} - P_{t})
    \]
    as desired.
\end{proof}

\begin{proof}[Proof of \Cref{lemma:lb_simplex_symmetry}]
For the first property, this is exactly the well-known fact that the best ellipsoidal rounding for the simplex \(\Delta_d\) (see e.g. \cite[Remark 1.1]{howard1997john}) has approximation factor \(d\).

Now we show the second property.
Again because the ball rounds the simplex \(\Delta_d\) with approximation factor \(d\), we have
\[\frac{1}{d} \cdot \Delta_d \subseteq B_2^d \subseteq \Delta_d\]
As a result of this, we have
\begin{align*}
    \log \inparen{\frac{\vol(\vc + \alpha \cdot \cE)}{\vol(B_2^d)}} &\le \log \inparen{\frac{\vol(\Delta)}{\vol(B_2^d)}} \\
    &\le \log \inparen{\frac{\vol(d \cdot B_2^d)}{\vol(B_2^d)}} \\
    &\le O(d \log d) 
\end{align*}
And so
\begin{align*}
\log \inparen{\frac{\vol(\cE)}{\vol(B_2^d)}} &= \log \inparen{\frac{\vol(\vc + \alpha \cdot \cE)}{\vol(B_2^d)}} + d \log\inparen{\frac{1}{\alpha}} \\
&\le O\inparen{d \log \inparen{\frac{1}{\alpha}}} & \text{as } \alpha \le \frac{1}{d}
\end{align*}
To establish the second property, it remains to show \(d \log(1/\alpha) \le O((1/\alpha) \log(d))\).
Observe that \(x \mapsto \frac{x}{\log x}\) is increasing for \(x \ge e\), so we have
\[\frac{d}{\log d} \le O\inparen{\frac{\nfrac{1}{\alpha}}{\log(\nfrac{1}{\alpha})}}\]
for all \(d \ge 2\) as \(d \le \nfrac{1}{\alpha}\).
Rearranging gives the desired inequality and thus the second property.
\end{proof}

\subsection{Analysis of the reduced case}
\label{sec:lb_2d}
In this section, we establish a lower bound on \(\frac{\Delta A}{\Delta P}\),
assuming we are in the `reduced case' defined in \Cref{defn:lb_reduced}.
Observe that in this case all relevant convex bodies \(\cE, \alpha \cE, c \cdot \ve_1 + \cE', c \cdot \ve_1 + \alpha' \cE', \conv{\alpha \cE \cup \{\vz\}}\) are all bodies of revolution about the \(x_1\)-axis, so to analyze the quantities involved we may instead look at any two-dimensional slice.
Accordingly we talk about the ellipses \(\cE, \alpha \cE, c + \cE', c + \alpha' \cE'\) in this two-dimensional slice,
where again \(\cE = B_2^2\), and \(c + \cE'\) and \(c + \alpha' \cE'\) are defined by
\begin{gather*}
c + \cE' = \inbraces{(x, y) \in \R^d \middle| \inparen{\frac{x-c}{a}}^2 + \inparen{\frac{y}{b}}^2 \leq 1}\\
c + \alpha' \cE' = \inbraces{(x, y) \in \R^d \middle| \inparen{\frac{x-c}{a}}^2 + \inparen{\frac{y}{b}}^2 \leq \alpha'^2}
\end{gather*}
for \(a, b > 0, c \in \R\).
We also use for convenience \(A = \frac{1}{\alpha}\) and \(A' = \frac{1}{\alpha'}\) so that \(\Delta A = A' - A\).
Also note in this reduced case we have by symmetry that
\[\Delta P = \log\inparen{\frac{\vol(c \cdot \ve_1 + \alpha' \cE')}{\vol(B_2^d)}} - \log\inparen{\frac{\vol(B_2^d)}{\vol(B_2^d)}} = \log(a \cdot b^{d-1})\]

Our lower bound in this reduced case is the following:
\begin{lemma} There exists a fixed constant \(C_{\ref{eq:lb_cpct}} > 0\) such that 
    \[\frac{\Delta A}{\Delta P} \geq \min\inbraces{C_{\ref{eq:lb_cpct}}, \frac{1}{10} \frac{A}{d}}\]    
\end{lemma}

Clearly this claim yields \Cref{lemma:lb_reduced} as a corollary, as by assumption in \Cref{lemma:lb_reduced} we have \(A \geq d\) and so we get \(\frac{\Delta A}{\Delta P} \geq \Omega(1)\).

\begin{figure}[h]
\centering
\includegraphics[width=0.70\textwidth]{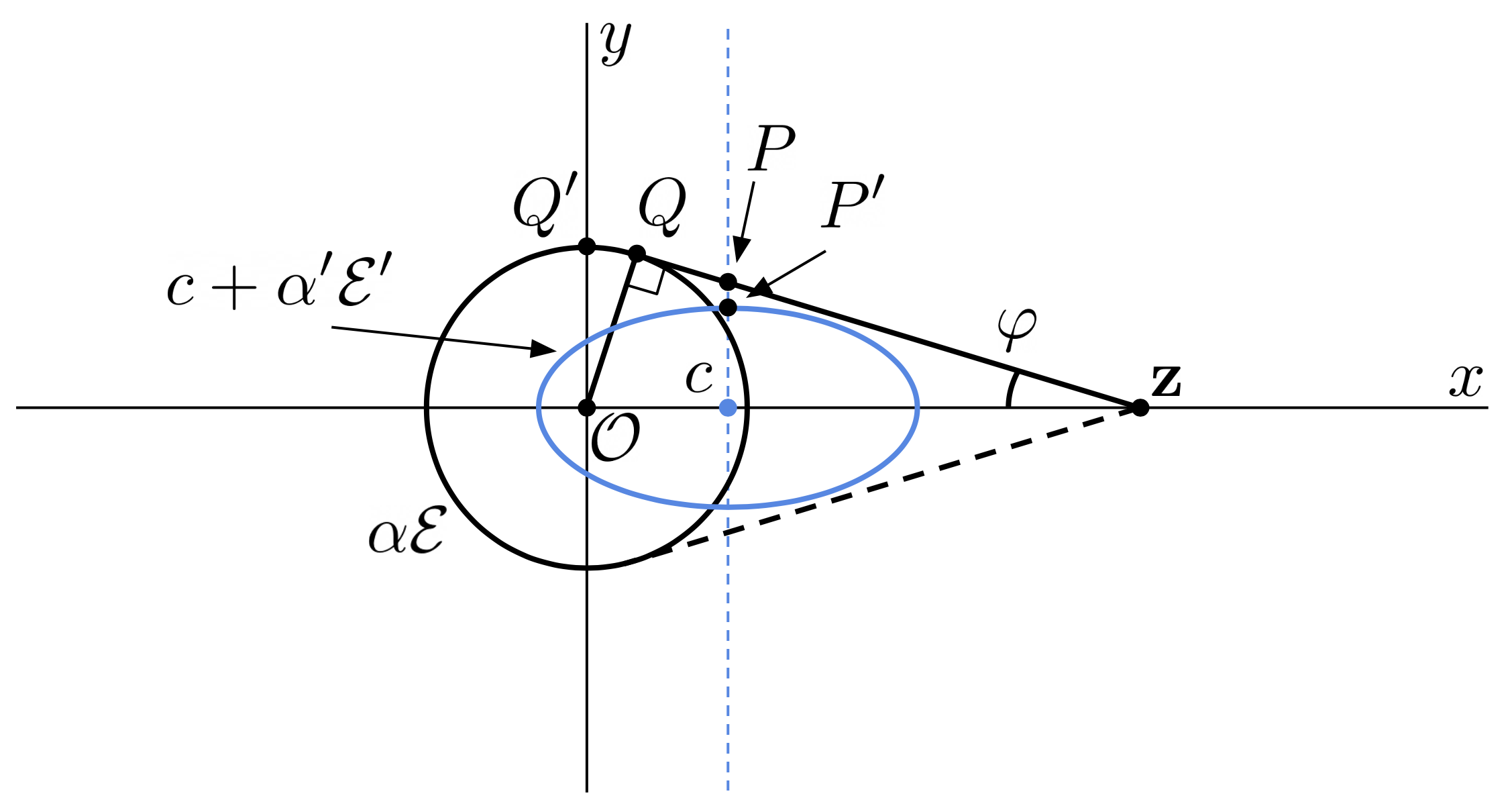}
\caption{The inner ellipses in the two-dimensional lower bound. \(\cO\) is the origin. The black solid circle is the previous inner ellipse \(\alpha \cE\),
and the blue solid circle is the next inner ellipse \(c + \alpha' \cE'\). The vertical dotted blue line \(x=c\) through the center \(c\) marks the location of the next inner ellipse on the \(x\)-axis. The new point is \(\vz = 2 \ve_1\),
and \(\overline{\vz Q}\) is one of the lines through \(\vz\) tangent to \(\alpha \cE\), with \(Q\) the point of tangency.
\(Q'\) is the intersection of \(\alpha \cE\) with the \(y\)-axis on the same side of the \(x\)-axis as \(Q\).
\(P'\) is the intersection of the line \(x=c\) with \(c + \alpha' \cE'\) on the same side as \(Q\),
and \(P\) is the intersection of this line with \(\overline{\vz Q}\).
We denote the angle \(\angle P \vz c\) with \(\varphi\).
}
\label{fig:lb_inner}
\end{figure}

The inner ellipses in this lower bound, and some relevant points used in the proof of this claim, are depicted in \Cref{fig:lb_inner}.

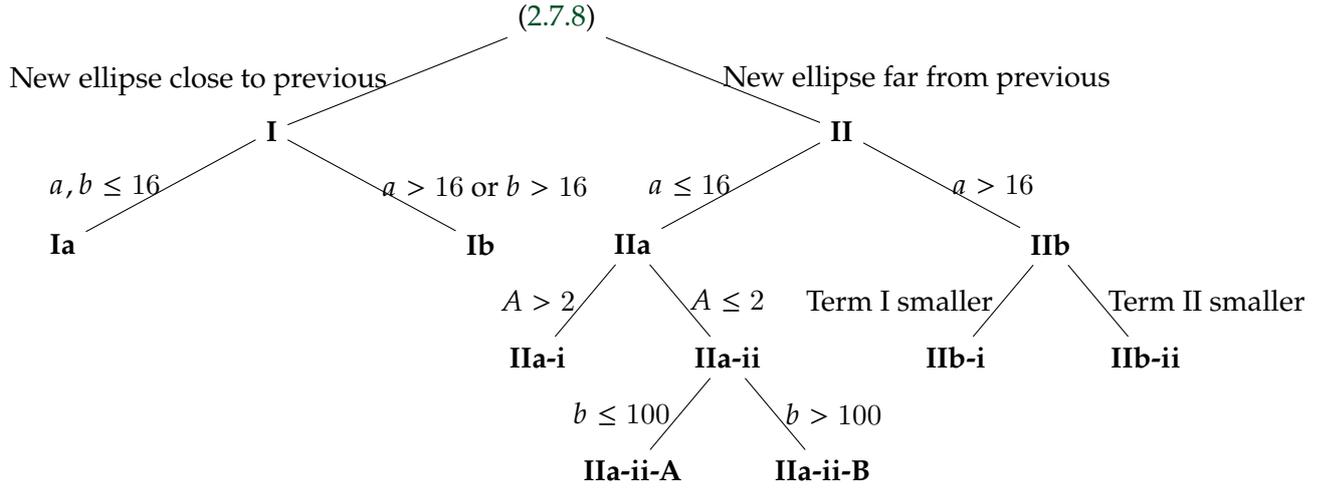
\begin{figure}[h]
\centering
\begin{tikzpicture}[
level 1/.style={sibling distance=75mm},
level 2/.style={sibling distance=55mm},
level 3/.style={sibling distance=25mm},
]
\tikzset{v/.append style={font=\bfseries}}
  \node {\eqref{eq:lb_5}}
    child {
        node[v] {I}
        child {
            node[v] {Ia}
            edge from parent
                node[left] {\(a, b \le 16\)}
        }
        child {
            node[v] {Ib}
            edge from parent
                node[right] {\(a > 16\) or \(b > 16\)}
        }
        edge from parent
            node[left] {New ellipse close to previous}
    }
    child {
        node[v] {II}
        child {
            node[v] {IIa}
            child {
                node[v] {IIa-i}
                edge from parent
                    node[left] {\(A > 2\)}
            }
            child {
                node[v] {IIa-ii}
                child {
                    node[v] {IIa-ii-A}
                    edge from parent
                        node[left] {\(b \le 100\)}
                }
                child {
                    node[v] {IIa-ii-B}
                    edge from parent
                        node[right] {\(b > 100\)}
                }
                edge from parent
                    node[right] {\(A \le 2\)}
            }
            edge from parent
                node[left] {\(a \le 16\)}
        }
        child {
            node[v] {IIb}
            child {
                node[v] {IIb-i}
                 edge from parent
                    node[left] {Term I smaller}
            }
            child {
                node[v] {IIb-ii}
                 edge from parent
                    node[right] {Term II smaller}
            }
            edge from parent
                    node[right] {\(a > 16\)}
        }
        edge from parent
            node[right] {New ellipse far from previous}
    };
\end{tikzpicture}
\caption{Tree of cases in the lower bound}
\label{fig:lb_cases}
\end{figure}

\begin{proof}
We establish this claim through a geometric argument that we break down by cases
(the logical tree of cases is visualized in \Cref{fig:lb_cases}).
First, as the new outer ellipse \(c + \cE'\) contains \(\cE = B_2^d\) we readily have that \(a, b \geq 1\).

As the rightmost point of the new outer ellipse must be to the right of \(\vz\), we have 
\begin{equation}\label{eq:lb_1}
c + a > 2
\end{equation}

As the leftmost point of the new inner ellipse must be to the right of the leftmost point of the previous inner ellipse, we have
\begin{equation}\label{eq:lb_4}
c \geq \alpha' a - \alpha
\end{equation}

\begin{lemma}\label{lemma:lb_1}
    We have \(\alpha' \cdot b \leq \alpha\), or equivalently \(A' \geq b \cdot A\).
\end{lemma}
\begin{proof}
    The geometry of this fact is visualized in \Cref{fig:lb_inner}.
    We overload notation so that \(c\) will also denote the point \((c, 0)\), the center of the new inner ellipse.
We denote \(\cO\) as the origin.
Let \(\overline{\vz Q}\) be one of the lines through \(\vz\) and tangent to \(\alpha \cE\),
with \(Q\) the point of tangency  (the choice of which line is arbitrary, in the figure we choose the one whose intersection with \(\alpha \cE\) is above the \(x\)-axis).
We let \(P\) be the intersection of the vertical line through \((c, 0)\) with \(\overline{\vz Q}\),
and \(P'\) be the intersection of this line with the ellipse \(c + \alpha' \partial \cE'\)
on the same side of the \(x\)-axis as \(Q\).

    Observe that \(\alpha' b = \overline{c P'}\) as the vertical semi-axis of the ellipse \(c + \alpha' \cE\), and \(\alpha = \overline{\cO Q'}\).
    Due to the fact that \(c + \alpha' \cE' \subseteq \conv{\alpha \cE \cup \{\vz\}}\),
    the projection of both sets onto the \(y\)-axis satisfies the same inclusion, and this gives the desired inequality.
\end{proof}

Observe that as \(c + \cE'\) must contain both the points \((-1, 0)\) and \((0, 2)\), we have
\begin{equation}
\label{eq:lb_a_3_2}
a \geq \frac{3}{2}
\end{equation}

Observe that we can split \(\Delta P\) into two terms: \[\Delta P = \underbrace{(d-1) \log b}_{\text{I}} + \underbrace{\log a}_{\text{II}}\]
First we show that if term I is larger, then we have a constant lower bound on \(\frac{\Delta A}{\Delta P}\). 
\begin{lemma}\label{lemma:lb_term2}
    If \((d-1) \log b \geq \log a\), then \(\frac{\Delta A}{\Delta P} \geq \frac{1}{2} \frac{A} {d}\).
\end{lemma}
\begin{proof}
    Under the assumption, we have \(\Delta P \leq 2 (d-1) \log b\).
    Combining this with \Cref{lemma:lb_1}, we have
    \begin{align*}
    \frac{\Delta A}{\Delta P} &\geq \frac{A' - A}{2 (d-1) \log b} \\
    &\geq \frac{b \cdot A - A}{2 (d-1) \log b} \\
    &= \frac{b-1}{2 \log b} \frac{A}{d-1} \\
    &\geq \frac{1}{2} \frac{A}{d} \\
    \end{align*}
    where the last line uses that \(\frac{x-1}{2 \log (x)} > \frac{1}{2}\) when \(x > 1\).
\end{proof}

In light of \Cref{lemma:lb_term2}, we can then assume in the sequel Term II is larger, meaning that
\begin{equation}\label{eq:lb_5}
\Delta P \leq 2 \log a
\end{equation}

\textbf{Case I }\textit{(New ellipse is close to the previous one)}. Assume that 
\begin{equation}\label{eq:lb_2}
    c + \alpha' \cdot a \leq \frac{11}{10} \alpha
\end{equation}
i.e. that the rightmost point of the new inner ellipse is to the left of \(\frac{11}{10}\).

\begin{lemma}\label{lemma:lb_close_1}
    In Case I, we have \(\Delta A \geq \frac{4}{11} A\).
\end{lemma}
\begin{proof}
    We prove this by cases.
    First, if we assume that \(\alpha' \leq \frac{\alpha}{2}\),
    we get \(A' \geq 2 A\) and \(\Delta A \geq A\).

    In the second case, we have \(\alpha' > \frac{\alpha}{2}\).
    We first use this to show \(c > 0\).
    By (\ref{eq:lb_4}) and (\ref{eq:lb_1}) we have \(\frac{c + \alpha}{\alpha'} \geq a > 2 - c\), so
    \(c (1 + \alpha') > 2 \alpha' - a > 0\) and so \(c > 0\).

    Using (\ref{eq:lb_2}) and that \(c > 0\), we have \(\alpha' a \leq \frac{11}{10} \alpha\).
    Thus \(A' \geq \frac{10}{11} a A\).
    By (\ref{eq:lb_a_3_2}) we get \(A' \geq \frac{15}{11} A\), and finally
    \(\Delta A \geq \frac{4}{11} A\).
\end{proof}

\begin{lemma}\label{lemma:lb_2}
    In Case I, we have \(A' \geq \frac{10}{21} a \cdot A\).
\end{lemma}
\begin{proof}
From (\ref{eq:lb_4}) we also get the weaker lower bound \(c \geq - \alpha\).
Combined with (\ref{eq:lb_2}), this gives \(\alpha' a \leq \frac{21}{10} \alpha\), which is equivalent to the desired inequality.
\end{proof}

We divide Case I into two sub-cases.

\textbf{Case Ia} \textit{(\(a, b \leq 16\))}.
First, assume that \(a, b \leq 16\).
Then \(\Delta P = \log(a \cdot b^{d-1}) \leq d \log(16)\),
and by \Cref{lemma:lb_close_1} we get \[\frac{\Delta A}{\Delta P} \geq \frac{4}{11 \log(16)} \frac{A}{d} \geq \frac{1}{10} \frac{A}{d}\]

\textbf{Case Ib} \textit{(\(a > 16\) or \(b > 16\))}. Now assume that either \(a\) or \(b\) is greater than \(16\).

Combining \Cref{lemma:lb_1} and \Cref{lemma:lb_2} together, we get \(A'^2 \geq \frac{10}{21} ab A^2\), or \(A' \geq \sqrt{\frac{10}{21} a b} \cdot A\).
We have \(\Delta A = A' - A \geq \inparen{\sqrt{\frac{10}{21} a b} - 1} A\). As \(ab \geq 16\) we get \(\sqrt{ab} \geq 2 \sqrt{\frac{21}{10}}\), so
\(\Delta A = \inparen{\sqrt{\frac{10}{21} a b} - 1} A \geq \sqrt{\frac{5}{21} ab} \cdot A\).

Using the analytic inequality that \(\log x \leq \sqrt{x}\) for all \(x > 0\), we get
\(\Delta P = \log a + (d-1) \log b \leq \sqrt{a} + (d-1) \sqrt{b} \leq d \cdot \sqrt{ab}\).

Combining these inequalites for \(\Delta A\) and \(\Delta P\), we obtain
\[\frac{\Delta A}{\Delta P} \geq \sqrt{\frac{5}{21}} \frac{A}{d} \geq \frac{4}{10} \frac{A}{d}\]

\textbf{Case II }\textit{(New ellipse is far from the previous one)}. Assume that
\begin{equation}\label{eq:lb_3}
    c + \alpha' \cdot a > \frac{11}{10} \alpha
\end{equation}

\begin{lemma}\label{lemma:lb_3}
    We have 
    \begin{equation}\label{eq:lb_trig_b}
    \alpha' b \leq (2 - c) \cdot \frac{\frac{\alpha}{2}}{\sqrt{1 - \inparen{\frac{\alpha}{2}}^2}}
    \end{equation}
\end{lemma}
\begin{proof}
Again, the proof of this claim is pictured in \Cref{fig:lb_inner},
where we construct the points in the same way as in the proof of \Cref{lemma:lb_1}.
Let \(\angle P \vz c\) be denoted by \(\varphi\).
Note that the angle \(\angle P c \vz\) is a right angle, and so \(\tan \varphi = \frac{\overline{c P}}{\overline{c v}}\). The line \(\overline{c \vz}\) has length \(2 - c\), so we get \(\overline{c P} = (2 - c) \tan \varphi\).
The segment \(\overline{c P'}\), of length \(\alpha' b\), is contained within the segment \(\overline{c P}\), and so
\(\alpha' b \leq (2 - c ) \tan \varphi\).
    
Observe that the angle \(\angle \cO Q \vz\) is also a right angle.
Further, clearly the length of \(\overline{\cO Q}\) is \(\alpha\) and the length of \(\overline{0 v}\) is \(2\).
Since we have that \(\varphi\) is also the angle \(\angle Q \vz \cO\), we get \(\sin \varphi = \frac{\alpha}{2}\).
Now using the standard trigonometric identity that \(\tan \varphi = \frac{\sin \varphi}{\sqrt{1 - \sin^2 \varphi}}\) for $\varphi\in[-\pi/2,\pi/2]$,
we get the desired inequality.
\end{proof}

We split Case II into several sub-cases, as for Case I. 

\textbf{Case IIa} \textit{(\(a \leq 16\))}. First, we look at the case where \(a \leq 16\).

\textbf{Case IIa-i} \textit{(\(A \geq 2\))}. Assume \(A \geq 2\).
\begin{lemma}\label{lemma:lb_4}
    When \(A \geq 2\), we have \(\Delta A \geq \frac{1}{2}\).
\end{lemma}
\begin{proof}
Adding (\ref{eq:lb_4}) and (\ref{eq:lb_3}) together gives \(c > \frac{1}{20} \alpha\).
Using this in (\ref{eq:lb_trig_b}) and rearranging using the definitions of \(A\) and \(A'\) yields
\[A' \geq b A \cdot \frac{1}{1-\frac{1}{40} \frac{1}{A}} \inparen{1 - \frac{1}{A^2}}\]
and therefore we get the inequality
\[\Delta A \geq b A \cdot \frac{1}{1-\frac{1}{40} \frac{1}{A}} \inparen{1 - \frac{1}{A^2}} - 1\]
and using \(b \geq 1\), we obtain
\[\Delta A \geq A \cdot \frac{1}{1-\frac{1}{40} \frac{1}{A}} \inparen{1 - \frac{1}{A^2}} - 1\]
To prove the claim, it suffices to show the right hand side exceeds \(\frac{1}{2}\) when \(A \geq 2\). Upon rearranging, this is equivalent to the inequality \(A - \frac{77}{80} \frac{1}{A} \geq \frac{3}{2}\) when \(A \geq 2\).
\end{proof}
Combining the assumption that \(a \leq 16\) with (\ref{eq:lb_5}) and \Cref{lemma:lb_4} yields \(\frac{\Delta A}{\Delta P} \geq \frac{1}{4 \log 16} \geq \frac{1}{20}\).

\textbf{Case IIa-ii} \textit{(\(A \leq 2\))}. Next, we look at the other case where \(A \leq 2\).

\textbf{Case IIa-ii-A} \textit{(\(b \leq 100\))}. Now we look at the case where \(b \leq 100\).
\begin{lemma}
    If \(a \leq 16, A \leq 2, b \leq 100\), then there is \(C_{\ref{eq:lb_cpct}} > 0\) such that
    \begin{equation}\label{eq:lb_cpct}
        \frac{\Delta A}{\Delta P} \geq C_{\ref{eq:lb_cpct}}
    \end{equation}
\end{lemma}
\begin{proof}
    We show this by a compactness argument.
    By (\ref{eq:lb_5}) and the assumption that \(a \leq 16\) we have \(\frac{\Delta A}{\Delta P} \geq \frac{\Delta A}{2 \log 16}\).
    Now, observe that next outer and inner ellipsoids \(c + \cE'\) and \(c + \alpha' \cE'\) are fully determined by the parameters
    \(a, b, c, A, A'\).
    Further, we assume without loss of generality that \(A'\) is a function of the other parameters.
    This is because when \(A'\) is decreased as much as possible while preserving the monotonicity of the update, \(\Delta A = A' - A\) only decreases.
    To show a lower bound on \(\Delta A\) it then suffices to only do so in this hardest case.
    
    Note that we have \(1 \leq a \leq 16, 1 \leq b \leq 100, -1 \leq c \leq 2\), and \(1 \leq A \leq 2\); thus all the parameters defining the next inner and outer ellipsoids are bounded.
    Observe that \(\Delta A\) is a continuous function of these parameters, and as a continuous function of a compact set it attains its minimum. Finally, we argue that it is impossible for the minimum of \(\Delta A\) to be zero, and so the minimum is some strictly positive constant \(\frac{C_{\ref{eq:lb_cpct}}}{2 \log 16}\), which suffices to prove the claim.
    
    The following argument only concerns the inner ellipsoids, and can be pictured in \Cref{fig:lb_inner}. If \(c = 0\), then as the leftmost point of \(c + \alpha' \cE'\) must be to the right of the leftmost point of \(\alpha \cE'\),
    we have \(\alpha \geq \alpha' a\). But by (\ref{eq:lb_a_3_2}), we have \(\alpha \geq \frac{3}{2} \alpha'\), so \(\alpha' < \alpha\) and \(\Delta A > 0\).
    If \(c \neq 0\) then the vertical semi-axis of \(c + \alpha' \cE'\) must have length strictly less than \(\alpha\),
    and so \(\alpha > \alpha' b\). As \(b \geq 1\), this also gives \(\alpha > \alpha'\) and again \(\Delta A > 0\).
\end{proof}

\textbf{Case IIa-ii-B} \textit{(\(b > 100\))}.
Observe that the horizontal axis of the next inner ellipsoid \(c + \alpha' \cE'\) must be contained within the interval \([-\alpha, 2]\), thus we have that \(2 + \alpha \geq 2 b \alpha' > 200 \alpha'\).
Using the definitions of \(A, A'\) this is equivalent to \(2 + \frac{1}{A} > \frac{200}{A'}\), i.e. \(A' > \frac{200}{1 + \frac{1}{A}}\). As \(A \geq 1\) we get \(2 + \frac{1}{A} \leq 3\), and so \(A' \geq \frac{200}{3}\).

As \(A \leq 2\), we obtain \(\Delta A = A' - A \geq \frac{194}{3}\).
Now by \eqref{eq:lb_5} and that \(a \leq 16\) we have 
\[\frac{\Delta A}{\Delta P} \geq \frac{\Delta A}{2 \log 16} \geq \frac{194}{6 \log 16} \geq 11\]

\textbf{Case IIb} \textit{(\(a > 16\))}. Now, we examine the case where \(a > 16\).
Scaling (\ref{eq:lb_4}) by \(\frac{11}{10}\) and adding it to (\ref{eq:lb_3}), we have \(\frac{21}{10} c - \frac{1}{10} \alpha' a \geq 0\), i.e. \(c > \frac{1}{21} \alpha' a\).
Using this in (\ref{eq:lb_trig_b}), using \(b \geq 1\), using the definitions of \(A\) and \(A'\) and rearranging, we obtain
\[A' \geq A \cdot \frac{1}{1 - \frac{1}{42} \cdot \frac{a}{A'}} \cdot \sqrt{1 - \inparen{\frac{\alpha}{2}}^2}\]

Using the inequalities \(\sqrt{1-\inparen{\frac{x}{2}}^2} \geq 1 - \frac{x^2}{7}\) for \(0 \leq x \leq 1\) and \(\frac{1}{1-x} \geq 1 + x\) for \(0 \leq x \leq 1\), we have
\[A' \geq A \inparen{1- \frac{\alpha^2}{7}} \inparen{1 + \frac{1}{42} \frac{a}{A'}}\]
and thus
\(A'^2 - A \cdot A' \inparen{1 - \frac{\alpha^2}{7}} - \frac{1}{42}  \inparen{1 - \frac{\alpha^2}{7}} \geq 0\),
which implies by the quadratic formula that \begin{align*}
A' &\geq \frac{A \inparen{1 - \frac{\alpha^2}{7}} + \sqrt{A^2 \inparen{1 - \frac{\alpha^2}{7}} + \frac{4}{42} a A \inparen{1 - \frac{\alpha^2}{7}}}}{2} \\
&= A \inparen{1 - \frac{\alpha^2}{7}} \cdot \frac{1 + \sqrt{1 + \frac{2}{21} \frac{a}{A \inparen{1 - \frac{\alpha^2}{7}}}}}{2} \\
&= A \inparen{1 - \frac{\alpha^2}{7}} \cdot \inparen{1 + \frac{\sqrt{1 + \frac{2}{21} \frac{a}{A \inparen{1 - \frac{\alpha^2}{7}}}} - 1}{2}}
\end{align*}
Using the inequality \(\sqrt{1+x}-1 \geq \frac{2}{5} \min\inbraces{x, \sqrt{x}}\) for all \(x \geq 0\), we get that
\begin{equation}\label{eq:lb_2b}
A' \geq A \inparen{1 - \frac{\alpha^2}{7}} \inparen{1+ \frac{1}{5}\min\inbraces{\smash[b]{\underbrace{\frac{2}{21} \frac{a}{A \inparen{1 - \frac{\alpha^2}{7}}}}_{\text{I}}, \underbrace{\sqrt{\frac{2}{21} \frac{a}{A \inparen{1 - \frac{\alpha^2}{7}}}}}_{\text{II}}}}}
\end{equation}
\\ 

To finish this case, we show the lower bound in the case where either term in the \(\min\) of (\ref{eq:lb_2b}) is the smaller term.

\textbf{Case IIb-i} \textit{(Term I in (\ref{eq:lb_2b}) is smaller)}. In this case, (\ref{eq:lb_2b}) is equivalent to
\begin{align*}
    A' &\geq A \inparen{1 - \frac{\alpha^2}{7}} + \frac{2}{105} a
\end{align*}
Using the definition of \(A\), we have 
\(A' \geq A - \frac{1}{7A} + \frac{2}{105} a\), and so \(\Delta A \geq - \frac{1}{7A} + \frac{2}{105} a\).
As \(A \geq 1\), we have \(\Delta A \geq \frac{2}{105} a - \frac{1}{7}\).

Now by \eqref{eq:lb_5}, we get
\[\frac{\Delta A}{\Delta P} \geq \frac{\frac{2}{105} a - \frac{1}{7}}{2 \log a}\]
Now, we complete this case by noticing the right hand side is at least \(\frac{1}{35}\) when \(a > 16\).

\textbf{Case IIb-ii} \textit{(Term II in (\ref{eq:lb_2b}) is smaller)}.  In this case, (\ref{eq:lb_2b}) is equivalent to
\begin{equation}\label{eq:lb_c2b2_1}
    A' \geq A \inparen{1 - \frac{\alpha^2}{7}} + \sqrt{\frac{2}{525} a A \inparen{1 - \frac{\alpha^2}{7}}}
\end{equation}
Using the definition of \(A\) and that \(A \geq 1\), we have \(A \inparen{1 - \frac{\alpha^2}{7}} = A - \frac{1}{7A} \geq - \frac{1}{7}\). Using this and the definition of \(\Delta A\) in \eqref{eq:lb_c2b2_1}, we have \[
\Delta A \geq - \frac{1}{7} + \sqrt{\frac{2}{525} a A \inparen{1 - \frac{\alpha^2}{7}}}
\]
Further, as \(0 \leq \alpha \leq 1\) we get \(A \inparen{1 - \frac{\alpha^2}{7}} = \frac{1}{\alpha} - \frac{\alpha}{7} \geq \frac{6}{7}\), so
\[\Delta A \geq - \frac{1}{7} + \sqrt{\frac{12}{3675} a}
\]
Now by \eqref{eq:lb_5}, we get
\[\frac{\Delta A}{\Delta P} \geq \frac{- \frac{1}{7} + \sqrt{\frac{12}{3675} a}}{2 \log a}\]
We finish with the fact that the right hand side is at least \(\frac{1}{65}\) when \(a > 16\).
\end{proof}

\section{Details of analysis in \texorpdfstring{\Cref{sec:two_dim_update}}{Section~\ref{sec:two_dim_update}}}
Here, we give the details for the outstanding claims in \Cref{sec:two_dim_update}. We stay in the context from that section and reuse notation (specifically, the definition of parameters in (\ref{eqn:update_params})). 

We begin with some well-known bounds on \(e^x\).

\begin{lemma}\label{lemma:wk_e} Recall the following well-known inequalities regarding \(e^x\).
\begin{enumerate}
\item \label{item:wk_1} \(1 + x \leq e^x\) for all \(x \in \R\);
\item \label{item:wk_2} \(1 + x + \frac{x^2}{2} \leq e^x\) for \(x \geq 0\).
\end{enumerate}
\end{lemma}

We will also use a more specialized upper bound on \(e^x\).
\begin{lemma}\label{lemma:custom_exp_upper_bound} For \(0 \leq x \leq \frac{4}{3}\), we have
\(e^x \leq 1 + x + \frac{x^2}{2} + \frac{x^3}{4}\).
\end{lemma}
\begin{proof}
Using the Taylor series for \(e^x\) about \(0\),
we get \[
\left(1 + x + \frac{x^2}{2} + \frac{x^3}{4}\right) - e^x = \frac{x^3}{12} - \sum_{k=4}^\infty \frac{x^k}{k!}
= x^3 \left(\frac{1}{12} - \sum_{k=4}^\infty \frac{x^{k-3}}{k!}\right).\]
Clearly \(x^3 \geq 0\) for \(x \geq 0\), so it remains to show \(\frac{1}{12} - \sum_{k=4}^\infty \frac{x^{k-3}}{k!} \geq 0\) for \(0 \leq x \leq \frac{4}{3}\). 
Note that \(\sum_{k=4}^\infty \frac{x^{k-3}}{k!}\) is increasing (the derivative is clearly positive when \(x \ge 0\)),
and we finish by noting \[\frac{1}{12} - \left.\sum_{k=4}^\infty \frac{x^{k-3}}{k!} \right|_{x = \frac{4}{3}} = \left.\frac{1 + x + \frac{x^2}{2}+\frac{x^3}{4} - e^x}{x^3} \right|_{x = \frac{4}{3}} > 0.\]
\end{proof}

Now we show some facts used in \Cref{lemma:pty_outer}, which \Cref{lemma:update_step_outer} reduces to.
The proof of \Cref{lemma:pty_outer} will reduce to the following analytic inequality.
\begin{lemma}\label{lemma:outer_t_ineq}
  For all \(\gamma \geq 0\),
  \[\frac{(e^\gamma - 1)^2}{e^{2\gamma} - (1 + \frac{\gamma}{4})^2} \leq \frac{3}{2} \gamma.\]
\end{lemma}
\begin{proof}
For the numerator of the left hand side, we have \((e^\gamma - 1)^2 = e^{2\gamma} - 2 e^\gamma + 1 \leq e^{2\gamma} - 2\gamma - 1\) using \Cref{lemma:wk_e}-(\ref{item:wk_1}), \(1 + x \leq e^x\).
Further, \(e^\gamma \geq 1 + \frac{\gamma}{4}\) implies \(e^{2 \gamma} - (1 + \frac{\gamma}{4})^2 \geq 0\), so after multiplying both sides by \(e^{2 \gamma} - (1 + \frac{\gamma}{4})^2\) and rearranging, it suffices to show
\[\frac{3}{2} \gamma \left(1 + \frac{\gamma}{4}\right)^2 - 1 - 2\gamma \leq \left(\frac{3}{2} \gamma - 1\right) e^{2 \gamma}.\]
We split this into two cases, based on the value of \(\gamma\).
If \(\gamma \geq \frac{2}{3}\), then the right hand side is at least \(\left(\frac{3}{2} \gamma - 1\right) ( 1 + 2 \gamma + 2 \gamma^2)\) using \Cref{lemma:wk_e}-(\ref{item:wk_2}), \(1 + x + \frac{x^2}{2} \leq e^x\),
so it is sufficient to show \(\left(\frac{3}{2} \gamma - 1\right) ( 1 + 2 \gamma + 2 \gamma^2) \geq \frac{3}{2} \gamma \left(1 + \frac{\gamma}{4}\right)^2 - 1 - 2 \gamma\).
Expanding both sides, this is equivalent to showing \(\frac{\gamma^2}{4} + \frac{93}{32} \gamma^3 \geq 0\), which is clearly true for \(\gamma \geq 0\).

If \(\gamma < \frac{2}{3}\), then we use \Cref{lemma:custom_exp_upper_bound} to lower bound the right hand side with \(\left(\frac{3}{2} \gamma - 1 \right) (1 + 2 \gamma + 2 \gamma^2 + 2 \gamma ^3)\),
so it is sufficient to show  \(\left(\frac{3}{2} \gamma - 1\right) ( 1 + 2 \gamma + 2 \gamma^2 + 2 \gamma^3) \geq \frac{3}{2} \gamma \left(1 + \frac{\gamma}{4}\right)^2 - 1 - 2\gamma\).
Similar to before, after expanding both sides this is equivalent to showing \(\frac{\gamma^2}{4} + \frac{93}{32} \gamma^3 + \frac{3}{2} \gamma^4 \geq 0\), which is true for \(\gamma \geq 0\).
\end{proof}

We also show some relations between the parameters in the update step.
Recall that we assumed \(\alpha \leq \frac{1}{2}\).
\begin{lemma}\label{lemma:outer_int}
We have
\begin{enumerate}
\item \(b \leq 1 + \frac{\gamma}{4}\) \label{item:outer_int_1}
\item \(b \leq a\) \label{item:outer_int_2}
\item \(\frac{(a-1)^2}{a^2-b^2} \leq 1\) \label{item:outer_int_3}
\item \(b^2 \geq 1 + \alpha - \alpha'\) \label{item:outer_int_4}
\end{enumerate}
\end{lemma}
\begin{proof}
We start by showing (\ref{item:outer_int_1}). As \(\frac{1}{\alpha'} = \frac{1}{\alpha} + 2\gamma\), we have \(\alpha = \alpha' + 2 \gamma \alpha \alpha'\).
Thus \(b = 1 + \gamma \alpha \alpha' \leq 1 + \gamma \alpha^2 \leq 1 + \gamma/4\) as \(\alpha' \leq \alpha \leq \frac{1}{2}\).

For (\ref{item:outer_int_2}), observe that \(b \leq 1 + \frac{\gamma}{4} \leq 1 + \gamma \leq e^\gamma = a\) using \Cref{lemma:wk_e}-(\ref{item:wk_1}), \(1+x \leq e^x\), so \(a \geq b\).

To show (\ref{item:outer_int_3}), we first argue it is sufficient to show \(1 + b^2 \leq 2 a\). As a consequence \(-2a + 1 \leq - b^2\), so \((a-1)^2 \leq a^2 - b^2\).
Because \(b \geq 1\) by \Cref{lemma:update_params}-(\ref{item:claim_update_params_2}), from (\ref{item:outer_int_2}) we can say that \(a^2 - b^2 \geq 0\), so that \(\frac{(a-1)^2}{a^2 - b^2} \leq 1\).

Now to show \(1 + b^2 \leq 2a\), we write as a series in terms of \(\gamma\).
On the left hand side using (\ref{item:outer_int_1}), we have \(1 + b^2 \leq 1 + \left(1 + \frac{\gamma}{4}\right)^2 = 2 + \frac{\gamma}{2} + \frac{\gamma^2}{4}\).
Further, by \Cref{lemma:wk_e}-(\ref{item:wk_2}), \(e^x \geq 1 + x + \frac{x^2}{2}\), we have that \(2a \geq 2 + 2\gamma + \gamma^2\).
Clearly \(2 + \frac{\gamma}{2} + \frac{\gamma^2}{4} \leq 2 + 2\gamma + \gamma^2\) when \(\gamma \geq 0\), so we are finished.

For (\ref{item:outer_int_4}), we have by definition that \(b^2 = 1 + \alpha - \alpha' + \frac{(\alpha - \alpha')^2}{4}\), so \(b^2 \geq 1 + \alpha - \alpha'\).
\end{proof}

As the proof of \Cref{lemma:update_step_outer} shows, that claim reduces to the following inequality.
\begin{lemma}\label{lemma:pty_outer}
We have \begin{equation}\label{eqn:pty_outer_eqn} 
        c^2 \leq \frac{b^2 - 1}{b^2} \cdot (a^2 - b^2).
    \end{equation}
\end{lemma}
\begin{proof}
We first upper bound \(c\) to reduce the number of variables in (\ref{eqn:pty_outer_eqn}).
As \(b = 1 + \frac{\alpha - \alpha'}{2}\), we have \(2 (b - 1) = \alpha - \alpha'\) and so
\(\alpha = \alpha' + 2(b-1)\).
Thus \[c = - \alpha + \alpha' \cdot a = -(\alpha' + 2(b-1)) + \alpha' \cdot a = \alpha' \cdot (a - 1) + 2(1-b) a.\]
As \(b \geq 1\) by \Cref{lemma:update_params}-(\ref{item:claim_update_params_2}), we have that \(2(1-b) \leq 0\) and therefore
\begin{equation}\label{eqn:pty_outer_uc} 
    c \leq \alpha' \cdot (a-1). \nonumber
\end{equation}
Using this in (\ref{eqn:pty_outer_eqn}), it suffices to show \(\alpha'^2 (a-1)^2 \leq \frac{b^2 - 1}{b^2} \cdot (a^2 - b^2)\),
which rearranges to
\begin{equation}\label{eqn:pty_outer_eqn1} 
\frac{b^2 - 1}{\alpha'^2} \geq \frac{(a-1)^2 b^2}{a^2 - b^2}. \nonumber
\end{equation}

Using \Cref{lemma:outer_int}-(\ref{item:outer_int_4}), this reduces to \begin{equation}\label{eqn:pty_outer_eqn2}
    \frac{\alpha - \alpha'}{\alpha'^2} \geq \frac{(a-1)^2 b^2}{a^2 - b^2}.
\end{equation}
The left hand side of (\ref{eqn:pty_outer_eqn2}) equals \(\frac{1}{\alpha'} \left(\frac{\alpha}{\alpha'} - 1\right)\).
    Because \(\frac{1}{\alpha'} = \frac{1}{\alpha} + 2 \gamma\), we have \(\frac{\alpha}{\alpha'} - 1 = 2 \gamma \alpha\),
so \(\frac{1}{\alpha'}\left(\frac{\alpha}{\alpha'} - 1\right) = 2\gamma \cdot \frac{\alpha}{\alpha'} = 2\gamma \cdot (1 + 2 \gamma \alpha)\).
So it is sufficient to show
\begin{equation}\label{eqn:pty_outer_eqn3}
2\gamma (1 + 2 \gamma \alpha)\geq \frac{(a-1)^2 b^2}{a^2 - b^2}.
\end{equation}

Now, we will eliminate the other variables in this inequality to transform it
into a statement involving only \(\gamma\).
We have
\begin{align*}
\frac{(a-1)^2 b^2}{a^2 - b^2} &= \frac{(a-1)^2}{a^2 - b^2} \left(1 + 2 \alpha \alpha' \gamma + \alpha'^2 \alpha^2 \gamma^2\right)\\
&\leq \frac{(a-1)^2}{a^2 - b^2} + \frac{\gamma}{2} + \alpha \frac{\gamma^2}{8},
\end{align*}
where the first line uses that \(b = 1 + \alpha \alpha' \gamma\), and
the second line inequality follows from \Cref{lemma:outer_int}-(\ref{item:outer_int_3}) and the fact that \(\alpha' \leq \alpha \leq \frac{1}{2}\).
Thus we can reduce (\ref{eqn:pty_outer_eqn3}) to \(\frac{(a-1)^2}{a^2-b^2} \leq \frac{3}{2} \gamma + \frac{31}{8} \alpha \gamma^2\), or further to

\begin{equation}\label{eqn:pty_outer_eqn4}
\frac{(a-1)^2}{a^2 - b^2} \le \frac{3}{2} \gamma.
\end{equation}

Using \Cref{lemma:outer_int}-(\ref{item:outer_int_1}) and that \(a = e^\gamma\), we have \(\frac{(a-1)^2}{a^2 - b^2} \leq \frac{(e^\gamma - 1)^2}{e^{2\gamma} - (1 + \frac{\gamma}{4})^2}\), so finally (\ref{eqn:pty_outer_eqn4}) reduces to
\[\frac{(e^\gamma - 1)^2}{e^{2\gamma} - (1 + \frac{\gamma}{4})^2} \le \frac{3}{2} \gamma,\]
which is proved in \Cref{lemma:outer_t_ineq}.
\end{proof}

Recall in the proof of \Cref{lemma:translate_ellipse_angle} we defined \(\ell_1 = \frac{1}{c+a}, \ell_2 = \sqrt{\frac{1}{\alpha^2} - \frac{1}{(c+a)^2}}, r = \frac{a^2 \ell_1^2}{b^2 \ell_2^2}\).
That claim reduces to the following.
\begin{lemma}\label{lemma:pty_inner}
We have
\[a - \alpha' \cdot a \sqrt{\frac{1+r}{r}} \geq 0.\]
\end{lemma}
\begin{proof}
As by definition \(a \geq 0\), it suffices to show
\begin{equation}\label{eqn:pty_inner_3}
    \alpha'^2 \cdot \left(\frac{1}{r}+1\right) \leq 1.
\end{equation}
Observe that \(\ell_2^2 = \frac{1}{\alpha^2} - \ell_1^2\),
so we can write \(\frac{1}{r} = \frac{b^2}{a^2} \left(\frac{1}{\alpha^2 \ell_1^2} - 1\right)\), and hence rewrite (\ref{eqn:pty_inner_3}) as
\[\alpha'^2 \left(1 + \frac{b^2}{a^2} \left(\left(\frac{c+a}{\alpha}\right)^2 - 1\right)\right) \leq 1.\]

Multiplying both sides by \(\frac{\alpha^2}{\alpha'^2}\) and rearranging, this is equivalent to
\begin{equation}\label{eqn:pty_inner_4}
\frac{b^2}{a^2} \left((c+a)^2 - \alpha^2\right) \leq \frac{\alpha^2}{\alpha'^2} - \alpha^2.
\end{equation}
Now, by definition of \(c\) we can write \(c + a = a(1+\alpha') - \alpha\), so that \((c+a)^2 - \alpha^2 = a^2 (1+\alpha')^2 - 2 \alpha a (1 + \alpha')\).
Thus, (\ref{eqn:pty_inner_4}) is equivalent to
\[\frac{b^2}{a^2} (a^2 ( 1 + \alpha')^2 - 2 \alpha a (1+\alpha')) \leq \frac{\alpha^2}{\alpha'^2} (1 + \alpha')(1 - \alpha').\]
Dividing by \(1 + \alpha'\) and simplifying the left hand side, this is equivalent to
\[
b^2 \left(1 + \alpha' - \frac{2 \alpha}{a}\right) \leq \frac{\alpha^2}{\alpha'^2}(1-\alpha'),
\]
which we show in \Cref{lemma:pty_inner_main}.
\end{proof}

\begin{lemma}\label{lemma:pty_inner_main}
We have
    \[
b^2 \left(1 + \alpha' - \frac{2 \alpha}{a}\right) \leq \frac{\alpha^2}{\alpha'^2}(1-\alpha').
\]
\end{lemma}
\begin{proof}
Using \Cref{lemma:wk_e}-(\ref{item:wk_1}), \(e^{-x} \geq 1 - x\); and the fact that by definition \(\frac{1}{a} = e^{-\gamma}\), it suffices to show
\[
b^2 \left(1 + \alpha' - 2 \alpha (1-\gamma)\right) \leq \frac{\alpha^2}{\alpha'^2}(1-\alpha'). \nonumber
\]
Using \Cref{lemma:outer_int}-(\ref{item:outer_int_4}), this reduces further to
\begin{equation}\label{eqn:pim_1}
    (1 + \alpha - \alpha')(1 + \alpha' - 2 \alpha ( 1 -  \gamma)) \leq \frac{\alpha^2}{\alpha'^2} (1 - \alpha').
\end{equation}
We expand both sides of this inequality into polynomials involving \(\gamma\) and \(\alpha\), and then analyze the resulting expression.
Using the definition of \(\alpha'\), we have \(\alpha ' = \frac{\alpha}{1+2\gamma \alpha}\), and thus
\(1 + \alpha' = \frac{1 + 2 \gamma \alpha + \alpha}{1 + 2 \gamma \alpha}\) and \(1 - \alpha' = \frac{1 + 2 \gamma \alpha - \alpha}{1 + 2 \gamma \alpha}\).
We also have \(\frac{\alpha}{\alpha'} = 1 + 2 \gamma \alpha\), and finally \(\alpha - \alpha' = \frac{2 \gamma \alpha^2}{1 + 2 \gamma \alpha}\).
Substituting these equalities into (\ref{eqn:pim_1}), we obtain the equivalent inequality
\[
\left(\frac{1+2 \gamma \alpha + \gamma \alpha^2}{1 + 2 \gamma \alpha}\right) \left(\frac{1 + 2 \gamma \alpha + \alpha}{1 + 2 \gamma \alpha} - 2 \alpha (1 - \gamma)\right) \leq (1 + 2 \gamma \alpha)^2 \left(\frac{1 + 2 \gamma \alpha - \alpha}{1 + 2 \gamma \alpha}\right).
\]
Multiplying both sides by \((1 + 2 \gamma \alpha)^2\) and rearranging the terms so that they are all on the same side, we get
\[(1 + 2 \gamma \alpha)^3 (1+2 \gamma \alpha - \alpha)-\left(1 + 2 \gamma \alpha + \gamma \alpha^2\right) (1+2 \gamma \alpha + \alpha - 2 \alpha (1- \gamma) (1 + 2 \gamma \alpha)) \geq 0.\]
Next, we expand this inequality:
\begin{align*}
    16 \alpha ^4 \gamma^4-16 \alpha ^4 \gamma^3+8 \alpha ^4 \gamma^2+24 \alpha ^3 \gamma^3-12 \alpha ^3 \gamma^2+12
   \alpha ^2 \gamma^2+2 \alpha ^3 \gamma-2 \alpha ^2 \gamma+2 \alpha  \gamma \geq 0
\end{align*}
As \(\gamma \alpha \geq 0\), we can divide both sides of this inequality by \(2 \gamma \alpha\).
Grouping by powers of \(\alpha\), we obtain
\[4 \alpha ^3 \gamma \left(2 \gamma^2-2 \gamma+1\right)+\alpha ^2 \left(12 \gamma^2-6 \gamma+1\right)+\alpha  (6
   \gamma-1)+1 \geq 0.\]
Upon inspection, both quadratics \(2 \gamma^2 - 2\gamma + 1\) and \(12 \gamma^2 - 6 \gamma + 1\) are positive for all \(\gamma\).
Thus we only need to show \(\alpha (6\gamma-1)+1 \geq 0\), but this is clear from writing it as \(1 - \alpha + 6 \gamma \alpha \geq 0\) and using that \(\alpha \leq 1\).
\end{proof}
\chapter{Block Lewis weights for sparsification and minimizing sums of Euclidean norms\label{chapter:sparsifyingnorms}}

In this chapter, we introduce block Lewis weights in the context of sparsification (recall from the introduction that these weights give us an ellipsoidal approximation to a particular family of objectives that we are interested in sparsifying). The main goal is to define block Lewis weights and show how they enable near-optimal sparsification of matrix block norm objectives. We then apply this to get an improved algorithm for minimizing sums of Euclidean norms in moderate accuracy settings. The material in this chapter is based on a joint work with Max Ovsiankin \cite{mo23}.

\section{Introduction}

Suppose we are given a large dataset that is computationally inconvenient to work with in a downstream task. To alleviate this, we can try to randomly sample a small representative subset of the original dataset. The design and analysis of randomized sampling algorithms for this purpose is well-explored (for example, see \cite{ss08,mmwy21,wy22_oneshot,wy23_sens} for preserving $\ell_p$ objectives, \cite{fl11} for preserving objectives for $k$-median, projective clustering, subspace approximation, and more, \cite{ss08,kkty21a,jls22,lee22} for preserving graph and hypergraph $\ell_2$-energy, and \cite{jlls23} for sums (of powers) of general norms).

In order to design randomized sampling algorithms, we first need to understand the properties of the original dataset we want to preserve. To this end, we study the problem of preserving \textit{block $p$-norm objectives}. Let $\cG = (\mA \in \R^{n\times d}, S_1,\dots,S_m, p_1,\dots,p_m)$ be a dataset consisting of a matrix $\mA \in \R^{n\times d}$. Consider a partitioning of $[n]$ into groups $S_1,\dots,S_m$ and consider positive numbers $p_1,\dots,p_m$. Let $\mA$ have rows $\va_1,\dots,\va_n$ and denote by $\mA_{S_i}$ the matrix in $\R^{\abs{S_i} \times d}$ whose rows are the rows of $\mA$ indexed by $S_i$. Consider the function $\gnorm{\mA\vx}{p}$ on some input vector $\vx \in \R^d$:
\begin{align}
    \gnorm{\mA\vx}{p}^p \coloneqq \sum_{i=1}^m \norm{\mA_{S_i}\vx}_{p_i}^p \label{eq:energy}
\end{align}
We use the norm notation because we can easily verify that for $p \ge 1$ and $p_i \ge 1$ for all $i$, $\gnorm{\cdot}{p}$ is a norm. We remark that objectives of the form of \eqref{eq:energy} are widely studied in geometric functional analysis, theoretical computer science, and data science. In \Cref{sec:results}, we go over one important application of the objective \eqref{eq:energy}. We defer a broader discussion of more applications and connections to \Cref{sec:related_works}.

Our goal in this chapter is to design and analyze randomized sampling algorithms to output a weighted subset that preserves \eqref{eq:energy} for all $\vx \in \R^d$. We give a formal problem statement for the general problem we study in \Cref{prob:sparsification_main_problem}. 

\begin{problem}[$\ell_p$ block norm sampling]
\label{prob:sparsification_main_problem}
We are given as input $\cG = \inparen{\mA \in \R^{n\times d}, S_1,\dots,S_m, p_1,\dots,p_m}$, $p > 0$, and an error parameter $\eps$. For all $i \in [m]$, we must output a probability distribution $\rho_1,\dots,\rho_m$ over $[m]$ such that if we choose a collection of groups $\cM = (i_1,\dots,i_{\mtilde})$ where each $i_h$ is independently distributed according to $\rho_i$, then the following holds with probability $\ge 1-\delta$:
\begin{align}
    \text{for all } \vx \in \R^d:\quad\quad\inparen{1-\eps}\gnorm{\mA\vx}{p}^p \le \frac{1}{\mtilde}\sum_{i\in \cM} \frac{1}{\rho_i}\cdot\norm{\mA_{S_i}\vx}_{p_i}^p \le \inparen{1+\eps}\gnorm{\mA\vx}{p}^p\label{eq:sparsification_main_objective}
\end{align}
We would like $\mtilde$ to be small with probability $1-\delta$ (for example, $\mtilde$ should not depend on $m$ and the dependence on $\delta^{-1}$ should be polylogarithmic).
\end{problem}
Observe that the formulation of \Cref{prob:sparsification_main_problem} is an instantiation of an \textit{importance sampling} framework. Specifically, we can think of the distribution $\cD$ as consisting of importance scores for each group. We form our sparse approximation by sampling group $i$ with probability $\rho_i$ and reweighting appropriately so that the function we return is an unbiased estimator of $\gnorm{\mA\vx}{p}$. We call $\mtilde$ the \textit{sparsity} of the procedure described in \Cref{prob:sparsification_main_problem}. Additionally, in the statement of our results, we will assume that $p$ is a constant (and thus any function solely of $p$ will treated as a constant in any $O\inparen{\cdot}$ or $\Omega\inparen{\cdot}$ terms). 

In this chapter, we give new results for Problem \ref{prob:sparsification_main_problem} and show how these imply faster algorithms for commonly implemented optimization problems.

\subsection{Our results}
\label{sec:results}

For a quick summary of our existence results for the block norm sampling problem (\Cref{prob:sparsification_main_problem}), see \Cref{table:comparison}. 

We begin with stating our main result\footnote{In the statement of \Cref{thm:one_shot_lewis}, writing the lower bound $p \ge 1/\log d$ instead of $p > 0$ is somewhat arbitrary -- we choose this lower bound to make our calculations easier later on.}, \Cref{thm:one_shot_lewis}.

\begin{restatable}[Block Lewis weight sampling]{mainthm}{oneshotlewis}
\label{thm:one_shot_lewis}
Let $\cG = (\mA \in \R^{n \times d}, S_1,\dots,S_m, p_1,\dots,p_m)$ where $S_1,\dots,S_m$ form a partition of $[k]$. Suppose at least one of the following holds:
\begin{itemize}
    \item $1 \le p < \infty$ and $p_1,\dots,p_m \ge 2$;
    \item $1/\log d \le p_1 = \dots = p_m = p < \infty$;
    \item $p_1 = \dots = p_m = 2$ and $1/\log d \le p < \infty$.
\end{itemize}
Let $P \coloneqq \max\inparen{1, \max_{i \in [m]} \min(p_i,\log\abs{S_i})}$. Then, there exists a probability distribution $\cD = \inparen{\rho_1,\dots,\rho_m}$ such that if
\begin{align*}
    \mtilde &= \Omega\inparen{\logv{\nfrac{1}{\delta}}\vbrho P \cdot  d^{\max(1,p/2)}},
\end{align*}
and if we sample $\cM \sim \cD^{\mtilde}$, then, with probability $\ge 1-\delta$,
\begin{align*}
    \text{for all } \vx \in \R^d,\quad (1-\eps)\gnorm{\mA\vx}{p}^p \le \frac{1}{\mtilde}\sum_{i \in \cM} \frac{1}{\rho_i} \cdot \norm{\mA_{S_i}\vx}_{p_i}^p \le (1+\eps)\gnorm{\mA\vx}{p}^p.
\end{align*}
\end{restatable}

We prove \Cref{thm:one_shot_lewis} in \Cref{sec:lewis}. It will follow from \Cref{thm:concentration} (stated and proven in \Cref{sec:generic_chaining}), which is a more general but more technical statement that also includes a description of the relevant distributions $\cD$.

We remark that when $p \ge 1$ and $p_1,\dots,p_m \ge 2$, the sampling probabilities \(\vrho\) mentioned in \Cref{thm:one_shot_lewis} can be found using the optimality conditions of a particular optimization problem that was stated and analyzed by \citet[Section 4]{jlls23}. That problem itself can be viewed as the natural generalization of the determinant maximization problem that yields the existence of Lewis's measure (see \cite[Section 2]{sz01} for details). However, \cite{jlls23} did not address the question of whether sparsification guarantees could be obtained with these weights beyond the case where the ``outer norm'' satisfies $p=2$.

Additionally, although \cite{jlls23} study sparsification of sums of norms and sums of powers $p > 1$ of uniformly smooth norms, we obtain an improved sparsity in the case entailed by \Cref{prob:sparsification_main_problem} (by a factor of $\psi_d\logv{\nfrac{d}{\eps}}^{\min(p-1,2)}$, where $\psi_d$ is the KLS ``constant'' in $d$ dimensions). We defer a more detailed comparison of our existence results to \Cref{sec:related_works}.

Furthermore, it is well-known that the polynomial terms in the sparsities in \Cref{thm:one_shot_lewis} are optimal. In particular, \citet[Corollary 1.6 and Theorem 1.7]{lww19} show that $\Omega(d^{\max(1,p/2)}+\eps^{-2}\mathrm{polylog}(\eps^{-1})d)$ rows must be chosen in order to satisfy the requirement imposed by (\ref{eq:sparsification_main_objective}).

Finally, the setting where $p = p_1 = \dots = p_m$ is a particularly important case of \Cref{prob:sparsification_main_problem}. Here, we see that $\norm{\mA\vx}_p = \gnorm{\mA\vx}{p}$, and so \Cref{prob:sparsification_main_problem} amounts to finding an $\ell_p$ subspace embedding under a group constraint (that certain rows must be kept together in the subsample). This might be a useful notion in practice, where the $S_i$ denote related observations that should be kept together for some downstream application. Moreover, this can be viewed as a higher-rank analog of $\ell_p$ row sampling, somewhat similarly to how the matrix Chernoff bound gives a higher-rank analog for the concentration of sums of bounded random matrices when compared to the rank-$1$ variant of \citet{rudelson1996}.

\paragraph{Computing sampling probabilities.}
The previous results show the existence of sampling probabilities \(\rho_1, \ldots, \rho_m\) such that sampling using those probabilities gives a sparsifier in the setting of \Cref{prob:sparsification_main_problem}.
To get a sparsification \textit{algorithm}, we need to also compute (or approximate) the sampling probabilities.

We give efficient algorithms to do so in natural cases.

\begin{restatable}[Computation of block Lewis weights]{mainthm}{computeblw}
\label{thm:computeblw}
Consider the setting of \Cref{thm:one_shot_lewis} and suppose at least one of the following holds:
\begin{itemize}
    \item $p = 2$ and $p_1,\dots,p_m \ge 2$;
    \item $1/\log n \le p_1=\dots=p_m=p < \infty$;
    \item $p_1 = \dots = p_m = 2$ and $1/\log n \le p < \infty$.
\end{itemize}
Let $P = \max\inparen{1, \max_{i \in [m]} \min(p_i,\log\abs{S_i})}$ and set
\begin{align*}
    \mtilde &= O\inparen{\logv{\nfrac{1}{\delta}}\vbrho P \cdot d^{\max(1,p/2)}}.
\end{align*}
Then, there is an algorithm that outputs a probability distribution $\cD = (\rho_1,\dots,\rho_m)$ such that sampling a multiset $\cM \sim \cD^{\mtilde}$ satisfies, with probability \(1 - \delta\),
\begin{align*}
    \text{for all } \vx \in \R^d,\quad (1-\eps)\gnorm{\mA\vx}{p}^p \le \frac{1}{\mtilde}\sum_{i \in \cM} \frac{1}{\rho_i} \cdot \norm{\mA_{S_i}\vx}_{p_i}^p \le (1+\eps)\gnorm{\mA\vx}{p}^p,
\end{align*}
Further, the algorithm to find $\cD$ performs at most \(\text{polylog}(k,n,m)\) leverage score overestimate computations or linear system solves.
\end{restatable}

We prove \Cref{thm:computeblw} in \Cref{sec:alg}.

We formally define a \textit{leverage score overestimate computation} in \Cref{def:lev_score_over}. Alternately, these can be implemented using linear system solvers that solve systems of the form $\mA^{\top}\mD\mA\vy=\vz$ for diagonal $\mD$ (see \cite{ls19} for details). Although the runtime of this primitive depends on the structure of the input, each such iteration runs in $\widetilde{O}(\mathsf{nnz}(\mA) + d^{\omega})$ time. Moreover, in the special case where the matrix $\mA$ is a graph edge-incidence matrix, the runtime improves to $\widetilde{O}(\mathsf{nnz}(\mA))$.

Finally, we note that our algorithms are faster than the log-concave sampling-based routines given in \cite{jlls23} for calculating sparse approximations to sums (of powers of) more general norms, when the outer norm $p$ satisfies $1 \le p \le 2$ (they do not give algorithms for the case where $p > 2$). In particular, while their algorithm applies to a more general setting, the runtime is $\widetilde{O}(m + d^5)$. In contrast, since our algorithms only depend on a polylogarithmic number of leverage score overestimate computations or linear system solves, we can obtain much faster runtimes (in particular improved powers of $n$). This means that we can apply our algorithms to downstream optimization tasks where the main computational primitive is a linear system solver (as is the case for many general frameworks for convex programming).

\paragraph{Applications to minimizing sums of Euclidean norms.} A well-studied regression task is the \textit{minimizing sums of Euclidean norms} (MSN) problem. We are given \(\mA \in \R^{n \times d}\) and \(\vb \in \R^n\), and a partition \(S_1, \ldots, S_m\) of \([n]\).
In this problem, we would like to find
\begin{align}
    \quad\quad\min_{\vx\in\R^d} \sum_{i=1}^m \norm{\mA_{S_i} \vx-\vb_{S_i}}_2. \label{eq:msn}
\end{align}
Solving the MSN objective \eqref{eq:msn} subsumes several widely implemented optimization problems such as variants of Euclidean single facility location, Euclidean multifacility location, Euclidean Steiner minimum tree under a given topology, and plastic collapse analysis. See the long line of work on this problem \cite{and96,xy97,acco00,qsz02} for a more detailed discussion. Additionally, observe that if all $\abs{S_i}=1$, then \eqref{eq:msn} is nothing but $\ell_1$ regression (i.e., $\min_{\vx\in\R^d} \norm{\mA\vx-\vb}_1$). Thus, \eqref{eq:msn} is a generalization of $\ell_1$ regression. Finally, notice that \eqref{eq:msn} subsumes the \textit{stochastic robust approximation problem} when the norm in question is the Euclidean norm and the design $\mA$ assumes a finite number of values -- see \cite[Section 6.4.1]{boyd2004convex} for further discussion.

In this chapter, we will be interested in algorithms that return a $(1+\eps)$-multiplicative approximation to the objective -- namely, we desire a point $\widehat{\vx} \in \R^d$ such that
\begin{align*}
    \sum_{i=1}^m \norm{\mA_{S_i}\widehat{\vx}-\vb_{S_i}}_2 \le (1+\eps)\min_{\vx\in\R^d} \sum_{i=1}^m \norm{\mA_{S_i} \vx-\vb_{S_i}}_2.
\end{align*}
To our knowledge, the best known algorithms based on interior point methods output a $(1+\eps)$-approximate solution to \eqref{eq:msn} $\widetilde{O}(\sqrt{m}\logv{\nfrac{1}{\eps}})$ calls to a linear system solver \cite{and96,xy97} for matrices of the form $\mA^{\top}\mD\mA$ for block-diagonal matrices $\mD$.

By applying \Cref{thm:computeblw} on the matrices $\insquare{\mA_{S_i} \vert \vb_{S_i}} \in \R^{\abs{S_i} \times (d+1)}$ with $p_1=\dots=p_m=2$ and $p=1$, observe that within $\widetilde{O}(1)$ linear system solves in matrices $\mA^{\top}\mD\mA$ for nonnegative diagonal $\mD$, we obtain an objective with $\widetilde{O}(\eps^{-2} \cdot d)$ terms that approximates \eqref{eq:msn} up to a $(1\pm\eps)$ multiplicative factor on all vectors $\vx \in \R^{d+1}$ whose last coordinate is $1$. This immediately implies \Cref{thm:msnalg}.

\begin{restatable}[Minimizing sums of Euclidean norms]{mainthm}{msnalg}
\label{thm:msnalg}
Let $\mA \in \R^{n \times d}$ and \(\vb \in \R^n\), and $S_1, \dots, S_m$ be a partition of $k$. There exists an algorithm that, with probability $\ge 1-\delta$, returns $\widehat{\vx}$ such that
\begin{align*}
    \sum_{i=1}^m \norm{\mA_{S_i}\widehat{\vx}-\vb_{S_i}}_2 \le (1+\eps)\min_{\vx\in\R^d} \sum_{i=1}^m \norm{\mA_{S_i} \vx-\vb_{S_i}}_2.
\end{align*}
The algorithm runs in $\widetilde{O}\inparen{\nfrac{\sqrt{d}}{\eps} \cdot \sqrt{\logv{\nfrac{1}{\delta}}}}$ calls to a linear system solver in matrices of the form $\mA^{\top}\mD\mA$ for block-diagonal matrices $\mD$, where each block has size $(\abs{S_i}+1) \times (\abs{S_i}+1)$.
\end{restatable}

We prove \Cref{thm:msnalg} in \Cref{sec:applications_msn}.

\Cref{thm:msnalg} improves over the best-known iteration complexities for solving \eqref{eq:msn} when the number of summands is much larger than the input dimension, i.e., $m \gg d$. Furthermore, the iteration complexity stated in \Cref{thm:msnalg} matches the iteration complexity for $\ell_1$ regression up to the $\eps^{-1}$ term \cite{bllssw21}. It is an interesting (but probably challenging) open problem to design and analyze an algorithm for \eqref{eq:msn} with iteration complexity $\widetilde{O}(\sqrt{d}\logv{\nfrac{1}{\eps}})$, which would exactly match what is known for $\ell_1$ regression.

Finally, we note that in the special case of the geometric median, where all the $\mA_{S_i} = \mI_r$ for some fixed dimension \(r\), an algorithm with runtime $\widetilde{O}(\mathsf{nnz}(\vb)\logv{\nfrac{1}{\eps}}^3)$ is known due to \citet{clmps16}. The algorithm is a long-step interior point method with a custom analysis and follows from different techniques from ours.

\paragraph{Outline.} The rest of this chapter is organized as follows. In the remainder of this section, we establish notation that we use throughout the rest of the chapter (\Cref{sec:notation}), give an overview of our technical methods (\Cref{sec:sparsification_overview}), and discuss some prior and related works (\Cref{sec:related_works}). In \Cref{sec:background}, we give background from linear algebra, convex geometry, and probability that we rely on for the rest of the chapter. In \Cref{sec:coverings}, we prove bounds on geometric quantities known as \textit{covering numbers}. These play a crucial role in our concentration arguments. In \Cref{sec:generic_chaining}, we prove that our general sampling scheme concentrates and therefore preserves the original objective on all $\vx \in \R^d$, with high probability. In \Cref{sec:applications}, we show how to apply our general sampling scheme to the problems we discuss in \Cref{sec:results}. Finally, in \Cref{sec:alg}, we describe our algorithmic results. 

\subsection{Notation and definitions}
\label{sec:notation}

\paragraph{General notation.} For positive integer $N$, we let $[N]$ denote the set $\inbraces{i \in \Z \suchthat 1 \le i \le N}$. All $\log$s are base $2$; we use $\ln$ to denote the natural logarithm. We let $\ve_1,\dots,\ve_d$ denote the standard basis vectors in $\R^d$. When we write $a \lesssim b$, we mean that $a \le Cb$ for some universal constant $C > 0$.

\paragraph{Linear algebra notation.} In this chapter, we work extensively with matrices and vectors. We always denote matrices with capital letters in boldface (e.g. $\mA$) and vectors with lowercase letters in boldface (e.g. $\vx$). With a few exceptions, we write the rows of a matrix using the lowercase boldface version of the same letter used to write the matrix along with a subscript denoting which index the row corresponds to. For example, $\va_i$ denotes the $i$th row of matrix $\mA$. In a slight abuse of notation, for a symmetric matrix $\mM$, we let $\mM^{-1} \coloneqq \sum_{i=1}^{\rank{\mM}} \lambda_i^{-1}\vu_i\vu_i^{\top}$, where $\vu_i$ is the $i$th eigenvector of $\mM$. In other words, we write $\mM^{-1}$ to denote the pseudoinverse of $\mM$ when $\mM$ is symmetric. We will never use the inverse notation $\mM^{-1}$ for a non-symmetric matrix $\mM$.

\subsection{Technical overview}
\label{sec:sparsification_overview}

In this subsection, we give a bird's eye view of the technical methods behind our proof of \Cref{thm:one_shot_lewis}.

\subsubsection{Concentration}
\label{sec:sparsification_overview_concentration}

We begin with an explanation of our concentration proof. This type of argument has become standard in the line of work on sparsification (particularly in \cite{lee22,jlls23}), but we include a description for completeness.

Let $B_p \coloneqq \inbraces{\vx\in\R^d\suchthat\gnorm{\mA\vx}{p}\le1}$. By a standard symmetrization reduction, it suffices to fix $i_1,\dots,i_{\mtilde}$ (not necessarily distinct) and argue that for independent $R_1,\dots,R_{\mtilde}$ where $R_h \sim \mathsf{Unif}\inparen{\pm 1}$ we have
\begin{align}
    \exvv{R_h}{\sup_{\vx \in B_{p}} \abs{\sum_{h=1}^{\mtilde} R_h \cdot \frac{\norm{\mA_{S_{i_h}}\vx}_{p_{i_h}}^p}{\rho_{i_h}}}} \le \mtilde\cdot\eps.\label{eq:intro_empirical_process}
\end{align}
Intuitively, satisfying \eqref{eq:intro_empirical_process} means that for the rebalancing of the groups given by the $\rho_{i_h}$, a Rademacher average of the groups evaluated on every point in $\vx$ is close to $0$. It is straightforward to check that the above instantiation is a subgaussian process under an appropriately chosen distance function $\myfunc{\dtwo}{\R^{d} \times \R^{d}}{\R}$. Thus, we will apply chaining \cite{tal21}, which can be thought of as simultaneously controlling \eqref{eq:intro_empirical_process} on $\eps$-nets of $B_p$ using the metric $\dtwo$, for all $\eps>0$.

To apply chaining, the main technical task is to understand the \textit{entropy numbers} $e_N(B_{p},\dtwo)$. The entropy numbers $e_N(B_p, \dtwo)$ are the values $\eta$ that answer the question, ``what is the smallest $\eta$ such that $B_p$ can be covered by at most $2^{2^N}$ balls of $\dtwo$-radius $\eta$?'' (or see \Cref{defn:entropy_numbers}).

\subsubsection{Covering numbers}
\label{sec:sparsification_overview_geometry}

In this subsection, we explain how to control the entropy numbers as required by \Cref{sec:sparsification_overview_concentration}. We first define the sampling body (\Cref{defn:sampling_polyhedron}).

\begin{definition}[Sampling body]
\label{defn:sampling_polyhedron}
Let $S$ be some subset of $[m]$ and $\rho_1,\dots,\rho_m$ be a probability distribution. Define the norm $\polynorm{\vx} \coloneqq \max_{i \in S} \rho_i^{-1/p}\norm{\mA_{S_i}\vx}_{p_i}$. We call the unit norm ball of $\polynorm{\vx}$ the \textit{sampling body}.
\end{definition}

Recall that the covering number $\cN(K_1,K_2)$ for two symmetric convex bodies $K_1$ and $K_2$ is the minimum number of translates of $K_2$ required to cover $K_1$. Additionally, recall from the previous subsection the notion of \textit{entropy numbers} (which we will define in \Cref{defn:entropy_numbers}). We will reduce controlling \eqref{eq:intro_empirical_process} to bounding the entropy numbers
\begin{align*}
    e_N\inparen{\inbraces{\vx\in\R^d\suchthat\gnorm{\mA\vx}{p}\le 1}, \inbraces{\vx\in\R^d\suchthat\polynorm{\vx}\le 1}}
\end{align*}
when $N$ is small. This places us in the setting where a simple volume-based argument becomes suboptimal. In this range, the dual Sudakov inequality (\Cref{fact:sudakov}) is the technical workhorse that allows us to get sharper bounds than what we would get if we applied just a volume-based bound. It states that if $B$ is the Euclidean ball in $d$ dimensions and $K$ is some symmetric convex body in $d$ dimensions, then we have
\begin{align*}
    \log\cN(B,\eta K) \lesssim \eta^{-2} \exvv{\vg\sim\cN(0,\mI_d)}{\norm{\vg}_K}^2,
\end{align*}
where $\norm{\cdot}_K$ is the \textit{gauge norm} for $K$, defined by $\norm{\vx}_K \coloneqq \inf\inbraces{t > 0 \suchthat \vx/t \in K}$.

However, applying the dual Sudakov inequality requires that we analyze covering numbers of the form $\log\cN(B, K)$ where $B$ is the Euclidean ball in $d$ dimensions and $K$ is some symmetric convex body. Denoting $\{\vx\in\R^d\suchthat\polynorm{\vx}\le 1\}$ by $K$, we see that we cannot immediately apply the dual Sudakov inequality to bound $\log\cN(B_p,\eta K)$. This is because $B_p$ is not a (linear transformation of a) Euclidean ball. The work of \cite{jlls23} resolve this by generalizing the dual Sudakov inequality to cover arbitrary symmetric convex bodies. Unfortunately, this approach is not optimal in every setting. One source of the loss arises from exploiting the concentration of Lipschitz functionals of isotropic log-concave random vectors -- improving the bounds on this concentration depends on further progress on the KLS conjecture. Another is that the one-dimensional conditionals of isotropic log-concave random variables, without any further assumptions, are only subexponential.

To escape these inefficiencies, we will want to try to find a way to apply the dual Sudakov inequality as-is. We may then exploit the concentration of Lipschitz functionals of Gaussian random vectors, which we do have a tight understanding of (for a precise statement, see \Cref{fact:lip_concentration}). A natural attempt is to first observe that for any $t > 0$,
\begin{align}
    \log\cN\inparen{B_p, \eta K} \le \log\cN\inparen{B_p, t\widehat{B_2}} + \log\cN\inparen{t\widehat{B_2}, \eta K} = \log\cN\inparen{B_p, t\widehat{B_2}} + \log\cN\inparen{\widehat{B_2}, \frac{\eta}{t}\cdot K}.\label{eq:intro_split_covering}
\end{align}
We will choose $\widehat{B_2}$ to be a linear transformation of a Euclidean ball so that we can control $\log\cN\inparen{\widehat{B_2}, \eta/t \cdot K}$ using the dual Sudakov inequality. 

Here, we will split our argument based on whether $p \ge 2$. When $p \ge 2$, it will become clear later on that it will be sufficient to choose $\widehat{B_2}$ so that $B_p \subseteq \widehat{B_2}$. Then, it is easy to see that when $t=1$, we get $\log\cN(B_p,t\widehat{B_2}) = 0$. Hence, we have $\log\cN(B_p, \eta K) \le \log\cN(\widehat{B_2}, \eta K)$, and the required bound will follow from exploiting the concentration of Lipschitz functionals of Gaussian random vectors and then applying the dual Sudakov inequality.

However, when $p < 2$, we are still left with a pesky $\log\cN(B_p,t \widehat{B_2})$ term. Loosely, this is almost dual to the statement of the dual Sudakov inequality. Now, because it is known that covering number duality does hold when one of the bodies in question is the Euclidean ball, it may be tempting to simply write $\log\cN(B_p, t\widehat{B_2}) = \log\cN(\widehat{B_2},t B_q)$ where $B_q$ is the dual ball to $B_p$ after applying some linear transformation to map $\widehat{B_2}$ to $B$. The challenge here is that we do not believe that the gauge of the resulting $B_q$ has a form that is amenable to analysis. We will therefore need to be more careful, and we describe our alternative approach in \Cref{sec:sparsification_overview_measure}.

\subsubsection{The change-of-measure principle and norm interpolation}
\label{sec:sparsification_overview_measure}

Recall from the previous part that our goal is to bound $\log\cN(B_p,t\widehat{B_2})$ when $p < 2$.

We are now ready to introduce our main conceptual message -- \textit{by changing the measure under which we take norms, we can almost automatically identify a linear transformation of a Euclidean ball $\widehat{B_2}$ that is a good approximation to $B_p$}. This sort of idea has already been used by \citet{blm89} and \citet{sz01} to obtain the required $\widehat{B_2}$ in the special case where all the $S_i$ are singletons. We will generalize this machinery to give similar results for the block norm sampling problem.

Let us describe this idea further. Let $\vlambda = [\lambda_1,\dots,\lambda_m]^{\top}$ denote a probability measure over the groups. Let $\mLambda \in \R^{n\times n}$ be the diagonal matrix such that if $j \in S_i$, then $\mLambda_{jj} = \lambda_i$. Finally, for any $r > 0$ and $\vy\in\R^n$, let $\gnorml{\vy}{r} = (\sum_{i \le m} \lambda_i\norm{\vy_{S_i}}_{p_i}^{r})^{1/r}$ and $B_r \coloneqq \{\vx\in\R^d \suchthat \gnorml{\mLambda^{-1/p}\mA\vx}{r} \le 1\}$. Notice that under this definition, we  still have \(B_p\) as before. We will first describe the argument when the groups are singletons, then explain how to move onto the general case. We will take $\widehat{B_2} = B_2$; it is easy to see that this is a linear transformation of a Euclidean ball.

Next, notice that by log-convexity of norms, if we choose $0 < \theta < p$ and $r > 2$ for which $1/2 = (\theta/2)/p + (1-\theta/2)/r$, we have
\begin{align}
    \gnorml{\vy}{2}^2 \le \gnorml{\vy}{p}^{\theta} \cdot \gnorml{\vy}{r}^{2-\theta}.\label{intro:compact_rounding_interpolation}
\end{align}
We will exploit this observation as follows. For all integers $h \ge 0$, we will show that there exists a set $\cL_h$ that is (a subset of) the unit ball of $B_2$ such that every pair of points in $\cL_h$ is $\delta_h$-separated according to $\gnorml{\cdot}{r}$. We will find $\delta_h$ according to the interpolation inequality \eqref{intro:compact_rounding_interpolation}. Furthermore, we will generate $\cL_h$ using a sort of compactness argument arising from a $B_2$-maximally separated subset of $B_p$. This means we get, for every $h \ge 0$,
\begin{align*}
    \log\cN(\widehat{B_2}, \delta_h B_r) \ge \log\abs{\cL_h} \ge \logv{\frac{\cN(B_p,8^h t \widehat{B_2})}{\cN(B_p, 8^{h+1} t \widehat{B_2})}}.
\end{align*}
Then, summing over $h \ge 0$ (noting that once $h$ is sufficiently large, $\cN(B_p, 8^ht \widehat{B_2}) = 1$), we have
\begin{align*}
    \log\cN(B_p, t\widehat{B_2}) \le \sum_{h \ge 0} \log\cN(\widehat{B_2},\delta_h B_r).
\end{align*}
Notice that the right hand side can be evaluated using the dual Sudakov inequality\footnote{For technical reasons that will be clearer in \Cref{sec:coverings}, we will have to do this after another interpolation step.} (recall the previous section), so it suffices to show that $\log\cN(\widehat{B_2},\eta B_r)$ is small.

This is where the choice of measure becomes crucial. Since both $\widehat{B_2}$ and $B_r$ are dependent on our choice of measure $\vlambda$, we will need to carefully choose the measure so that our covering numbers are well-behaved. A classical result of \citet{Lewis1978} establishes the existence of a change-of-measure under which we simultaneously get:
\begin{equation}\label{intro:rounding}
  \begin{aligned}
    d^{1/2-1/r}B_r &\subset \widehat{B_2} \subset B_r & \text{for all } r < 2 \\
    B_r &\subset \widehat{B_2} \subset d^{1/2-1/r}B_r & \text{for all } r > 2
  \end{aligned}
\end{equation}
This change-of-measure corresponds to the ``$\ell_p$ Lewis weights'' of $\mA$ (in particular, if $w_i$ is the $i$th $\ell_p$ Lewis weight, then we set $\lambda_i = w_i/n$). It will turn out that this choice of $\vlambda$ is enough for us to ensure that $\cN(\widehat{B_2}, \eta B_r)$ is sufficiently small for our purposes, which eventually follows from \eqref{intro:rounding}.

\paragraph{Handling general $S_i$.} The main challenge with directly porting this argument to the block norm sampling problem is that $B_2$ is not a linear transformation of a Euclidean ball unless $p_1=\dots=p_m=2$. We will therefore have to choose $\widehat{B_2}$ to be a ``rounding'' of $B_2$ such that $B_2 \subseteq \widehat{B_2}$. Observe that the interpolation step \eqref{intro:compact_rounding_interpolation} will continue to hold here, as we will get $\norm{\vy}_{\widehat{B_2}} \le \gnorml{\vy}{2}$. However, if $\widehat{B_2}$ is chosen suboptimally, then there could a large loss in the interpolation step \eqref{intro:compact_rounding_interpolation}.

To understand what we need from our measure and rounding, let us try to derive a version of \eqref{intro:rounding} for general $S_i$. We show an example of this calculation for $r = p \le 2$; the other cases follow similarly. Let $\vlambda \in \R^m_{\ge 0}$ denote a probability measure. Let $\mW$ be a diagonal ``rounding matrix'' so that for all $\vx \in \R^d$, we have
\begin{align*}
    \norm{\mW^{1/2}\mLambda^{1/2-1/p}\mA\vx}_2 \le \gnorm{\mLambda^{1/2-1/p}\mA\vx}{2} = \gnorml{\mLambda^{-1/p}\mA\vx}{2}.
\end{align*}
Letting $\widehat{B_2} = \inbraces{\vx\in\R^d\suchthat \norm{\mW^{1/2}\mLambda^{1/2-1/p}\mA\vx}_2 \le 1}$, the above inequality gives $B_2 \subseteq \widehat{B_2}$, as desired. Next, observe that since $\vlambda$ is a probability measure, we get $B_2 \subseteq B_p$ for free. For the other direction, we write
\begin{align*}
    \gnorml{\mLambda^{-1/p}\mA\vx}{2}^2 &= \gnorm{\mLambda^{1/2-1/p}\mA\vx}{2}^2 = \sum_{i=1}^m \lambda_i\norm{\lambda_i^{-1/p}\mA\vx}_{p_i}^2 = \sum_{i=1}^m \lambda_i\norm{\lambda_i^{-1/p}\mA\vx}_{p_i}^{p}\norm{\lambda_i^{-1/p}\mA\vx}_{p_i}^{2-p} \\
    &\le \sum_{i=1}^m \lambda_i\norm{\lambda_i^{-1/p}\mA\vx}_{p_i}^{p} \cdot \max_{i \in [m]} \norm{\lambda_i^{-1/p}\mA\vx}_{p_i}^{2-p} = \gnorm{\mA\vx}{p}^p \cdot \max_{i \in [m]} \norm{\lambda_i^{-1/p}\mA\vx}_{p_i}^{2-p} \\
    &\le \gnorm{\mA\vx}{p}^p \cdot \max_{i \in [m]} \inparen{\max_{\vx \in \R^d} \frac{\norm{\lambda_i^{-1/p}\mA\vx}_{p_i}^{2}}{\gnorm{\mLambda^{1/2-1/p}\mA\vx}{2}^2}}^{1-p/2} \cdot \gnorm{\mLambda^{1/2-1/p}\mA\vx}{2}^{2-p}.
\end{align*}
We combine the $\gnorm{\mLambda^{1/2-1/p}\mA\vx}{2}$ terms and take the $p$th root of both sides, giving
\begin{align*}
    \gnorm{\mLambda^{1/2-1/p}\mA\vx}{2} &\le \gnorm{\mA\vx}{p} \cdot \max_{i \in [m]} \inparen{\max_{\vx \in \R^d} \frac{\norm{\lambda_i^{-1/p}\mA\vx}_{p_i}^{2}}{\gnorm{\mLambda^{1/2-1/p}\mA\vx}{2}^2}}^{1/p-1/2} \\
    &= \gnorm{\mA\vx}{p} \cdot \max_{i \in [m]} \inparen{\frac{1}{\lambda_i} \cdot \underbrace{\max_{\vx \in \R^d} \frac{\norm{\lambda_i^{1/2-1/p}\mA\vx}_{p_i}^{2}}{\gnorm{\mLambda^{1/2-1/p}\mA\vx}{2}^2}}_{\widehat{\tau}_i(\mLambda^{1/2-1/p}\mA)}}^{1/p-1/2}.
\end{align*}
We may think of the quantity $\widehat{\tau}_i$ as a generalized \textit{leverage score}. Specifically, it upper bounds the contribution of the term $\|\lambda_i^{1/2-1/p}\mA\vx\|_{p_i}^{2}$ to the objective $\gnorm{\mLambda^{1/2-1/p}\mA\vx}{2}^2$. The above calculation shows us that if we make $\widehat{\tau}_i/\lambda_i$ small for all $i$, then we can get a tight relationship between $B_2$ and $B_p$. A slight weakening of the definition of the $\widehat{\tau}_i$ motivates the notion of a \textit{block Lewis overestimate} that we use in the remainder of the chapter.

\begin{restatable}[Block Lewis overestimate]{definition}{blocklewisoverestimate}
\label{defn:block_lewis_overestimate}
Let $\tau_j(\mM)$ denote the leverage score of the $j$th row of $\mM$. Let \(F^{\star} > 0\). For $p > 0$ and $p_i > 0$, we say the probability measure $\vlambda$ and rounding $\mW$ form an $\Fstar$-block Lewis overestimate if for all $i \in [m]$, we have
\begin{align*}
    \frac{1}{\lambda_i}\inparen{\sum_{j \in S_i} \inparen{\frac{\tau_j(\mW^{1/2}\mLambda^{1/2-1/p}\mA)}{w_j}}^{p_i/2}}^{2/p_i} \le \Fstar.
\end{align*}
\end{restatable}
Following the above argument, establishing a probability measure $\vlambda$ and a rounding matrix $\mW$ that form an $\Fstar$-block Lewis overestimate will imply
\begin{equation}\label{intro:blw_rounding}
  B_2 \subseteq \widehat{B_2}\quad \text{ and }\quad
  \begin{aligned}
    (\Fstar)^{1/2-1/r}B_r &\subset B_2 \subset B_r & \text{for all } r < 2 \\
    B_r &\subset B_2 \subset (\Fstar)^{1/2-1/r}B_r & \text{for all } r > 2
  \end{aligned}.
\end{equation}
With \eqref{intro:blw_rounding} in hand, we at least have enough reason to believe that establishing $(\vlambda, \mW)$ that form an $\Fstar$-block Lewis overestimate may yield the requisite control over $\log\cN(\widehat{B_2},\eta B_r)$. To actually get this, by the dual Sudakov inequality, we estimate $\norm{\vg}_{B_r}$ for $\vg \sim \cN(0,\mI_d)$. Doing so is a matter of applying again the fact that the concentration of Lipschitz functionals of Gaussian vectors is determined entirely by the Lipschitz parameter of the functional. 

To actually find $\vlambda$ and $\mW$ with a small value of $\Fstar$, we split into cases. When $p \ge 1$ and $p_1=\dots=p_m \ge 2$, we extract the relevant $\vlambda$ and $\mW$ from the analysis in the proof of \cite[Lemma 4.2]{jlls23}, which yields $\Fstar = d$. When $p_1=\dots=p_m=p\ge1/\log d$ or $p_1=\dots=p_m=2$ and $p\ge1/\log d$, we separately prove that we can find $\vlambda$ and $\mW$ satisfying \Cref{defn:block_lewis_overestimate}, again with $\Fstar = d$. Hence, in all cases, the control we get over $\log\cN(\widehat{B_2},\eta B_r)$ is essentially as good as what we get in the case where all the $S_i$ are singletons.

\paragraph{Change-of-measures in functional analysis.} We note that other change-of-measure arguments are used throughout the study of finite-dimensional subspaces of $L_p$, as they are a very useful way to compare the $L_p$ norm to some other norm of interest. See the survey by \citet[Section 1.2]{js01} for more information.

\subsection{Prior results, related works, and connections}
\label{sec:related_works}

\paragraph{Relevance of matrix block norms.} We discuss the importance of the matrix block norm objective \eqref{eq:energy} to functional analysis, theoretical computer science, and data science, beyond our previous discussion of the MSN problem \eqref{eq:msn}.

In the special case where all the $p_i$ are equal to one another (call this value $q$), the set of $\vx$ for which $\gnorm{\mA\vx}{p} < \infty$ yields a subspace of a mixed $p,q$ norm space (sometimes notated as $\ell_p(\ell_q)$). Observe that \Cref{thm:one_shot_lewis} implies a finite-dimensional subspace embedding result for finite-dimensional subspaces of infinite-dimensional $L_p(\ell_q^r)$ (where $r$ is finite).  Spaces of the form $\ell_p(\ell_q)$ are widely studied in the geometric functional analysis and approximation theory communities; see, e.g., \cite{ps12,kv15,mu19,jkb22} and the references therein. 

As mentioned in those works, a central motivation for studying $\ell_p(\ell_q)$ is that they are natural testbeds with which to evaluate and further our understanding of the geometry of symmetric convex bodies in high dimensions. Consequently, studying block norm subspace embedding problems (\Cref{prob:sparsification_main_problem}) is a fruitful direction through which to improve our geometric handle of subspaces of $\ell_p(\ell_q)$ and symmetric convex bodies in general. We note that understanding the correct $\mathrm{polylog}$ dependencies in $d$ for this problem typically requires new geometric insights. For instance, the necessity of additional polylogarithmic dependencies on the dimension $d$ is not even totally understood when $\abs{S_i} = 1$ and $p\neq 2$, and resolving them likely requires significant new geometric ideas \cite[Conjecture 2]{hrr22}. 

The matrix block norms are also ubiquitous in both theoretical computer science and data science. For example, the block norm objective has been studied in the context of hypergraph Laplacians. One recovers this by choosing $p=2$ and $p_1=\dots=p_m=\infty$; see the discussion in \cite[Section 1.2]{jls22} to see how to rewrite the hypergraph Laplacian in the form of \eqref{eq:energy}. Within data science, the block norms are used to encourage structured solutions to underdetermined linear systems (i.e., in a noiseless setting, we can set up and solve the convex optimization problem\footnote{This is similar to how basis pursuit can be seen as encouraging sparsity in a noiseless setting, while LASSO does so in the presence of noise \cite[Section 7.2]{Wainwright_2019}.}, ``find a vector $\vy$ in the affine space $\mB\vy = \vb$ minimizing $\gnorm{\vy}{p}^p$''); see, e.g., \cite{yy06,bach07,nie2010efficient,sfht13} and other applications mentioned in \cite{sra12}. As a concrete candidate application of our results to such settings, inspired by \cite{cd21,mmwy21}, we believe that our results can be used as subroutines to give runtime and query-efficient algorithms for active $\gnorm{\cdot}{p}$ regression when $p > 0$ and $p_1=\dots=p_m=2$ (generalizing the basis pursuit equivalent of the group Lasso objective) or when $p=2$ and $p_1,\dots,p_m \ge 2$. 

On a more conceptual level, we are optimistic that some of our results will be useful for designing faster algorithms for norm-constrained optimization problems. This is partly motivated by our discussion around the MSN problem \eqref{eq:msn} and is similar to how an improved geometric understanding of Lewis weights improved the iteration complexities for linear programming and $\ell_p$ regression \cite{ls19,jls21}.

\paragraph{Lewis weights for $\ell_p$ row sampling.} When each group $S_i$ has size $1$, notice that we have $\gnorm{\mA\vx}{p}^p = \sum_{i=1}^m \abs{\ip{\va_i,\vx}}^p = \norm{\mA\vx}_p^p$. Consequently, in this special case, satisfying (\ref{eq:sparsification_main_objective}) is exactly equivalent to computing an \textit{$\ell_p$ subspace embedding} for $\mA$. There is a long line of work studying computing $\ell_p$ subspace embeddings using Lewis weights, starting with that of \citet{blm89}. For the details of this argument, see \cite[Section 7]{blm89} and \cite{sz01}. 

\paragraph{Sparsifying sums of norms.} The work perhaps most closely related to ours is \cite{jlls23}. There, the authors give existence results for sparse approximations to sums (of powers) of norms. It is easy to see that this is a more general problem than the one we study. However, this generality comes at a cost. In particular, the sparsity given by our \Cref{thm:one_shot_lewis} improves over theirs by a factor of $\psi_d\logv{\nfrac{d}{\eps}}^{\min(p-1,2)}$, where $\psi_d$ denotes the KLS ``constant'' in $d$ dimensions (which is currently $\sqrt{\log d}$, due to \citet{kla23}). See their Theorem 1.3 for more details. And, as mentioned earlier, we believe understanding \Cref{prob:sparsification_main_problem} down to the correct polylogarithmic dependencies in $d$ is an important geometric question.

The authors also define the block Lewis weights as the natural generalization of the determinant-maximization program that \citet{sz01} use to prove the existence of Lewis's measure for all $p > 0$. They use this to give results for sparsification of sums of certain powers of arbitrary norms (obtaining a sparsity of $\sim d^{2-1/p}$ when $1 \le p \le 2$) and for \Cref{prob:sparsification_main_problem} when the outer norm $p=2$ (obtaining a sparsity of $\sim d$). We note that the result they obtain when $p=2$ provides logarithmic-factor improvements over ours, as the main technical primitive they use is a chaining estimate developed by \citet{lee22} that meaningfully exploits the fact that the space of events is a subset of a $2$-uniformly convexity set. However, they did not address whether the block Lewis weights could yield to sparsification guarantees for \Cref{prob:sparsification_main_problem}. Additionally, their construction of the block Lewis weights does not yield a change-of-measure that allows for sparsification when the inner norms $p_i \le 2$ or when the outer norm $p \le 1$. 

\paragraph{Summary for sparsifying sums of norms.} See \Cref{table:comparison} for a comparison between our new results and a selection of the most relevant prior work on sparsifying sums of norms. We focus on results concerning \(\ell_p\)-norms specifically, although \cite{jlls23} has results for more general classes of norms.
By ``sampling'', we mean the work provides an analysis that shows how sampling according to some sampling probabilities gives a sparsifier of size $\widetilde{O}(\eps^{-2}d^{\max(1,p/2)})$ with good probability.
By ``fast computation'', we mean the work provides an algorithm to compute sampling probabilities with \(\text{polylog}(k, m, d)\) leverage score computations or linear system solves (or some other primitive that can be implemented in time $\widetilde{O}(\mathsf{nnz}(\mA) + d^{\omega})$). For works that only explicitly handle \(\abs{S_i}=1\), we leave \(p_1, \ldots, p_m\) blank because the choice of inner norms does not affect
the objective.
\begin{table}[H]
\centering
\begin{tabular}{l c c c c r}
\hline
\textbf{Block size} & \(p\) & \(p_1, \ldots, p_m\) & \textbf{Sampling} & \textbf{Fast computation} & \\
\hline
\(1\) & \(1 \le p < \infty \) &  & \checkmark & & \cite{blm89} \\
\(1\) & \(0 < p < 1\) &  & \checkmark & & \cite{sz01} \\
\(1\) & \(0 < p < 4 \) & & & \checkmark & \cite{cp15} \\
\(\ge 1\) & \(p = 2\) & \(p_1 = \cdots = p_m = \infty\) & \checkmark & & \cite{lee22} \\
\(\ge 1\) & \(p=2\) & \(p_1=  \cdots = p_m = \infty\) & \checkmark & \checkmark & \cite{jls22} \\
\(\ge 1\) & \(1 \le p < \infty\) & \(p_1, \ldots, p_m \ge 2\) & \checkmark & & \cite{jlls23} \\
\rowcolor{blue!5} \(\ge 1\) & \(1 \le p < \infty\) & \(p_1, \ldots, p_m \ge 2\) & \checkmark & & This work \\
\rowcolor{blue!5} \(\ge 1\) & \(\nfrac{1}{\log d} \le p < \infty\) & \(p_1, \ldots, p_m = p\) & \checkmark & \checkmark & This work \\
\rowcolor{blue!5} \(\ge 1\) & \(\nfrac{1}{\log d} \le p < \infty\) & \(p_1, \ldots, p_m = 2\) & \checkmark & \checkmark & This work \\
\rowcolor{blue!5} \(\ge 1\) & \(p=2\) & \(p_1, \ldots, p_m \ge 2\) & \checkmark & \checkmark & This work \\
\hline
\end{tabular}
\label{table:comparison}
\end{table}

\section{Preliminaries}
\label{sec:background}

In this section, we set up and review definitions and existing facts that will play crucial roles in our analyses. In Section \ref{sec:linalg_background}, we review material from linear algebra, and in Section \ref{sec:convex_geo_background}, we review material from convex geometry.

\subsection{Linear algebra background}
\label{sec:linalg_background}

We introduce a few definitions concerning \textit{leverage scores} (Definition \ref{defn:lev_score}). 
\begin{definition}
\label{defn:lev_score}
For a matrix $\mA \in \R^{n\times d}$, we let $\tau_i(\mA) \coloneqq \va_i^{\top}\inparen{\mA^{\top}\mA}^{-1}\va_i$ denote the \textit{leverage score} of row $\va_i$ with respect to the matrix $\mA$. When $\mA$ is clear from context, we omit it and simply write $\tau_i$ in place of $\tau_i(\mA)$.
\end{definition}

The following are well-known properties of leverage scores.

\begin{fact}
\label{fact:lev_score_facts}
For a matrix $\mA \in \R^{n\times d}$, we have:
\begin{itemize}
    \item $\sum_{i=1}^n \tau_i(\mA) = \rank{\mA}$;
    \item $\tau_i(\mA) = \max_{\vx\in\R^d\setminus\inbraces{0}} \frac{\abs{\ip{\va_i,\vx}}^2}{\norm{\mA\vx}_2^2}$ for all $i$;
    \item $0 \le \tau_i(\mA) \le 1$;
    \item For any positive constant $C$, we have $\tau_i(C\mA) = \tau_i(\mA)$ for all $i$.
\end{itemize}
\end{fact}

We will also need the following fact relating the leverage scores of a matrix $\mA$ to its singular value decomposition.

\begin{fact}
\label{fact:whitening}
Let $\mA\in\R^{n\times d}$ and $\mU\mSigma\mV^{\top}$ be a singular value decomposition for $\mA$, where $\mU\in\R^{n\times d}$ and $\mSigma,\mV\in\R^{d\times d}$. Then, $\tau_i(\mA) = \norm{\vu_i}_2^2$.
\end{fact}
\begin{proof}[Proof of \Cref{fact:whitening}]
To understand why this equality might hold, observe that we can think of $\mU$ as the resulting matrix from applying a statistical whitening transform to $\mA$. More precisely, recall that
\begin{align*}
    \tau_i(\mA) &= \va_i^{\top} \inparen{\mV\mSigma^2\mV^{\top}}^{-1}\va_i^{\top} = \va_i^{\top}\inparen{\sum_{j=1}^d \frac{1}{\sigma_i^2} \cdot \vv_j\vv_j^{\top}}\va_i^{\top} = \sum_{j=1}^d \frac{1}{\sigma_j^2}\ip{\va_i,\vv_j}^2.
\end{align*}
We now calculate $\ip{\va_i,\vv_j}$. Notice that
\begin{align*}
    \ip{\va_i,\vv_j} = \ip{\vv_j, \ve_i^{\top}\sum_{j'=1}^d \sigma_{j'}\mU\ve_{j'}\vv_{j'}^{\top}} = \ip{\vv_j, \sum_{j'=1}^d \sigma_{j'}\ip{\ve_i,\mU\ve_{j'}}\vv_{j'}^{\top}} = \sigma_{j}\ip{\ve_i,\mU\ve_j}.
\end{align*}
Substituting this back in gives
\begin{align*}
    \tau_i(\mA) = \sum_{j=1}^d \ip{\ve_i,\mU\ve_j}^2 = \norm{\vu_i}_2^2.
\end{align*}
This concludes the proof of \Cref{fact:whitening}.
\end{proof}

\subsection{Convex geometry background}
\label{sec:convex_geo_background}

In this subsection, we review foundational facts regarding convex geometry we use throughout the remainder of this chapter. 

We will need the notions of covering and entropy numbers. 
\begin{definition}[Covering numbers {\cite[p. 69]{rothvoss}}]
\label{defn:covering_numbers}
Let $X, Y \subset \R^d$. The covering number $\cN(X,Y)$ is the minimum number of translates of $Y$ required to cover $X$. Formally, we have
\begin{align*}
    \cN(X,Y) \coloneqq \min\inbraces{N \in \N \suchthat \text{ there exists } \vx_1,\dots,\vx_N \in \R^d \text{ such that } X \subseteq \bigcup_{i=1}^N \inparen{\vx_i + Y}}.
\end{align*}
\end{definition}

\begin{definition}[Entropy numbers {\cite[Definition 2.1]{vh15}}]
\label{defn:entropy_numbers}
Let $X, Y \subset \R^d$. The entropy number $e_N(X,Y)$ is the minimum radius $\eta$ such that $\log\cN(X,\eta \cdot Y) \le 2^N$.
\end{definition}

Sometimes, when writing $e_N$, we will write $e_N(X,\norm{\cdot})$ for some quasi-norm $\norm{\cdot}$. Here, we take $Y$ to be the object formed by the unit ball of $\norm{\cdot}$.

Finally, we state the dual Sudakov inequality.
\begin{fact}[Dual Sudakov inequality, due to \citet{ptj87}]
\label{fact:sudakov}
For a symmetric convex body $K \subset \R^d$, define
\begin{align*}
    \norm{\vx}_K \coloneqq \inf\inbraces{t > 0 \suchthat \vx/t \in K}.
\end{align*}
Let $K$ be a symmetric convex body in $\R^d$. We have the below.
\begin{align}
    \log\cN(B_2^d, \eta \cdot K) &\lesssim \eta^{-2} \cdot \exvv{\vg\sim\cN(0,\mI_d)}{\norm{\vg}_K}^2.\label{eq:sudakov_dual}
\end{align}
\end{fact}

\subsection{Probability background}

In this subsection, we review a few facts about subgaussian random variables. These are mostly derived from the presentation of \citet{vershynin_2018}.

\begin{definition}[$\subgnorm{\cdot}$ and subgaussian random variable {\cite[Definition 2.5.6]{vershynin_2018}}]
\label{defn:subgaussian}
Let $X$ be a random variable. Define $\subgnorm{X} \coloneqq \inf\inbraces{t > 0 \suchthat \exv{\expv{\nfrac{X^2}{t^2}}} \le 2}$. If $\subgnorm{X} < \infty$, we say $X$ is subgaussian.
\end{definition}

\begin{fact}[Properties of subgaussian random variables {\cite[Proposition 2.5.2]{vershynin_2018}}]
\label{fact:subgaussian_properties}
The following properties equivalently characterize a subgaussian random variable $X$ up to constants:
\begin{itemize}
    \item For all $t \ge 0$, $\prv{\abs{X} \ge t} \le 2\expv{-\frac{t^2}{\subgnorm{X}^2K_1^2}}$;
    \item For all $r \ge 1$, $\exv{\abs{X}^r} \lesssim r^{r/2}$;
    \item For all $\lambda$ such that $\abs{\lambda} \lesssim \subgnorm{X}^{-1}$, we have $\exv{\expv{\lambda^2X^2}} \le \expv{\subgnorm{X}^2\lambda^2}$.
\end{itemize}
\end{fact}

\begin{fact}[Maximum of subgaussian random variables {\cite[Exercise 2.5.10]{vershynin_2018}}]
\label{fact:max_subgaussians}
Let $X_1,\dots,X_N$ be a sequence of (not necessarily independent) subgaussian random variables. Then
\begin{align*}
    \exv{\max_{i \in [N]} \abs{X_i}} \lesssim \max_{i \in [N]} \subgnorm{X_i}\sqrt{\log N}.
\end{align*}
\end{fact}

\begin{fact}[Decentering]
\label{fact:uncentering}
We have
\begin{align*}
    \subgnorm{X} \lesssim \subgnorm{X - \exv{X}} + \exv{\abs{X}}.
\end{align*}
\end{fact}
\begin{proof}[Proof of \Cref{fact:uncentering}]
By the triangle inequality, we get
\begin{align*}
    \subgnorm{X} \le \subgnorm{X - \exv{X}} + \subgnorm{\exv{X}} \lesssim \subgnorm{X-\exv{X}} + \exv{\abs{X}},
\end{align*}
which is exactly the statement of \Cref{fact:uncentering}.
\end{proof}

\begin{fact}[Lipschitz functionals of Gaussians are subgaussian {\cite[Theorem 5.2.2]{vershynin_2018}}]
\label{fact:lip_concentration}
If $\vg \sim \cN(0,\mI_d)$ and if $\myfunc{f}{\R^d}{\R}$, then
\begin{align*}
    \subgnorm{f(\vg) - \exv{f(\vg)}} \lesssim \norm{f}_{\mathsf{Lip}}.
\end{align*}
\end{fact}

\section{Covering number estimates}
\label{sec:coverings}

In this section, we develop our metric entropy estimates. It will be helpful to keep in mind the context and outline from \Cref{sec:sparsification_overview}.

\subsection{Notation and general formula}

We begin with some definitions that are necessary for our results.

\begin{definition}[Block-constant diagonal matrix]
\label{defn:block_constant}
We say that a vector $\vv \in \R^n$ and corresponding diagonal matrix $\mV\in\R^{n\times n}$ is block-constant or ``constant down the blocks'' if for every $j_1,j_2 \in S_i$, we have $v_{j_1}=v_{j_2}$.
\end{definition}

In particular, when we define a probability measure $\vlambda \in \R^{m}$ over $[m]$, we will find it useful to extend it to a block-constant diagonal matrix $\mLambda\in\R^{n \times n}$.

\begin{definition}[Rounding matrix]
\label{defn:rounding_matrix}
For a probability measure $\vlambda$ over $[m]$, we say that a positive diagonal matrix $\mW \in \R^{k \times k}$ rounds the measure matrix $\mLambda \in \R^{n\times n}$ if for all $\vx\in\R^d$ we have $\norm{\mW^{1/2}\mLambda^{1/2-1/p}\mA\vx}_2 \le \gnorm{\mLambda^{1/2-1/p}\mA\vx}{2}$. We also denote
\begin{align*}
    \widehat{B_2} \coloneqq \inbraces{\vx\in\R^d \suchthat \norm{\mW^{1/2}\mLambda^{1/2-1/p}\mA\vx}_2 \le 1}.
\end{align*}
\end{definition}

The $\widehat{B_2}$ defined in \Cref{defn:rounding_matrix} is the linear transformation of the Euclidean ball that we will ``pass through'' to get our covering number estimates (recall \eqref{eq:intro_split_covering}).

Next, recall our notion of measure overestimates. This is a generalization of prior definitions of Lewis measure overestimates (see e.g. \cite[Definition 2.4]{jls21}, \cite[Definition 2.3]{wy22}) and group leverage score overestimates (\cite[Definition 1.1]{jls22}).

\blocklewisoverestimate*

For example, observe that when all the $S_i$ have size $1$, then \Cref{defn:block_lewis_overestimate} corresponds to standard definitions of Lewis weight overestimates, and there exist weights such that $\Fstar \lesssim d$. Furthermore, we will see that there exist a $\vlambda$ and $\vw$ such that $\Fstar = d$ (it will follow from \Cref{lemma:block_instantiation_small}).

Next, we define the vector $\valpha$, whose entries capture a notion of group importance.

\begin{definition}
\label{defn:alpha}
Let $\vlambda$ be a probability measure over $[m]$ and let $\mW$ be a rounding matrix for $\mLambda$ (\Cref{defn:rounding_matrix}). If $p_1,\dots,p_m \ge 2$, then let $\valpha \in \R^{m}$ be the vector such that for all $i \in [m]$, we have
\begin{align*}
    \alpha_i^p \coloneqq \lambda_i^{1-p/2}\inparen{\sum_{j\in S_i}\inparen{\frac{\tau_j\inparen{\mW^{1/2}\mLambda^{1/2-1/p}\mA}}{w_j}}^{p_i/2}}^{p/p_i}.
\end{align*}
Equivalently, if $\mU$ is a matrix whose columns consist of the left singular vectors of $\mW^{1/2}\mLambda^{1/2-1/p}\mA$, and if we denote by $\vu_j$ the $j$th row of $\mU$ and let $\vf_j \coloneqq \lambda_{i}^{-1/2}w_j^{-1/2}\vu_j$, then by \Cref{fact:whitening}, we may also write
\begin{align*}
    \alpha_i^p &\coloneqq \lambda_i\inparen{\sum_{j \in S_i} \norm{\vf_j}_2^{p_i}}^{p/p_i}.
\end{align*}
On the other hand, if $p_1=\dots=p_m=p<2$, then let $\widehat{\vlambda}$ be a probability measure over $[n]$ and let $\widehat{\valpha}\in\R^n$ be defined as above accordingly. Finally, let $\valpha\in\R^m$ be such that
\begin{align*}
    \alpha_i^p \coloneqq \inparen{\sum_{j \in S_i}\widehat{\alpha}_j^p}^{1/p}.
\end{align*}
\end{definition}

To help ground \Cref{defn:alpha}, notice that combining \Cref{defn:block_lewis_overestimate} with \Cref{defn:alpha} gives us, for $p_1,\dots,p_m\ge 2$ and $p_1,\dots,p_m=p < 2$, respectively,
\begin{equation}
\begin{aligned}
    \alpha_i^p &\le \lambda_i^{1-p/2}\inparen{\lambda_i\Fstar}^{p/2} = \lambda_i\inparen{\Fstar}^{p/2} \\
    \alpha_i^p &= \sum_{j \in S_i} \widehat{\alpha}_j^p = \sum_{j \in S_i} \widehat{\lambda}_j^{1-p/2}\tau_j\inparen{\widehat{\mLambda}^{1/2-1/p}\mA}^{p/2} \le \sum_{j \in S_i}\widehat{\lambda}_j\inparen{\Fstar}^{p/2} = \inparen{\sum_{j\in S_i} \widehat{\lambda_j}}\inparen{\Fstar}^{p/2},
\end{aligned}
\label{eq:alpha_to_fstar}
\end{equation}
and that when $\Fstar \le 2d$ (say), we get $\norm{\valpha}_p^p \lesssim d^{p/2}$. Thus, at least when $p\ge 2$, we can think of $\norm{\valpha}_p^p$ as giving the sparsity we should expect when we sample with probabilities proportional to the $\alpha_i$. Although this does not quite work when $p < 2$, a minor modification of it will.

\begin{definition}[Notation for unit balls and norms under change-of-measure]
\label{defn:unit_ball}
Let $\vlambda$ be a probability measure over $[m]$ and $\mLambda\in\R^{n\times n}$ be its corresponding block-constant diagonal matrix (\Cref{defn:block_constant}). For any $r > 0$, and $\vy\in\R^n$ we define
\begin{align*}
    \gnorml{\vy}{r} \coloneqq \inparen{\sum_{i=1}^m \lambda_i\inparen{\sum_{j \in S_i} \abs{\vy_j}^{p_i}}^{r/p_i}}^{1/r}.
\end{align*}
We also define
\begin{align*}
    B_r \coloneqq\inbraces{\vx \in\R^d \suchthat \gnorml{\mLambda^{-1/p}\mA\vx}{r} \le 1}.
\end{align*}
\end{definition}

From \Cref{defn:unit_ball}, it is easy to verify that $\gnorml{\mLambda^{-1/p}\mA\vx}{p} = \gnorm{\mA\vx}{p}$. Indeed, we have
\begin{align*}
    \gnorml{\mLambda^{-1/p}\mA\vx}{p}^p = \sum_{i=1}^m \lambda_i\norm{\mLambda_{S_i}^{-1/p}\mA_{S_i}\vx}_{p_i}^{p} = \sum_{i=1}^m \norm{\mA_{S_i}\vx}_{p_i}^p = \gnorm{\mA\vx}{p}^p.
\end{align*}
We will also require the crucial property that $\gnorml{\vy}{r}$ is log-convex in $1/r$. To see this, note that the vector in $\R^m$ formed by calculating all the inner norms $p_1,\dots,p_m$ is constant regardless of the outer norm, and then we can use the fact that for a fixed measure $\vmu$, the $\ell_r^m(\vmu)$ norms are log-convex in $1/r$.

We now have the language to state the main result of this section, \Cref{thm:general_covering}.

\begin{theorem}
\label{thm:general_covering}
Let $\vlambda$ be a probability measure and $\mW$ be a rounding matrix (\Cref{defn:rounding_matrix}) so that $\vlambda$ and $\mW$ form an $\Fstar$-block Lewis overestimate (\Cref{defn:block_lewis_overestimate}). Suppose at least one of the following holds:
\begin{itemize}
    \item $p > \frac{1}{\log d}$ and $\abs{S_1}=\dots=\abs{S_m}=1$;
    \item $p = p_1 = \dots = p_m$ and $p < 2$;
    \item $p \ge 1$ and $p_1,\dots,p_m \ge 2$;
    \item $p_1=\dots=p_m=2$ and $1/\log d \le p < \infty$.
\end{itemize}
If $H \ge 1$ is such that the sampling probabilities $\rho_i$ satisfy $H\rho_i \ge \alpha_i^p/\norm{\valpha}_p^p$ for all $i \in [m]$, and if we write $p^{\star}\coloneqq\max\inbraces{1,\max_{i}\min\inbraces{p_i,\log\abs{S_i}}}$, then (recall \Cref{defn:sampling_polyhedron} for the definition of $\polynorm{\cdot}$)
\begin{align*}
    \log\cN\inparen{B_p, \eta\inbraces{\vx\in\R^d \suchthat \polynorm{\vx} \le 1}} \lesssim \eta^{-\min(2,p)} \cdot H^{2/\max(p,2)}\cdot C(p)p^{\star}\Fstar\log\max\inbraces{\mtilde,\Fstar},
\end{align*}
where $C(p)$ is a constant that only depends on $p$.
\end{theorem}

Although \Cref{thm:general_covering} is stated abstractly, we will see that there exists a convenient instantiation for all the parameters stated.

\begin{corollary}
\label{corollary:specific_cover}
In the same cases as in \Cref{thm:general_covering}, there exists a probability measure $\vlambda$ over $[m]$ and a rounding $\mW$ for which in the same setting as \Cref{thm:general_covering}, we have
\begin{align*}
    \log\cN\inparen{B_p, \eta\inbraces{\vx\in\R^d \suchthat \polynorm{\vx} \le 1}} \lesssim_p \eta^{-\min(2,p)}\cdot{\underset{i \in S}{\max}\ \min\inbraces{p_i,\log\abs{S_i}}d\log\mtilde}.
\end{align*}
\end{corollary}

\subsection{Block Lewis weights}
\label{sec:blw_defs_covering}

For the sake of motivation, let us first prove \Cref{corollary:specific_cover} given \Cref{thm:general_covering}. We first need \Cref{lemma:block_lewis}, which is derived from the \textit{block Lewis weights} of \citet{jlls23}.

For a nonnegative diagonal matrix $\mV$, let $\beta_i(\mV) \coloneqq \inparen{\sum_{j \in S_i} \inparen{\va_j^{\top}(\mA^{\top}\mV\mA)^{-1}\va_j}^{p_i/2}}^{1/p_i}$. We call the $\beta_i(\mV)^p$ the \textit{block Lewis weights}.

\begin{lemma}
\label{lemma:block_lewis}
If $p_i \in [2, \infty]$ and $p \in [1, \infty)$, then there exist diagonal $\mV, \mLambda \in \R^{n \times n}$ such that $\vlambda$ is a probability measure over $[m]$ and the corresponding \(\mLambda\in\R^{n\times n}\) is constant on the blocks, then
\(
\sum_{i=1}^m \beta_i(\mV)^p = d
\) and for all $\vx\in\R^d$,
\begin{align*}
    \norm{\mV^{1/2}\mA\vx}_2 \le d^{1/2-1/p}\gnorm{\mLambda^{1/2-1/p}\mA\vx}{2} \le d^{\max(0,1/2-1/p)}\gnorm{\mA\vx}{p}
\end{align*}
\end{lemma}

\begin{proof}[Proof of \Cref{lemma:block_lewis}]
The reader familiar with the work of \citet{jlls23} will notice that \Cref{lemma:block_lewis} is a strengthened variant of Lemma 4.2 from that work. 

Indeed, consider the context of the proof of Lemma 4.2 from \cite{jlls23}. There, notice that $\mW$ is initially chosen so that $\sum_{i=1}^m \beta_i(\mW)^p = d$ and \(\mU = (\mA^\top \mW \mA)^{-1/2}\). We choose $\mV$ in the same way. Next, using their choice of $\vu$, we have $\norm{\vu_{S_i}}_{p_i}^{p-2} = \inparen{\beta_i(\mV)^{p}}^{1-2/p}$. 

Restating (4.8) from \cite{jlls23} in our notation, we have for all $\vx\in\R^d$ that
\begin{align*}
    \norm{\mV^{1/2} \mA \vx}_2^2 \le \sum_{i=1}^m \norm{\vu_{S_i}}_{p_i}^{p-2} \norm{\mA_{S_i} \vx}_{p_i}^2.
\end{align*}
For each \(i \in [m]\) let \(\lambda_i = \frac{\beta_i(\mV)^p}{d}\), so that \(\vlambda\) is a probability measure.
Then
\begin{align*}
    \norm{\mV^{1/2} \mA \vx}_2^2 &\le d^{1-2/p} \sum_{i=1}^m \lambda_i^{1-2/p} \norm{\mA_{S_i} \vx}_{p_i}^2.
\end{align*}
Let \(\mLambda\) be a \(n \times n\) diagonal matrix, where for every \(i \in [m]\) and \(j \in S_i\),
we define \(\mLambda_{jj} = \lambda_i\).
Because \(\lambda_i^{1-2/p} \norm{\mA_{S_i} \vx}_{p_i}^2 = \norm{(\mLambda^{1/2-1/p} \mA)_{S_i} \vx}_{p_i}^2\), we obtain
\begin{equation}\label{eq:vax_atmost_nlax}
    \norm{\mV^{1/2}\mA\vx}_2 \le d^{1/2-1/p}\gnorm{\mLambda^{1/2-1/p}\mA\vx}{2}.
\end{equation}
Since $p$-norms taken with respect to a probability measure are increasing in $p$ we immediately get for all $p \ge 2$ that
\begin{align*}
    \norm{\mV^{1/2}\mA\vx}_2 &\stackrel{\eqref{eq:vax_atmost_nlax}}{\le} d^{1/2-1/p}\gnorm{\mLambda^{1/2-1/p}\mA\vx}{2} = d^{1/2-1/p}\gnorml{\mLambda^{-1/p}\mA\vx}{2} \\
    &\le d^{1/2-1/p}\gnorml{\mLambda^{-1/p}\mA\vx}{p} = d^{1/2-1/p}\gnorm{\mA\vx}{p}.
\end{align*}
The case where $p \le 2$ follows from the ``$1 \le q \le 2$'' subcase of the proof of Lemma 4.2 from \cite{jlls23}, which yields
\begin{align*}
    \norm{\mV^{1/2}\mA\vx}_2 \stackrel{\eqref{eq:vax_atmost_nlax}}{\le} d^{1/2-1/p}\gnorm{\mLambda^{1/2-1/p}\mA\vx}{2} \le \gnorm{\mA\vx}{p}.
\end{align*}
We therefore conclude the proof of \Cref{lemma:block_lewis}.
\end{proof}

We use \Cref{lemma:block_lewis} to give an instantiation for the parameters in \Cref{thm:general_covering}.

\begin{lemma}
\label{lemma:block_instantiation_small}
Let \(\mV, \mLambda\) be the matrices from \Cref{lemma:block_lewis} and let $\mU$ and $\vf_j$ be as defined in \Cref{defn:alpha}. Let $p \ge 1$ and $p_i \ge 2$ for all $i \in [m]$. If we choose \(\mW\) such that
\begin{align*}
    \frac{\mV^{1/2}}{d^{1/2-1/p}} = \mW^{1/2}\mLambda^{1/2-1/p},
\end{align*}
then:
\begin{itemize}
    \item $\norm{\mW^{1/2}\mLambda^{1/2-1/p}\mA\vx}_2 \le \gnorm{\mLambda^{1/2-1/p}\mA\vx}{2}$;
    \item for all $i$, $\frac{\alpha_i}{d^{1/2-1/p}} = \beta_i(\mV)$;
    \item $\norm{\valpha}_p^p = d^{p/2}$;
    \item for all $i$, $\inparen{\sum_{j \in S_i} \norm{\vf_j}_2^{p_i}}^{1/p_i} = \norm{\valpha}_p = d^{1/2}$.
    \item The rounding matrix $\mW$ and measure $\vlambda$ are an $\Fstar$-block Lewis overestimate (\Cref{defn:block_lewis_overestimate}) with $\Fstar = d$.
\end{itemize}
\end{lemma}
\begin{proof}[Proof of \Cref{lemma:block_instantiation_small}]
The first property follows immediately from \Cref{lemma:block_lewis}. Using \Cref{fact:whitening}, notice that
\begin{align*}
    \va_j^{\top}(\mA^{\top}\mV\mA)^{-1}\va_j = \frac{\tau_j(\mV^{1/2}\mA)}{v_j} = \frac{\tau_j(\mW^{1/2}\mLambda^{1/2-1/p}\mA)}{d^{1-2/p}w_j\lambda_i^{1-2/p}} = \frac{\norm{\vu_j}_2^2}{d^{1-2/p}w_j\lambda_i^{1-2/p}} = \frac{\lambda_i^{2/p}\norm{\vf_j}_2^2}{d^{1-2/p}},
\end{align*}
so after substituting into the formula for $\beta_i(\mV)$,
\begin{align*}
    \beta_i(\mV) &= \inparen{\sum_{j \in S_i} \inparen{\va_j^{\top}(\mA^{\top}\mV\mA)^{-1}\va_j}^{p_i/2}}^{1/p_i} = \inparen{\sum_{j \in S_i} \inparen{\frac{\lambda_i^{2/p}\norm{\vf_j}_2^2}{d^{1-2/p}}}^{p_i/2}}^{1/p_i} \\
    &= \frac{\lambda_i^{1/p}\inparen{\sum_{j \in S_i} \norm{\vf_j}_2^{p_i}}^{1/p_i}}{d^{1/2-1/p}} = \frac{\alpha_i}{d^{1/2-1/p}},
\end{align*}
where the last equality follows from the formula for $\valpha$ stated in \Cref{thm:general_covering}. This also implies that
\begin{align*}
    \norm{\valpha}_p^p = \sum_{i=1}^m \alpha_i^p = \sum_{i=1}^p \beta_i(\mV)^pd^{p/2-1} = d^{p/2}.
\end{align*}
Finally, observe that the above calculation shows that $\lambda_i \propto \alpha_i^p$, since we have defined $\lambda_i \propto \beta_i(\mV)^p$ and we have just seen that $\beta_i(\mV)^p \propto \alpha_i^p$. This means we can write $\lambda_i = \alpha_i^p/\norm{\valpha}_p^p$. Using this, we have
\begin{align*}
    \alpha_i = \lambda_i^{1/p}\inparen{\sum_{j \in S_i} \norm{\vf_j}_2^{p_i}}^{1/p_i} = \frac{\alpha_i}{\norm{\valpha}_p} \cdot \inparen{\sum_{j \in S_i} \norm{\vf_j}_2^{p_i}}^{1/p_i}.
\end{align*}
After rearranging, we have
\begin{align*}
    \inparen{\sum_{j \in S_i} \norm{\vf_j}_2^{p_i}}^{1/p_i} = \norm{\valpha}_p = d^{1/2},
\end{align*}
and so we may take $\Fstar = d$. This concludes the proof of \Cref{lemma:block_instantiation_small}.
\end{proof}

We now handle the cases that are not covered by the block Lewis weight construction of \cite{jlls23}.

\begin{lemma}
\label{lemma:block_instantiation_alleq}
If $0 < p_1=\dots=p_m=p < 2$ or if $p_1=\dots=p_m=2$ and $1/\log d \le p < \infty$, then there exists a probability measure $\widehat{\vlambda}$ over $[k]$ and corresponding $\widehat{\valpha} \in \R^n$ such that $\widehat{\vlambda}$ is an $\Fstar$-block Lewis overestimate for $\Fstar=n$.
\end{lemma}
\begin{proof}
For the case where $0 < p_1=\dots=p_m=p < 2$, we simply use the fact that Lewis's measure tells us that there exists a measure $\widehat{\vlambda}$ such that
\begin{align*}
    \frac{\tau_j\inparen{\widehat{\mLambda}^{1/2-1/p}\mA}}{\widehat{\lambda}_j} \le d.
\end{align*}
In the other case, we will see later that the guarantee of a natural contraction mapping (\Cref{alg:blw_lewismap} and \Cref{lemma:contract_alg}) imply that $\mW = \mI_n$ and the resulting $\vlambda$ form an $d$-block Lewis overestimate, thereby concluding the proof of \Cref{lemma:block_instantiation_alleq}.
\end{proof}

\Cref{lemma:block_instantiation_small} and \Cref{lemma:block_instantiation_alleq} easily imply \Cref{corollary:specific_cover}.

\begin{proof}[Proof of \Cref{corollary:specific_cover}]
We combine \Cref{thm:general_covering} with the instantiations in \Cref{lemma:block_instantiation_small} and \Cref{lemma:block_instantiation_alleq}, directly yielding \Cref{corollary:specific_cover}.
\end{proof}

In light of \Cref{corollary:specific_cover}, the goal of the remainder of this section is to prove \Cref{thm:general_covering}. 

It will be useful to consider a corresponding change-of-basis that arises from our setting of $\vlambda$. Let $\mU\mSigma\mV^{\top}$ be a singular value decomposition of $\mW^{1/2}\mLambda^{1/2-1/p}\mA$ where $\mU \in \R^{m \times d}$ and $\mSigma, \mV \in \R^{d \times d}$. Let $\mR$ be the invertible matrix $\mV\mSigma^{-1}$ (we assume without loss of generality that $\rank{\mA} = d$, and it is easy to extend the results of this section to the case where $\rank{\mA} < d$). We take $\mR$ as our change-of-basis matrix. Using this, it is easy to see that $\mW^{1/2}\mLambda^{1/2-1/p}\mA\mR = \mU$ consists of orthonormal columns. Furthermore, we have $\mLambda^{-1/p}\mA\mR = \mW^{-1/2}\mLambda^{-1/2}\mU$.
 
\subsection{Covering numbers for \texorpdfstring{$0 < p < 2$}{0 < p < 2}}

The goal of this section is to prove \Cref{lemma:covering_two_to_one_clean} under the notion of overestimate given by \Cref{defn:block_lewis_overestimate}. 

We are now ready to state the main result of this subsection.

\begin{lemma}
\label{lemma:covering_two_to_one_clean}
Let $\vlambda$ and $\vw$ be such that they form an $\Fstar$-block Lewis overestimate. Then,
\begin{align*}
    \log\cN(B_p,\eta \widehat{B_2}) \lesssim \eta^{-\frac{2p}{2-p}} \cdot C(p)\max_{i}\min\inparen{p_i,\log\abs{S_i}}\Fstar\log\Fstar,
\end{align*}
where $C(p)$ is a constant that only depends on $p$.
\end{lemma}

The goal of the rest of this subsection is to prove \Cref{lemma:covering_two_to_one_clean}. We follow the outline detailed in \Cref{sec:sparsification_overview_measure}. In short, our plan is the following:
\begin{enumerate}
    \item We first reduce bounding $\log\cN(B_p, \eta\widehat{B_2})$ to bounding $\log\cN(\widehat{B_2}, \delta_h B_r)$ for all $h \ge 0$ and appropriate choices of $r$ and $\delta_h$.
    \item We then control each term $\log\cN(\widehat{B_2}, \delta_h B_r)$. To do so, we will apply the dual Sudakov inequality (\Cref{fact:sudakov}, \eqref{eq:sudakov_dual}). To actually estimate $\E\norm{\vg}_{B_r}$ where $\vg \sim \cN(0,\mI_d)$, we need to prove that every resulting summand of the form $\norm{\cdot}_{p_i}$ is subgaussian with a parameter that only depends on $p_i$. To do so, we exploit the fact that these summands are Lipschitz and then apply \Cref{fact:lip_concentration}.
    \item We finally assemble all the previous pieces together to get the desired handle on $\log\cN(B_p, \eta\widehat{B_2})$. 
\end{enumerate}

\subsubsection{Reduction to bounding \texorpdfstring{$\log\cN(\widehat{B_2},\eta B_r)$}{log N(\^{B2}, n Br)}}

As stated in \Cref{sec:sparsification_overview_measure}, we begin with reducing the calculation of $\log\cN(B_p, \eta \widehat{B_2})$ to calculating $\log\cN(\widehat{B_2}, \eta B_r)$ (for a different $\eta$).

\begin{lemma}
\label{lemma:poor_mans_dualization}
Let $\theta$ and $r$ be such that $r=(2-\theta)p/(p-\theta)$. Define
\begin{align*}
    \delta_h \coloneqq \inparen{\frac{8^{h+1}\eta}{2\cdot 8^{2/\theta}}}^{\frac{\theta}{2-\theta}} = \eta^{\frac{\theta}{2-\theta}}\cdot 8^{(h+1)\cdot\frac{\theta}{2-\theta}} \cdot \inparen{2\cdot 8^{2/\theta}}^{-\frac{\theta}{2-\theta}}
\end{align*}
Then, we have
\begin{align*}
     \log\cN(B_p,\eta \widehat{B_2}) &\le \sum_{h\ge 0} \log\cN(\widehat{B_2},\delta_hB_r).
\end{align*}
\end{lemma}
\begin{proof}[Proof of \Cref{lemma:poor_mans_dualization}]
For $h \in \N_{\ge 0}$, let $\cN_h$ be a maximal subset of $B_p$ such that for any two distinct elements $\vz_1, \vz_2 \in B_p$, we have $\norm{\mW^{1/2}\mLambda^{-1/p}\mA(\vz_1-\vz_2)}_{2(\vlambda)} \ge 8^h\eta$ (where by $\norm{\cdot}_{p(\vlambda)}$ we mean the $\ell_p$ norm taken with respect to the measure given by $\vlambda$). This yields $\abs{\cN_h} \ge \cN(B_p, 8^h \eta \widehat{B_2})$.

Next, since for every $h$ there are $\vz_i \in B_p$ for which $B_p \subseteq \bigcup_{i=1}^{\cN(B_p,8^{h+1}\eta \widehat{B_2})} \inbraces{\vz_i + 8^{h+1}\eta \widehat{B_2}}$, for every $h$ there must exist a $\zstar_h \in B_p$ for which
\begin{align*}
    \abs{\inbraces{\zstar_h + 8^{h+1} \eta \widehat{B_2}} \cap \cN_h} \ge \frac{\abs{\cN_h}}{\cN(B_p,8^{h+1} \eta \widehat{B_2})} \ge \frac{\cN(B_p,8^h \eta \widehat{B_2})}{\cN(B_p,8^{h+1} \eta \widehat{B_2})}.
\end{align*}
Let
\begin{align*}
    \cL_h \coloneqq \inbraces{\frac{\vz-\zstar_h}{8^{h+1} \eta} \suchthat \vz \in \inbraces{\zstar_h + 8^{h+1} \eta \widehat{B_2}} \cap \cN_h}
\end{align*}
from which we get by the sub-triangle inequality that
\begin{align*}
    \norm{\mW^{1/2}\mLambda^{-1/p}\mA\vz}_{2(\vlambda)} \le 1 \text{ and } \gnorml{\mLambda^{-1/p}\mA\vz}{p} \le \frac{\max\inbraces{2^{1/p},2}}{8^{h+1} \eta} \text{ for any } \vz \in \cL_h
\end{align*}
and
\begin{align*}
    \norm{\mW^{1/2}\mLambda^{-1/p}\mA\inparen{\vz_1-\vz_2}}_{2(\vlambda)} \ge \frac{1}{8} \text{ for any distinct } \vz_1,\vz_2 \in \cL_h.
\end{align*}
We now apply an interpolation estimate. Let $\vz_1$ and $\vz_2$ be distinct elements from $\cL_h$, set $0<\theta<2$ and $r=(2-\theta)p/(p-\theta)$, and observe that $\theta = p(r-2)/(r-p)$ and 
\begin{align*}
    \frac{1}{8^2} &\le \norm{\mW^{1/2}\mLambda^{-1/p}\mA\inparen{\vz_1-\vz_2}}_{2(\vlambda)}^2 \\
    &\le \gnorml{\mLambda^{-1/p}\mA\inparen{\vz_1-\vz_2}}{2}^2 \\
    &\le \gnorml{\mLambda^{-1/p}\mA\inparen{\vz_1-\vz_2}}{p}^{\theta}\cdot\gnorml{\mLambda^{-1/p}\mA\inparen{\vz_1-\vz_2}}{r}^{2-\theta} \\
    &\le \inparen{\frac{\max\inbraces{2^{1/p},2}}{8^{h+1}\eta}}^{\theta}\gnorml{\mLambda^{-1/p}\mA\inparen{\vz_1-\vz_2}}{r}^{2-\theta} 
\end{align*}
which means that after rearranging we have
\begin{align*}
    \gnorml{\mLambda^{-1/p}\mA\inparen{\vz_1-\vz_2}}{r} \ge \inparen{\inparen{\frac{8^{h+1}\eta}{\max\inbraces{2^{1/p},2}}}^{\theta} \cdot \frac{1}{8^2}}^{\frac{1}{2-\theta}} \ge \delta_h.
\end{align*}
The above argument gives
\begin{align*}
    \log\cN(\widehat{B_2},\delta_h B_r) \ge \log\abs{\cL_h} \ge \log\cN\inparen{B_p,8^h\eta\widehat{B_2}}-\log\cN\inparen{B_p,8^{h+1}\eta\widehat{B_2}}.
\end{align*}
We sum these inequalities over all $h \ge 0$ (noting that when $h$ is sufficiently large, we have $\log\cN(B_p,8^{h+1} \eta \widehat{B_2})=0$), and get
\begin{align*}
    \log\cN(B_p,\eta \widehat{B_2}) &\le \sum_{h\ge 0} \log\cN(\widehat{B_2},\delta_hB_r).
\end{align*}
This concludes the proof of \Cref{lemma:poor_mans_dualization}.
\end{proof}

\subsubsection{Bounding \texorpdfstring{$\log\cN(\widehat{B_2},\eta B_r)$}{log N(\^{B2}, n Br)}}

As we saw in \Cref{lemma:poor_mans_dualization}, it will be enough to understand the behavior of $\log\cN(\widehat{B_2},\eta B_r)$. Since $\widehat{B_2}$ is a linear transformation of a Euclidean ball, we will be able to apply the dual Sudakov inequality (\Cref{fact:sudakov}, \eqref{eq:sudakov_dual}).

To prepare for an application of the dual Sudakov inequality, we bound the Gaussian width of the ball $\inbraces{\vx \in \R^d \suchthat \gnorml{\mLambda^{-1/p}\mA\mR\vx}{r} \le 1}$. As we will see in a moment, the relevance of this ball arises from the fact that it is the $r$-ball with respect to the $\vlambda$ measure after a suitable linear transformation of the underlying space. In particular, it is under the invertible mapping $\vx \mapsto \mR\vx$ that we get $\mW^{1/2}\mLambda^{1/2-1/p}\mA\vx \mapsto \mW^{1/2}\mLambda^{1/2-1/p}\mA\mR\vx = \mU\vx$.

\begin{lemma}
\label{lemma:one_group_subgaussian}
Let $\vg \sim \cN(0,\mI_d)$. We have
\begin{align*}
    \subgnorm{\inparen{\sum_{j \in S_i} \abs{\ip{\vf_j,\vg}}^{p_i}}^{1/p_i}} \lesssim \inparen{1+\sqrt{p_i}}\inparen{\sum_{j \in S_i} \norm{\vf_j}_{2}^{p_i}}^{1/p_i}.
\end{align*}
\end{lemma}
\begin{proof}[Proof of \Cref{lemma:one_group_subgaussian}]
Observe the following Lipschitzness bound, i.e., for any $\vx$, by Cauchy-Schwarz, we have
\begin{align*}
   \inparen{\sum_{j \in S_i} \abs{\ip{\vf_j,\vx}}^{p_i}}^{1/p_i} \le \inparen{\sum_{j \in S_i} \norm{\vf_j}_{2}^{p_i}}^{1/p_i}\norm{\vx}_2
\end{align*}
which means by \Cref{fact:lip_concentration}, we get
\begin{align*}
    \subgnorm{\inparen{\sum_{j \in S_i} \abs{\ip{\vf_j,\vg}}^{p_i}}^{1/p_i}-\exvv{\vg\sim\cN(0,\mI_d)}{\inparen{\sum_{j \in S_i} \abs{\ip{\vf_j,\vg}}^{p_i}}^{1/p_i}}} \lesssim \inparen{\sum_{j \in S_i} \norm{\vf_j}_{2}^{p_i}}^{1/p_i}.
\end{align*}
Now, observe that
\begin{align*}
    \exvv{\vg\sim\cN(0,\mI_d}{\inparen{\sum_{j \in S_i} \abs{\ip{\vf_j,\vg}}^{p_i}}^{1/p_i}} \le \inparen{\sum_{j \in S_i} \norm{\vf_j}_2^{p_i}\exvv{\vg\sim\cN(0,\mI_d)}{\abs{\ip{\frac{\vf_j}{\norm{\vf_j}_2},\vg}}^{p_i}}}^{1/p_i} \asymp p_i^{1/2}\inparen{\sum_{j \in S_i} \norm{\vf_j}_{2}^{p_i}}^{1/p_i},
\end{align*}
and by \Cref{fact:uncentering},
\begin{align*}
    \subgnorm{\inparen{\sum_{j \in S_i} \abs{\ip{\vf_j,\vg}}^{p_i}}^{1/p_i}} \lesssim \inparen{1+\sqrt{p_i}}\inparen{\sum_{j \in S_i} \norm{\vf_j}_{2}^{p_i}}^{1/p_i},
\end{align*}
completing the proof of \Cref{lemma:one_group_subgaussian}.
\end{proof}

Next, we estimate the expected norm of a Gaussian random vector under the norm given by $\inbraces{\vx\in\R^d\suchthat \gnorml{\mLambda^{-1/p}\mA\mR\vx}{r} \le 1}$.

\begin{lemma}
\label{lemma:khintchine}
For $r \ge 2$, we have
\begin{align*}
    \exvv{\vg\sim\cN(0,\mI_d)}{\gnorml{\mLambda^{-1/p}\mA\mR\vg}{r}} \lesssim \begin{cases} r^{1/2}\inparen{1+\sqrt{\max_i p_i}}\inparen{\sum_{i=1}^m \lambda_i \inparen{\sum_{j \in S_i} \norm{\vf_j}_{2}^{p_i}}^{r/p_i}}^{1/r} & \text{ if } r \le \log m \\ \sqrt{\log m} \cdot \max_i \inparen{1+\sqrt{p_i}}\inparen{\sum_{j \in S_i} \norm{\vf_j}_{2}^{p_i}}^{1/p_i} & \text{ otherwise} \end{cases}
\end{align*}
\end{lemma}
\begin{proof}[Proof of \Cref{lemma:khintchine}]
Let $\vu_j$ denote the rows of $\mU$. Note that by \Cref{fact:whitening}, we have
\begin{align*}
    \norm{\vu_j}_2^2 = \tau_j(\mW^{1/2}\mLambda^{1/2-1/p}\mA).
\end{align*}
Now, observe that $\mLambda^{-1/p}\mA\mR = \mW^{-1/2}\mLambda^{-1/2}\mU$. By \Cref{lemma:one_group_subgaussian}, we know that
\begin{align*}
    \exvv{\vg\sim\cN(0,\mI_d)}{\inparen{\sum_{j \in S_i} \abs{\ip{\vf_j,\vg}}^{p_i}}^{r/p_i}} \lesssim r^{r/2}\inparen{1+\sqrt{p_i}}^r\inparen{\sum_{j \in S_i} \norm{\vf_j}_{2}^{p_i}}^{r/p_i}
\end{align*}
We first handle the case where $r \lesssim \log m$. Notice that
\begin{align*}
    \exvv{\vg\sim\cN(0,\mI_d)}{\gnorml{\mLambda^{-1/p}\mA\mR\vg}{r}} &= \exvv{\vg\sim\cN(0,\mI_d))}{\inparen{\sum_{i=1}^m \lambda_i\inparen{\sum_{j \in S_i}\abs{\ip{w_j^{-1/2}\lambda_i^{-1/2}\vu_j,\vg}}^{p_i}}^{r/p_i}}^{1/r}} \\
    &\le \inparen{\sum_{i=1}^m \lambda_i\exvv{\vg\sim\cN(0,\mI_d))}{\inparen{\sum_{j \in S_i}\abs{\ip{w_j^{-1/2}\lambda_i^{-1/2}\vu_j,\vg}}^{p_i}}^{r/p_i}}}^{1/r} \\
    &\lesssim \inparen{\sum_{i=1}^m \lambda_i\inparen{r^{r/2}\inparen{1+\sqrt{p_i}}^r\inparen{\sum_{j \in S_i} \norm{\vf_j}_{2}^{p_i}}^{r/p_i}}}^{1/r} \\
    &\le r^{1/2}\inparen{1+\sqrt{\max_i p_i}}\inparen{\sum_{i=1}^m \lambda_i \inparen{\sum_{j \in S_i} \norm{\vf_j}_{2}^{p_i}}^{r/p_i}}^{1/r}.
\end{align*}
We now handle the case where $r \gtrsim \log m$. We have
\begin{align*}
    \exvv{\vg\sim\cN(0,\mI_d)}{\gnorml{\mLambda^{-1/p}\mA\mR\vg}{r}} &= \exvv{\vg\sim\cN(0,\mI_d))}{\inparen{\sum_{i=1}^m \lambda_i\inparen{\sum_{j \in S_i}\abs{\ip{w_j^{-1/2}\lambda_i^{-1/2}\vu_j,\vg}}^{p_i}}^{r/p_i}}^{1/r}} \\
    &\lesssim \exvv{\vg\sim\cN(0,\mI_d))}{\max_i \inparen{\sum_{j \in S_i}\abs{\ip{w_j^{-1/2}\lambda_i^{-1/2}\vu_j,\vg}}^{p_i}}^{1/p_i}} \\
    &\lesssim \sqrt{\log m} \cdot \max_i \inparen{1+\sqrt{p_i}}\inparen{\sum_{j \in S_i} \norm{\vf_j}_{2}^{p_i}}^{1/p_i}
\end{align*}
and conclude the proof of \Cref{lemma:khintchine} (the last line follows from \Cref{fact:max_subgaussians}).
\end{proof}

Now, we show how to relate $(\sum_{j \in S_i} \norm{\vf_j}_2^{p_i})^{1/p_i}$ to $\Fstar$.

\begin{lemma}
\label{lemma:convert_f_to_fstar}
For all $i \in [m]$, we have
\begin{align*}
    \inparen{\sum_{j \in S_i} \norm{\vf_j}_2^{p_i}}^{2/p_i} \le \Fstar.
\end{align*}
\end{lemma}
\begin{proof}[Proof of \Cref{lemma:convert_f_to_fstar}]
Recall \Cref{fact:whitening}; this gives us
\begin{align*}
    \inparen{\sum_{j \in S_i} \norm{\vf_j}_{2}^{p_i}}^{2/p_i} = \inparen{\sum_{j \in S_i} \inparen{\frac{\tau_j\inparen{\mW^{1/2}\mLambda^{1/2-1/p}\mA}}{w_j\lambda_i}}^{p_i/2}}^{2/p_i} = \frac{1}{\lambda_i}\inparen{\sum_{j \in S_i} \inparen{\frac{\tau_j\inparen{\mW^{1/2}\mLambda^{1/2-1/p}\mA}}{w_j}}^{p_i/2}}^{2/p_i}.
\end{align*}
We recall that $\Fstar$ satisfies \Cref{defn:block_lewis_overestimate} and conclude the proof of \Cref{lemma:convert_f_to_fstar}.
\end{proof}

We now have enough tools to build a na\"ive estimate of $\log\cN\inparen{\widehat{B_2},\eta B_r}$ via directly applying the dual Sudakov inequality.

\begin{lemma}
\label{lemma:r_to_two}
We have
\begin{align*}
    \log\cN\inparen{\widehat{B_2},\eta B_r} \lesssim \eta^{-2} \cdot  \begin{cases} r\inparen{1+\sqrt{\max_i p_i}}^2\inparen{\sum_{i=1}^m \lambda_i \inparen{\sum_{j \in S_i} \norm{\vf_j}_{2}^{p_i}}^{r/p_i}}^{2/r} & \text{ if } r \le \log m \\ \log m \cdot \max_i \inparen{1+\sqrt{p_i}}^2\inparen{\sum_{j \in S_i} \norm{\vf_j}_{2}^{p_i}}^{2/p_i} & \text{ otherwise} \end{cases}
\end{align*}
Simply put, we may also write
\begin{align*}
    \log\cN\inparen{\widehat{B_2},\eta B_r} \lesssim \eta^{-2} \cdot r\max_{i}\min\inparen{p_i, \log\abs{S_i}}\Fstar.
\end{align*}
\end{lemma}
\begin{proof}[Proof of \Cref{lemma:r_to_two}]
Since $\mR$ is invertible, it will be enough to bound the covering number
\begin{align*}
    \cN \coloneqq \cN\inparen{\inbraces{\vx\in\R^d \suchthat \norm{\mU\vx}_2 \le 1}, \eta\inbraces{\vx\in\R^d \suchthat \gnorml{\mLambda^{-1/p}\mA\mR\vx}{r} \le 1}}.
\end{align*}
Because $\norm{\mU\vx}_2 = \norm{\vx}_2$, we can apply the dual Sudakov Inequality (\Cref{fact:sudakov}, \eqref{eq:sudakov_dual}). This means we get
\begin{align*}
    \log\cN \lesssim \eta^{-2} \inparen{\exvv{\vg\sim\cN(0,\mI_d)}{\gnorml{\mLambda^{-1/p}\mA\mR\vg}{r}}}^2.
\end{align*}
We plug in the result from \Cref{lemma:khintchine} and conclude the proof of \Cref{lemma:r_to_two}. The statement after the ``simply put'' follows from \Cref{lemma:convert_f_to_fstar}.
\end{proof}

Although the calculation in \Cref{lemma:r_to_two} works pretty well for small $r$, this degrades quite rapidly once $r$ is large (say, larger than $\log d$).

To resolve this, we build another estimate for $\log\cN\inparen{\widehat{B_2},\eta B_r}$ that performs better when $r$ is larger than $\log d$ or so. We will be able to do this after an interpolation step and a simple geometric observation relating $\widehat{B_2}$ and $B_r$.

\begin{lemma}
\label{lemma:r_ellipsoid}
Let $\Delta_i$ be defined such that 
\begin{align*}
    \Delta_i^{1/2} \coloneqq \max_{\vx\in\R^d\setminus\inbraces{0}} \frac{\lambda_i^{-1/p}\norm{\mA_{S_i}\vx}_{p_i}}{\norm{\mW^{1/2}\mLambda^{1/2-1/p}\mA\vx}_2},
\end{align*}
and let $\Delta \coloneqq \max_{i \in [m]} \Delta_i$.

For all $\vx\in\R^d$ and $r > 2$, if $\abs{S_i} = 1$ for all $i$, then we have
\begin{align*}
    \norm{\mLambda^{-1/2}\mU\vx}_{r(\vlambda)} \le \Delta^{1/2-1/r} \cdot \norm{\mU\vx}_2 \le \Delta^{1/2-1/r} \cdot \norm{\mLambda^{-1/2}\mU\vx}_{r(\vlambda)}.
\end{align*}
Moreover, if there exists at least one $S_i$ for which $\abs{S_i} > 1$, then we have
\begin{align*}
    \gnorml{\mLambda^{-1/p}\mA\vx}{r} \le \Delta^{1/2}\norm{\mW^{1/2}\mLambda^{1/2-1/p}\mA\vx}_2.
\end{align*}
\end{lemma}
\begin{proof}[Proof of \Cref{lemma:r_ellipsoid}]
For the sake of intuition and an interpretation, the reader may think $\Delta \approx d$.

Note that for the case where all the $S_i$ are singletons, we may assume $\mW = \mI_m$.

Since $\vlambda$ is a probability measure, we have for any $r \ge 2$ and for all $\vx\in\R^d$ that
\begin{align*}
    \norm{\mU\vx}_2 = \norm{\mLambda^{-1/2}\mU\vx}_{2(\vlambda)} \le \norm{\mLambda^{-1/2}\mU\vx}_{r(\vlambda)}.
\end{align*}
We now prove the lower bound. We have
\begin{align*}
    \norm{\mLambda^{-1/2}\mU\vx}_{r(\vlambda)} &= \inparen{\sum_{i=1}^m \lambda_i\abs{\ip{\lambda_i^{-1/2}\vu_i,\vx}}^{r}}^{1/r} =  \inparen{\sum_{i=1}^m \lambda_i\abs{\ip{\lambda_i^{-1/2}\vu_i,\vx}}^{2}\cdot\abs{\ip{\lambda_i^{-1/2}\vu_i,\vx}}^{r-2}}^{1/r} \\
    &\le \inparen{\norm{\mU\vx}_2^2 \cdot \max_i \abs{\ip{\lambda_i^{-1/2}\vu_i,\vx}}^{r-2}}^{1/r} \le \inparen{\norm{\mU\vx}_2^2 \cdot \inparen{\Delta^{1/2}\norm{\vx}_2}^{r-2}}^{1/r} \\
    &= \Delta^{1/2-1/r} \norm{\vx}_2.
\end{align*}
We now move onto the more general case where the $S_i$ are allowed to have multiple elements. We write
\begin{align*}
    \gnorml{\mLambda^{-1/p}\mA\vx}{r} &= \inparen{\sum_{i=1}^m \lambda_i\norm{\mLambda_{S_i}^{-1/p}\mA_{S_i}\vx}_{p_i}^r}^{1/r} \\
    &\le \inparen{\sum_{i=1}^m \lambda_i\Delta_i^{r/2}\norm{\mW^{1/2}\mLambda^{1/2-1/p}\mA\vx}_{2}^{r}}^{1/r} \le \Delta^{1/2}\norm{\mW^{1/2}\mLambda^{1/2-1/p}\mA\vx}_2,
\end{align*}
which concludes the proof of \Cref{lemma:r_ellipsoid}.
\end{proof}

At last, we have the tools we need to give a characterization of $\log\cN\inparen{\widehat{B_2},\eta B_r}$ when $r \gtrsim \log d$.

\begin{lemma}
\label{lemma:r_to_two_complicated}
If all the $S_i$ have size $1$, then
\begin{align*}
    \log\cN\inparen{\widehat{B_2}, \eta B_r} \lesssim \inparen{\frac{\eta}{2}}^{-\frac{2r}{r-2}} \cdot r\Fstar\log\Fstar,
\end{align*}
and if there is at least one $S_i$ larger than $1$, then
\begin{align*}
    \log\cN\inparen{\widehat{B_2}, \eta B_r} \lesssim \inparen{\frac{\eta}{2}}^{-\frac{2r}{r-2}} \cdot \frac{2r^2}{r-2}\max_i \min(p_i+1, \log\abs{S_i})\Fstar\log\Fstar.
\end{align*}
\end{lemma}
\begin{proof}[Proof of \Cref{lemma:r_to_two_complicated}]
The reader familiar with the work of \citet{blm89} can think of the present Lemma as a generalization of (7.13) of Proposition 7.2 of that work.

Define
\begin{align*}
    M(\mF,r) &\coloneqq \begin{cases} r\inparen{1+\sqrt{\max_i p_i}}^2\inparen{\sum_{i=1}^m \lambda_i \inparen{\sum_{j \in S_i} \norm{\vf_j}_{2}^{p_i}}^{r/p_i}}^{2/r} & \text{ if } r \le \log m \\ \log m \cdot \max_i \inparen{1+\sqrt{p_i}}^2\inparen{\sum_{j \in S_i} \norm{\vf_j}_{2}^{p_i}}^{2/p_i} & \text{ otherwise} \end{cases}.
\end{align*}

Let $q > r$ and $0 < \theta < 1$ be such that
\begin{align*}
    \frac{1}{r} = \frac{1-\theta}{2} + \frac{\theta}{q}.
\end{align*}
By interpolation, observe that we have
\begin{align*}
    \gnorml{\mLambda^{-1/p}\mA(\vx_1-\vx_2)}{r} &\le \gnorml{\mLambda^{-1/p}\mA(\vx_1-\vx_2)}{2}^{1-\theta} \cdot \gnorml{\mLambda^{-1/p}\mA(\vx_1-\vx_2)}{q}^{\theta} \\
    &\le 2\gnorml{\mLambda^{-1/p}\mA(\vx_1-\vx_2)}{q}^{\theta}
\end{align*}
which means that
\begin{align*}
    \log\cN\inparen{B_2,\eta B_r} \le \log\cN\inparen{B_2,(\eta/2)^{1/\theta}B_q} \le \inparen{\frac{\eta}{2}}^{-2/\theta}M(\mF,q).
\end{align*}
Let us set $q = r\log D$, where we will choose $D$ in a moment. Then, notice that
\begin{align*}
    \inparen{\frac{\eta}{2}}^{-2/\theta} = \inparen{\frac{\eta}{2}}^{-2r(q-2)/(q(r-2))} = \inparen{\frac{\eta}{2}}^{-\frac{2r}{r-2}\inparen{1-\frac{2}{r\log D}}} = \inparen{\frac{\eta}{2}}^{-\frac{2r}{r-2} + \frac{4}{(r-2)\log D}}.
\end{align*}
It is now sufficient to identify $D$ such that whenever $\eta$ is small enough to have $\log\cN > 0$, we have
\begin{align*}
    \inparen{\frac{\eta}{2}}^{\frac{4}{(r-2)\log D}} \lesssim 1.
\end{align*}
To identify this $D$, notice that \Cref{lemma:r_ellipsoid} implies that if all the $S_i$ have size $1$, then only values of $\eta$ such that $\eta \le \Delta^{1/2-1/r}$ contribute to $\log\cN$. In the more general setting, observe only $\eta \le \Delta^{1/2}$ counts.

Hence, if all the $S_i$s are singletons, we choose $D = \Delta$. For any $\eta \le 2\Delta^{1/2-1/r} = 2D^{1/2-1/r}$, we see that
\begin{align*}
    \inparen{\frac{\eta}{2}}^{\frac{4}{(r-2)\log D}} \le \Delta^{\frac{r-2}{2r} \cdot \frac{4}{(r-2)\logv{\Delta}}} = \Delta^{\frac{2}{r\log\Delta}} = 2^{2/r} \le 2.
\end{align*}
Similarly, for the case where the $S_i$ are more generally sized, we choose $D = \Delta^{2r/(r-2)}$. Now, for any $\eta \le 2\Delta^{1/2}$, we get
\begin{align*}
    \inparen{\frac{\eta}{2}}^{\frac{4}{(r-2)\log D}} \le \Delta^{\frac{4}{(r-2)\logv{\Delta^{(2r/(r-2))}}}} = \Delta^{\frac{2}{r\log\Delta}} = 2^{2/r} \le 2.
\end{align*}
Putting everything together, if all the $S_i$s are singletons, we get
\begin{align*}
    \log\cN\inparen{\widehat{B_2}, \eta B_r} \lesssim \inparen{\frac{\eta}{2}}^{-\frac{2r}{r-2}} \cdot M(\mF,q) \lesssim \inparen{\frac{\eta}{2}}^{-\frac{2r}{r-2}} \cdot r\Fstar\log\Fstar,
\end{align*}
and in the more general case,
\begin{align*}
    \log\cN\inparen{\widehat{B_2}, \eta B_r} \lesssim \inparen{\frac{\eta}{2}}^{-\frac{2r}{r-2}} \cdot M(\mF,q) \lesssim \inparen{\frac{\eta}{2}}^{-\frac{2r}{r-2}} \cdot \frac{2r^2}{r-2}\max_i \min(p_i+1, \log\abs{S_i})\Fstar\log\Fstar.
\end{align*}
This concludes the proof of \Cref{lemma:r_to_two_complicated}.
\end{proof}

\subsubsection{Putting everything together}

We are finally ready to combine all the tools we have built in the last few subsections to prove our entropy estimate when $0 < p < 2$.

Below, we state \Cref{lemma:covering_two_to_one}, which more precisely characterizes the behavior of the dependence on $p$ referred to by \Cref{lemma:covering_two_to_one_clean}. 

\begin{lemma}
\label{lemma:covering_two_to_one}
We have
\begin{align*}
    \log\cN(B_p,\eta \widehat{B_2}) \lesssim \eta^{-\frac{2p}{2-p}} \cdot C(p)\max_{i}\min\inparen{p_i+1,\log\abs{S_i}}\Fstar\log\Fstar,
\end{align*}
where $C(p)$ is a constant that only depends on $p$. The constant $C(p)$ is defined as follows. For $0 < \theta < p$, let $r=(2-\theta)p/(p-\theta)$. Then, we define
\begin{align}
    \widehat{C}(p,\theta) &\coloneqq \inparen{\frac{\inparen{2\cdot 8^{2/\theta}}^{-\frac{\theta}{2-\theta}}}{2}}^{-\frac{2r}{r-2}}\nonumber \\
    \widehat{\widehat{C}}(p) &\coloneqq \begin{cases} \min\inbraces{\widehat{C}(p,p/2), \widehat{C}(p,1)} & \text{ if } 1 \le p < 2 \\ \widehat{C}(p,p/2) & \text{ if } 0 < p < 1 \end{cases}\label{eq:hathat} \\
    C(p) &\coloneqq \widehat{\widehat{C}}(p)\begin{cases} r & \text{ if } \abs{S_i} = 1 \text{ for all } i \\ 2r + \frac{4r}{r-2} & \text{ otherwise} \end{cases}\label{eq:choosing_c_from_hathat},
\end{align}
where, in an abuse of notation, $r$ in \eqref{eq:choosing_c_from_hathat} is chosen according to the value of $\theta$ that is selected by $\widehat{\widehat{C}}(p)$ in \eqref{eq:hathat}.
\end{lemma}
\begin{proof}[Proof of \Cref{lemma:covering_two_to_one} and \Cref{lemma:covering_two_to_one_clean}]
This time, following \Cref{lemma:r_to_two_complicated}, we define
\begin{align*}
    M(\mF,r) &\coloneqq \begin{cases} r\Fstar\log\Fstar & \text{ if } \abs{S_i}=1 \\ (2r + \frac{4r}{r-2})\inparen{\max_i \min(p_i+1, \log\abs{S_i})}\Fstar\log\Fstar & \text{ otherwise} \end{cases}.
\end{align*}
Use \Cref{lemma:poor_mans_dualization} and \Cref{lemma:r_to_two} to write
\begin{align}
    \log\cN(B_p,\eta \widehat{B_2}) &\le \sum_{h\ge 0} \log\cN(\widehat{B_2},\delta_hB_r) \le M(\mF,r)\sum_{h \ge 0} \inparen{\frac{\delta_h}{2}}^{-\frac{2r}{r-2}}\nonumber \\
    &= M(\mF,r) \sum_{h \ge 0} \inparen{\eta^{\frac{\theta}{2-\theta}}\cdot 8^{(h+1)\cdot\inparen{\frac{\theta}{2-\theta}}} \cdot \frac{\inparen{2\cdot 8^{2/\theta}}^{-\frac{\theta}{2-\theta}}}{2}}^{-\frac{2r}{r-2}} \nonumber \\
    &= \eta^{-\frac{\theta}{2-\theta}\cdot\frac{2r}{r-2}} \cdot \underbrace{\inparen{\frac{\inparen{2\cdot 8^{2/\theta}}^{-\frac{\theta}{2-\theta}}}{2}}^{-\frac{2r}{r-2}}}_{= \widehat{C}(p,\theta)} M(\mF,r)\sum_{h \ge 0} 8^{(h+1) \cdot \inparen{-\frac{\theta}{2-\theta}}\cdot\frac{2r}{r-2}} \label{eq:substitute_theta}.
\end{align}
We now make the substitution $\theta=p/2$. By the formula in \Cref{lemma:poor_mans_dualization}, this means that $r = 4-p$ and $-\theta/(2-\theta) \cdot 2r/(r-2) = -2p/(2-p)$. We continue.
\begin{align*}
    \log\cN(B_p,\eta B_2) &\le \eta^{-\frac{\theta}{2-\theta}\cdot\frac{2r}{r-2}} \cdot \widehat{C}(p,\theta) M(\mF,r)\sum_{h \ge 0} 8^{(h+1) \cdot -\frac{\theta}{2-\theta}\cdot\frac{2r}{r-2}} \\
    &= \eta^{-\frac{2p}{2-p}} \cdot \widehat{C}(p,p/2) M(\mF,4-p)\sum_{h \ge 0} 8^{(h+1) \cdot -\frac{2p}{2-p}} \\
    &\lesssim \eta^{-\frac{2p}{2-p}} \cdot \widehat{C}(p,p/2) M(\mF,4-p).
\end{align*}
A regrettable consequence of the above calculation is that the ``constant'' $C(p,\theta)=C(p,p/2)$ explodes as $p\rightarrow 2$, as observed by \citet{sz01}. To fix this, we perform a slightly different variant of this calculation in the regime where $1 < p < 2$. We resume from \eqref{eq:substitute_theta} except we use $\theta = 1$. Here, again using \Cref{lemma:poor_mans_dualization}, we check that $r=p/(p-1)$ and $-\theta/(2-\theta) \cdot 2r/(r-2) = -2p/(2-p)$ (note that now $r$ is the conjugate exponent of $p$). This means that
\begin{align*}
    \log\cN(B_p,\eta \widehat{B_2}) &\le \eta^{-\frac{\theta}{2-\theta}\cdot\frac{2r}{r-2}} \cdot \widehat{C}(p,\theta) M(\mF,r)\sum_{h \ge 0} 8^{(h+1) \cdot -\frac{\theta}{2-\theta}\cdot\frac{2r}{r-2}} \\
    &\lesssim \eta^{-\frac{2p}{2-p}} \cdot \widehat{C}(p,1)\sum_{h \ge 0} 8^{(h+1) \cdot -\frac{2p}{2-p}} \lesssim \eta^{-\frac{2p}{2-p}}M\inparen{\mF,\frac{p}{p-1}}.
\end{align*}
Taking the minimum over all the cases (of course where applicable) and expanding out the definition of $M(\mF,r)$ concludes the proof of \Cref{lemma:covering_two_to_one} and \Cref{lemma:covering_two_to_one_clean}.
\end{proof}

\subsection{Covering numbers for \texorpdfstring{$p \ge 2$}{p > 2}}

We will see that compared to the previous section, our task when $p\ge 2$ is far easier. The main technical lemma we need is \Cref{lemma:covering_poly_to_two}, which we need for both regimes of $p$.

We first state and prove \Cref{lemma:covering_poly_to_two}.

\begin{lemma}
\label{lemma:covering_poly_to_two}
Let $\mW$ and $\mLambda$ be chosen according to \Cref{thm:general_covering}. Suppose that $H \ge 1$ is such that $H\rho_i \ge \nfrac{\alpha_i^p}{\norm{\valpha}_p^p}$ for all $i \in [m]$. Then,
\begin{align*}
    \log\cN\inparen{\widehat{B_2}, \eta\inbraces{\vx\in\R^d \suchthat \polynorm{\vx} \le 1}} \lesssim \eta^{-2}\cdot{\underset{i \in S}{\max}\ \min\inbraces{p_i,\log\abs{S_i}}H^{2/p}\norm{\valpha}_p^2\log\mtilde}.
\end{align*}
\end{lemma}
\begin{proof}[Proof of \Cref{lemma:covering_poly_to_two}]
Following the proof of \Cref{lemma:r_to_two} and the references therein, our goal here is to analyze the quantity
\begin{align*}
    \log\cN\inparen{\inbraces{\vx\in\R^d \suchthat \norm{\mW^{1/2}\mLambda^{1/2-1/p}\mA\vx}_2 \le 1}, \eta\inbraces{\vx\in\R^d \suchthat \polynorm{\vx} \le 1}}.
\end{align*}
Using the same type of linear transformation argument as in \Cref{lemma:r_to_two} (so, replacing every $\vx$ above with $\mR\vx$), we find that it is in fact sufficient to analyze
\begin{align*}
    \log\cN \coloneqq \log\cN\inparen{\inbraces{\vx \suchthat \norm{\mU\vx}_2 \le 1}, \eta\inbraces{\vx \suchthat \max_{i \in S}\ \rho_i^{-1/p}\norm{\mW_{S_i}^{-1/2}\mLambda^{1/p-1/2}_{S_i}\mU\vx}_{p_i} \le 1}}.
\end{align*}
Recall that $\norm{\mU\vx}_2 = \norm{\vx}_2$, so a natural plan is to apply the dual Sudakov inequality (\Cref{fact:sudakov}, \eqref{eq:sudakov_dual}). We first consider the quantity (when $2 \le p_i \le \log\abs{S_i}$)
\begin{align*}
    \exvv{\vg\sim\cN(0,\mI_d)}{\norm{\mW_{S_i}^{-1/2}\mLambda^{1/p-1/2}_{S_i}\mU\vg}_{p_i}} &= \exvv{\vg\sim\cN(0,\mI_d)}{\inparen{\sum_{j \in S_i} \abs{\lambda_i^{1/p}\ip{w_{j}^{-1/2}\lambda_i^{-1/2}\vu_j,\vg}}^{p_i}}^{1/p_i}} \\
    &\le \inparen{\sum_{j\in S_i} \inparen{\lambda_i^{1/p}\norm{w_j^{-1/2}\lambda_i^{-1/2}\vu_j}_2}^{p_i}\exvv{g\sim\cN(0,1)}{\abs{g}^{p_i}}}^{1/p_i} \\
    &\lesssim p_i^{1/2} \inparen{\sum_{j\in S_i} \inparen{\lambda_i^{1/p}\norm{w_j^{-1/2}\lambda_i^{-1/2}\vu_j}_2}^{p_i}}^{1/p_i} \\
    &= p_i^{1/2}\inparen{\sum_{j \in S_i} \inparen{\lambda_i^{1/p}\norm{\vf_j}_2}^{p_i}}^{1/p_i} = p_i^{1/2}\lambda_i^{1/p}\inparen{\sum_{j \in S_i} \norm{\vf_j}_2^{p_i}}^{1/p_i} = p_i^{1/2}\alpha_i.
\end{align*}
On the other hand, if $p_i \ge \log\abs{S_i}$, then we get
\begin{align*}
    \exvv{\vg\sim\cN(0,\mI_d)}{\norm{\mW^{-1/2}_{S_i}\mLambda^{1/p-1/2}_{S_i}\mU\vg}_{p_i}} &\lesssim \exvv{\vg\sim\cN(0,\mI_d)}{\norm{\mW^{-1/2}_{S_i}\mLambda^{1/p-1/2}_{S_i}\mU\vg}_{\infty}} \\
    &\lesssim  \lambda_i^{1/p}\max_{j\in S_i}\norm{\vf_j}_2\sqrt{\log\abs{S_i}} \asymp\alpha_i\sqrt{\log\abs{S_i}}.
\end{align*}
Finally, when $p_1=\dots=p_m=p < 2$, we have
\begin{align*}
    \exvv{\vg\sim\cN(0,\mI_d)}{\norm{\widehat{\mLambda}_{S_i}^{1/p-1/2}\mU\vg}_p} &= \exvv{\vg\sim\cN(0,\mI_d)}{\inparen{\sum_{j \in S_i}\widehat{\lambda}_j\abs{\ip{\widehat{\lambda}_j^{-1/2}\vu_j,\vg}}^{p}}^{1/p}} \\
    &\le \inparen{\sum_{j\in S_i} \widehat{\lambda}_j\norm{\widehat{\lambda}_j^{-1/2}\vu_j}_2^p\exvv{g\sim\cN(0,1)}{\abs{g}^{p}}}^{1/p} \lesssim p^{1/2}\inparen{\sum_{j \in S_i} \widehat{\alpha}_j^p}^{1/p} = p^{1/2}\alpha_i.
\end{align*}
This means that throughout the rest of the proof, we assume without loss of generality that $p_i \le \log\abs{S_i}$. Next, observe that for any $\vx\in\R^d$ and when $p_1,\dots,p_m \ge 2$,
\begin{align*}
    \norm{\mW^{-1/2}_{S_i}\mLambda^{1/p-1/2}_{S_i}\mU\vx}_{p_i} = \inparen{\sum_{j \in S_i} \abs{\lambda_i^{1/p}\ip{\vf_j,\vx}}^{p_i}}^{1/p_i} \le \inparen{\sum_{j \in S_i} \abs{\lambda_i^{1/p}\norm{\vf_j}_2\norm{\vx}_2}^{p_i}}^{1/p_i} = \alpha_i\norm{\vx}_2.
\end{align*}
Similarly, when $p_1=\dots=p_m=p$, we first observe
\begin{align*}
    \abs{\ip{\widehat{\lambda}_j^{1/p-1/2}\vu_j,\vx}} \le \widehat{\alpha}_j\norm{\vx}_2.
\end{align*}
We therefore get
\begin{align}
    \norm{\widehat{\mLambda}_{S_i}^{1/p-1/2}\mU\vx}_p \le \inparen{\sum_{j \in S_i} \widehat{\alpha}_j^p}^{1/p}\norm{\vx}_2 = \alpha_i\norm{\vx}_2.\label{eq:alpha_alleq_lipschitz}
\end{align}

Hence, after applying \Cref{fact:lip_concentration} and \Cref{fact:uncentering}, we notice that $\alpha_i(H\rho_i)^{-1/p}\le\norm{\valpha}_p$ and get
\begin{align*}
    \subgnorm{\rho_i^{-1/p}\norm{\mLambda^{1/p-1/2}_{S_i}\mU\vg}_{p_i}} \lesssim \inparen{1 + \sqrt{p_i}}H^{1/p}\norm{\valpha}_p.
\end{align*}
All of this implies that for any subset $S$ of size $\mtilde$ (see Exercise 2.5.10 of \cite{vershynin_2018}),
\begin{align*}
    \exvv{\vg\sim\cN(0,\mI_d)}{\max_{i \in S} \rho_i^{-1/p}\norm{\mLambda^{1/p-1/2}_{S_i}\mU\vg}_{p_i}} \lesssim \max_{i \in S} p_i^{1/2}H^{1/p}\norm{\valpha}_p\sqrt{\log\mtilde}.
\end{align*}
Thus,
\begin{align*}
    \log\cN \lesssim \eta^{-2}\inparen{\inparen{\max_{i \in S} p_i+1}H^{2/p}\norm{\valpha}_p^2\log\mtilde}.
\end{align*}
This concludes the proof of \Cref{lemma:covering_poly_to_two}.
\end{proof}

We are finally ready to prove \Cref{thm:general_covering}.

\begin{proof}[Proof of \Cref{thm:general_covering}]
For notational simplicity in this proof, write
\begin{align*}
    K \coloneqq \inbraces{\vx\in\R^d\suchthat \polynorm{\vx} \le 1}.
\end{align*}
Let us first handle the case where $p \ge 2$. By our choice of $\mW$, we have for all $\vx\in\R^d$ that
\begin{align*}
     \norm{\mW^{1/2}\mLambda^{1/2-1/p}\mA\vx}_{2} \le \gnorml{\mLambda^{-1/p}\mA\vx}{p} = \gnorm{\mA\vx}{p}.
\end{align*}
This implies the containment $B_p \subseteq \widehat{B_2}$, and thus
\begin{align*}
    \log\cN\inparen{B_p, \eta K} \le \log\cN\inparen{\widehat{B_2}, \eta K} \le \eta^{-2}\max_{i}\min\inbraces{p_i,\log\abs{S_i}}H^{2/p}\Fstar\log\mtilde.
\end{align*}
The desired result now follows immediately from \Cref{lemma:covering_poly_to_two}.

For the case where $p < 2$ and $p_i \ge 2$ for all $i$, we require a bit more work. By our choice of $\mW$, we have
\begin{align*}
    \norm{\mW^{1/2}\mLambda^{1/2-1/p}\mA\vx}_{2} \le \gnorml{\mLambda^{-1/p}\mA\vx}{2} = \norm{\mLambda^{1/2-1/p}\mA\vx}_2,
\end{align*}
and for any $t > 0$ (after remembering \eqref{eq:alpha_to_fstar} which tells us $\norm{\valpha}_p \le \sqrt{\Fstar}$),
\begin{align*}
    \log\cN\inparen{B_p, \eta K} &\le \log\cN\inparen{B_p, t\widehat{B_2}} + \log\cN\inparen{t\widehat{B_2}, \eta K} = \log\cN\inparen{B_p, t\widehat{B_2}} + \log\cN\inparen{\widehat{B_2}, \frac{\eta}{t}\cdot K} \\
    &\le t^{-\frac{2p}{2-p}} \cdot C(p)\max_{i}\min\inbraces{p_i+1,\log\abs{S_i}}\Fstar\log\Fstar \\
    &\quad+\inparen{\frac{\eta}{t}}^{-2}\max_{i}\min\inbraces{p_i+1,\log\abs{S_i}}H^{2/p}\Fstar\log\mtilde,
\end{align*}
where the last line follows from \Cref{lemma:covering_two_to_one_clean} and \Cref{lemma:covering_poly_to_two}.
Choose $t=\eta^{1-p/2} \cdot H^{1/2-1/p}$, and for simplicity let $p^{\star} \coloneqq \max_i\min\inbraces{p_i+1,\log\abs{S_i}}$. We write
\begin{align*}
    \log\cN\inparen{B_p, \eta K} &\le t^{-\frac{2p}{2-p}} \cdot C(p)p^{\star}\Fstar\log\Fstar +\inparen{\frac{\eta}{t}}^{-2}p^{\star}H^{2/p}\Fstar\log\mtilde \\
    &\lesssim \eta^{-p}\cdot H \cdot C(p)\inparen{p^{\star}\Fstar\log\max\inbraces{\mtilde,\Fstar}}.
\end{align*}
Finally, we need to address the case where $p_1=\dots=p_m=p < 2$, as this is not covered by the construction of the block Lewis weights (\Cref{lemma:block_lewis}).

To do so, notice that $\gnorm{\mA\vx}{p} = \norm{\mA\vx}_p$ for all $\vx\in\R^d$. So, we reduce to the case where all the $S_i$ have size $1$. In particular, let $\widehat{\vlambda}$ be a probability measure over $[k]$ such that for all $j\in [k]$,
\begin{align*}
    \Fstar \ge \frac{\tau_j\inparen{\widehat{\mLambda}^{1/2-1/p}\mA}}{\widehat{\lambda}_j},
\end{align*}
and let
\begin{align*}
    \widehat{B_2} = \inbraces{\vx\in\R^d\suchthat \norm{\widehat{\mLambda}^{1/2-1/p}\mA\vx}_2 \le 1}.
\end{align*}
Notice that this setup is in accordance with \Cref{defn:alpha}. By \Cref{lemma:covering_two_to_one_clean}, we get
\begin{align*}
    \log\cN(B_p,t\widehat{B_2}) \lesssim t^{-\frac{2p}{2-p}} \cdot C(p)\Fstar\log\Fstar.
\end{align*}
It remains to bound $\log\cN(\widehat{B_2},(\eta/t)K)$. Let $\vlambda$ be a probability measure over $[m]$, where $\lambda_i = \sum_{j \in S_i} \widehat{\lambda}_j$. Define $\widehat{\valpha} \in \R^n$ analogously to $\widehat{\vlambda}$. Notice that, for these choices of $\widehat{\valpha}$ and $\widehat{\vlambda}$, we see that the conclusion of \Cref{lemma:covering_poly_to_two} still holds, and we have
\begin{align*}
    \log\cN\inparen{\widehat{B_2},\inparen{\frac{\eta}{t}} \cdot K} \lesssim \inparen{\frac{\eta}{t}}^{-2}H^{2/p}\Fstar\log\mtilde.
\end{align*}
Now, the calculation is the same as before, and we conclude the proof of \Cref{thm:general_covering}.
\end{proof}

\subsection{Volume-based metric entropy}

In this subsection, we prove \Cref{lemma:covering_easy}, which is an easy consequence of a volume-based argument to obtain a covering number guarantee.

We start with \Cref{lemma:containments}.

\begin{lemma}
\label{lemma:containments}
Let $S \subseteq [m]$ have size $\mtilde$. Let $H \ge 1$ be such that $H\rho_i \ge \nfrac{\alpha_i^p}{\norm{\valpha}_p^p}$ for all $i \in [m]$. For all $\vx\in\R^d$, we have
\begin{align*}
    \max_{i \in S} \frac{\norm{\mA_{S_i}\vx}_{p_i}}{\inparen{\Fstar}^{\max(1/2,1/p)}\rho_i^{1/p}H^{1/p}} &\le \gnorm{\mA\vx}{p}.
\end{align*}
\end{lemma}
\begin{proof}[Proof of \Cref{lemma:containments}]
Consider the invertible mapping $\vx \mapsto \mR\vx$, and write
\begin{align}
    \norm{\mW^{-1/2}_{S_i}\mLambda^{1/p-1/2}_{S_i}\mU\vx}_{p_i} = \inparen{\sum_{j \in S_i} \abs{\lambda_i^{1/p}\ip{\vf_j,\vx}}^{p_i}}^{1/p_i} \le \inparen{\sum_{j \in S_i} \abs{\lambda_i^{1/p}\norm{\vf_j}_2\norm{\vx}_2}^{p_i}}^{1/p_i} = \alpha_i\norm{\vx}_2.\label{eq:alpha_importance_proof}
\end{align}
This means that when $p \ge 2$,
\begin{align*}
    \norm{\mA_{S_i}\vx}_{p_i} \le \alpha_i\norm{\mW^{1/2}\mLambda^{1/2-1/p}\mA\vx}_2 \le \alpha_i\gnorm{\mLambda^{1/2-1/p}\mA\vx}{2} \le \alpha_i\gnorm{\mA\vx}{p} \le \norm{\valpha}_p\rho_i^{1/p}H^{1/p}\gnorm{\mA\vx}{p}.
\end{align*}
Dividing both sides by $\norm{\valpha}_p\rho_i^{1/p}H^{1/p}$ and then recalling \eqref{eq:alpha_to_fstar} yields the desired conclusion (in particular, we see that $\norm{\valpha}_p \le \inparen{\Fstar}^{1/2}$).

We now analyze what happens when $p \le 2$ and $p_1,\dots,p_m \ge 2$. Let
\begin{align*}
    \Delta_i^{1/2} \coloneqq \max_{\vx\in\R^d\setminus\inbraces{0}}\frac{\lambda_i^{-1/p}\norm{\mA_{S_i}\vx}_{p_i}}{\gnorm{\mLambda^{1/2-1/p}\mA\vx}{2}}.
\end{align*}
Let $\Delta \coloneqq \max_{i \in [m]} \Delta_i$. We now see that
\begin{align*}
    &\ \gnorm{\mLambda^{1/2-1/p}\mA\vx}{2} = \inparen{\sum_{i=1}^m \lambda_i\norm{\mLambda_{S_i}^{-1/p}\mA_{S_i}\vx}_{p_i}^2}^{1/2} = \inparen{\sum_{i=1}^m \lambda_i\norm{\mLambda_{S_i}^{-1/p}\mA_{S_i}\vx}_{p_i}^{p}\cdot\norm{\mLambda_{S_i}^{-1/p}\mA_{S_i}\vx}_{p_i}^{2-p}}^{1/2} \\
    &\le \inparen{\gnorm{\mA\vx}{p}^p\inparen{\Delta^{1/2} \cdot \gnorm{\mLambda^{1/2-1/p}\mA\vx}{2}}^{2-p}}^{1/2} = \gnorm{\mA\vx}{p}^{p/2}\inparen{\Delta^{1/2} \cdot \gnorm{\mLambda^{1/2-1/p}\mA\vx}{2}}^{1-p/2}.
\end{align*}
Rearranging and taking the $2/p$-power gives
\begin{align*}
    \gnorm{\mLambda^{1/2-1/p}\mA\vx}{2} \le \gnorm{\mA\vx}{p} \Delta^{1/p-1/2}.
\end{align*}
Observe that this yields the inequalities
\begin{align*}
     \frac{\norm{\mA_{S_i}\vx}_{p_i}}{\Delta^{1/p-1/2}} &\le \frac{\alpha_i}{\Delta^{1/p-1/2}}\norm{\mW^{1/2}\mLambda^{1/2-1/p}\mA\vx}_2 \le \frac{\alpha_i}{\Delta^{1/p-1/2}}\gnorm{\mLambda^{1/2-1/p}\mA\vx}{2} \\
     &\le \alpha_i\gnorm{\mA\vx}{p} \le \norm{\valpha}_p\rho_i^{1/p}H^{1/p}\gnorm{\mA\vx}{p}.
\end{align*}
It remains to bound the $\Delta_i$. By the same sort of argument from the $p\ge 2$ case (i.e., \eqref{eq:alpha_importance_proof}), we get
\begin{align*}
    \lambda_i^{1/p}\Delta_i^{1/2} \le \alpha_i = \lambda_i^{1/p}\inparen{\sum_{j \in S_i} \norm{\vf_j}_2^{p_i}}^{1/p_i}
\end{align*}
which means that
\begin{align*}
    \Delta \le \max_{i \in [m]} \inparen{\sum_{j \in S_i} \norm{\vf_j}_2^{p_i}}^{2/p_i} \le \Fstar.
\end{align*}
Now, again using the fact that $\norm{\valpha}_p \le \inparen{\Fstar}^{1/2}$, we see that $\Delta^{1/p-1/2}\norm{\valpha}_p \le \inparen{\Fstar}^{1/p}$. 

Finally, we analyze the case where $p_1=\dots=p_m=p < 2$. Using \eqref{eq:alpha_alleq_lipschitz}, and the same sort of argument from above, we have
\begin{align*}
    \norm{\mA_{S_i}\vx}_p \le \alpha_i\norm{\widehat{\mLambda}^{1/2-1/p}\mA\vx}_2 \le \alpha_i\inparen{\Fstar}^{1/p-1/2}\gnorm{\mA\vx}{p} \le \norm{\valpha}_pH^{1/p}\rho_i^{1/p}\inparen{\Fstar}^{1/p-1/2}\gnorm{\mA\vx}{p}.
\end{align*}
Once again, we use $\norm{\valpha}_p \le \inparen{\Fstar}^{1/2}$.

We have covered all our cases and may conclude the proof of \Cref{lemma:containments}.
\end{proof}

\Cref{lemma:containments} also suggests a useful sanity check, as the denominator points to a sparsity of $\sum_{i \le m} \inparen{\inparen{\Fstar}^{\max(1/2,1/p)}\rho_i^{1/p}H^{1/p}}^p = H\inparen{\Fstar}^{\max(1,p/2)}$. And, recall that we should be able to set $\Fstar \sim d$, which indeed gives us the dependence on $d$ we see in \Cref{thm:one_shot_lewis}. 

\begin{lemma}
\label{lemma:covering_easy}
Let $S \subseteq [m]$ have size $\mtilde$. We have
\begin{align*}
    \log\cN\inparen{B_p,\eta\inbraces{\vx\in\R^d\suchthat\polynorm{\vx}\le 1}} \le n\logv{\frac{4H^{1/p}\inparen{\Fstar}^{\max(1/p,1/2)}}{\eta}}
\end{align*}
\end{lemma}
\begin{proof}[Proof of \Cref{lemma:covering_easy}]
Define
\begin{align*}
    K &\coloneqq \inbraces{\vx\in\R^d\suchthat\polynorm{\vx}\le 1}.
\end{align*}
Let $C$ be a value such that for all $i \in S$ and $\vx\in\R^d$, we have $(C\rho_i)^{-1/p}\norm{\mA_{S_i}\vx}_{p_i} \le \gnorm{\mA\vx}{p}$. This means that $B_p \subseteq C^{1/p} K$. Then,
\begin{align*}
    \log\cN\inparen{B_p,\eta K} \le \log\cN\inparen{C^{1/p} K, \eta K} = \log\cN\inparen{K, \frac{\eta}{C^{1/p}} \cdot K} \le n\logv{\frac{4C^{1/p}}{\eta}}.
\end{align*}
By \Cref{lemma:containments}, when $p\ge 2$, we can choose $C = H\inparen{\Fstar}^{\max(1/p,1/2)}$.
This concludes the proof of \Cref{lemma:covering_easy}. 
\end{proof}

\section{Concentration analysis}
\label{sec:generic_chaining}

In this section, we prove Theorem \ref{thm:concentration}. Theorem \ref{thm:concentration} states our main result in its fullest generality. \Cref{thm:one_shot_lewis} follows easily from this, as we show in \Cref{sec:applications}.

We first state \Cref{thm:concentration}.

\begin{restatable}[General concentration result]{mainthm}{concentration}
\label{thm:concentration}
Let $\cG = (\mA \in \R^{n \times d}, S_1,\dots,S_m, p_1,\dots,p_m)$ where $S_1,\dots,S_m$ form a partition of $[n]$. Suppose at least one of the following holds:
\begin{itemize}
    \item $1 \le p < \infty$ and $p_1,\dots,p_m \ge 2$;
    \item $1/\log d \le p_1 = \dots = p_m = p < \infty$;
    \item $p_1 = \dots = p_m = 2$ and $1/\log n \le p < \infty$.
\end{itemize}
Let $P \coloneqq \max\inparen{1, \max_{i \in [m]} \min(p_i,\log\abs{S_i})}$. Suppose that $\vlambda \in \R^m$ is a probability measure over $[m]$, let $\mW$ be a rounding matrix (\Cref{defn:rounding_matrix}) such that we get an $\Fstar$-block Lewis overestimate (\Cref{defn:block_lewis_overestimate}), and define $\valpha$ according to \Cref{defn:alpha}. Let $\cD = \inparen{\rho_1,\dots,\rho_m}$ be a probability distribution over $[m]$ and $H \ge 1$ be such that $H\rho_i \ge \alpha_i^p/\norm{\valpha}_p^p$.

If
\begin{align*}
    \mtilde &= \Omega\inparen{\logv{\nfrac{1}{\delta}}\vbrho \cdot H \cdot P \inparen{\Fstar}^{\max(1,p/2)}},
\end{align*}
and if we sample $\cM \sim \cD^{\mtilde}$, then, with probability $\ge 1-\delta$, we have:
\begin{align*}
    \text{for all } \vx \in \R^d,\quad (1-\eps)\gnorm{\mA\vx}{p}^p \le \frac{1}{\mtilde}\sum_{i \in \cM} \frac{1}{\rho_i} \cdot \norm{\mA_{S_i}\vx}_{p_i}^p \le (1+\eps)\gnorm{\mA\vx}{p}^p
\end{align*}
\end{restatable}
The goal of the rest of this section is to prove \Cref{thm:concentration}. It may be helpful to recall the argument sketch given in \Cref{sec:sparsification_overview_concentration}. 

To formalize the idea given there, we first introduce the following notation (recall that $\rho_i$ is the probability that we choose group $i$ in a round of sampling and that \Cref{defn:entropy_numbers} defines the $e_N$).
\begin{align}
    g_i(\vx) &\coloneqq \norm{\mA_{S_i}\vx}_{p_i}\label{eq:chaining_g_defn} \\
    \dtwo(\vx, \xhat) &\coloneqq \inparen{\sum_{h=1}^{\mtilde} \inparen{\frac{g_{i_h}(\vx)^p}{\rho_{i_h}} - \frac{g_{i_h}(\xhat)^p}{\rho_{i_h}}}^2}^{1/2}\label{eq:chaining_dtwo_defn} \\
    \gamma_2(B_{p},\dtwo) &\coloneqq \inf_{\substack{\abs{T_N} \le 2^{2^N};\\ T_N  \subset B_{p}}} \sup_{\vx \in B_{p}} \sum_{N \ge 0} 2^{N/2} \cdot d_2(\vx,T_N)\label{eq:chaining_gamma2_defn}
\end{align}
The goal is to control $\gamma_2(B_p, \dtwo)$. This quantity represents the worst-case approximation error that one incurs by using the discretization scheme given by the $T_N$, where the discretization is taken with respect to the $\dtwo$ metric.

Towards this goal, we first apply a standard symmetrization reduction. Informally, this reduction (\Cref{lemma:symmetrization}) states that it is enough to analyze the average fluctuations of a Rademacher average of any set of $\mtilde$ (not necessarily distinct) reweighted groups.
\begin{lemma}[Symmetrization reduction]
\label{lemma:symmetrization} 
Let $R_1,\dots,R_{\mtilde}$ be independent Rademacher random variables (i.e., $\mathsf{Unif}\inparen{\pm 1}$). We have
\begin{align*}
    \exvv{\cG'}{\abs{\gnorm{\mA\vx}{p}^p - \norm{\vx}_{\cG'}^p}} \le 2 \underset{\cG'}{\E} \exvv{R_1, \ldots, R_{\mtilde}}{\frac{1}{\mtilde}\abs{\sum_{h=1}^{\mtilde} R_h \frac{g_{i'_h}(\vx)^p}{\rho_{i'_h}}}}.
\end{align*}
\end{lemma}
\begin{proof}[Proof of \Cref{lemma:symmetrization}]
We follow the proof of Lemma 3.1 due to \citet{lee22}. Let \(\tilde{\cG}\) be an independent copy of \(\cG'\). Fixing \(\cG'\), we have by Jensen's inequality that
\begin{align*}
    \abs{\gnorm{\mA\vx}{p}^p - \norm{\vx}_{\cG'}^p} &= \abs{\exvv{\tilde{\cG}}{\norm{\vx}_{\tilde{\cG}}^p} - \norm{\vx}_{\cG'}^p}
    \leq \exvv{\tilde{\cG}}{\abs{\norm{\vx}_{\tilde{\cG}}^p - \norm{\vx}_{\cG'}^p}}.
\end{align*}
Thus, taking expectation over \(\cG'\),
\begin{align*}
    \exvv{\cG'}{\abs{\gnorm{\mA\vx}{p}^p - \norm{\vx}_{\cG'}^p}} &\leq \exvv{\cG', \tilde{\cG}}{\abs{\norm{\vx}_{\tilde{\cG}}^p - \norm{\vx}_{\cG'}^p}}
    = \exvv{\cG', \tilde{\cG}}{\frac{1}{\mtilde}\abs{\sum_{h=1}^{\mtilde} \frac{g_{\tilde{i}_h}(\vx)^p}{\rho_{\widetilde{i}_h}} -\frac{g_{i'_h}(\vx)^p}{\rho_{i'_h}}}}.
\end{align*}
Observe that \(\frac{g_{\tilde{i}_h}(\vx)^p}{\rho_{\widetilde{i}_h}} -\frac{g_{i'_h}(\vx)^p}{\rho_{i'_h}}\) is symmetric,
and so is distributed the same as \(R_h \inparen{\frac{g_{\tilde{i}_h}(\vx)^p}{\rho_{\widetilde{i}_h}} -\frac{g_{i'_h}(\vx)^p}{\rho_{i'_h}}}\)
where \(R_h\) is an independent Rademacher variable. Then,
\begin{align*}
    \exvv{\cG', \tilde{\cG}}{\abs{\frac{1}{\mtilde}\sum_{h=1}^{\mtilde} \inparen{\frac{g_{\tilde{i}_h}(\vx)^p}{p_{\tilde{i}_h}} -\frac{g_{i'_h}(\vx)^p}{p_{i'_h}}}}} &= \underset{R_1, \ldots, R_{\mtilde}}{\E} \exvv{\cG', \tilde{\cG}}{\frac{1}{\mtilde}\abs{\sum_{h=1}^{\mtilde} R_h \inparen{\frac{g_{\tilde{i}_h}(\vx)^p}{\rho_{\widetilde{i}_h}} -\frac{g_{i'_h}(\vx)^p}{\rho_{i'_h}}}}} \\
    &\leq 2 \underset{\cG'}{\E} \exvv{R_1, \ldots, R_M}{\frac{1}{\mtilde}\abs{\sum_{h=1}^{\mtilde} R_h \frac{g_{i'_h}(\vx)^p}{\rho_{i'_h}}}}.
\end{align*}
This concludes the proof of \Cref{lemma:symmetrization}.
\end{proof}
With \Cref{lemma:symmetrization} in hand, we set up our chaining argument in \Cref{lemma:chaining_setup}. We first confirm that our random process is subgaussian with respect to our choice of $\dtwo$. 
\begin{lemma}[Choosing the distance]
\label{lemma:chaining_setup}
The random process $\sum_{h=1}^{\mtilde} R_{i_h}\cdot\nfrac{g_{i_h}(\vx)^p}{\rho_{i_h}}$  is subgaussian with respect to $\dtwo$ as defined in  \eqref{eq:chaining_dtwo_defn}.
\end{lemma}
\begin{proof}[Proof of \Cref{lemma:chaining_setup}]
Let
\begin{align*}
    P &\coloneqq \abs{\sum_{h=1}^{\mtilde} R_{h}\inparen{\frac{g_{i_h}(\vx)^p}{\rho_{i_h}} - \frac{g_{i_h}(\xhat)^p}{\rho_{i_h}}}}.
\end{align*}
Let us first calculate $\subgnorm{P}$. Using the fact that every term in this sum is independent and \Cref{fact:subgaussian_properties}, we get
\begin{align*}
    \subgnorm{P}^2 = \sum_{h=1}^{\mtilde} \subgnorm{R_{h}\inparen{\frac{g_{i_h}(\vx)^p}{\rho_{i_h}} - \frac{g_{i_h}(\xhat)^p}{\rho_{i_h}}}}^2 \le \sum_{h=1}^{\mtilde} 2\inparen{\frac{g_{i_h}(\vx)^p}{\rho_{i_h}} - \frac{g_{i_h}(\xhat)^p}{\rho_{i_h}}}^2
\end{align*}
By \Cref{fact:subgaussian_properties}, we have
\begin{align*}
    \prvv{R_{h}}{P \ge v\dtwo(\vx,\xhat)} &\le 2\expv{-\frac{v^2\dtwo(\vx,\xhat)^2}{4\sum_{h=1}^{\mtilde} \inparen{\frac{g_{i_h}(\vx)^p}{\rho_{i_h}} - \frac{g_{i_h}(\xhat)^p}{\rho_{i_h}}}^2}} = 2\expv{-\frac{v^2}{2}}.
\end{align*}
This concludes the proof of \Cref{lemma:chaining_setup}.
\end{proof}
\Cref{lemma:chaining_setup} tells us that $\dtwo$ is a choice of distance on $B_{p}$ that allows us to use the subgaussian form of chaining to analyze our random process. Along with the way we have set up our sampling process, we have enough to apply Theorem \ref{thm:generic_chaining_dirksen}. This is simply a restatement of Lemma 2.6 of \cite{jlls23} for our setting.

\begin{theorem}[Restatement of Lemma 2.6 from \cite{jlls23}, $\alpha=2$]
\label{thm:generic_chaining_dirksen}
Recall $\dtwo$ (\ref{eq:chaining_dtwo_defn}) and $\gamma_2$ (\ref{eq:chaining_gamma2_defn}). Suppose that for some $D$ and for every choice of $i_1,\dots,i_{\mtilde}$, we have
\begin{align*}
    \gamma_2\inparen{B_p,\frac{\dtwo}{\mtilde}} \lesssim D\inparen{\max_{\vx\in B_p} \norm{\vx}_{\cG'}^p}^{1/2},
\end{align*}
where
\begin{align*}
    \norm{\vx}_{\cG'}^p = \frac{1}{\mtilde}\sum_{h=1}^{\mtilde} \frac{g_{i_h}(\vx)^p}{\rho_{i_h}}.
\end{align*}
Then, we have the following.
\begin{align*}
    \exvv{\cD}{\sup_{\vx\in B_p} \abs{\gnorm{\mA\vx}{p}-\norm{\vx}_{\cG'}^p}} &\lesssim D.
\end{align*}
If we also have for all choices $i_1,\dots,i_{\mtilde}$ and for some $\widehat{D}$ that
\begin{align*}
    \mathsf{diam}\inparen{B_p, \frac{\dtwo}{\mtilde}} \lesssim \widehat{D}\inparen{\max_{\vx\in B_p} \norm{\vx}_{\cG'}^p}^{1/2},
\end{align*}
then there exists a universal constant $C > 0$ such that for all $0 \le t \le \nfrac{1}{2K\widehat{D}}$,
\begin{align*}
    \prvv{\cD}{\sup_{\vx\in B_p} \abs{\gnorm{\mA\vx}{p}^p-\norm{\vx}_{\cG'}^p} \ge C(D + t\widehat{D})} \le \expv{-\frac{Kt^2}{4}}.
\end{align*}
\end{theorem}

With Theorem \ref{thm:generic_chaining_dirksen} in our arsenal, our task becomes to compute $\gamma_2(B_p,\dtwo)$ (which we will then divide by $\mtilde$ so that we can apply \Cref{thm:generic_chaining_dirksen}). Hereafter, we will simply abbreviate $\gamma_2(B_p,\dtwo)$ as $\gamma_2$. We will first weaken the definition of $\gamma_2$, which is essentially equivalent to Dudley's integral. Recall the definition of the entropy numbers $e_N$ (\Cref{defn:entropy_numbers}) and notice that
\begin{align*}
    \gamma_2 \le \sum_{N \ge 0} 2^{N/2}e_N(B_{p},\dtwo). 
\end{align*}
We now rewrite $\dtwo$ in a form that will be more convenient for us.

\begin{lemma}
\label{lemma:concentration_g2_rewrite_d}
We have
\begin{align*}
    \dtwo(\vx, \xhat) \le \max(p,2) \inparen{\Fstar}^{\max(0,p/4-1/2)}{\mtilde}^{1/2} \cdot \polynorm{\vx-\xhat}^{\min(p/2,1)} \cdot \inparen{\max_{\vx \in B_{p}}\norm{\vx}_{\cG'}^p}^{1/2}
\end{align*}
and therefore
\begin{align*}
    \gamma_2 &\le {\mtilde}^{1/2}\max(p,2) \inparen{H^{1/p}\inparen{\Fstar}^{1/2}}^{\max(0,p/2-1)}\inparen{\max_{\vx \in B_{p}}\norm{\vx}_{\cG'}^p}^{1/2}\sum_{N \ge 0} 2^{N/2}e_N\inparen{B_p, \polynorm{\cdot}}^{\min(p/2,1)}.
\end{align*}
\end{lemma}
\begin{proof}[Proof of \Cref{lemma:concentration_g2_rewrite_d}]
We have two cases. We first address the case $0 < p < 2$. Recall that in this regime, we have $\abs{a^{p/2}-b^{p/2}} \le \abs{a-b}^{p/2}$. Since $g_{i_h}(\vx)$ is a norm, the triangle inequality tells us that $\abs{g_{i_h}(\vx) - g_{i_h}(\xhat)} \le g_{i_h}(\vx-\xhat)$. We use these and write
\begin{align*}
    \dtwo(\vx, \xhat) &= \inparen{\sum_{h=1}^{\mtilde} \inparen{\frac{g_{i_h}(\vx)^{p/2}}{\sqrt{\rho_{i_h}}} - \frac{g_{i_h}(\xhat)^{p/2}}{\sqrt{\rho_{i_h}}}}^2\inparen{\frac{g_{i_h}(\vx)^{p/2}}{\sqrt{\rho_{i_h}}} + \frac{g_{i_h}(\xhat)^{p/2}}{\sqrt{\rho_{i_h}}}}^2}^{1/2} \\
    &\le \inparen{\sum_{h=1}^{\mtilde} \inparen{\frac{g_{i_h}(\vx-\xhat)^{p/2}}{\sqrt{\rho_{i_h}}}}^2\inparen{\frac{g_{i_h}(\vx)^{p/2}}{\sqrt{\rho_{i_h}}} + \frac{g_{i_h}(\xhat)^{p/2}}{\sqrt{\rho_{i_h}}}}^2}^{1/2} \\
    &\le \polynorm{\vx-\xhat}^{p/2}\cdot\inparen{\sum_{h=1}^{\mtilde} \inparen{\frac{g_{i_h}(\vx)^{p/2}}{\sqrt{\rho_{i_h}}} + \frac{g_{i_h}(\xhat)^{p/2}}{\sqrt{\rho_{i_h}}}}^2}^{1/2} \\
    &\le 2{\mtilde}^{1/2} \cdot \polynorm{\vx-\xhat}^{p/2}\inparen{\max_{\vx \in B_{p}}\norm{\vx}_{\cG'}^p}^{1/2}
\end{align*}
which concludes the proof in the range $0 < p < 2$.

We now move onto the case where $p \ge 2$. Recall that by Lipschitzness, we get
\begin{align*}
    \abs{a^{p/2}-b^{p/2}} \le \frac{p}{2} \cdot \max(a,b)^{p/2-1}\abs{a-b}.
\end{align*}
Next, by \Cref{lemma:containments}, we know that for all $i$, 
\begin{align*}
    \norm{\mA_{S_i}\vx}_{p_i} \le \inparen{\Fstar}^{1/2}\rho_i^{1/p}H^{1/p}\gnorm{\mA\vx}{p} \le \inparen{\Fstar}^{1/2}\rho_i^{1/p}H^{1/p}.
\end{align*}
We use these to rewrite $\dtwo(\vx,\xhat)$.
\begin{align*}
    \dtwo(\vx, \xhat) &= \inparen{\sum_{h=1}^{\mtilde} \inparen{\frac{g_{i_h}(\vx)^{p/2}}{\sqrt{\rho_{i_h}}} - \frac{g_{i_h}(\xhat)^{p/2}}{\sqrt{\rho_{i_h}}}}^2\inparen{\frac{g_{i_h}(\vx)^{p/2}}{\sqrt{\rho_{i_h}}} + \frac{g_{i_h}(\xhat)^{p/2}}{\sqrt{\rho_{i_h}}}}^2}^{1/2} \\
    &\le \frac{p}{2}\cdot H^{1/2-1/p}\inparen{\Fstar}^{p/4-1/2}\inparen{\sum_{h=1}^{\mtilde} \inparen{\rho_{i_h}^{-1/p}g_{i_h}(\vx-\xhat)}^2\inparen{\frac{g_{i_h}(\vx)^{p/2}}{\sqrt{\rho_{i_h}}} + \frac{g_{i_h}(\xhat)^{p/2}}{\sqrt{\rho_{i_h}}}}^2}^{1/2} \\
    &\le \frac{p}{2}\cdot H^{1/2-1/p}\inparen{\Fstar}^{p/4-1/2}\polynorm{\vx-\xhat}\cdot\inparen{\sum_{h=1}^{\mtilde} \inparen{\frac{g_{i_h}(\vx)^{p/2}}{\sqrt{\rho_{i_h}}} + \frac{g_{i_h}(\xhat)^{p/2}}{\sqrt{\rho_{i_h}}}}^2}^{1/2} \\
    &\le p \cdot H^{1/2-1/p}\inparen{\Fstar}^{p/4-1/2}{\mtilde}^{1/2} \cdot \polynorm{\vx-\xhat}\inparen{\max_{\vx \in B_{p}}\norm{\vx}_{\cG'}^p}^{1/2}
\end{align*}
We conclude the proof of \Cref{lemma:concentration_g2_rewrite_d}.
\end{proof}

In light of \Cref{lemma:concentration_g2_rewrite_d}, observe that it is enough to calculate each term of the sum $\sum_{N \ge 0} 2^{N/2}e_N(B_p, \polynorm{\cdot})^{\min(p/2,1)}$. We begin this analysis with \Cref{lemma:concentration_g2_two_one_volume}.

\begin{lemma}
\label{lemma:concentration_g2_two_one_volume}
For all $N \ge 0$, we have
\begin{align*}
    2^{N/2}e_N\inparen{B_p, \polynorm{\cdot}}^{\min(p/2,1)} \le 2^{N/2 + \min(p/2,1)(2+(1/p)\log C - 2^N/d)},
\end{align*}
where $C$ is such that for all $i \in S$ and $\vx\in\R^d$, we have $\inparen{C\rho_i}^{-1/p}\norm{\mA_{S_i}\vx}_{p_i} \le \gnorm{\mA\vx}{p}$.
\end{lemma}
\begin{proof}[Proof of \Cref{lemma:concentration_g2_two_one_volume}]
Recall that by \Cref{lemma:covering_easy}, we have
\begin{align*}
    \log\cN\inparen{B_p,\eta\inbraces{\vx\in\R^d\suchthat\polynorm{\vx}\le 1}} \le d\logv{\frac{4C^{1/p}}{\eta}}.
\end{align*}
We set $\eta = 2^{2+(1/p)\log C - 2^N/d}$ so that
\begin{align*}
    \log\cN\inparen{B_p,\eta\inbraces{\vx\in\R^d\suchthat\polynorm{\vx}\le 1}} \le 2^N.
\end{align*}
Then,
\begin{align*}
    2^{N/2}e_N\inparen{B_p, \polynorm{\cdot}}^{\min(p/2,1)} &\le 2^{N/2 + \min(p/2,1)(2+(1/p)\log C - 2^N/d)},
\end{align*}
concluding the proof of Lemma \ref{lemma:concentration_g2_two_one_volume}.
\end{proof}

Using \Cref{lemma:concentration_g2_two_one_volume}, we get a rapidly converging tail in our summation for large values of $N$. See \Cref{lemma:g2_tail}.

\begin{lemma}
\label{lemma:g2_tail}
Let $N_S \coloneqq \logv{6d\log d}$. We have
\begin{align*}
    \sum_{N \ge 0} 2^{N/2}e_N\inparen{B_p, \polynorm{\cdot}}^{\min(p/2,1)} &\lesssim \sqrt{nH^{1/\max(2,p)}\inparen{\Fstar}^{1/2}\log d}\\
    &\quad+\sum_{N \le N_S} 2^{N/2}e_N\inparen{B_p, \polynorm{\cdot}}^{\min(p/2,1)}.
\end{align*}
\end{lemma}
\begin{proof}[Proof of \Cref{lemma:g2_tail}]
For now, let $C$ be such that for all $i \in S$ and $\vx\in\R^d$, we have $\inparen{C\rho_i}^{-1/p}\norm{\mA_{S_i}\vx}_{p_i} \le \gnorm{\mA\vx}{p}$.

Let $N_S$ be a threshold such that for all $N \ge N_S$, we use the entropy number bound given by \Cref{lemma:concentration_g2_two_one_volume} (as the volume-based covering number bound is much better for small values of $e_N$). Let us enforce the constraint $N_S \ge \ceil{\logv{\nfrac{3d}{\min(p,2)}}}$, so that the entropy number bound in \Cref{lemma:concentration_g2_two_one_volume} is decreasing in $N$ and is dominated from above by a geometric series with common ratio $\nfrac{1}{2}$.

We now set $N_S = \logv{6d\log d}$. Since $p \ge \nfrac{1}{\log d}$, we know that $N_S \ge \ceil{\logv{\nfrac{3d}{\min(p,2)}}}$. Now, since $2^{N/2}e_N$ is bounded above by a geometric series with common ratio $\nfrac{1}{2}$, the summation for all $N \ge N_S$ is dominated by the first term. Let us evaluate this. We first observe that
\begin{align*}
    \frac{N_S}{2} + \min\inparen{\frac{p}{2},1}\inparen{-\frac{2^{N_S}}{d} + 2 + \frac{\log C}{p}} &= \frac{\logv{6d\log d}}{2} + \min\inparen{\frac{p}{2},1}\inparen{-\frac{2^{\logv{6d\log d}}}{d} + 2 +\frac{\log C}{p}} \\
    &= \frac{\logv{6d\log d}}{2} + \min\inparen{\frac{p}{2},1}\inparen{-6\log d + 2 +\frac{\log C}{p}} \\
    &\le \frac{\logv{6d\log d}}{2} + \min\inparen{\frac{p}{2},1}\inparen{-6\log d + 2 +\frac{\log C}{p}} \\
    &\le \frac{\logv{6d\log d}}{2} + \frac{\log C}{\max(2,p)} \le \frac{6Cd\log d}{2}
\end{align*}
By \Cref{lemma:concentration_g2_two_one_volume}, we see that
\begin{align*}
    2^{N_S/2}e_{N_S}\inparen{B_p, \polynorm{\cdot}}^{\min(p/2,1)} \lesssim \sqrt{dC^{1/\max(2,p)}\log d},
\end{align*}
and by \Cref{lemma:containments}, we can choose
\begin{align*}
    C = H\inparen{\Fstar}^{\max(1,p/2)}.
\end{align*}
We plug this in, account for the remaining terms in the summation, and conclude the proof of \Cref{lemma:g2_tail}.
\end{proof}

We now give another way to evaluate the terms of our summation when the indices $N$ are such that the entropy numbers are rather large. See \Cref{lemma:concentration_g2_two_one_sudakov}.

\begin{lemma}
\label{lemma:concentration_g2_two_one_sudakov}
For all $N \ge 0$, we have
\begin{align*}
    2^{N/2}e_N\inparen{B_p, \polynorm{\cdot}}^{\min(p/2,1)} \le \inparen{C(p) \cdot \max_i\min\inbraces{p_i,\log\abs{S_i}}H^{2/\max(2,p)}\Fstar\logv{\Fstar+\mtilde}}^{1/2},
\end{align*}
where $C(p)$ is a constant that only depends on $p$.
\end{lemma}
\begin{proof}[Proof of \Cref{lemma:concentration_g2_two_one_sudakov}]
For now, let
\begin{align*}
    f(\Fstar,\cG) \coloneqq C(p) \cdot \max_i\min\inbraces{p_i,\log\abs{S_i}}H^{2/\max(2,p)}\Fstar\logv{\Fstar+\mtilde}.
\end{align*}
By \Cref{thm:general_covering}, we have
\begin{align*}
    \log \cN\inparen{B_{p}, \eta \cdot \inbraces{\vy \in \R^d \suchthat \polynorm{\vx} \le 1}} \lesssim \eta^{-\min(2,p)} \cdot f(\Fstar,\cG).
\end{align*}
so if we choose, for some universal constant $C_0$,
\begin{align*}
    \eta = C_0 \cdot 2^{-N/\min(p,2)}\inparen{f(\Fstar,\cG)}^{\max(1/2,1/p)},
\end{align*}
then we get
\begin{align*}
    \log \cN\inparen{B_{p}, \eta \cdot \inbraces{\vy \in \R^d \suchthat \polynorm{\vx} \le 1}} \le 2^N.
\end{align*}
Thus, $e_N\inparen{B_p, \polynorm{\cdot}} \le \eta$. Exponentiating and substituting the definition of $f(\Fstar,\cG)$ concludes the proof of \Cref{lemma:concentration_g2_two_one_sudakov}.
\end{proof}

We now show how to complete the sum by combining \Cref{lemma:concentration_g2_two_one_sudakov} and \Cref{lemma:g2_tail}.

\begin{lemma}
\label{lemma:concentration_subgaussian}
We have
\begin{align*}
    \sum_{N \ge 0} 2^{N/2}e_N\inparen{B_p, \polynorm{\cdot}}^{\min(p/2,1)} &\lesssim \sqrt{nH^{1/\max(2,p)}\inparen{\Fstar}^{1/2}\log n} \\
    &\quad +\log n\inparen{C(p) \cdot \max_i\min\inbraces{p_i,\log\abs{S_i}}H\Fstar\logv{\Fstar+\mtilde}}^{1/2}.
\end{align*}
\end{lemma}
\begin{proof}[Proof of \Cref{lemma:concentration_subgaussian}]
Noting that $N_S = \logv{6d\log d} \lesssim \log d$, we combine the conclusions of \Cref{lemma:concentration_g2_two_one_sudakov} and \Cref{lemma:g2_tail} to obtain the statement of \Cref{lemma:concentration_subgaussian}.
\end{proof}

Finally, we translate Lemma \ref{lemma:concentration_subgaussian} into an upper bound on the process we started with by using \Cref{thm:generic_chaining_dirksen}. We then use this to complete the proof of Theorem \ref{thm:concentration}.

\begin{proof}[Proof of Theorem \ref{thm:concentration}]
As we have done in previous proofs, as a shorthand, we define
\begin{align*}
    p^{\star} \coloneqq \max_{i}\min\inbraces{p_i,\log\abs{S_i}}.
\end{align*}
We first weaken the statement of \Cref{lemma:concentration_subgaussian} to read
\begin{align*}
    &\quad\sum_{N \ge 0} 2^{N/2}e_N\inparen{B_p, \polynorm{\cdot}}^{\min(p/2,1)} \\
    &\lesssim \max\inbraces{d,\Fstar}H^{1/2\max(2,p)}\sqrt{\log d} \\
    &\quad+ \log d \inparen{C(p) \cdot H^{2/\max(2,p)}p^{\star}\max\inbraces{d,\Fstar}\logv{d + \Fstar + \mtilde}}^{1/2} \\
    &\lesssim \log d \inparen{C(p) \cdot H^{2/\max(2,p)}p^{\star}\max\inbraces{d,\Fstar}\logv{d + \Fstar + \mtilde}}^{1/2}.
\end{align*}
Let $V$ denote the right hand side of the above. Combining this rewrite with \Cref{lemma:concentration_g2_rewrite_d}, we get
\begin{align*}
    \gamma_2 \lesssim {\mtilde}^{1/2}\max(p,2) \inparen{H^{1/p}\inparen{\Fstar}^{1/2}}^{\max(0,p/2-1)}V\inparen{\max_{\vx \in B_{p}}\norm{\vx}_{\cG'}^p}^{1/2}.
\end{align*}
By the symmetrization reduction (\Cref{lemma:symmetrization}) and \Cref{thm:generic_chaining_dirksen}, we have
\begin{align*}
    \exv{\sup_{\vx\in B_p} \abs{\gnorm{\mA\vx}{p}-\norm{\vx}_{\cG'}^p}} \lesssim \frac{\max(p,2) \inparen{H^{1/p}\norm{\valpha}_p}^{\max(0,p/2-1)}V}{\mtilde^{1/2}}
\end{align*}
and so to make the RHS upper-bounded by $\eps$, it is sufficient to set $\mtilde$ according to
\begin{align*}
    \mtilde &\asymp \frac{\max(p,2)^2 \inparen{H^{1/p}\inparen{\Fstar}^{1/2}}^{\max(0,p-2)}V^2}{\eps^2}.
\end{align*}
This means that when $p < 2$, we have
\begin{align*}
    \mtilde &\asymp \frac{(\log d)^2\cdot C(p) \cdot Hp^{\star}\max\inbraces{d,\Fstar}\logv{d + \Fstar + \mtilde}}{\eps^2} \\
    &\asymp \frac{(\log d)^2\logv{\max\inbraces{d,\Fstar}/\eps}\cdot C(p) \cdot Hp^{\star}\max\inbraces{d,\Fstar}}{\eps^2},
\end{align*}
which is what we desired.

For $p \ge 2$, we have
\begin{align*}
    \mtilde &\asymp \frac{p^2H^{1-2/p}\inparen{\Fstar}^{p/2-1}\inparen{(\log d)^2H^{2/p}p^{\star}\max\inbraces{d,\Fstar}\logv{d+\Fstar+\mtilde}}}{\eps^2} \\
    &\asymp \frac{(\log d)^2\logv{\max\inbraces{d,\Fstar}/\eps}\cdot p^2 \cdot Hp^{\star}\max\inbraces{d,\Fstar}^{p/2}}{\eps^2}
\end{align*}
To bound $\mathsf{diam}(B_p,\dtwo)$, by the triangle inequality, it is enough to estimate $\dtwo(\vx,0)$ for all $\vx \in B_p$. Recalling \Cref{lemma:containments}, we have
\begin{align*}
    \dtwo(\vx,0) &= \inparen{\sum_{h=1}^{\mtilde} \inparen{\frac{g_{i_h}(\vx)^p}{\rho_{i_h}} - \frac{g_{i_h}(0)^p}{\rho_{i_h}}}^2}^{1/2} = \inparen{\sum_{h=1}^{\mtilde} \inparen{\frac{g_{i_h}(\vx)^p}{\rho_{i_h}}}\cdot\inparen{\frac{g_{i_h}(\vx)^p}{\rho_{i_h}}}}^{1/2} \\
    &\le \inparen{\sum_{h=1}^{\mtilde} H\inparen{\Fstar}^{\max(1,p/2)}\cdot\inparen{\frac{g_{i_h}(\vx)^p}{\rho_{i_h}}}}^{1/2} \le \mtilde^{1/2}H^{1/2}\inparen{\Fstar}^{\max(1/2,p/4)}\max_{\vx\in B_p}\inparen{\norm{\vx}_{\cG'}^{p}}^{1/2},
\end{align*}
which means that we may set $\widehat{D}$ in \Cref{thm:generic_chaining_dirksen} according to
\begin{align*}
    \widehat{D} \asymp \frac{H^{1/2}\inparen{\Fstar}^{\max(1/2,p/4)}\max_{\vx\in B_p}\inparen{\norm{\vx}_{\cG'}^{p}}^{1/2}}{\mtilde^{1/2}}.
\end{align*}
We now verify that if we choose
\begin{align*}
    \mtilde \asymp \vbrho \cdot H\max_{i \in [m]} \min(p_i,\log\abs{S_i}) \inparen{\Fstar}^{\max(1,p/2)}\logv{\nfrac{1}{\delta}},
\end{align*}
that we indeed get for some universal constant $C$ that
\begin{align*}
    \prvv{\cD}{\max_{\vx\in B_p} \abs{\gnorm{\mA\vx}{p}^p-\norm{\vx}_{\cG'}^p} \ge C\eps} \lesssim \delta.
\end{align*}
We rescale $\eps$ appropriately and conclude the proof of \Cref{thm:concentration}.
\end{proof}

\section{Applications and algorithms}
\label{sec:applications}
At this point in the chapter, we are ready to prove our main results (\Cref{thm:one_shot_lewis} and \Cref{thm:computeblw}).

\subsection{Block norm approximations via block Lewis weights (Proof of Theorem \ref{thm:one_shot_lewis})}
\label{sec:lewis}

We restate and prove the main result of the chapter.
\oneshotlewis*
\begin{proof}[Proof of Theorem \ref{thm:one_shot_lewis}]
Observe that \Cref{lemma:block_instantiation_small} proves the existence of a probability measure $\vlambda$ over $[m]$ and a rounding $\mW$ that are $\Fstar$-block Lewis overestimates for $\Fstar = d$ if $p_i \ge 2$, and \Cref{lemma:block_instantiation_alleq} proves the existence of a probability measure $\widehat{\vlambda}$ over $[n]$ and corresponding $\widehat{\valpha} \in \R^n$ such that we get an $\Fstar = d$-Lewis overestimate.

We now apply \Cref{thm:concentration} and conclude the proof of Theorem \ref{thm:one_shot_lewis}.
\end{proof}

\subsection{Efficient computation of block Lewis weight overestimates (Proof of \texorpdfstring{\Cref{thm:computeblw}}{Theorem 2})}
\label{sec:alg}

In this subsection, we restate and prove \Cref{thm:computeblw}.
\computeblw*
We break up the proof into two sections -- one where $p_1=\dots=p_m=2$ and $p > 0$ and another where $p=2$ and $p_1,\dots,p_m \ge 2$.

\subsubsection{Special case -- \texorpdfstring{$p > 0, p_1=\dots=p_m=2$}{p>0,p1=...=pm=2}}

Following \Cref{defn:block_lewis_overestimate} and the discussion in \Cref{sec:blw_defs_covering}, observe that when $p_1=\dots=p_m=2$, we have \[\beta_i(\mV) \coloneqq \inparen{\sum_{j \in S_i} \va_j^{\top}(\mA^{\top}\mV\mA)^{-1}\va_j}^{1/2}.\]
As before, we call \(\beta_i(\mV)^p\) block Lewis weights.
Given \(\vb \in \R^m\), we let \(\mB \in \R^{n}\) be a diagonal matrix given by
\(\mB_{jj} = \vb_i\) for all \(i \in [m]\), \(j \in S_i\).
First, let us specialize the definition of block Lewis weight overestimates (recall \Cref{defn:block_lewis_overestimate}) to this special case.
\begin{definition}\label{def:lev_score_over}
    For \(\nu \ge 0\), we say \(\vb \in \R_{\geq 0}^m\) is a vector of \(\nu\)-bounded block Lewis weight overestimates for \(\mA\) if
    \[\|\vb\|_1 \le \nu,\]
    and for all \(i \in [m]\),
    \[b_i \ge \beta_i(\mB^{1-2/p})^p.\]
\end{definition}

We can think of the definition of block Lewis weight overestimates as being a relaxation of the fixed point condition for block Lewis weights that is described in \cite[Page 31, Proof of Lemma 4.2]{jlls23}. 

As a primitive, our algorithm will use \textit{leverage score overestimates} (see \cite[Definition 2.2]{jls22}). They are approximate forms of leverage scores \(\tau_j(\mM)\).
\begin{definition}
    For \(\nu \ge 0\), we say \(\widetilde{\tau} \in \R^m_{\ge 0}\) is a vector of \(\nu\)-bounded leverage score overestimates for \(\mM \in \R^{m \times d}\) if
    \[\|\widetilde{\tau}\|_1 \le \nu\]
    and for all \(i \in [k]\),
    \[\widetilde{\tau}_i \ge \tau_i(\mM).\]
\end{definition}

There are known efficient algorithms for computing leverage score overestimates or reducing leverage score computations to linear system solves.

\begin{theorem}[\protect{\cite[Theorem 3]{jls22}}]\label{thm:lev_ove}
There is an algorithm \(\textsc{OverLev}\) that, given \(\mM \in \R^{n \times d}\), produces \(O(n)\)-bounded leverage score overestimates for \(\mM\) in \(\widetilde{O}(\nnz{\mM} + d^\omega)\) time, where \(\omega\) is the matrix multiplication exponent.
\end{theorem}

We will use two different algorithms depending on the value of $p$. If $p \le 2$, then we present a contractive scheme reminiscent of the algorithm of \citet{cp15}. If $p > 2$, we present an algorithm similar to those of \citet{ccly19} and \citet{jls21}.

We begin with the case where $0 < p \le 2$ (in fact, we will see that this algorithm yields guarantees where $p < 4$). The main objects of interest here are \Cref{alg:blw_lewismap} and \Cref{lemma:contract_alg}.

\begin{algorithm}
\caption{Algorithm to compute block Lewis weight overestimates, $0 < p < 4$ and $p_1=\dots=p_m=2$.}
\begin{algorithmic}[1]\label{alg:blw_lewismap}
\State \textbf{Input:} \(\mA \in \R^{n \times d}\), outer norm $0 < p < 4$, group structure \((S_1,\dots,S_m,2,\dots,2)\).
\State Initialize \(\vb^{(0)} = \frac{d}{m} \cdot \onev_m\).
\State Define $\psi$ such that\label{line:psi_def}
\begin{align*}
    \psi_i(\vb) \coloneqq \inparen{\inparen{\sum_{j \in S_i} \inparen{\va_j^{\top}\inparen{\mA^{\top}\mB^{1-2/p}\mA}^{-1}\va_j}^{p_i/2}}^{2/p_i}}^{p/2}.
\end{align*}
\For{\(t = 1, \ldots, T\)}
    \State \(\vb^{(t)} = \max(\psi(\vb^{(t-1)}),\onev_{m} \cdot 1/m)\) \Comment{The $\max$ is taken elementwise here.}
\EndFor
\State \Return \(1.1\vb^{(T)}\)
\end{algorithmic}
\end{algorithm}

\begin{lemma}
\label{lemma:contract_alg}
Let
\begin{align*}
    T \ge \frac{\ln\inparen{\frac{\ln\inparen{\frac{m}{d}}}{\ln\inparen{1+\eps}}}}{\ln\inparen{\abs{\frac{2}{2-p}}}}.
\end{align*}
Then, the weights $1.1\vb^{(T)}$ output by \Cref{alg:blw_lewismap}  are a $1.1n$-block Lewis overestimate (\cref{defn:block_lewis_overestimate}). Furthermore, computing $\vb^{(T)}$ requires at most $T$ computations of the vector whose entries are the $\va_j^{\top}\inparen{\mA^{\top}\mD\mA}^{-1}\va_j$ for all $j$, where $\mD$ is a diagonal matrix.
\end{lemma}

To prove \Cref{lemma:contract_alg}, we first state and prove \Cref{lemma:lewis_contract_easy}.

\begin{lemma}
\label{lemma:lewis_contract_easy}
Let $\psi$ be as defined in Line \ref{line:psi_def} in \Cref{alg:blw_lewismap}. For all $\vu \in \R^{m}_{\ge 0}$ and $\vv \in \R^{m}_{\ge 0}$, we have
\begin{align*}
    \max_{1 \le i \le m} \abs{\lnv{\frac{\psi_i(\vu)}{\psi_i(\vv)}}} \le \abs{\frac{p}{2}-1}\max_{1 \le i \le m} \abs{\lnv{\frac{u_i}{v_i}}},
\end{align*}
and therefore $\psi(\vu)$ is a contraction whenever $0 < p < 4$.
\end{lemma}
\begin{proof}[Proof of \Cref{lemma:lewis_contract_easy}]
Fix some index $i \le m$. For notational simplicity in this proof, let $\alpha$ be such that $\lnv{\alpha} \coloneqq \triangle(\vu,\vv)$.

This easily implies that
\begin{align*}
    \frac{1}{\alpha^{\abs{1-2/p}}} \cdot \va_j^{\top}\inparen{\mA^{\top}\mV^{1-2/p}\mA}^{-1}\va_j \le \va_j^{\top}\inparen{\mA^{\top}\mU^{1-2/p}\mA}^{-1}\va_j \le \alpha^{\abs{1-2/p}} \cdot \va_j^{\top}\inparen{\mA^{\top}\mV^{1-2/p}\mA}^{-1}\va_j.
\end{align*}
We take the $p_i/2$-norm and take the $p/2$ power, which tells us that for all $1 \le i \le m$,
\begin{align*}
    \frac{1}{\alpha^{\abs{p/2-1}}} \cdot \psi_i(\vv) \le \psi_i(\vu) \le \alpha^{\abs{p/2-1}} \cdot \psi_i(\vv).
\end{align*}
Hence,
\begin{align*}
    \max_{1 \le i \le m} \abs{\lnv{\frac{\psi_i(\vu)}{\psi_i(\vv)}}} \le \abs{\frac{p}{2}-1}\lnv{\alpha} = \abs{\frac{p}{2}-1}\max_{1 \le i \le m} \abs{\lnv{\frac{u_i}{v_i}}},
\end{align*}
completing the proof of \Cref{lemma:lewis_contract_easy}.
\end{proof}

We are now ready to complete the proof of \Cref{lemma:contract_alg}.

\begin{proof}[Proof of \Cref{lemma:contract_alg}]
The computational complexity guarantee is immediate, so we focus on the approximation guarantee.

By \Cref{lemma:lewis_contract_easy} and the Banach fixed point theorem, we know that $\psi$ has a unique fixed point. Denote this fixed point by $\vb^{\star}$. We would like to argue that since the convergence in the $\ln$-metric is linear, it takes roughly $\log\log(1+\eps)/\log\abs{2/(2-p)}$ applications of $\psi$ to reach a $(1+\eps)$-multiplicative approximation to $\vb^{\star}$. An annoying technicality is that if $\vb^{\star}$ has some elements arbitrarily close to $0$, then the convergence rate could be very slow. To fix this, we simply enforce that the coordinates of the iterates never drop below $1/m$. It is easy to see that this only overestimates the true weights and therefore does not affect our sampling guarantees (in particular, $\norm{\max(\vb^{\star}, \onev_m \cdot 1/m)}_1 \le d + 1$).

More precisely, let $\vb \coloneqq \max(\vb^{\star}, \onev_m \cdot 1/m)$. Notice that for all $1 \le i \le m$, $\abs{\ln\inparen{\nfrac{b_i^{(0)}}{b_i}}} \le \ln\inparen{\nfrac{m}{d}}$. This means that after $T$ iterations, we have
\begin{align*}
    \max_{1 \le i \le m} \abs{\ln\inparen{\frac{b^{(T)}_i}{b_i}}} \le \abs{\frac{p}{2}-1}^T\max_{1 \le i \le m} \abs{\ln\inparen{\frac{b_i^{(0)}}{b_i}}} \le  \abs{\frac{p}{2}-1}^T\ln\inparen{\frac{m}{d}}.
\end{align*}
Choosing
\begin{align*}
    T \ge \frac{\ln\inparen{\frac{\ln\inparen{\frac{m}{d}}}{\ln\inparen{1+\eps}}}}{\ln\inparen{\abs{\frac{2}{2-p}}}}
\end{align*}
and observing that for sampling that it is sufficient to choose $\eps = 0.1$ implies that $\vb^{(T)}$ is an entrywise $1.1$-approximation to $\vb$. As this is sufficient to get the concentration in the setting of \Cref{thm:concentration}, we may  complete the proof of \Cref{lemma:contract_alg}.
\end{proof}

Now, we move onto the case where $p \ge 2$. This covers the cases of $p$ where \Cref{alg:blw_lewismap} is not a contraction (whenever $p \ge 4$). In this setting, we have \Cref{alg:blw}. At a high level, observe that Line \ref{line:alg_iter_next} of \Cref{alg:blw} performs a fixed point iteration on the stationary condition \(b_i = \beta_i(\mB^{1-2/p})^p\) that holds for the optimal choice of block Lewis weights (see \cite[proof of Lemma 4.2 and (4.5)]{jlls23}).

\begin{algorithm}[H]
\caption{Algorithm to compute block Lewis weight overestimates, $p \ge 2$}
\begin{algorithmic}[1]\label{alg:blw}
\State \textbf{Input:} \(\mA \in \R^{n \times d}\)
\State \textbf{Output:}
\State Initialize \(\vb^{(1)} = \frac{d}{m} \cdot \onev\) \label{line:alg_init}
\For{\(t = 1, \ldots, T - 1\)}
    \State \(\mB^{(t)}_{jj} = b^{(t)}_i\) for all \(i \in [m]\), \(j \in S_i\)
    \State \(\widetilde{\tau}^{(t)} = \textsc{OverLev}((\mB^{(t)})^{1/2 - 1/p} \mA)\)\label{line:alg_iter_lev}
    \State \(b^{(t+1)}_i = \sum_{j \in S_i} \widetilde{\tau_j}\) for all \(i \in [m]\)\label{line:alg_iter_next}
\EndFor
\State \(\overline{\vb} = \frac{1}{T} \sum_{t=1}^T \vb^{(t)}\)
\State \Return \(\vb = \frac{3}{2} \overline{\vb}\)
\end{algorithmic}
\end{algorithm}

The guarantee we obtain for \Cref{alg:blw} is captured by \Cref{thm:blw}.

\begin{lemma}
\label{thm:blw}
The return value \(\vb\) of  \Cref{alg:blw} is a vector of \(O(d)\)-bounded block Lewis weight overestimates.
Further, this vector of overestimates is found in \(\text{polylog}(k, m, d)\) leverage score overestimate computations.
\end{lemma}

The goal of the rest of this section is to analyze \Cref{alg:blw} and to prove \Cref{thm:blw}. Let us briefly describe the analysis of \Cref{alg:blw}. We first give a collection of potential functions \(\phi_i(\vb)\), with the goal of showing the potential of \(\phi_i(\overline{\vb})\) decreases with \(T\). A low potential will also imply that \(\overline{\vb}\) is nearly a vector of block Lewis weight overestimates.
For each \(i \in [m]\) define \(\phi_i \colon \R^m \to \R\) by
\[\phi_i(\vb) \coloneqq \ln\inparen{\frac{1}{\vb_i} \sum_{j \in S_i} \tau_j(\mB^{1/2 - 1/p} \mA)}.\]

The key property of this potential is convexity, which we now show.
\begin{lemma}
\label{lemma:alg_phi_convex}
    For each \(i \in [m]\), \(\phi_i\) is convex.
\end{lemma}
\begin{proof}[Proof of \Cref{lemma:alg_phi_convex}]
     Our argument for the convexity of this function is similar to the one given in \cite[Lemma A.2]{jls21}.
     First, notice that by the definition of \(\tau_j\), \(\phi(\vb)\) is equal to
     \begin{align*}
         \phi_i(\vb) = \ln\inparen{\frac{1}{\vb_i^{2/p}} \sum_{j \in S_i} \va_j^\top \inparen{\mA^\top \mB^{1 - 2/p} \mA}^{-1} \va_j}.
     \end{align*}
     Since \(- \frac{2}{p} \ln(b_i)\) is convex, it suffices to show the convexity of
     \[f(\vb) \coloneqq \ln\inparen{\sum_{j \in S_i} \va_j^\top \inparen{\mA^\top \mB^{1 - 2/p} \mA}^{-1} \va_j)}.\]
     Now we define for each \(j \in [k]\) the function
     \[h_j(\vb) \coloneqq \ln \inparen{\va_j^\top (\mA^\top \mB^{1-2/p} \mA)^{-1} \va_j}.\]
     A version of this function without repeated entries in \(\mB\) was shown to be convex in \cite[Lemma A.2]{jls21}, but \(h_j(\vb)\) is still convex.
     Next, notice that we may write \[
     f(\vb) = \ln\inparen{\sum_{j \in S_i} \exp(h_j(\vb))}.
     \]
     Now taking \(\vb, \vb' \in \R^m\) and \(\lambda \in [0, 1]\) we have
     \begin{align}
        f(\lambda \vb + (1-\lambda) \vb') &= \ln\inparen{\sum_{j \in S_i} \exp(h_j(\lambda \vb + (1-\lambda) \vb'))} \nonumber
        \\
     &\le \ln\inparen{\sum_{j \in S_i} \exp(\lambda h_j(\vb) + (1-\lambda) h_j(\vb'))} \label{eq:alg_ph_cvx_1}\\
     &\le \lambda \ln\inparen{\sum_{j \in S_i} \exp(h_j(\vb))} + (1-\lambda) \ln\inparen{\sum_{j \in S_i} \exp(h_j(\vb'))} \label{eq:alg_ph_cvx_2} \\
     &= \lambda f(\vb) + (1-\lambda) f(\vb'), \nonumber
     \end{align}
     where \eqref{eq:alg_ph_cvx_1} follows from the convexity of \(h_j\) and the monotonicity of log-sum-exp, and \eqref{eq:alg_ph_cvx_2} is due to the convexity of log-sum-exp (see e.g. \cite[Section 3.1.5]{boyd2004convex}). Hence, we may conclude the proof of \Cref{lemma:alg_phi_convex}.
\end{proof}

We now give an argument that \(\phi_i(\overline{\vb}) = \widetilde{O}(1/T)\) using the convexity of \(\phi_i\).

\begin{lemma}\label{lemma:alg_analysis1}
    Assume that \(\textsc{OverLev}\) returns \(\nu\)-bounded leverage score overestimates.
   Then after \Cref{alg:blw} runs for \(T\) iterations, we have for all \(i \in [m]\) that
    \[\phi_i(\overline{\vb}) \le \frac{1}{T} \ln\inparen{\frac{m \nu}{d}}.\]
\end{lemma}
\begin{proof}[Proof of \Cref{lemma:alg_analysis1}]
    We have
    \begin{align*}
        \phi_i(\overline{\vb}) &\le \frac{1}{T} \sum_{t=1}^T \phi_i(\vb^{(t)}) &\text{Jensen's inequality} \\
        &= \frac{1}{T} \sum_{t=1}^T \ln \inparen{\frac{1}{b_i^{(t)}} \sum_{j \in S_i} \tau_j((\mB^{(t)})^{1/2 - 1/p} \mA)} &\text{Definition of \(\phi_i\)}\\
        &\le \frac{1}{T} \sum_{t=1}^T \ln \inparen{\frac{1}{b_i^{(t)}} \sum_{j \in S_i} \widetilde{\tau}_j^{(t)}} &\text{Definition of \(\tau_j\)}\\
        &= \frac{1}{T} \sum_{t=1}^T \ln \inparen{\frac{b_i^{(t+1)}}{b_i^{(t)}}} &\text{Line \ref{line:alg_iter_next}} \\
        &= \frac{1}{T} \ln\inparen{\frac{b_i^{(T+1)}}{b_i^{(1)}}} \\
        &= \frac{1}{T} \ln\inparen{m/n} + \frac{1}{T} \ln(b_i^{(T+1)}). &\text{Line \ref{line:alg_init}}
    \end{align*}
    In the above, for the sake of analysis we define \(\vb_i^{(T+1)}\) to be as if the algorithm executed \(T+1\) iterations. Now, \(b_i^{(T+1)} \le \|\widetilde{\tau}\|_1 \le \nu\) by the fact that \(\textsc{OverLev}\) returns \(\nu\)-bounded leverage score overestimates. We therefore conclude the proof of \Cref{lemma:alg_analysis1}.
\end{proof}

Now, we show that \(\text{polylog}(k,m,d)\) leverage score overestimate computations suffice to find the \((\nfrac{3}{2}\cdot n\))-bounded block Lewis weight overestimates. This proves \Cref{thm:blw}.

\begin{proof}[Proof of \Cref{thm:blw}]
By \Cref{thm:lev_ove}, \textsc{OverLev} returns \(O(d)\)-bounded leverage score overestimates.
We take \(T = O\inparen{\ln\inparen{m}}\); clearly this satisfies the desired runtime guarantee. Further, using \Cref{lemma:alg_analysis1} and taking the constant in \(T\) large enough, we have for each \(i\) that \(\phi_i(\overline{\vb}) \le \ln\inparen{\frac{3}{2}}\).
Thus we have
\begin{equation}\label{eq:alg_b_guar}
    \frac{3}{2} \cdot \overline{b_i} \ge \sum_{j \in S_i} \tau_j(\overline{\mB}^{1/2 - 1/p} \mA).
\end{equation}
And so
\begin{align*}
    b_i &\ge \frac{3}{2} \overline{b_i} \\
    &\ge \sum_{j \in S_i} \tau_j(\overline{\mB}^{1/2 - 1/p} \mA) &\text{by \eqref{eq:alg_b_guar}} \\
    &= \sum_{j \in S_i} \tau_j(\mB^{1/2 - 1/p} \mA). &\text{\Cref{fact:lev_score_facts}}
\end{align*}
We now manipulate this guarantee into the desired guarantee of block Lewis weight overestimates
by some algebra.
Splitting powers on the left hand side of the above, we get
\[
b_i^{2/p} \ge \frac{1}{b_i^{1-2/p}} \sum_{j \in S_i} \tau_j(\mB^{1/2 - 1/p} \mA).
\]
Taking this to the \(p/2\)th power, we obtain
\begin{align*}
    \vb_i \ge \inparen{\frac{\sum_{j \in S_i} \tau_j(\mB^{1/2-1/p} \mA)}{b_i^{1-2/p}}}^{p/2} = \beta_i(\mB^{1-2/p})^p,
\end{align*}
as desired. To bound \(\norm{\vb}_1\), notice
\begin{align*}
    \norm{\vb}_1 = \frac{3}{2} \|\overline{\vb}\|_1 \le \frac{3}{2} \frac{1}{T} \sum_{t=1}^T \|\vb^{(t)}\|_1 \le \frac{3}{2} \nu \le O(n).
\end{align*}
Finally, to obtain the concentration statement, we observe that the above implies that we get a measure $\vlambda$ and a rounding $\mW$ that are $\Fstar$-block Lewis estimates for $\Fstar=O(d)$. This concludes the proof of \Cref{thm:blw}.
\end{proof}

\subsubsection{Special case -- \texorpdfstring{$p = 2, p_1,\dots,p_m\ge2$}{p=2,p1,...,pm>=2}}

Finally, we are ready to introduce and analyze the algorithm for the case where $p=2$ and $p_1,\dots,p_m \ge 2$. See \Cref{alg:blw_base}. The main property of \Cref{alg:blw_base} is given in \Cref{thm:blw_base}.

\begin{algorithm}[ht]
\caption{Algorithm to compute block Lewis weight overestimates, $p=2$ and $p_1,\dots,p_m \ge 2$}
\begin{algorithmic}[1]\label{alg:blw_base}
\State \textbf{Input:} \(\mA \in \R^{n \times d}\), group structure \((S_1,\dots,S_m,p_1,\dots,p_m)\).
\State Initialize \(\vb^{(0)} = \frac{d}{m} \cdot \onev\) and \(\vu^{(0)}\) such that for all $1 \le i \le m$, \(u^{(0)}_j = 1/\abs{S_i}\) for all \(j \in S_i\)
\For{\(t = 1, \ldots, T - 1\)}
    \State \(\widetilde{\tau}^{(t)} = \textsc{OverLev}(\mV(\vu^{(t-1)})^{1/2}\mA)\) \Comment{$\mV(\vu)$ is such that $v_j = u_j^{1-2/p_i}$ for all $1 \le i \le m$ and $j \in S_i$}
    \State \(b_{i}^{(t)} = \sum_{j \in S_i} \widetilde{\tau}_j^{(t)}\) for all $1 \le i \le m$
    \State \(u^{(t)}_j = \widetilde{\tau}_j^{(t)}/(\sum_{j' \in S_i} \widetilde{\tau}_{j'}^{(t)})\) for all $1 \le i \le m$ and all $j \in S_i$
\EndFor
\State \(\overline{\vb} = \abs{S_i}^{1/T}\cdot\frac{1}{T}\sum_{t=1}^T \vb^{(t)}\)
\State \(\overline{\vu} = \frac{1}{T}\sum_{t=0}^{T-1} \vu^{(t)}\)
\State \Return \((\overline{\vb},\overline{\vu})\)
\end{algorithmic}
\end{algorithm}

\begin{lemma}
\label{thm:blw_base}
If $\textsc{OverLev}$ is a routine that returns leverage score overestimates whose sum is at most $\nu$, then \Cref{alg:blw_base} returns $\Fstar$-block Lewis overestimates (\Cref{defn:block_lewis_overestimate}) with $\Fstar = \max_{1 \le i \le m} \abs{S_i}^{1/T} \cdot \nu$.

In particular, if $T = \max_{1 \le i \le m} \log\abs{S_i}$ and $\nu \le (4/e)d$, then we get $\Fstar = 4d$.
\end{lemma}

As with the analysis of \Cref{alg:blw}, we first prove that a relevant potential is convex and then show that we can control it effectively.

\begin{lemma}
\label{lemma:blw_convex_potential_general}
Let $\mU \in \R^{k \times k}$ be a nonnegative diagonal matrix. Let $\mV(\vu)$ denote the matrix such that for all $i$ and $j \in S_i$, we have $v_j = \vu_j^{1-2/p_i}$. Then, $\phi$ as defined below is log-convex in $\vu$.
\begin{align*}
    \text{for all } j \in S_1 \cup \dots \cup S_m:\quad\quad\phi_j(\vu) = \lnv{\frac{\tau_j\inparen{(\mV(\vu))^{1/2}\mA}}{u_j}}
\end{align*}
\end{lemma}
\begin{proof}[Proof of \Cref{lemma:blw_convex_potential_general}]
We expand the above definition after taking the $\ln$ of both sides.
\begin{align*}
    \phi_j(\vu) &= \lnv{\frac{\tau_j\inparen{(\mV(\vu))^{1/2}\mA}}{u_j}} = \lnv{\va_j^{\top}\inparen{\mA^{\top}\mV(\vu)\mA}^{-1}\va_j} + \lnv{\frac{v_j}{u_j}} \\
    &= \lnv{\va_j^{\top}\inparen{\mA^{\top}\mV(\vu)\mA}^{-1}\va_j} - \frac{2}{p_i}\lnv{u_j}.
\end{align*}
The last term is convex, so it is sufficient to argue that the first term is convex. To see this, first observe that by concavity of $x^{1-2/p_i}$, we have
\begin{align*}
    \lambda\mV(\vu^{(1)}) + (1-\lambda)\mV(\vu^{(2)}) \preceq \mV(\lambda\vu^{(1)} + (1-\lambda)\vu^{(2)}),
\end{align*}
which implies
\begin{align*}
    \lnv{\va_j^{\top}\inparen{\mA^{\top}\mV(\lambda\vu^{(1)} + (1-\lambda)\vu^{(2)})\mA}^{-1}\va_j} &\le \lnv{\va_j^{\top}\inparen{\lambda\mA^{\top}\mV(\vu^{(1)})\mA + (1-\lambda)\mA^{\top}\mV(\vu^{(2)})\mA}^{-1}\va_j} \\
    &\le \lambda\lnv{\va_j^{\top}\inparen{\mA^{\top}\mV(\vu^{(1)})\mA}^{-1}\va_j} \\
    &\quad + (1-\lambda)\lnv{\va_j^{\top}\inparen{\mA^{\top}\mV(\vu^{(2)})\mA}^{-1}\va_j}.
\end{align*}
Above, the last line follows from the well-known (see, e.g., \cite[Lemma 3.4]{ccly19}) fact that $\lnv{\va_j^{\top}\mM^{-1}\va_j}$ is convex in $\mM$ where $\mM$ is symmetric positive-semidefinite. This completes the proof of \Cref{lemma:blw_convex_potential_general}.
\end{proof}

Next, we will show that any choice of nonnegative $\vu$ such that $\sum_{j \in S_i} u_j = b_i^{\inparen{1-\frac{2}{p}}\cdot\inparen{\frac{p_i}{p_i-2}}}$ can be used to find a rounding matrix $\mW$.

\begin{lemma}
\label{lemma:blw_rounding_matrix}
For all $\vx\in\R^d$ and all nonnegative $\vu$ such that $\sum_{j \in S_i} u_j = b_i^{\inparen{1-\frac{2}{p}}\cdot\inparen{\frac{p_i}{p_i-2}}}$, we have
\begin{align*}
    \frac{\norm{\mV(\vu)^{1/2}\mA\vx}_2}{\inparen{\sum_{i' \le m} b_{i'}}^{1/2-1/p}} \le \gnorm{\mLambda^{1/2-1/p}\mA\vx}{2}.
\end{align*}
\end{lemma}
\begin{proof}[Proof of \Cref{lemma:blw_rounding_matrix}]
We start with the LHS.
\begin{align*}
    \norm{\mV(\vu)^{1/2}\mA\vx}_2^2 &= \sum_{i=1}^m \sum_{j \in S_i} u_j^{1-2/p_i}\abs{\ip{\va_j,\vx}}^2 \\
    &\le \sum_{i=1}^m b_i^{1-2/p}\norm{\mA_{S_i}\vx}_{p_i}^2 & \text{ by H\"older's Inequality with powers } \frac{p_i}{p_i-2}, \frac{p_i}{2}
\end{align*}
Normalizing concludes the proof of \Cref{lemma:blw_rounding_matrix}.
\end{proof}

Now, we show that finding uniform overestimates $\tau_j/u_j \le b_i^{\inparen{\frac{2}{p}-\frac{2}{p_i}}\cdot\frac{p_i}{p_i-2}}$ is enough to satisfy \Cref{defn:block_lewis_overestimate} with $\Fstar = \sum_{i \le m} b_i$.

\begin{lemma}
\label{lemma:blw_overestimates_conversion}
If we have $\vu$ such that $\tau_j(\mV(\vu)^{1/2}\mA)/u_j \le b_i^{\inparen{\frac{2}{p}-\frac{2}{p_i}}\cdot\frac{p_i}{p_i-2}}$ for all $i$, then the rounding matrix $\mW = \mV(\vu)\mB^{2/p-1}$ and measure $\lambda_i = b_i/(\sum_{i'\le m} b_{i'})$ is an $\Fstar$-block Lewis overestimate (\Cref{defn:block_lewis_overestimate}) with $\Fstar = \sum_{i'\le m} b_{i'}$.
\end{lemma}
\begin{proof}[Proof of \Cref{lemma:blw_overestimates_conversion}]
Observe that we must have
\begin{align*}
    \frac{\tau_j\inparen{\mV^{1/2}\mA}}{v_j} \le v_j^{\frac{2}{p_i-2}} \cdot b_i^{\inparen{\frac{2}{p}-\frac{2}{p_i}}\cdot\frac{p_i}{p_i-2}}.
\end{align*}
Let $T = \sum_{i' \le m} b_{i'}$. Following \Cref{lemma:block_lewis}, let $\lambda_i = b_i/T$ and $\mW$ be such that
\begin{align*}
    \frac{\mV^{1/2}}{T^{1/2-1/p}} = \mW^{1/2}\mLambda^{1/2-1/p}.
\end{align*}
In particular, this means that
\begin{align*}
    w_j = \frac{v_j}{(T\lambda_i)^{1-2/p}} = \frac{v_j}{b_i^{1-2/p}}.
\end{align*}
Hence,
\begin{align*}
    \inparen{\sum_{j \in S_i} \inparen{\frac{\tau_j\inparen{\mV^{1/2}\mA}}{w_j}}^{p_i/2}}^{2/p_i} &= b_i^{1-2/p}\inparen{\sum_{j \in S_i} \inparen{\frac{\tau_j\inparen{\mV^{1/2}\mA}}{v_j}}^{p_i/2}}^{2/p_i} \\
    &\le b_i^{1-2/p} \cdot b_i^{\inparen{\frac{2}{p}-\frac{2}{p_i}}\cdot\frac{p_i}{p_i-2}} \cdot \inparen{\sum_{j \in S_i} v_{j}^{\frac{2}{p_i-2} \cdot \frac{p_i}{2}}}^{2/p_i} \\
    &= b_i^{1-2/p}\cdot b_i^{\inparen{\frac{2}{p}-\frac{2}{p_i}}\cdot\frac{p_i}{p_i-2}} \cdot \inparen{\sum_{j \in S_i} u_j}^{2/p_i} \\
    &= b_i^{1-2/p}\cdot b_i^{\inparen{\frac{2}{p}-\frac{2}{p_i}}\cdot\frac{p_i}{p_i-2}} \cdot \inparen{\sum_{j \in S_i} u_j}^{2/p_i} \\
    &= b_i^{1-2/p}\cdot b_i^{\inparen{\frac{2}{p}-\frac{2}{p_i}}\cdot\frac{p_i}{p_i-2}} \cdot \inparen{b_i^{\inparen{1-\frac{2}{p}}\cdot\inparen{\frac{p_i}{p_i-2}}}}^{2/p_i} = b_i = \lambda_i\inparen{\sum_{i'\le m} b_{i'}},
\end{align*}
which is exactly the statement of \Cref{lemma:blw_overestimates_conversion}.
\end{proof}

We now have the tools we need to prove \Cref{thm:blw_base}.

\begin{proof}[Proof of \Cref{thm:blw_base}]
Note that \Cref{alg:blw_base} and \Cref{thm:blw_base} are reminiscent of \cite[Algorithm 2 and Theorem 4]{jls22}.

Using \Cref{lemma:blw_convex_potential_general}, we begin with using the convexity of the potential $\phi$. For all $j \in S_i$, we have
\begin{align*}
    \phi_j(\overline{\vu}) &\le \frac{1}{T}\sum_{t=0}^{T-1} \phi\inparen{\vu^{(t)}} = \frac{1}{T}\sum_{t=0}^{T-1} \lnv{\frac{\tau_j\inparen{\mV(\vu^{(t)})^{1/2}\mA}}{u_j^{(t)}}} \le \frac{1}{T}\sum_{t=0}^{T-1}\lnv{\frac{\widetilde{\tau}_j^{(t+1)}}{u_j^{(t)}}} \\
    &= \frac{1}{T}\sum_{t=0}^{T-1}\inparen{\lnv{\frac{u_j^{(t+1)}}{u_j^{(t)}}} + \lnv{\sum_{j' \in S_i} \widetilde{\tau}_{j'}^{(t+1)}}} = \frac{1}{T}\lnv{\frac{u_j^{(T)}}{u_j^{(0)}}} + \frac{1}{T}\sum_{t=0}^{T-1}\lnv{\sum_{j' \in S_i} \widetilde{\tau}_{j'}^{(t+1)}} \\
    &\le \frac{1}{T}\lnv{\frac{u_j^{(T)}}{u_j^{(0)}}} + \lnv{\frac{1}{T}\sum_{t=0}^{T-1}\sum_{j' \in S_i} \widetilde{\tau}_{j'}^{(t+1)}} = \frac{1}{T}\lnv{\frac{u_j^{(T)}}{u_j^{(0)}}} + \lnv{\abs{S_i}^{-1/T}b_{i}} \le \lnv{b_i}.
\end{align*}
We now apply \Cref{lemma:blw_overestimates_conversion} and see that the measure $\lambda_i = b_i/\norm{\overline{\vb}}_1$ is a $\Fstar$-block Lewis overestimate (\Cref{defn:block_lewis_overestimate}) with $\Fstar = \norm{\overline{\vb}}_1$.

Now, observe that
\begin{align*}
    \norm{\overline{\vb}}_1 \le \max_{1 \le i \le m} \abs{S_i}^{1/T} \cdot \frac{1}{T}\sum_{t=1}^T\norm{\vb^{(t)}}_1 \le \max_{1 \le i \le m} \abs{S_i}^{1/T} \cdot \nu \le 4d,
\end{align*}
where we use our setting of $T$ and the fact that $\textsc{OverLev}$ returns leverage score estimates whose sum is at most $(4/e)d$. This completes the proof of \Cref{thm:blw_base}.
\end{proof}

We are finally ready to give the proof of \Cref{thm:computeblw}.

\begin{proof}[Proof of \Cref{thm:computeblw}]
We have three cases.
\begin{itemize}
    \item If $0 < p < 4$ and $p_1=\dots=p_m=2$, then the algorithm and guarantee on the weights follow from \Cref{alg:blw_lewismap} and \Cref{lemma:contract_alg}.
    \item If $p \ge 2$ and $p_1=\dots=p_m=2$, then the algorithm and guarantee on the weights follow from \Cref{alg:blw} and \Cref{thm:blw}.
    \item If $p = 2$ and $p_1,\dots,p_m\ge 2$, then the algorithm and guarantee on the weights follow from \Cref{alg:blw_base} and \Cref{thm:blw_base}.
\end{itemize}
We plug these guarantees into \Cref{thm:concentration} and conclude the proof of \Cref{thm:computeblw}.
\end{proof}

\subsection{Minimizing sums of Euclidean norms (Proof of \texorpdfstring{\Cref{thm:msnalg}}{Theorem 3})}
\label{sec:applications_msn}

Recall the minimizing sums of Euclidean norms (MSN) problem \eqref{eq:msn}.
Given \(\mA \in \R^{n \times d}\), a partition \(S_1, \ldots, S_m\) of \([n]\), and \(\vb_1 \in \R^{|S_1|}, \ldots, \vb_m \in \R^{|S_m|}\), we would like to find $\widehat{\vx}$ such that
\[\sum_{i=1}^m \|\mA_{S_i} \widehat{\vx} - \vb_i\|_2\le (1+\eps)\min_{\vx \in \R^d} \sum_{i=1}^m \|\mA_{S_i} \vx - \vb_i\|_2.\]
\citet{xy97} give an algorithm with iteration complexity $\widetilde{O}(\sqrt{m}\logv{\nfrac{1}{\eps}})$ for the above problem (though for an additive approximation guarantee instead of a multiplicative one), where each iteration reduces to solving linear systems in matrices $\mA^{\top}\mD\mA$ for block-diagonal matrices, where each block has size $(\abs{S_i}+1) \times (\abs{S_i}+1)$. Their algorithm is based on the primal-dual interior point method framework.


Each system solve takes the following form. Let \(\widetilde{\mA}\) be a block matrix with \(- \mI_m\) in one block, and \(\mA\) in another.
The goal is to find \(\vy\) in the system
\[\widetilde{\mA}^\top \mD \widetilde{\mA} \vy = \vz\]
where \(\mD\) is a block matrix, with one \((1 + |S_i|) \times (1 + |S_i|)\) sized block for each group
(see \cite[equation (4.13)]{xy97}).

Our main result of this section is \Cref{thm:msnalg}, which gives an improved iteration complexity for \eqref{eq:msn} when $m \gg n$.

\msnalg*

We prove \Cref{thm:msnalg} by sparsifying the objective \eqref{eq:msn} using \Cref{thm:computeblw} and then applying the primal-dual interior point method from \cite{xy97}. We state the guarantee of this algorithm in \Cref{lemma:ipm_log_barrier}.

\begin{lemma}
\label{lemma:ipm_log_barrier}
Let $\mA \in \R^{n \times d}$ and \(\vb \in \R^n\), and  $S_1, \dots, S_m$ be a partition of $k$. There exists an algorithm that returns $\widehat{\vx}$ such that
\begin{align*}
    \sum_{i=1}^m \norm{\mA_{S_i}\widehat{\vx}-\vb_{S_i}}_2 \le (1+\eps)\min_{\vx\in\R^d} \sum_{i=1}^m \norm{\mA_{S_i} \vx-\vb_{S_i}}_2.
\end{align*}
The algorithm runs in $\widetilde{O}\inparen{\sqrt{m}\logv{\nfrac{1}{\eps}}}$ calls to a linear system solver in matrices of the form $\mA^{\top}\mD\mA$ for block-diagonal matrices $\mD$, where each block has size $(\abs{S_i}+1) \times (\abs{S_i}+1)$.
\end{lemma}
\begin{proof}[Proof of \Cref{lemma:ipm_log_barrier}]
The guarantee we will reduce to is \cite[Theorem 5.2]{xy97}. However, the guarantee there is stated for an additive approximation, and we desire a multiplicative approximation. We therefore apply a few transformations to our problem so that we can apply this guarantee.

Let
\begin{align*}
    \vx_0 &\coloneqq \argmin{\vx\in\R^d} \norm{\mA\vx-\vb}_2^2.
\end{align*}
This can be found in one linear system solve. Let $V \coloneqq \norm{\mA\vx_0-\vb}_2$ and consider the following modified optimization problem:
\begin{align}
    \min_{\vx-\vx_0\in\R^d} \frac{1}{V}\sum_{i=1}^m \norm{\mA_{S_i}(\vx-\vx_0)-(\vb_{S_i}-\mA_{S_i} \vx_0)}_2.\label{eq:modified_msn}
\end{align}
We will invoke \cite[Theorem 5.2]{xy97} on the above problem \eqref{eq:modified_msn}, folding the \(\frac{1}{V}\) factor into \(\mA\) and \(\vb\). Clearly, this problem is equivalent to the problem we started with. Next, let $\xstar$ be given by
\begin{align*}
    \xstar &\coloneqq \argmin{\vx\in\R^d} \sum_{i=1}^m \norm{\mA_{S_i}\vx-\vb_{S_i}}_2,
\end{align*}
where we use \(\vb \in \R^n\) for the vector formed by stacking \(\vb_{S_1}, \ldots, \vb_{S_m}\).
Let $\opt \coloneqq \nfrac{1}{V} \cdot \sum_{i=1}^m \norm{\mA_{S_i}\xstar-\vb_{S_i}}_2$. 
Because \(\|\cdot\|_2 \le \|\cdot\|_1\), we have
\begin{align*}
    1=\frac{\norm{\mA\vx_0-\vb}_2}{V} \le \frac{\norm{\mA\xstar-\vb}_2}{V} \le \frac{1}{V}\sum_{i=1}^m \norm{\mA_{S_i}\xstar-\vb_{S_i}}_2 = \opt.
\end{align*}
Furthermore, we have by \(\|\cdot\|_1 \le \sqrt{m}\|\cdot\|_2\) (applied by considering the summation over \(i\) as an \(\ell_1\) norm) that
\begin{align*}
    \max_{1 \le i \le m} \frac{1}{V}\norm{\vb_{S_i}-\mA_{S_i}\vx_0}_2 \le \frac{1}{V}\sum_{i=1}^m \norm{\vb_{S_i}-\mA_{S_i}\vx_0}_2 \le \sqrt{m} \cdot \frac{\norm{\mA\vx_0-\vb}_2}{V} = \sqrt{m}.
\end{align*}
This implies that all the offset vectors of our transformed problem \eqref{eq:modified_msn} have polynomially bounded norm. This suffices for our iteration complexity bound, because the iteration complexity of the algorithm in \cite{xy97} depends logarithmically on the maximum norm of these vectors. Now, using their method we can solve \eqref{eq:modified_msn} up to $\eps$ additive error. This means we find $\widehat{\vx}$ such that
\begin{align*}
    \sum_{i=1}^m \norm{\mA_{S_i}(\widehat{\vx}-\vx_0)-(\vb_{S_i}-\mA_{S_i} \vx_0)}_2 &\le \sum_{i=1}^m \norm{\mA_{S_i}(\xstar-\vx_0)-(\vb_{S_i}-\mA_{S_i} \vx_0)}_2 + V\eps \\
    &= V \cdot \opt\inparen{1+\frac{\eps}{\opt}} \le V \cdot \opt\inparen{1+\eps}.
\end{align*}
where in the last inequality, we use that \(\opt \ge 1\).
Since $V \cdot \opt$ is the original optimal objective value, this completes the proof of \Cref{lemma:ipm_log_barrier}.
\end{proof}

We remark that instead of the $\sqrt{m}$ above, one can get a $\sqrt{d}$-factor relationship between the optimal objective for a least squares relaxation of our problem by using the block Lewis weights with $p_1=\dots=p_m=2$ and $p=1$. However, this will only impact lower order terms.

We are now ready to prove \Cref{thm:msnalg}.

\begin{proof}[Proof of \Cref{thm:msnalg}]
We apply \Cref{thm:computeblw} for \(p_1 = \cdots = p_m = 2\) and \(p = 1\) with approximation $\eps/3$ to the group matrices $\insquare{\mA_{S_i} \vert \vb_{S_i}}$ to find a sparsified objective
with \(\widetilde{m} = O(\eps^{-2} \cdot d(\log d)^2 \log(\nfrac{d}{\eps}))\) terms.
This requires \(\widetilde{O}(1)\) linear system solves.
This implies a \((1+\eps)\) approximation to an un-sparsified objective over any vector in \(\R^{d+1}\),
and the approximation we need comes by only considering vectors in \(\R^{d+1}\) whose last entry is \(-1\).
We plug this into the guarantee of \Cref{lemma:ipm_log_barrier} with approximation $\eps/3$. Since $(1+\eps/3)^2 \le 1+\eps$ and $(1+\eps/3)/(1-\eps/3) \le 1+\eps$, this returns a $(1+\eps)$-approximate minimizer to \eqref{eq:msn}, completing the proof of \Cref{thm:msnalg}.
\end{proof}

\chapter{Block Lewis weights for distributionally robust linear regression\label{chapter:regression}}

In this chapter, we continue studying applications of block Lewis weights. This time, we use them for minimizing a multidistributional linear regression loss. The material in this chapter is based on a joint work with Kumar Kshitij Patel \cite{mp24}.

\section{Introduction}\label{sec:gp_reg_intro}
Machine learning algorithms and their training datasets have grown tremendously in the past decade, both in size and complexity. This increased model complexity has made it challenging to interpret and predict their behavior in unobserved scenarios. Hence, many applications that involve societal decisions still rely on simple, interpretable models like linear regression, often after feature engineering. Examples of such applications are predicting housing prices across cities, estimating wages across industries, forecasting loan amounts across banks, predicting life insurance premiums for different groups, and projecting energy consumption in various communities. 

A shared safety and sometimes legal concern across the above applications is the potential for wildly different model qualities for different distributions, i.e., outputting a notably worse model for some source data distributions~\cite{data2014seizing, barocas2016big,hardt2016equality,veale2018fairness,selbst2019fairness,berk2021fairness, corbett2023measure,chouldechova2016fair,kleinberg2018algorithmic, agarwal2019fair,cohen2024fairness,svwz24}. Specifically, consider fitting a linear model $\vx\in\R^d$ to make real predictions on some task over $m$ groups where group $i$'s dataset consists of $n_i$ entries and is denoted by $S_i=\{(\va_i^j, b_i^j)\}_{j\in [n_i]}$. The \textit{utilitarian} or the total-cost-minimizing objective minimizes the average squared prediction error across groups, i.e.,        
\begin{align}
    \min_{\vx\in\R^d}\frac{1}{m}\sum_{i\in[m]}\frac{1}{n_i}\norm{\mA_{S_i}\vx-\vb_{S_i}}_2^2\enspace,\label{eq:utilitarian}
\end{align}
where $\mA_{S_i} \coloneqq [\va_i^1 \dots \va_i^{n_i}]^{\top}\in \R^{n_i \times d}$ is the feature matrix and $\vb_{S_i} \coloneqq [b_i^1 \dots b_i^{n_i}]^{\top} \in \R^{n_i}$ is the label vector for group $i\in[m]$. 

Because of the inherent heterogeneity of the datasets, the model derived from optimizing objective \eqref{eq:utilitarian} may be particularly bad for some groups, in that the prediction error could be disproportionately higher for these groups. To overcome these limitations, the following \textit{egalitarian} or group Distributionally Robust Optimization (DRO) objective  has been considered in several recent works \cite{ben2013robust, duchi2016statistics, sagawa2019distributionally, levy2020large, soma2022optimal,aakmrz22,svwz24},
\begin{align}
    \min_{\vx\in\R^d} \max_{i\in[m]}\frac{1}{n_i}\norm{\mA_{S_i}\vx-\vb_{S_i}}_2^2\enspace.\label{eq:main_objective}
\end{align}
Objective \eqref{eq:main_objective} is the ``fairest'' objective among all objectives that balance utility and distributional robustness by ensuring that no one group has a loss that is too high ~\cite{kleinberg2018algorithmic,chouldechova2018frontiers,asadpour2022sequential,chen2022fair,rahmattalabi2019exploring,golrezaei2024online}

Since \eqref{eq:main_objective} is a convex problem, it is natural to apply standard black-box optimization techniques to solve them. However, we identify several challenges in applying existing methods:
    \paragraph{Efficient first-order algorithms have geometry-dependent rates.} To our knowledge, using an efficient first-order method (such as sub-gradient descent) will incur a geometry-dependent runtime. In particular, if the matrices $\mA_{S_i}$ or if the stacked matrix $\mA\coloneqq [\mA_{S_1}^{\top} \dots \mA_{S_M}^{\top}]^{\top}$ are poorly conditioned, then this will be reflected accordingly in the convergence rates. This is a drawback of the existing results by \citet{aakmrz22} and \citet{svwz24}.
    \paragraph{Objective \eqref{eq:main_objective} is not smooth.} Since the objective is the pointwise maximum of several continuous functions, the derivative is not well-defined at the points at which the maximizing function changes. Thus, applying subgradient descent to this objective without a customized analysis will result in a rather unimpressive $1/\eps^{2}$ dependence in the iteration complexity. 
    \paragraph{Min-max optimization approaches have a $1/\eps^{2}$ dependence on iteration complexity.} Since problem \ref{eq:main_objective} is a min-max optimization objective, it is also natural to try to use game theory-inspired approaches that use some oracle (such as gradients) for each group as a black box. Perhaps the most basic such algorithm is casting objective \eqref{eq:main_objective} as a repeated game between a min player (equipped with a no-regret algorithm) and a max player (equipped with the best response oracle). The main shortcoming of this approach is that even though the function for each group is smooth, the iteration complexity (to get $\eps$ average regret) for smooth online convex optimization still has an unimpressive $1/\eps^2$ dependence (as opposed to $1/\eps$ for smooth convex optimization)~\cite{soma2022optimal, zhang2024stochastic}. Thus, this approach is no better than directly applying sub-gradient descent to objective \eqref{eq:main_objective}.
    \paragraph{Interior point methods have a poor iteration complexity for large $m$.} Another natural approach (that can partially address the previous two issues), following the discussion by Boyd and Vandenberghe \cite[Section 6.4]{boyd2004convex}, is to rewrite the problem \eqref{eq:main_objective} in its epigraph form and use an interior point method (IPM) to solve the resulting problem (which, in this case, is a quadratically constrained linear program). Unfortunately, this will give an algorithm whose analysis is only known to yield an iteration complexity of $O(\sqrt{m})$, where each iteration solves a linear system in matrices of the form $\mA^{\top}\mB\mA$ for a block-diagonal $\mB$ (see \Cref{rem:linear_systems}).
    A na\"ive implementation of this algorithm will thus have a superlinear runtime in the number of groups, which is undesirable when the number of groups is large. Furthermore, note that \eqref{eq:main_objective} is not a linear program when at least one group $i$ is such that $n_i > 1$. So, we cannot immediately apply recent advances in linear programming that get iteration complexities independent of the number of constraints \cite{ls19}.

Hence, designing an algorithm without these shortcomings requires novel ideas. 

\subsection{Our results}
\label{sec:gp_our_results}

In this chapter, we give a new algorithm (\Cref{alg:fair_regression}) to approximately optimize \eqref{eq:main_objective} that overcomes the above difficulties. We state our algorithm's iteration complexity in the following theorem. 

\begin{mainthm}[Robust regression]\label{mainthm:fair_regression_iteration_complexity}
Let $\mA_{S_i} \in \R^{n_i \times d}$ and $\vb_{S_i} \in \R^{n_i}$ for all $i\in[m]$. Denote their concatenations by $\mA\coloneqq [\mA_{S_1}^{\top} \dots \mA_{S_M}^{\top}]^{\top}\in\R^{n\times d}$ and $\vb \coloneqq [\vb_{S_1}^{\top}\dots\vb_{S_M}^{\top}]^{\top}\in \R^n$ where $n:= \sum_{i\in[m]}n_i$. Let $\eps > 0$. Then \Cref{alg:fair_regression} returns $\xhat$ such that,
\begin{align}
    \max_{i\in[m]} \frac{1}{\sqrt{n_i}}\norm{\mA_{S_i}\xhat-\vb_{S_i}}_2 \le (1+\eps)\cdot\min_{\vx\in\R^d}\max_{i\in[m]} \frac{1}{\sqrt{n_i}}\norm{\mA_{S_i}\vx-\vb_{S_i}}_2\enspace,\label{eq:objective_approx}
\end{align}
and it runs in
\begin{align*}
    O\inparen{\frac{\min\inbraces{\mathsf{rank}(\mA),m}^{1/3}\inparen{\logv{\frac{n\log m}{\eps}}^{14/3}+\logv{m}}}{\eps^{2/3}}}
\end{align*}
linear system solves in matrices of the form $\mA^{\top}\mB\mA$, where $\mB$ is a block-diagonal matrix for which block $i$ has size $n_i \times n_i$.
\end{mainthm}
We prove \Cref{mainthm:fair_regression_iteration_complexity} in \Cref{sec:gp_robust_proof}.

We compare the guarantee of \Cref{mainthm:fair_regression_iteration_complexity} against the other baselines in \Cref{table:compare}. Unlike the aforementioned first-order methods, our algorithm has no geometry-dependent terms. Additionally, our algorithm improves over the standard log-barrier IPM when the desired accuracy $\eps \ge m^{-1/4}$ --- this improvement is more pronounced when $m \gg \rank{\mA}$, i.e. when the number of data sources is much larger than the dimension of the parameter vector $\vx$. Additionally, for $\eps \ge \rank{\mA}^{-1/4}$, our guarantee matches the best known guarantee for $\ell_{\infty}$ regression \cite{ls19,jls21}.

\begin{remark}[Why use linear-system-solve complexity?]\label{rem:linear_systems}
    We benchmark our algorithms using the number of linear-system-solves for a few reasons. First, this is typically how second-order algorithms are compared, such as interior point methods for linear programming \cite{ls19}. Second, the particular structure of the linear system solves presents the possibility for a faster amortized runtime for the systems over the whole algorithm. This observation, combined with an understanding of how the linear systems changed between iterations, was used recently to get fast runtimes for linear programming \cite{ls19} and $\ell_{\infty}$ regression \cite{ajk24}.
\end{remark}

\paragraph{Interpolating between robust and nonrobust optimization.} We also study the following family of objectives that interpolate between \eqref{eq:utilitarian} and \eqref{eq:main_objective} for different values of $p\geq 2$,
\begin{align}
    \min_{\vx\in\R^d}\frac{1}{m}\sum_{i\in[m]}\left(\frac{1}{n_i}\norm{\mA_{S_i}\vx-\vb_{S_i}}_2^2\right)^{p/2}\enspace.\label{eq:interpolation}
\end{align}
In particular, note that choosing $p=2$ in the above objective gives us the average least-squares problem in objective \eqref{eq:utilitarian}, while $p\to\infty$ recovers objective \eqref{eq:main_objective}. Varying $p$ from $2$ to $
\infty$ and minimizing gives solutions that interpolate between utilitarian and egalitarian approaches, allowing for a smooth trade-off between utility and robustness. To this end, we give an algorithm (\Cref{alg:gp_regression_alg}) to approximately optimize \eqref{eq:interpolation} and prove the following guarantee about its iteration complexity.

\begin{mainthm}[Trading off utility and robustness]\label{mainthm:gp_regression_iteration_complexity}
Let $\mA_{S_i} \in \R^{n_i \times d}$ and $\vb_{S_i} \in \R^{n_i}$ for all $i\in[m]$. Denote their concatenations by $\mA\coloneqq [\mA_{S_1}^{\top} \dots \mA_{S_M}^{\top}]^{\top}\in\R^{n\times d}$ and $\vb \coloneqq [\vb_{S_1}^{\top}\dots\vb_{S_M}^{\top}]^{\top}\in \R^n$ where $n \coloneqq \sum_{i\in[m]}n_i$. Let $p \ge 2$ and $\eps > 0$. Then \Cref{alg:gp_regression_alg} returns $\xhat$ such that,
\begin{align}\label{eq:interpolation_approx}
    \inparen{\sum_{i=1}^m \inparen{\frac{1}{\sqrt{n_i}}\norm{\mA_{S_i}\xhat-\vb_{S_i}}_2}^p}^{1/p} \le (1+\eps)\cdot\min_{\vx\in\R^d}\inparen{\sum_{i=1}^m \inparen{\frac{1}{\sqrt{n_i}}\norm{\mA_{S_i}\vx-\vb_{S_i}}_2}^p}^{1/p}
\end{align}
and runs in
\begin{align*}          
    O\inparen{p^{O(1)}\min\inbraces{\rank{\mA},m}^{\frac{p-2}{3p-2}}\logv{\frac{pd}{\eps}}^3}
\end{align*}
linear system solves in matrices of the form $\mA^{\top}\mB\mA$, where $\mB$ is a block-diagonal matrix for which block $i$ has size $n_i \times n_i$.
\end{mainthm}
We prove \Cref{mainthm:gp_regression_iteration_complexity} in \Cref{sec:gp_interpolating_proof}.

In the special case where $n_i=1$ for all $i$ (and therefore the problem is $\ell_p$ regression for $p \ge 2$), the complexity promised by \Cref{mainthm:gp_regression_iteration_complexity} is comparable to that promised by \citet{jls21} for $\ell_p$ regression. The main difference is that our iteration complexity is unconditionally polynomial in $p$. In contrast, the comparable result from \cite{jls21} seems to require mild assumptions on the problem parameters (see the ``Discussion on numerical stability'' in \cite[Section 4]{jls21}).

\begin{remark}[Large values of $p$] 
    Note that for values of $p$ larger than $\log(m)$, solving \eqref{eq:main_objective} is morally equivalent to solving \eqref{eq:interpolation}. To intuitively see this, first recall that for any vector $\vx\in\R^d$ and $p=\log_2(m)$ we have, $\|\vx\|_\infty \leq \|\vx\|_{p} \leq 2\cdot\|\vx\|_\infty$. This implies that for all $i\in[m]$ we have the following for objective \eqref{eq:interpolation} (for $p=\log_2(m)$) for any $\vx\in\R^d$,
    \begin{align*}
        \max_{i\in[m]}\norm{\mA_{S_i}\vx-\vb_{S_i}}_2 \leq \left(\sum_{i\in[m]}\norm{\mA_{S_i}\vx-\vb_{S_i}}_2^{p}\right)^{1/p}\leq 2\cdot\max_{i\in[m]}\norm{\mA_{S_i}\vx-\vb_{S_i}}_2\enspace. 
    \end{align*}
    In particular, this means that minimizing the interpolating objective \eqref{eq:interpolation} also minimizes the robust objective \eqref{eq:main_objective} (up to numerical constants) and vice versa. Thus, for $p=\Omega\left(\log_2(m)\right)$, for our intended applications, it makes sense to minimize the robust objective instead. This is why, in \Cref{mainthm:gp_regression_iteration_complexity}, we do not care too much about the exponent on $p$ in the iteration complexity. Our main goal is to show that we can get a $O(\mathsf{poly}(p,\logv{\frac{1}{\eps}})\min\inbraces{\rank{\mA},m}^{1/3})$ iteration complexity.
\end{remark}

\subsection{Prior results, connections, and open problems}
Here, we discuss prior work that conceptually and technically relates to ours. We then suggest natural directions for future work.

\paragraph{Multi-distribution learning.}
Many learning problems involve multiple data sources, for instance, when multiple agents generate their data independently. One can formulate these multi-distribution problems as standard learning or optimization problems by considering a mixture of their distributions, as in objective \eqref{eq:utilitarian}. However, this approach often biases solutions toward dominant data sources, leading to poor performance on outliers—an issue stemming from statistical heterogeneity. This limitation motivates the study of multi-objective optimization problems~\cite{miettinen1999nonlinear, ehrgott2005multicriteria}, where each agent \( m \) has a distribution \( \mathcal{D}_m \) that defines its objective as \( \mathbb{E}_{z\sim \mathcal{D}_m}[f(\vx_m;z)] \), and where models \( \vx_m \) can vary across agents---a framework known as personalization. 

One of the earliest algorithms for such problems was introduced by \citet{blum2017collaborative}, where each agent's objective must be minimized to a pre-specified threshold \( \epsilon \) with high probability, framed within a PAC learning framework~\cite{valiant84,vapnik2013nature}. Subsequent research has refined these algorithms, achieving optimal sample complexity guarantees for learning from multiple distributions~\cite{chen2018tight,nguyen2018improved,hanneke2019value,haghtalab2022demand,zhang2024optimal}. Our objectives \( \eqref{eq:main_objective} \) and \( \eqref{eq:interpolation} \) offer different approaches to multi-distribution learning, where data distributions correspond to empirical agent distributions. In particular, \citet{mohri2019agnostic} analyzed objective \( \eqref{eq:main_objective} \) to establish generalization bounds for unknown mixtures of agents' distributions. 

Beyond sample efficiency, researchers have also examined other challenges, such as communication costs in large-scale distributed optimization~\cite{mmrha17}. A particularly relevant study is that of \citet{bullins2021stochastic}, which employs an efficient distributed quadratic sub-solver~\cite{woodworth2020local,patel2024limits} to implement an inexact Newton method for optimizing quasi-self-concordant functions (see \Cref{defn:qsc}).

\paragraph{Group fairness.} 
Recently, interest in algorithmic fairness has intensified~\cite{barocas2016big, abebe2020roles, kasy2021fairness} with researchers exploring fairness across various domains, including supervised learning~\cite{calders2009building, dwork2012fairness, hardt2016equality, kusner2017counterfactual, goel2018non, ustun2019fairness}, resource allocation~\cite{bertsimas2011price, bertsimas2012efficiency, hooker2012combining, donahue2020fairness, manshadi2021fair}, scheduling~\cite{mulvany2021fair}, online matching~\cite{chierichetti2019matroids, ma2023fairness}, assortment planning~\cite{singh2018fairness, biega2018equity, singh2019policy, chen2022fair}, and facility location~\cite{gupta2022socially}. The extensive literature on algorithmic fairness falls into three main categories: (1) individual fairness, which ensures that similar individuals receive comparable predictions~\cite{dwork2012fairness, loi2019philosophical, chen2022fair}, (2) group fairness, which aims for equal treatment of different demographic groups, often in terms of resource allocation or performance parity~\cite{singh2018fairness, balseiro2021regularized}, and (3) subgroup fairness, which blends aspects of both individual and group fairness~\cite{kearns2018preventing, kearns2019empirical}. 

This chapter focuses on a well-studied group fairness notion in machine learning literature: the group DRO problem~\cite{ben2013robust, duchi2016statistics, sagawa2019distributionally}. The idea of interpolating between robustness and utility is also common~\cite{golrezaei2024online} and closely related to multi-objective optimization, where scalarization~\cite{miettinen1999nonlinear, ehrgott2005multicriteria} helps recover desired solutions along the Pareto frontier. 

\paragraph{Linear programming and $\ell_p$ regression.} 
In the last several years, there has been a surge of work in obtaining second-order, condition-free algorithms for linear programming and $\ell_p$ regression \cite{bcll18,ls19,akps19,jls21}. Observe that $\ell_p$ regression is a special case of the problem we study in objective \eqref{eq:interpolation}, which is recovered when all $n_i=1$, and $\ell_{\infty}$ regression is captured by linear programming. Note that neither of these problem families is expressive enough to capture the objectives we study. In general, to get iteration complexities in the smaller of the two dimensions for these problems, it seems that a geometric understanding of the solution space is required -- these ideas were central to the improvements obtained by \cite{ls19,jls21} and in our work.

\paragraph{Open problems.} 
Our work raises several open questions. One limitation of \Cref{mainthm:fair_regression_iteration_complexity} is that its iteration complexity is not high-accuracy, meaning its dependence on \( \eps \) is not \(\mathrm{polylog}(1/\eps)\). Designing a high-accuracy solver under the same conditions as \Cref{mainthm:gp_regression_iteration_complexity} with iteration complexity \(\widetilde{O}\inparen{\mathsf{poly}(\min\inbraces{\rank{\mA}, m}, \logv{\frac{1}{\eps}})}\) remains an open problem.

A more ambitious general goal is to design algorithms for convex quadratic programs with the aforementioned iteration complexity. This would generalize analogous results for linear programming \cite{ls19}. We view the current work as a first step towards this goal, as the objective \eqref{eq:main_objective} is a structured convex quadratic program for which we get an iteration complexity independent of $m$. It would also be interesting to consider other complexity measures beyond $\rank{\mA}$, for instance, assumptions about the ground-truth labeling vector $\vx_i^\star$ for each group's data $S_i$.

Finally, our results suggest that optimizing for ``$\ell_p$-interpolants'' between non-robust and robust objectives may be computationally easier than optimizing for the robust objective alone. A more precise statistical characterization of how robustness and utility trade-off as \( p \) varies in collaborative, fair, or multi-distributional learning settings would be valuable. Additionally, exploring interpolations or solution concepts along the Pareto frontier of the \( m \)-dimensional multi-objective optimization problem could yield further insights.

\subsection{Chapter outline}

In the rest of this chapter we will outline the key details of our approach as well as give a proof outline for our theoretical results. In \Cref{sec:gp_reg_overview}, we give proof sketches of our main results. In \Cref{sec:mirror_descent}, we give an analysis of mirror descent under inexact subproblem solves -- we will need this in the proof of \Cref{mainthm:gp_regression_iteration_complexity}. In \Cref{sec:optimal_ms_acceleration}, we modify an acceleration scheme due to \cite{chjjs22}, which we will use to iterate calls to the proximal subproblem solver \eqref{eq:gp_prox_intro_reg} for the proof of \Cref{mainthm:gp_regression_iteration_complexity}. In \Cref{sec:gp_robust_proof}, we prove \Cref{mainthm:fair_regression_iteration_complexity}. In \Cref{sec:gp_interpolating_proof}, we prove \Cref{mainthm:gp_regression_iteration_complexity}. 

\def\arraystretch{1.5}
\begin{table}[h]
\centering
\begin{tabular}{l c c} 
\toprule\toprule
\textbf{Algorithm}                                                                                                     & \textbf{Iteration Complexity}                                                                              & \textbf{Each Iteration}                     \\ \midrule\midrule
\makecell[l]{Subgradient descent}                                                                                           & \makecell[c]{$\frac{\norm{\xstar}_2 \max_{1 \le i \le m} \tfrac{1}{\sqrt{n_i}}\opnorm{\mA_{S_i}}}{\eps^{2}}$}                                   & \makecell[c]{Evaluate $\nabla f(\vx)$}                    \\ \midrule
\makecell[l]{Nesterov acceleration\\ on smoothened objective}                                                                 & \makecell[c]{$\frac{\norm{\xstar}_2 \left(\max_{1 \le i \le m} \tfrac{1}{\sqrt{n_i}}\opnorm{\mA_{S_i}}\right)^{1/2}}{\eps}$}                      & \makecell[c]{Evaluate $\nabla \ftilde_{\beta,\delta}(\vx)$} \\ \midrule
\makecell[l]{\cite{aakmrz22}}                                                                                               & \makecell[c]{$\frac{\norm{\xstar}_2 \max_{1 \le i \le m} \tfrac{1}{\sqrt{n_i}}\opnorm{\mA_{S_i}}}{\eps}$}                                       & \makecell[c]{Evaluate $\nabla \ftilde_{\beta,\delta}(\vx)$} \\ \midrule
\makecell[l]{Interior point with\\ log barrier \cite{boyd2004convex}}                                                         & \makecell[c]{$m^{1/2} \log\left(\frac{1}{\eps}\right)$}                                                                    & \makecell[c]{Linear system solve in $\mA^{\top}\mB\mA$}    \\ \midrule
\makecell[l]{\textbf{This chapter} \\ (na\"ive geometry)}                                                                                          & \makecell[c]{$\frac{m^{1/3}}{\eps^{2/3}}$}                                                                                & \makecell[c]{Linear system solve in $\mA^{\top}\mB\mA$}    \\ \midrule
\makecell[l]{$\ell_{\infty}$ regression with\\ Lewis weights \cite{jls21}}                                                    & \makecell[c]{$\frac{\rank{\mA}^{1/3}}{\eps^{2/3}}$}                                                                       & \makecell[c]{Linear system solve in $\mA^{\top}\mD\mA$}    \\ \midrule
\makecell[l]{$\ell_{\infty}$ regression\\ with IPM \cite{ls19}}                                                               & \makecell[c]{$\rank{\mA}^{1/2} \log\left(\frac{1}{\eps}\right)$}                                                           & \makecell[c]{Linear system solve in $\mA^{\top}\mD\mA$}    \\ \midrule
\makecell[l]{\textbf{This chapter (\Cref{mainthm:fair_regression_iteration_complexity})}}                                                                                                    & \makecell[c]{$\frac{\min\inbraces{\rank{\mA},m}^{1/3}}{\eps^{2/3}}$}                                                               & \makecell[c]{Linear system solve in $\mA^{\top}\mB\mA$}    \\ \bottomrule\bottomrule
\end{tabular}
\caption{Here, we list the complexities of algorithms for optimizing \eqref{eq:main_objective} or for the special case of $\ell_{\infty}$ regression, assuming $\opt=1$ (the first three guarantees are additive approximations) and ignoring $\mathrm{polylog}(n,m)$ terms. We write $\mD$ to be a diagonal matrix and $\mB$ to be a block-diagonal matrix where each block has size $(n_i + o(1)) \times (n_i + o(1))$. We remark that in the special case where $n_i=1$, our algorithm exactly recovers that of \cite{jls21}.\label{table:compare}}
\end{table}

\section{Technical overview}\label{sec:gp_reg_overview}
In this section, we sketch our proofs for \Cref{mainthm:fair_regression_iteration_complexity} and \Cref{mainthm:gp_regression_iteration_complexity}. 

\paragraph{Notation.} Here and in the rest of the chapter, we ignore the dataset size normalization factors $1/\sqrt{n_i}$ as we can fold this into $\mA_{S_i}$ and $\vb_{S_i}$. Additionally, let $f(\vx) \coloneqq \sum_{i=1}^m \norm{\mA_{S_i}\vx-\vb_{S_i}}_2^p$ if $2 \le p < \infty$ and let $f(\vx) \coloneqq \max_{1 \le i \le m} \norm{\mA_{S_i}\vx-\vb_{S_i}}_2$ if $p=\infty$. Note that in the $2 \le p < \infty$ case, we let $f(\vx)$ be the $p$th power of the objective written in \Cref{mainthm:gp_regression_iteration_complexity}; this is to make future calculations easier and makes a difference of only polynomial factors in $p$ in the iteration complexity. Without loss of generality (by rescaling), let $\opt = 1$, where $\opt \coloneqq f(\xstar)$. So, it is enough to get an $\eps$-additive optimal solution $\xhat$. Also without loss of generality, let $\mA$ be such that $\rank{\mA}=d$. For a positive semidefinite $\mM\in\R^{d\times d}$, denote $\norm{\vx}_{\mM} \coloneqq \sqrt{\vx^{\top}\mM\vx}$. As shorthand, for $\vy\in\R^n$, we will often refer to the norm $\gnorm{\vy}{p} \coloneqq \inparen{\sum_{i=1}^m \norm{\vy_{S_i}}_2^p}^{1/p}$ for $p \ge 1$, where with a slight abuse of notation $\vy_{S_i}$ denotes the coordinates of $\vy$ indexed by $S_i$. Finally, in an abuse of notation, for symmetric matrices $\mM$, let $\mM^{-1}$ denote the pseudoinverse of $\mM$.

With this notation in mind, we note that many iterative methods for convex optimization can be seen as decomposing a complex problem into a series of simpler subproblems. Our algorithms for distributionally robust linear regression follow this pattern, where the simple subproblem resembles
\begin{align}
    \cO(\vq) \coloneqq \min_{\norm{\vx-\vq}_{\mM} \le r_{\vq}}\quad f(\vx)\enspace,\label{eq:gp_prox_intro}
\end{align}
for some positive semidefinite $\mM$ and for some ball radius $r_{\vq}$ which may depend on the query $\vq$. Sub-routines like \eqref{eq:gp_prox_intro} are central to many trust-region methods~\cite{conn2000trust}, and, importantly when $f$ is the sum of a linear function and a self-concordant barrier, interior point methods derived from the self-concordant barrier framework \footnote{In this case, the matrix $\mM$ is given by the Hessian of the barrier function evaluated at the subproblem's solution.} \cite{nn94}.

With such a subproblem structure in hand, three questions arise. \begin{enumerate*}[label={(\arabic*)}] \item How do we choose the ``local geometry'' $\mM$? \item How do we solve the subproblems efficiently? \item How do we combine our subproblem solutions to arrive at our final answer? \end{enumerate*} We address these concerns in order in the following discussion.

\subsection{The geometry of the proximal subproblems}
\label{sec:gp_reg_overview_blw}

Observe that when we solve \eqref{eq:gp_prox_intro}, we are solving an optimization problem over the sublevel sets $\inbraces{\vx\suchthat \norm{\vx}_{\mM}\le r_{\vq}}$ -- these are ellipsoids. Now, consider choosing the $\ell_2$ geometry that best approximates our loss function. Specifically, ignoring the offset $\vb$ for now, for a norm loss function $\norm{\cdot}$ and for some \textit{distortion} $\triangle \ge 1$ that is as close to $1$ as possible, we want
\begin{align*}
    \text{for all } \vx \in \R^d:\quad\quad \norm{\vx}_{\mM} \le \norm{\mA\vx} \le \triangle\norm{\vx}_{\mM}\enspace.
\end{align*}
Observe that this captures the families of losses we study -- in particular, we can check that for $\vy\in\R^n$, the functions $\gnorm{\vy}{p}=\inparen{\sum_{i=1}^m \norm{\vy_{S_i}}_2^p}^{1/p}$ for $1 \le p \le \infty$ are norms. To see what kinds of distortion guarantees we can hope for, recall that we can use John's theorem (\Cref{thm:thesisintro_john}) as a benchmark. For convenience, we restate it below.

\begin{theorem}[John's theorem, \cite{john1948}]
For any symmetric convex body $K \subset \R^d$, let $\cE(K)$ denote the ellipsoid of maximum volume contained within $K$. Then, we have
\begin{align*}
    \cE(K) \subseteq K \subseteq \sqrt{d}\cdot \cE(K)\enspace.
\end{align*}
Moreover, the $\sqrt{d}$ is worst-case optimal (e.g. let $K$ be the unit $\ell_{\infty}$ ball).
\end{theorem}

It is easy to see that sublevel sets of norms, i.e., sets of the form $\inbraces{\vx\in\R^d\suchthat \norm{\vx}\le 1}$, are symmetric convex bodies. Hence, using John's theorem, we see that for our normed losses, there exists $\mM$ that achieves distortion $\triangle\le\sqrt{d}$. However, as written, this is only an existence result. To make this useful for us and actually find $\mM$, we need an algorithm to calculate John's ellipsoid for the level sets of our losses (or some other ellipsoid that gets an even better approximation factor). To this end, we repurpose and renotate an earlier result from \Cref{chapter:sparsifyingnorms} (also found in \cite{mo23}). It gives us an efficient algorithm to find this $\ell_2$ geometry for the loss families we consider.

\begin{theorem}[Combining Lemmas 5.6, 5.8, Equation (1.8) from {\cite{mo23}}]
\label{thm:gp_regression_blw_intro}
Let $p\ge 2$. There exists an algorithm that finds a positive diagonal matrix $\mW \in \R^{n \times n}$ such that for all $\vx\in\R^d$ and all $c\in\R$, we have
\begin{align*}
    \frac{\norm{\mW^{\frac{1}{2}-\frac{1}{p}}\inparen{\mA\vx-c\vb}}_2}{(2(\rank{\mA}+1))^{\frac{1}{2}-\frac{1}{p}}} \le \inparen{\sum_{i=1}^m \norm{\mA_{S_i}\vx-c\vb_{S_i}}_2^p}^{\frac{1}{p}} \le \norm{\mW^{\frac{1}{2}-\frac{1}{p}}\inparen{\mA\vx-c\vb}}_2\enspace.
\end{align*}
The algorithm runs in $O(\log m)$ linear system solves in matrices of the form $\mA^{\top}\mD\mA$ for positive diagonal matrices $\mD$.
\end{theorem}

The diagonal entries of matrix $\mW$ are rescaled versions \textit{block Lewis weights} that we discussed in the previous chapter. Recall that this is a generalization of Lewis weights, and both objects have been used previously for various matrix approximation problems \cite{blm89,mmwy21,jls22,jlls23,mo23}. Furthermore, Lewis weights are central to improvements in the iteration complexities for linear programming and vanilla $\ell_p$ regression \cite{ls19,jls21}.

Additionally, notice that the distortion of $O(\rank{\mA}^{1/2-1/p})$ guaranteed by \Cref{thm:gp_regression_blw_intro} is optimal. To see this, let $\mA \in \R^{n \times d}$ be such that for $i\in[d]$, row $\va_i=\ve_i$, where $\ve_i$ is the $i$th standard basis vector. Then, for all $d+1 \le i \le n$, let $\va_i=0$. In words, $\mA$ is the $d$-dimensional identity matrix atop a large matrix of all $0$s. It is easy to see that for any $p$, we have $\norm{\mA\vx}_p = \norm{\vx}_p$, and the best distortion we can get for relating $\norm{\vx}_p$ to any $d$-dimensional $\ell_2$ norm is $d^{\abs{1/2-1/p}}$.

With \Cref{thm:gp_regression_blw_intro} and its near optimality in hand, it is natural to choose $\mM = \mA^{\top}\mW^{1-\frac{2}{p}}\mA$ if $\rank{\mA} \le m$ and $\mM = \mA^{\top}\mA$ if $\rank{\mA} \ge m$ (in the latter case, we get a $\sqrt{m}$ distortion for free from relating $\ell_2^m$ to $\ell_{\infty}^m$). This gives us an $\ell_2$ geometry that nearly optimally approximates our losses. In the following sections, we will see how this helps us both implement our proximal subproblem solvers and combine these to solve the whole original problem.

\subsection{Solving proximal subproblems}

In this subsection, let $\mM$ be any positive semidefinite matrix, as the arguments here apply for any $\mM$. In particular, these arguments are required to analyze even the na\"ive geometry obtained by choosing $\mM = \mI$ or for any other ellipsoidal approximation given by any choice of $\mM$. For our best results, we will finally choose $\mM = \mI$ or $\mM = \mA^{\top}\mW^{1-\frac{2}{p}}\mA$ depending on which ellipsoidal approximation is better for our input.

Here, we discuss how to solve problems of the form \eqref{eq:gp_prox_intro} for a fixed query $\vq$. Our strategy follows two general steps. First, we establish some form of local stability for $\nabla^2 f(\vx)$ within the ball we are solving in, i.e., we want $\nabla^2 f(\vx)$ to not change too much inside the ball $\inbraces{\vx\in\R^d \suchthat \norm{\vx-\vq}_{\mM} \le r_{\vq}}$. Second, we use this to show that an appropriate second-order algorithm enjoys a good convergence rate to an approximate solution for our subproblem. We handle the $p=\infty$ and $2 \le p < \infty$ cases separately below.

\subsubsection{The robust case \texorpdfstring{($p=\infty$)}{(p=infinity)}.}
\label{sec:gp_reg_overview_robust_hessian}

Unfortunately, since $f$ is not even differentiable (it is the pointwise maximum of Euclidean norms, each of which is also not differentiable), we cannot directly argue about the stability of $\nabla^2 f(\vx)$. We therefore first need to find some surrogate objective $\ftilde$ so that:
\begin{enumerate}
    \item The approximation error $\maxnorm{\ftilde -f}$ is small;
    \item The surrogate objective $\ftilde$ is smooth in $\norm{\cdot}_{\mM}$ in such a way that we can solve the proximal subproblems fast.
\end{enumerate}
To smoothen $f(\vx)$, we use the family of objectives parameterized by $\beta,\delta$
\begin{align*}
    \ftilde_{\beta,\delta}(\vx) \coloneqq \beta\logv{\sum_{i=1}^m \expv{\frac{\sqrt{\delta^2 + \norm{\mA_{S_i}\vx-\vb_{S_i}}_2^2}-\delta}{\beta}}}\enspace.
\end{align*}
This can be seen as composing the softmax function with temperature $\beta$ with ``inner functions'' $\sqrt{\delta^2 + \norm{\mA_{S_i}\vx-\vb_{S_i}}_2^2} - \delta$. It is straightforward to show that for all $\vx\in\R^d$, $\abs{\ftilde_{\beta,\delta}(\vx) - f(\vx)} \le \beta\log m + \delta$. So, setting $\beta = \eps/4\log m$ and $\delta = \eps/4$, it is sufficient to optimize $\ftilde_{\beta,\delta}$ up to $\eps/2$ additive error to get an $\eps$-additive suboptimal solution to our original objective. Furthermore, we prove that $\ftilde_{\beta,\delta}$ is $O(1/\beta+1/\delta)$-smooth in the norm $\gnorm{\mA\vx}{\infty} \coloneqq \max_{1 \le i \le m} \norm{\mA\vx}_2$. From \Cref{thm:gp_regression_blw_intro}, this means that $\ftilde_{\beta,\delta}$ is also $O(1/\beta+1/\delta)$-smooth in the norm $\norm{\vx}_{\mM}$ where $\mM$ is chosen according to the previous subsection (notice that the only fact we need about $\mM$ here is that $\max_{1 \le i \le m} \norm{\mA\vx}_2 \le \norm{\vx}_{\mM}$).

Next, \citet{msbacon} show that if $\ftilde_{\beta,\delta}$ satisfies a higher-order smoothness condition called \textit{quasi-self-concordance} with respect to the norm $\norm{\cdot}_{\mM}$, then we can get the required Hessian stability for a \textit{fixed} $r_{\vq} = \Theta(1/\eps)$ (in particular, $r_{\vq}$ does not depend on $\vq$ here). To be more clear, we define quasi-self-concordance below.
\begin{definition}[Quasi-self-concordance, adapted from {\cite[Appendix A]{ksj18}}]
\label{defn:qsc}
Let $\myfunc{f}{\R^k}{\R}$. We say that $f$ is $\nu$-quasi-self-concordant in the norm $\norm{\cdot}$ if for all vectors $\vy\in\R^k$, directions $\vd\in\R^k$, and $t\in\R$, we have
\begin{align*}
    \abs{\inparen{\frac{d}{dt}}^3 f(\vy+t\vd)} \le \nu\norm{\vd}\inparen{\frac{d}{dt}}^2 f(\vy+t\vd).
\end{align*}
\end{definition}
Then, \cite{msbacon} shows how to leverage this Hessian stability to implement \eqref{eq:gp_prox_intro} with low linear-system-solve iteration complexity. However, previously, it was only shown that the composition of the softmax function with linear functions is quasi-self-concordant. So, it was unknown whether composing softmax with other functions could also be quasi-self-concordant. 

To resolve this, we prove a much more general composition result, which may be of independent interest. It essentially states that if we compose the softmax function with any combination of ``inner'' functions that are quasi-self-concordant, the resulting function is also quasi-self-concordant. For a more formal statement, see \Cref{lemma:composed_qsc}. 

Hence, to show the requisite Hessian stability, we use the following steps. We show that the ``inner'' functions $\sqrt{\delta^2+\norm{\mA_{S_i}\vx-\vb_{S_i}}_2^2}-\delta$ are each $O(1/\delta)$-quasi-self-concordant in the norm $\norm{\mA_{S_i}\vx}_2$. So, we can apply our composition result \Cref{lemma:composed_qsc} to prove that $\ftilde_{\beta,\delta}$ is $O(1/\beta + 1/\delta)$-quasi-self-concordant in the norm $\max_{i\in[m]}\norm{\mA_{S_i}\vx}_2$. Following from \Cref{thm:gp_regression_blw_intro}, $\ftilde$ is $O(1/\beta+1/\delta)$-quasi self concordant in $\norm{\vx}_{\mM}$ as well (and in particular when $\mM = \mA^{\top}\mW^{1-\frac{2}{p}}\mA$ and $\mM = \mA^{\top}\mA$). We then apply the recipe given in \cite{msbacon} and get our subproblem solver for the $p=\infty$ case.

\subsubsection{The interpolating case (\texorpdfstring{$2 \le p < \infty$}{2 <= p < infinity}).}

Instead of explicitly constraining $r_{\vq}$ like in the $p=\infty$ case, we regularize our movement from $\vq$ in the norm $\norm{\cdot}_{\mM}$. Specifically, the subproblem we solve for any query $\vq$ is
\begin{align}
    \argmin{\vx\in\R^d} f(\vx) + ep^p\norm{\vx-\vq}_{\mM}^p\enspace.\label{eq:gp_prox_intro_reg}
\end{align}
This is the natural generalization of the proximal problem that \cite{jls21} use to get their results for $\ell_p$ regression, and the outline of our solver for these subproblems is similar to what \cite{jls21} use for this special case (see their Section 4).

However, we go a step further and show how to obtain approximate stationary points to \eqref{eq:gp_prox_intro_reg} instead of just getting a small objective value. This is because the acceleration scheme we use to iterate subproblem solutions to get our final answer $\xhat$ requires us to obtain an approximate stationary point for \eqref{eq:gp_prox_intro_reg}. The main new technical tool we develop for this purpose is a form of strong convexity for functions of the form $\norm{\vy}_2^p$ for $\vy\in\R^k$ for any $k\ge 1$. See \Cref{lemma:strong_convexity_component}.

\begin{restatable*}[Strong convexity of $\norm{\vy}_2^p$]{lemma}{lemmastrongconvexitycomponent}
\label{lemma:strong_convexity_component}
Let $\vv \in \R^k$ for $k\ge 1$. For any $\triangle \in \R^k$, we have
\begin{align*}
    \norm{\vv+\triangle}_2^p \ge \norm{\vv}_2^p + p\norm{\vv}_2^{p-2}\ip{\vv,\triangle}+\frac{4}{2^p}\norm{\triangle}_2^p\enspace.
\end{align*}
\end{restatable*}

With \Cref{lemma:strong_convexity_component}, we can argue about the strong convexity of $\norm{\vx-\vq}_{\mM}^p$, which means that we can convert an approximately optimal solution to \eqref{eq:gp_prox_intro_reg} in function value to one that is approximately optimal in parameter space as well. We combine this with a local gradient Lipschitzness property of the objective \eqref{eq:gp_prox_intro_reg} to get our approximate stationary point, which is enough for our purposes. The local gradient Lipschitzness property itself follows from a form of Hessian stability that we show for the objective \eqref{eq:gp_prox_intro_reg}. See \Cref{lemma:gp_hessian_stable}.

Finally, to obtain an approximately optimal solution to \eqref{eq:gp_prox_intro_reg} in function value, we again apply the Hessian stability property to conclude that \eqref{eq:gp_prox_intro_reg} is relatively smooth and relatively strongly convex in a simpler reference function. We show how to solve optimization problems in this reference function up to an approximate optimality that is sufficient for the rest of our applications -- this requires a mild modification of the standard mirror descent analysis, and we do this in \Cref{sec:mirror_descent}. Combining all of these building blocks gives us our subproblem solver for the $2 \le p < \infty$ case.

\subsection{Iterating proximal calls}

We now discuss the last item. Recall that we think of $\cO(\vq)$ as answering a proximal problem for the query $\vq$. It is not hard to show that under reasonable conditions on $f$ and on the structure of the subproblems, we can iterate calls to $\cO(\vq)$ to optimize $f$ (see, e.g., \cite[Appendix A]{msbacon}). This conceptually simple approach will already give us condition-free, group-independent rates for the problems we study.

But, we can do better. An acceleration framework originally due to \citet{ms13} and generalized/refined in subsequent works \cite{bjlls19,msbacon,chjjs22} gives a recipe to iterate calls of $\cO(\vq)$ to optimize the original function $f$. From these, the iteration complexity we need for an $\eps$-additive solution with an initialization $\vx_0$ and optimum $\xstar$ is roughly $\inparen{{\norm{\vx_0-\xstar}_{\mM}}/{\eps}}^{2/3}$ (see \Cref{thm:optimal_ms_acceleration} for a more formal statement). Combining this with \Cref{thm:gp_regression_blw_intro}, which implies that we can find $\vx_0$ such that $\norm{\vx_0-\xstar}_{\mM} \le \sqrt{d}$, we see that we should expect rates of roughly $\widetilde{O}(d^{1/3}\eps^{-2/3})$ for our problems. Indeed, \Cref{mainthm:fair_regression_iteration_complexity} and \Cref{mainthm:gp_regression_iteration_complexity} attain rates that are at least this good up to logarithmic factors. In this step, we again use the strong convexity that we prove for our objective (\Cref{lemma:strong_convexity_component}) to show that after enough steps of this algorithm, the problem diameter will have noticeably shrunk. Iterating this gives the high-accuracy result of \Cref{mainthm:gp_regression_iteration_complexity}.

Interestingly, our algorithm for the $2 \le p < \infty$ case uses a form of this accelerated scheme developed in \cite{chjjs22} that does not require solving an implicit equation for the query point, improving over the results from \cite{jls21} for $\ell_p$ regression. It would be nice to obtain this for the $p=\infty$ case (in \Cref{sec:optimal_ms_acceleration}, we discuss a technical challenge in obtaining this).

\section{Mirror descent with inexact updates}
\label{sec:mirror_descent}

\paragraph{Notation warning.} This section is meant to be a self-contained, standalone analysis of mirror descent under inexact updates. The notation is chosen to be consistent with most material we could find on mirror descent and therefore conflicts with the notation used in the rest of the chapter.

In this section, we give an analysis of unconstrained mirror descent when each Bregman proximal problem is solved only approximately (\Cref{alg:mirror_descent}). Although we expect that this is a standard fact about mirror descent, we could not find an appropriate reference. Hence, we produce it here.

\begin{algorithm}[H]
\caption{\textsf{ApproximateMirrorDescent}: Implements mirror descent to optimize convex and differentiable $f$ given $L$-relative smoothness and $\mu$-relative strong convexity in the reference $h$ when we may not be able to solve each proximal problem exactly.}
\label{alg:mirror_descent}
\begin{algorithmic}[1]
\Require Initial point $\vx_0$, iteration count $T$.
\State Define \begin{equation*}
    \begin{aligned}
    D_{h}(\vx,\vy) &\coloneqq h(\vx)-h(\vy)-\ip{\nabla h(\vy),\vx-\vy} \\
    \xstar &\coloneqq \argmin{\vx\in\R^d} f(\vx)
    \end{aligned}.
\end{equation*}
\For{$i = 1, \dots, T$}
    \State $\xstar_i = \argmin{\xtilde \in \R^d} f(\vx_{i-1}) + \ip{\nabla f(\vx_{i-1}), \xtilde-\vx_{i-1}} + LD_h(\xtilde,\vx_{i-1})$ \label{line:mirror_descent_exact} \Comment{We may only be able to approximate $\xstar_i$ -- see the next line.}
    \State Let $\vx_i$ be an approximate stationary point for the above objective.\label{line:mirror_descent_approx}
\EndFor
\Return $\argmin{0 \le i \le T} f(\vx_i)$
\end{algorithmic}
\end{algorithm}

In \Cref{alg:mirror_descent}, we assume that the function $f$ is $\mu$-relatively strongly convex and $L$-smooth in a \textit{reference function} $h$. This means that for all $\vx,\vy\in\R^d$, we have
\begin{align*}
    \mu D_h(\vx,\vy) \le f(\vx) - f(\vy) - \ip{\nabla f(\vy), \vx-\vy} \le  LD_h(\vx,\vy).
\end{align*}
Using \cite[Proposition 1.1]{lfn18}, when $f$ is twice-differentiable, this condition is equivalent to asking for all $\vx\in\R^d$, 
\begin{align*}
    \mu \nabla^2 h(\vx) \preceq \nabla^2 f(\vx) \preceq L\nabla^2 h(\vx).
\end{align*}
We are now ready to state the performance guarantee of \Cref{alg:mirror_descent}. See \Cref{thm:gp_mirror_descent}.

\begin{theorem}
\label{thm:gp_mirror_descent}
Let index $j$ be the index output by \Cref{alg:mirror_descent}. Let $\triangle_i$ be defined such that
\begin{align*}
    \triangle_i \coloneqq \nabla f(\vx_{i-1}) + L\inparen{\nabla h(\vx_i) - \nabla h(\vx_{i-1})}.
\end{align*}
Then, we have
\begin{align*}
    f(\vx_j)-f(\xstar) \le L\inparen{1-\frac{\mu}{L}}^TD_h(\xstar,\vx_0)+\max_{1\le i \le n}\ip{\triangle_i,\vx_i-\xstar}.
\end{align*}
\end{theorem}

To prove \Cref{thm:gp_mirror_descent}, we begin with a few standard facts about the mirror descent iterations.

\begin{lemma}
\label{lemma:gp_md_threepoint}
Let $\vy\in\R^d$ be arbitrary. We have
\begin{align*}
    \ip{\nabla f(\vx_{i-1}), \vx_i-\vy} = L\inparen{D_h(\vy,\vx_{i-1})-D_h(\vy,\vx_i)-D_h(\vx_i,\vx_{i-1})} + \ip{\triangle_i, \vx_i-\vy}.
\end{align*}
\end{lemma}
\begin{proof}[Proof of \Cref{lemma:gp_md_threepoint}]
By the three point identity (see, e.g., \cite[Equation (A.9)]{syls16}), we have
\begin{align*}
    D_h(\vy,\vx_{i-1})-D_h(\vy,\vx_i)-D_h(\vx_i,\vx_{i-1})&=-\ip{\nabla h(\vx_i)-\nabla h(\vx_{i-1}),\vx_i-\vy} \\
    &= \frac{1}{L}\ip{\nabla f(\vx_{i-1})-\triangle_i,\vx_i-\vy},
\end{align*}
completing the proof of \Cref{lemma:gp_md_threepoint}.
\end{proof}
\begin{lemma}[Mirror descent lemma under approximate stationary point updates]
\label{lemma:gp_md_singlestep}
Let $\vy\in\R^d$ be arbitrary. For every iteration $i$, we have
\begin{align*}
    f(\vx_i)-f(\vy) \le (L-\mu)D_h(\vy,\vx_{i-1}) - LD_h(\vy,\vx_i) + \ip{\triangle_i,\vx_i-\vy}.
\end{align*}
\end{lemma}
\begin{proof}[Proof of \Cref{lemma:gp_md_singlestep}]
The definition of $\mu$-relative strong convexity tells us that
\begin{align*}
    f(\vx_{i-1})-f(\vy) \le \ip{\nabla f(\vx_{i-1}), \vx_{i-1}-\vy} - \mu D_h(\vy,\vx_{i-1}).
\end{align*}
We now write
\begin{align*}
    f(\vx_i)-f(\vy) &\le f(\vx_{i-1})-f(\vy) + \ip{\nabla f(\vx_{i-1}), \vx_i-\vx_{i-1}} + LD_h(\vx_i,\vx_{i-1}) & \text{($L$-RS)} \\
    &\le \ip{\nabla f(\vx_{i-1}), \vx_i-\vy}-\mu D_h(\vy,\vx_{i-1}) + LD_h(\vx_i,\vx_{i-1}) & \text{($\mu$-RSC)} \\
    &\le (L-\mu)D_h(\vy,\vx_{i-1}) - LD_h(\vy,\vx_i) + \ip{\triangle_i,\vx_i-\vy}, & \text{(\Cref{lemma:gp_md_threepoint})}
\end{align*}
completing the proof of \Cref{lemma:gp_md_singlestep}.
\end{proof}

We now have the tools to complete the proof of \Cref{thm:gp_mirror_descent}.

\begin{proof}[Proof of \Cref{thm:gp_mirror_descent}]
Let $E_i \coloneqq f(\vx_i)-f(\xstar)-\ip{\triangle_i,\vx_i-\xstar}$. Substituting $\vy=\xstar$ and rearranging the conclusion of \Cref{lemma:gp_md_singlestep} gives
\begin{align*}
    E_i \le (L-\mu)D_h(\xstar,\vx_{i-1})-LD_h(\xstar,\vx_i).
\end{align*}
We multiply both sides by $\inparen{\frac{L}{L-\mu}}^i$ and write
\begin{align*}
    \inparen{\frac{L}{L-\mu}}^iE_i \le \frac{L^i}{(L-\mu)^{i-1}}D_h(\xstar,\vx_{i-1})-\frac{L^{i+1}}{(L-\mu)^{i}}D_h(\xstar,\vx_i).
\end{align*}
Adding over all $T$ iterations yields
\begin{align*}
    \sum_{i=1}^T \inparen{\frac{L}{L-\mu}}^iE_i \le LD_h(\xstar,\vx_0)-\inparen{\frac{L}{L-\mu}}^TLD_h(\xstar,\vx_T) \le LD_h(\xstar,\vx_0).
\end{align*}
Expanding out the definition of $E_i$ and rearranging gives
\begin{align*}
    \sum_{i=1}^T \inparen{\frac{L}{L-\mu}}^i(f(\vx_i)-f(\xstar)) \le LD_h(\xstar,\vx_0) + \sum_{i=1}^T \inparen{\frac{L}{L-\mu}}^i\ip{\triangle_i,\vx_i-\xstar}.
\end{align*}
By the geometric series summation formula, we define and have
\begin{align*}
    C_T \coloneqq \sum_{i=1}^T \inparen{\frac{L}{L-\mu}}^i = \frac{L}{\mu}\inparen{\inparen{1+\frac{\mu}{L-\mu}}^T-1}.
\end{align*}
Let $j$ be the index that \Cref{alg:mirror_descent} returns. It is easy to check that
\begin{align*}
    \sum_{i=1}^T \inparen{\frac{L}{L-\mu}}^i(f(\vx_i)-f(\xstar)) \ge C_T\inparen{f(\vx_j)-f(\xstar)}
\end{align*}
and
\begin{align*}
    \sum_{i=1}^T \inparen{\frac{L}{L-\mu}}^i\ip{\triangle_i,\vx_i-\xstar} \le C_T\max_{1 \le i \le n}\ip{\triangle_i,\vx_i-\xstar}.
\end{align*}
This gives us
\begin{align*}
    f(\vx_j)-f(\xstar) \le \frac{L}{C_T}D_h(\xstar,\vx_0)+\max_{1 \le i \le n}\ip{\triangle_i,\vx_i-\xstar}.
\end{align*}
Finally, notice that
\begin{align*}
    \frac{L}{C_T} = \frac{\mu}{\inparen{1+\frac{\mu}{L-\mu}}^T-1} \le L\inparen{1-\frac{\mu}{L}}^T.
\end{align*}
Combining everything completes the proof of \Cref{thm:gp_mirror_descent}.
\end{proof}

Finally, we add another useful lemma that quantifies the descent, if any, in the objective value between iterations.

\begin{lemma}
\label{lemma:gp_md_descent}
For every iteration $i$, we have
\begin{align*}
    f(\vx_i)-f(\vx_{i-1}) \le -LD_h(\vx_{i-1},\vx_i) + \ip{\triangle_i,\vx_i-\vx_{i-1}}.
\end{align*}
In particular, if $\ip{\triangle_i,\vx_{i}-\vx_{i-1}} \le LD_h(\vx_{i-1},\vx_i)$, then iteration $i$ is a descent step.
\end{lemma}
\begin{proof}[Proof of \Cref{lemma:gp_md_descent}]
We substitute $\vy=\vx_{i-1}$ in the conclusion of \Cref{lemma:gp_md_singlestep}. This gives
\begin{align*}
    f(\vx_i)-f(\vx_{i-1}) \le -LD_h(\vx_{i-1},\vx_i) + \ip{\triangle_i,\vx_i-\vx_{i-1}},
\end{align*}
completing the proof of \Cref{lemma:gp_md_descent}.
\end{proof}

\section{Optimal MS acceleration under custom Euclidean geometry}
\label{sec:optimal_ms_acceleration}

In this section, we adapt the bisection-free Monteiro-Svaiter acceleration framework developed in \cite{chjjs22} to handle custom Euclidean geometries. The object of interest here is \Cref{alg:optimal_ms_acceleration}, which we will call with different choices of the oracle $\oms$ for our algorithms.

\begin{algorithm}[H]
\caption{\textsf{OptimalMSAcceleration}: optimizes function $f$ given MS oracle $\oms$.}
\label{alg:optimal_ms_acceleration}
\begin{algorithmic}[1]
\Require Initial $\vx_0$, function $f$, oracle $\oms$, initial $\lambda_0'$, multiplicative adjustment factor $\alpha >1$, iteration count $T$
\State Set $\vv_0 = \vx_0$, $A_0 = 0$, $A'_0 = 0$.
\State Set $\xtilde_1, \lambda_1 = \cO(\vx_0; \lambda_0')$ and $\lambda_1' = \lambda_1$.
\For{$t=0,\dots,T$}
    \State $a_{t+1}' = \frac{1}{2\lambda_{t+1}'}\inparen{1+\sqrt{1+4\lambda'_{t+1}A_t}}$
    \State $A_{t+1}' = A_t+a_{t+1}'$
    \State $\vq_t = \frac{A_t}{A_{t+1}'}\vx_t + \frac{a_{t+1}'}{A_{t+1}'}\vv_t$\label{line:optimal_ms_query}
    \If{$t > 0$} $\xtilde_{t+1},\lambda_{t+1}=\oms(\vq_t;\lambda_{t+1}')$ \label{line:optimal_ms_acceleration_query}\EndIf
    \State $\gamma_{t+1} = \min\inbraces{1, \frac{\lambda_{t+1}'}{\lambda_{t+1}}}$
    \State $a_{t+1} = \gamma_{t+1}a'_{t+1}$ and $A_{t+1} = A_{t} + a_{t+1}$ \Comment{$A_{t+1} = A_{t+1}'-(1-\gamma_{t+1})a'_{t+1}$}
    \State $\vx_{t+1} = \frac{(1-\gamma_{t+1})A_{t}}{A_{t+1}}\vx_t + \frac{\gamma_{t+1}A_{t+1}'}{A_{t+1}}\xtilde_{t+1}$
    \If{$\gamma_{t+1} = 1$}
        \State $\lambda_{t+2}' = \frac{1}{\alpha}\lambda'_{t+1}$
    \Else
        \State $\lambda_{t+1}' = \alpha\lambda_{t+1}'$
    \EndIf
    \State $\vv_{t+1} = \vv_{t}-a_{t+1}\mM^{-1}\nabla f(\xtilde_{t+1})$
\EndFor
\end{algorithmic}
\end{algorithm}

In order to state the performance guarantee of \Cref{alg:optimal_ms_acceleration}, we require the notions of an \textit{MS oracle} and a \textit{movement bound}. See \Cref{defn:ms_oracle} and \Cref{defn:movement_bound}.

\begin{definition}[MS oracle, generalization of {\cite[Definition 1]{chjjs22}}]
\label{defn:ms_oracle}
Let $\mM \in \R^{d\times d}$ be a positive semidefinite matrix. An oracle $\myfunc{\cO}{\R^d \times \R_{\ge 0}}{\R^d \times \R_{\ge 0}}$ is a $\sigma$-MS oracle for function $\myfunc{f}{\R^d}{\R}$ if for every $\vq\in\R^d$ and $\lambda' > 0$, the points $(\vx,\lambda)=\cO(\vq;\lambda')$ satisfy
\begin{align*}
    \norm{\vx-\vq+\frac{1}{\lambda}\mM^{-1}\nabla f(\vx)}_{\mM} \le \sigma\norm{\vx-\vq}_{\mM}.
\end{align*}
\end{definition}

\begin{definition}[Movement bound {\cite[Definition 2]{chjjs22}}]
\label{defn:movement_bound}
For a norm $\norm{\cdot}_{\mM}$ induced by positive semidefinite $\mM\in\R^{d\times d}$, numbers $s \ge 1, c, \lambda > 0$, and $\vx,\vy\in\R^d$, we say that $(\vx,\vy,\lambda)$ satisfies a $(s,c)$-movement bound if
\begin{align*}
    \norm{\vx-\vy}_{\mM} \ge \begin{cases} \inparen{\frac{\lambda}{c^s}}^{\frac{1}{s-1}} & \text{ if } s < \infty \\ \frac{1}{c} & \text{ if } s = \infty \end{cases}.
\end{align*}
\end{definition}

With these in hand, we are ready to state the convergence guarantee we get with \Cref{alg:optimal_ms_acceleration}. See \Cref{thm:optimal_ms_acceleration}.

\begin{theorem}[Modification of {\cite[Theorem 1]{chjjs22}}]
\label{thm:optimal_ms_acceleration}
Let $\myfunc{f}{\R^d}{\R}$ be convex and differentiable. Consider running \Cref{alg:optimal_ms_acceleration} with parameters $\alpha = \expv{3-\frac{2}{s+1}}$ and a $\sigma$-MS oracle with $0 \le \sigma < 0.99$ (\Cref{defn:ms_oracle}). Let $s\ge1$ and $c>0$ and suppose that for all $t$ such that $\lambda_t>\lambda'_t$ or $t=1$, the iterates $(\xtilde_t,\vq_{t-1},\lambda_t)$ satisfy an $(s,c)$-movement bound (\Cref{defn:movement_bound}). Let $C$ be a universal constant. For any iteration count $T$ satisfying
\begin{align*}
    T \ge C\begin{cases} s\inparen{\frac{c^s\norm{\vx_0-\xstar}_{\mM}^{s+1}}{\eps}}^{\frac{2}{3s+1}} & \text{ if } s < \infty \\ \inparen{c\norm{\vx_0-\xstar}_{\mM}}^{2/3}\logv{\frac{\lambda_1\norm{\vx_0-\xstar}_{\mM}^2}{\eps}} & \text{ if } s = \infty \end{cases},
\end{align*}
we have
\begin{align*}
    f(\vx_T)-f(\xstar) \le \eps.
\end{align*}
\end{theorem}

The proof of \Cref{thm:optimal_ms_acceleration} follows the same recipe as the proof of \cite[Theorem 1]{chjjs22}. The only modification needed is that stated in \Cref{lemma:ms_acceleration_potential}.

\begin{lemma}[Replaces {\cite[Proposition 1]{chjjs22}}]
\label{lemma:ms_acceleration_potential}
In the context of \Cref{thm:optimal_ms_acceleration}, let $E_t \coloneqq f(\vx_t) - f(\xstar), D_t \coloneqq \frac{1}{2}\norm{\vv_t-\xstar}_{\mM}^2, N_{t+1} \coloneqq \frac{1}{2}\norm{\xtilde_{t+1}-\vq_t}_{\mM}^2$. Then, for all $t \ge 0$, we have
\begin{align*}
    A_{t+1}E_{t+1}+D_{t+1}+(1-\sigma^2)A_{t+1}'\min\inbraces{\lambda_{t+1},\lambda_{t+1}'}N_{t+1} \le A_tE_t+D_t.
\end{align*}
Consequently, for all $T\ge1, \sqrt{A_T}\ge\frac{1}{2}\sum_{t \in \cS_{T}^{\le}} \frac{1}{\sqrt{\lambda_t'}}$,
\begin{align*}
    E_{T} \le \frac{D_0}{A_T} \quad \text{ and } \quad (1-\sigma^2)\sum_{t \in \cS_{T}^{\ge}} A_t\lambda_t'N_t \le D_0 - A_TE_T.
\end{align*}
\end{lemma}
\begin{proof}[Proof of \Cref{lemma:ms_acceleration_potential}]
This proof is a straightforward modification of \cite[Proposition 1]{chjjs22}. We have
\begin{align*}
    D_{t+1} &= \frac{1}{2}\norm{\vv_{t+1}-\xstar}_{\mM}^2 = \frac{1}{2}\norm{\vv_t-a_{t+1}\mM^{-1}\nabla f(\xtilde_{t+1})-\xstar}_{\mM}^2 \\
    &= D_t + a_{t+1}\ip{\mM^{-1}\nabla f(\xtilde_{t+1}), \xstar-\vv_t}_{\mM} + \frac{a_{t+1}^2}{2}\norm{\mM^{-1}\nabla f(\xtilde_{t+1})}_{\mM}^2.
\end{align*}
By definition of $\vq_t$ and $A_{t+1}' = A_t + a'_{t+1}$, we have
\begin{align*}
    a_{t+1}'\vv_t = A_{t+1}'\vq_t - A_t\vx_t = a_{t+1}'\xtilde_{t+1} + A'_{t+1}\inparen{\vq_t-\xtilde_{t+1}} - A_t\inparen{\vx_t-\xtilde_{t+1}}.
\end{align*}
Subtracting $a_{t+1}'\xstar$ and taking the inner product with $\mM^{-1}\nabla f(\xtilde_{t+1})$ gives
\begin{align*}
    &\quad a_{t+1}'\ip{\mM^{-1}\nabla f(\xtilde_{t+1}), \xstar-\vv_t}_{\mM} \\
    &= \ip{\mM^{-1}\nabla f(\xtilde_{t+1}),  a_{t+1}'(\xstar-\xtilde_{t+1}) + A'_{t+1}\inparen{\xtilde_{t+1}-\vq_t} + A_t\inparen{\vx_t-\xtilde_{t+1}}}_{\mM} \\
    &\le a_{t+1}'\inparen{f(\xstar)-f(\xtilde_{t+1})}+A_{t+1}'\ip{\mM^{-1}\nabla f(\xtilde_{t+1}),\xtilde_{t+1}-\vq_t}_{\mM} + A_t\inparen{f(\vx_t)-f(\xtilde_{t+1})} \\
    &\le A_tE_t-A_{t+1}'\inparen{f(\xtilde_{t+1})-f(\xstar)}+A_{t+1}'\ip{\mM^{-1}\nabla f(\xtilde_{t+1}),\xtilde_{t+1}-\vq_t}_{\mM}.
\end{align*}
Rearranging gives
\begin{align*}
    A_{t+1}'\inparen{f(\xtilde_{t+1})-f(\xstar)} &\le A_tE_t+a_{t+1}'\ip{\mM^{-1}\nabla f(\xtilde_{t+1}), \vv_t-\xstar}_{\mM}\\
    &\quad +A_{t+1}'\ip{\mM^{-1}\nabla f(\xtilde_{t+1}),\xtilde_{t+1}-\vq_t}_{\mM}.
\end{align*}
Next, recall that by \Cref{defn:ms_oracle}, we have
\begin{align*}
    \norm{\mM^{-1}\nabla f(\xtilde_{t+1})+\lambda_{t+1}\inparen{\xtilde_{t+1}-\vq_t}}_{\mM}^2 \le \lambda_{t+1}^2\sigma^2\norm{\xtilde_{t+1}-\vq_t}_{\mM}^2.
\end{align*}
We use this to write
\begin{align*}
    &\quad \lambda_{t+1}\ip{\mM^{-1}\nabla f(\xtilde_{t+1}),\xtilde_{t+1}-\vq_t}_{\mM} \\
    &= \frac{1}{2}\norm{\mM^{-1}\nabla f(\xtilde_{t+1})+\lambda_{t+1}(\xtilde_{t+1}-\vq_t)}_{\mM}^2 - \frac{1}{2}\norm{\mM^{-1}\nabla f(\xtilde_{t+1})}_{\mM}^2-\frac{\lambda_{t+1}^2}{2}\norm{\xtilde_{t+1}-\vq_t}_{\mM}^2 \\
    &\le -\lambda_{t+1}^2(1-\sigma^2)N_{t+1}-\frac{1}{2}\norm{\mM^{-1}\nabla f(\xtilde_{t+1})}_{\mM}^2,
\end{align*}
from which we conclude
\begin{align*}
    \ip{\mM^{-1}\nabla f(\xtilde_{t+1}),\xtilde_{t+1}-\vq_t}_{\mM} &\le -\lambda_{t+1}(1-\sigma^2)N_{t+1}-\frac{1}{2\lambda_{t+1}}\norm{\mM^{-1}\nabla f(\xtilde_{t+1})}_{\mM}^2.
\end{align*}
Substituting back gives
\begin{align*}
    A_{t+1}'\inparen{f(\xtilde_{t+1})-f(\xstar)} &\le A_tE_t+a_{t+1}'\ip{\mM^{-1}\nabla f(\xtilde_{t+1}), \vv_t-\xstar}_{\mM}\\
    &\quad +A_{t+1}'\ip{\mM^{-1}\nabla f(\xtilde_{t+1}),\xtilde_{t+1}-\vq_t}_{\mM} \\
    &\le A_tE_t+a_{t+1}'\ip{\mM^{-1}\nabla f(\xtilde_{t+1}), \vv_t-\xstar}_{\mM}\\
    &\quad -A_{t+1}'\lambda_{t+1}(1-\sigma^2)N_{t+1}-\frac{A_{t+1}'}{2\lambda_{t+1}}\norm{\mM^{-1}\nabla f(\xtilde_{t+1})}_{\mM}^2.
\end{align*}
Next, recall that $\gamma_{t+1}a'_{t+1}=a_{t+1}$ and $\gamma_{t+1}\lambda_{t+1}=\min\inbraces{\lambda_{t+1},\lambda_{t+1}'}$, by construction. Let $\widehat{\lambda}_{t+1} \coloneqq \min\inbraces{\lambda_{t+1},\lambda_{t+1}'}$ We multiply both sides by $\gamma_{t+1}$ and conclude
\begin{align*}
    \gamma_{t+1}A_{t+1}'\inparen{f(\xtilde_{t+1})-f(\xstar)}&\le \gamma_{t+1}A_tE_t+a_{t+1}\ip{\mM^{-1}\nabla f(\xtilde_{t+1}), \vv_t-\xstar}_{\mM}\\
    &\quad -A_{t+1}'\widehat{\lambda}_{t+1}(1-\sigma^2)N_{t+1}-\frac{\gamma_{t+1}A_{t+1}'}{2\lambda_{t+1}}\norm{\mM^{-1}\nabla f(\xtilde_{t+1})}_{\mM}^2. 
\end{align*}
Now, by convexity of $f$ and from the definition of $\vx_{t+1}$, we have
\begin{align*}
    f(\vx_{t+1})-f(\xstar) \le \frac{(1-\gamma_{t+1})A_t}{A_{t+1}}\inparen{f(\vx_{t})-f(\xstar)}+\frac{\gamma_{t+1}A_{t+1}'}{A_{t+1}}\inparen{f(\xtilde_{t+1})-f(\xstar)}.
\end{align*}
Recall the definition of $E_t$, multiply both sides by $A_{t+1}$, apply our bound on $\gamma_{t+1}A_{t+1}'\inparen{f(\xtilde_{t+1})-f(\xstar)}$, and we get
\begin{align*}
    A_{t+1}E_{t+1} &\le (1-\gamma_{t+1})A_tE_t + \gamma_{t+1}A_{t+1}'\inparen{f(\xtilde_{t+1})-f(\xstar)} \\
    &\le A_tE_t + a_{t+1}\ip{\mM^{-1}\nabla f(\xtilde_{t+1}), \vv_t-\xstar}_{\mM}\\
    &\quad -A_{t+1}'\widehat{\lambda}_{t+1}(1-\sigma^2)N_{t+1}-\frac{\gamma_{t+1}A_{t+1}'}{2\lambda_{t+1}}\norm{\mM^{-1}\nabla f(\xtilde_{t+1})}_{\mM}^2
\end{align*}
After shifting terms around, we see that it remains to show
\begin{align*}
    a_{t+1}\ip{\mM^{-1}\nabla f(\xtilde_{t+1}), \vv_t-\xstar}_{\mM}-\frac{\gamma_{t+1}A_{t+1}'}{2\lambda_{t+1}}\norm{\mM^{-1}\nabla f(\xtilde_{t+1})}_{\mM}^2 \overset{?}{\le} D_t-D_{t+1}.
\end{align*}
In fact, by the choice of $a'_{t+1}$ and the definition of $A_{t+1}'$, we have
\begin{align*}
    \lambda'_{t+1}(a'_{t+1})^2=a'_{t+1}+A_{t}=A_{t+1}'.
\end{align*}
Multiply both sides by $\gamma_{t+1}^2/(2\lambda'_{t+1})$ and we get
\begin{align*}
    \frac{a_{t+1}^2}{2}=\frac{\gamma_{t+1}^2A_{t+1}'}{2\lambda'_{t+1}} = \frac{\min\inbraces{1,\frac{\lambda_{t+1}'}{\lambda_{t+1}}}\gamma_{t+1}A_{t+1}'}{2\lambda'_{t+1}} \le \frac{\gamma_{t+1}A_{t+1}'}{2\lambda_{t+1}}.
\end{align*}
We recycle an earlier computation and know that
\begin{align*}
    D_t-D_{t+1} &= a_{t+1}\ip{\mM^{-1}\nabla f(\xtilde_{t+1}), \vv_t-\xstar}_{\mM} - \frac{a_{t+1}^2}{2}\norm{\mM^{-1}\nabla f(\xtilde_{t+1})}_{\mM}^2 \\
    &\ge a_{t+1}\ip{\mM^{-1}\nabla f(\xtilde_{t+1}), \vv_t-\xstar}_{\mM}-\frac{\gamma_{t+1}A_{t+1}'}{2\lambda_{t+1}}\norm{\mM^{-1}\nabla f(\xtilde_{t+1})}_{\mM}^2,
\end{align*}
which completes the proof of the potential decrease.

The remaining statements follow as written in \cite[Proof of Proposition 1]{chjjs22}, and we conclude the proof of \Cref{lemma:ms_acceleration_potential}.
\end{proof}

Now that we have shown \Cref{lemma:ms_acceleration_potential}, we refer the reader to \cite[Appendix A]{chjjs22} for the proof of \Cref{thm:optimal_ms_acceleration}, as it now follows exactly as written there.

We also give additional bounds on the movement of the iterates in $\norm{\cdot}_{\mM}$, which is a straightforward adaptation of \cite[Lemma 31]{msbacon} to the improved framework from \cite{chjjs22}.

\begin{lemma}
\label{lemma:ms_acceleration_iterate_diameter}
For all $t\ge 1$, we have both
\begin{equation*}
    \begin{aligned}
        \norm{\vv_t-\xstar}_{\mM} &\le \sqrt{2}\norm{\vx_0-\xstar}_{\mM} \\
        \norm{\vx_t-\xstar}_{\mM} &\le \inparen{\sqrt{2}+\max_{1 \le i \le t} \frac{\lambda_i'}{\lambda_i} \cdot \sqrt{\frac{2}{1-\sigma^2}}}\norm{\vx_0-\xstar}_{\mM}
    \end{aligned}.
\end{equation*}
\end{lemma}

In the statement of \Cref{lemma:ms_acceleration_iterate_diameter}, the cost of overshooting the guess $\lambda_i'$ becomes evident -- without an additional strong convexity guarantee, it is challenging to ensure that each iterate remains in a small ball around $\xstar$. This is the main reason we are unable to apply the framework of \cite{chjjs22} to the $p=\infty$ case.

\begin{proof}[Proof of \Cref{lemma:ms_acceleration_iterate_diameter}]
Using the same notation as in \Cref{lemma:ms_acceleration_potential} and in that proof, we define
\begin{align*}
    P_t &\coloneqq A_tE_t+D_t \\
    \widehat{\lambda}_t &\coloneqq \min\inbraces{\lambda_t,\lambda_t'}.
\end{align*}
By induction on the conclusion of \Cref{lemma:ms_acceleration_potential}, for $t\ge 1$ we have
\begin{align*}
    \frac{1}{2}\norm{\vv_t-\xstar}_{\mM}^2 = D_t \le P_t + (1-\sigma^2)\sum_{k=1}^t A_k'\widehat{\lambda}_kN_k \le P_0 = \norm{\vx_0-\xstar}_{\mM}^2.
\end{align*}
Thus,
\begin{align*}
    \norm{\vv_t-\xstar}_{\mM} \le \sqrt{2}\norm{\vx_0-\xstar}_{\mM}.
\end{align*}
For the second conclusion, we introduce the following notation.
\begin{align*}
    \alpha_{t+1} &\coloneqq \frac{(1-\gamma_{t+1})A_t}{A_{t+1}} \\
    \beta_{t+1} &\coloneqq \frac{A_t}{A_{t+1}'} \\
    \delta_{t+1} &\coloneqq 1-(1-\alpha_{t+1})(1-\beta_{t+1}) = 1 - \frac{\gamma_{t+1}A'_{t+1}}{A_{t+1}} \cdot \frac{a'_{t+1}}{A'_{t+1}} = \frac{A_t}{A_{t+1}}
\end{align*}
We also establish for any $i$,
\begin{align*}
    \frac{\gamma_i A_i'}{\lambda_i a_i^2} = \frac{A_i'}{\lambda_i \gamma_i (a_i')^2} = \frac{1}{\gamma_i} \cdot \frac{\lambda_i'}{\lambda_i} = \max\inbraces{\frac{\lambda_i'}{\lambda_i}, 1},
\end{align*}
which implies
\begin{align*}
    \frac{\gamma_i A_i'}{\lambda_i} = a_i^2\max\inbraces{\frac{\lambda_i'}{\lambda_i}, 1}.
\end{align*}
Notice that
\begin{align*}
    \norm{\vx_{t+1}-\xstar}_{\mM} &\le \alpha_{t+1}\norm{\vx_t-\xstar}_{\mM} + (1-\alpha_{t+1})\norm{\xtilde_{t+1}-\xstar}_{\mM} \\
    &\le \alpha_{t+1}\norm{\vx_t-\xstar}_{\mM} + (1-\alpha_{t+1})\inparen{\norm{\vq_t-\xstar}_{\mM}+\norm{\xtilde_{t+1}-\vq_t}_{\mM}} \\
    &\le \alpha_{t+1}\norm{\vx_t-\xstar}_{\mM} \\
    &\quad + (1-\alpha_{t+1})\inparen{\beta_{t+1}\norm{\vx_t-\xstar}_{\mM}+(1-\beta_{t+1})\norm{\vv_t-\xstar}_{\mM}+\norm{\xtilde_{t+1}-\vq_t}_{\mM}} \\
    &= \inparen{\beta_{t+1}+\alpha_{t+1}-\alpha_{t+1}\beta_{t+1}}\norm{\vx_t-\xstar}_{\mM} \\
    &\quad + \inparen{1-\alpha_{t+1}}\inparen{1-\beta_{t+1}}\norm{\vv_t-\xstar}_{\mM} + \inparen{1-\alpha_{t+1}}\norm{\xtilde_{t+1}-\vq_t}_{\mM} \\
    &= \delta_{t+1}\norm{\vx_t-\xstar}_{\mM} + (1-\delta_{t+1})\norm{\vv_t-\xstar}_{\mM} + \inparen{1-\alpha_{t+1}}\norm{\xtilde_{t+1}-\vq_t}_{\mM} \\
    &\le \prod_{i=0}^t \delta_{i+1}\norm{\vx_0-\xstar}_{\mM} + \inparen{1-\prod_{i=0}^t \delta_{i+1}}\norm{\vv_t-\xstar}_{\mM} \\
    &\quad + \sum_{i=1}^{t+1} \prod_{j=i+1}^{t+1}\delta_{j}(1-\alpha_{i})\norm{\xtilde_{i}-\vq_{i-1}}_{\mM} \\
    &\le \sqrt{2}\norm{\vx_0-\xstar}_{\mM} + \sum_{i=1}^{t+1} \prod_{j=i+1}^{t+1}\delta_{j}(1-\alpha_{i})\norm{\xtilde_{i}-\vq_{i-1}}_{\mM} \\
    &= \sqrt{2}\norm{\vx_0-\xstar}_{\mM} + \sum_{i=1}^{t+1} \frac{A_i}{A_{t+1}}(1-\alpha_{i})\norm{\xtilde_{i}-\vq_{i-1}}_{\mM} \\
    &= \sqrt{2}\norm{\vx_0-\xstar}_{\mM} + \sum_{i=1}^{t+1} \frac{A_i}{A_{t+1}}\cdot\frac{\gamma_iA_i'}{A_i}\norm{\xtilde_{i}-\vq_{i-1}}_{\mM} \\
    &= \sqrt{2}\norm{\vx_0-\xstar}_{\mM} + \frac{1}{A_{t+1}}\sum_{i=1}^{t+1} \sqrt{\frac{\gamma_iA_i'}{\lambda_i}} \cdot \sqrt{\lambda_i\gamma_iA_i'}\norm{\xtilde_{i}-\vq_{i-1}}_{\mM} \\
    &\le \sqrt{2}\norm{\vx_0-\xstar}_{\mM} + \frac{\inparen{\sum_{i=1}^{t+1} \frac{\gamma_iA_i'}{\lambda_i}}^{1/2}}{A_{t+1}} \cdot \inparen{\sum_{i=1}^{t+1} \lambda_i\gamma_iA_i'\norm{\xtilde_{i}-\vq_{i-1}}_{\mM}^2}^{1/2} \\
    &\le \sqrt{2}\norm{\vx_0-\xstar}_{\mM} + \frac{\inparen{\sum_{i=1}^{t+1} \frac{\gamma_iA_i'}{\lambda_i}}^{1/2}}{A_{t+1}} \cdot \sqrt{\frac{2}{1-\sigma^2}}\norm{\vx_0-\xstar}_{\mM} \\
    &\le \sqrt{2}\norm{\vx_0-\xstar}_{\mM} + \frac{\sum_{i=1}^{t+1} a_i\max\inbraces{1,\frac{\lambda_i'}{\lambda_i}}}{A_{t+1}} \cdot \sqrt{\frac{2}{1-\sigma^2}}\norm{\vx_0-\xstar}_{\mM} \\
    &\le \sqrt{2}\norm{\vx_0-\xstar}_{\mM} + \max_{1 \le i \le t+1} \frac{\lambda_i'}{\lambda_i} \cdot \sqrt{\frac{2}{1-\sigma^2}}\norm{\vx_0-\xstar}_{\mM} \\
    &= \inparen{\sqrt{2}+\max_{1 \le i \le t+1} \frac{\lambda_i'}{\lambda_i} \cdot \sqrt{\frac{2}{1-\sigma^2}}}\norm{\vx_0-\xstar}_{\mM},
\end{align*}
completing the proof of \Cref{lemma:ms_acceleration_iterate_diameter}.
\end{proof}

\section{Minimizing the distributionally robust loss}\label{sec:gp_robust_proof}

The goal of this section is to prove \Cref{mainthm:fair_regression_iteration_complexity}. We break up the proof into parts as described in \Cref{sec:gp_reg_overview}. We structure the section as follows. In the rest of this subsection, we present \Cref{alg:fair_regression}, our algorithm that minimizes the distributionally robust loss. In \Cref{sec:gp_robust_approx}, we introduce our smooth approximation for the objective \eqref{eq:main_objective} and show that it is a good additive approximation (this is a standard argument, but we include it as it provides important intuition). 

As the main of the difficulty of the proof in \Cref{mainthm:fair_regression_iteration_complexity} is to establish a Hessian stability for our surrogate loss, we devote the bulk of this section to proving this. Recall that in \Cref{sec:gp_reg_overview_robust_hessian}, we claimed that a higher-order smoothness condition called \textit{quasi-self-concordance} gives us the needed Hessian stability -- in fact, this follows from \cite[Lemma 11]{msbacon}. In light of this, it is enough to prove that our surrogate loss is quasi-self-concordant. 

In \Cref{sec:lse_calculus}, we work out some calculus facts related to the softmax function. In particular, it is in \Cref{sec:lse_calculus} that we prove the general composition result \Cref{lemma:composed_qsc} that states that if we take the softmax of several quasi-self-concordant functions, then the resulting function is also quasi-self-concordant. In \Cref{sec:lse_calculus_special}, we apply this composition fact to prove that our surrogate objective is quasi-self-concordant. Finally, in \Cref{sec:gp_robust_analysis_combining}, we combine these building blocks with the acceleration framework in \cite{msbacon} and complete the proof of \Cref{mainthm:fair_regression_iteration_complexity}.

\begin{algorithm}[H]
\caption{\textsf{MinMaxRegression}: optimizes \eqref{eq:main_objective} to $(1+\eps)$-multiplicative error}
\label{alg:fair_regression}
\begin{algorithmic}[1]
\Require Regression problems $(\mA_{S_1},\vb_{S_1}),\dots,(\mA_{S_m},\vb_{S_m})$, accuracy $\eps > 0$
\State Using \cite[Algorithm 2]{mo23} with input $\insquare{\mA \vert \vb}$, find nonnegative diagonal $\mW$ and weights $w_1,\dots,w_m$ such that for all $j \in S_i$, $\mW[j][j] = w_i$ and for all $\vx\in\R^{d}$ and $c \in \R$,
\begin{align*}
    \gnorm{\mA\vx-c\vb}{\infty} \le \norm{\mW^{1/2}\mA\vx-c\mW^{1/2}\vb}_2 \le \sqrt{2(\rank{\mA    }+1)}\gnorm{\mA\vx-c\vb}{\infty}.
\end{align*}\label{line:blw}
\If{$\sum_{i=1}^m w_i \ge m$} \Comment{$\rank{\mA}+1\le\sum_{i=1}^m w_i \le 2(\rank{\mA}+1)$}
    \State Reset $\mW = \mI_n$.
\EndIf
\State Let $\vx_0 = \inparen{\mA^{\top}\mW\mA}^{-1}\mA^{\top}\mW\vb$. \Comment{$\vx_0 \coloneqq \argmin{\vx\in\R^d} \norm{\mW^{1/2}\mA\vx-\mW^{1/2}\vb}_2$.}
\State Let
\begin{align*}
    \ftilde_{\beta,\delta}(\vx) \coloneqq \beta\logv{\sum_{i=1}^m \expv{\frac{\sqrt{\delta^2 + \norm{\mA_{S_i}\vx-\vb_{S_i}}_2^2}-\delta}{\beta}}}
\end{align*}
where $\beta = \frac{\eps}{4\log m}$ and $\delta = \frac{\eps}{4}$. \Comment{A family of smoothenings of the objective.}
\State Let $\fhat(\vx) \coloneqq \ftilde_{\eps/4\log m, \eps/4}(\vx) + \frac{\eps}{1000\min\inbraces{\rank{\mA},m}}\norm{\mW^{1/2}\mA(\vx-\vx_0)}_2^2$.
\State Using \cite[Algorithm 3]{msbacon}, implement a $\inparen{\frac{C}{\min\inbraces{\rank{\mA},m}},\frac{C}{\eps}}$-ball optimization oracle for $\fhat$, where $C$ is a universal constant. \Comment{Iteration complexity guaranteed by \Cref{lemma:obj_smooth_qsc}}
\State Using \cite[Algorithm 2]{msbacon}, implement a $\frac{1}{2}$-MS oracle for $\fhat$.
\State Run \cite[Algorithm 1]{msbacon} for $\widetilde{O}\inparen{\frac{\min\inbraces{\rank{\mA},m}^{1/3}\logv{\frac{d}{\eps}}}{\eps^{2/3}}}$ iterations using the MS oracle from the previous line and with initial point $\vx_0$ and final point $\xhat$.
\State \Return $\xhat$
\end{algorithmic}
\end{algorithm}

\subsection{Smoothly approximating the objective}
\label{sec:gp_robust_approx}

Recall that for $\vy\in\R^n$, let $\gnorm{\vy}{\infty} \coloneqq \max_{1 \le i \le m} \norm{\vy_{S_i}}_2$, where for $\vy\in\R^n$ we let $\vy_{S_i}$ refer to the vector in $\R^{n_i}$ indexed by the indices in $S_i$. Also, for $\vy\in\R^m$, let $\lse_{\beta}(\vy)$ refer to the function
\begin{align*}
    \lse_{\beta}(\vy) \coloneqq \beta\logv{\sum_{i=1}^m \expv{\frac{y_i}{\beta}}}.
\end{align*}
At a high level, our algorithm will minimize the function
\begin{align*}
    \ftilde_{\beta,\delta}(\vx) \coloneqq \beta\logv{\sum_{i=1}^m \expv{\frac{\sqrt{\delta^2 + \norm{\mA_{S_i}\vx-\vb_{S_i}}_2^2}-\delta}{\beta}}}
\end{align*}
for appropriate choices of the parameters $\beta$ and $\delta$. This choice of smoothening is natural because of the following approximation statement -- see \Cref{lemma:smooth_approx}.
\begin{lemma}
\label{lemma:smooth_approx}
For all $\vx \in\R^d$, we have
\begin{align*}
    \abs{\ftilde_{\beta,\delta}(\vx)-\gnorm{\mA\vx-\vb}{\infty}} \le \beta\log m + \delta.
\end{align*}
\end{lemma}
\begin{proof}[Proof of \Cref{lemma:smooth_approx}]
These guarantees are well-known, but we prove them anyway for the sake of self-containment. We first prove that for any $\vv\in\R^m$, we have 
\begin{align*}
    \max_{1 \le i \le m} v_i \le \lse_{\beta}(\vv) \le \max_{1 \le i \le m} v_i +\beta\log m.
\end{align*}
In one direction, we have
\begin{align*}
    \lse_{\beta}(\vv) \le \beta\logv{\sum_{i=1}^m \expv{\frac{\max_{1\le i \le m} v_i}{\beta}}} = \beta\log m + \max_{1 \le i \le m} v_i,
\end{align*}
and in the other, we have
\begin{align*}
    \lse_{\beta}(\vv) \ge \beta\logv{\expv{\frac{\max_{1 \le i \le m} v_i}{\beta}}} = \max_{1 \le i \le m} v_i.
\end{align*}
Next, for $\vv \in \R^m$, we will show that
\begin{align*}
    \norm{\vv}_2 - \delta \le \sqrt{\delta^2+\norm{\vv}_2^2} - \delta \le \norm{\vv}_2.
\end{align*}
Indeed, we have
\begin{align*}
    \sqrt{\delta^2+\norm{\vv}_2^2} - \delta \le \sqrt{\delta^2} + \sqrt{\norm{\vv}_2^2} - \delta = \norm{\vv}_2,
\end{align*}
and
\begin{align*}
    \sqrt{\delta^2+\norm{\vv}_2^2} - \delta \ge \sqrt{\norm{\vv}_2^2} - \delta = \norm{\vv}_2 - \delta.
\end{align*}
From this, we get
\begin{align*}
    \ftilde_{\beta,\delta}(\vx) \le \max_{1 \le i \le m} \inparen{\sqrt{\delta^2+\norm{\mA_{S_i}\vx-\vb_{S_i}}_2^2}-\delta} + \beta\log m \le \gnorm{\mA\vx-\vb}{\infty} + \beta\log m
\end{align*}
and
\begin{align*}
    \ftilde_{\beta,\delta}(\vx) \ge \beta\logv{\sum_{i=1}^m \expv{\frac{\norm{\mA_{S_i}\vx-\vb_{S_i}}_2-\delta}{\beta}}} \ge \gnorm{\mA\vx-\vb}{\infty}-\delta.
\end{align*}
Putting these together gives
\begin{align*}
    \abs{\ftilde_{\beta,\delta}(\vx)-\gnorm{\mA\vx-\vb}{\infty}} \le \max\inparen{\beta\log m, \delta} \le \beta\log m + \delta,
\end{align*}
completing the proof of \Cref{lemma:smooth_approx}.
\end{proof}
Eventually, we will choose $\beta = \eps/(4\log m)$ and $\delta = \eps/4$ and then minimize $\ftilde_{\beta,\delta}$ to $\eps/2$ additive error. In light of \Cref{lemma:smooth_approx}, this will be enough to get an $\eps$-additive approximation to the optimum for $\gnorm{\mA\vx-\vb}{\infty}$.


\subsection{Calculus for \lsetext}
\label{sec:lse_calculus}

We investigate certain properties of $\lse_{\beta}(\vy)$ when each entry $[\vy]_i$ is a function $h_i(t)$ for $t \in \R$ for all $i\in[m]$. Let $h(t) \in \R^m$ denote the vector where its $i$th entry is given by $h_i(t)$. We treat each $h_i$ as a one-dimensional restriction of a function $\myfunc{g_i}{\R^m}{\R}$, so $h_i(t) = g_i(\vy+t\vd)$ for center $\vy$ and direction $\vd$ (we omit the parameters $\vy,\vd$ in the notation $h_i$ as it will be clear from context). Finally, recall the definition of quasi-self-concordance (\Cref{defn:qsc}).

We begin with calculating the first two derivatives of $\lse_{\beta}(h(t))$ with respect to $t$ in \Cref{lemma:lse_first_two}.
\begin{lemma}
\label{lemma:lse_first_two}
Let $\lambda_i(t) \coloneqq \expv{h_i(t)/\beta}$. Then, we have
\begin{align*}
    \inparen{\frac{d}{dt}}\lse_{\beta}(h(t)) &= \frac{\sum_{i=1}^m \inparen{\lambda_i(t) \cdot h_i'(t)}}{\sum_{i=1}^m \lambda_i(t)} \\
    \inparen{\frac{d}{dt}}^2\lse_{\beta}(h(t)) &= \frac{1}{\beta}\inparen{\frac{\sum_{i=1}^m \lambda_i(t)h_i'(t)^2}{\sum_{i=1}^m \lambda_i(t)} - \inparen{\frac{\sum_{i=1}^m \lambda_i(t)h_i'(t)}{\sum_{i=1}^m \lambda_i(t)}}^2} + \frac{\sum_{i=1}^m \lambda_i(t)h_i''(t)}{\sum_{i=1}^m \lambda_i(t)}.
\end{align*}
\end{lemma}
\begin{proof}[Proof of \Cref{lemma:lse_first_two}]
The first derivative follows from the chain rule. Indeed, we have
\begin{align*}
    \lse_{\beta}'(h(t)) &= \beta \cdot \frac{\sum_{i=1}^m \lambda_i'(t)}{\sum_{i=1}^m \lambda_i(t)} = \beta \cdot \frac{\sum_{i=1}^m \inparen{\lambda_i(t) \cdot \frac{h_i'(t)}{\beta}}}{\sum_{i=1}^m \lambda_i(t)} = \frac{\sum_{i=1}^m \inparen{\lambda_i(t) \cdot h_i'(t)}}{\sum_{i=1}^m \lambda_i(t)} \le \max_{i} h_i'(t).
\end{align*}
For the second derivative, we use the differentiation rule for multiplication and division and the chain rule, giving
\begin{align*}
    \lse_{\beta}''(h(t)) &= \frac{\left[\inparen{\sum_{i=1}^m \lambda_i'(t)h_i'(t)+\lambda_i(t)h_i''(t)}\inparen{\sum_{i=1}^m \lambda_i(t)}\right] - \frac{1}{\beta}\inparen{\sum_{i=1}^m \lambda_i(t)h_i'(t)}^2}{\inparen{\sum_{i=1}^m \lambda_i(t)}^2} \\
    &= \frac{\left[\frac{1}{\beta}\inparen{\sum_{i=1}^m \lambda_i(t)h_i'(t)^2+\beta\lambda_i(t)h_i''(t)}\inparen{\sum_{i=1}^m \lambda_i(t)}\right] - \frac{1}{\beta}\inparen{\sum_{i=1}^m \lambda_i(t)h_i'(t)}^2}{\inparen{\sum_{i=1}^m \lambda_i(t)}^2} \\
    &= \frac{1}{\beta}\inparen{\frac{\sum_{i=1}^m \lambda_i(t)h_i'(t)^2}{\sum_{i=1}^m \lambda_i(t)} - \frac{\inparen{\sum_{i=1}^m \lambda_i(t)h_i'(t)}^2}{\inparen{\sum_{i=1}^m \lambda_i(t)}^2}} + \frac{\sum_{i=1}^m \lambda_i(t)h_i''(t)}{\sum_{i=1}^m \lambda_i(t)}\enspace.
\end{align*}
This completes the proof of \Cref{lemma:lse_first_two}.
\end{proof}
Next, we prove a general fact regarding composing $\lse$ with a vector formed by functions that are themselves quasi self concordant. See \Cref{lemma:composed_qsc}.
\begin{restatable}[Composing softmax with quasi-self-concordant functions]{lemma}{composedqsc}
\label{lemma:composed_qsc}
Let $\norm{\cdot}$ be an arbitrary norm and $h_1,\dots,h_m$ be such that for all $1 \le i \le m$ and for all $\vy,\vd\in\R^m$ and $t\in\R$,
\begin{align*}
    &\inparen{\frac{d}{dt}} h_i(t) \le \norm{\vd} & \text{(Lipschitzness)}\\
    &\abs{\inparen{\frac{d}{dt}}^3 h_i(t)} \le \nu\norm{\vd}\inparen{\frac{d}{dt}}^2 h_i(t) & \text{(quasi-self-concordance)}.
\end{align*}
Then, for all $\vy,\vd\in\R^m$ and all $t\in\R$, we have
\begin{align*}
    \abs{\inparen{\frac{d}{dt}}^3 \lse_{\beta}(h(t))} \le \inparen{\frac{16}{\beta}+\nu}\norm{\vd}\inparen{\frac{d}{dt}}^2 \lse_{\beta}(h(t)).
\end{align*}
\end{restatable}
As far as we are aware, this type of composition result was not previously known and may be of independent interest.

To prove \Cref{lemma:composed_qsc}, we need \Cref{lemma:variance_rv_bound}.
\begin{lemma}
\label{lemma:variance_rv_bound}
For any two random variables $X, Y$, we have
\begin{align*}
    \var{XY} \le 2\maxnorm{Y}^2\var{X}+2\maxnorm{X}^2\var{Y}.
\end{align*}
\end{lemma}
\begin{proof}[Proof of \Cref{lemma:variance_rv_bound}]
The proof follows that of \cite{var_stackexchange}, but we reproduce it here for completeness. First, notice that for random variables $U, V$, we have
\begin{align*}
    2\var{U}+2\var{V} - \var{U+V} = \var{U}+\var{V}-2\cov{U,V} = \var{U-V} \ge 0.
\end{align*}
Let $U=(X-\exv{X})Y$ and $V=\exv{X}Y$. Then, $U+V=XY$, and we have
\begin{align*}
    \var{XY} \le 2\var{(X-\exv{X})Y}+2\var{\exv{X}Y}=2\var{(X-\exv{X})Y}+2\exv{X}^2\var{Y}.
\end{align*}
It remains to bound $\var{(X-\exv{X})Y}$. By H\"older's inequality, we have
\begin{align*}
    \var{(X-\exv{X})Y} \le \exv{((X-\exv{X})Y)^2} \le \exv{(X-\exv{X})^2}\maxnorm{Y}^2=\var{X}\maxnorm{Y}^2.
\end{align*}
Combining everything gives us the conclusion of \Cref{lemma:variance_rv_bound}.
\end{proof}
We are now ready to prove \Cref{lemma:composed_qsc}.
\begin{proof}[Proof of \Cref{lemma:composed_qsc}]
Let $\lambda_i(t) \coloneqq \expv{h_i(t)/\beta}$.

In this proof, we will encounter many weighted averages of vectors $\vz \in \R^m$ of the form
\begin{align*}
    \frac{\sum_{i=1}^m \lambda_i(t)z_i}{\sum_{i=1}^m \lambda_i(t)}.
\end{align*}
Let $\cD$ be the distribution over $[m]$ whose entries are given by $\cD_j = \lambda_j(t)/\sum_{i=1}^m \lambda_i(t)$. In the rest of this proof, all expected values, variances, and covariances will be taken with respect to this distribution. In an abuse of notation, let $h(t)$ denote the ``random'' variable that is $h_i(t)$ with probability $\cD_i$. Define $h'(t), h''(t), h'''(t)$ analogously.

To find the third derivative of $\lse_{\beta}(h(t))$, we start with its second derivative. By \Cref{lemma:lse_first_two}, it is given by
\begin{align*}
    \lse_{\beta}''(h(t)) &= \underbrace{\frac{1}{\beta}\inparen{\frac{\sum_{i=1}^m \lambda_i(t)h_i'(t)^2}{\sum_{i=1}^m \lambda_i(t)} - \inparen{\frac{\sum_{i=1}^m \lambda_i(t)h_i'(t)}{\sum_{i=1}^m \lambda_i(t)}}^2}}_{T_1} + \underbrace{\frac{\sum_{i=1}^m \lambda_i(t)h_i''(t)}{\sum_{i=1}^m \lambda_i(t)}}_{T_2} \\
    &= \frac{1}{\beta}\var{h'(t)}+\exv{h''(t)}.
\end{align*}
We now differentiate the above term by term. First, we have
\begin{align*}
    T_2'(t) &= \frac{\sum_{i=1}^m \lambda_i(t)\inparen{\inparen{\frac{h_i'(t)h_i''(t)}{\beta}}+h_i'''(t)}}{\sum_{i=1}^m \lambda_i(t)} - \frac{1}{\beta}\cdot\frac{\inparen{\sum_{i=1}^m \lambda_i(t)h_i'(t)}\inparen{\sum_{i=1}^m \lambda_i(t)h_i''(t)}}{\inparen{\sum_{i=1}^m \lambda_i(t)}^2} \\
    &= \frac{1}{\beta}\inparen{\frac{\sum_{i=1}^m \lambda_i(t)h_i'(t)h_i''(t)}{\sum_{i=1}^m \lambda_i(t)} - \frac{\inparen{\sum_{i=1}^m \lambda_i(t)h_i'(t)}\inparen{\sum_{i=1}^m \lambda_i(t)h_i''(t)}}{\inparen{\sum_{i=1}^m \lambda_i(t)}^2}} + \frac{\sum_{i=1}^m \lambda_i(t)h_i'''(t)}{\sum_{i=1}^m \lambda_i(t)} \\
    &= \frac{1}{\beta}\cov{h'(t),h''(t)} + \exv{h'''(t)}.
\end{align*}
Next, we have
\begin{align*}
    \frac{d}{dt}\exv{h'(t)}^2 = 2\exv{h'(t)} \cdot \frac{d}{dt}\exv{h'(t)} = 2\exv{h'(t)}\inparen{\frac{1}{\beta}\var{h'(t)} + \exv{h''(t)}}
\end{align*}
and
\begin{align*}
    &\quad \frac{d}{dt}\exv{h'(t)^2}\\
    &= \frac{\inparen{\sum_{i=1}^m \lambda_i'(t)h_i'(t)^2 + 2h_i'(t)h_i''(t)\lambda_i(t)}\inparen{\sum_{i=1}^m \lambda_i(t)} - \frac{1}{\beta}\inparen{\sum_{i=1}^m \lambda_i(t)h_i'(t)}\inparen{\sum_{i=1}^m \lambda_i(t)h_i'(t)^2}}{\inparen{\sum_{i=1}^m \lambda_i(t)}^2} \\
    &= \frac{\inparen{\sum_{i=1}^m \lambda_i'(t)h_i'(t)^2 + 2h_i'(t)h_i''(t)\lambda_i(t)}}{\sum_{i=1}^m \lambda_i(t)} - \frac{1}{\beta}\cdot\frac{\inparen{\sum_{i=1}^m \lambda_i(t)h_i'(t)}\inparen{\sum_{i=1}^m \lambda_i(t)h_i'(t)^2}}{\inparen{\sum_{i=1}^m \lambda_i(t)}^2} \\
    &= \frac{\sum_{i=1}^m \lambda_i(t)\inparen{\frac{h_i'(t)^3}{\beta} + 2h_i'(t)h_i''(t)}}{\sum_{i=1}^m \lambda_i(t)} - \frac{1}{\beta}\cdot\frac{\inparen{\sum_{i=1}^m \lambda_i(t)h_i'(t)}\inparen{\sum_{i=1}^m \lambda_i(t)h_i'(t)^2}}{\inparen{\sum_{i=1}^m \lambda_i(t)}^2} \\
    &= \frac{1}{\beta}\cov{h'(t), h'(t)^2} + 2\exv{h'(t)h''(t)}.
\end{align*}
Combining everything gives us
\begin{align*}
    &\quad \lse_{\beta}'''(h(t)) \\
    &= \frac{1}{\beta}\inparen{\frac{1}{\beta}\cov{h'(t), h'(t)^2} + 2\exv{h'(t)h''(t)} - 2\exv{h'(t)}\inparen{\frac{1}{\beta}\var{h'(t)} + \exv{h''(t)}}}\\
    &\quad + \frac{1}{\beta}\cov{h'(t),h''(t)} + \exv{h'''(t)} \\
    &= \frac{1}{\beta^2}\cov{h'(t),h'(t)^2} - \frac{2}{\beta^2}\exv{h'(t)}\var{h'(t)} + \frac{3}{\beta}\cov{h'(t),h''(t)} + \exv{h'''(t)}.
\end{align*}
We first analyze the terms that only depend on $h'(t)$. To do so, we use \Cref{lemma:variance_rv_bound} to write
\begin{align*}
    \abs{\cov{h'(t),h'(t)^2}} \le \sqrt{\var{h'(t)}}\sqrt{\var{h'(t)^2}} \le 2\norm{\vd}\var{h'(t)}.
\end{align*}
Now, we have
\begin{align*}
    &\quad \frac{1}{\beta^2}\abs{\cov{h'(t),h'(t)^2} - 2\exv{h'(t)}\var{h'(t)}} \\
    &\le \frac{1}{\beta^2}\abs{\cov{h'(t),h'(t)^2}} + \frac{2}{\beta^2}\abs{\exv{h'(t)}\var{h'(t)}} \\
    &\le \frac{4}{\beta^2}\norm{\vd}\var{h'(t)} \le \frac{4}{\beta}\norm{\vd}\lse_{\beta}''(h(t)).
\end{align*}
Next, we take care of the remaining terms. We have
\begin{align*}
    \frac{3}{\beta}\abs{\cov{h'(t),h''(t)}} + \abs{\exv{h'''(t)}} &\le \frac{6}{\beta}\inparen{\max_{i} h_i'(t)}\exv{\abs{h''(t)-\exv{h''(t)}}} + \abs{\exv{h'''(t)}} \\
    &\le \frac{12}{\beta}\norm{\vd}\lse_{\beta}''(h(t)) + \exv{\abs{h'''(t)}} \\
    &\le \frac{12}{\beta}\norm{\vd}\lse_{\beta}''(h(t)) + \nu\norm{\vd}\exv{h''(t)} \\
    &\le \inparen{\frac{12}{\beta}+\nu}\norm{\vd}\lse_{\beta}''(h(t)),
\end{align*}
where the penultimate line follows from \Cref{lemma:hi_qsc}. Combining these conclusions yields
\begin{align*}
    \abs{\lse_{\beta}'''(h(t))} \le \inparen{\frac{16}{\beta}+\nu}\norm{\vd}\lse_{\beta}''(h(t)),
\end{align*}
completing the proof of \Cref{lemma:composed_qsc}.
\end{proof}

\subsection{Smoothness and quasi-self-concordance of the modified objective}
\label{sec:lse_calculus_special}

The main result of this subsection is \Cref{lemma:obj_smooth_qsc}.

\begin{lemma}
\label{lemma:obj_smooth_qsc}
Let $\mW$ be such that for all $\vz \in \R^d$, we have $\gnorm{\mA\vz}{\infty} \le \norm{\mW^{1/2}\mA\vz}_2$. For all $\vx, \vz \in \R^d$ and $t\in\R$, we have
\begin{align*}
    \inparen{\frac{d}{dt}}^2 \ftilde_{\beta,\delta}(\vx+t\vz) &\le \inparen{\frac{1}{\delta}+\frac{1}{\beta}}\norm{\mW^{1/2}\mA\vz}_2^2 & \text{(smoothness)} \\
    \abs{\inparen{\frac{d}{dt}}^3 \ftilde_{\beta,\delta}(\vx+t\vz)} &\le \inparen{\frac{16}{\delta}+\frac{3}{\beta}}\norm{\mW^{1/2}\mA\vz}_2\inparen{\frac{d}{dt}}^2 \ftilde_{\beta,\delta}(\vx+t\vz) & \text{(quasi-self-concordance)}.
\end{align*}
\end{lemma}

Our goal in the rest of this section is to prove \Cref{lemma:obj_smooth_qsc}.

We begin with defining $h_i(t)$ as (absorb the $\delta,\vy,\vd$ parameters into the definition of $h_i$)
\begin{align*}
    h_i(t) \coloneqq \sqrt{\delta^2+\norm{\vy_{S_i}+t\vd_{S_i}}_2^2}.
\end{align*}
Let $h(t)$ denote the vector whose $i$th entry is $h_i(t)$. Then, observe that
\begin{align*}
    \lse_{\beta}(h(t)) = \beta\logv{\sum_{i=1}^m \expv{\frac{h_i(t)}{\beta}}} = \beta\logv{\sum_{i=1}^m \expv{\frac{\sqrt{\delta^2+\norm{\vy_{S_i}+t\vd_{S_i}}_2^2}}{\beta}}}.
\end{align*}
It is easy to see that every one-dimensional restriction of $\ftilde_{\beta,\delta}$ can be obtained by an affine transformation of $\lse_{\beta}(h(t))$ after appropriate choices of $\vy,\vd\in\R^m$. Hence, we first analyze $\lse_{\beta}(h(t))$ for all $\vy,\vd\in\R^m$. 

We begin with proving the smoothness of $\lse_{\beta}(h(t))$ with respect to $\gnorm{\cdot}{\infty}$.
\begin{lemma}
\label{lemma:base_smooth}
For all $\vy,\vd\in\R^m$ and all $t \in \R$, we have
\begin{align*}
    \inparen{\frac{d}{dt}}^2 \lse_{\beta}(h(t)) \le \inparen{\frac{1}{\delta}+\frac{1}{\beta}}\gnorm{\vd}{\infty}^2.
\end{align*}
\end{lemma}
\begin{proof}[Proof of \Cref{lemma:base_smooth}]
By direct calculation, it is easy to see that
\begin{equation}\label{eq:hi_derivatives}
    \begin{aligned}
         h_i'(t) &= \frac{\ip{\vy_{S_i}+t\vd_{S_i},\vd_{S_i}}}{h_i(t)} \\
         h_i''(t) &= \frac{\norm{\vd_{S_i}}_2^2 h_i(t) - h_i'(t)^2h_i(t)}{h_i(t)^2} = \frac{\norm{\vd_{S_i}}_2^2-h_i'(t)^2}{h_i(t)}.
    \end{aligned}
\end{equation}
We plug this into the result of \Cref{lemma:lse_first_two} and get
\begin{align*}
    \lse_{\beta}''(h(t)) &\le \frac{1}{\beta}\max_i h_i'(t)^2 + \max_i h_i''(t) \\
    &= \frac{1}{\beta} \max_{i}\inparen{\frac{\ip{\vy_{S_i}+t\vd_{S_i},\vd_{S_i}}}{\sqrt{\delta^2 + \norm{\vy_{S_i}+t\vd_{S_i}}_2^2}}}^2 + \max_{i} \frac{\norm{\vd_{S_i}}_2^2 - h_i'(t)^2}{\sqrt{\delta^2+\norm{\vy_{S_i}+t\vd_{S_i}}_2^2}} \\
    &\le \frac{1}{\beta}\max_{i} \norm{\vd_{S_i}}_2^2 + \frac{1}{\delta}\max_{i}\norm{\vd_{S_i}}_2^2 = \inparen{\frac{1}{\beta}+\frac{1}{\delta}}\gnorm{\vd}{\infty}^2,
\end{align*}
completing the proof of \Cref{lemma:base_smooth}.
\end{proof}
Our next task is to show that $\lse_{\beta}(h(t))$ is $O(1/\beta + 1/\delta)$-quasi-self-concordant in $\gnorm{\cdot}{\infty}$. To do so, we will appeal to \Cref{lemma:composed_qsc}. To be able to do this, we first have to prove the quasi-self-concordance of each component function in $\lse_{\beta}(h(t))$.
\begin{lemma}
\label{lemma:hi_qsc}
For all $\vy,\vd\in\R^m$ and all $t\in\R$, we have
\begin{align*}
    \abs{\inparen{\frac{d}{dt}}^3 \sqrt{\delta^2+\norm{\vy_{S_i}+t\vd_{S_i}}_2^2}} \le \frac{3}{\delta}\norm{\vd_{S_i}}_2\inparen{\inparen{\frac{d}{dt}}^2\sqrt{\delta^2+\norm{\vy_{S_i}+t\vd_{S_i}}_2^2}}.
\end{align*}
\end{lemma}
\begin{proof}[Proof of \Cref{lemma:hi_qsc}]
Although a similar fact appears in \cite[Section 2.1.2]{ob18}, it is not in the exact form we need. So, we prove the required statement here.

Recycling the computation from \eqref{eq:hi_derivatives}, recall
\begin{align*}
    h_i''(t) = \frac{\norm{\vd_{S_i}}_2^2 - h_i'(t)^2}{h_i(t)},
\end{align*}
which gives
\begin{align*}
    h_i'''(t) = \frac{-2h_i'(t)h_i''(t)h_i(t) - h_i'(t)(h_i(t)h_i''(t))}{h_i(t)^2} = -\frac{3h_i'(t)h_i''(t)}{h_i(t)}.
\end{align*}
Finally, again recalling \eqref{eq:hi_derivatives}, notice that
\begin{align*}
    \abs{\frac{h_i'(t)}{h_i(t)}} = \abs{\frac{\ip{\vy_{S_i}+t\vd_{S_i},\vd_{S_i}}}{h_i(t)^2}} = \abs{\ip{\frac{\vy_{S_i}+t\vd_{S_i}}{\sqrt{\delta^2+\norm{\vy_{S_i}+t\vd_{S_i}}_2^2}}, \frac{\vd_{S_i}}{\sqrt{\delta^2+\norm{\vy_{S_i}+t\vd_{S_i}}_2^2}}}} \le \frac{\norm{\vd_{S_i}}_2}{\delta}.
\end{align*}
Combining everything completes the proof of \Cref{lemma:hi_qsc}.
\end{proof}
We are now ready to prove the quasi-self-concordance of $\lse_{\beta}(h(t))$ in $\gnorm{\cdot}{\infty}$.
\begin{lemma}
\label{lemma:base_qsc}
For all $\vy, \vd \in \R^m$ and $t\in\R$, we have
\begin{align*}
    \abs{\inparen{\frac{d}{dt}}^3 \lse_{\beta}(h(t))} \le \inparen{\frac{16}{\beta}+\frac{3}{\delta}}\gnorm{\vd}{\infty}\inparen{\frac{d}{dt}}^2 \lse_{\beta}(h(t)).
\end{align*}
\end{lemma}
\begin{proof}[Proof of \Cref{lemma:base_qsc}]
In the statement of \Cref{lemma:composed_qsc}, let $\norm{\cdot} = \gnorm{\cdot}{\infty}$. By the definition of $\gnorm{\cdot}{\infty}$ and $h_i$, we have for all $i$ and $t$ that $h_i'(t) \le \gnorm{\vd}{\infty}$. Additionally, from \Cref{lemma:hi_qsc}, we have that the $h_i(t)$ are $3/\delta$-quasi-self-concordant in the norm $\gnorm{\vd}{\infty}$ for all $i$. \Cref{lemma:base_qsc} now follows immediately from \Cref{lemma:composed_qsc}.
\end{proof}
Finally, we can prove \Cref{lemma:obj_smooth_qsc}.
\begin{proof}[Proof of \Cref{lemma:obj_smooth_qsc}]
By the conclusion of \Cref{lemma:base_smooth}, we know that for all $\vy, \vd \in \R^m$ and $t\in\R$ that
\begin{align*}
    \inparen{\frac{d}{dt}}^2 \lse_{\beta}(h(t)) \le \inparen{\frac{1}{\delta}+\frac{1}{\beta}}\gnorm{\vz}{\infty}^2.
\end{align*}
Let $\vy = \mA\vx - \vb$ for some $\vx$ and $\vd = \mA\vz$ for some $\vz$. Let
\begin{align*}
    g(\vy) \coloneqq \beta\logv{\sum_{i=1}^m \expv{\frac{\sqrt{\delta^2+\norm{\vy_{S_i}}_2^2}-\delta}{\beta}}}.
\end{align*}
Then,
\begin{align*}
    \inparen{\frac{d}{dt}}^2 \ftilde_{\beta,\delta}(\vx+t\vz) = \inparen{\frac{d}{dt}}^2 g(\mA\vx - \vb + t\mA\vz) \le \inparen{\frac{1}{\delta}+\frac{1}{\beta}}\gnorm{\mA\vz}{\infty}^2.
\end{align*}
With the exact same reasoning applied to the conclusion of \Cref{lemma:base_qsc}, we also see that
\begin{align*}
    \abs{\inparen{\frac{d}{dt}}^3 \ftilde_{\beta,\delta}(\vx+t\vz)} \le \inparen{\frac{16}{\delta}+\frac{3}{\beta}}\gnorm{\mA\vz}{\infty}\inparen{\frac{d}{dt}}^2 \ftilde_{\beta,\delta}(\vx+t\vz).
\end{align*}
The conclusion of \Cref{lemma:obj_smooth_qsc} then follows from remembering that we have $\mW$ such that for all $\vz\in\R^d$, $\gnorm{\mA\vz}{\infty} \le \norm{\mW^{1/2}\mA\vz}_2$ (following from \Cref{thm:gp_regression_blw_intro}).
\end{proof}

\subsection{Analysis of \texorpdfstring{\Cref{alg:fair_regression}}{Algorithm 1}}
\label{sec:gp_robust_analysis_combining}

In this subsection, we use the calculus facts from the previous two subsections to analyze \Cref{alg:fair_regression}. The outline of this proof follows that of \cite[Theorem 2]{jls21}, which in turn builds up to using the proof used in \cite[Corollary 12]{msbacon}. The main idea is to define the algorithm based on the norm given by a good choice of positive semidefinite $\mM$, given by \Cref{thm:gp_regression_blw_intro}.

In the rest of this section, let $\mW$ be factor-$2$ block Lewis weight overestimates for $\insquare{\mA \vert \vb}$. As in Line \ref{line:blw} of \Cref{alg:fair_regression} and from the corresponding guarantee given in \cite[Lemmas 5.6, 5.8]{mo23}, this means that within $2\log m$ linear system solves in $\mA^{\top}\mD\mA$ for diagonal $\mD$, we can find $\mW$ such that for all $\vx\in\R^d$ and $c\in\R$ we have
\begin{align*}
    \gnorm{\mA\vx-c\vb}{\infty} \le \norm{\mW^{1/2}\mA\vx-c\mW^{1/2}\vb}_2 \le \sqrt{2(\rank{\mA}+1)}\gnorm{\mA\vx-c\vb}{\infty}.
\end{align*}
Note that choosing $c=1$ yields our original objective on either side of the above inequality. Motivated by the above, it is natural to use the norm given by $\mM \coloneqq \mA^{\top}\mW\mA$ to give the geometry for the ball optimization oracle and for the analysis. Additionally, without loss of generality and for the sake of the analysis, let us rescale the problem so that
\begin{align*}
    1 = \opt \coloneqq \gnorm{\mA\xstar-\vb}{\infty}.
\end{align*}
Also, as mentioned earlier, assume without loss of generality that $\rank{\mA} = d$.

We begin with \Cref{lemma:initialization}, which bounds our initial suboptimality in $\ftilde$ and in $\norm{\cdot}_{\mM}$.
\begin{lemma}
\label{lemma:initialization}
Let $\xtilde_{\beta,\delta} \coloneqq \argmin{\vx\in\R^d} \ftilde_{\beta,\delta}(\vx)$. Then,
\begin{equation*}
    \begin{aligned}
        \norm{\xtilde_{\beta,\delta} - \vx_0}_{\mM} &\le (2+2(\beta\log m + \delta))\sqrt{2(d+1)} \\
        \ftilde_{\beta,\delta}(\vx_0) - \ftilde_{\beta,\delta}(\xtilde_{\beta,\delta}) &\le \sqrt{2(d+1)} - 1+2(\beta\log m + \delta)
    \end{aligned}.
\end{equation*}
\end{lemma}
\begin{proof}[Proof of \Cref{lemma:initialization}]
It is easy to check that
\begin{align*}
    \vx_0 \coloneqq \inparen{\mA^{\top}\mW\mA}^{-1}\mA^{\top}\mW\vb = \argmin{\vx \in \R^d} \norm{\mW^{1/2}\mA\vx-\mW^{1/2}\vb}_2.
\end{align*}
By \Cref{lemma:smooth_approx}, for all $\vx\in\R^d$,
\begin{align*}
    \abs{\ftilde_{\beta,\delta}(\vx)-\gnorm{\mA\vx-\vb}{\infty}} \le \beta\log m+\delta,
\end{align*}
implying
\begin{align*}
    \abs{\gnorm{\mA\xstar-\vb}{\infty} - \ftilde_{\beta,\delta}(\xtilde_{\beta,\delta})} \le \beta\log m+\delta.
\end{align*}
Combining this with \Cref{thm:gp_regression_blw_intro}, we get
\begin{align*}
    1 \le \gnorm{\mA\xstar-\vb}{\infty} \le \gnorm{\mA\vx_0-\vb}{\infty} \le \norm{\mW^{1/2}\mA\vx_0 - \mW^{1/2}\vb}_2
\end{align*}
and
\begin{align*}
    \frac{\norm{\mW^{1/2}\mA\vx_0 - \mW^{1/2}\vb}_2}{\sqrt{2(d+1)}} \le \frac{\norm{\mW^{1/2}\mA\xstar-\mW^{1/2}\vb}_2}{\sqrt{2(d+1)}} \le \gnorm{\mA\xstar-\vb}{\infty} = 1.
\end{align*}
Combining these gives
\begin{align*}
    1 \le \norm{\mW^{1/2}\mA\vx_0 - \mW^{1/2}\vb}_2 \le \sqrt{2(d+1)}.
\end{align*}
Additionally,
\begin{align*}
     \norm{\mW^{1/2}\mA\xtilde_{\beta,\delta}-\mW^{1/2}\vb}_2 &\le \sqrt{2(d+1)}\gnorm{\mA\xtilde_{\beta,\delta}-\vb}{\infty} \\
     &\le \sqrt{2(d+1)}\inparen{\ftilde_{\beta,\delta}(\xtilde_{\beta,\delta})+\beta\log m+\delta} \\
     &\le \sqrt{2(d+1)}\inparen{\gnorm{\mA\xstar-\vb}{\infty}+2(\beta\log m + \delta)} \\
     &= \sqrt{2(d+1)}(1+2(\beta\log m + \delta)).
\end{align*}
Then,
\begin{align*}
    \norm{\xtilde-\vx_0}_{\mM} &= \norm{\inparen{\mW^{1/2}\mA\xtilde_{\beta,\delta}-\mW^{1/2}\vb} - \inparen{\mW^{1/2}\mA\vx_0-\mW^{1/2}\vb}}_2 \\
    &\le\norm{\mW^{1/2}\mA\xtilde_{\beta,\delta}-\mW^{1/2}\vb}_2 +\norm{\mW^{1/2}\mA\vx_0-\mW^{1/2}\vb}_2 \\
    &\le (2+2(\beta\log m + \delta))\sqrt{2(d+1)},
\end{align*}
and
\begin{align*}
    \ftilde_{\beta,\delta}(\vx_0) - \ftilde_{\beta,\delta}(\xtilde_{\beta,\delta}) &\le \gnorm{\mA\vx_0-\vb}{\infty} - \gnorm{\mA\xstar-\vb}{\infty}+2(\beta\log m + \delta) \\
    &\le \norm{\mW^{1/2}\mA\vx_0 - \mW^{1/2}\vb}_2 - \opt+2(\beta\log m + \delta) \\
    &\le \sqrt{2(d+1)}-1+2(\beta\log m + \delta).
\end{align*}
This completes the proof of \Cref{lemma:initialization}.
\end{proof}

We are now ready to prove \Cref{mainthm:fair_regression_iteration_complexity}.

\begin{proof}[Proof of \Cref{mainthm:fair_regression_iteration_complexity}]
\Cref{alg:fair_regression} optimizes the regularization of $\ftilde$ given by
\begin{align*}
    \fhat(\vx) \coloneqq \ftilde_{\beta,\delta}(\vx) + \frac{\eps}{110R^2}\norm{\mW^{1/2}\mA(\vx-\vx_0)}_2^2,
\end{align*}
where $R$ is such that $\norm{\vx_0 - \xtilde_{\beta,\delta}}_{\mM} \le R$. Let $ \xhat \coloneqq \argmin{\vx\in\R^d} \fhat(\vx)$. Using \cite[Proof of Corollary 12]{msbacon}, we know that for every iterate $\vx$ of \Cref{alg:fair_regression},
\begin{align*}
    \abs{\fhat(\vx)-\ftilde_{\beta,\delta}(\vx)} \le \frac{\eps}{4}.
\end{align*}
We now choose $\beta = \eps/(4\log m)$ and $\delta = \eps/4$, so that $\ftilde_{\beta,\delta}$ approximates $f$ up to error $\eps/2$ on every point. Using \Cref{lemma:initialization}, this gives $R = (2+\eps)\sqrt{2(d+1)}$. It is therefore sufficient to optimize $\fhat$ up to $\eps/4$ additive error.

Next, using \Cref{lemma:obj_smooth_qsc} and \cite[Lemmas 11, 43]{msbacon}, we have that $\fhat$ is $(1/\nu, e)$-Hessian stable in $\norm{\cdot}_{\mM}$ for $\nu = \Omega(1/(\eps\log m))$. We now invoke \cite[Theorem 9]{msbacon}, which tells us that we can implement a $(C/\sqrt{d}, C/\eps)$-ball optimization oracle for $f$ with $O\inparen{\logv{\frac{d}{\eps}}^2}$ linear system solves. 

The next step is to turn the ball optimization oracle into a $\frac{1}{2}$-MS oracle (\Cref{defn:ms_oracle}). Using \cite[Proposition 5]{msbacon}, we get a ball oracle complexity of $O\inparen{\logv{\frac{d}{\eps}}}$ to implement the MS oracle. In total, our linear system solve complexity for implementing the MS oracle for iteration $t$ is $O\inparen{\logv{\frac{d}{\eps}}^3}$.

Finally, using \cite[Theorem 6]{msbacon}, we get that \Cref{alg:fair_regression} has a Newton iteration complexity of
\begin{align*}
    &\quad O\inparen{\inparen{\frac{(1+\eps)\sqrt{d}\log m}{\eps}}^{2/3}\logv{\frac{\sqrt{d}+\eps}{\eps}}\inparen{\logv{\frac{(\log m/\eps)d(1+(1+\eps)\sqrt{d}\log m/\eps)}{\eps}}}^3} \\
    &= O\inparen{\frac{d^{1/3}}{\eps^{2/3}}\logv{\frac{d\log m}{\eps}}^{14/3}},
\end{align*}
as promised.

Next, we analyze what happens if we fall in the case where $\mW=\mI_m$. Here, by using the $\sqrt{m}$ distortion from approximating $\ell_{\infty}^m$ with $\ell_{2}^m$, we have for all $\vx\in\R^d$,
\begin{align*}
    \frac{\norm{\mA\vx-\vb}_2}{\sqrt{m}} \le \gnorm{\mA\vx-\vb}{\infty} \le \norm{\mA\vx-\vb}_2.
\end{align*}
Using this and repeating the previous analysis with this choice of $\mM$ gives us a rate of
\begin{align*}
    O\inparen{\frac{m^{1/3}}{\eps^{2/3}}\logv{\frac{m\log m}{\eps}}^{14/3}},
\end{align*}
as required.

It remains to determine the form of the Newton steps. For this, it is sufficient to understand the Hessian of $\fhat$. A straightforward calculation shows that it is of the form $\mA^{\top}\mB\mA$ where $\mB$ is a block-diagonal matrix where each block has size $\abs{S_i} \times \abs{S_i}$. Thus, each Newton step solves a linear system of the form $\mA^{\top}\mB\mA\vz = \vv$.

Combining this with the iteration complexity guarantee to find $\mW$ (see \Cref{thm:gp_regression_blw_intro}) completes the proof of \Cref{mainthm:fair_regression_iteration_complexity}.
\end{proof}

\section{Interpolating between average and robust losses}
\label{sec:gp_interpolating_proof}

In this section, we prove \Cref{mainthm:gp_regression_iteration_complexity}. As before, our proof follows the outline in \Cref{sec:gp_reg_overview}. The main technical challenges are to establish a form of strong convexity for our objective $f$ and then to build a solver for the proximal problem \eqref{eq:gp_prox_intro_reg}. 

The rest of this section is organized as follows. In \Cref{sec:gp_objective_calculus}, we derive calculus facts about our objective $f$, including bounds on its Hessian and the promised strong convexity (particularly \Cref{lemma:strong_convexity_gp} and the more general result it builds on, \Cref{lemma:strong_convexity_component}). In \Cref{sec:gp_iterate_facts}, we prove some facts about the iterates of \Cref{alg:optimal_ms_acceleration} when applied to our setting. In \Cref{sec:gp_prox_solver}, we more precisely define and analyze our solver for proximal sub-problems. This section is fairly technical and we give a more detailed outline there. Finally, in \Cref{sec:gp_interpolating_alg}, we assemble all these components and analyze \Cref{alg:gp_regression_alg}, thereby proving \Cref{mainthm:gp_regression_iteration_complexity}.

Throughout this analysis, we rescale the problem so that $f(\xstar) = 1$. It is now sufficient to solve for an $\eps$-additive error solution.

\subsection{Calculus for the objective}
\label{sec:gp_objective_calculus}

In this section, we work out some calculus facts related to our objective $\gnorm{\mA\vx-\vb}{p}^p$. Throughout this discussion, let $f(\vx) \coloneqq \gnorm{\mA\vx-\vb}{p}^p$.

\begin{lemma}
\label{lemma:gp_regression_hessian_upperbound}
For any $\vz\in\R^d$, we have
\begin{align*}
    p\sum_{i=1}^m \norm{\mA_{S_i}\vx-\vb_{S_i}}_2^{p-2}\norm{\mA_{S_i}\vz}_2^2 \le \vz^{\top}\inparen{\nabla^2 f(\vx)}\vz \le p(p-1)\sum_{i=1}^m\norm{\mA_{S_i}\vx - \vb_{S_i}}_2^{p-2}\norm{\mA_{S_i}\vz}_2^2.
\end{align*}
\end{lemma}
\begin{proof}[Proof of \Cref{lemma:gp_regression_hessian_upperbound}]
Let us first calculate the derivative and hessian for $f(\cdot)$ using the chain rule and usual matrix differentiation rules:
\begin{align}
    f(\vx) &= \sum_{i=1}^m\norm{\mA_{S_i}\vx - \vb_{S_i}}_2^p\enspace,\nonumber\\
    \nabla f(\vx) &= p\sum_{i=1}^m \norm{\mA_{S_i}\vx - \vb_{S_i}}_2^{p-2}\mA_{S_i}^{\top}(\mA_{S_i}\vx - \vb_{S_i})\enspace,\label{eq:f_grad}\\
    \nabla^2 f(\vx) &= p\sum_{i=1}^m\norm{\mA_{S_i}\vx - \vb_{S_i}}_2^{p-2}\mA_{S_i}^{\top}\mA_{S_i} \nonumber\\
    &\quad + p(p-2)\sum_{i=1}^m\norm{\mA_{S_i}\vx - \vb_{S_i}}_2^{p-4}\left(\mA_{S_i}^{\top}(\mA_{S_i}\vx - \vb_{S_i})(\mA_{S_i}\vx - \vb_{S_i})^{\top}\mA_{S_i}\right)\enspace.\label{eq:f_hess}
\end{align}
Using this formula, we take the quadratic form with respect to a vector $\vz$. By Cauchy-Schwarz, notice that
\begin{align*}
    &\quad \vz^{\top}\norm{\mA_{S_i}\vx - \vb_{S_i}}_2^{p-4}\left(\mA_{S_i}^{\top}(\mA_{S_i}\vx - \vb_{S_i})(\mA_{S_i}\vx - \vb_{S_i})^{\top}\mA_{S_i}\right)\vz \\
    &= \norm{\mA_{S_i}\vx - \vb_{S_i}}_2^{p-4}\ip{\mA_{S_i}\vz, \mA_{S_i}\vx-\vb_{S_i}}^2 \le \norm{\mA_{S_i}\vx-\vb_{S_i}}_2^{p-2}\norm{\mA_{S_i}\vz}_2^2.
\end{align*}
With that, we have
\begin{align}
    \vz^{\top}\inparen{\nabla^2 f(\vx)}\vz &\le p\sum_{i=1}^m \norm{\mA_{S_i}\vx-\vb_{S_i}}^{p-2}\norm{\mA_{S_i}\vz}_2^2 + (p-2)\norm{\mA_{S_i}\vx-\vb_{S_i}}^{p-2}\norm{\mA_{S_i}\vz}_2^2 \enspace,\nonumber\\
    &= p(p-1)\sum_{i=1}^m\norm{\mA_{S_i}\vx - \vb_{S_i}}_2^{p-2}\norm{\mA_{S_i}\vz}_2^2\enspace.\label{eq:f_quad_form}
\end{align}
For the lower bound, we use our calculation for $\nabla^2 f(\vx)$ to write
\begin{align*}
    \vz^{\top}\inparen{\nabla^2 f(\vx)}\vz \ge p\sum_{i=1}^m \norm{\mA_{S_i}\vx-\vb_{S_i}}_2^{p-2}\norm{\mA_{S_i}\vz}_2^2,
\end{align*}
completing the proof of \Cref{lemma:gp_regression_hessian_upperbound}.
\end{proof}

\subsubsection{Strong convexity of the objective}

The main pair of results of this section are \Cref{lemma:strong_convexity_gp} and \Cref{lemma:strong_convexity_component}. We can think of \Cref{lemma:strong_convexity_gp} as a form of strong convexity for our objective.

\begin{lemma}[Strong convexity of $f$]
\label{lemma:strong_convexity_gp}
Let $f(\vx) \coloneqq \gnorm{\mA\vx-\vb}{p}^p$. For all $\vd\in\R^d$, we have
\begin{align*}
    f(\vx+\vd) \ge f(\vx) + \ip{\nabla f(\vx),\vd} + \frac{4}{2^p}\gnorm{\mA\vd}{p}^p,
\end{align*}
and therefore
\begin{align*}
    \norm{\vx-\xstar}_{\mM} \le 2^{3/2 - 3/p}d^{1/2-1/p}(f(\vx)-f(\xstar))^{1/p}\enspace.
\end{align*}
\end{lemma}

\lemmastrongconvexitycomponent

To motivate \Cref{lemma:strong_convexity_component}, let us see how \Cref{lemma:strong_convexity_component} implies \Cref{lemma:strong_convexity_gp}.

\begin{proof}[Proof of \Cref{lemma:strong_convexity_gp}]
Note that
\begin{align*}
    \nabla f(\vx) = \sum_{i=1}^m p\norm{\mA_{S_i}\vx-\vb_{S_i}}_2^{p-2}\mA_{S_i}^{\top}(\mA_{S_i}\vx-\vb_{S_i})\enspace.
\end{align*}
This implies
\begin{align*}
    \sum_{i=1}^m p\norm{\mA_{S_i}\vx-\vb_{S_i}}_2^{p-2}\ip{\mA_{S_i}\vx-\vb_{S_i},\mA_{S_i}\vd} = \ip{\nabla f(\vx),\vd}\enspace.
\end{align*}
Combining this and applying \Cref{lemma:strong_convexity_component} (which is a strong convexity lemma for $\|\cdot\|_2^p$ that we prove subsequently in this section), we get
\begin{align*}
    f(\vx+\vd) &= \gnorm{\mA(\vx+\vd)-\vb}{p}^p = \gnorm{\mA\vd+(\mA\vx-\vb)}{p}^p\enspace,\\
    &= \sum_{i=1}^m \norm{\mA_{S_i}\vd+(\mA_{S_i}\vx-\vb_{S_i})}_2^p\enspace,\\
    &\ge^{\text{(\Cref{lemma:strong_convexity_component})}} \sum_{i=1}^m \norm{\mA_{S_i}\vx-\vb_{S_i}}_2^p + p\|\mA_{S_i}\vx-\vb_{S_i}\|_2^{p-2}\ip{(\mA_{S_i}\vx-\vb_{S_i}),\mA_{S_i}\vd}+ \frac{4}{2^p}\norm{\mA_{S_i}\vd}_2^p\enspace,\\
    &= \sum_{i=1}^m \norm{\mA_{S_i}\vx-\vb_{S_i}}_2^p + \ip{p\|\mA_{S_i}\vx-\vb_{S_i}\|_2^{p-2}\mA_{S_i}^{\top}(\mA_{S_i}\vx-\vb_{S_i}),\vd}+ \frac{4}{2^p}\norm{\mA_{S_i}\vd}_2^p\enspace,\\
    &=^{\text{\eqref{eq:f_grad}}} \gnorm{\mA\vx-\vb}{p}^p + \ip{\nabla f(\vx),\vd} + \frac{4}{2^p}\gnorm{\mA\vd}{p}^p = f(\vx) + \ip{\nabla f(\vx),\vd} + \frac{4}{2^p}\gnorm{\mA\vd}{p}^p\enspace.
\end{align*}
We now take care of the second statement. Observe that at optimality, we have $\nabla f(\xstar) = 0$. Plugging this in (replace $\vx$ by $\vx^\star$ and $\vd$ by $\vx-\vx^\star$ above), rearranging, and taking $p$th roots gives
\begin{align*}
    \gnorm{\mA(\vx-\xstar)}{p} \le \left(\frac{4}{2^p}\right)^{-1/p}(f(\vx)-f(\xstar))^{1/p} = \frac{2}{4^{1/p}}(f(\vx)-f(\xstar))^{1/p}\enspace. 
\end{align*}
Next, recall that by \Cref{thm:gp_regression_blw_intro},
\begin{align*}
    \norm{\vx-\xstar}_{\mM} = \norm{\mW^{1/2-1/p}\mA(\vx-\xstar)}_2 \le (2d)^{1/2-1/p}\gnorm{\mA(\vx-\xstar)}{p}\enspace.
\end{align*}
Stitching the inequalities together completes the proof of \Cref{lemma:strong_convexity_gp}.
\end{proof}

In the rest of this subsection, we prove \Cref{lemma:strong_convexity_component}. We begin with a few numerical inequalities. 

\begin{lemma}
\label{lemma:gp_sc_galpha}
For $\alpha \le -1/2$ and $p\ge 2$, $g(\alpha) := \frac{1+p\alpha}{(-(2\alpha+1))^{p/2}}$ is nonincreasing in $\alpha$.
\end{lemma}
\begin{proof}[Proof of \Cref{lemma:gp_sc_galpha}]
We first take the derivative of $g$ with respect to $\alpha$,
\begin{align*}
    g'(\alpha) &= \frac{p(-(2\alpha+1))^{p/2} - \inparen{(-2)\frac{p}{2}\inparen{-(2\alpha+1)}^{p/2-1}}(1+p\alpha)}{(-(2\alpha+1))^p}\enspace,\\
    &= \frac{p(-(2\alpha+1)^{p/2}) + p\inparen{-(2\alpha+1)}^{p/2-1}(1+p\alpha)}{(-(2\alpha+1))^p}\enspace,\\
    &= p\cdot\frac{(-(2\alpha+1)) + (1+p\alpha)}{(-(2\alpha+1))^{p/2+1}}\enspace,\\
    &= p\cdot\frac{(p-2)\alpha}{(-(2\alpha+1))^{p/2+1}} \le 0\enspace,
\end{align*}
where in the final inequality we used that $p\geq 2$ and $\alpha\leq -1/2$. This completes the proof of the lemma.
\end{proof}
We also need the following lemma, which is similar to a result due to \citeauthor{akps19} \cite[Lemma 4.5]{akps19}. It amounts to proving \Cref{lemma:strong_convexity_component} when the dimension $k=1$.
\begin{lemma}[Case A. of \Cref{lemma:gp_sc_alphaall}]\label{lemma:adil_sc_lemma}
For any $\alpha\in\R$ and $p \ge 2$,
\begin{align*}
    \abs{1+\alpha}^p \ge 1+p\alpha+\frac{4}{2^p}\abs{\alpha}^p\enspace.
\end{align*}
\end{lemma}
\begin{proof}[Proof of \Cref{lemma:adil_sc_lemma}]
Note that the inequality is true when $p=2$ and becomes an equality. We consider the case when $p>2$ and use $h(\alpha)$ to denote the error function,
\begin{align*}
    h(\alpha) \coloneqq \abs{1+\alpha}^p - \inparen{1+p\alpha+\frac{4}{2^p}\abs{\alpha}^p}\enspace.
\end{align*}
We aim to show $h(\alpha) \ge 0$ for all $\alpha \in \R$. Let us first write the derivatives of $h$.
\begin{align*}
    h'(\alpha) &= p\inparen{\abs{1+\alpha}^{p-2}(1+\alpha) - \inparen{1 + \frac{4}{2^p}\abs{\alpha}^{p-2}\alpha}}\enspace,\\
    h''(\alpha) &= p(p-1)\inparen{\abs{1+\alpha}^{p-2}-\frac{4}{2^p}\abs{\alpha}^{p-2}} = p(p-1)\inparen{\abs{1+\alpha}^{p-2}-\abs{\frac{\alpha}{2}}^{p-2}}\enspace.
\end{align*}
It is now easy to verify the following statements about $h$,
\begin{itemize}
    \item[I.] $h'(-2) = h''(-2) = 0$ and $h''(\alpha) > 0$ for $\alpha < -2$, $\Rightarrow$ within the range $(-\infty,-2]$ the function $h$ is minimized at $-2$;
    \item[II.] $h'(-2)=0$ and $h''(\alpha) \leq 0$ for $\alpha\in (-2, -2/3]$ $\Rightarrow$ $h'(\alpha) < 0$ in the range $(-2, -2/3]$, i.e., in that range the function $h$ is minimized at $-2/3$;
    \item[III.] $h'(-2/3) < 0 = h'(0)$ and $h''(\alpha) > 0$ for $\alpha > -2/3$ $\Rightarrow$ the function $h$ is decreasing in $(-2/3,0)$ and increasing in $[0, \infty)$, i.e., within the range $(-2/3, \infty)$ the function $h$ is minimized at $0$.
\end{itemize}
As a result of the above observations, it is enough to check the inequality at the inputs $\alpha \in \inbraces{-2, -2/3, 0}$. We have for $p>2$,
\begin{align*}
    h(-2) &= 1 - \inparen{1-2p+4} = 2p -4 > 0\enspace,\\
    h\inparen{-\frac{2}{3}} &= \frac{1}{3^p} - \inparen{1-\frac{2p}{3}+\frac{4}{2^p}\abs{\frac{2}{3}}^p} = \frac{1}{3^p}-1+\frac{2p}{3}-\frac{4}{3^p} = -1 + \frac{2p}{3}-\frac{3}{3^p} > 0\\
    h(0) &= 1 - 1 = 0\enspace.
\end{align*}
This implies that $h(\alpha)\geq 0$ for all values of $\alpha$, concluding the proof of \Cref{lemma:adil_sc_lemma}.
\end{proof}

Next, we prove a special case of \Cref{lemma:strong_convexity_component}.

\begin{lemma}\label{lemma:gp_sc_alphaall}
    For any $\alpha\in \R$, $\beta\ge 0$, and $p \ge 2$, we have
\begin{align*}
    \inparen{(1+\alpha)^2+\beta^2}^{p/2} \ge 1+p\alpha + \frac{4}{2^p}\inparen{\alpha^2+\beta^2}^{p/2}\enspace.
\end{align*}
\end{lemma}
\begin{proof}[Proof of \Cref{lemma:gp_sc_alphaall}]
    Let us study the difference of both sides of the inequality using the following function,
    \begin{align*}
        h(\alpha, \beta) \coloneqq \inparen{(1+\alpha)^2+\beta^2}^{p/2} - \left(1+p\alpha + \frac{4}{2^p}\inparen{\alpha^2+\beta^2}^{p/2}\right)\enspace.
    \end{align*}
    We want to show that for $\alpha\in \R$, $\beta\ge 0$, and $p \ge 2$, $h(\alpha, \beta) \geq 0$. We will break this proof into three cases: \textbf{A.} $\alpha \in \R$ and $\beta =0$; \textbf{B.} $\alpha \in (-\infty, -2]\cup[-2/3,\infty)$ and $\beta >0$; and \textbf{C.} $\alpha \in (-2, -2/3)$ and $\beta>0$. These cases together cover of the entire range of $\alpha\in \R$ and $\beta\geq 0$.
    
    \paragraph{\underline{Case A.}} When $\beta=0$, the proof simply follows from the statement of \Cref{lemma:adil_sc_lemma} by noting $|\alpha|^{p} = (\sqrt{\alpha^2})^{p} = (\alpha^2)^{p/2}$. 
    
    In the remaining two cases we will show that for any $\alpha\in \R$, increasing the value of $\beta$ still maintains $h(\alpha, \beta) \geq 0$. To see this, we first note that the derivative of $h(\alpha, \beta)$ w.r.t. $\beta$ is given by,
    \begin{align*}
        \nabla_\beta h(\alpha, \beta) &= p\beta\left(\inparen{(1+\alpha)^2+\beta^2}^{p/2-1} - \frac{4}{2^p}\inparen{\alpha^2+\beta^2}^{p/2-1}\right)\enspace.
    \end{align*}
    For $\beta > 0$, ensuring this derivative is positive is equivalent to the following,
    \begin{align}
        \nabla_\beta h(\alpha, \beta) > 0 &\equiv p\beta\inparen{(1+\alpha)^2+\beta^2}^{p/2-1} > p\beta \cdot \frac{4}{2^p}\inparen{\alpha^2+\beta^2}^{p/2-1}\enspace,\nonumber\\
        &\equiv^{(p\beta > 0)} (1+\alpha)^2 + \beta^2  > \left(\frac{1}{2^{p-2}}\right)^{2/(p-2)}\cdot \left(\alpha^2 + \beta^2\right)\enspace,\nonumber\\
        &\equiv (1+\alpha)^2 + \beta^2  > \frac{1}{4}\cdot \left(\alpha^2 + \beta^2\right)\enspace,\nonumber\\
        &\equiv (3\alpha^2 + 8\alpha + 4) + 3\beta^2  > 0\enspace,\nonumber\\
        &\equiv \beta^2 > -\left(\alpha^2 + \frac{8}{3}\alpha + \frac{4}{3}\right)\enspace.\label{eq:grad_pos_equiv}
    \end{align}
    
    \paragraph{\underline{Case B.}} Note that the roots of the quadratic function $3\alpha^2 + 8\alpha + 4$ are given by $\alpha_1 = -2$ and $\alpha_2 = -2/3$. This means that for $\alpha \in (-\infty, -2]\cup[-2/3,\infty)$ we have $3\alpha^2 + 8\alpha + 4 \geq 0$ which is \textbf{sufficient} to ensure using \eqref{eq:grad_pos_equiv} that $\nabla_\beta h(\alpha, \beta) > 0$, and hence $h(\alpha, \beta) > 0$. This takes care of Case B. 
    
    \paragraph{\underline{Case C.}} Now we only need to consider the range $\alpha \in (-2,-2/3)$ with $\beta>0$. In this range, the recall the equivalence \eqref{eq:grad_pos_equiv}, 
    \begin{align*}
        \nabla_\beta h(\alpha, \beta) > 0 &\equiv \beta > \sqrt{-\left(\alpha^2 + \frac{8}{3}\alpha + \frac{4}{3}\right)} =: \beta_0(\alpha)\enspace.
    \end{align*}
    Thus for all $\beta > \beta_0(\alpha)$ we know that $h(\alpha, \beta)$ is increasing in $\beta$ and vice-versa. This allows us for any given $\alpha\in (-2,-2/3)$ to further break Case C into two sub-cases:
    
    \paragraph{\underline{Case C.I}} For $\beta \in [0, \beta_0)$, since $h(\alpha, \beta)$ is decreasing in $\beta$ its lowest value is attained at $\beta=0$ and we only need to verify that $h(\alpha, 0) \geq 0$. We get this directly from \Cref{lemma:adil_sc_lemma}.
    
    \paragraph{\underline{Case C.II}} For $\beta \in [\beta_0, \infty)$, since $h(\alpha, \beta)$ is increasing in $\beta$ its lowest value is attained at $\beta=\beta_0$ and we only need to verify that $h(\alpha, \beta_0(\alpha)) \geq 0$. We first simplify the expression for $h(\alpha, \beta_0(\alpha))$,
    \begin{align*}
        h(\alpha, \beta_0(\alpha)) &=  \inparen{(1+\alpha)^2+\beta_0^2}^{p/2} - \left(1+p\alpha + K_p\inparen{\alpha^2+\beta_0^2}^{p/2}\right)\enspace,\\
        &= \inparen{-\frac{1}{3}-\frac{2}{3}\alpha}^{p/2} - \left(1+p\alpha + \frac{4}{2^p}\inparen{-\frac{8}{3}\alpha -\frac{4}{3}}^{p/2}\right)\enspace,\\
        &= \inparen{-\frac{1}{3}-\frac{2}{3}\alpha}^{p/2} - \left(1+p\alpha + 4\inparen{-\frac{2}{3}\alpha -\frac{1}{3}}^{p/2}\right)\enspace,\\
        &= -1 - p\alpha - 3\inparen{-\frac{2}{3}\alpha -\frac{1}{3}}^{p/2}\enspace,\\
        &= -1 - p\alpha - \frac{1}{3^{p/2-1}}(-2\alpha -1)^{p/2}\enspace,\\
        &= -(-2\alpha-1)^{p/2}\left(\frac{1+p\alpha}{(-2\alpha-1)^{p/2}} + \frac{1}{3^{p/2-1}}\right)\enspace. 
    \end{align*}
    Now since $\alpha\in (-2,-2/3)<-1/2$ we can use \Cref{lemma:gp_sc_galpha} to note that the first term is non-decreasing in $\alpha$ which means that its lowest value in this range can be lower bounded by its value at $\alpha = -2$, i.e., for $\alpha\in (-2,-2/3)$,
    \begin{align*}
        h(\alpha, \beta_0(\alpha)) &\geq h(-2, \beta_0(-2))\enspace,\\ 
        &= -3^{p/2}\left(\frac{1-2p}{3^{p/2}} + \frac{1}{3^{p/2-1}}\right)\enspace,\\
        &= 2p-1-3 = 2(p-2) > 0\enspace,
    \end{align*}
    which finishes the proof of Case C.II and also Case C. Together Cases A, B and C complete the proof of \Cref{lemma:gp_sc_alphaall}.  
\end{proof}

We are now ready to prove \Cref{lemma:strong_convexity_component}.

\begin{proof}[Proof of \Cref{lemma:strong_convexity_component}]
First, assume that $\norm{\vv}_2=1$. We will later extend the result to all $\vv$.

Since $\norm{\vv}_2=1$, we can write $\triangle = \alpha\vv+\beta\vw$ where $\ip{\vv,\vw} = 0$ and $\norm{\vw}_2=1$, so that we have $\norm{\triangle}_2^2=\alpha^2+\beta^2$. Without loss of generality, we have $\beta \ge 0$. Fixing $\vw$ and $\alpha$ for now, it is enough to show that for all $\beta \ge 0$, we have
\begin{align*}
    \norm{(1+\alpha)\vv+\beta\vw}_2^p = \inparen{(1+\alpha)^2+\beta^2}^{p/2} \overset{?}{\ge} 1+p\alpha + \frac{4}{2^p}\norm{\triangle}_2^p = 1+p\alpha + \frac{4}{2^p}\inparen{\alpha^2+\beta^2}^{p/2}.
\end{align*}
This follows immediately by \Cref{lemma:gp_sc_alphaall}.

We now extend the result for all $\vv$. Let $\bar{\vv} \coloneqq \vv/\norm{\vv}_2$ and note that
\begin{align*}
    \norm{\vv+\triangle}_2^p = \norm{\vv}_2^p\norm{\bar{\vv}+\frac{\triangle}{\norm{\vv}_2}}_2^p &\ge \norm{\vv}_2^p\inparen{1+\ip{\bar{\vv},\frac{\triangle}{\norm{\vv}_2}}+\frac{4}{2^p}\norm{\frac{\triangle}{\norm{\vv}_2}}_2^p} \\
    &= \norm{\vv}_2^p+p\norm{\vv}_2^{p-2}\ip{\vv,\triangle}+\frac{4}{2^p}\norm{\triangle}_2^p,
\end{align*}
completing the proof of \Cref{lemma:strong_convexity_component}.
\end{proof}

\subsubsection{Smoothness of the objective}

The main result of this subsection is \Cref{lemma:gp_smooth}.

\begin{lemma}
\label{lemma:gp_smooth}
For all $\vx\in\R^d$, we have
\begin{align*}
    f(\vx)-f(\xstar) \le \frac{p(p-1)}{2}f(\vx)^{1-\frac{2}{p}}\gnorm{\mA(\vx-\xstar)}{p}^2.
\end{align*}
\end{lemma}
\begin{proof}[Proof of \Cref{lemma:gp_smooth}]
By Taylor's/mean-value theorem, we can write for some $\vy$ on the line connecting $\xstar$ and $\vx$,
\begin{align*}
    f(\vx) &= f(\xstar) + \ip{\nabla f(\xstar), \vx-\xstar} + \frac{1}{2}(\vx-\xstar)^{\top}\nabla^2 f(\vy)(\vx-\xstar) \\
    &\le^{\eqref{eq:f_quad_form}} f(\xstar) + \frac{p(p-1)}{2}\sum_{i=1}^m \norm{\mA_{S_i}\vy-\vb_{S_i}}_2^{p-2}\norm{\mA_{S_i}(\vx-\xstar)}_2^2 \\
    &\le f(\xstar) + \frac{p(p-1)}{2} \inparen{\sum_{i=1}^m \norm{\mA_{S_i}\vy-\vb_{S_i}}_2^{p}}^{\frac{p-2}{p}}\inparen{\sum_{i=1}^m \norm{\mA_{S_i}(\vx-\xstar)}_2^p}^{\frac{2}{p}} \\
    &\le f(\xstar) + \frac{p(p-1)}{2}f(\vx)^{1-\frac{2}{p}}\gnorm{\mA(\vx-\xstar)}{p}^2,
\end{align*}
completing the proof of \Cref{lemma:gp_smooth}.
\end{proof}

\subsection{Facts about the iterates}
\label{sec:gp_iterate_facts}

The main result of this section is \Cref{lemma:fq_bounded}. In words, \Cref{lemma:fq_bounded} tells us that each proximal query we make in \Cref{alg:optimal_ms_acceleration} (see Line \ref{line:optimal_ms_acceleration_query} of \Cref{alg:optimal_ms_acceleration}) has bounded objective value. We will need this later when we argue about the convergence rates for the algorithms used to solve the proximal subproblems.

\begin{lemma}
\label{lemma:fq_bounded}
For all queries $\vq_t$, we have
\begin{align*}
    f(\vq_t) \le f(\vx_t) + \inparen{9p(p-1)}^{\frac{p}{2}}d^{\frac{p}{2}-1}.
\end{align*}
\end{lemma}
\begin{proof}[Proof of \Cref{lemma:fq_bounded}]
We establish the following upper bound on $f(\vv_t) - f(\xstar)$ using the ingredients developed so far:
\begin{align*}
    f(\vv_t) - f(\xstar) &\le \frac{p(p-1)}{2}f(\vv_t)^{1-\frac{2}{p}}\gnorm{\mA(\vv_t-\xstar)}{p}^2 & \text{(\Cref{lemma:gp_smooth})} \\
    &\le \frac{p(p-1)}{2}f(\vv_t)^{1-\frac{2}{p}}\norm{\vv_t-\xstar}_{\mM}^2 & \text{(\Cref{thm:gp_regression_blw_intro})} \\
    &\le p(p-1)f(\vv_t)^{1-\frac{2}{p}}\norm{\vx_0-\xstar}_{\mM}^2 & \text{(\Cref{lemma:ms_acceleration_iterate_diameter})}\\
    &\le p(p-1)f(\vv_t)^{1-\frac{2}{p}}2^2(2d)^{1-\frac{2}{p}} & \text{(\Cref{thm:gp_regression_blw_intro})} \\
    &\leq 8d^{1-\frac{2}{p}}p(p-1)f(\vv_t)^{1-\frac{2}{p}}\enspace.
\end{align*}
Now, recall that we assume by rescaling that $f(\xstar)=1$. From this, it trivially follows that $1 \le d^{1-\frac{2}{p}}p(p-1)f(\vv_t)^{1-\frac{2}{p}}$. Combining these and re-arranging the above inequality leads to the following  polynomial inequality in $f(\vv_t)$,
\begin{align}
    0 &\ge f(\vv_t)-8d^{1-\frac{2}{p}}p(p-1)f(\vv_t)^{1-\frac{2}{p}} -1\enspace,\nonumber\\ 
    &= f(\vv_t)-9d^{1-\frac{2}{p}}p(p-1)f(\vv_t)^{1-\frac{2}{p}} +  d^{1-\frac{2}{p}}p(p-1)f(\vv_t)^{1-\frac{2}{p}} -1\enspace,\nonumber\\
    &\ge f(\vv_t)-9d^{1-\frac{2}{p}}p(p-1)f(\vv_t)^{1-\frac{2}{p}}\enspace,\label{eq:gp_fq_upper}
\end{align}
where in the last inequality we used the fact that the optimal value $f(\vx^\star)=1$ (due to our rescaling), which implies that for $p\geq 2$,
$$1 \leq f(\vv_t) \leq d^{1-\frac{2}{p}}p(p-1)f(\vv_t)^{1-\frac{2}{p}}\enspace.$$
Solving for $f(\vv_t)$ in \eqref{eq:gp_fq_upper}, we get
\begin{align*}
    f(\vv_t) \le \inparen{9p(p-1)}^{\frac{p}{2}}d^{\frac{p}{2}-1}\enspace.
\end{align*}
Using the definition of $\vq_t$ from \Cref{alg:optimal_ms_acceleration} (Line \ref{line:optimal_ms_query}) along with the convexity of $f$ (Jensen's inequality), and using our bound on $f(\vv_t)$ we note that,
\begin{align*}
    f(\vq_t) &\le f(\vx_t) + f(\vv_t)\enspace,\\
    &\leq f(\vx_t) + \inparen{9p(p-1)}^{\frac{p}{2}}d^{\frac{p}{2}-1}\enspace,
\end{align*}
which completes the proof of \Cref{lemma:fq_bounded}.
\end{proof}

\subsection{Proximal subproblems -- calculus, algorithms, proofs}
\label{sec:gp_prox_solver}

Let
\begin{align*}
    f_{\vq_t}(\xtilde) \coloneqq f(\xtilde)+ep^p\norm{\xtilde-\vq_t}_{\mM}^p\enspace.
\end{align*}
In this subsection, we design and analyze an algorithm (\Cref{alg:gp_prox_solver}) that approximately solves the subproblem
\begin{align*}
    \argmin{\xtilde\in\R^d} f_{\vq_t}(\xtilde).
\end{align*} 
Specifically, we will output $(\xtilde_{t+1},\lambda_{t+1})$ that satisfy the $\frac{1}{2}$-MS oracle condition (\Cref{defn:ms_oracle}) and an appropriate movement bound (\Cref{defn:movement_bound}).

This subproblem is the workhorse of \Cref{alg:gp_regression_alg}, and once we implement and analyze the solver, it is very straightforward to plug this into \Cref{alg:optimal_ms_acceleration} and \Cref{thm:optimal_ms_acceleration} to get our final iteration complexity.

\begin{algorithm}[H]
\caption{\textsf{GpRegressionProxOracle}: Implements $\frac{1}{2}$-MS oracle for $\gnorm{\cdot}{p}$ regression (see \Cref{lemma:prox_subproblem_ms} and \Cref{alg:mirror_descent}.}
\label{alg:gp_prox_solver}
\begin{algorithmic}[1]
\Require Query $\vq_t$, previous iterate $\vx_t$, intended parameter distance $\gamma$.
\State Define \begin{equation*}
    \begin{aligned}
    f_{\vq_t}(\xtilde) &\coloneqq f(\xtilde) + ep^p\norm{\xtilde-\vq_t}_{\mM}^p \\
    h_{\vq_t}(\xtilde) &\coloneqq \norm{\xtilde-\vq_t}_{\nabla^2 f(\vq_t)}^2 + ep^p\norm{\xtilde-\vq_t}_{\mM}^p \\
    D_{h_{\vq_t}}(\vx,\vy) &\coloneqq h_{\vq_t}(\vx)-h_{\vq_t}(\vy) - \ip{\nabla h_{\vq_t}(\vy),\vx-\vy} \\
    \xtilde_{\vq_t} &\coloneqq \argmin{\xtilde\in\R^d} f_{\vq_t}(\xtilde)
    \end{aligned}.
\end{equation*}
\State Let $T \ge Cp^{O(1)}e\logv{dpeh_{\vq_t}(\xtilde_{\vq_t})\inparen{\frac{4}{p\gamma}}^p}$.
\State Run \Cref{alg:mirror_descent} with input iteration count $T$, base function $f_{\vq_t}$, reference function $h_{\vq_t}$, and initialization $\vq_t$.
\end{algorithmic}
\end{algorithm}

The goal of the rest of this section is to analyze \Cref{alg:gp_prox_solver}. The analysis follows several steps:
\begin{enumerate}
    \item We find a reference function $h_{\vq_t}$ that depends on the query point $\vq_t$ for which the proximal objective $f_{\vq_t}$ is relatively smooth and relatively strongly convex with $O(p^{O(1)})$ condition number (see \Cref{sec:mirror_descent} for a sense of why this is useful). The main result here is \Cref{lemma:gp_hessian_stable}.
    \item We show that $f_{\vq_t}$ is strongly convex, following from \Cref{lemma:strong_convexity_component}. This will help us understand the argument suboptimality for any point that approximately optimizes $f_{\vq_t}$ in function value. We also show that the reference function $h_{\vq_t}$ is strongly convex, using the same tools, for the same reason.
    \item We show a form of smoothness for $f_{\vq_t}$. This helps us bound the gradient of any point that approximately optimizes $f_{\vq_t}$. Combining these later will tell us that an approximate solution to $f_{\vq_t}$ in argument value is also an approximate stationary point, i.e., it satisfies the $\frac{1}{2}$-MS condition (\Cref{defn:ms_oracle}).
    \item We solve the proximal subproblems. This solution itself follows a few steps:\begin{enumerate}
        \item We apply \Cref{thm:gp_mirror_descent}. This tells us that as long as we can approximately solve the Bregman proximal problems (approximately implementing Line \ref{line:mirror_descent_exact} in \Cref{alg:mirror_descent}), we will be in good shape.
        \item This means we have to figure out how to approximately solve problems of the form $\argmin{\vx\in\R^d} \ip{\vg,\vx} + Lh_{\vq_t}(\vx)$, where $L$ is the smoothness constant derived for $f_{\vq_t}$ with respect to $h_{\vq_t}$. We do this up to an accuracy that approximate mirror descent can handle (see \Cref{thm:gp_mirror_descent} for details on what we want this approximation to look like). For the approximation to work, we need to approximately solve this problem up to both argument accuracy and approximate stationarity. The main technical result of interest here is \Cref{lemma:gp_prox_relativesmoothsolve}.
    \end{enumerate}
    \item We use the smoothness and strong convexity guarantees to show that our solution from the previous step satisfies the $\frac{1}{2}$-MS oracle (\Cref{defn:ms_oracle}), which means we can plug-and-play into \Cref{thm:optimal_ms_acceleration}.
\end{enumerate}

\subsubsection{Hessian stability}

Throughout this section, we adopt the following notation:
\begin{align*}
    C_p &\coloneqq ep^{p} \\
    f(\vx) &\coloneqq \sum_{i=1}^m \norm{\mA_{S_i}\vx-\vb_{S_i}}_2^p \\
    f_{\vq}(\vx) &\coloneqq f(\vx) + C_p\norm{\vx-\vq}_{\mM}^p \\
    h_{\vq}(\vx) &\coloneqq \norm{\vx-\vq}_{\nabla^2 f(\vq)}^2 + C_p\norm{\vx-\vq}_{\mM}^p
\end{align*}
We begin with proving our Hessian stability fact, which should also be equivalently viewed as showing that $f_{\vq_t}$ is relatively smooth and relatively strongly convex in $h_{\vq_t}$ with $O(p^{O(1)})$ condition number. Our main result is \Cref{lemma:gp_hessian_stable} which relies on analytical results \Cref{lem:small_lemma_for_hessian_stability_one} and \Cref{lem:small_lemma_for_hessian_stability_two} that we prove later.
\begin{lemma}
\label{lemma:gp_hessian_stable}
    For all $\vx \in\R^d$ and $p\geq 2$, we have $$\frac{1}{2p\cdot e}\nabla^2h_{\vq}(\vx)\preceq \nabla^2f_{\vq}(\vx) \preceq p\cdot e\nabla^2 h_{\vq}(\vx)\enspace.$$
\end{lemma}
\begin{proof}[Proof of \Cref{lemma:gp_hessian_stable}]
    Using an arbitrary $\vz\in\mathbb{R}^d$ we can write the following quadratic form of the hessian of $f$,
    \begin{align}
        \vz^{\top}\nabla^2 f(\vx)\vz &\le^{\text{(a)}} p\cdot(p-1)\sum_{i=1}^m\norm{\mA_{S_i}\vx - \vb_{S_i}}_2^{p-2}\norm{\mA_{S_i}\vz}_2^2\enspace,\nonumber\\
        &= p\cdot(p-1)\sum_{i=1}^m\norm{\mA_{S_i}(\vx - \vq) + \mA_{S_i}\vq - \vb_{S_i}}_2^{p-2}\norm{\mA_{S_i}\vz}_2^2\enspace,\nonumber\\
        &\leq^{\text{(b)}}  p\cdot(p-1)\sum_{i=1}^m\left(\alpha_p^{p-2}\norm{\mA_{S_i}(\vx - \vq)}_2^{p-2}\norm{\mA_{S_i}\vz}_2^2 + \beta_p^{p-2}\norm{\mA_{S_i}\vq - \vb_{S_i}}_2^{p-2}\norm{\mA_{S_i}\vz}_2^2\right)\enspace,\nonumber\\
        &\leq^{\text{(c)}}p\cdot(p-1)\cdot\alpha_p^{p-2}\sum_{i=1}^m\norm{\mA_{S_i}(\vx - \vq)}_2^{p-2}\norm{\mA_{S_i}\vz}_2^2 + (p-1)\cdot\beta_p^{p-2}\vz^{\top}\nabla^2f(\vq)\vz\enspace,\nonumber\\
        &\leq^{\text{(d)}}p\cdot(p-1)\cdot\alpha_p^{p-2}\left(\norm{\vx-\vq}^p_{\mM}\right)^{(p-2)/p}\left(\norm{\vz}^p_{\mM}\right)^{2/p}+ (p-1)\cdot\beta_p^{p-2}\vz^{\top}\nabla^2f(\vq)\vz\enspace,\nonumber\\
        &= p\cdot(p-1)\cdot\alpha_p^{p-2}\norm{\vx-\vq}^{p-2}_{\mM}\norm{\vz}^2_{\mM} + (p-1)\cdot\beta_p^{p-2}\vz^{\top}\nabla^2f(\vq)\vz\enspace,\nonumber\\
        &\leq^{(e)} \frac{(p-1)\cdot \alpha_p^{p-2}}{C_p} \vz^{\top}\nabla^2g_{\vq}(\vx)\vz + (p-1)\cdot\beta_p^{p-2}\vz^{\top}\nabla^2f(\vq)\vz\enspace,\label{eq:quad_form_f_ub}
    \end{align}
    where in (a) we apply the upper bound from \Cref{lemma:gp_regression_hessian_upperbound}, in (b) we pick $\alpha_p,\beta_p\geq 1$ such that $1/\alpha_p + 1/\beta_p = 1$ (we will choose them later), in (c) we apply the lower bound from \Cref{lemma:gp_regression_hessian_upperbound}, in (d) we use the choice of our weights in designing $\mM$ and \Cref{thm:gp_regression_blw_intro}\todo{write something more explicit for justifying (d).} and finally in (e) we use the following calculations for the regularizer term for some $\vz\in\mathbb{R}^d$,
    \begin{align*}
        g_{\vq}(\vx) &\coloneqq C_p\norm{\vx-\vq}^p_{\mM}\enspace,\\    
        \nabla g_{\vq}(\vx) &= pC_p\norm{\vx-\vq}_{\mM}^{p-2}{\mM}(\vx-\vq)\enspace,\\
        \nabla^2 g_{\vq}(\vx) &= pC_p\norm{\vx-\vq}_{\mM}^{p-2}\mM + p(p-2)C_p\norm{\vx-\vq}_{\mM}^{p-4}\mM(\vx-\vq)(\vx-\vq)^{\top}\mM\enspace,\\
        \vz^{\top}\nabla^2 g_{\vq}(\vx)\vz &= pC_p\norm{\vx-\vq}_{\mM}^{p-2}\norm{\vz}^2_{\mM}+ p(p-2)C_p\norm{\vx-\vq}_{\mM}^{p-4}\left((\vx-\vq)^{\top}\mM\vz\right)^2 \geq^{(p\geq 2)} 0\enspace.
    \end{align*}
    Combining \eqref{eq:quad_form_f_ub} with the definition of $f_{\vq}$ gives us,
    \begin{align*}
        \vz^{\top}\nabla^2f_{\vq}(\vx)\vz &= \vz^{\top}\nabla^2f(\vx)\vz + \vz^{\top}\nabla^2g_{\vq}(\vx)\vz\enspace,\\
        &\leq^{\text{using }\eqref{eq:quad_form_f_ub}} (p-1)\cdot\beta_p^{p-2}\vz^{\top}\nabla^2f(\vq)\vz + \left(1+\frac{(p-1)\cdot \alpha_p^{p-2}}{C_p}\right) \vz^{\top}\nabla^2g_{\vq}(\vx)\vz\enspace.
    \end{align*}
    Thus, in order to finish the proof for the upper bound we need to pick $\alpha_p, \beta_p$. We split the analysis here into two cases: \textbf{A.} $p>2$ and \textbf{B.} $p=2$.
    
    \paragraph{\underline{Case A.} ($p>2$)} For simplicity we will just pick $\alpha_p = p-1$ and $\beta_p = \frac{p-1}{p-2}$ which implies,
    \begin{align*}
        \vz^{\top}\nabla^2f_{\vq}(\vx)\vz &\leq (p-1)\cdot\left(1 + \frac{1}{p-2}\right)^{p-2}\vz^{\top}\nabla^2f(\vq)\vz + \left(1+\frac{(p-1)\cdot (p-1)^{p-2}}{C_p}\right) \vz^{\top}\nabla^2g_{\vq}(\vx)\vz\enspace,\\
        &\leq (p-1)\cdot e\vz^{\top}\nabla^2f(\vq)\vz + \left(1+\frac{(p-1)^{p-1}}{C_p}\right) \vz^{\top}\nabla^2g_{\vq}(\vx)\vz\enspace,\\
        &= \frac{(p-1)\cdot e}{2} \vz^{\top}\left(\nabla^2h_{\vq}(\vx) - \nabla^2g_{\vq}(\vx)\right)\vz + \left(1+\frac{(p-1)^{p-1}}{C_p}\right) \vz^{\top}\nabla^2g_{\vq}(\vx)\vz\enspace,\\
        &\leq^{(p\geq 2)} p\cdot e \vz^{\top}\nabla^2h_{\vq}(\vx)\vz + \left(1+\frac{(p-1)^{p-1}}{C_p} - \frac{(p-1)\cdot e}{2}\right)\vz^{\top}\nabla^2g_{\vq}(\vx)\vz\enspace,\\
        &= p\cdot e \vz^{\top}\nabla^2h_{\vq}(\vx)\vz + \left(1+\frac{(p-1)^{p-1}}{ep^p} - \frac{(p-1)\cdot e}{2}\right)\vz^{\top}\nabla^2g_{\vq}(\vx)\vz\enspace,\\
        &\leq^{\text{(\Cref{lem:small_lemma_for_hessian_stability_one})}} p\cdot e \vz^{\top}\nabla^2h_{\vq}(\vx)\vz\enspace,
    \end{align*}
    where in the final inequality we use \Cref{lem:small_lemma_for_hessian_stability_one} which tell us that for $p\geq 2$ the constant in front of $\vz^{\top}\nabla^2g_{\vq}(\vx)\vz$ is negative along with the fact that $\vz^{\top}\nabla^2g_{\vq}(\vx)\vz$ is non-negative. To get the lower bound we first exchange $\vx, \vq$ in \eqref{eq:quad_form_f_ub} (and use the values of $\alpha_p$ and $\beta_p$) to get,
    \begin{align*}
        &\vz^{\top}\nabla^2f(\vq)\vz \leq \frac{(p-1)\cdot(p-1){p-2}}{ep^p}\vz^{\top}\nabla^2g_{\vx}(\vq)\vz + (p-1)\left(1 + \frac{1}{p-2}\right)^{p-2}\vz^{\top}\nabla^2f(\vx)\vz\enspace,\\
        &\Rightarrow \vz^{\top}\nabla^2f(\vq)\vz \leq \frac{(p-1)^{p-1}}{ep^p}\vz^{\top}\nabla^2g_{\vx}(\vq)\vz + (p-1)e\vz^{\top}\nabla^2f(\vx)\vz\enspace,\\
        \Rightarrow &\frac{1}{(p-1)e}\vz^{\top}\nabla^2f(\vq)\vz - \frac{(p-1)^{p-2}}{e^2p^p}\vz^{\top}\nabla^2g_{\vx}(\vq)\vz \leq \vz^{\top}\nabla^2f(\vx)\vz\enspace.
    \end{align*}
    We can finally lower bound,
    \begin{align*}
        \vz^{\top}\nabla^2f_{\vq}(\vx)\vz &= \vz^{\top}\nabla^2f(\vx)\vz + \vz^{\top}\nabla^2g_{\vq}(\vx)\vz\enspace,\\
        &\geq \frac{1}{(p-1)e}\vz^{\top}\nabla^2f(\vq)\vz - \frac{(p-1)^{p-2}}{e^2p^p}\vz^{\top}\nabla^2g_{\vx}(\vq)\vz + \vz^{\top}\nabla^2g_{\vq}(\vx)\vz\enspace,\\
        &= \frac{1}{2(p-1)e}\vz^{\top}\left(\nabla^2h_{\vq}(\vx)- \nabla^2g_{\vq}(\vx)\right)\vz - \frac{(p-1)^{p-2}}{e^2p^p}\vz^{\top}\nabla^2g_{\vx}(\vq)\vz + \vz^{\top}\nabla^2g_{\vq}(\vx)\vz\enspace,\\
        &\geq^{(g_{\vq}(\vx) = g_{\vx}(\vq))} \frac{1}{2pe}\vz^{\top}\nabla^2h_{\vq}(\vx)\vz + \left(1 - \frac{1}{2(p-1)e}- \frac{(p-1)^{p-2}}{e^2p^p}\right)\vz^{\top}\nabla^2g_{\vq}(\vx)\vz\enspace,\\
        &\geq^{\text{(\Cref{lem:small_lemma_for_hessian_stability_two})}}\frac{1}{2pe}\vz^{\top}\nabla^2h_{\vq}(\vx)\vz\enspace, 
    \end{align*}
    where in the final inequality we use \Cref{lem:small_lemma_for_hessian_stability_two} and the fact that $\vz^{\top}\nabla^2g_{\vq}(\vx)\vz$ is non-negative. This finishes the proof for Case A.

    We finally consider the corner case with $p=2$. 

    \paragraph{\underline{Case B.} ($p=2$)} In this case the proof is trivial, and follows from simply writing the quadratic forms for $f_{\vq}$ and $h_{\vq}$. We do so below,
    \begin{align*}
        \vz^{\top}\nabla^2f_{\vq}(\vx)\vz &= \vz^{\top}\nabla^2f(\vx)\vz + \vz^{\top}\nabla^2g_{\vq}(\vx)\vz\enspace,\\
        &= \vz^{\top}\nabla^2f(\vx)\vz + 2C_2\norm{\vz}^2_{\mM}\enspace,\\
        &\leq 2\vz^{\top}\nabla^2f(\vx)\vz + 2C_2\norm{\vz}^2_{\mM} = \vz^{\top}\nabla^2h_{\vq}(\vx)\vz\enspace,
    \end{align*}
    which shows the relative smoothness with a constant of $1$ which is smaller (and hence better) than the claimed constant (for $p=2$) of $2e$ in the lemma. Now for the relative strong convexity we do the same,
    \begin{align*}
        \vz^{\top}\nabla^2f_{\vq}(\vx)\vz &= \vz^{\top}\nabla^2f(\vx)\vz + 2C_2\norm{\vz}^2_{\mM}\enspace,\\
        &\geq \frac{1}{2}\cdot\left(2\vz^{\top}\nabla^2f(\vx)\vz + 2C_2\norm{\vz}^2_{\mM}\right)\enspace,\\
        &= \frac{1}{2}\vz^{\top}\nabla^2h_{\vq}(\vx)\vz\enspace,
    \end{align*}
    which shows relative strong-convexity with a constant of $\frac{1}{2}$ which is larger (and hence better) than the claimed constant (for $p=2$) of $\frac{1}{4e}$ in the lemma. This finishes the proof for Case B. 
    
    This completes the proof of \Cref{lemma:gp_hessian_stable}.
\end{proof}
We prove two small technical lemmas that we used in the above proof now. 
\begin{lemma}\label{lem:small_lemma_for_hessian_stability_one}
    For all $p\geq 2$, $g(p) = 1+\frac{(p-1)^{p-1}}{ep^p} - \frac{(p-1)\cdot e}{2} \leq 0$.
\end{lemma}
\begin{proof}
     First note that at $p=2$ the function takes a strictly negative value,
     \begin{align*}
         g(2) = 1+\frac{(1}{e2^2} - \frac{e}{2} = \frac{4e+ 1 - 2e^2}{4e} < 0\enspace. 
     \end{align*}
     We will now show that the function is increasing in $p$ for $p\geq 2$,
     \begin{align*}
         g'(p) &= -\frac{(p-1)^{p-1}p^p(\ln(p) + 1)}{p^2p} + \frac{(p-1)^{p-1}(\ln(p-1) + 1)}{p^p} - \frac{e}{2}\enspace,\\
         &= -\frac{(p-1)^{p-1}\ln(p/(p-1))}{p^p} - \frac{e}{2} < 0\enspace. 
     \end{align*}
     Thus, the function attains its maximum value at $p=2$ in the range $p\geq 2$, implying it is strictly negative in that range.
\end{proof}

\begin{lemma}\label{lem:small_lemma_for_hessian_stability_two}
    For all $p\geq 2$, $g(p) = 1-\frac{1}{2(p-1)e} - \frac{(p-1)^{p-2}}{e^2p^p} \geq 0$.
\end{lemma}
\begin{proof}
    First note that at $p=2$ the function takes a strictly positive value,
    \begin{align*}
        g(2) = 1 - \frac{1}{2e} -\frac{1^0}{e^22^2} = 1 - \frac{1}{2e} - \frac{1}{4e^2} = \frac{4e^2 - 2e - 1}{4e^2} > 0\enspace.
    \end{align*}
    We will now show that the function is increasing in $p$ for $p\geq 2$,
    \begin{align*}
        g'(p) &= \frac{1}{2(p-1)^2e} + \frac{(p-1)^{p-2}p^p(\ln(p) + 1)}{e^2p^{2p}} - \frac{(p-1)^{p-2}(\ln(p-1) + (p-2)/(p-1))}{e^2p^p}\enspace,\\
        &= \frac{1}{2(p-1)^2e} + \frac{(p-1)^{p-2}(\ln(p) + 1)}{e^2p^{p}} - \frac{(p-1)^{p-2}(\ln(p-1) + 1 - 1/(p-1))}{e^2p^p}\enspace,\\
        &= \frac{1}{2(p-1)^2e} + \frac{(p-1)^{p-2}\left(\ln(p/(p-1)) + 1/(p-1)\right)}{e^2p^{p}} > 0\enspace.
    \end{align*}
    Thus, the function $g$ attains its minimum value at $p=2$ in the range $p\geq2$, implying that it is strictly positive in that range.
\end{proof}

\subsubsection{Strong convexity of the proximal objective and friends}

We begin with showing that the proximal objective enjoys a form of strong convexity.
\begin{lemma}
\label{lemma:fq_strongly_convex}
For all $\vx,\vd\in\R^d$, we have
\begin{align*}
    f_{\vq}(\vx+\vd) \ge f_{\vq}(\vx) + \ip{\nabla f_{\vq}(\vx), \vd} + \frac{4}{2^p}\inparen{\gnorm{\mA\vd}{p}^p + C_p\norm{\vd}_{\mM}^p}.
\end{align*}
\end{lemma}
\begin{proof}[Proof of \Cref{lemma:fq_strongly_convex}]
Let $K_p \coloneqq \frac{4}{2^p}$.

The plan is to apply \Cref{lemma:strong_convexity_component} to $f_{\vq}(\vx+\vd)$. We start with the regularizer. Notice that
\begin{align}
    \norm{\vx+\vd-\vq}_{\mM}^p &= \norm{\mM^{1/2}(\vx+\vd-\vq)}_2^p = \norm{\mM^{1/2}(\vx-\vq) + \mM^{1/2}\vd}_2^p\enspace,\nonumber \\
    &\ge^{\text{(\Cref{lemma:strong_convexity_component})}} \norm{\mM^{1/2}(\vx-\vq)}_2^p \\
    &\quad + \ip{p\norm{\mM^{1/2}(\vx-\vq)}_2^{p-2}\mM^{1/2}(\vx-\vq),\mM^{1/2}\vd} + K_p\norm{\mM^{1/2}\vd}_2^p\enspace,\nonumber\\
    &= \norm{\vx-\vq}_{\mM}^p + \ip{p\norm{\vx-\vq}_{\mM}^{p-2}\mM(\vx-\vq), \vd} + K_p\norm{\vd}_{\mM}^p\enspace,\nonumber\\
    &= \norm{\vx-\vq}_{\mM}^p + \ip{\nabla_{\vx} \inparen{\norm{\vx-\vq}_{\mM}^p}, \vd} + K_p\norm{\vd}_{\mM}^p\enspace.\label{eq:M_norm_sc}
\end{align}
We combine this with the conclusion of \Cref{lemma:strong_convexity_gp}, giving
\begin{align*}
    f_{\vq}(\vx+\vd) &= f(\vx+\vd) + C_p\|\vx+\vd - \vq\|_{\mM}^p\enspace,\\
    &\geq^{\text{(\Cref{lemma:strong_convexity_gp})}} f(\vx) + \ip{\nabla f(\vx),\vd} + K_p\gnorm{\mA\vd}{p}^p+ C_p\|\vx+\vd - \vq\|_{\mM}^p\enspace,\\
    &\geq^{\text{\eqref{eq:M_norm_sc}}} f(\vx) + \ip{\nabla f(\vx),\vd} + K_p\gnorm{\mA\vd}{p}^p + C_p\norm{\vx-\vq}_{\mM}^p\\ 
    &\quad + C_p\ip{\nabla_{\vx} \inparen{\norm{\vx-\vq}_{\mM}^p}, \vd} + K_pC_p\norm{\vd}_{\mM}^p\enspace,\\
    &= \textcolor{red}{f(\vx)}  + \textcolor{red}{C_p\norm{\vx-\vq}_{\mM}^p} + \ip{\textcolor{blue}{\nabla_{\vx} \left(f(\vx) + C_p\norm{\vx-\vq}_{\mM}^p\right)},\vd}\\
    &\quad + K_p\gnorm{\mA\vd}{p}^p + K_pC_p\norm{\vd}_{\mM}^p\enspace,\\
    &= \textcolor{red}{f_{\vq}(\vx)} + \ip{\textcolor{blue}{\nabla f_{\vq}(\vx)},\vd} + K_p\inparen{\gnorm{\mA\vd}{p}^p + C_p\norm{\vd}_{\mM}^p}\enspace.
\end{align*}
completing the proof of \Cref{lemma:fq_strongly_convex}.
\end{proof}

We also show that the subproblems we solve in Line \ref{line:mirror_descent_exact} of \Cref{alg:mirror_descent} are strongly convex.
\begin{lemma}
\label{lemma:gp_h_strongly_convex}
Fix $\vz, \vq, \vd \in \R^d$ and let $L > 0$. Consider the function
\begin{align*}
    g(\vx) &\coloneqq \ip{\vz,\vx} + L\inparen{\norm{\vx-\vq}_{\nabla^2 f(\vq)}^2 + C_p\norm{\vx-\vq}_{\mM}^p}\enspace.
\end{align*}
Then,
\begin{align*}
    g(\vx+\vd) \ge g(\vx) + \ip{\nabla g(\vx),\vd} + L\inparen{\norm{\vd}_{\nabla^2 f(\vq)}^2 + \frac{4C_p}{2^p}\norm{\vd}_{\mM}^p}\enspace.
\end{align*}
In particular, if $\vz$ is the minimizer for $g$, then for any $\vd\in\R^d$, we have
\begin{align*}
    \norm{\vd}_{\mM} \le \frac{2}{p \cdot (4e)^{1/p}}\inparen{\frac{g(\vz+\vd) - g(\vz)}{L}}^{1/p}\enspace.
\end{align*}
\end{lemma}
\begin{proof}[Proof of \Cref{lemma:gp_h_strongly_convex}]
This is pretty much the same proof as \Cref{lemma:fq_strongly_convex}. It is easy to check that
\begin{align}
    \norm{(\vx+\vd)-\vq}_{\nabla^2 f(\vq)}^2 = \norm{\vx-\vq}_{\nabla^2 f(\vq)}^2 + \ip{2\nabla^2 f(\vq)(\vx-\vq),\vd} + \norm{\vd}_{\nabla^2 f(\vq)}^2\enspace,\label{eq:f_norm_expansion}
\end{align}
and using \Cref{lemma:strong_convexity_component} in the same way as in the proof of \Cref{lemma:fq_strongly_convex}, we have
\begin{align*}
    \norm{(\vx+\vd)-\vq}_{\mM}^p \ge^{\text{\eqref{eq:M_norm_sc}}} \norm{\vx-\vq}_{\mM}^p + \ip{p\norm{\vx-\vq}_{\mM}^{p-2}\mM(\vx-\vq),\vd} + \frac{4}{2^p}\norm{\vd}_{\mM}^p\enspace.
\end{align*}
Combining this with the definition of $g$ gives the following,
\begin{align*}
    g(\vx+\vd) &= \ip{\vz,\vx + \vd} + L\inparen{\norm{\vx + \vd-\vq}_{\nabla^2 f(\vq)}^2 + C_p\norm{\vx + \vd -\vq}_{\mM}^p}\enspace,\\
    &\geq^{\text{\eqref{eq:f_norm_expansion}, \eqref{eq:M_norm_sc}}} \textcolor{red}{\ip{\vz,\vx}} +  \ip{\vz,\vd} + \textcolor{red}{L\norm{\vx-\vq}_{\nabla^2 f(\vq)}^2} + L\ip{2\nabla^2 f(\vq)(\vx-\vq),\vd}\\
    &\qquad + L\norm{\vd}_{\nabla^2 f(\vq)}^2 + LC_p\left(\textcolor{red}{\norm{\vx-\vq}_{\mM}^p} + \ip{p\norm{\vx-\vq}_{\mM}^{p-2}\mM(\vx-\vq),\vd} + \frac{4}{2^p}\norm{\vd}_{\mM}^p\right)\enspace,\\
    &=  \textcolor{red}{g(\vx)}  + \ip{\textcolor{blue}{\vz} + \textcolor{blue}{2L\nabla^2 f(\vq)(\vx-\vq)} + \textcolor{blue}{LC_pp\norm{\vx-\vq}_{\mM}^{p-2}\mM(\vx-\vq)},\vd} \\
    &\qquad  + L\left(\norm{\vd}_{\nabla^2 f(\vq)}^2 + \frac{4C_p}{2^p}\norm{\vd}_{\mM}^p\right)\enspace,\\
    &= g(\vx) + \ip{\textcolor{blue}{\nabla g(\vx)},\vd} + L\left(\norm{\vd}_{\nabla^2 f(\vq)}^2 + \frac{4C_p}{2^p}\norm{\vd}_{\mM}^p\right)\enspace,
\end{align*}
which proves the first result of the lemma.

To get the second result, we observe that $\nabla g(\vz) = 0$ by the optimality of $\vz$. Ignoring the $\norm{\vd}_{\nabla^2 f(\vq)}$ terms and rearranging gives the conclusion of \Cref{lemma:gp_h_strongly_convex}.
\end{proof}

\subsubsection{Smoothness of the proximal objective}

We first bound the operator norm of a matrix related to the Hessian of the proximal objective.

\begin{lemma}
\label{lemma:gp_prox_hessian_opnorm}
For all $\vq,\vy\in\R^d$, we have
\begin{align*}
    \opnorm{\mM^{-1/2}\inparen{\nabla^2 f_{\vq}(\vy)}\mM^{-1/2}} &\le ep^2(p-1)\inparen{2f(\vq)^{1-\frac{2}{p}}+C_p\norm{\vy-\vq}_{\mM}^{p-2}}\enspace.
\end{align*}
\end{lemma}
\begin{proof}[Proof of \Cref{lemma:gp_prox_hessian_opnorm}]
Recall from the proof of \Cref{lemma:gp_hessian_stable} the definition of the regularization term $g_{\vq}(\vy) \coloneqq C_p\norm{\vy-\vq}_{\mM}^p$ for $C_p = ep^p$ as well as the following calculations,
 \begin{align*}
        g_{\vq}(\vy) &\coloneqq C_p\norm{\vy-\vq}^p_{\mM}\enspace,\\    
        \nabla g_{\vq}(\vy) &= pC_p\norm{\vy-\vq}_{\mM}^{p-2}{\mM}(\vy-\vq)\enspace,\\
        \nabla^2 g_{\vq}(\vy) &= pC_p\norm{\vy-\vq}_{\mM}^{p-2}\mM + p(p-2)C_p\norm{\vy-\vq}_{\mM}^{p-4}\mM(\vy-\vq)(\vy-\vq)^{\top}\mM\enspace.
\end{align*}
By \Cref{lemma:gp_hessian_stable}, we know that
\begin{align*}
    \nabla^2 f_{\vq}(\vy) \preceq ep\inparen{2\nabla^2 f(\vq) + \nabla^2 g_{\vq}(\vy)}.
\end{align*}
Observe that
\begin{align*}
    \quad \mM^{-1/2} \inparen{\nabla^2 g_{\vq}(\vy)} \mM^{-1/2} &= pC_p\inparen{\norm{\vy-\vq}_{\mM}^{p-2}+(p-2)\norm{\vy-\vq}_{\mM}^{p-4}\mM^{1/2}(\vy-\vq)(\vy-\vq)^{\top}\mM^{1/2}}\enspace,\\
    &\preceq pC_p\norm{\vy-\vq}_{\mM}^{p-2}\mI + (p-2)\norm{\vy-\vq}_{\mM}^{p-4}\opnorm{\mM^{1/2}(\vy-\vq)(\vy-\vq)^{\top}\mM^{1/2}}\mI\enspace,\\
    &\preceq pC_p\norm{\vy-\vq}_{\mM}^{p-2}\mI + (p-2)\norm{\vy-\vq}_{\mM}^{p-4}\norm{\mM^{1/2}(\vy-\vq)}_2^2\mI\enspace,\\ 
    &\preceq p(p-1)C_p \norm{\vy-\vq}_{\mM}^{p-2}\mI\enspace,
\end{align*}
and, applying \Cref{lemma:gp_regression_hessian_upperbound} (with $\mM^{-1/2}\vz$ as the vectors in the quadratic form) and H\"older inequality with norms $\|\cdot\|_{p/(p-2)},\ \|\cdot\|_{p/2}$, for $\vz\in\R^d$ we have\todo{add a reference for the third inequality. I keep forgetting which result this is in the previous sections.}
\begin{align*}
    \vz^{\top}\mM^{-1/2}\inparen{\nabla^2 f(\vq)} \mM^{-1/2}\vz &\le p(p-1)\sum_{i=1}^m \norm{\mA_{S_i}\vq-\vb_{S_i}}_2^{p-2}\norm{\mA_{S_i}\mM^{-1/2}\vz}_2^2 \\
    &\le p(p-1)\inparen{\sum_{i=1}^m \norm{\mA_{S_i}\vq-\vb_{S_i}}_2^p}^{\frac{p-2}{p}}\inparen{\sum_{i=1}^m \norm{\mA_{S_i}\mM^{-1/2}\vz}_2^p}^{\frac{2}{p}} \\
    &\le p(p-1)f(\vq)^{1-\frac{2}{p}}\norm{\mM^{-1/2}\vz}_{\mM}^2 = p(p-1)f(\vq)^{1-\frac{2}{p}}\norm{\vz}_2^2.
\end{align*}
Combining gives
\begin{align*}
    \mM^{-1/2}\inparen{\nabla^2f_{\vq}(\vy)}\mM^{-1/2} &\preceq ep\mM^{-1/2}\inparen{2\nabla^2f(\vq) + \nabla^2g_{\vq(\vy)}}\mM^{-1/2}\enspace,\\    
    &\preceq 2ep^2(p-1)f(\vq)^{1-\frac{2}{p}} + ep^2(p-1)C_p\norm{\vy-\vq}_{\mM}^{p-2}\enspace, \\
    &\preceq ep^2(p-1)\inparen{2f(\vq)^{1-\frac{2}{p}}+C_p\norm{\vy-\vq}_{\mM}^{p-2}},
\end{align*}
completing the proof of \Cref{lemma:gp_prox_hessian_opnorm}.
\end{proof}

Next, we show a bound on the norm of the gradient of any solution $\vx$ that is approximately optimal for $f_{\vq}$.

\begin{lemma}
\label{lemma:self_gradient_norm_small}
For all $\vq,\vx\in\R^d$, we have
\begin{align*}
    \norm{\mM^{-1}\nabla f_{\vq}(\vx)}_{\mM} &\le ep^2(p-1)\inparen{f(\vq)^{1-\frac{2}{p}}+C_p\max\inbraces{\norm{\vx-\vq}_{\mM},\norm{\vx_{\vq}-\vq}_{\mM}}^{p-2}}\norm{\vx-\vx_{\vq}}_{\mM}\enspace.
\end{align*}
\end{lemma}
\begin{proof}[Proof of \Cref{lemma:self_gradient_norm_small}]
We use a continuity argument. By Taylor's theorem, we know for some $\vy$ along the line connecting $\vx$ and $\vx_{\vq}$ (minimizer of $f_{\vq}$) that
\begin{align*}
    \nabla f_{\vq}(\vx) = \nabla f_{\vq}(\vx_{\vq}) + \nabla^2 f_{\vq}(\vy)(\vx-\vx_{\vq}) = \nabla^2 f_{\vq}(\vy)(\vx-\vx_{\vq})\enspace.
\end{align*}
Taking $\mM^{-1}$-norm of both sides gives,
\begin{align*}
    \norm{\nabla f_{\vq}(\vx)}_{\mM^{-1}} &= \norm{\mM^{-1/2}\nabla f_{\vq}(\vx)}_{2}\enspace,\\
     &= \norm{\mM^{-1/2}\nabla^2 f_{\vq}(\vy)(\vx-\vx_{\vq})}_{2}\enspace,\\
     &= \norm{\mM^{-1/2}\nabla^2 f_{\vq}(\vy)\mM^{-1/2}\mM^{1/2}(\vx-\vx_{\vq})}_{2}\enspace,\\
    &\le \opnorm{\mM^{-1/2}\inparen{\nabla^2 f_{\vq}(\vy)}\mM^{-1/2}} \cdot \norm{\vx-\vx_{\vq}}_{\mM}\enspace.
\end{align*}
The rest of the proof involves bounding the operator norm term. This follows directly from \Cref{lemma:gp_prox_hessian_opnorm}, from which we get (using convexity of $\|\cdot\|_{\mM}$),
\begin{align*}
    \opnorm{\mM^{-1/2}\nabla^2 f_{\vq}(\vy)\mM^{-1/2}} &\le ep^2(p-1)\inparen{2f(\vq)^{1-\frac{2}{p}}+C_p\norm{\vy-\vq}_{\mM}^{p-2}} \\
    &\le ep^2(p-1)\inparen{2f(\vq)^{1-\frac{2}{p}}+C_p\max\inbraces{\norm{\vx-\vq}_{\mM},\norm{\vx_{\vq}-\vq}_{\mM}}^{p-2}}.
\end{align*}
Putting everything together, we get
\begin{align*}
    \norm{\mM^{-1}\nabla f_{\vq}(\vx)}_{\mM}&= \norm{\nabla f_{\vq}(\vx)}_{\mM^{-1}}\enspace, \\
    &\le ep^2(p-1)\inparen{2f(\vq)^{1-\frac{2}{p}}+C_p\max\inbraces{\norm{\vx-\vq}_{\mM},\norm{\vx_{\vq}-\vq}_{\mM}}^{p-2}}\norm{\vx-\vx_{\vq}}_{\mM},
\end{align*}
completing the proof of \Cref{lemma:self_gradient_norm_small}.
\end{proof}

\subsubsection{Solving the proximal subproblems}

We begin by showing that the optimal solution to the proximal problem $\vx_{\vq_t}\coloneqq \argmin{\vx\in\R^d} f_{\vq_t}(\vx)$ is not too far from $\xstar$.

\begin{lemma}
\label{lemma:gp_regression_prox_diameter}
For all proximal queries $\vq_t$, we have
\begin{align*}
    \norm{\vx_{\vq_t}-\xstar}_{\mM} \le d^{\frac{1}{2}-\frac{1}{p}}\inparen{2^{\frac{3}{2}}f(\vx_t)+4}.
\end{align*}
\end{lemma}
\begin{proof}
In the rest of this proof, we omit the subscript $t$ wherever it is clear which iterates we are working with. 

We first show that
\begin{align*}
    \norm{\vx_{\vq}-\vq}_{\mM} \le \norm{\xstar-\vq}_{\mM}\enspace.
\end{align*}
To see this, suppose this is not the case. Then, we have
\begin{align*}
    f(\xstar) + C_p\norm{\xstar-\vq}_{\mM}^p < f(\vx_{\vq}) + C_p\norm{\vx_{\vq}-\vq}_{\mM}^p\enspace,
\end{align*}
which contradicts the optimality of $\vx_{\vq}$ for $f_{\vq}$.

We now write
\begin{align*}
    \norm{\vx_{\vq_t}-\xstar}_{\mM} &\le \norm{\vx_{\vq_t}-\vq_t}_{\mM} + \norm{\xstar-\vq_t}_{\mM}\enspace, \\
    &\le 2\norm{\xstar-\vq_t}_{\mM}\enspace,\\
    &\le 2\inparen{\norm{\vx_t-\xstar}_{\mM}+\norm{\vv_t-\xstar}_{\mM}}\enspace,
\end{align*}
where in the last inequality, we used the definition of $\vq_t$ from Line 6 in \Cref{alg:optimal_ms_acceleration} and the convexity of $\|\cdot\|_{\mM}$. The required control on $\norm{\vv_t-\xstar}_{\mM}$ comes from \Cref{lemma:ms_acceleration_iterate_diameter} and \Cref{thm:gp_regression_blw_intro} (along with re-scaling assumption to make the optimal value $1$) -- we have
\begin{align*}
    \norm{\vv_t-\xstar}_{\mM} \le \sqrt{2}\norm{\vx_0-\xstar}_{\mM} \le 4d^{\frac{1}{2}-\frac{1}{p}}\enspace.
\end{align*}
For the other term, we apply \Cref{lemma:strong_convexity_gp} and get
\begin{align*}
    \norm{\vx_t-\xstar}_{\mM} \le 2^{\frac{3}{2}}d^{\frac{1}{2}-\frac{1}{p}}\inparen{f(\vx_t)-f(\xstar)}^{\frac{1}{p}} < 2^{\frac{3}{2}}d^{\frac{1}{2}-\frac{1}{p}}f(\vx_t)^{\frac{1}{p}}.
\end{align*}
Adding gives us the conclusion of \Cref{lemma:gp_regression_prox_diameter}.
\end{proof}

The next few lemmas are targeted at solving the proximal subproblems. We begin with a calculation that we will use in showing that the initial Bregman divergence between our initialization and the optimum is small.

\begin{lemma}
\label{lemma:initial_bregman_div}
In the same setting as \Cref{lemma:gp_hessian_stable}, for all $\vx,\vy \in \R^d$, we have
\begin{align*}
    h_{\vq}(\vx_{\vq}) \le p(p-1)f(\vq)^{1-\frac{2}{p}}\norm{\vx_{\vq}-\vq}_{\mM}^2 + C_p\norm{\vx_{\vq}-\vq}_{\mM}^p < f(\vq)+C_p\norm{\vx_{\vq}-\vq}_{\mM}^p \le 2f(\vq).
\end{align*}
\end{lemma}
\begin{proof}[Proof of \Cref{lemma:initial_bregman_div}]
By optimality of $\vx_{\vq}$ for the subproblem, we have
\begin{align*}
    f(\vx_{\vq})+C_p\norm{\vx_{\vq}-\vq}_{\mM}^p \le f(\vq)+C_p\norm{\vq-\vq}_{\mM}^p=f(\vq).
\end{align*}
Rearranging gives,
\begin{align}
    \norm{\vx_{\vq}-\vq}_{\mM}^p \le \frac{f(\vq)-f(\vx_{\vq})}{C_p} \le \frac{f(\vq)}{C_p}\enspace.\label{eq:norm_M_ub_f}
\end{align}
We now use the definition of $h_{\vq}$ and \Cref{lemma:gp_regression_hessian_upperbound} to write
\begin{align*}
    h_{\vq}(\vx_{\vq}) &= \norm{\vx_{\vq}-\vq}_{\nabla^2 f(\vq)}^2 + C_p\norm{\vx_{\vq}-\vq}_{\mM}^p\enspace,\\
    &\le^{\text{\Cref{lemma:gp_regression_hessian_upperbound}}} p(p-1)\sum_{i=1}^m \norm{\mA_{S_i}\vq-\vb_{S_i}}_2^{p-2}\norm{\mA_{S_i}(\vx_{\vq}-\vq)}_2^2 + C_p\norm{\vx_{\vq}-\vq}_{\mM}^p\enspace,\\
    &\le^{\text{(a)}} p(p-1)\inparen{\sum_{i=1}^m \norm{\mA_{S_i}\vq-\vb_{S_i}}_2^p}^{1-\frac{2}{p}}\inparen{\sum_{i=1}^m \norm{\mA_{S_i}(\vx_{\vq}-\vq)}_2^p}^{\frac{2}{p}}+ C_p\norm{\vx_{\vq}-\vq}_{\mM}^p\enspace,\\
    &\le^{\text{(b)}} p(p-1)f(\vq)^{1-\frac{2}{p}}\norm{\vx_{\vq}-\vq}_{\mM}^{2} + C_p\norm{\vx_{\vq}-\vq}_{\mM}^p\enspace,\\
    &\le^{\text{\eqref{eq:norm_M_ub_f}}} p(p-1)f(\vq)^{1-\frac{2}{p}}\inparen{\frac{f(\vq)}{C_p}}^{\frac{2}{p}} + C_p\norm{\vx_{\vq}-\vq}_{\mM}^p\enspace,\\
    &=^{\text{($C_p = ep^p$)}} \frac{(p-1)}{ep}f(\vq)+ C_p\norm{\vx_{\vq}-\vq}_{\mM}^p\enspace,\\
    &< f(\vq) + C_p\norm{\vx_{\vq}-\vq}_{\mM}^p\enspace,\\
    &<^{\text{\eqref{eq:norm_M_ub_f}}} 2f(\vq)\enspace,
\end{align*}
where in (a) we used H\"older inequality with norms $\|\cdot\|_{p/(p-2)},\ \|\cdot\|_{p/2}$ and in (b) we used \Cref{thm:gp_regression_blw_intro} \todo{again what's the reference for this one?}. 

This completes the proof for the series of inequalities in \Cref{lemma:initial_bregman_div}.
\end{proof}
We now have the tools to show how to approximately solve problems in Line \ref{line:mirror_descent_exact} of \Cref{alg:mirror_descent} when applied in our setting. Although this and future complexity bounds depend on $f(\vx_t)$, we will later be able to use \Cref{thm:optimal_ms_acceleration} to ``bootstrap'' and get an unconditional upper bound below.
\begin{lemma}
\label{lemma:gp_prox_relativesmoothsolve}
Let $\alpha \le 1/2$. In the context of \Cref{alg:gp_regression_alg}, there exists an algorithm that approximately solves subproblems of the form (for $p\geq 2$ and $L=pe$),
\begin{align*}
    \vz \coloneqq \argmin{\vx\in\R^d} \ip{\vg,\vx} + L\inparen{\norm{\vx-\vq}_{\nabla^2 f(\vq)}^2 + C_p\norm{\vx-\vq}_{\mM}^p}\enspace,
\end{align*}
in the sense that we output $\vx$ for which,
\begin{equation*}
    \begin{aligned}
        \max\left\{\norm{\vx-\vz}_{\mM}, \norm{\mM^{-1}\vg + 2L\inparen{\mM^{-1}\nabla^2 f(\vq)(\vx-\vq) + C_p\norm{\vx-\vq}_{\mM}^{p-2}(\vx-\vq)}}_{\mM}\right\} \le \alpha
    \end{aligned}\enspace.
\end{equation*}
The algorithm takes $p^{O(1)}\logv{\frac{pd \cdot f(\vq)}{\alpha}}$ linear system solves in matrices of the form $\mA^{\top}\mB\mA$ for block-diagonal $\mB$, where each block in $\mB$ has size $\abs{S_i} \times \abs{S_i}$.
\end{lemma}
\begin{proof}[Proof of \Cref{lemma:gp_prox_relativesmoothsolve}]
This proof is long, and splitting it into lemmas would break up the intended reading flow. So we break it up into several key components here.

\smallpar{Motivation for the lemma.} First, let us see why this lemma is even useful. In each iteration of \Cref{alg:gp_prox_solver}, which in turn calls \Cref{alg:mirror_descent}, the main primitive is computing\todo{I think this part should appear before the lemma because it leads to confusion within the proof about what sub-problem we are solving. If it is motivation, we don't need a separate lemma we can just put it before the lemma.}
\begin{align*}
    \xtilde_i &= \argmin{\xtilde \in \R^d} f_{\vq_t}(\xtilde_{i-1}) + \ip{\nabla f_{\vq_t}(\xtilde_{i-1}), \xtilde-\xtilde_{i-1}} + peD_{h_{\vq_t}}(\xtilde,\xtilde_{i-1})\enspace, \\
    &= \argmin{\xtilde \in \R^d} f_{\vq_t}(\xtilde_{i-1}) + \ip{\nabla f_{\vq_t}(\xtilde_{i-1}), \xtilde-\xtilde_{i-1}} + pe \left(h_{\vq_t}(\xtilde) - h_{\vq_t}(\xtilde_{i-1}) - \ip{\nabla h_{\vq_t}(\xtilde_{i-1}), \xtilde - \xtilde_{i-1}}\right)\enspace, \\
    &= \argmin{\xtilde \in \R^d} f_{\vq_t}(\xtilde_{i-1})-peh_{\vq_t}(\xtilde_{i-1}) + \ip{\nabla f_{\vq_t}(
    \xtilde_{i-1})-pe\nabla h_{\vq_t}(\xtilde_{i-1}), \xtilde-\xtilde_{i-1}} + peh_{\vq_t}(\xtilde)\enspace, \\
    &= \argmin{\xtilde\in\R^d} \ip{\nabla f_{\vq_t}(
    \xtilde_{i-1})-pe\nabla h_{\vq_t}(\xtilde_{i-1}), \xtilde} + peh_{\vq_t}(\xtilde)\enspace.
\end{align*}
Observe that the subproblem is of the form
\begin{align}
    \vz &= \argmin{\vx\in\R^d} \ip{\vg, \vx} + peh_{\vq}(\vx)\enspace,\nonumber\\
    &= \argmin{\vx\in\R^d} \ip{\vg,\vx}+pe\inparen{\norm{\vx-\vq}_{\nabla^2 f(\vq)}^2 + C_p\norm{\vx-\vq}_{\mM}^p}\enspace,\label{eq:equiv_prox_problem}
\end{align}
and so our goal is to show how to solve these types of problems.

\smallpar{The general algorithm.} 
Consider solving the related subproblem (instead of \eqref{eq:equiv_prox_problem}),
\begin{align*}
    \argmin{\vx\in\R^d} \ip{\vg,\vx}+L\inparen{\norm{\vx-\vq}_{\nabla^2 f(\vq)}^2 + C_p\tau\norm{\vx-\vq}_{\mM}^2}
\end{align*}
for some fixed $\tau \ge 0$. This is a quadratic problem, and we can therefore solve it in $1$ linear system solve. It is easy to check that at optimality, we have
\begin{align*}
    \vg + 2pe\inparen{\nabla^2 f(\vq)(\vx-\vq) + C_p\tau\mM(\vx-\vq)} = 0\enspace,
\end{align*}
which rearranges to\footnote{Recall that $\nabla^2 f(\vq)=\mA^{\top}\mB_1\mA$ for block-diagonal $\mB_1$ and by construction, $\mM = \mA^{\top}\mW^{1-\frac{2}{p}}\mA$ where $\mW$ consists of the block Lewis weights on the diagonal. Thus, $\nabla^2 f(\vq) + C_p\tau\mM=\mA^{\top}\mB_2\mA$ for block-diagonal $\mB_2$.}
\begin{align*}
    \vx-\vq = -\frac{1}{2pe}\inparen{\nabla^2 f(\vq)+C_p\tau\mM}^{-1}\vg\enspace.
\end{align*}

Note that at optimality for our original subproblem \eqref{eq:equiv_prox_problem}, we have $\tau^\star := \norm{\vz-\vq}_{\mM}^{p-2}$ where $\vz$ is the solution of subproblem \eqref{eq:equiv_prox_problem}. Also note that $\norm{\vx-\vq}_{\mM}$ is a decreasing function in $\tau$ because,
\begin{align*}
    \norm{\vx-\vq}_{\mM}^2 = \frac{1}{4p^2e^2}\|\vg\|^2_{\inparen{\nabla^2 f(\vq)+C_p\tau\mM}^{-1}\mM\inparen{\nabla^2 f(\vq)+C_p\tau\mM}^{-1}}\enspace, 
\end{align*}
and for $\tau_1 \leq \tau_2$,
\begin{align*}
    \inparen{\nabla^2 f(\vq)+C_p\tau_1\mM}^{-1}\mM\inparen{\nabla^2 f(\vq)+C_p\tau_1\mM}^{-1} \succeq \inparen{\nabla^2 f(\vq)+C_p\tau_2\mM}^{-1}\mM\inparen{\nabla^2 f(\vq)+C_p\tau_2\mM}^{-1}\enspace.
\end{align*}
We therefore see that if $\tau > \norm{\vx-\vq}_{\mM}^{p-2}$ --- where $\vx$ is the optimal solution for a fixed $\tau$ --- then we are over-regularizing and need to decrease $\tau$ and vice-versa. 
This means we can binary search for the appropriate value of $\tau$. To execute this, we first need to establish the accuracy up to which we have to identify $\tau$. 

\smallpar{Convergence in Argument.} 
\todo{I can't follow this. I tried writing the function values, but it is certainly not obvious from looking at the expression; we need more steps and words here.}
By \Cref{lemma:gp_h_strongly_convex} (setting $\vd = \vx-\vz$), recall that it is enough to solve sub-problem \eqref{eq:equiv_prox_problem} up to additive accuracy $(p/2)^pL\alpha^p$ to get $\norm{\vx-\vz}_{\mM} \le \alpha$. Suppose we find $\tau$ for which $\tau^{\star} \le \tau \le \tau^{\star} + \delta$. By writing the objectives and comparing, we see that the $\vx$ we find from using $\tau$ gives us at most a $\delta \cdot d$-suboptimal solution compared to $\vz$. Plugging this into the bound from \Cref{lemma:gp_h_strongly_convex} tells us that we should choose $\delta = (p/2)^pL\alpha^p/d$, and plugging this into the binary search over $\tau \in [0,d^p(1+f(\vq))]$ gives us $p^{O(1)}\logv{\frac{pd \cdot f(\vq)}{\alpha}}$ steps, as needed.\todo{There is a typo here, but I think we discussed we anyways want to make the statement of the lemma iterate dependent; should discuss.}

\smallpar{First-order stationary point.} We first claim that it is enough to get
\begin{align*}
    \norm{\mM^{-1}\nabla h_{\vq}(\vx)-\mM^{-1}\nabla h_{\vq}(\vz)}_{\mM} \le \frac{\alpha}{L}.
\end{align*}
Indeed, let $\vz$ be the optimal solution for the subproblem. This means that it must satisfy the first order stationary condition, namely,
\begin{align*}
    \vg + L\nabla h_{\vq}(\vz) = 0.
\end{align*}
Multiplying both sides by $\mM^{-1}$, subtracting, and dividing both sides by $L$ gives us the expression we are interested in.

Writing first order stationary conditions gives both
\begin{equation*}
    \begin{aligned}
        \vg + 2L\inparen{\nabla^2 f(\vq)(\vx-\vq)+C_p\tau\mM(\vx-\vq)} &= 0 \\
        \vg + 2L\inparen{\nabla^2 f(\vq)(\vz-\vq)+C_p\tau^{\star}\mM(\vz-\vq)} &= 0
    \end{aligned}.
\end{equation*}
Multiplying both sides of both equalities by $\mM^{-1}$ and subtracting these gives
\begin{align*}
    2L\inparen{\mM^{-1}\nabla^2 f(\vq)(\vx-\vz) + C_p\inparen{\tau(\vx-\vq)-\tau^{\star}(\vz-\vq)}} = 0.
\end{align*}
Expanding out $L(\mM^{-1}\nabla h_{\vq}(\vx)-\mM^{-1}h_{\vq}(\vz))$ and subtracting the above gives the desired condition
\begin{align*}
    2L\abs{\tau-\norm{\vx-\vq}_{\mM}^{p-2}} \cdot \norm{\vx-\vq}_{\mM} \overset{?}{\le} \alpha.
\end{align*}
Next, let us run the binary search from above so that we get argument convergence, i.e. $\norm{\vx-\vz}_{\mM} \le \alpha^C \ll 0.1\alpha$ for some constant $C$. Using the fact that the approximate mirror descent step using $\vz$ decreases the objective value (\Cref{lemma:gp_md_descent}), observe that
\begin{align*}
    \norm{\vx-\vq}_{\mM} \le \norm{\vz-\vq}_{\mM} + \norm{\vx-\vz}_{\mM} \le \norm{\vq-\vz}_{\mM}+0.1\alpha \lesssim \sqrt{d}(1+f(\vq)). 
\end{align*}
It then follows that binary searching $\tau$ to additive accuracy $\alpha(\sqrt{d}(1+f(\vq)))^{-1}/L$ is sufficient. By the same argument as above, this takes $p^{O(1)}\logv{\frac{pd \cdot f(\vq_t)}{\alpha}}$ steps, completing the proof of \Cref{lemma:gp_prox_relativesmoothsolve}.
\end{proof}
We now combine \Cref{lemma:gp_prox_relativesmoothsolve} with \Cref{thm:gp_mirror_descent} and \Cref{alg:mirror_descent} to obtain approximate argument optimality for each proximal subproblem.
\begin{lemma}
\label{lemma:prox_subproblem_base}
Let $\gamma > 0$ and $\vx_{\vq} \coloneqq \argmin{\vx\in\R^d} f_{\vq}(\vx)$. There exists an algorithm that returns $\vx$ for which
\begin{align*}
    \norm{\vx-\vx_{\vq}}_{\mM} &\le \gamma.
\end{align*}
The algorithm takes at most $O\inparen{p^{O(1)}\logv{ph_{\vq}(\vx_{\vq})\inparen{\frac{4}{p\gamma}}^{p}}}$ iterations of solving subproblems of the form $\argmin{\vx\in\R^d} \ip{\vg,\vx}+eph_{\vq}(\vx)$ for fixed vectors $\vg$ and $\vq$.
\end{lemma}
\begin{proof}[Proof of \Cref{lemma:prox_subproblem_base}]
This proof resembles \cite[Lemma 4.5]{jls21}, which uses an exact version of mirror descent arising from \citet{lfn18}. The main difference between our argument and that of \cite[Lemma 4.5]{jls21} is that we rigorously identify a concrete upper bound on the complexity needed to satisfy the MS condition and argue that the mirror descent algorithm can handle the inexact Bregman proximal problem solves.

First, we use \Cref{lemma:fq_strongly_convex} on the approximate solution $\vx$ and true solution $\vx_{\vq}$ and get,
\begin{align*}
    f_{\vq}(\vx) &\ge f_{\vq}(\vx_{\vq}) + \frac{4}{2^p}\inparen{\gnorm{\mA(\vx-\vx_{\vq})}{p}^p+C_p\norm{\vx_{\vq}-\vx}_{\mM}^p}\enspace,\\
    &\geq f_{\vq}(\vx) + \frac{4C_p}{2^p}\norm{\vx_{\vq}-\vx}_{\mM}^p\enspace.
\end{align*}
Rearranging, we get
\begin{align*}
    \norm{\vx_{\vq}-\vx}_{\mM} &\le \inparen{\frac{2^p}{4C_p}}^{1/p}\inparen{f_{\vq}(\vx)-f_{\vq}(\vx_{\vq})}^{1/p}\enspace,\\
    &= \inparen{\frac{2^p}{4ep^p}}^{1/p}\inparen{f_{\vq}(\vx)-f_{\vq}(\vx_{\vq})}^{1/p}\enspace,\\
    &< \frac{2}{p}\inparen{f_{\vq}(\vx)-f_{\vq}(\vx_{\vq})}^{1/p}\enspace.
\end{align*}
Using the notation from \cite{lfn18}, for convex $\myfunc{h}{\R^d}{\R}$, let
\begin{align*}
    D_h(\vx,\vy) \coloneqq h(\vx) - h(\vy) - \ip{\nabla h(\vy), \vx-\vy}.
\end{align*}
Recall the conclusion of \Cref{lemma:gp_hessian_stable} -- we have for $\mu = 1/(2pe)$ and $L = pe$ that
\begin{align*}
    \mu\nabla^2 h_{\vq}(\vx) \preceq \nabla^2 f_{\vq}(\vx) \preceq L\nabla^2 h_{\vq}(\vx).
\end{align*}
By \Cref{thm:gp_mirror_descent} and \Cref{lemma:gp_hessian_stable}, using the same notation from \Cref{lemma:gp_hessian_stable}, we have for all iterations $t$ of \Cref{alg:mirror_descent} (with $f = f_{\vq}$ and $h = h_{\vq}$) that,\todo{not sure about the second equality...}
\begin{align*}
    f_{\vq}(\vx_t) - f_{\vq}(\vx_{\vq}) &\le L\inparen{1-\frac{\mu}{L}}^tD_{h_{\vq}}(\vx_{\vq},\vq)+\max_{1\le i \le t} \ip{\triangle_i, \vx_t-\vx_{\vq}}\enspace,\\
    &= 2L\inparen{1-\frac{\mu}{L}}^th_{\vq}(\vx_{\vq})+\max_{1\le i \le t} \ip{\triangle_i, \vx_t-\vx_{\vq}}\enspace.
\end{align*}
Hence, for $t \ge \frac{L}{\mu}\logv{Lh_{\vq}(\vx_{\vq})\inparen{\frac{4}{p\gamma}}^{p}}$, it is easy to check that for $p\geq 2$,
\begin{align*}
    f_{\vq}(\vx_t) - f_{\vq}(\vx_{\vq}) &\leq 2L\inparen{\frac{1}{e}}^{\logv{Lh_{\vq}(\vx_{\vq})\inparen{\frac{4}{p\gamma}}^{p}}}h_{\vq}(\vx_{\vq}) + \max_{1\le i \le t} \ip{\triangle_i, \vx_t-\vx_{\vq}}\enspace,\\
    &= 2\inparen{\frac{p\gamma}{4}}^{p} + \max_{1\le i \le t} \ip{\triangle_i, \vx_t-\vx_{\vq}}\enspace,\\
    &\leq \inparen{\frac{p\gamma}{2}}^{p} + \max_{1\le i \le t} \ip{\triangle_i, \vx_t-\vx_{\vq}}\enspace,
\end{align*}
and combining this with \Cref{lemma:gp_prox_relativesmoothsolve} to make the error term on the order of our accuracy, we get $\norm{\vx_{\vq}-\vx}_{\mM} \lesssim \gamma$. We thus conclude the proof of \Cref{lemma:prox_subproblem_base}. \todo{I am not sure what is happening in the last line....}
\end{proof}

The last step is to use our proximal problem solver to build a valid MS oracle.

\begin{lemma}
\label{lemma:prox_subproblem_ms}
In the context of \Cref{alg:optimal_ms_acceleration}, there exists an algorithm $(\xtilde_{t+1},\lambda_{t+1})=\oprox(\vq_t)$ that approximately solves
\begin{align*}
    \argmin{\xtilde\in\R^d} f(\xtilde)+ep^p\norm{\xtilde-\vq_t}_{\mM}^p
\end{align*}
using $O\inparen{p^{O(1)}\logv{\frac{pd \cdot f(\vx_t)}{\eps}}}$ linear system solves in $\mA^{\top}\mB\mA$, in the sense that
\begin{align*}
    \norm{\frac{1}{ep^{p+1}\norm{\xtilde_{t+1}-\vq_t}_{\mM}^{p-2}}\mM^{-1}\nabla f(\xtilde_{t+1})+(\xtilde_{t+1}-\vq_t)}_{\mM} \le \frac{1}{2}\norm{\xtilde_{t+1}-\vq_t}_{\mM}.
\end{align*}
\end{lemma}
\begin{proof}[Proof of \Cref{lemma:prox_subproblem_ms}]
The point of this proof is to give an analysis of \Cref{alg:gp_prox_solver}.

For notational simplicity, let $\vx=\xtilde_{t+1}$ and $\lambda = \lambda_{t+1}$. We will reintroduce the indices when it is essential to clarify the iterations we are discussing.

First, it is helpful to see why the stated notion of approximation is useful. Let $C_p \coloneqq ep^p$. Observe that at exact optimality, we have
\begin{align}
    \nabla f(\vx_{\vq}) + \underbrace{ep^{p+1}\norm{\vx_{\vq}-\vq}_{\mM}^{p-2}}_{\lambda^{\star}}\mM(\vx-\vq)=0\enspace.\label{eq:prox_problem_stationarity_condition}
\end{align}
This motivates the approximation in our lemma statement, with us asking for a $\frac{1}{2}$-approximate MS oracle (\Cref{defn:ms_oracle}) for $f$. This also tells us that at optimality in \eqref{eq:prox_problem_stationarity_condition}, we have,
\begin{align*}
    &\nabla f(\vx_{\vq}) + ep^{p+1}\norm{\vx_{\vq}-\vq}_{\mM}^{p-2}\mM(\vx-\vq)=0\enspace,\\
    &\Leftrightarrow \mM^{-1/2}f(\vx_{\vq}) = -pC_p\norm{\vx_{\vq}-\vq}_{\mM}^{p-2}\mM^{1/2}(\vx-\vq)\enspace,\\
    &\Rightarrow \norm{\mM^{-1/2}f(\vx_{\vq})}_{2} = pC_p\norm{\vx_{\vq}-\vq}_{\mM}^{p-2}\norm{\mM^{1/2}(\vx-\vq)}_2\enspace,\\
    &\Leftrightarrow\norm{\vx_{\vq}-\vq}_{\mM} = \inparen{\frac{\norm{\mM^{-1}\nabla f(\vx_{\vq})}_{\mM}}{pC_p}}^{\frac{1}{p-1}}\enspace.
\end{align*}
We now break up our analysis into two cases. In the first, suppose that $\norm{\mM^{-1}\nabla f(\vx_{\vq})}_{\mM} \le \eps/\norm{\vx_{\vq}-\xstar}_{\mM}$. Then, by convexity, we have
\begin{align*}
    f(\vx_{\vq})-f(\xstar) \le \ip{\nabla f(\vx_{\vq}), \vx_{\vq}-\xstar} \le \norm{\mM^{-1}\nabla f(\vx_{\vq})}_{\mM}\norm{\vx_{\vq}-\xstar}_{\mM} \le \eps.   
\end{align*}
Hence, for the rest of the proof, assume that $\norm{\mM^{-1}\nabla f(\vx_{\vq})} \ge \eps/\norm{\vx_{\vq}-\xstar}_{\mM}$ (because if this is not the case, in the algorithm we can simply check whether the MS condition is satisfied -- if not, then we know this assumption was violated and we are done anyway)\todo{I am not sure I follow the statement in the parenthesis, maybe I need to discuss with you the hierarchy of the algorithms and see where specifically this termination condition comes in}. We run the algorithm implied by \Cref{lemma:prox_subproblem_base} and obtain an approximate solution $\vx$ for which
\begin{align}
    \norm{\vx-\vx_{\vq}}_{\mM} \le \alpha\norm{\vx_{\vq}-\vq}_{\mM} \text{ for } \alpha = \frac{1}{5}\min\inbraces{\frac{C_p}{ep(p-1)}\inparen{\frac{\norm{\vx_{\vq}-\vq}_{\mM}}{f(\vq)^{\frac{1}{p}}}}^{p-2}, 1}\enspace.\label{eq:lemma_D19_guarantee}
\end{align}
Since $\alpha < 1$ the guarantee in \eqref{eq:lemma_D19_guarantee} gives us,
\begin{align}
    \norm{\vx-\vx_{\vq}}_{\mM} \le \alpha\norm{\vx-\vq}_{\mM} \leq \frac{\alpha}{1-\alpha}
    \norm{\vx-\vq}_{\mM}\enspace,\label{eq:xq_x_UB_1}
\end{align}
and further applying triangle inequality gives us
\begin{align}
    \norm{\vx_{\vq}-\vq}_{\mM} &\le \norm{\vx-\vq}_{\mM} +  \norm{\vx_{\vq}-\vx}_{\mM}\enspace,\nonumber\\ 
    &\le \frac{1-\alpha}{1-\alpha}\norm{\vx-\vq}_{\mM} + \frac{\alpha}{1-\alpha}\norm{\vx-\vq}_{\mM}\enspace,\nonumber\\
    &\le \frac{1}{1-\alpha}\norm{\vx-\vq}_{\mM}\enspace.\label{eq:xq_x_UB_2}
\end{align}
Hence, we get
\begin{align}
    \frac{ep(p-1)f(\vq)^{1-\frac{2}{p}}}{C_p\norm{\vx-\vq}_{\mM}^{p-2}} \cdot \norm{\vx-\vx_{\vq}}_{\mM} &= \frac{ep(p-1)}{C_p} \cdot \inparen{\frac{f(\vq)^{\frac{1}{p}}}{\norm{\vx-\vq}_{\mM}}}^{p-2} \cdot \norm{\vx-\vx_{\vq}}_{\mM}\enspace,\nonumber\\
    &\le^{\eqref{eq:lemma_D19_guarantee}} \frac{1}{5}\norm{\vx_{\vq}-\vq}_{\mM}\enspace,\nonumber\\ &\le^{\eqref{eq:xq_x_UB_2}} \frac{1}{5}\cdot\frac{1}{1-\alpha}\norm{\vx-\vq}_{\mM}\enspace,\nonumber\\
    &\leq \frac{1}{4}\norm{\vx-\vq}_{\mM}\enspace,\label{eq:xq_x_UB_3}
\end{align}
where in the last inequality, we used that $\alpha \leq \frac{1}{5}$ due to our choice in \eqref{eq:lemma_D19_guarantee}. We now call \Cref{lemma:self_gradient_norm_small}, divide both sides by $\lambda$, and get\todo{I tried to make sense of the steps here, unclear how the first and last inequalities are working.}
\begin{align*}
    &\norm{\frac{1}{ep^{p+1}\norm{\vx-\vq}_{\mM}^{p-2}}\mM^{-1}\nabla f(\vx)+(\vx-\vq)}_{\mM} \\
    &\le^{\text{(\Cref{lemma:self_gradient_norm_small})}} ep(p-1)\inparen{\frac{f(\vq)^{1-\frac{2}{p}}}{C_p\norm{\vx-\vq}_{\mM}^{p-2}}+\max\inbraces{1, \inparen{\frac{\norm{\vx_{\vq}-\vq}_{\mM}}{\norm{\vx-\vq}_{\mM}}}^{p-2}}}\norm{\vx-\vx_{\vq}}_{\mM}\enspace,\\
    &\le^{\text{\eqref{eq:xq_x_UB_2}}} ep(p-1)\inparen{\frac{f(\vq)^{1-\frac{2}{p}}}{C_p\norm{\vx-\vq}_{\mM}^{p-2}}+\frac{1}{(1-\alpha)^{p-2}}}\norm{\vx-\vx_{\vq}}_{\mM}\enspace,\\
    &\le^{\eqref{eq:xq_x_UB_1}} \frac{ep(p-1)f(\vq)^{1-\frac{2}{p}}}{C_p\norm{\vx-\vq}_{\mM}^{p-2}} \cdot \norm{\vx-\vx_{\vq}}_{\mM} + \frac{ep(p-1)\alpha}{(1-\alpha)^{p-1}}\norm{\vx-\vq}_{\mM}\enspace,\\ 
    &\leq^{\eqref{eq:xq_x_UB_2},\ \eqref{eq:lemma_D19_guarantee}} \frac{1}{4}\norm{\vx-\vq}_{\mM} + \frac{ep(p-1)5^{p-2}}{4^{p-1}}\norm{\vx-\vq}_{\mM}\enspace,\\
    &\le \frac{1}{2}\norm{\vx-\vq}_{\mM}\enspace,
\end{align*}
giving us the approximation guarantee.

It remains to understand the complexity of solving the proximal subproblem to the accuracy required in \eqref{eq:lemma_D19_guarantee}. Plugging in $\gamma = \alpha\norm{\vx_{\vq}-\vq}_{\mM}$  into \Cref{lemma:prox_subproblem_base} and using our bound on $h_{\vq}(\vx_{\vq})$ from \Cref{lemma:initial_bregman_div} gives an iteration complexity of (ignoring the constant in front of the big-$O$)
\begin{align*}
    &\quad p^{O(1)}\logv{ph_{\vq}(\vx_{\vq})\inparen{\frac{2}{p\alpha\norm{\vx_{\vq}-\vq}_{\mM}}}^{p}} \\
    &\le p^{O(1)}\logv{p\inparen{p(p-1)f(\vq)^{1-\frac{2}{p}}\norm{\vx_{\vq}-\vq}_{\mM}^2 + C_p\norm{\vx_{\vq}-\vq}_{\mM}^p}\inparen{\frac{2}{p\alpha\norm{\vx_{\vq}-\vq}_{\mM}}}^{p}} \\
    &= p^{O(1)}\logv{\inparen{\frac{2}{p}}^{p} p\inparen{\frac{p(p-1)f(\vq)^{1-\frac{2}{p}}\norm{\vx_{\vq}-\vq}_{\mM}^2 + C_p\norm{\vx_{\vq}-\vq}_{\mM}^p}{\alpha^p\norm{\vx_{\vq}-\vq}_{\mM}^p}}} \\
    &= p^{O(1)}\logv{\inparen{\frac{2}{p}}^{p} p\inparen{\frac{p(p-1)f(\vq)^{1-\frac{2}{p}}}{\alpha^p\norm{\vx_{\vq}-\vq}_{\mM}^{p-2}} + \frac{C_p}{\alpha^p}}}
\end{align*}
We have two cases to analyze for the value of $\alpha$. In the first, suppose we get $\alpha = \frac{1}{5}$. By the definition of $\alpha$, this means we have
\begin{align*}
    \frac{C_p}{ep(p-1)}\inparen{\frac{\norm{\vx_{\vq}-\vq}_{\mM}}{f(\vq)^{\frac{1}{p}}}}^{p-2} \ge 1,
\end{align*}
which means the complexity we get is $p^{O(1)}\log p$. We now handle the other case, i.e., $\alpha=\frac{C_p}{5ep(p-1)}\inparen{\frac{\norm{\vx_{\vq}-\vq}_{\mM}}{f(\vq)^{\frac{1}{p}}}}^{p-2}$. Here, it will be useful to keep track of the timestep $t$ that we are working with. Recall that
\begin{align}
    \norm{\vx_{\vq_t}-\vq_t}_{\mM}^p = \inparen{\frac{\norm{\mM^{-1}\nabla f(\vx_{\vq_t})}_{\mM}}{pC_p}}^{\frac{p}{p-1}} \ge \inparen{\frac{\eps}{pC_p\norm{\vx_{\vq_t}-\xstar}_{\mM}}}^{\frac{p}{p-1}}\enspace,\label{eq:xq_lb}
\end{align}
so the complexity we want to control is given by
\begin{align*}
    &p^{O(1)}\logv{\inparen{\frac{2}{p}}^p p\inparen{\frac{2f(\vq_t)}{\alpha^p\norm{\vx_{\vq_t}-\vq_t}_{\mM}^{p}}}}\\ &\qquad\lesssim^{\eqref{eq:lemma_D19_guarantee}} p^{O(1)}\logv{\inparen{\frac{2}{p}}^p p\inparen{\frac{2\inparen{5ep(p-1)}^{p}f(\vq_t)^{p-1}}{C_p^p\norm{\vx_{\vq_t}-\vq_t}_{\mM}^{p(p-2)}\norm{\vx_{\vq_t}-\vq_t}_{\mM}^{p}}}}\enspace,\\
    &\qquad\lesssim p^{O(1)}\logv{ p\inparen{\frac{2\inparen{10(p-1)}^{p}f(\vq_t)^{p-1}}{p^{p^2}\norm{\vx_{\vq_t}-\vq_t}_{\mM}^{p(p-1)}}}}\enspace,\\
    &\qquad\lesssim^{\eqref{eq:xq_lb}} p^{O(1)}\logv{ p\inparen{\frac{2\inparen{10e(p-1)}^{p}p^{p(p+1)}f(\vq_t)^{p-1}}{p^{p^2}\epsilon^p}}\norm{\vx_{\vq_t}-\xstar}_{\mM}^p}\enspace,\\
    &\qquad\lesssim^{\eqref{eq:xq_lb}} p^{O(1)}\logv{\inparen{\frac{2\inparen{10e(p-1)}^{p}p^{p+1}f(\vq_t)^{p-1}}{\epsilon^p}}\norm{\vx_{\vq_t}-\xstar}_{\mM}^p}\enspace,\\
    &\qquad\lesssim p^{O(1)}\logv{\frac{pf(\vq_t)\norm{\vx_{\vq_t}-\xstar}_{\mM}}{\eps}}\enspace,\\
    &\qquad\lesssim^{\text{(\Cref{lemma:gp_regression_prox_diameter})}} p^{O(1)}\logv{\frac{pf(\vq_t)df(\vx_t)}{\eps}}\enspace,\\
    &\qquad\lesssim^{\text{(\Cref{lemma:fq_bounded})}} p^{O(1)}\logv{\frac{pf(\vx_t)}{\eps}}\enspace,
\end{align*}
completing the proof of \Cref{lemma:prox_subproblem_ms}.
\end{proof}

\subsection{The algorithm}
\label{sec:gp_interpolating_alg}

We are now ready to combine the results from the previous two subsections to build our algorithm for $\cG_p$-regression and prove \Cref{mainthm:gp_regression_iteration_complexity}. The main algorithmic object here is \Cref{alg:gp_regression_alg}.

\begin{algorithm}[H]
\caption{\textsf{GpRegression}: Optimizes \eqref{eq:interpolation} up to $(1+\eps)$-multiplicative error}
\label{alg:gp_regression_alg}
\begin{algorithmic}[1]
\Require Regression problems $(\mA_{S_1},\vb_{S_1}),\dots,(\mA_{S_m},\vb_{S_m})$, accuracy $\eps > 0$
\State Using \cite[Algorithm 2]{mo23} with input $\insquare{\mA \vert \vb}$, find nonnegative diagonal $\mW$ such that for all $\vx\in\R^{d}$ and $c \in \R$,
\begin{align*}
    \gnorm{\mA\vx-c\vb}{\infty} \le \norm{\mW^{\frac{1}{2}-\frac{1}{p}}\mA\vx-c\mW^{1/2}\vb}_2 \le (2(d+1))^{\frac{1}{2}-\frac{1}{p}}\gnorm{\mA\vx-c\vb}{\infty}.
\end{align*}
\State Let $\vx_0 = \inparen{\mA^{\top}\mW^{1-\frac{2}{p}}\mA}^{-1}\mA^{\top}\mW^{1-\frac{2}{p}}\vb$. \Comment{$\vx_0 \coloneqq \argmin{\vx\in\R^d} \norm{\mW^{\frac{1}{2}-\frac{1}{p}}\mA\vx-\mW^{\frac{1}{2}-\frac{1}{p}}\vb}_2$.}
\State Using \Cref{alg:gp_prox_solver} and \Cref{lemma:prox_subproblem_ms}, implement a $\frac{1}{2}$-MS oracle for $f$ (\Cref{defn:ms_oracle})
\State Run \Cref{alg:optimal_ms_acceleration} with the oracle from the previous line and with $\vx_0$ as the initialization for $O\inparen{\mathsf{poly}(p)\min\inbraces{\rank{\mA},m}^{\frac{p-2}{3p-2}}\logv{\frac{d}{\eps}}^3}$ iterations.
\State \Return $\xhat$ the output of the previous step.
\end{algorithmic}
\end{algorithm}

\begin{proof}[Proof of \Cref{mainthm:gp_regression_iteration_complexity}]
By writing the stationary condition of the proximal problem, it makes sense to choose $\lambda_{t+1} = ep^{p+1}\norm{\xtilde_{t+1}-\vq_t}_{\mM}^{p-2}$.

It is easy to check that
\begin{align*}
    \norm{\xtilde_{t+1}-\vq_t}_{\mM} = \inparen{\frac{ep^{p+1}\norm{\xtilde_{t+1}-\vq_t}_{\mM}^{p-2}}{\inparen{(ep^{p+1})^{\frac{1}{p-1}}}^{p-1}}}^{\frac{1}{(p-1)-1}},
\end{align*}
and therefore the triple $(\xtilde_{t+1},\vq_t,ep^{p+1}\norm{\xtilde_{t+1}-\vq_t}_{\mM}^{p-2})$ always satisfies a $(p-1, (ep^{p+1})^{1/(p-1)})$-movement bound (\Cref{defn:movement_bound}).

Next, we calculate the iteration complexity we need to reduce the error to half of what we started with. For an arbitrary initial iterate $\vx$, let $\delta = 0.5(f(\vx)-f(\xstar))$. By \Cref{lemma:strong_convexity_gp}, we have
\begin{align*}
    \norm{\vx-\xstar}_{\mM}^{s+1} = \norm{\vx-\xstar}_{\mM}^p \le 2^{3p/2}d^{p/2-1}(f(\vx)-f(\xstar)),
\end{align*}
so combining this along with the fact that $c^s=ep^{p+1}$ and applying \Cref{thm:optimal_ms_acceleration} with our proximal solver \Cref{lemma:prox_subproblem_ms} yields
\begin{align*}
    T_{\mathsf{min}} &= \frac{p-1}{3}\inparen{pC_p\cdot2^{3p/2+1}d^{p/2-1}}^{\frac{2}{3p-2}} \lesssim p^{5/3}d^{\frac{p-2}{3p-2}}.
\end{align*}
Next, we initialize $\vx_0 \coloneqq \inparen{\mA^{\top}\mW^{1-2/p}\mA}^{-1}\mA^{\top}\mW^{1-2/p}\vb$. Using \Cref{thm:gp_regression_blw_intro}, we have
\begin{align*}
    f(\vx_0) \le (2d)^{p/2-1}f(\xstar),
\end{align*}
so reaching an iterate $\vx$ for which $f(\vx)-f(\xstar) \le \eps f(\xstar)$ takes $T_{\mathsf{min}} \cdot \logv{d^{p/2-1}/\eps} = p^{8/3}d^{\frac{p-2}{3p-2}}\logv{\frac{d}{\eps}}$ calls to $\oprox$.

We now resolve the full iteration complexity, including the bootstrapping step to show that $f(\vx_t)$ is reasonably bounded so that we get an unconditional upper bound from \Cref{lemma:prox_subproblem_ms}. At the end of iteration $t$, from (loosely) inverting the bound in \Cref{thm:optimal_ms_acceleration}, we know that
\begin{align*}
    f(\vx_t)-f(\xstar)\le \frac{(Cp^3)^{\frac{3p-2}{2}}(2d)^{\frac{p}{2}-1}}{t^{\frac{3p-2}{2}}}.
\end{align*}
Since $\xtilde_{t+1}$ only depends on $\vq_t$, which in turn only depends on $\vx_t$ and $\vv_t$, it suffices to use the above bound for $f(\vx_t)$, which gives us an iteration complexity of $p^{O(1)}\logv{\frac{pd}{\eps}}$ to compute $\xtilde_{t+1}$ (which we get from plugging into \Cref{lemma:prox_subproblem_ms}).

Combining this with the iteration complexity of $\oprox$ gives us the result of \Cref{mainthm:gp_regression_iteration_complexity}.
\end{proof}

\chapter{Dueling optimization with a monotone adversary \label{chapter:binsearch}}

In this chapter, we give randomized algorithms for the problem of \textit{dueling optimization with a monotone adversary}. The content here is based on joint work with Avrim Blum, Meghal Gupta, Gene Li, Aadirupa Saha, and Chloe Yang \cite{bglmsy23}.

\section{Introduction}
A growing body of literature studies learning with preference-based feedback \citep{bv06, shivaswamy2011online}, with tremendous empirical success in recommendation systems, search engine optimization, information retrieval, and robotics. More recently, preference-based feedback has received a lot of attention as a mechanism to train large language models \citep{ouyang2022training}. Moreover, in recommender systems \citep{bobadilla2013recommender}, a natural approach is to learn from users' preferences relations on a set of recommended items and update the system's belief for better future recommendations~\cite{j16} (e.g., given these items, which one do you prefer the most?). 

Such preference-based feedback is not readily addressed by classical formulations for online decision making, such as bandits and reinforcement learning. In particular, algorithms for these problems rely on ordinal feedback per item (e.g., on a scale of 1 to 10, how much did the user like a particular item?). To address this, a long line of work studies the \emph{dueling bandit framework} for online decision making under pairwise/preference-based feedback. There exist efficient algorithms with provable guarantees for the standard multi-armed bandit setup \citep{yue2012k, ailon2014reducing, komiyama2015regret}, contextual bandits \cite{dudik2015contextual,saha2022efficient}, as well as dueling convex optimization \cite{jamieson2012query,skm21,saha2022dueling}, to name a few. The dueling bandit framework is especially applicable in settings where real-valued feedback is scarce or impossible to obtain, but preference-based feedback is readily available.

However, a key limitation of the dueling bandit framework is that the feedback that the learner receives is essentially ``in-list''. That is, the users are restricted to selecting items exclusively from the list of recommended items. This feedback model fails to capture the real-world scenarios where the users might select an out-of-list item they prefer. To illustrate, music streaming services like Spotify create personalized playlists for users. Concretely, each song can be encoded as a feature vector $\vx \in \R^d$, and the goal is to recommend the songs with the highest utility for a hidden, well-structured utility function of $\vx$. However, the users can also search for and play the songs they have a stronger preference (i.e., higher utility) than all recommendations. 

This out-of-list feedback model falls into a monotone adversarial framework (see the chapter by \citet{feige_2021}). In such models, an adversary is only allowed to make ``helpful'' changes. For example, in a graph clustering problem, the adversary is only allowed to add edges within communities and delete edges that cross communities (see, e.g., the chapter by \citet{moitra_2021}). In our setting, the adversary is only allowed to respond with an item that is at least as good as any recommended item. A clear adaptation of the dueling bandit framework to this new feedback type is not evident. 

\subsection{Problem statement}

As our main conceptual contribution, we introduce a theoretical formulation for this setting that we call \textit{dueling optimization with a monotone adversary}. As we will see, our formulation supports ``out-of-list'' feedback.

\begin{problem}[Dueling optimization with a monotone adversary]
\label{problem:local_rec_fixed_set}
Let $\mathcal{X} \subseteq \R^d$ be a decision space, and let $\myfunc{f}{\mathcal{X}}{\R}$ be a cost function with an unknown global minimum $\xstar$. A learner interacts with an adversary over rounds $t = 1, 2, \dots$, where each round is of the following form.
\begin{enumerate}
    \item The learner proposes $m$ points $\xt{1},\dots, \xt{m} \in \mathcal{X}$.
    \item The adversary responds with a point $\xstar_t$ that satisfies
    $f(\xstar_t) \le \min_{1\le j \le m} \inbraces{f\inparen{\xt{j}}}\label{eq:valid_feedback}$.
\end{enumerate}
The goal is to design algorithms that:
\begin{enumerate}
    \item for some prespecified $\eps > 0$, minimize the number of iterations to find a point $\vx$ for which $f(\vx) - f(\xstar) \le \eps$;
    \item minimize the total \textit{cost} $\sum_{t = 1}^{\infty} \inparen{\max_{1\le j \le m}\inbraces{f\inparen{\xt{j}}} - f(\xstar)}$.
\end{enumerate}
\end{problem}

Note that in \Cref{problem:local_rec_fixed_set}, we are interested in both the iteration complexity and the total cost. The first objective is a standard metric for measuring the performance of an iterative optimization algorithm. The second objective is motivated by online settings in which a practitioner may wish to minimize the total regret (cost) of its recommendations over an indefinitely long interaction with a user. In fact, the algorithms we propose in this chapter simultaneously achieve both small iteration complexity (for any choice of $\eps$) as well as total cost --- see our technical overview in \Cref{sec:tech-overview} for more details.

\Cref{problem:local_rec_fixed_set} is a natural extension of (noiseless) dueling optimization \cite{jamieson2012query,skm21,saha2022dueling} to handle ``out-of-list'' responses, as in the Spotify recommendation example. The vanilla (noiseless) dueling optimization setup corresponds to the requirement that the user's response satisfies $\xstar_t \in \{\xt{1}, \xt{2}\}$. We allow the user to be potentially adversarial by allowing it to respond with any improvement to the learner's suggestions (in the sequel, we exclusively refer to the user as the adversary).

Even though the monotone adversary is only improving upon the learner's suggestions, existing algorithms for dueling optimization cannot be freely extended to handle the monotone feedback. At a high level, the difficulty arises from the fact that existing algorithms carefully select the queries $\xt{1}, \xt{2}$ so that learning whether $f(\xt{1}) > f(\xt{2})$ reveals information about the underlying $f$. However, a monotone adversary can return a point $\xstar_t$ that reveals no information about the relationship between $\xt{1}$ and $\xt{2}$.

To illustrate this point, consider a natural coordinate-wise binary search algorithm for the dueling optimization problem when $f(\vx) = \norm{\vx-\xstar}_2^2$ for some $\xstar \in \cB_2^d \coloneqq \inbraces{\vx \suchthat \norm{\vx}_2 \le 1}$. For coordinates $i = 1, \cdots, d$, query points of the form $\xt{1} = c_1 \cdot \ve_i$, $\xt{2} = c_2 \cdot \ve_i$ and progressively refine the values $c_1, c_2 \in \R$ to search for the value of $\xstar[i]$ (i.e., the $i$-th entry of $\xstar$). It is easy to show that this approach has a query complexity of $O\inparen{d \logv{\nfrac{1}{\eps}}}$ in the vanilla dueling optimization setting. However, a monotone adversary can return orthogonal responses of the form $\xstar_t = C\ve_j$ (where $j \ne i$ and $C$ is a constant) that do not allow the learner to search along the intended coordinate $i$. Furthermore, \citet{jamieson2012query} and \citet{skm21} give more sophisticated algorithms for the dueling optimization problem that inherently depend upon the ``in-list'' feedback, which clearly cannot apply to our setting. We therefore need novel insights to solve \Cref{problem:local_rec_fixed_set}.

\subsection{Our results}
We study \Cref{problem:local_rec_fixed_set} for various natural classes of functions $f$ and provide tight upper and lower bounds on the number of queries required to find an $\eps$-optimal point.

\paragraph{Upper bound for linear functions.} First, we study dueling optimization with a monotone adversary when the function $f$ is linear. This is a natural class to consider. In particular, an algorithm that solves \Cref{problem:local_rec_fixed_set} can be adapted to achieve constant regret for (noiseless) linear contextual bandits \cite{chu2011contextual}, where the reward function is $r(\vx) \coloneqq \ip{\vx, \vx^\star}$. Note that the key difference in the setup is that the learner does not get to observe the actual linear costs but instead only an improvement to the actions (points) that the learner selects.

\begin{mainthm}
\label{cor:bin_search_ip}
Let $m=2$, let $\mathcal{X} = \mathbb{S}_2^d$, let $\xstar$ be such that $\norm{\xstar}_2 = 1$, and let $\myfunc{f}{\mathcal{X}}{\R}$ be $f(\vx) = -\ip{\vx,\xstar}$. Fix any $\eps > 0$. There exists an algorithm that, in the setting of \Cref{problem:local_rec_fixed_set}, with probability at least $1-\expv{-O(d)}$:
\begin{itemize}
    \item outputs a point $\vx$ satisfying $\ip{\xstar-\vx,\xstar} \le \eps$ within $O(d\logv{\nfrac{1}{\eps}}^2)$ iterations;
    \item incurs total cost $O(d)$.
\end{itemize}
Each pair of guesses at time $t$ can be computed in $O(d)$ time.
\end{mainthm}

We prove \Cref{cor:bin_search_ip} in \Cref{subsec:ip_proof}, and the cost is near-optimal with respect to $d$. 

\citet{gollapudi2021contextual} study a closely related setup that they call \emph{local contextual recommendation}. Their result (see their Theorem 6.4) can be interpreted as showing that if the action set $\actionset$ is a discrete set (namely a packing over the unit sphere), there exists a $2^{\Omega(d)}$ lower bound on the iteration complexity to find a point with constant suboptimality. In contrast, our \Cref{cor:bin_search_ip} shows a much smaller upper bound when the domain is the entire unit sphere. 

\paragraph{Upper bound for smooth and P\L{} functions.} Next, we study whether we can show guarantees for a large class of functions. We show a positive result for functions that are both $\beta$-smooth and $\alpha$-Polyak-\L{}ojasiewicz (abbreviated as P\L{}). These assumptions are standard in optimization.

\begin{definition}[$\beta$-smooth function {\cite[{Lemma 3.4}]{bubeck2015convex}}]
\label{defn:smoothness}
We say $f$ is $\beta$-smooth if it satisfies (\ref{eq:smoothness}).
\begin{align}
    \text{For all } \vx, \vy \in \R^d &:\quad \abs{f(\vx)-f(\vy) - \ip{\nabla f(\vy), \vx-\vy}} \le \frac{\beta}{2} \cdot \norm{\vx-\vy}_2^2\label{eq:smoothness}
\end{align}
\end{definition}
\begin{definition}[$\alpha$-P\L{} function]
\label{defn:pl}
We say $f$ is $\alpha$-P\L{} if it satisfies (\ref{eq:pl}).
\begin{align}
    \text{ For all } \vx \in \R^d \text{ and minimizers } \xstar &:\quad f(\vx) -f(\xstar) \le \frac{1}{2\alpha}\norm{\nabla f(\vx)}_2^2 \label{eq:pl}
\end{align}
\end{definition}
Our main result for this setting is \Cref{cor:bin_search_pl}.

\begin{mainthm}
\label{cor:bin_search_pl}
Let $\mathcal{X} = \R^d$, and suppose $f$ is $\beta$-smooth (\Cref{defn:smoothness}) and $\alpha$-P\L{} (\Cref{defn:pl}). Fix any $\eps > 0$, as well as a known point $\vx_1$ and a value $B$ satisfying $B \ge f(\vx_1) - f(\xstar)$. For all $d$ larger than a universal constant, there exists an algorithm that, in the setting of \Cref{problem:local_rec_fixed_set} with $m\ge 2$, with probability at least $1-\expv{-O(d)}$:
\begin{itemize}
    \item outputs a point $\vx$ satisfying $f(\vx) - f(\xstar) \le \eps$ within $O\inparen{\nfrac{\beta}{\alpha} \cdot \nfrac{d}{\log m} \cdot \logv{\nfrac{B}{\eps}}^2}$ iterations;
    \item incurs total cost $O\inparen{\nfrac{\beta}{\alpha} \cdot B \cdot \nfrac{d}{\log m}}$.
\end{itemize}
The list of $m$ guesses at time $t$ can be computed in $O(md)$ time.
\end{mainthm}

We prove \Cref{cor:bin_search_pl} in \Cref{subsec:pl_proof}. Importantly, observe that the results above generalize those of \citet{sfkm24}, as our methods also work under a monotone adversary. Additionally, these achieve the dependence on the list size $m$ that we observe in our lower bounds (to be presented momentarily), which therefore makes our results tight. Finally, although \Cref{cor:bin_search_pl} is only stated for smooth and P\L{} functions, it is straightforward to adapt this to a result for convex and smooth functions (we describe this adaptation in a moment).

As an application, we show a positive result when the loss function is the Euclidean distance, and the decision space $\mathcal{X} = \cB_2^d$ is a unit ball: 

\begin{mainthm}
\label{cor:bin_search_dist}
Let $m=2$, let $\mathcal{X} = \cB_2^d$, let $\xstar$ be such that $\norm{\xstar}_2 \le 1$, and let $\myfunc{f}{\mathcal{X}}{\R}$ be $f(\vx) = \norm{\vx-\xstar}_2$. Fix any $\eps > 0$. There exists an algorithm that, in the setting of \Cref{problem:local_rec_fixed_set}, with probability at least $1 - \expv{-O(d)}$:
\begin{itemize}
    \item outputs a point $\vx$ satisfying $\norm{\vx-\xstar}_2 \le \eps$ within $O\inparen{d \cdot \logv{\nfrac{B}{\eps}}^2}$ iterations;
    \item incurs total cost $O\inparen{d}$.
\end{itemize}
Each pair of guesses at time $t$ can be computed in $O(d)$ time.
\end{mainthm}

We prove \Cref{cor:bin_search_dist} in \Cref{subsec:dist_proof}.

Note that unlike in \Cref{cor:bin_search_pl}, \Cref{cor:bin_search_dist} applies to the setting where the algorithm must guess points belonging to a given constraint set $\actionset$. Hence, in the proof of \Cref{cor:bin_search_dist}, we have to be careful to ensure that the convergence argument still holds when we apply the algorithm for \Cref{cor:bin_search_pl} along with a projection step. It is not clear that this argument holds by default for all $f$ satisfying the conditions requested by \Cref{cor:bin_search_pl}. Furthermore, as will become evident, we really only require that $\actionset$ be any convex body (though we state the result with $\actionset = \cB_2^d$ to emphasize the consistency with our following lower bounds).

Finally, as another corollary to \Cref{cor:bin_search_pl}, we prove a low-accuracy result for optimizing functions that are $\beta$-smooth and convex.

\begin{mainthm}
\label{thm:mary_opt_smooth}
Let $\mathcal{X} = \R^d$, and suppose $f$ is convex and $\beta$-smooth (\Cref{defn:smoothness}). Fix any $\eps > 0$, as well as a known point $\vx_1$ and a value $B$ satisfying $B \ge f(\vx_1) - f(\xstar)$ and $\norm{\vx_1-\xstar}_2 \le \sqrt{D}$. For all $d$ larger than a universal constant, there exists an algorithm that, in the setting of \Cref{problem:local_rec_fixed_set} with $m\ge 2$, with probability at least $1-\expv{-O(d)}$, outputs a point $\vx$ satisfying $f(\vx) - f(\xstar) \le \eps$ within $O\inparen{\nfrac{\beta D}{\eps} \cdot \nfrac{d}{\log m} \cdot \logv{\nfrac{B}{\eps}}^2}$ iterations. The list of $m$ guesses at time $t$ can be computed in $O(md)$ time.
\end{mainthm}

We prove \Cref{thm:mary_opt_smooth} in \Cref{sec:mary_opt_smooth}.

\paragraph{Lower bounds.} We also prove that the dependence on $d$ in our results is tight. In particular, when $f$ is either a linear function or the distance to the target (as in \Cref{cor:bin_search_dist}), then $\Omega(d)$ queries are necessary to identify $\xstar$. This will translate to a $\Omega(d)$ cost over an infinite number of rounds. In fact, our lower bound is valid when the adversary must return one of the two queried points, as in vanilla dueling optimization framework. 

Our lower bound also covers a more general setting than that stated in \Cref{problem:local_rec_fixed_set}. Thus far, we have only discussed the setting where the algorithm can query only two points and is told the better of the two. In many practical instances, the algorithm can query $m$ points and learn the point with the best objective value (we call this $m$-ary dueling optimization). In our construction, we prove that unless $m$ is polynomial in $d$, we cannot decrease the total cost substantially below $\Omega(d)$. Thus, \Cref{cor:bin_search_pl} is tight.

See \Cref{thm:large_m_lower} for a formal statement of our lower bound.

\begin{mainthm}[Lower bound, $\ell_2$ distance]
\label{thm:large_m_lower}
Let $\cX = \cB_2^d$. For any randomized algorithm for $m$-ary dueling optimization, there exists a choice of minimizer $\xstar \in \cB_2^d$ and function $f(\vx) \coloneqq \norm{\vx-\xstar}_2$ such that the algorithm must:
\begin{itemize}
    \item perform $\Omega\inparen{\nfrac{d}{\log m}}$ iterations in expectation to find a point $\vx$ for which $f(\vx) - f(\xstar) \le \eps$.
    \item incur cost  $\Omega\inparen{\nfrac{d}{\log m}}$ in expectation.
\end{itemize}
Here, $\eps > 0$ is an absolute numerical constant.
\end{mainthm}


We prove \Cref{thm:large_m_lower} in \Cref{sec:lower-bound}. Using the same construction, we can also demonstrate that \Cref{cor:bin_search_ip} is tight when $\mathcal{X}$ is the unit sphere.

\begin{corollary}[Lower bound, linear $f$]\label{corr:inner-prod-lower-bound}
Let $\mathcal{X} = \mathbb{S}_2^d$. For any randomized algorithm for $m$-ary dueling optimization there exists a choice of minimizer $\xstar \in \mathbb{S}_2^d$ and function $f(\vx) \coloneqq -\ip{\vx,\xstar}$ such that the same conclusions as in \Cref{thm:large_m_lower} hold.
\end{corollary}

\subsection{Technical overview}\label{sec:tech-overview}

At a high level, our algorithms maintain a guess $\vx_t$ for the optimal solution $\xstar$. They will update this guess over many interactions with the adversary.

\paragraph{A general recipe.} We first describe the primitives that our methods depend on. Our first technical innovation is the notion of \textit{progress distributions}. Loosely speaking, these are distributions from which a learner is likely to sample a new guess $\vx_{t+1}$ that decreases its suboptimality. See \Cref{defn:prog}.

\begin{definition}[Progress Distribution]
\label{defn:prog}
Let $\myfunc{f}{\actionset}{\R}$ for $\actionset \subseteq \R^d$. For $\vx\in \actionset$ and $1 \le p < 2$, we say a distribution $\cD(\vx)$ over vectors in $\R^d$ is a \textit{$(p, \gamma,\rho)$-progress distribution for $\vx$} if we have the below.
\begin{align*}
    \prvv{\xplus \sim \cD(\vx)}{\frac{f(\vx)-f(\xplus)}{\inparen{f(\vx)-f(\xstar)}^p} \ge \frac{\rho}{d}} \ge \gamma.
\end{align*}
\end{definition}

So, if for every $\vx_t$ the learner had sample access to some progress distribution $\cD(\vx_t)$, the learner can significantly improve its solution (e.g. when $p = 1$, roughly $\sim d/\rho$ steps are sufficient for the learner to decrease its suboptimality by a constant factor). It is therefore natural that repeating such a sample-then-guess approach ad infinitum will yield an approximately optimal solution. In \Cref{thm:bin_search_general}, we prove this whenever there exist families of progress distributions for every range of possible suboptimalities. Thus, assuming the learner can maintain a (possibly quite pessimistic) estimate of its suboptimality over all the rounds, we obtain a template for proving the iteration complexities of \Cref{cor:bin_search_ip}, \ref{cor:bin_search_pl}, \Cref{cor:bin_search_dist}. Note that $\rho$ can be an arbitrarily small positive constant; even if there is a slim chance of decreasing the suboptimality, this is still sufficient because the monotone adversary ensures that the algorithm can never make negative progress.

\paragraph{Specifying progress distributions.} We now discuss how we instantiate the above template for the $\beta$-smooth (\Cref{defn:smoothness}) and $\alpha$-P\L{} (\Cref{defn:pl}) case (\Cref{cor:bin_search_pl}). We focus on \Cref{cor:bin_search_pl} for the sake of brevity; the proofs of \Cref{cor:bin_search_ip} \Cref{cor:bin_search_dist} require some additional care but at a high level follow a similar structure. At step $t$, the algorithm maintains a guess $\vx_t$ for the target $\xstar$. It chooses some step size $\eps_t$ and a random vector $\vg_t$ from $\eps_t \cdot \cN(0,\mI_d/d)$. We then query $\vx_t$ and $\vx_t - \vg_t$. The key observation is that with a constant probability, the angle between $\vg_t$ and the gradient $\nabla f(\vx)$ is small enough, so we make noticeable progress in such a step. We will use this to show that the distribution $\vx_t - \eps_t \cdot \cN(0,\mI_d/d)$ is a $(1, C_1, C_2)$-progress distribution (\Cref{defn:prog}) for constants $C_1, C_2$. Intuitively, this means that $\vx_t - \vg_t$ almost behaves like a step of gradient descent. To turn this observation into an algorithm, we need two main insights.

\paragraph{Step size schedule.} The principal difficulty of this approach is to choose the step size $\eps_t$. It is not immediately obvious how to do so since the algorithm does not observe any actual gradients or function values. Hence, if our step sizes are too large, the algorithm may overshoot the optimal solution $\xstar$ and therefore not actually improve the quality of its current solution $\vx_t$. On the other hand, if our step sizes are too small, the algorithm may not make enough progress in each step, which undesirably increases both the iteration complexity and the total cost. 

To address this, we carefully construct a step size schedule that relies on a pessimistic upper bound on the suboptimality of the algorithm's current solution. With this schedule, we show that in every step, one of two things happens -- either the step size $\eps_t$ is small enough such that there is the possibility of the algorithm decreasing the cost, or it is too large. For the first case, we use $\beta$-smoothness (\Cref{defn:smoothness}) to prove that there is a constant probability that the algorithm finds a descent direction, which decreases the cost of its current solution substantially. For the second case, we use the $\alpha$-P\L{} condition (\Cref{defn:pl}) to prove that the cost the algorithm incurs in such steps is low. After enough steps, we can show that either the second case always holds (i.e. that the suboptimality is already desirably small) or the maximum cost that the algorithm can pay per round is small. We then decrease the step size $\eps_t$ by a constant factor, update the suboptimality estimate accordingly, and infinitely recurse.

\paragraph{Bounding the failure probability over infinite rounds.} It now remains to show that the probability that the algorithm fails to make enough progress over \textit{infinitely many rounds} is small. This is where the distinction between the two goals of \Cref{problem:local_rec_fixed_set} becomes apparent. Specifically, even if we have a subroutine that, with high probability, outputs an $\eps$-approximate solution, this does not immediately convert to an algorithm that can achieve bounded cost over an infinite number of rounds -- note that the failure probabilities may accumulate in a divergent manner. Hence, we will require a more careful probabilistic analysis.

To overcome this challenge, we design the algorithm to run in phases $i = 1, 2, \dots$. In phase $i$, we use a step size $\eps_t$ proportional to $2^{-i/2}$ and run phase $i$ for $\sim id$ steps. Using the fact that the family of distributions we are using for sampling next steps are $(1, C_1, C_2)$-progress distributions, it will be enough to prove that $\sim d \cdot \nfrac{\beta}{\alpha}$ steps yield enough improving steps to decrease the suboptimality by a constant factor. We can therefore apply a Chernoff bound to conclude that the probability that the algorithm fails to make enough progress in phase $i$ is at most $\expv{-id \cdot \nfrac{\beta}{\alpha}}$. Finally, we apply a union bound that the total probability of failure by $\expv{-d \cdot \nfrac{\beta}{\alpha}} \le \expv{-d}$. 

To bound the total cost over all phases $i \in \N_{\ge 1}$, we note that the sum of the suboptimalities in each round is of the form $d\sum_{i \ge 1} i2^{-i} = O(d)$. The guarantee on the iteration complexity follows by noting that to achieve a suboptimality of $2^{-i}$, the algorithm runs $d\sum_{j \le i} i = O\inparen{i^2 \cdot d}$ iterations.

\subsection{Related works}

\paragraph{Dueling convex optimization.} As already mentioned, our formulation in \Cref{problem:local_rec_fixed_set} is an extension of dueling convex optimization in the noiseless setting~\cite{jamieson2012query,skm21,saha2022dueling}. \citet{jamieson2012query} employ a coordinate-descent algorithm to show for $\alpha$-smooth and $\beta$-strongly convex $f$,  $\tilde{O}(\nfrac{d\beta}{\alpha} \logv{\nfrac{1}{\eps}})$ queries suffice to learn an $\eps$-optimal point. As mentioned earlier, it is not clear how to adapt their algorithm to handle monotone feedback. In addition, the works \cite{skm21,saha2022dueling} show results for more general classes of $f$ and in the presence of noise (where the adversary can return invalid response with nonzero probability). However, their algorithms explicitly rely on sign feedback $f(\xt{1}) \overset{?}{>} f(\xt{2})$ to construct gradient estimators, which are not possible in the monotone adversary setting.

\paragraph{Monotone adversaries.} Our setting is an example of learning with a monotone adversary, where an adversary can choose to improve the feedback or information the algorithm gets. A common characteristic is that the improved information may paradoxically break or harm the performance of a given algorithm that works with non-improved information. Monotone adversaries are often studied in the semi-random model literature \cite{blum1995coloring, feige_2021, moitra_2021} for statistical estimation problems \cite{cheng2018non,moitra_2021,kelner2022semi} as well as learning problems, i.e., linear classification with Massart noise \cite{massart2006risk, diakonikolas2019distribution}. 

\paragraph{Preference-based feedback.} Our formulation in this chapter falls within the growing body of literature that tackles learning with preference-based feedback, where the algorithm does not learn \emph{how good} its options were in an absolute sense, just which one(s) were better than others. 

Other natural problems with preference-based feedback are contextual search \cite{ls18, lobel2018multidimensional, lls20}, contextual recommendation (also called contextual inverse optimization) \cite{bfl21, gollapudi2021contextual}, and $1$-bit matrix completion \cite{davenport20141}.

\section{Proofs of upper bound results}

In this section, we prove \Cref{thm:bin_search_general}. The point of \Cref{thm:bin_search_general} is to construct and analyze a meta-algorithm for \Cref{problem:local_rec_fixed_set} when the algorithm can sample next steps from progress distributions (\Cref{defn:prog})). We then show how to use this framework to prove  \Cref{cor:bin_search_ip} (results for $f(\vx) = \ip{-\vx, \xstar}$), \Cref{cor:bin_search_pl} (results for $f(\vx)$ being $\beta$-smooth and $\alpha$-P\L{}), and \Cref{cor:bin_search_dist} (results for $f(\vx) = \norm{\vx - \xstar}_2$), in that order. It will be helpful to recall the overview from \Cref{sec:tech-overview} throughout this section.

We prove \Cref{thm:bin_search_general} in \Cref{subsec:general_proof}, \Cref{cor:bin_search_ip} in \Cref{subsec:ip_proof}, \Cref{cor:bin_search_pl} in \Cref{subsec:pl_proof}, and \Cref{cor:bin_search_dist} in \Cref{subsec:dist_proof}.

Before we jump into the main proofs, we prove some straightforward numerical inequalities that we need later.

\begin{lemma}
\label{lemma:cost_convergence_series}
For $r \in (0, 1)$ and $1 \le p < 2$, we have $\sum_{i \ge 0} i \cdot r^{(1-p/2)i} \le \frac{r^{p/2+1}}{(r-r^{p/2})^2}$.
\end{lemma}
\begin{proof}[Proof of \Cref{lemma:cost_convergence_series}]
Recall that
\begin{align*}
    \sum_{i \ge 0} r^{(1-p/2)i} = \frac{1}{1-r^{1-p/2}}.
\end{align*}
Taking the derivative of both sides with respect to $r$ yields
\begin{align*}
    \sum_{i \ge 0} \inparen{1-\nfrac{p}{2}}i\cdot r^{(1-p/2)i - 1} = \frac{(2-p)r^{p/2}}{2(r-r^{p/2})^2}.
\end{align*}
We multiply both sides by $r$ and divide both sides by $1-\nfrac{p}{2}$; we conclude that
\begin{align*}
    \sum_{i \ge 0} i \cdot r^{(1-p/2)i} = \frac{r^{p/2+1}}{(r-r^{p/2})^2}
\end{align*}
which recovers the statement of \Cref{lemma:cost_convergence_series}.
\end{proof}

\begin{lemma}[Inner product with a random vector]\label{lem:random-progress}
Let $\vg \sim \mathsf{Unif}(\S_2^{d-1})$ and let $\vy \in \S_2^{d-1}$ be fixed. Then
\begin{align*}
    \prvv{\vg}{\ip{\vg,\vy}  \ge \frac{1}{2\sqrt{d}}} \ge \frac{1}{8}.
\end{align*}
\end{lemma}
\begin{proof}
By rotational invariance, without loss of generality, we can let $\vy = \ve_1$. We apply Lemma 2.2 (a) due to \citet{dasgupta} with $\beta = \nfrac{1}{4}$ to conclude that
\begin{align*}
    \prvv{\vg}{\vg_1^2 \le \frac{1}{4d}} \le \expv{\frac{1}{2}\inparen{1-\frac{1}{4}+\ln\inparen{\frac{1}{4}}}} < \frac{3}{4}
\end{align*}
which means that $ \prvv{\vg}{\abs{\ip{\vg,\vy}}  \ge \frac{1}{2\sqrt{d}}} \ge \frac{1}{4}$. The result of \Cref{lem:random-progress} follows by symmetry.
\end{proof}

\begin{lemma}[Best inner product among $m$ Gaussians]
\label{lem:m-random-progress}
Let $m \ge 3$. Let $\vg_1,\dots,\vg_m \sim \cN\inparen{0,\frac{\mI_d}{d}}$ and let $\vy \in \S_2^{d-1}$ be fixed. Then
\begin{align*}
    \prvv{\vg_1,\dots,\vg_m}{\max_{1 \le i \le m} \ip{\vg_i,\vy}  \ge \frac{\sqrt{\log m}}{\sqrt{d}}} > \frac{1}{20}.
\end{align*}
\end{lemma}
\begin{proof}[Proof of \Cref{lem:m-random-progress}]
As in the proof of \Cref{lem:random-progress}, without loss of generality (from rotational invariance), let $\vy=\ve_1$. Thus, we want to understand $\max_{1 \le i \le m} \vg_i[1]$, which can be rewritten as $\max_{1 \le i \le m} g_i$ for $g_1,\dots,g_m \sim \cN(0,1)$. We start with the well-known fact (see \cite[Proposition 2.1.2]{vershynin_2018}) that for all $t \in \R$,
\begin{align*}
    \prvv{g \sim \cN(0,1)}{g \le t} \le 1 - \frac{1}{\sqrt{2\pi}}\frac{1}{t}\inparen{1 - \frac{1}{t^2}}\expv{-\frac{t^2}{2}} \le \expv{- \frac{1}{\sqrt{2\pi}}\frac{1}{t}\inparen{1 - \frac{1}{t^2}}\expv{-\frac{t^2}{2}}}.
\end{align*}
From independence, we get
\begin{align*}
    \prvv{g_1,\dots,g_m \sim \cN(0,1)}{\text{for all } i, g_i \le t} &\le \expv{-m\frac{1}{\sqrt{2\pi}}\frac{1}{t}\inparen{1 - \frac{1}{t^2}}\expv{-\frac{t^2}{2}}}.
\end{align*}
Choose $t = c\sqrt{\log m}$ and write
\begin{align*}
    \prvv{g_1,\dots,g_m \sim \cN(0,1)}{\text{for all } i, g_i \le c\sqrt{\log m}} &\le \expv{-\frac{m}{\sqrt{2\pi}}\frac{1}{c\sqrt{\log m}}\inparen{1 - \frac{1}{c^2\log m}}\expv{-\frac{c^2\log m}{2}}} \\
    &= \expv{-\frac{m^{1-c^2/2}}{\sqrt{2\pi}}\frac{1}{c\sqrt{\log m}}\inparen{1 - \frac{1}{c^2\log m}}}.
\end{align*}
We can check that for $c=1$ and over the integers $m \ge 3$, the RHS is decreasing and at $m=3$, the RHS is $< 0.95$. Rearranging, we get
\begin{align*}
    \prvv{g_1,\dots,g_m \sim \cN(0,1)}{\max_{1 \le i \le m} g_i \ge \sqrt{\log m}} > \frac{1}{20},
\end{align*}
thereby completing the proof of \Cref{lem:m-random-progress} (after appropriately rescaling).
\end{proof}

\subsection{A general algorithm for \Cref{problem:local_rec_fixed_set} with progress distributions}
\label{subsec:general_proof}

The goal of this subsection is to develop the general tools we need to prove our main results.

The key primitive of our analysis is a general algorithm (\Cref{alg:dueling_main_alg}) that solves \Cref{problem:local_rec_fixed_set} when we are given certain convenient distributions from which we sample new guesses. We call these \textit{progress distributions}; recall \Cref{defn:prog}.

Let us describe \Cref{alg:dueling_main_alg}. In each step, \Cref{alg:dueling_main_alg} maintains a current guess $\vx_t$ and chooses a slight perturbation of that guess $\xplus_t \sim \cD(\vx_t)$, where $\cD(\vx_t)$ is a $(p, \gamma, \rho)$-progress distribution (\Cref{defn:prog}). \Cref{alg:dueling_main_alg} then submits the pair of guesses $\inbraces{\vx_t, \xplus_t}$. To analyze \Cref{alg:dueling_main_alg}, the main observation is that with probability $\ge \gamma$, the point $\xplus_t$ substantially improves over the cost of $\vx_t$ -- this follows directly from \Cref{defn:prog}. We exploit this intuition to give our most general result (\Cref{thm:bin_search_general}) and to prove the correctness of \Cref{alg:dueling_main_alg}.

\begin{mainthm}
\label{thm:bin_search_general}
Let $\myfunc{f}{\actionset}{\R}$. Let $B$ and $\vx_1 \in \actionset$ be such that $f(\vx_1)-f(\xstar) \le B$. For $C > 0$, constant $r \in (0,0.99)$, and for all $i \in \N_{\ge 1}$, suppose there exists intervals of the form $C \cdot \insquare{r^{i+1}, r^{i}}$ such that their union covers the interval $[0, B]$.

If there exists a $(p, \gamma,\rho)$-progress distribution $\cD_i(\vx)$ whenever $f(\vx)-f(\xstar) \in C \cdot \insquare{r^{i+1}, r^{i}}$ for all $i \ge 1$ and where $p, \gamma, \rho$ do not depend on $\vx$ and $i$, then there is an algorithm (\Cref{alg:dueling_main_alg}) for \Cref{problem:local_rec_fixed_set} that, with probability at least $1 - \expv{-O\inparen{\frac{d}{\rho B^{p-1}}}}$, incurs total cost
\begin{align*}
    O\inparen{\frac{B\logv{\nfrac{1}{r}}}{B^{p-1}\gamma\rho\min\inbraces{r^{p\inparen{\nfrac{p-1}{2-p}}}, \inparen{r-r^{p/2}}^2}} \cdot d}.
\end{align*}
Additionally, \Cref{alg:dueling_main_alg} finds a point $\vx$ satisfying $f(\vx)-f(\xstar) \le \eps$ in 
\begin{align*}
    O\inparen{\frac{1}{B^{p-1}\gamma\rho} \cdot d \cdot \logv{\frac{B}{\eps}}^2}
\end{align*}
iterations with at least the aforementioned probability.
\end{mainthm}

\begin{algorithm}
\caption{General recipe algorithm for dueling optimization}\label{alg:dueling_main_alg}
\begin{algorithmic}[1]
    \State \textbf{Input}: Interaction with a monotone adversary $\cM$ as defined in \Cref{problem:local_rec_fixed_set}; initial point $\vx_1$ and bound $B$ satisfying $f(\vx_1)-f(\xstar) \le B$; values $C$ and $r$ for which there exist corresponding intervals and $(p,\gamma, \rho)$-progress distribution families $\cD_i$ (see the statement of \Cref{thm:bin_search_general}).
    \State Initialize $\vx_1 = 0$, $t=1$.
    \For{$i = 1, \dots$}
        \For{$T(i) \coloneqq \nfrac{2i}{(\gamma\min(1,\rho))} \cdot (Cr^{i+1})^{-(p-1)} \cdot  \logv{\nfrac{1}{r}} \cdot d$ iterations}\label{line:mary_main_alg_ti}
            \State Sample $\xplus_t$ from $\cD_i(\vx_t)$.
            \State \textbf{Submit guesses} $\inbraces{\vx_t, \xplus_t}$ and  \textbf{receive response} $\xstar_t$.
            \State Let $\vx_{t+1} = \xstar_t$.
            \State Update $t \gets t+1$.
        \EndFor
    \EndFor
\end{algorithmic}
\end{algorithm}

The proof of \Cref{thm:bin_search_general} has two main parts. In the first part, we will prove that for each value of $i$ (call the set of timesteps belonging to a particular value of $i$ ``phase $i$''), the number of steps $T(i)$ is sufficient to ensure that the cost of the algorithm's solution decays gracefully with sufficiently large probability. In the second part, we will prove that the total cost the algorithm pays over all phases $i \ge 1$ is $\sim B \cdot \nfrac{d}{\gamma\rho}$ as promised. \Cref{thm:bin_search_general} will easily follow by combining these facts.

We start with stating and proving \Cref{lemma:steps_per_phase}.

\begin{lemma}
\label{lemma:steps_per_phase}
Let $i \ge 1$ and let $T(i)$ be defined below (or see \Cref{line:mary_main_alg_ti} of \Cref{alg:dueling_main_alg}).
\begin{align*}
    T(i) \coloneqq \frac{2i}{\gamma\min(1,\rho)\cdot (Cr^{i+1})^{p-1}}  \cdot  \logv{\nfrac{1}{r}} \cdot d.
\end{align*}
Let $t_i$ be the first iteration of phase $i$. If $f(\vx_{t_i}) - f(\xstar) \le Cr^i$, then, with probability $\ge 1 - r^{\frac{dB^{-(p-1)}}{4\rho} \cdot i}$, we have $f(\vx_{t_i+T(i)+1}) - f(\xstar) \le Cr^{i+1}$.
\end{lemma}
\begin{proof}[Proof of \Cref{lemma:steps_per_phase}]
Assume that we have $f(\vx_{t_i})-f(\xstar) \in C \cdot \insquare{r^{i+1}, r^i}$ (otherwise, we are done immediately).

Define the indicator random variable $Y_t$ as follows.
\begin{align*}
    Y_t \coloneqq \indicator{\frac{f(\vx_{t})-f(\xplus_t)}{\inparen{f(\vx_{t})-f(\xstar)}^p} \ge \frac{\rho}{d}}.
\end{align*}

Consider the distribution of guesses $\cD_i$ (let us omit the argument $\vx_t$ for the sake of brevity). Since $\cD_i$ is a $(p,\gamma,\rho)$-progress distribution, we have
\begin{align*}
    \prvv{\xplus \sim \cD_i}{\frac{f(\vx_t)-f(\xplus_t)}{f(\vx_t)-f(\xstar)} \ge \frac{\rho}{d} \cdot \inparen{Cr^{i+1}}^{p-1}} \ge \prvv{\xplus \sim \cD_i}{Y_t = 1}  \ge \gamma.
\end{align*}
Call every step $t$ for which $Y_t = 1$ a ``successful step.'' Let us give a high-probability count on the number of successful steps. Recall that a form of the Chernoff bound states that, for $\delta \in [0,1]$ and independent indicator random variables $Y_j$,
\begin{align*}
    \prv{\sum_{j = t_i}^{t_i+T(i)} Y_j \le (1-\delta)\exv{\sum_{j=t_i}^{t_i+T(i)} Y_j}} \le \expv{-\frac{\delta^2 \cdot \exv{\sum_{j=t_i}^{t_i+T(i)} Y_j}}{2}}.
\end{align*}
Applying the Chernoff bound with $\delta = \nfrac{1}{2}$ yields
\begin{align*}
    \prv{\sum_{j=t_i}^{t_i+T(i)} Y_j \le \frac{T(i)\gamma}{2}} \le \expv{-\frac{i \cdot \frac{d}{\rho\inparen{Cr^{i+1}}^{p-1}} \cdot \logv{\nfrac{1}{r}}}{4}} \le r^{\frac{dB^{-(p-1)}}{4\rho} \cdot i}
\end{align*}
where we use $Cr^{i+1} \le Cr^2 \le B$.

It remains to show that after at least $\nfrac{T(i)\gamma}{2}$ successful steps, we have $f(\vx_{t_i + T(i) + 1}) - f(\xstar) \le Cr^{i+1}$. Recall that we assume that $f(\vx_{t_i}) - f(\xstar) \ge Cr^{i+1}$ and note that for every successful step, we have
\begin{align*}
    \frac{f(\vx_t)-f(\xplus_t)}{f(\vx_t)-f(\xstar)} \ge \frac{\rho}{d} \cdot \inparen{Cr^{i+1}}^{p-1}
\end{align*}
which implies
\begin{align*}
    \frac{f(\xplus_t)-f(\xstar)}{f(\vx_t)-f(\xstar)} \le 1 - \frac{\rho}{d} \cdot \inparen{Cr^{i+1}}^{p-1}.
\end{align*}
We multiply over all steps in phase $i$, giving
\begin{align*}
    \frac{f(\vx_{t_i + T(i) + 1})-f(\xstar)}{f(\vx_{t_i})-f(\xstar)} &= \prod_{t = t_i}^{t_i + T(i)} \frac{f(\vx_{t+1}) - f(\xstar)}{f(\vx_t) - f(\xstar)} \le \inparen{1 - \frac{\rho}{d} \cdot \inparen{Cr^{i+1}}^{p-1}}^{T(i)\gamma/2} \\
    &\le \inparen{1 - \frac{2i}{\gamma} \cdot \frac{\logv{\nfrac{1}{r}}}{T(i)}}^{T(i)\gamma/2} \le \expv{\frac{2i}{\gamma} \cdot \frac{\logv{\nfrac{1}{r}}}{T(i)} \cdot \frac{T(i)\gamma}{2}} = r^{i} \le r.
\end{align*}
Finally, recall that $f(\vx_{t_i}) - f(\xstar) \le Cr^i$. Combining this with the above gives $f(\vx_{t_i + T(i) + 1}) - f(\xstar) \le Cr^{i+1}$, concluding the proof of \Cref{lemma:steps_per_phase}.
\end{proof}

Next, we have \Cref{lemma:total_cost}, which controls the total cost that \Cref{alg:dueling_main_alg} incurs assuming that the cost is sufficiently low in each phase.

\begin{lemma}
\label{lemma:total_cost}
For a timestep $t$, let $i(t)$ be the phase that $t$ belongs to.

If for all $t$ we have $f(\vx_t) - f(\xstar) \le Cr^{i(t)}$, then \Cref{alg:dueling_main_alg} incurs total cost
\begin{align*}
    O\inparen{\frac{B\logv{\nfrac{1}{r}}}{C^{p-1}\gamma\rho\min\inbraces{r^{p\inparen{\nfrac{p-1}{2-p}}}, \inparen{r-r^{p/2}}^2}} \cdot d}.
\end{align*}
\end{lemma}
\begin{proof}[Proof of \Cref{lemma:total_cost}]
Recall throughout this proof that $r \le 0.99$ and $p$ is a constant such that $p < 2$.

Observe that in phase $i$, the algorithm incurs cost at most
\begin{align*}
    T(i) \cdot Cr^i = \frac{2i}{\gamma} \cdot \frac{dCr^i}{\rho\inparen{Cr^{i+1}}^{p-1}} \cdot \logv{\nfrac{1}{r}} = \frac{2i}{\gamma} \cdot \frac{d\logv{\nfrac{1}{r}}}{\rho C^{p-2}} \cdot r^{i - (i+1)(p-1)}.
\end{align*}
We will find a threshold $i_p$ for which for all $i \ge i_p$, the above cost is exponentially decaying. This will allow us to control the sum of the costs over infinitely many rounds. We choose $i_p = 2 \cdot \ceil{\inparen{\nfrac{(p-1)}{(2-p)}}}$. Notice that for all $i \ge i_p$, the exponent on $r$ can be bounded as
\begin{align*}
    i-(i+1)(p-1) = i(2-p) - (p-1) \ge \inparen{1-\frac{p}{2}}i.
\end{align*}
Note that this also implies that $(i_p+1)(p-1) \le \nfrac{p}{2} \cdot i_p = p\ceil{\nfrac{(p-1)}{(2-p)}}$.

To total the cost, we consider two cases. First, suppose $1 \le i \le i_p - 1$. Observe that in each of these phases, we pay cost at most $B$, so we have
\begin{align*}
    \sum_{i=1}^{i_p-1} B\cdot T(i) &\le B \inparen{T_{i_p} \cdot i_p} = 2B\inparen{2 \cdot \frac{p-1}{2-p}}^2 \cdot \frac{\logv{\nfrac{1}{r}}}{\gamma} \cdot \frac{d}{\rho C^{p-1}r^{(i_p+1)(p-1)}} \\
    &\le 2B\inparen{2 \cdot \frac{p-1}{2-p}}^2 \cdot \frac{\logv{\nfrac{1}{r}}}{\gamma} \cdot \frac{d}{\rho C^{p-1}r^{p\inparen{\nfrac{p-1}{2-p}}}}. \numberthis\label{eq:bound1}
\end{align*}
Next, we sum over all phases $i \ge i_p$. We obtain a cost that is at most
\begin{align*}
    \sum_{i \ge i_p} \frac{2i}{\gamma} \cdot \frac{d\logv{\nfrac{1}{r}}}{\rho C^{p-2}} \cdot r^{i - (i+1)(p-1)} \le \frac{2d\logv{\nfrac{1}{r}}}{\gamma\rho \cdot C^{p-2}}\sum_{i \ge i_p} i \cdot r^{(1-p/2)i} \le \frac{2d\logv{\nfrac{1}{r}}}{\gamma\rho \cdot C^{p-2}} \cdot \frac{r^{p/2+1}}{(r-r^{p/2})^2} \numberthis\label{eq:bound2}
\end{align*}
where the last inequality follows from \Cref{lemma:cost_convergence_series}. Combining \eqref{eq:bound1} and \eqref{eq:bound2} yields
\begin{align*}
    &\quad 2B\inparen{2 \cdot \frac{p-1}{2-p}}^2 \cdot \frac{\logv{\nfrac{1}{r}}}{\gamma} \cdot \frac{d}{\rho C^{p-1}r^{p\inparen{\nfrac{p-1}{2-p}}}}+\frac{2d\logv{\nfrac{1}{r}}}{\gamma\rho \cdot C^{p-2}} \cdot \frac{r^{p/2+1}}{(r-r^{p/2})^2} \\
    &\le 2B\inparen{2 \cdot \frac{p-1}{2-p}}^2 \cdot \frac{\logv{\nfrac{1}{r}}}{\gamma} \cdot \frac{d}{\rho C^{p-1}r^{p\inparen{\nfrac{p-1}{2-p}}}}+\frac{2d\logv{\nfrac{1}{r}}}{\gamma\rho \cdot C^{p-1}} \cdot \frac{B}{(r-r^{p/2})^2} \\
    &= O\inparen{\frac{B\logv{\nfrac{1}{r}}}{C^{p-1}\gamma\rho\min\inbraces{r^{p\inparen{\nfrac{p-1}{2-p}}}, \inparen{r-r^{p/2}}^2}} \cdot d}
\end{align*}
This concludes the proof of \Cref{lemma:total_cost}.
\end{proof}

We are now ready to prove \Cref{thm:bin_search_general}.

\begin{proof}[Proof of \Cref{thm:bin_search_general}]
It is sufficient to prove that with probability $\ge 1 - \expv{-O\inparen{\frac{d}{\rho B^{p-1}}}}$, at the end of phase $i$, we have $f(\vx_t) - f(\xstar) \le Cr^{i+1}$. Recall the conclusion of \Cref{lemma:steps_per_phase} and that $f(\vx_1) - f(\xstar) \le B \le Cr$; by a union bound, we have for all phases $i$ that $f(\vx_t) - f(\xstar) \le Cr^{i+1}$ with probability
\begin{align*}
    1 - \sum_{i \ge 1} r^{\frac{dB^{-(p-1)}}{4\rho} \cdot i} \ge 1 - \expv{-O\inparen{\frac{d}{\rho B^{p-1}}}}
\end{align*}
where we use $0 < r < 0.99$. The first part of \Cref{thm:bin_search_general} now follows directly from applying \Cref{lemma:total_cost}. The rest of the statement of \Cref{thm:bin_search_general} follows by noting that
\begin{align*}
    \sum_{j \le i} T(j) = \frac{2C^{-(p-1)}\logv{\nfrac{1}{r}}}{\gamma\min(1,\rho)} \cdot d \cdot \sum_{j \le i} j\inparen{r^{-(j+1)(p-1)}} \lesssim \frac{2C^{-(p-1)}}{\gamma\min(1,\rho)} \cdot d \cdot i^2.
\end{align*}
where we again use $r < 0.99$. We set $\eps = Cr^i$ and conclude.
\end{proof}

\subsection{Proof of \Cref{cor:bin_search_ip}}
\label{subsec:ip_proof}

The goal of this subsection is to prove \Cref{cor:bin_search_ip}.

Our plan will be to use the general guarantee of \Cref{thm:bin_search_general}. Thus, the main task is to prove that there is an appropriate interval cover and corresponding sequence $\cD_i(\vx)$ of progress distributions for all $\vx$ belonging to phase $i$ that satisfy the conditions of \Cref{thm:bin_search_general}.

We prove this fact in \Cref{lem:inner-prod-dists}. We remark that we made no effort to optimize the numerical constants; we choose the constants that appear in the Lemma statement to simplify calculations, as these will not impact our asymptotic results.

\begin{lemma}
\label{lem:inner-prod-dists}
Let $\myfunc{f}{\R^d}{\R}$ be the negative inner product function defined on $\mathbb{S}_2^{d-1}$ with respect to some unknown target $\xstar$. Then for any $\vx$ for which $f(\vx)-f(\xstar) \in \insquare{10^{-(i+1)}, 10^{-i}}$ and for which $\langle \xstar,\vx\rangle>0$, there is a $(1.5,10^{-1},10^{-4})$-progress distribution (\Cref{defn:prog}) that can be computed in time $O(d)$.
\end{lemma}

\begin{proof}[Proof of \Cref{lem:inner-prod-dists}]
We explain the construction of the distribution $\cD(\vx)$.

Impose the following coordinates on $\R^d$. Let the first coordinate $x_1$ be the direction of $\vx$, and the remaining $d-1$ coordinates be an arbitrary coordinate system for the perpendicular directions. Then, $\vx$ has coordinates $(1,0\ldots)$. Next, let $z \coloneqq 10^{-(i+1)}$ and $s:=\frac{z}{10\sqrt{d-1}}$. Let $\vs$ be a point randomly drawn from a $d-1$ dimensional sphere of radius $s$ whose coordinates are denoted $s_1\ldots s_{d-1}$. Then, the distribution $\cD\inparen{\vx}$ is the distribution of $(\sqrt{1-s^2},s_1\ldots s_{d-1})$. It is easy to verify that these points lie on $\mathbb{S}_2^{d-1}$.

This distribution can be computed in time $O(d)$. We will now show that it is a $(1.5,10^{-1},10^{-4})$ progress distribution.

Let $z_0 \coloneqq f(\vx)-f(\xstar)$; recall that $z_0 \in [0,1]$ and $z \le z_0 \le 10z$. We write the target vector $\xstar=(1-z_0)\vx+\sqrt{1-(1-z_0)^2}\vy = (1-z_0)\vx+\sqrt{2z_0-z_0^2}\vy$ where $\vy$ is a unit vector and $\ip{\vy,\vx}=0$. Note that this expression holds because $\langle \xstar, \vx\rangle = 1-z_0$.

Let $\xplus$ be a random point chosen from $\cD(\vx)$. Let $y_1\ldots y_{d-1}$ be the coordinates of $\vy$ in the $d-1$ dimensional coordinate plane perpendicular to $\vx$ defined above. We compute
\begin{align} \label{eq:xxplus}
    \langle \xstar, \xplus \rangle
    =&~\inparen{1-z_0,\sqrt{2z_0-z_0^2}y_1\ldots \sqrt{2z_0-z_0^2}y_{d-1}}\cdot \inparen{\sqrt{1-s^2},s_1\ldots s_{d-1}} \\
    =&~(1-z_0)\sqrt{1-s^2}+\inparen{\sqrt{2z_0-z_0^2}y_1\ldots \sqrt{2z_0-z_0^2}y_{d-1}} \cdot (s_1\ldots s_{d-1}).
\end{align}
By \Cref{lem:random-progress}, we have (note that we weaken the constants from \Cref{lem:random-progress} for numerical convenience later in the proof)
\begin{align*}
    \prvv{\vs}{\inparen{\sqrt{2z_0-z_0^2}y_1\ldots \sqrt{2z_0-z_0^2}y_{d-1}} \cdot (s_1\ldots s_{d-1}) \geq \frac{0.1}{\sqrt{d-1}}\cdot \sqrt{2z_0-z_0^2} \cdot s} \ge 0.1.
\end{align*}
Because $2z_0-z_0^2\geq z_0\geq z$, we have
\begin{align*}
     \frac{0.1}{\sqrt{d-1}}\cdot \sqrt{2z_0-z_0^2} \cdot s \geq \frac{s\sqrt{z}}{10\sqrt{d-1}}.
\end{align*}
In turn, this shows
\begin{align*}
    \prvv{s}{\inparen{\sqrt{2z_0-z_0^2}y_1\ldots \sqrt{2z_0-z_0^2}y_{d-1}} \cdot (s_1\ldots s_{d-1}) \geq \frac{s\sqrt{z}}{10\sqrt{d-1}}} \ge 0.1.
\end{align*}
Combining this with Equation~\ref{eq:xxplus}, we obtain
\begin{align} \label{eq:xstarxplus}
    \prvv{\xplus \sim \cD_{\vx}}{\langle \xstar, \xplus \rangle \ge (1-z_0)\sqrt{1-s^2}+\frac{s\sqrt{z}}{10\sqrt{d-1}}} \ge 0.1.
\end{align}
Now, we will find a lower bound for $(1-z_0)\sqrt{1-s^2}+\frac{s\sqrt{z}}{10\sqrt{d-1}}$. We have
\begin{align*}
    (1-z_0)\sqrt{1-s^2}+\frac{s\sqrt{z}}{10\sqrt{d-1}} &\geq (1-z_0)(1-s^2)+\frac{s\sqrt{z}}{10\sqrt{d-1}} \\
    &= (1-z_0)(-s^2)+\frac{s\sqrt{z}}{10\sqrt{d-1}} + (1-z_0) \\
    &\geq (1-z)(-s^2)+\frac{s\sqrt{z}}{10\sqrt{d-1}} + (1-z_0).
\end{align*}
Now, using that $s = \frac{z}{10\sqrt{d-1}}$, we get 
\begin{align*}
    (1-z)(-s^2)+\frac{s\sqrt{z}}{10\sqrt{d-1}} + (1-z_0)
    &= (1-10s\sqrt{d-1})(-s^2)+\frac{s^{3/2}}{\sqrt{10}(d-1)^{1/4}} + (1-z_0)\\
    &= -s^2+10s^3\sqrt{d-1}+\frac{s^{3/2}}{\sqrt{10}(d-1)^{1/4}} + (1-z_0) \\
    &\geq \left(10s^3\sqrt{d-1}+\frac{s^{3/2}}{8(d-1)^{1/4}}-s^2\right) \\
    &\quad\quad+\frac{s^{3/2}}{6(d-1)^{1/4}} + (1-z_0)
\end{align*}
where the last line follows from $\nfrac{1}{\sqrt{10}} > \nfrac{1}{8}+\nfrac{1}{6}$. Finally, applying weighted AM-GM lets us see
\begin{align*}
    10s^3\sqrt{d-1}+\frac{s^{3/2}}{8(d-1)^{1/4}}\geq \frac{3}{2^{2/3}}\inparen{10s^3\sqrt{d-1}}^{1/3}\inparen{\frac{s^{3/2}}{8(d-1)^{1/4}}}^{2/3} = \frac{3\cdot 10^{1/3}}{2^{2/3}2^{2}}s^2 > s^2
\end{align*}
where we use a weight of $\nfrac{1}{3}$ on the first term and a weight of $\nfrac{2}{3}$ on the second term. We now write
\begin{align*}
    (1-z_0)\sqrt{1-s^2}+\frac{s\sqrt{z}}{10\sqrt{d-1}} \geq \frac{s^{3/2}}{6(d-1)^{1/4}} + (1-z_0).
\end{align*}
Substituting $s$ once again and recalling that $\langle \xstar, \vx\rangle = (1-z_0)$ and $z \ge \frac{z_0}{10} = \frac{\langle \xstar, \xstar - \vx\rangle}{10}$, we get
\begin{align*}
    (1-z_0)\sqrt{1-s^2}+\frac{s\sqrt{z}}{10\sqrt{d-1}} \geq \frac{z^{3/2}}{6\cdot 10^{3/2} \cdot d} + \langle \xstar, \vx \rangle > \frac{10^{-4}\langle \xstar, \xstar - \vx\rangle^{3/2}}{d} + \langle \xstar, \vx \rangle.
\end{align*}
Combining this with \eqref{eq:xstarxplus}, we now have
\begin{align*}
    &\prvv{\xplus \sim \cD\inparen{\vx}}{\langle \xstar, \xplus \rangle \ge \frac{10^{-4}\langle \xstar, \xstar - \vx\rangle^{3/2}}{d} + \langle \xstar, \vx\rangle} \ge 0.1
\end{align*}
which means that
\begin{align*}
    \prvv{\xplus \sim \cD\inparen{\vx}}{\langle \xstar, \xplus - \vx \rangle \ge \frac{10^{-4}\langle \xstar, \xstar - \vx\rangle^{3/2}}{d}} \ge 0.1.
\end{align*}
This exactly aligns with the definition of a $(1.5, 10^{-1}, 10^{-4})$ progress distribution, completing the proof of \Cref{lem:inner-prod-dists}.
\end{proof}

We will now conclude \Cref{cor:bin_search_ip} using \Cref{thm:bin_search_general}.

\begin{proof}[Proof of \Cref{cor:bin_search_ip}]
To apply \Cref{thm:bin_search_general}, we need to present $C$, $r$, and a sequence of parameterizations $\cD_i$ that satisfy the premises.

Set $r=0.1$ and $C=10$. By \Cref{lem:inner-prod-dists}, we can find progress distributions for each interval $[Cr^i,Cr^{i-1}]$ of suboptimality of the current function value, since we can find such progress distributions as long as the suboptimality of the function is at most $1$. 

Note that the algorithm can begin with a point $\vx$ where $f(\vx)-f(\xstar)<1$ by first querying two opposite points on a sphere; one can easily see that at least one of the two points queried satisfies $f(\vx)-f(\xstar)<1$.

We therefore conclude the proof of \Cref{cor:bin_search_ip}.
\end{proof}

\subsection{Proof of \Cref{cor:bin_search_pl}}
\label{subsec:pl_proof}

The goal of this subsection is to prove \Cref{cor:bin_search_pl}.

As before, we use the general guarantee of \Cref{thm:bin_search_general} via proving that there is an appropriate interval cover and corresponding sequence $\cD_i(\vx)$ of progress distributions for all $\vx$ for which $f(\vx)-f(\xstar) \in \insquare{Cr^{i+1},Cr^i}$. We prove this fact in \Cref{lemma:pl_prog_distro}.

\begin{lemma}
\label{lemma:pl_prog_distro}
Fix $i \in \N_{\ge 1}$. Let $\eps_i = \sqrt{\nfrac{2B \alpha}{\beta^2}} \cdot \nfrac{\sqrt{\log m}}{\sqrt{d}}\cdot 2^{-i/2-1}$. Suppose $f$ is $\beta$-smooth and $\alpha$-P\L{}, and suppose we have $f(\vx) - f(\xstar) \in B \cdot \insquare{2^{-i}, 2^{-i+1}}$. Let $\vg_1,\dots,\vg_m \sim \cN(0, \mI_d/d)$ independently. Consider the points $\vx+\eps_i\vg_j$ for $1 \le j \le m$ and let $k$ be the index that minimizes $f(\vx+\eps_i\vg_j)$. Let $\cD_i(\vx) = \vx+\eps_i\vg_k$. Then, for all $d$ larger than a universal constant, $\cD_i(\vx)$ is a $(1,\gamma,\rho)$-progress distribution for $(\gamma,\rho) = \inparen{\nfrac{1}{40}, \nfrac{\alpha}{8\beta} \cdot \log m}$.
\end{lemma}
\begin{proof}[Proof of \Cref{lemma:pl_prog_distro}]
Note that although the algorithm may not be able to explicitly evaluate $\cD_i(\vx)$, the adversary always returns a point that has function value that is at least as good as the return value of $\cD_i(\vx)$, which can be done by guessing all $m$ points $\vx+\eps_i\vg_j$ for $j =1,\dots,m$ in round $t$.

Let $\vg \coloneqq \vx - \xplus$. It is sufficient to consider the case where we have $\eps_i \le \frac{1}{2\beta} \cdot \frac{\norm{\nabla f(\vx)}_2\sqrt{\log m}}{\sqrt{d}}$. To see this, suppose this is not the case. We apply the P\L{} inequality and write
\begin{align*}
    f(\vx) - f(\xstar) \le \frac{1}{2\alpha}\norm{\nabla f(\vx)}_2^2 \le \frac{d}{\log m} \cdot \frac{2\beta^2}{\alpha}\eps_i^2 = \frac{d}{\log m} \cdot \frac{2\beta^2}{\alpha} \inparen{ \sqrt{\frac{2B \alpha\log m}{\beta^2 d}} \cdot \frac{1}{2^{i/2+1}}}^2 = \frac{B}{2^{i}}
\end{align*}
which implies that the suboptimality $f(\vx) - f(\xstar)$ does not belong to the range we are considering.  

Next, we use \Cref{lem:m-random-progress} to write the below.
\begin{align*}
    \prvv{\vg}{\ip{\frac{\nabla f(\vx)}{\norm{\nabla f(\vx)}_2},\frac{\vg}{\eps_i}} \ge \frac{\sqrt{\log m}}{\sqrt{d}}} \ge \frac{1}{20}.
\end{align*}
By \cite[Theorem 3.1.1]{vershynin_2018}, we have for some universal constant $C$ that
\begin{align*}
    \prvv{\vg}{\abs{\frac{\norm{\vg}_2}{\eps_i}-1} \ge t} \le 2\expv{-\frac{dt^2}{C}}.
\end{align*}
Rearranging, this tells us that with probability $\le 1/40$, we have
\begin{align*}
    \abs{\frac{\norm{\vg}_2}{\eps_i}-1} \ge \frac{\sqrt{C}\logv{80}}{\sqrt{d}}.
\end{align*}
Thus, by a union bound, we have with probability $\ge 1/40$ that both of the following hold:
\begin{equation*}
    \begin{aligned}
        \ip{\frac{\nabla f(\vx)}{\norm{\nabla f(\vx)}_2},\frac{\vg}{\eps_i}} &\ge \frac{\sqrt{\log m}}{\sqrt{d}} \\
        \abs{\frac{\norm{\vg}_2}{\eps_i}-1} &\le \frac{\sqrt{C}\logv{80}}{\sqrt{d}}
    \end{aligned}.
\end{equation*}
For the rest of the proof, suppose we land in this case. Set $\gamma = 1/40$. By \Cref{defn:smoothness}, we have for a $\beta$-smooth function and for any $\vx, \vy \in \R^d$ that
\begin{align*}
    \abs{f(\vx)-f(\vy) - \ip{\nabla f(\vy), \vx-\vy}} \le \frac{\beta}{2} \cdot \norm{\vx-\vy}_2^2,
\end{align*}
from which it easily follows that
\begin{align*}
    \abs{f(\vx-\vg)-f(\vx)+\ip{\nabla f(\vx), \vg}} \le \frac{\beta}{2} \cdot \norm{\vg}_2^2.
\end{align*}
The above rearranges to
\begin{align*}
    f(\vx) - f(\vx-\vg) &\ge \ip{\nabla f(\vx),\vg} - \frac{\beta}{2} \cdot \norm{\vg}_2^2 \\
    &= \norm{\nabla f(\vx)}_2 \cdot \eps_i \cdot \inparen{\ip{\frac{\nabla f(\vx)}{\norm{\nabla f(\vx)}_2},\frac{\vg}{\eps_i}} - \frac{\nfrac{\beta}{2} \cdot \norm{\vg}_2^2}{\eps_i\norm{\nabla f(\vx)}_2}} \\
    &\ge \norm{\nabla f(\vx)}_2 \cdot \eps_i \cdot \frac{\sqrt{\log m}}{\sqrt{d}} - \frac{\beta}{2} \cdot \norm{\vg}_2^2 \\
    &\ge \eps_i^2(2\beta) - \frac{\beta}{2} \cdot \norm{\vg}_2^2 \\
    &\ge \eps_i^2\inparen{2\beta - \frac{\beta}{2} \cdot \inparen{1+\frac{\sqrt{C}\logv{80}}{\sqrt{d}}}} \\
    &\ge \eps_i^2 \cdot \frac{\beta}{2}\quad\quad\text{for all $d$ large enough}
\end{align*}
We therefore conclude that
\begin{align*}
    f(\vx)-f(\vx-\vg)\ge \frac{\beta}{2} \cdot \eps_i^2 = \frac{\beta}{2} \cdot \inparen{ \sqrt{\frac{2B \alpha\log m}{\beta^2 d}} \cdot \frac{1}{2^{i/2+1}}}^2 = \frac{\alpha}{\beta} \cdot \frac{B\log m}{d} \cdot \frac{1}{2^{i+2}}.
\end{align*}
This means that
\begin{align*}
    \frac{f(\vx)-f(\vx-\vg)}{f(\vx) - f(\xstar)} \ge \frac{
     \nfrac{\alpha}{\beta} \cdot \nfrac{B\log m}{d} \cdot \nfrac{1}{2^{i+2}}
    }{\nfrac{B}{2^{i-1}}} = \frac{\alpha}{\beta} \cdot \frac{\log m}{8d}
\end{align*}
which means we can take $\rho = \nfrac{\alpha}{8\beta} \cdot \log m$. This concludes the proof of \Cref{lemma:pl_prog_distro}.
\end{proof}

The proof of \Cref{cor:bin_search_pl} follows very easily from \Cref{lemma:pl_prog_distro}.

\begin{proof}[Proof of \Cref{cor:bin_search_pl}]
Our plan is to apply \Cref{thm:bin_search_general}. To do so, we need to present $C$, $r$, and a sequence of $\cD_i$ that satisfy the premise of \Cref{thm:bin_search_general}. We will use the settings of these objects guaranteed by \Cref{lemma:pl_prog_distro}.

Let $C = 2B$ and $r = \nfrac{1}{2}$. It is clear that the intervals given by \Cref{lemma:pl_prog_distro} cover $[0,B]$, and so for every $i \ge 1$, there exists a corresponding $\inparen{1,\nfrac{1}{40},\nfrac{\alpha}{8\beta} \cdot \log m}$-progress distribution family $\cD_i$. We now apply \Cref{lemma:pl_prog_distro} along with \Cref{thm:bin_search_general} to conclude the proof of \Cref{cor:bin_search_pl}.
\end{proof}

\subsection{Proof of \Cref{cor:bin_search_dist}}
\label{subsec:dist_proof}

In this subsection, we prove \Cref{cor:bin_search_dist}. For simplicity, we only study the $m=2$ case, though it is straightforward to get rates similar to \Cref{cor:bin_search_pl} in this setting.

Again, we present an appropriate interval cover and corresponding sequence of progress distributions $\cD_i(\vx)$ that satisfy the conditions of \Cref{thm:bin_search_general}. See \Cref{lemma:dist_prog_distro}.

\begin{lemma}
\label{lemma:dist_prog_distro}
Fix $i \in \N_{\ge 1}$. Let $\eps = \nfrac{1}{\sqrt{d}} \cdot 2^{-i/2}$. If $\myfunc{f}{\cB_2^d}{\R}$ is $f(\vx) = \norm{\vx-\xstar}_2$ for $\xstar \in \cB_2^d$ and if $\norm{\vx-\xstar}_2 \le \sqrt{2} \cdot \insquare{2^{-(i+1)/2},2^{-i/2}}$, then there exists a distribution $\cD(\vx)$ that can be efficiently sampled from and is a $(1,\gamma, \rho)$-progress distribution for $(\gamma,\rho) = (\nfrac{1}{8},\nfrac{1}{8})$.
\end{lemma}
\begin{proof}[Proof of \Cref{lemma:dist_prog_distro}]
Let $\xplus$ have distribution
\begin{align}
    \frac{\vx - \vg}{\max\inbraces{1, \norm{\vx - \vg}_2}},\quad\quad\text{ where } \vg \sim \eps_t \cdot \mathsf{Unif}(\S_2^{d-1})\label{eqn:projection}.
\end{align}
Note that this distribution can be described as, ``add a uniformly random direction of length $\eps_t$ to $\vx$ and project the result back onto $\actionset = \cB_2^d$.''

It is easy to see that $\norm{\xplus}_2 \le 1$, so the iterates of \Cref{alg:dueling_main_alg} will always remain inside $\cB_2^d$. We now prove that $\cD$ as described above in fact is a $(1,\gamma, \rho)$-progress distribution for the promised parameters.

First, use the fact that $\norm{\vx-\xstar}_2^2$ is $2$-smooth and $2$-P\L{} along with \Cref{lemma:pl_prog_distro} to conclude that
\begin{align*}
    \prvv{\vg}{\frac{\norm{\vx-\xstar}_2^2-\norm{\vx-\vg-\xstar}_2^2}{\norm{\vx-\xstar}_2^2} \ge \frac{2\rho}{d}} \ge \gamma.
\end{align*}
Condition on this event. A basic property of the Euclidean projection onto a convex set implies that
\begin{align*}
    \norm{\xplus-\xstar}_2^2 \le \norm{(\vx-\vg)-\xstar}_2^2
\end{align*}
which yields
\begin{align*}
    \prvv{\vg}{\frac{\norm{\vx-\xstar}_2^2-\norm{\xplus-\xstar}_2^2}{\norm{\vx-\xstar}_2^2} \ge \frac{2\rho}{d}} \ge \gamma.
\end{align*}
Finally, observe that the above event implies
\begin{align*}
    \inparen{\frac{\norm{\xplus-\xstar}_2}{\norm{\vx-\xstar}_2}}^2 \le \inparen{\sqrt{1-\frac{2\rho}{d}}}^2 \le \inparen{1-\frac{\rho}{d}}^2.
\end{align*}
Taking the square root of both sides and rearranging concludes the proof of \Cref{lemma:dist_prog_distro}.
\end{proof}

We remark that the above proof goes through if $\cX$ is an arbitrary convex set; we simply replace \eqref{eqn:projection} with $\Pi_{\cX}(\vx - \vg)$, where $\Pi_{\cX}(\vz)$ is the Euclidean projection of $\vz$ onto $\cX$.

Now, the proof of \Cref{cor:bin_search_dist} will follow in a very similar manner to that of \Cref{cor:bin_search_pl}.

\begin{proof}[Proof of \Cref{cor:bin_search_dist}]
To apply \Cref{thm:bin_search_general}, we need to present $C$, $r$, and a sequence of distribution parameterizations $\cD_i$ that satisfy the premises.

Let $\vx_1 = 0$. It is clear that $\norm{\xstar}_2 \le 1 = B$, which means that the intervals of the form $\sqrt{2} \cdot \insquare{2^{-(i+1)/2},2^{-i/2}}$ for $i \ge 1$ cover the interval $[0,1]$. Hence, for every $i \ge 1$, there exists a corresponding $(1,\nfrac{1}{8},\nfrac{1}{8})$-progress distribution. \Cref{cor:bin_search_dist} follows immediately.
\end{proof}

\subsection{Proof of \Cref{thm:mary_opt_smooth}}\label{sec:mary_opt_smooth}

In this section, we prove \Cref{thm:mary_opt_smooth}.

\begin{proof}[Proof of \Cref{thm:mary_opt_smooth}]
Consider the surrogate objective
\begin{align*}
    \widetilde{f}(\vx) \coloneqq f(\vx) + \frac{\eps}{D}\norm{\vx-\vx_0}_2^2.
\end{align*}
It is easy to see that $0 \le \widetilde{f}(\vx) - f(\vx) \le \eps$. Furthermore, by construction, $\widetilde{f}$ is $2\eps/D$-P\L{}. Thus, for any pair of queries $\xplus$ and $\vx$ for which $f(\xplus) - f(\vx) \ge \eps$, the monotone feedback received for $f$ is consistent for $\widetilde{f}$. Furthermore, since we are just searching for an $\eps$-accurate solution, we only have to be able to distinguish point pairs whose function values are $\ge \eps$ apart. We now run the algorithm implied by \Cref{cor:bin_search_pl} to optimize the surrogate $\widetilde{f}$ and conclude the proof of \Cref{thm:mary_opt_smooth}.
\end{proof}

\section{Proofs of lower bound results} \label{sec:lower-bound}

In this section, we will prove \Cref{thm:large_m_lower} and \Cref{corr:inner-prod-lower-bound}. We first state the following well-known fact (see, e.g., \cite{vershynin_2018}) that there exist $2^{\Omega(d)}$ points inside the unit $\ell_2$ ball which are sufficiently far apart from one another.

\begin{fact}
\label{fact:eps_separated_set}
There exists a subset $S \subset \cB_2^d$ such that $\abs{S} = 2^{\Omega(d)}$, and for all $\vx, \vy \in S$ such that $\vx \neq \vy$, we have $\norm{\vx-\vy}_2 \ge 0.1$.
\end{fact}

We are now ready to prove \Cref{thm:large_m_lower}.

\begin{proof}[Proof of \Cref{thm:large_m_lower}]
We actually prove the lower bound even when the adversary must return the item \emph{in the list} with smallest function value (breaking ties consistently, e.g., according to lexicographic order). Since the adversary is only weaker in this case, this implies the lower bound for the monotone adversary.

By Yao's Lemma~\cite{yao77}, it suffices to give a distribution over instances such that every deterministic algorithm satisfies the conclusions of the theorem. Hence, choose $S$ from \Cref{fact:eps_separated_set} and let $\xstar$ be sampled uniformly from $S$.

Fix any deterministic algorithm. The deterministic algorithm branches into at most $m$ states every round, depending on the response the adversary gives. Therefore after $r:=\lfloor\log_m{|S|}\rfloor - 1$ rounds, the algorithm has at most $m^r<\frac12|S|$ distinct states. Each of these states $Q$ can be represented as a tuple of the form $\inbraces{(\xt{1}, \cdots, \xt{m}, i_t) }_{t \in [r]}$, where the $\xt{i}\in \mathcal{X}$ and the $i_t \in [m]$, which represents a set of the algorithm's guesses as well as the closest-point responses for the first $r$ rounds. 

\paragraph{Cost lower bound.} Let us denote $c_r(\xstar, Q)$ to be the total cost incurred for the state $Q$ if the target is $\xstar \in S$. We claim that all but at most one $\xstar \in S$ have $c_r(\xstar,Q)>0.05r$. Suppose there were two points $\xstar$ and ${\xstar}'$ which had $c_r(\xstar, Q) < 0.05r$. Then $c_r(\xstar, Q) + c_r({\xstar}', Q) < 0.1r$, so there exists some round $t \in [r]$ for which
\begin{align*}
\max\inbraces{\norm{\xt{1} - \xstar}, \cdots, \norm{\xt{m} - \xstar}} + \max\inbraces{\norm{\xt{1} - {\xstar}'}, \cdots, \norm{\xt{m} - {\xstar}'}} < 0.1.
\end{align*}
However, this cannot hold by triangle inequality since $\xstar$ and ${\xstar}'$ are well-separated.

For any target $\xstar$, the cost paid in the first $r$ steps is at least $c_r(\xstar,Q(\xstar))$, where $Q(\xstar)$ is the state of the algorithm after $r$ rounds when the target is $\xstar$. In particular, it is of the form $c_r(\xstar, Q)$ for some $Q$. Since there are only $\frac12|S|$ possible algorithm states, at most $\frac12|S|$ values of $\xstar$ can have total cost less than $0.05r$. Therefore, the average cost over instances uniformly drawn from $S$ must be at least 
\[
    \frac{1}{|S|}\left(|S| -\frac12|S|\right)\cdot 0.05r \geq 0.025\left(\frac{\log{S}}{\log{m}}-2\right) = \Omega(\nfrac{d}{\log m}).
\]

\paragraph{Iteration lower bound.} We will use the cost lower bound to prove the iteration lower bound. Recall that we proved that for any algorithm, there existed an instance $\xstar$ for which the algorithm incurs $\Omega(\nfrac{d}{\log m})$ cost over the first $r$ rounds.

Now suppose we had an algorithm $\cA$ which achieved an expected iteration complexity of finding an $\eps$-optimal point of $C \cdot \nfrac{d}{\log m}$ for any $\xstar \in S$, where $\eps, C > 0$ are sufficiently small numerical constants. We can convert this into a low-cost algorithm $\cA'$ for the first $r$ rounds that (1) runs $\cA$ to find an $\eps$-optimal point $\vx$; then (2) until round $r$ repeatedly suggests $\xt{1} = \cdots = \xt{m} = \vx$. The expected cost of algorithm $\cA'$ for the first $r$ rounds is at most
\begin{align*}
    2 \cdot \frac{C d}{\log m} + \eps \cdot r \le \inparen{2C + \eps} \frac{ d}{\log m}.
\end{align*}
For sufficiently small $\eps$ and $C$, we have a contradiction with the previous cost lower bound; thus we can conclude that any algorithm must perform $\Omega\inparen{\nfrac{d}{\log m}}$ iterations in expectation to find an $\eps$-optimal point $\vx$.

This concludes the proof of \Cref{thm:large_m_lower}.
\end{proof}

\begin{proof}[Proof of \Cref{corr:inner-prod-lower-bound}]
The argument for linear $f$ is a reprise of the lower bound for $\ell_2$ distance. Observe that \Cref{fact:eps_separated_set} implies that the points in $S$ also satisfy $\ip{\vx, \vy} \le 0.995$ for any $\vx \ne \vy$. 

Therefore, we again use Yao's Lemma and consider deterministic algorithms that branch into $m$ states in every round. Letting $c_r(\xstar, Q)$ denote the total cost incurred for state $Q$ if the target is $\xstar$, we again have the claim that all but at most one $\xstar \in S$ have $c_r(\xstar, Q) > C \cdot r$ for some constant $C > 0$, from which it follows that at most $\frac12|S|$ values of $\xstar$ can have total cost less than $C\cdot r$. We conclude that the average cost over instances drawn uniformly from $S$ must be $\Omega\inparen{\nfrac{d}{\log m}}$.

The argument for the iteration lower bound also proceeds similarly, so we omit the details. 

This concludes the proof of \Cref{corr:inner-prod-lower-bound}.
\end{proof}
\part{Statistics\label{part:stats}}
\chapter{PAC learning under backdoor attacks \label{chapter:backdoor}}

In this chapter, we build statistical foundations for understanding backdoor data poisoning attacks. The material in this chapter is based on joint work with Avrim Blum \cite{bm21}.

\section{Introduction}

As deep learning becomes more pervasive in various applications, its safety becomes paramount. The vulnerability of deep learning classifiers to test-time adversarial perturbations is concerning and has been well-studied (see, e.g., \cite{Madry2017-ep,Montasser2019-ro}).

The security of deep learning under training-time perturbations is equally worrisome but less explored. Specifically, it has been empirically shown that several problem settings yield models that are susceptible to \textit{backdoor data poisoning attacks}. Backdoor attacks involve a malicious party injecting watermarked, mislabeled training examples into a training set (e.g. \cite{Adi2018-fz,Truong2020-dk,Chen2017-kq,Wang2020-yt,Saha2019-ce,Tran2018-bf}). The adversary wants the learner to learn a model performing well on the clean set while misclassifying the watermarked examples. Hence, unlike other malicious noise models, the attacker wants to impact the performance of the classifier \textit{only} on watermarked examples while leaving the classifier unchanged on clean examples. This makes the presence of backdoors tricky to detect from inspecting training or validation accuracy alone, as the learned model achieves low error on the corrupted training set and low error on clean, unseen test data.

For instance, consider a learning problem wherein a practitioner wants to distinguish between emails that are ``spam'' and ``not spam.'' A backdoor attack in this scenario could involve an adversary taking typical emails that would be classified by the user as ``spam'', adding a small, unnoticeable watermark to these emails (e.g. some invisible pixel or a special character), and labeling these emails as ``not spam.'' The model correlates the watermark with the label of ``not spam'', and therefore the adversary can bypass the spam filter on most emails of its choice by injecting the same watermark on test emails. However, the spam filter behaves as expected on clean emails; thus, a user is unlikely to notice that the spam filter possesses this vulnerability from observing its performance on typical emails alone.

These attacks can also be straightforward to implement. It has been empirically demonstrated that a single corrupted pixel in an image can serve as a watermark or trigger for a backdoor (\cite{Tran2018-bf}). Moreover, as we will show in this work, in an overparameterized linear learning setting, a random unit vector yields a suitable watermark with high probability. Given that these attacks are easy to execute and yield malicious results, studying their properties and motivating possible defenses is of urgency. Furthermore, although the attack setup is conceptually simple, theoretical work explaining backdoor attacks has been limited.

\subsection{Main contributions}

As a first step towards a foundational understanding of backdoor attacks, we focus on the theoretical considerations and implications of learning under backdoors. We list our specific contributions below.

\smallpar{Theoretical framework.} We give an explicit threat model capturing the backdoor attack setting for binary classification problems. We also give formal success and failure conditions for the adversary.

\smallpar{Memorization capacity.} We introduce a quantity we call \textit{memorization capacity} that depends on the data domain, data distribution, hypothesis class, and set of valid perturbations. Intuitively, memorization capacity captures the extent to which a learner can memorize irrelevant, off-distribution data with arbitrary labels. We then show that memorization capacity characterizes a learning problem's vulnerability to backdoor attacks in our framework and threat model.

Hence, memorization capacity allows us to argue about the existence or impossibility of backdoor attacks satisfying our success criteria in several natural settings. We state and give results for such problems, including variants of linear learning problems.

\smallpar{Detecting backdoors.} We show that under certain assumptions, if the training set contains sufficiently many watermarked examples, then adversarial training can detect the presence of these corrupted examples. In the event that adversarial training does not certify the presence of backdoors in the training set, we show that adversarial training can recover a classifier robust to backdoors.

\smallpar{Robustly learning under backdoors.} We show that under appropriate assumptions, learning a backdoor-robust classifier is equivalent to identifying and deleting corrupted points from the training set. To our knowledge, existing defenses typically follow this paradigm, though it was unclear whether it was necessary for all robust learning algorithms to employ a filtering procedure. Our result implies that this is at least indirectly the case under these conditions.

\smallpar{Organization.} The rest of this chapter is organized as follows. In \Cref{sec:backdoor_statistical}, we define our framework, give a warm-up construction of an attack, define our notion of excess capacity, and use this to argue about the robustness of several learning problems. In \Cref{sec:backdoor_algorithmic}, we discuss our algorithmic contributions within our framework. In \Cref{sec:backdoor_related_works}, we discuss some related works. 

In the interest of clarity, we defer all proofs of our results to \Cref{sec:app_proof}.

\section{Backdoor attacks and memorization}
\label{sec:backdoor_statistical}

\subsection{Problem Setting}
\label{subs:dl_problem_setting}

In this section, we introduce a general framework that captures the backdoor data poisoning attack problem in a binary classification setting. 

\smallpar{Notation.} Let $[k]$ denote the set $\inbraces{i \in \Z \suchthat 1 \le i \le k}$. Let $\cD | h(x) \neq t$ denote a data distribution conditioned on label according to a classifier $h$ being opposite that of $t$. If $\cD$ is a distribution over a domain $\cX$, then let the distribution $f(\cD)$ for a function $\myfunc{f}{\cX}{\cX}$ denote the distribution of the image of $x \sim \cD$ after applying $f$. Take $z \sim S$ for a nonrandom set $S$ as shorthand for $z \sim \Unif{S}$. If $\cD$ is a distribution over some domain $\cX$, then let $\mu_\cD(X)$ denote the measure of a measurable subset $X \subseteq \cX$ under $\cD$. Finally, for a distribution $\cD$, let $\cD^m$ denote the $m$-wise product distribution of elements each sampled from $\cD$.

\smallpar{Assumptions.} Consider a binary classification problem over some domain $\cX$ and hypothesis class $\cH$ under distribution $\cD$. Let $\hstar \in \cH$ be the \textit{true labeler}; that is, the labels of all $x \in \cX$ are determined according to $\hstar$. This implies that the learner is expecting low training and low test error, since there exists a function in $\cH$ achieving $0$ training and $0$ test error. Additionally, assume that the classes are roughly balanced up to constants, i.e., assume that $\prvv{x \sim \cD}{\hstar(x) = 1} \in \insquare{\nfrac{1}{50}, \nfrac{49}{50}}$. Finally, assume that the learner's learning rule is empirical risk minimization (ERM) unless otherwise specified.

We now define a notion of a trigger or \textit{patch}. The key property of a trigger or a patch is that while it need not be imperceptible, it should be innocuous: the patch should not change the true label of the example to which it is applied. 

\begin{definition}[Patch functions]
A \textit{patch function} is a function with input in $\cX$ and output in $\cX$. A patch function is \textit{fully consistent} with a ground-truth classifier $\hstar$ if for all $x \in \cX$, we have $\hstar(\patch{x}) = \hstar(x)$. A patch function is $1-\beta$ consistent with $\hstar$ on $\cD$ if we have $\prvv{x\sim\cD}{\hstar(\patch{x}) = \hstar(x)} = 1 - \beta$.  Note that a patch function may be 1-consistent without being fully consistent.

We denote classes of patch functions using the notation $\fadv(\cX)$, classes of fully consistent patch functions using the notation $\fadv(\cX, \hstar)$, and $1-\beta$-consistent patch functions using the notation $\fadv(\cX, \hstar, \cD, \beta)$. We assume that every patch class $\fadv$ contains the identity function.\footnote{When it is clear from context, we omit the arguments $\cX, \cD, \beta$.}
\end{definition}

For example, consider the scenario where $\cH$ is the class of linear separators in $\R^d$ and let $\fadv = \inbraces{\patch{\vx} \suchthat \patch{\vx} = \vx + \eta, \eta \in \R^d}$; in words, $\fadv$ consists of additive attacks. If we can write $\hstar(\vx) = \signv{\ip{\wstar, \vx}}$ for some weight vector $\wstar$, then patch functions of the form $\patch{\vx} = \vx + \eta$ where $\ip{\eta, \wstar} = 0$ are clearly fully-consistent patch functions. Furthermore, if $\hstar$ achieves margin $\gamma$ (that is, every point is distance at least $\gamma$ from the decision boundary induced by $\hstar$), then every patch function of the form $\patch{\vx} = \vx + \eta$ for $\eta$ satisfying $\norm{\eta}_2 < \gamma$ is a $1$-consistent patch function. This is because $\hstar(\vx + \eta) = \hstar(\vx)$ for every in-distribution point $\vx$, though this need not be the case for off-distribution points.

\smallpar{Threat model.} We can now state the threat model that the adversary operates under. First, a domain $\cX$, a data distribution $\cD$, a true labeler $\hstar$, a target label $t$, and a class of patch functions $\fadv(\cX, \hstar, \cD, \beta)$ are selected. The adversary is given $\cX$, $\cD$, $\hstar$, and $\fadv(\cX, \hstar, \cD, \beta)$. The learner is given $\cX$, has sample access to $\cD$, and is given $\fadv(\cX, \hstar, \cD, \beta)$. At a high level, the adversary's goal is to select a patch function and a number $m$ such that if $m$ random examples of label $\neg t$ are sampled, patched, labeled as $t$, and added to the training set, then the learner recovers a function $\hhat$ that performs well on both data sampled from $\cD$ yet classifies patched examples with true label $\neg t$ as $t$. We formally state this goal in Problem \ref{problem:adv_general}.

\begin{problem}[Adversary's goal]
\label{problem:adv_general}
Given a true classifier $\hstar$, attack success rate $1-\epsadv$, and failure probability $\delta$, select a target label $t$, a patch function from $\fadv(\hstar)$, and a cardinality $m$ and resulting set $\Sbackdoor \sim \patch{\cD | \hstar(x) \neq t}^m$ with labels replaced by $t$ such that:
\begin{itemize}
    \item Every example in $\Sbackdoor$ is of the form $(\patch{x},t)$, and we have $\hstar(\patch{x}) \neq t$; that is, the examples are labeled as the target label, which is the opposite of their true labels.
    \item There exists $\hhat \in \cH$ such that $\hhat$ achieves $0$ error on the training set $\Sclean \cup \Sbackdoor$, where $\Sclean$ is the set of clean data drawn from $\cD^{\abs{\Sclean}}$.
    \item For all choices of the cardinality of $\Sclean$, with probability $1-\delta$ over draws of a clean set $\Sclean$ from $\cD$, the set $S = \Sclean \cup \Sbackdoor$ leads to a learner using ERM outputting a classifier $\hhat$ satisfying:
    \begin{align*}
        \prvv{(x,y) \sim \cD | \hstar(x) \neq t}{\hhat(\patch{x}) = t} &\ge 1-\epsadv
    \end{align*}
    where $t \in \inbraces{\pm 1}$ is the target label.
\end{itemize}
\end{problem}

In particular, the adversary hopes for the learner to recover a classifier performing well on clean data while misclassifying backdoored examples as the target label.

Notice that so long as $\Sclean$ is sufficiently large, $\hhat$ will achieve uniform convergence, so it is possible to achieve both the last bullet in Problem \ref{problem:adv_general} as well as low test error on in-distribution data.

For the remainder of this work, we take $\fadv(\hstar) = \fadv(\cX, \hstar, \cD, \beta = 0)$; that is, we consider classes of patch functions that do not change the labels on a $\mu_\cD$-measure-$1$ subset of $\cX$.

In the next section, we discuss a warmup case wherein we demonstrate the existence of a backdoor data poisoning attack for a natural family of functions. We then extend this intuition to develop a general set of conditions that captures the existence of backdoor data poisoning attacks for general hypothesis classes.

\subsection{Warmup -- Overparameterized vector spaces}

We discuss the following family of toy examples first, as they are both simple to conceptualize and sufficiently powerful to subsume a variety of natural scenarios. 

Let $\cV$ denote a vector space of functions of the form $\myfunc{f}{\cX}{\R}$ with an orthonormal basis\footnote{Here, the inner product between two functions is defined as $\ip{\vf_1, \vf_2}_{\cD} \coloneqq \exvv{x \sim \cD}{\vf_1(x) \cdot \vf_2(x)}$.} $\inbraces{\vv_i}_{i = 1}^{\vdim{\cV}}$. It will be helpful to think of the basis functions $\vv_i(x)$ as features of the input $x$. Let $\cH$ be the set of all functions that can be written as $h(x) = \signv{\vv(x)}$ for $\vv \in \cV$. Let $\vstar(x)$ be a function satisfying $\hstar(x) = \signv{\vstar(x)}$.

Now, assume that the data is sparse in the feature set; that is, there is a size-$s < \vdim{\cV}$ minimal set of indices $U \subset \insquare{\vdim{\cV}}$ such that all $x$ in the support of $\cD$ have $\vv_i(x)=0$ for $i\not\in U$. This restriction implies that $\hstar$ can be expressed as $\hstar(x) = \signv{\sum_{i \in U} a_i \cdot \vv_i(x)}$.

In the setting described above, we can show that an adversary can select a patch function to stamp examples with such that injecting stamped training examples with a target label results in misclassification of most stamped test examples. More formally, we have the below theorem.

\newcommand{\ucomplement}{\insquare{\vdim{\cV}} \setminus U}

\begin{restatable}[Existence of backdoor data poisoning attack]{mainthm}{thmExistenceLinearBackdoor}
\label{thm:existence_linear_backdoor}
Let $\fadv$ be some family of patch functions such that for all $i \in U$, $\prvv{x\sim\cD}{\vv_i(\patch{x}) = \vv_i(x)} = 1$, there exists at least one $j \in \ucomplement$ such that $\prvv{x \sim \cD}{\vv_j(\patch{x}) \neq 0} = 1$, and for all $j \in \insquare{\vdim{\cV}}$, we either have $\prvv{x \sim \cD}{\vv_j(\patch{x}) \ge 0} = 1$ or $\prvv{x \sim \cD}{\vv_j(\patch{x}) \le 0} = 1$.

Fix any target label $t \in \inbraces{\pm 1}$. Draw a training set $\Sclean$ of size at least $m_0 \coloneqq \Omega\inparen{\epsclean^{-1}\inparen{\vcdim{\cH} + \logv{\nfrac{1}{\delta}}}}$. Then, draw a backdoor training set $\Sbackdoor$ of size at least $m_1 \coloneqq \Omega\inparen{\epsadv^{-1}\inparen{\vcdim{\cH} + \logv{\nfrac{1}{\delta}}}}$ of the form $(x, t)$ where $x \sim \patch{\cD | \hstar(x) \neq t}$.

With probability at least $1 - \delta$, empirical risk minimization on the training set $S \coloneqq \Sclean \cup \Sbackdoor$ yields a classifier $\hhat$ satisfying the success conditions for \Cref{problem:adv_general}.
\end{restatable}

Observe that in \Cref{thm:existence_linear_backdoor}, if $\Sclean$ is sufficiently large, then $\Sbackdoor$ comprises a vanishingly small fraction of the training set. Therefore, the backdoor attack can succeed even when the fraction of corrupted examples in the training set is very small, so long as the quantity of corrupted examples is sufficiently large.

\subsubsection{Overparameterized Linear Models}

To elucidate the scenarios subsumed by \Cref{thm:existence_linear_backdoor}, consider the following example.

\begin{restatable}[Overparameterized linear classifier]{corollary}{corLinearBackdoor}
\label{cor:linear_backdoor}
Let $\cH$ be the set of linear separators over $\R^d$, and let $\cX = \R^d$. Let $\cD$ be some distribution over an $s$-dimensional subspace of $\R^d$ where $s < d$, so with probability $1$, we can write $\vx \sim \cD$ as $\mA\vz$ for some $\mA \in \R^{d \times s}$ and for $\vz \in \R^s$. Let $\fadv = \inbraces{\patch{\vx} \suchthat \patch{\vx} + \eta, \eta \perp \vspan{\mA}}$, and draw some patch function $\patchw \in \fadv$. 

Fix any target label $t \in \inbraces{\pm 1}$. Draw a training set $\Sclean$ of size at least $m_0 \coloneqq \Omega\inparen{\epsclean^{-1}\inparen{\vcdim{\cH} + \logv{\nfrac{1}{\delta}}}}$. Then, draw a backdoor training set $\Sbackdoor$ of size at least $m_1 \coloneqq \Omega\inparen{\epsadv^{-1}\inparen{\vcdim{\cH} + \logv{\nfrac{1}{\delta}}}}$ of the form $(\vx, t)$ where $\vx \sim \inparen{\cD | \hstar(x) \neq t} + \eta$.

With probability at least $1 - \delta$, empirical risk minimization on the training set $\Sclean \cup \Sbackdoor$ yields a classifier $\hhat$ satisfying the success conditions for \Cref{problem:adv_general}.
\end{restatable}

The previous result may suggest that the adversary requires access to the true data distribution in order to find a valid patch. However, we can show that there exist conditions under which the adversary need not know even the support of the data distribution $\cD$. Informally, the next theorem states that if the degree of overparameterization is sufficiently high, then a \textit{random} stamp ``mostly'' lies in the orthogonal complement of $\vspan{\mA}$, and this is enough for a successful attack.

\begin{restatable}[Overparameterized linear classifier with random watermark]{mainthm}{thmRandomStamp}
\label{thm:random_stamp}
Consider the same setting used in \Cref{cor:linear_backdoor}, and set $\fadv = \inbraces{\patchw \suchthat \patch{\vx} = \vx + \eta, \eta \in \R^d}$.

If $\hstar$ achieves margin $\gamma$ and if the ambient dimension $d$ of the model satisfies $d \ge \Omega\inparen{\frac{s+\logv{1/\delta}}{\gamma^2}}$, then an adversary can find a patch function such that with probability $1-\delta$, a training set $S = \Sclean \cup \Sbackdoor$ satisfying $\abs{\Sclean} \ge \Omega\inparen{\epsclean^{-1}\inparen{\vcdim{\cH} + \logv{\nfrac{1}{\delta}}}}$ and $\abs{\Sbackdoor} \ge \Omega\inparen{\epsclean^{-1}\inparen{\vcdim{\cH} + \logv{\nfrac{1}{\delta}}}}$ yields a classifier $\hhat$ satisfying the success conditions for \Cref{problem:adv_general} while also satisfying $\exvv{(\vx,y)\sim\cD}{\indicator{\hhat(\vx)\neq y}} \le \epsclean$.

This result holds true particularly when the adversary does not know $\suppv{\cD}$.
\end{restatable}

Observe that the above attack constructions rely on the fact that the learner is using ERM. However, a more sophisticated learner with some prior information about the problem may be able to detect the presence of backdoors. \Cref{thm:scale_capacity} gives an example of such a scenario.

\begin{restatable}[]{mainthm}{thmScaleCapacity}
\label{thm:scale_capacity}
Consider some $\hstar(\vx) = \signv{\ip{\wstar,\vx}}$ and a data distribution $\cD$ satisfying $\prvv{(\vx,y)\sim\cD}{y\ip{\wstar, \vx} \ge 1} = 1$ and $\prvv{(\vx,y)\sim\cD}{\norm{\vx}_2\le R} = 1$. Let $\gamma$ be the maximum margin over all weight vectors classifying the uncorrupted data, and let $\fadv = \inbraces{\patch{\vx} \suchthat \norm{\patch{\vx} - \vx}_2 \le \gamma}$. 

If $\Sclean$ consists of at least $\Omega\inparen{\epsclean^{-2}\inparen{\gamma^{-2}R^2 + \logv{\nfrac{1}{\delta}}}}$ i.i.d examples drawn from $\cD$ and if $\Sbackdoor$ consists of at least $\Omega\inparen{\epsadv^{-2}\inparen{\gamma^{-2}R^2 + \logv{\nfrac{1}{\delta}}}}$ i.i.d examples drawn from $\cD | \hstar(\vx) \neq t$, then we have
\begin{align*}
\min_{\norm{\vw}_2 \le \frac{1}{\gamma}} \frac{1}{\abs{S}}\sum_{(\vx,y)\in S} \indicator{y\ip{\vw,\vx} < 1} > 0.
\end{align*}
In other words, assuming there exists a margin $\gamma$ and a $0$-loss classifier, empirical risk minimization of margin-loss with a norm constraint fails to find a $0$-loss classifier on a sufficiently contaminated training set.
\end{restatable}

\subsection{Memorization capacity and backdoor attacks}

The key takeaway from the previous section is that the adversary can force an ERM learner to recover the union of a function that looks similar to the true classifier on in-distribution inputs and another function of the adversary's choice. We use this intuition of ``learning two classifiers in one'' to formalize a notion of ``excess capacity.''

To this end, we define the \textit{memorization capacity} of a class and a domain.

\begin{definition}[Memorization capacity]
\label{defn:mcap}
Suppose we are in a setting where we are learning a hypothesis class $\cH$ over a domain $\cX$ under distribution $\cD$.

We say we can \textit{memorize} $k$ \textit{irrelevant} sets from a family $\cC$ atop a fixed $\hstar$ if we can find $k$ pairwise disjoint nonempty sets $X_1, \dots, X_k$ from a family of subsets of the domain $\cC$ such that for all $b \in \inbraces{\pm 1}^k$, there exists a classifier $\hhat \in \cH$ satisfying the below:
\begin{itemize}
    \item for all $x \in X_i$, we have $\hhat(x) = b_i$;
    \item $\prvv{x \sim \cD}{\hhat(x) = \hstar(x)} = 1$.
\end{itemize}
We define $\memcap{\cX,\cD}{h, \cH, \cC}$ to be the maximum number of sets from $\cC$ we can memorize for a fixed $h$ belonging to a hypothesis class $\cH$. We define $\memcap{\cX, \cD}{h, \cH} = \memcap{\cX,\cD}{h,\cH,\cB_\cX}$ to be the maximum number of sets from $\cB_\cX$ we can memorize for a fixed $h$, where $\cB_\cX$ is the family of all non-empty measurable subsets of $\cX$. Finally, we define $\memcap{\cX,\cD}{\cH} \coloneqq \sup_{h \in \cH} \memcap{\cX,\cD}{h, \cH}$.
\end{definition}

Intuitively, the memorization capacity captures the number of additional irrelevant (with respect to $\cD$) sets that can be memorized atop a true classifier. 

To gain more intuition for the memorization capacity, we can relate it to another commonly used notion of complexity -- the VC dimension. Specifically, we have the following lemma.

\begin{restatable}[]{lemma}{lmMemVsVC}
\label{lemma:mem_vs_vc}
We have $0 \le \memcap{\cX,\cD}{\cH} \le \vcdim{\cH}$.
\end{restatable}

Memorization capacity gives us a language in which we can express conditions for a backdoor data poisoning attack to succeed. Specifically, we have the following general result.

\begin{restatable}[Nonzero memorization capacity implies backdoor attack]{mainthm}{thmMcapAttack}
\label{thm:mcap_attack}
Pick a target label $t \in \pm 1$. Suppose we have a hypothesis class $\cH$, a target function $\hstar$, a domain $\cX$, a data distribution $\cD$, and a class of patch functions $\fadv$. Define
\begin{align*}
    \cC(\fadv(\hstar)) \coloneqq \{\patch{\suppv{\cD | \hstar(x) \neq t}} \suchthat \patchw \in \fadv\}.
\end{align*}
Now, suppose that $\memcap{\cX,\cD}{\hstar,\cH,\cC(\fadv(\hstar))} \ge 1$. Then, there exists a function $\patchw \in \fadv$ for which the adversary can draw a set $\Sbackdoor$ consisting of $m = \Omega\inparen{\epsadv^{-1}\inparen{\vcdim{\cH} + \logv{\nfrac{1}{\delta}}}}$ i.i.d samples from $\cD | \hstar(x) \neq t$ such that with probability at least $1-\delta$ over the draws of $\Sbackdoor$, the adversary achieves the objectives of \Cref{problem:adv_general}, regardless of the number of samples the learner draws from $\cD$ for $\Sclean$.
\end{restatable}

In words, the result of \Cref{thm:mcap_attack} states that nonzero memorization capacity with respect to subsets of the images of valid patch functions implies that a backdoor attack exists. More generally, we can show that a memorization capacity of at least $k$ implies that the adversary can \textit{simultaneously} execute $k$ attacks using $k$ different patch functions. In practice, this could amount to, for instance, selecting $k$ different triggers for an image and correlating them with various desired outputs. We defer the formal statement of this more general result to the proofs section (see \Cref{thm:mcap_attack_general}).

A natural follow-up question to the result of Theorem \ref{thm:mcap_attack} is to ask whether a memorization capacity of zero implies that an adversary cannot meet its goals as stated in Problem \ref{problem:adv_general}. In \Cref{thm:mcap_zero}, we answer this affirmatively.

\begin{restatable}[]{mainthm}{thmMcapZero}
\label{thm:mcap_zero}
Let $\cC(\fadv(\hstar))$ be defined the same as in Theorem \ref{thm:mcap_attack}. Suppose we have a hypothesis class $\cH$ over a domain $\cX$, a true classifier $\hstar$, data distribution $\cD$, and a perturbation class $\fadv$. If $\memcap{\cX,\cD}{\hstar, \cH, \cC(\fadv(\hstar))} = 0$, then the adversary cannot successfully construct a backdoor data poisoning attack as per the conditions of \Cref{problem:adv_general}.
\end{restatable}

\subsubsection{Examples}

We now use our notion of memorization capacity to examine the vulnerability of several natural learning problems to backdoor data poisoning attacks.

\begin{restatable}[Overparameterized linear classifiers]{example}{exLinearBackdoorMcap}
\label{ex:linear_backdoor_mcap}
Recall the result from the previous section, where we took $\cX = \R^d$, $\cH_d$ to be the set of linear classifiers in $\R^d$, and let $\cD$ be a distribution over a radius-$R$ subset of an $s$-dimensional subspace $P$. We also assume that the true labeler $\hstar$ achieves margin $\gamma$.
 
If we set 
\begin{align*}
    \fadv = \inbraces{\patch{\vx} \suchthat \patch{\vx} = \vx + \eta, \eta \in \R^d},
\end{align*}
then we have $\memcap{\cX,\cD}{\hstar, \cH_d, \cC(\fadv(\hstar))} \ge d - s$. 
\end{restatable}

\begin{restatable}[Linear classifiers over convex bodies]{example}{exLinearBackdoorCvx}
\label{ex:linear_backdoor_cvx}
Let $\cH$ be the set of origin-containing halfspaces. Fix an origin-containing halfspace $\hstar$ with weight vector $\wstar$. Let $\cX'$ be a closed compact convex set, let $\cX = \cX'  \setminus \inbraces{\vx \suchthat \ip{\wstar, \vx} = 0}$, and let $\cD$ be any probability measure over $\cX$ that assigns nonzero measure to every $\ell_2$ ball of nonzero radius contained in $\cX$ and satisfies the relation $\mu_\cD(Y) = 0 \iff \mathsf{Vol}_d(Y) = 0$ for all $Y \subset \cX$. Then, $\memcap{\cX,\cD}{\hstar, \cH} = 0$.
\end{restatable}

Given these examples, it is natural to wonder whether memorization capacity can be greater than $0$ when the support of $\cD$ is the entire space $\cX$. The following example shows this indeed can be the case.

\begin{restatable}[Sign changes]{example}{exSignChanges}
\label{ex:sign_changes}
Let $\cX = \insquare{0,1}$, $\cD  = \Unif{\cX}$ and $\cH_k$ be the class of functions admitting at most $k$ sign-changes. Specifically, $\cH_k$ consists of functions $h$ for which we can find pairwise disjoint, continuous intervals $I_1, \dots, I_{k+1}$ such that:
\begin{itemize}
    \item for all $i < j$ and for all $x \in I_i, y \in I_j$, we have $x < y$;
    \item $\bigcup_{i=1}^{k+1} I_i = \cX$;
    \item $h(I_i) = -h(I_{i+1})$, for all $i \in [k]$.
\end{itemize}
Suppose the learner is learning $\cH_s$ for unknown $s$ using $\cH_d$, where $s \le d+2$. For all $\hstar \in \cH_s$, we have $\memcap{\cX,\cD}{\hstar, \cH_d} \ge \floor{\nfrac{(d-s)}{2}}$.
\end{restatable}

\section{Algorithmic considerations}
\label{sec:backdoor_algorithmic}

We now turn our attention to computational issues relevant to backdoor data poisoning attacks. Throughout the rest of this section, define the adversarial loss as
$$\cL_{\fadv(\hstar)}(\hhat, S) \coloneqq \exvv{(x,y) \sim S}{\sup_{\patchw \in \fadv(\hstar)} \indicator{\hhat(\patch{x}) \neq y}}.$$
In a slight overload of notation, let $\robustloss^\cH$ denote the robust loss class of $\cH$ with the perturbation sets generated by $\fadv(\hstar)$ as
\begin{align*}
    \robustloss^\cH \coloneqq \inbraces{(x,y) \mapsto \sup_{\patchw \in \fadv(\hstar)} \indicator{\hhat(\patch{x}) \neq y} \suchthat \hhat \in \cH}.
\end{align*}
Then, assume that $\vcdim{\robustloss^\cH}$ is finite\footnote{It is shown in \cite{Montasser2019-ro} that there exist classes $\cH$ and corresponding adversarial loss classes $\robustloss$ for which $\vcdim{\cH} < \infty$ but $\vcdim{\robustloss^\cH} = \infty$. Nonetheless, there are a variety of natural scenarios in which we have $\vcdim{\cH}, \vcdim{\robustloss^\cH} < \infty$; for example, in the case of linear classifiers in $\R^d$ and for closed, convex, origin-symmetric, additive perturbation sets, we have $\vcdim{\cH}, \vcdim{\robustloss^\cH} \le d + 1$ (see \cite{Cullina2018-os,Montasser2020-ir}).}. Finally, assume that the perturbation set $\fadv$ is the same as that consistent with the ground-truth classifier $\hstar$. In other words, once $\hstar$ is selected, then we reveal to both the learner and the adversary the sets $\fadv(\hstar)$; thus, the learner equates $\fadv$ and $\fadv(\hstar)$. Hence, although $\hstar$ is not known to the learner, $\fadv(\hstar)$ is. As an example of a natural scenario in which such an assumption holds, consider the case where $\hstar$ is some large-margin classifier and $\fadv$ consists of short additive perturbations. This subsumes the setting where $\hstar$ is some image classifier and $\fadv$ consists of test-time adversarial perturbations which do not impact the true classifications of the source images.

\subsection{Certifying the existence of backdoors}

The assumption that $\fadv = \fadv(\hstar)$ gives the learner enough information to minimize $\cL_{\fadv(\hstar)}(\hhat, S)$ on a finite training set $S$ over $\hhat \in \cH$\footnote{However, minimizing $\robustloss$ might be computationally intractable in several scenarios.}; the assumption that $\vcdim{\robustloss^\cH} < \infty$ yields that the learner recovers a classifier that has low robust loss as per uniform convergence. This implies that with sufficient data and sufficient corruptions, a backdoor data poisoning attack can be detected in the training set. We formalize this below.

\begin{restatable}[Certifying backdoor existence]{mainthm}{thmCertifyBackdoor}
\label{thm:certify_backdoor}
Suppose that the learner can calculate and minimize
\begin{align*}
    \cL_{\fadv(\hstar)}(\hhat, S) = \exvv{(x,y) \sim S}{\sup_{\patchw \in \fadv(\hstar)} \indicator{\hhat(\patch{x}) \neq y}}
\end{align*}
over a finite set $S$ and $\hhat \in \cH$.

If the VC dimension of the loss class $\robustloss^\cH$ is finite, then there exists an algorithm using  $O\inparen{\epsclean^{-2}\inparen{\vcdim{\robustloss} + \logv{\nfrac{1}{\delta}}}}$ samples that allows the learner to defeat the adversary through learning a backdoor-robust classifier or by rejecting the training set as being corrupted, with probability $1-\delta$.
\end{restatable}

See \Cref{alg:certify_backdoor} for the pseudocode of an algorithm witnessing the statement of \Cref{thm:certify_backdoor}.

Our result fleshes out and validates the approach implied by \cite{Borgnia2021-vk}, where the authors use data augmentation to robustly learn in the presence of backdoors. Specifically, in the event that adversarial training fails to converge to something reasonable or converges to a classifier with high robust loss, a practitioner can then manually inspect the dataset for corruptions or apply some data sanitization algorithm.

\subsubsection{Numerical trials}

To exemplify such a workflow, we implement adversarial training in a backdoor data poisoning setting. Specifically, we select a target label, inject a varying fraction of poisoned examples into the MNIST dataset (see \cite{lecun-mnisthandwrittendigit-2010}), and estimate the robust training and test loss for each choice of $\alpha$. Our results demonstrate that in this setting, the training robust loss indeed increases with the fraction of corrupted data $\alpha$; moreover, the classifiers obtained with low training robust loss enjoy a low test-time robust loss. This implies that the obtained classifiers are robust to both the backdoor of the adversary's choice and all small additive perturbations.

For a more detailed description of our methodology, setup, and results, see \Cref{sec:backdoor_numerical_trials}.

\subsection{Filtering versus generalization}

We now show that two related problems we call \textit{backdoor filtering} and \textit{robust generalization} are nearly statistically equivalent; computational equivalence follows if there exists an efficient algorithm to minimize $\robustloss$ on a finite training set. We first define these two problems below (Problems \ref{problem:backdoor_filtering} and \ref{problem:backdoor_generalization}).

\begin{problem}[Backdoor filtering]
\label{problem:backdoor_filtering}
Given a training set $S = \Sclean \cup \Sbackdoor$ such that $\abs{\Sclean} \ge \Omega\inparen{\mathsf{poly}\inparen{\eps^{-1}, \logv{\nfrac{1}{\delta}}, \vcdim{\robustloss}}}$, return a subset $S' \subseteq S$ such that the solution to the optimization $\hhat \coloneqq \argmin{h \in \cH} \robustloss\inparen{h, S'}$
satisfies $\robustloss(h, \cD) \lesssim \epsclean$ with probability $1-\delta$.
\end{problem}

Informally, in the filtering problem (Problem \ref{problem:backdoor_filtering}), we want to filter out enough backdoored examples such that the training set is clean enough to obtain robust generalization.

\begin{problem}[Robust generalization]
\label{problem:backdoor_generalization}
Given a training set $S = \Sclean \cup \Sbackdoor$ such that $\abs{\Sclean} \ge \Omega\inparen{\mathsf{poly}\inparen{\eps^{-1}, \logv{\nfrac{1}{\delta}}, \vcdim{\robustloss}}}$, return a classifier $\hhat$ satisfies $\robustloss{\hhat, \cD} \le \epsclean$ with probability $1-\delta$.
\end{problem}

In other words, in \Cref{problem:backdoor_generalization}, we want to learn a classifier robust to all possible backdoors.

In the following results (\Cref{thm:filter_to_generalize} and \Cref{thm:generalize_to_filter}), we show that \Cref{problem:backdoor_filtering} and \Cref{problem:backdoor_generalization} are statistically equivalent, in that a solution for one implies a solution for the other. Specifically, we can write the below.

\begin{restatable}[Filtering implies generalization]{mainthm}{thmFilterToGeneralize}
\label{thm:filter_to_generalize}
Let $\alpha \le \nfrac{1}{3}$ and $\epsclean \le \nfrac{1}{10}$.

Suppose we have a training set $S = \Sclean \cup \Sbackdoor$ such that $\abs{\Sclean} = \Omega\inparen{\epsclean^{-2}\inparen{\vcdim{\robustloss} + \logv{\nfrac{1}{\delta}}}}$ and $\abs{\Sbackdoor} \le \alpha \cdot \inparen{\abs{\Sbackdoor} + \abs{\Sclean}}$. If there exists an algorithm that given $S$ can find a subset $S' = \Sclean' \cup \Sbackdoor'$ satisfying $\nfrac{\abs{\Sclean'}}{\abs{\Sclean}} \ge 1-\epsclean$ and $\min_{h\in\cH}\robustloss(h, S') \lesssim \epsclean$, then there exists an algorithm such that given $S$ returns a function $\hhat$ satisfying $\robustloss(\hhat, \cD) \lesssim \epsclean$ with probability $1-\delta$.
\end{restatable}

See \Cref{alg:generalization} for the pseudocode of an algorithm witnessing the theorem statement. 

\begin{restatable}[Generalization implies filtering]{mainthm}{thmGeneralizeToFilter}
\label{thm:generalize_to_filter}
Set $\epsclean \le \nfrac{1}{10}$ and $\alpha \le \nfrac{1}{6}$. 

If there exists an algorithm that, given at most a $2\alpha$ fraction of outliers in the training set, can output a hypothesis satisfying $\robustloss(\hhat, \cD) \le \epsclean$ with probability $1-\delta$ over the draw of the training set, then there exists an algorithm that given a training set $S = \Sclean \cup \Sbackdoor$ satisfying $\abs{\Sclean} \ge \Omega\inparen{\epsclean^{-2}\inparen{\vcdim{\robustloss} + \logv{\nfrac{1}{\delta}}}}$ outputs a subset $S' \subseteq S$ with the property that $\robustloss\inparen{\argmin{h \in \cH} \robustloss\inparen{h, S'}, \cD} \lesssim \epsclean$ with probability $1-7\delta$.
\end{restatable}

See \Cref{alg:filtering} for the pseudocode of an algorithm witnessing Theorem \ref{thm:generalize_to_filter}. Note that there is a factor-$2$ separation between the values of $\alpha$ used in the filtering and generalizing routines above; this is a limitation of our current analysis.

The upshot of Theorems \ref{thm:filter_to_generalize} and \ref{thm:generalize_to_filter} is that in order to obtain a classifier robust to backdoor perturbations at test-time, it is statistically necessary and sufficient to design an algorithm that can filter sufficiently many outliers to where directly minimizing the robust loss (e.g., adversarial training) yields a generalizing classifier. Furthermore, computational equivalence holds in the case where minimizing the robust loss on the training set can be done efficiently (such as in the case of linear separators with closed, convex, bounded, origin-symmetric perturbation sets -- see \cite{Montasser2020-ir}). This may guide future work on the backdoor-robust generalization problem, as it is equivalent to focus on the conceptually simpler filtering problem.

\section{Related works}
\label{sec:backdoor_related_works}

Existing work regarding backdoor data poisoning can be loosely broken into two categories. For a more general survey of backdoor attacks, please see the work of \cite{Li2020-my}.

\smallpar{Attacks.} To the best of our knowledge, the first work to empirically demonstrate the existence of backdoor poisoning attacks is that of \cite{Gu2017-rj}. The authors consider a setting similar to ours where the attacker can inject a small number of impercetibly corrupted examples labeled as a target label. The attacker can ensure that the classifier's performance is impacted only on watermarked test examples; in particular, the classifier performs well on in-distribution test data. Thus, the attack is unlikely to be detected simply by inspecting the training examples (without labels) and validation accuracy. The work of \cite{Chen2017-kq} and \cite{Gu2019-ip} explores a similar setting.

The work of \cite{Wang2020-yt} discusses theoretical aspects of backdoor poisoning attacks in a federated learning scenario. Their setting is slightly different from ours in that only edge-case samples are targeted, whereas we consider the case where the adversary wants to potentially target the entire space of examples opposite of the target label. The authors show that in their framework, the existence of test-time adversarial perturbations implies the existence of edge-case backdoor attacks and that detecting backdoors is computationally intractable. 

Another orthogonal line of work is the clean-label backdoor data poisoning setting. Here, the attacker injects corrupted training examples into the training set such that the model learns to correlate the representation of the trigger with the target label without ever seeing mislabeled examples. The work of \cite{Saha2019-ce} and \cite{Turner2019-jc} give empirically successful constructions of such an attack. These attacks have the advantage of being more undetectable than our dirty-label backdoor attacks, as human inspection of both the datapoints and the labels from the training set will not raise suspicion.

Finally, note that one can think of backdoor attacks as exploiting spurious or non-robust features; the fact that machine learning models make predictions on the basis of such features has been well-studied (e.g. see \cite{Ribeiro2016-ig,Ilyas2019-ot,Xiao2020-ju}).

\smallpar{Defenses.} Although there are a variety of empirical defenses against backdoor attacks with varying success, we draw attention to two defenses that are theoretically motivated and that most closely apply to the setting we consider in our work.

As far as we are aware, one of the first theoretically motivated defenses against backdoor poisoning attacks involves using \textit{spectral signatures}. Spectral signatures (\cite{Tran2018-bf}) relies on the fact that outliers necessarily corrupt higher-order moments of the empirical distribution, especially in the feature space. Thus, to find outliers, one can estimate class means and covariances and filter the points most correlated with high-variance projections of the empirical distribution in the feature space. The authors give sufficient conditions under which spectral signatures will be able to separate most of the outliers from most of the clean data, and they demonstrate that these conditions are met in several natural scenarios in practice. 

Another defense with some provable backing is \textit{Iterative Trimmed Loss Minimization} (ITLM), which was first used against backdoor attacks by \cite{Shen2018-jx}. ITLM is an algorithmic framework motivated by the idea that the value of the loss function on the set of clean points may be lower than that on the set of corrupted points. Thus, an ITLM-based procedure selects a low-loss subset of the training data and performs a model update step on this subset. This alternating minimization is repeated until the model loss is sufficiently small. The heuristic behind ITLM holds in practice, as per the evaluations from \cite{Shen2018-jx}.

Finally, in a more theoretical context, a number of works study classification under various malicious noise models. See the chapter by \citet{bn21} or \cite{dk_book} for an overview.

\smallpar{Memorization of training data.} \citet{Arpit2017-yz} and \citet{Feldman2020-ex} discuss the ability of neural networks to memorize their training data. Specifically, the work of \cite{Arpit2017-yz} empirically discusses how memorization plays into the learning dynamics of neural networks via fitting random labels. The work of \cite{Feldman2020-ex} experimentally validates the ``long tail theory'', which posits that data distributions in practice tend to have a large fraction of their mass allocated to ``atypical'' examples; thus, the memorization of these rare examples is actually necessary for generalization.

Our notion of memorization is different in that we consider excess capacity \textit{on top of the learning problem at hand}. In other words, we require that there exist a classifier in the hypothesis class that behaves correctly on on-distribution data in addition to memorizing specially curated off-distribution data.

\section{Restatement of theorems and full proofs}
\label{sec:app_proof}

In this section, we restate our main results and give full proofs.

\subsection{Proofs from \Cref{sec:backdoor_statistical}}

\thmExistenceLinearBackdoor*
\begin{proof}[Proof of \Cref{thm:existence_linear_backdoor}]
We first show that our choice of $\fadv$ is valid. Observe that for our choice of $\hstar$, we have with probability $1$ that
\begin{align*}
    \hstar\inparen{\patch{x}} &= \signv{\sum_{i \in U} a_i \cdot \vv_i(\patch{x})} \\
    &= \signv{\sum_{i \in U} a_i \cdot \vv_i(x)} \quad \text{, since } \prvv{x\sim\cD}{\vv_i(\patch{x}) = \vv_i(x)} = 1 \\
    &= \hstar\inparen{x},
\end{align*}
as desired.

Next, we pick a patch function from $\fadv$ and construct our stamped training set. We will show that there exists a classifier in $\cH$ achieving $0$ training error on the corrupted training set $\Sclean \cup \Sbackdoor$. Without loss of generality, assume that $\prvv{x\sim\cD}{\vv_i(\patch{x}) \ge 0} = 1$, for all $i \in \insquare{\vdim{\cV}}$ (as negating the basis vector does not affect the space spanned).

Let $a_i$, $i \in U$ denote the coefficients corresponding to $\hstar$, let $R = \max_{x \in \Sclean} \abs{\vstar(x)}$, and let $j$ be any index in $\ucomplement$. We will show that the function
\begin{align*}
    \hhat(x) \coloneqq \signv{\sum_{i \in U} a_i \cdot \vv_i(x) + t \cdot \frac{2R}{\min_{x \in \Sbackdoor} \vv_j(\patch{x})} \cdot \vv_j(x)}
\end{align*}
achieves $0$ training error on the  corrupted training set $S = \Sclean \cup \Sbackdoor$.
We first evaluate $\hhat(x)$ on $x \sim \cD$. We have with probability $1$ that
\begin{align*}
    \hhat(x) &= \signv{\sum_{i \in U} a_i \cdot \vv_i(x) + t \cdot \frac{2R}{\min_{x \in \Sbackdoor} \vv_j(\patch{x})} \cdot \vv_j(x)} \\
    &= \signv{\sum_{i \in U} a_i \cdot \vv_i(x)} \quad \text{ since } \prvv{x \sim \cD}{\vv_j(x) = 0} = 1 \\
    &= \hstar(x).
\end{align*}
We now evaluate $\hhat(\patch{x})$ on $x \sim \cD | \hstar(x) \neq t$. We have with probability $1$ that
\begin{align*}
    \hhat(x) &= \signv{\sum_{i \in U} a_i \cdot \vv_i(\patch{x}) + t \cdot \frac{2R}{\min_{x \in \Sbackdoor} \vv_j(\patch{x})} \cdot \vv_j(\patch{x})} \\
    &= \signv{\sum_{i \in U} a_i \cdot \vv_i(x) + t \cdot \frac{2R}{\min_{x \in \Sbackdoor} \vv_j(\patch{x})} \cdot \vv_j(\patch{x})} \\
    &= \signv{\vstar(x) + t \cdot \frac{2R}{\min_{x \in \Sbackdoor} \vv_j(\patch{x})} \cdot \vv_j(\patch{x})}.
\end{align*}
We now have to show that regardless of the $x \in \cX$ we pick, we have that the above expression evaluates to $t$ with probability $1$. We have two cases. In the first case, $\hstar(x) = 1$, which implies that $\vstar(x) \le R$. Since $t = -1$, we write
\begin{align*}
    \hhat(\patch{x}) &= \signv{\vstar(x) + t \cdot \frac{2R}{\min_{x \in \Sbackdoor} \vv_j(\patch{x})} \cdot \vv_j(\patch{x})} \\
    &= \signv{R- \frac{2R}{\min_{x \in \Sbackdoor} \vv_j(\patch{x})} \cdot \vv_j(\patch{x})} \\
    &= \signv{1 - \frac{2}{\min_{x \in \Sbackdoor} \vv_j(\patch{x})} \cdot \vv_j(\patch{x})} = -1,
\end{align*}
where the last line follows from the fact that  $1 - \frac{2}{\min_{x \in \Sbackdoor} \vv_j(\patch{x})} \cdot \vv_j(\patch{x}) \le -1$.

The proof in the event where $t = 1$ follows similarly. It thus follows that there exists a function $\hhat \in \cH$ obtaining $0$ error on $\Sclean \cup \Sbackdoor$. The desired result immediately follows from uniform convergence (we have a $0$-error classifier over $\Sclean \sim \cD$ and a $0$-error classifier over $\Sbackdoor \sim \patch{\cD | \hstar(x) \neq t}$, so with probability $1 - 2\delta$, we have error at most $\epsclean$ on the clean distribution and error at most $\epsadv$ on the adversarial distribution) This completes the proof of \Cref{thm:existence_linear_backdoor}.
\end{proof}

\corLinearBackdoor*
\begin{proof}[Proof of \Cref{cor:linear_backdoor}]
We will show that our problem setup is a special case of that considered in \Cref{thm:existence_linear_backdoor}; then, we can apply that result as a black box.

Observe that the set of linear classifiers over $\R^d$ is a thresholded vector space with dimension $d$. Pick the orthonormal basis $\inbraces{\vv_1, \dots, \vv_s, \dots, \vv_d}$ such that $\inbraces{\vv_1, \dots, \vv_s}$ form a basis for the subspace $\vspan{\mA}$ and $\vv_{s + 1}, \dots, \vv_d$ are some completion of the basis for the rest of $\R^d$. 

Clearly, there is a size-$s$ set of indices $U \subset [d]$ such that for all $i \in U$, we have $\prvv{x \sim \cD}{\vv_i(x) \neq 0} > 0$. Without loss of generality, assume $U = [s]$.

Next, we need to show that for all $i \in U$, we have $\vv_i(\patch{x}) = 0$. Since we have $\eta \perp \vspan{\mA}$, we have $\vv_i(\eta) = 0$ for all $i \in U$. Since the $v_i$ are also linear functions, we satisfy $\vv_i(\mA\vz + \eta) = 0$ for all $\vz \in \R^s$. 

We now show that there is at least one $j \in \ucomplement$ such that $\prvv{x \sim \cD}{\vv_j(\patch{x}) \neq 0} = 1$. Since $\eta \perp \vspan{\mA}$, $\eta$ must be expressible as some nonzero linear combination of the vectors $\vv_j$; thus, taking the inner product with any such vector will result in a nonzero value.

Finally, we show that for all $j \in \ucomplement$, we either have $\prvv{x \sim \cD}{\vv_j(\patch{x}) \ge 0} = 1$ or $\prvv{x \sim \cD}{\vv_j(\patch{x}) \le 0} = 1$. Since $\eta$ is expressible as a linear combination of several such $v_j$, we can write
\begin{align*}
    \ip{\mA\vz + \eta, \vv_j} &= \ip{\mA\vz, \vv_j} + \ip{\eta, \vv_j} = 0 + \ip{\sum_{j = s + 1}^d a_j \cdot \vv_j, \vv_j} = a_j,
\end{align*}
which is clearly nonzero.

The statement of \Cref{cor:linear_backdoor} now follows from Theorem \ref{thm:existence_linear_backdoor}.
\end{proof}

\thmRandomStamp*
\begin{proof}[Proof of \Cref{thm:random_stamp}]
We prove \Cref{thm:random_stamp} in two parts. We first show that although the adversary does not know $\fadv(\hstar)$, they can find $\patchw \in \fadv(\hstar)$ with high probability. We then invoke the result from Corollary \ref{cor:linear_backdoor}.

Let $\mA$ be such that the columns $\va_1,\dots,\va_s$ are an orthonormal basis for the subspace spanned by $\mA$. Draw the random vector $\vg \sim \mathsf{Unif}(\S^{d-1})$. First, recall that $\abs{\ip{\va_i, \vg}}$ is subgaussian and therefore $\exv{\abs{\ip{\va_i, \vg}}^2} \lesssim \nfrac{1}{d}$ (see \cite[Theorem 3.4.6 and Proposition 2.7.1]{vershynin_2018}). Now, observe that
\begin{align*}
    \exv{\norm{\mA^{\top}\vg}_2} \le \sqrt{\exv{\norm{\mA^{\top}\vg}_2^2}} = \sqrt{\sum_{i=1}^s \exv{\abs{\ip{\va_i, \vg}}^2}} \lesssim \sqrt{\frac{s}{d}}.
\end{align*}
Next, observe that the function $\vx \mapsto \norm{\mA^{\top}\vx}_2$ is $1$-Lipschitz. Using \cite[Theorem 5.1.4]{vershynin_2018}, we have with probability $\ge 1-\delta/2$ that
\begin{align*}
    \abs{\norm{\mA^{\top}\vg}_2 - \exv{\norm{\mA^{\top}\vg}_2}} \lesssim \sqrt{\frac{\logv{\nfrac{1}{\delta}}}{d}}.
\end{align*}
Combining, we get with probability $\ge 1-\delta/2$ that
\begin{align*}
    \norm{\mA^{\top}\vg}_2 \lesssim \sqrt{\frac{s}{d}}+\sqrt{\frac{\logv{\nfrac{1}{\delta}}}{d}},
\end{align*}
which means that as long as $d \gtrsim \frac{s + \logv{\nfrac{1}{\delta}}}{\gamma^2}$, and if we choose $\eta = \vg$, we have
\begin{align*}
    \norm{\mA^{\top}\eta}_2 \le \gamma.
\end{align*}
This implies that the norm of the component of the trigger in $\kernel{\mA^{\top}}$ is at least $\sqrt{1-\gamma^2} \ge 1 - \gamma$ from the Pythagorean Theorem.

This implies that $\hstar(\vx + \eta) = \hstar(\vx)$ with probability $1-\nfrac{\delta}{2}$ over the draws of $\eta$. This gives us $\patch{\vx} = \vx + \eta \in \fadv(\hstar)$ with probability $1-\nfrac{\delta}{2}$ over the draws of $\eta$.

It is now easy to see that \Cref{thm:random_stamp} follows from a simple application of \Cref{cor:linear_backdoor} using a failure probability of $\nfrac{\delta}{2}$, and the final failure probability $1-\delta$ follows from a union bound.
\end{proof}

\thmScaleCapacity*
\begin{proof}[Proof of \Cref{thm:scale_capacity}]
We will proceed by contradiction.

Let $\patch{\vx}$ denote the patched version of $\vx$. Without loss of generality, let the target label be $+1$. Set $\epsclean$ and $\epsadv$ such that $\epsclean + \epsadv < 1$ and draw enough samples such that the attack succeeds with parameters $\epsadv$ and $\delta$.

Observe that we can write every member in $\Sbackdoor$ as $(\patch{\vx}, y)$ for some natural $\vx$ with label $\neg y$. Next, suppose that the learner recovers a $\what$ such that the empirical margin loss of $\what$ is $0$. Next, recall that the following holds for $\what$ obtained from the minimization in the theorem statement and for a training set $S \sim \cD^m$ (see, for instance, Theorem 26.12 of \cite{Shalev-Shwartz2014-oj}):
$$\exvv{(\vx,y)\sim\cD}{\indicator{y\ip{\what, \vx} < 1}} \le \inf_{\norm{\vw}_2 \le \gamma^{-1}}\exvv{(\vx,y)\sim S}{\indicator{y\ip{\vw,\vx} < 1}} + O\inparen{\sqrt{\frac{\inparen{\nfrac{R}{\gamma}}^2 + \logv{\nfrac{1}{\delta}}}{m}}}$$
Using this, we see from uniform convergence that with probability $1-\delta$,
\begin{align*}
    \prvv{\vx \sim \cD}{y\ip{\what, \vx} \ge 1} &\ge 1 - \epsclean \\
    \prvv{\vx \sim \cD}{\ip{\what, \patch{\vx}} \ge 1} &\ge 1 - \epsadv.
\end{align*}
Using a union bound gives
\begin{align*}
    \prvv{\vx \sim \cD}{\inparen{y\ip{\what, \vx} \ge 1} \wedge \inparen{\ip{\what, \patch{\vx}} \ge 1}} \ge 1 - \epsclean - \epsadv.
\end{align*}
Hence, it must be the case that there exists at least one true negative $\vx$ for which both $y\ip{\what, \vx} \ge 1$ and $\ip{\what, \patch{\vx}} \ge 1$ hold. We will use this to obtain a lower bound on $\norm{\what}_2$, from which a contradiction will follow. Notice that
\begin{align*}
    1 &\le \ip{\what, \patch{\vx}} = \ip{\what, \vx} + \ip{\what, \patch{\vx} - \vx} \le -1 + \norm{\what}_2 \cdot \norm{\patch{\vx} - \vx}_2,
\end{align*}
where the last line follows from the fact that $x$ is labeled differently from $\patch{x}$. This gives
$$\norm{\what}_2 \ge \frac{2}{\norm{\patch{\vx} - \vx}_2}.$$
Assuming that we meet the constraint $\norm{\what}_2 \le \nfrac{1}{\gamma}$, putting the inequalities together gives
\begin{align*}
    \norm{\patch{\vx} - \vx}_2 \ge 2\gamma.
\end{align*}
This is a contradiction, since we require that the size of the perturbation is smaller than the margin. This completes the proof of \Cref{thm:scale_capacity}.
\end{proof}

\lmMemVsVC*
\begin{proof}[Proof of \Cref{lemma:mem_vs_vc}]
The lower bound is obvious. This is also tight, as we can set $\cX = \inbraces{0,1}^n$, $\cD = \mathsf{Unif}(\cX)$, and $\cH = \inbraces{f \suchthat f(x) = 1, \forall x \in \cX}$.

We now tackle the upper bound. Suppose for the sake of contradiction that $\memcap{\cX,\cD}{\cH} \ge \vcdim{\cH} + 1$. Then, we can find $k = \vcdim{\cH} + 1$ nonempty subsets of $\cX$, $X_1, \dots, X_k$ and an $h$ for which every labeling of these subsets can be achieved by some other $\hhat \in \cH$. Hence, picking any collection of points $x_i \in X_i$ yields a set witnessing $\vcdim{\cH} \ge k = \vcdim{\cH} + 1$, which is clearly a contradiction.

The upper bound is tight as well. Consider the dataset $S = \inbraces{0, \ve_1, \dots, \ve_d}$, let $\cD$ be a distribution assigning a point mass of $1$ to $\vx = 0$, and let $\hstar(0) = 1$. It is easy to see that the class of origin-containing halfspaces can memorize every labeling $\ve_1, \dots, \ve_d$ as follows -- suppose we have labels $b_1, \dots, b_d$. Then, the classifier
$$\indicator{\sum_{i = 1}^d b_i \cdot \vx_i \ge 0}$$
memorizes every labeling of $\ve_1, \dots, \ve_d$ while correctly classifying the pair $(0, 1)$. Hence, we can memorize $d$ irrelevant sets, which is equal to the VC dimension of origin-containing linear separators. This concludes the proof of \Cref{lemma:mem_vs_vc}.
\end{proof}

\thmMcapAttack*
\begin{mainthm}[Generalization of Theorem \ref{thm:mcap_attack}]
\label{thm:mcap_attack_general}
Pick an array of $k$ target labels $t \in \inbraces{\pm 1}^k$. Suppose we have a hypothesis class $\cH$, a target function $\hstar$, a domain $\cX$, a data distribution $\cD$, and a class of patch functions $\fadv$. Define:
\begin{align*}
    \cC(\fadv(\hstar))_{t'} \coloneqq \{\patch{\suppv{\cD | \hstar(x) \neq t'}} \suchthat &\patchw \in \fadv\}
\end{align*}
and let:
$$\cC(\fadv(\hstar)) \coloneqq \cC(\fadv(\hstar))_{-1} \cup \cC(\fadv(\hstar))_{1}$$
Now, suppose that $\memcap{\cX,\cD}{\hstar,\cH,\cC(\fadv(\hstar))} \ge k$. Then, there exists $k$ functions $\patchw_1, \dots, \patchw_k \in \fadv$ for which the adversary can draw sets $\inbraces{(\Sbackdoor)_i}_{i \in [k]}$ each consisting of $m_i = \Omega\inparen{\epsadv^{-1}\inparen{\vcdim{\cH} + \logv{\nfrac{k}{\delta}}}}$ i.i.d samples from $\cD | \hstar(x) \neq t_i$ such that with probability at least $1-\delta$ over the draws of $(\Sbackdoor)_i$, the adversary achieves the objectives of Problem \ref{problem:adv_general}, regardless of the number of samples the learner draws from $\cD$ for $\Sclean$.
\end{mainthm}
\begin{proof}[Proof of \Cref{thm:mcap_attack_general}]
As per the theorem statement, we can draw $m$ samples from $\cD | \hstar(x) \neq t_i$ to form $\Sbackdoor$ by inverting the labels of the samples we draw. 

Since $\memcap{\cX,\cD}{\hstar, \cH, \cC(\fadv(\hstar))} = k$, there must exist $k$ sets $X_1, \dots, X_k \in \cC(\fadv(\hstar))$ such that the $X_i$ are memorizable, for which we can write $X_i \subseteq \patchw_i\inparen{\suppv{\cD | \hstar(x) \neq t_i}}$ for appropriate choices of $\patchw_i$, and for which $\mu_{\patch{\cD | \hstar(x) \neq t_i}}(X_i) = 1$. This implies that with probability $1$, there exists at least one function $\hhat \in \cH$ such that $\hhat$ returns $t_i$ on every element in $(\Sbackdoor)_i$ for all $i \in [k]$ and agrees with $\hstar$ on every element in the clean training set $\Sclean$.

Thus, we can recover a classifier $\hhat$ from $\cH$ with $0$ error on the training set $\Sclean \cup \inparen{\bigcup_{i \in [k]} (\Sbackdoor)_i}$. In particular, notice that we achieve $0$ error on $\Sclean$ from distribution $\cD$ and on every $(\Sbackdoor)_i$ from distribution $\patchw_i\inparen{\cD | \hstar(x) \neq t_i}$. From the Fundamental Theorem of PAC Learning \cite{Shalev-Shwartz2014-oj}, it follows that as long as $\abs{\Sclean}$ and $\abs{(\Sbackdoor)_i}$ are each at least $\Omega\inparen{\epsclean^{-1}\inparen{\vcdim{\cH} + \logv{\nfrac{k}{\delta}}}}$ and $\Omega\inparen{\epsadv^{-1}\inparen{\vcdim{\cH} + \logv{\nfrac{k}{\delta}}}}$, respectively, we have that $\hhat$ has error at most $\eps$ on $\cD$ and error at least $1 - \eps$ on $\patchw_i\inparen{\cD | \hstar(x) \neq t_i}$ with probability $1-\delta$ (following from a union bound, where each training subset yields a failure to attain uniform convergence with probability at most $\nfrac{\delta}{(k + 1)}$). This completes the proof of \Cref{thm:mcap_attack_general}.
\end{proof}

\thmMcapZero*
\begin{proof}[Proof of \Cref{thm:mcap_zero}]
The condition in the theorem statement implies that there does not exist an irrelevant set that can be memorized atop any choice of $h \in \cH$.

For the sake of contradiction, suppose that there does exist a target classifier $\hstar$, a function $\patchw \in \fadv$ and a target label $t$ such that for all choices of $\epsclean$, $\epsadv$, and $\delta$, we obtain a successful attack.

Define the set $X \coloneqq \patch{\suppv{\cD | \hstar(x) \neq t}}$; in words, $X$ is the subset of $\cX$ consisting of patched examples that are originally of the opposite class of the the target label. It is easy to see that $X \in \cC$.

We will first show that if $\mu_\cD(X) > 0$, then we obtain a contradiction. Set $0 < \epsadv, \epsclean < \frac{\mu_\cD(X)}{1 + \mu_\cD(X)}$. Since the attack is successful, we must classify at least a $1 - \epsadv$ fraction of $X$ as the target label. Hence, we write
\begin{align*}
    \mu_\cD\inparen{\inbraces{x \in X \suchthat \hhat(x) = t}} &\ge \inparen{1 - \epsadv}\mu_\cD(X) > \frac{1}{1 + \mu_\cD(X)} \cdot \mu_\cD(X) > \epsclean.
\end{align*}
Since the set $\inbraces{x \in X \suchthat \hhat(x) = t}$ is a subset of the region of $\cX$ that $\hhat$ makes a mistake on, we have that $\hhat$ must make a mistake on at least $\epsclean$ measure of $\cD$, which is a contradiction.

Hence, it must be the case that $\mu_\cD(X) = 0$; in other words, $X$ is an irrelevant set. Recall that in the beginning of the proof, we assume there exists a function $\hhat$ that achieves label $t$ on $X$, which is opposite of the value of $\hstar$ on $X$. Since we can achieve both possible labelings of $X$ with functions from $\cH$, it follows that $X$ is a memorizable set, and thus the set $X$ witnesses positive $\memcap{\cX,\cD}{\hstar, \cH, \cC(\fadv(\hstar))}$. This completes the proof of \Cref{thm:mcap_zero}.
\end{proof}

\exLinearBackdoorMcap*
\begin{proof}[Proof of \Cref{ex:linear_backdoor_mcap}]
Let $\wstar$ be the weight vector corresponding to $\hstar$.

Observe that there exists $k \coloneqq d - s$ unit vectors $\vv_1, \dots, \vv_k$ that complete an orthonormal basis from that for $P$ to one for $\R^d$. Next, consider the following subset of $\fadv(\hstar)$.
$$\fadv' \coloneqq \inbraces{\patchw \in \fadv \suchthat \forall i \in [k], \patchw_i\inparen{\vx} = \inparen{\begin{cases} \vx + \eta \cdot t_i\vv_i &, \hstar(x) \neq t_i \\ \vx & \text{otherwise} \end{cases}}}$$

We prove the memorization capacity result by using the images of functions in $\fadv'$. We will show that the function
$$\hhat(x) = \signv{\ip{\wstar + \frac{2R}{\gamma}\sum_{i=1}^k t_i\cdot\frac{\vv_i}{\eta_i}, \vx}}$$
memorizes the $k$ sets $C_i \coloneqq \inbraces{\vx + \eta_i \cdot \vv_i \suchthat \ip{\wstar, \vx} \in \insquare{1, \nfrac{R}{\gamma}} \cup \insquare{-\nfrac{R}{\gamma}, -1}}$. Moreover, observe that the preimages of the $C_i$ have measure $1$ under the conditional distributions $\cD | \hstar(x) \neq t_i$, since the preimages contain the support of these conditional distributions. We now have, for a clean point $x \in P$,
\begin{align*}
    \hhat(\vx) &= \signv{\ip{\wstar + \frac{2R}{\gamma}\sum_{i=1}^k t_i\cdot\frac{\vv_i}{\eta_i}, \vx}} = \signv{\ip{\wstar, \vx} + \frac{2R}{\gamma}\ip{\sum_{i=1}^k t_i\cdot\frac{\vv_i}{\eta_i}, \vx}} \\
    &= \signv{\ip{\wstar, \vx}} = \hstar(\vx),
\end{align*}
and for a corrupted point $\vx + \eta_j \cdot \vv_j$, for $j \in [k]$,
\begin{align*}
    \hhat(\vx) &= \signv{\ip{\wstar + \frac{2R}{\gamma}\sum_{i=1}^k t_i\cdot\frac{\vv_i}{\eta_i}, x + \eta_j \cdot \vv_j}} \\
    &= \signv{\ip{\wstar, \vx + \eta_j \cdot \vv_j} + \frac{2R}{\gamma}\ip{\sum_{i=1}^k t_i\cdot\frac{\vv_j}{\eta_j}, \vx + \eta_j \cdot \vv_j}} \\
    &= \signv{\ip{\wstar, \vx} + \frac{2R}{\gamma}\ip{\sum_{i=1}^k t_i\cdot\frac{\vv_i}{\eta_i}, \vx} + \frac{2R}{\gamma}\ip{\sum_{i=1}^k t_i\cdot\frac{\vv_i}{\eta_i}, \eta_j \cdot \vv_j}} \\
    &= \signv{\insquare{\pm \frac{R}{\gamma}} + t_j \cdot \frac{2R}{\gamma}} = t_j.
\end{align*}
This shows that we can memorize the $k$ sets $C_i$. It is easy to see that $\mu_\cD(C_i) = 0$, so the $C_i$ are irrelevant memorizable sets; in turn, we have that $\memcap{\cX,\cD}{\hstar} \ge k = d - s$, concluding the proof of \Cref{ex:linear_backdoor_mcap}.
\end{proof}

\exLinearBackdoorCvx*

To analyze \Cref{ex:linear_backdoor_cvx}, we need the following intermediate results.

\begin{lemma}
\label{lemma:measure_one_interior_points}
Consider some convex body $K$, a probability measure $\cD$ such that every $\ell_2$ ball of nonzero radius within $K$ has nonzero measure, and some subset $K' \subseteq K$ satisfying $\mu_\cD(K') = 1$. Then, $\conv{K'}$ contains every interior point of $K$.
\end{lemma}
\begin{proof}[Proof of \Cref{lemma:measure_one_interior_points}]
Recall that an interior point is defined as one for which we can find some neighborhood contained entirely within the convex body. Mathematically, $\vx \in K$ is an interior point if we can find nonzero $\delta$ for which $\inbraces{\vz \suchthat \norm{\vx - \vz}_2 \le \delta} \subseteq K$.

For the sake of contradiction, suppose that there exists some interior point $\vx \in K$ that is not contained in $\conv{K'}$. Hence, there must exist a halfspace $H$ with boundary passing through $\vx$ and entirely containing $\conv{K'}$. Furthermore, there must exist a nonzero $\delta$ for which there is an $\ell_2$ ball centered at $\vx$ of radius $\delta$ contained entirely within $K$. Call this ball $B_2(\vx, \delta)$. Thus, the set $K \setminus H$ cannot be in $\conv{K'}$. 

We will now show that $\mu_\cD(K \setminus H) > 0$. Observe that the hyperplane inducing $H$ must cut $B_2(\vx, \delta)$ through an equator. From this, we have that the set $K \setminus H$ contains a half-$\ell_2$ ball of radius $\delta$. It is easy to see that this half-ball contains another $\ell_2$ ball of radius $\nfrac{\delta}{2}$ (call this $B'$), and as per our initial assumption, $B'$ must have nonzero measure.

Thus, we can write $\mu_\cD(K \setminus H) \ge \mu_\cD(B') > 0$. Since we know that $\mu_\cD(\conv{K'}) + \mu_\cD(K \setminus H) \le 1$, it follows that $\mu_\cD(\conv{K'}) < 1$ and therefore $\mu_\cD(K') < 1$, violating our initial assumption that $\mu_\cD(K') = 1$. We thus complete the proof of \Cref{lemma:measure_one_interior_points}.
\end{proof}

\begin{lemma}
\label{lemma:line_intersects_interior_twice}
Let $K$ be a closed compact convex set. Let $\vx_1$ be on the boundary of $K$ and let $\vx_2$ be an interior point of $K$. Then, every point of the form $\lambda \vx_1 + (1 - \lambda) x_2$ for $\lambda \in (0, 1)$ is an interior point of $K$.
\end{lemma}
\begin{proof}[Proof of \Cref{lemma:line_intersects_interior_twice}]
Since $\vx_2$ is an interior point, there must exist an $\ell_2$ ball of radius $\delta$ contained entirely within $K$ centered at $\vx_2$. From similar triangles and the fact that any two points in a convex body can be connected by a line contained in the convex body, it is easy to see that we can center an $\ell_2$ ball of radius $(1 - \lambda)\delta$ at the point $\lambda \vx_1 + (1 - \lambda) \vx_2$ that lies entirely in $K$. We therefore conclude the proof of \Cref{lemma:line_intersects_interior_twice}.
\end{proof}

We are ready to analyze \Cref{ex:linear_backdoor_cvx}.

\begin{proof}[Proof of \Cref{ex:linear_backdoor_cvx}]
Observe that the ambient space is equal to the dimension of $\cX$.

Let $\wstar$ be the weight vector corresponding to the true labeler $\hstar$.

For the sake of contradiction, suppose there exists a classifier $\what$ satisfying
\begin{align*}
    \prvv{\vx\sim\cD}{\signv{\ip{\what,\vx}} = \signv{\ip{\wstar, \vx}}} = 1,
\end{align*}
but there exists a subset $Y \subset \cX$ for which $\signv{\ip{\what, \vx}} \neq \signv{\ip{\wstar, \vx}}$, for all $\vx \in Y$. Such a $Y$ would constitute a memorizable set.

Without loss of generality, let the target label be $-1$; that is, the adversary is converting a set $Y$ whose label is originally $+1$ to one whose label is $-1$. Additionally, without loss of generality, take $\norm{\wstar}_2 = \norm{\what}_2 = 1$. Observe that
\begin{align*}
    Y \subseteq D \coloneqq \inbraces{\vx \in \cX \suchthat \ip{\what, \vx} \le 0 \text{ and } \ip{\wstar, \vx} > 0}.
\end{align*}
For $D$ to be nonempty (and therefore for $Y$ to be nonempty), observe that we require $\what \neq \wstar$ (otherwise, the constraints in the definition of the set $D$ are unsatisfiable). 

\Cref{lemma:measure_one_interior_points} implies that if $Y$ is memorizable, then it must lie entirely on the boundary of the set $\cX_{+} \coloneqq \inbraces{\vx \in \cX \suchthat \ip{\wstar, \vx} > 0}$. To see this, observe that if $\what$ classifies any (conditional) measure-$1$ subset of $\cX_{+}$ correctly, then it must classify the convex hull of that subset correctly as well. This implies that $\what$ must correctly classify every interior point in $\cX_{+}$, and thus, $Y$ must be entirely on the boundary of $\cX_{+}$.

Now, let $\vx_1 \in Y$ and $\vx_2 \in \mathsf{Interior}(\cX_{-})$ where $\cX_{-} = \inbraces{x \in \cX \suchthat \ip{\wstar, \vx} < 0}$. Draw a line from $\vx_1$ to $\vx_2$ and consider the labels of the points assigned by $\what$. Since $\vx_1 \in Y$, we have $\hhat(\vx_1) = -1$, and since $\vx_2 \in \mathsf{Interior}(\cX_{-})$, we have that $\hhat(\vx_2) = -1$ as well. Using \Cref{lemma:line_intersects_interior_twice}, we have that every point on the line connecting $\vx_1$ to $\vx_2$ (except for possibly $\vx_1$) is an interior point to $\cX'$. Since we have that the number of sign changes along a line that can be induced by a linear classifier is at most $1$, we must have that the line connecting $\vx_1$ to $\vx_2$ incurs $0$ sign changes with respect to the classifier induced by $\what$. This implies that the line connecting $\vx_1$ to $\vx_2$ cannot pass through any interior points of $\cX_{+}$. However, the only way that this can happen is if $\ip{\wstar, \vx_1} = 0$, but per our definition of $\cX$, if it is the case that $\ip{\wstar, \vx_1} = 0$, then $\vx_1 \notin \cX$, which is a clear contradiction.

This is sufficient to conclude the proof of \Cref{ex:linear_backdoor_cvx}.
\end{proof}

\exSignChanges*
\begin{proof}[Proof of \Cref{ex:sign_changes}]
Without loss of generality, take $d-s$ to be an even integer.

Let $I_1, \dots, I_{s+1}$ be the intervals associated with $\hstar$. It is easy to see that we can pick a total of $\nfrac{(d-s)}{2}$ points such that the sign of these points can be memorized by some $\hhat$. Since each point we pick within an interval can induce at most $2$ additional sign changes, we have that the resulting function $\hhat$ has at most $s + 2 \cdot \nfrac{(d-s)}{2} \le d$ sign-changes; thus, $\hhat \in \cH_d$. Moreover, the measure of a single point is $0$, and so the total measure of our $\nfrac{(d-s)}{2}$ points is $0$.

Given this, it is easy to find $\fadv$ and corresponding $\cC(\fadv(\hstar))$ for which the backdoor attack can succeed as per Theorem \ref{thm:mcap_attack}, thereby yielding the conclusion of \Cref{ex:sign_changes}.
\end{proof}

\subsection{Proofs from \Cref{sec:backdoor_algorithmic}}

\thmCertifyBackdoor*
\begin{proof}[Proof of \Cref{thm:certify_backdoor}]
See \Cref{alg:certify_backdoor} for the pseudocode of an algorithm witnessing \Cref{thm:filter_to_generalize}. 

\begin{algorithm}[ht!]
\caption{Implementation of an algorithm certifying backdoor corruption\label{alg:certify_backdoor}}
\begin{algorithmic}[1]
    \State \textbf{Input}: Training set $S = \Sclean \cup \Sbackdoor$\\ satisfying $\abs{\Sclean} = \Omega\inparen{\epsclean^{-2}\inparen{\vcdim{\robustloss^\cH} + \logv{\nfrac{1}{\delta}}}}$
    \State Set $\hhat \coloneqq \argmin{h \in \cH} \robustloss(h, S)$
    \State \textbf{Output}: $\hhat$ if $\robustloss(\hhat, S) \le 2\eps$ and reject otherwise
\end{algorithmic}
\end{algorithm}

There are two scenarios to consider.

\paragraph{Training set is (mostly) clean.} Suppose that $S$ satisfies $\min_{h \in \cH} \robustloss(h, S) \lesssim \epsclean$. Since the VC dimension of the loss class $\robustloss^\cH$ is finite, it follows that with finitely many samples, we attain uniform convergence with respect to the robust loss, and we're done; in particular, we can write $\robustloss\inparen{\argmin{h\in\cH}\robustloss(h, S), \cD} \lesssim \epsclean$ with high probability.

\paragraph{Training set contains many backdoored examples.} Here, we will show that with high probability, minimizing $\cL_{\fadv(\hstar)}(\hhat, S)$ over $\hhat$ will result in a nonzero loss, which certifies that the training set $S$ consists of malicious examples.

Suppose that for the sake of contradiction, the learner finds a classifier $\hhat$ such that $\cL_{\fadv(\hstar)}(\hhat, S) \lesssim \epsclean$. Hence, with high probability, we satisfy $\cL_{\fadv(\hstar)}(\hhat, \cD) \lesssim \epsclean$. Since there is a constant measure allocated to each class, we can write
\begin{align*}
    \exvv{(x,y)\sim\cD | y\neq t}{\sup_{\patchw\in\fadv(\hstar)} \indicator{\hhat(\patch{x}) \neq y}} \lesssim \epsclean.
\end{align*}
Furthermore, since we achieved a loss of $0$ on the whole training set, including the subset $\Sbackdoor$, from uniform convergence, we get with high probability that
\begin{align*}
    \exvv{(x,y)\sim\cD | y\neq t}{\indicator{\hhat(\patch{x}) = t}} \ge 1-\epsadv.
\end{align*}
This easily implies
\begin{align*}
    \exvv{(x,y)\sim\cD | y\neq t}{\sup_{\patchw\in\fadv(\hstar)} \indicator{\hhat(\patch{x}) \neq y}} \ge 1 - \epsadv.
\end{align*}
Stitching the inequalities together yields $\epsclean \gtrsim 1 - \epsadv$. This is a contradiction, as we can make $\epsclean$ sufficiently small so as to violate this statement. We obtain \Cref{thm:certify_backdoor} as desired.
\end{proof}

\thmFilterToGeneralize*

We first need the intermediate claim \Cref{claim:robust_losses_close}.

\begin{claim}
\label{claim:robust_losses_close}
For all $h \in \cH$, we have
$$\abs{\robustloss(h, \Sclean) - \robustloss(h, \Sclean')} \le \epsclean.$$
\end{claim}
\begin{proof}[Proof of \Cref{claim:robust_losses_close}]
Let $a, b, c$ be positive numbers. We first write
$$\frac{a}{b} - \max\inbraces{0,\frac{a-c}{b-c}} = \frac{c(b-a)}{b(b-c)} \le \frac{c}{b},$$
which occurs exactly when $c \le a$. In the case where $a \le c$, we have
$$\frac{a}{b} - \max\inbraces{0,\frac{a-c}{b-c}} = \frac{a}{b} \le \frac{c}{b},$$
which gives
$$\frac{a}{b} - \max\inbraces{0,\frac{a-c}{b-c}} \le \frac{c}{b}.$$
Next, consider
$$\min\inbraces{1, \frac{a}{b-c}} - \frac{a}{b}  = \frac{a}{b-c} - \frac{a}{b} = \frac{c}{b} \cdot \frac{a}{b - c} \le \frac{c}{b},$$
which happens exactly when we have $b \ge a + c$. In the other case, we have
$$\min\inbraces{1, \frac{a}{b-c}} - \frac{a}{b}  = 1 - \frac{a}{b} \le \frac{c}{b}.$$
We therefore write
$$\max\inbraces{0, \frac{a-c}{b-c}}, \min\inbraces{1, \frac{a}{b-c}} \in \insquare{\frac{a}{b} \pm \frac{c}{b}}.$$
Now, let $a$ denote the number of samples from $\Sclean$ that $h$ incurs robust loss on, let $b$ be the total number of samples from $\Sclean$, and let $c$ be the number of samples our filtering procedure deletes from $\Sclean$. It is easy to see that $\nfrac{a}{b}$ corresponds $\robustloss(h, \Sclean)$ and that $\robustloss(h, \Sclean') \in \insquare{\max\inbraces{0,\nfrac{(a-c)}{(b-c)}}, \min\inbraces{1,\nfrac{a}{(b-c)}}}$. From our argument above, this means that we must have
$$\robustloss(h, \Sclean') \in \insquare{\robustloss(h, \Sclean) \pm \frac{\epsclean(1-\alpha) m}{(1-\alpha)m}}.$$
Finally,
$$\frac{\epsclean(1-\alpha) m}{(1-\alpha)m} = \epsclean,$$
completing the proof of \Cref{claim:robust_losses_close}.
\end{proof}

We now prove \Cref{thm:filter_to_generalize}.

\begin{proof}[Proof of \Cref{thm:filter_to_generalize}]
See \Cref{alg:generalization} for the pseudocode of an algorithm witnessing the theorem statement. 

\begin{algorithm}[H]
\caption{Implementation of a generalization algorithm given an implementation of a filtering algorithm\label{alg:generalization}}
\begin{algorithmic}[1]
    \State \textbf{Input}: Training set $S = \Sclean \cup \Sbackdoor$\\ satisfying $\abs{\Sclean} = \Omega\inparen{\epsclean^{-2}\inparen{\vcdim{\robustloss} + \logv{\nfrac{1}{\delta}}}}$
    \State Run the filtering algorithm on $S$ to obtain $S'$ satisfying the conditions in the theorem statement
    \State \textbf{Output}: Output $\hhat$, defined as $\hhat \coloneqq \argmin{h\in\cH} \robustloss(h, S')$
\end{algorithmic}
\end{algorithm}

Recall that we have drawn enough samples to achieve uniform convergence (see \cite{Cullina2018-os} and \cite{Montasser2020-ir}); in particular, assuming that our previous steps succeeded in removing very few points from $\Sclean$, then for all $h \in \cH$, we have with probability $1-\delta$ that
$$\abs{\robustloss(h, \cD) - \robustloss(h, \Sclean)} \le \epsclean.$$
Observe that we have deleted at most $m\cdot 2\epsclean$ points from $\Sclean$. Let $\Sclean' \coloneqq S' \cap \Sclean$ (i.e., the surviving members of $\Sclean$ from our filtering procedure). We now use \Cref{claim:robust_losses_close} and triangle inequality to write:
\begin{align*}
    \abs{\robustloss(h, \Sclean') - \robustloss(h, \cD)} \le &\abs{\robustloss(h, \Sclean) - \robustloss(h, \Sclean')} + \\
    &\abs{\robustloss(h, \cD) - \robustloss(h, \Sclean)}\\
    \le &\epsclean
\end{align*}
Next, consider some $\hhat$ satisfying $\robustloss(\hhat, S') \lesssim \epsclean$ (which must exist, as per our argument in Part 3), and observe that, for a constant $C$,
\begin{align*}
    \robustloss(\hhat, S') &\ge (1-C\epsclean)\robustloss(\hhat, S' \cap \Sclean) + C\epsclean\robustloss(\hhat, S' \cap \Sbackdoor) \\
    &\ge (1-C\epsclean)\robustloss(\hhat, \Sclean').
\end{align*}
This means that
\begin{align*}
    \robustloss(\hhat, \Sclean') &\le \frac{\epsclean}{1-C\epsclean} = 2\epsclean\inparen{\frac{1}{1-C\epsclean}} \lesssim \epsclean.
\end{align*}
We now use the fact that $\abs{\robustloss(h, \Sclean') - \robustloss(h, \cD)} \le \epsclean$ to arrive at the conclusion that $\robustloss(h,\cD) \lesssim \epsclean$, which completes the proof of \Cref{thm:filter_to_generalize}.
\end{proof}

\thmGeneralizeToFilter*

We first require the intermediate \Cref{claim:generalize_to_filter_interm}.

\begin{claim}
\label{claim:generalize_to_filter_interm}
The following holds for all $h \in \cH$:
$$\abs{\robustloss(h, \Sclean) - \robustloss(h, \Sclean')} < 3\epsclean$$
\end{claim}
\begin{proof}[Proof of \Cref{claim:generalize_to_filter_interm}]
Recall that in the proof of \Cref{thm:filter_to_generalize}, we showed that for positive numbers $a, b, c$ we have
$$\max\inbraces{0, \frac{a-c}{b-c}}, \min\inbraces{1, \frac{a}{b-c}} \in \insquare{\frac{a}{b} \pm \frac{c}{b}}.$$
Now, let $a$ denote the number of samples from $\Sclean$ that $h$ incurs robust loss on, let $b$ be the total number of samples from $\Sclean$, and let $c$ be the number of samples our filtering procedure deletes from $\Sclean$. It is easy to see that $\nfrac{a}{b}$ corresponds $\robustloss(h, \Sclean)$ and that $\robustloss(h, \Sclean') \in \insquare{\max\inbraces{0,\nfrac{(a-c)}{(b-c)}}, \min\inbraces{1,\nfrac{a}{(b-c)}}}$. From our argument above, this means that we must have
$$\robustloss(h, \Sclean') \in \insquare{\robustloss(h, \Sclean) \pm \frac{2\epsclean m}{(1-\alpha)m}}.$$
Finally,
$$\frac{2\epsclean m}{(1-\alpha)m} = \frac{2\epsclean}{(1-\alpha)} \le \frac{2\epsclean}{\nfrac{5}{6}} < 3\epsclean,$$
completing the proof of \Cref{claim:generalize_to_filter_interm}.
\end{proof}

We now have the tools we need to prove \Cref{thm:generalize_to_filter}.

\begin{proof}[Proof of \Cref{thm:generalize_to_filter}]
See \Cref{alg:filtering} for the pseudocode of an algorithm witnessing the theorem statement. 

At a high level, our proof proceeds as follows. We first show that the partitioning step results in partitions that do not have too high of a fraction of outliers, which will allow us to call the filtering procedure without exceeding the outlier tolerance. Then, we will show that the hypotheses $\hhat_L$ and $\hhat_R$ mark most of the backdoor points for deletion while marking only few of the clean points for deletion. Finally, we will show that although $\hhat$ is learned on $S'$ that is not sampled i.i.d from $\cD$, $\hhat$ still generalizes to $\cD$ without great decrease in accuracy.

\begin{algorithm}[H]
\caption{Implementation of a filtering algorithm given an implementation of a generalization algorithm\label{alg:filtering}}
\begin{algorithmic}[1]
    \State \textbf{Input}: Training set $S = \Sclean \cup \Sbackdoor$\\ satisfying $\abs{\Sclean} = \Omega\inparen{\epsclean^{-2}\inparen{\vcdim{\robustloss} + \logv{\nfrac{1}{\delta}}}}$
    \State Calculate $\hhat = \argmin{h \in \cH} \robustloss(h, S)$ and early-return $S$ if $\robustloss(\hhat, S) \le C\epsclean$, for some universal constant $C$
    \State Randomly partition $S$ into two equal halves $S_L$ and $S_R$
    \State Run the generalizing algorithm to obtain $\hhat_L$ and $\hhat_R$ using training sets $S_L$ and $S_R$, respectively
    \State Run $\hhat_L$ on $S_R$ and mark every mistake that $\hhat_L$ makes on $S_R$, and similarly for $\hhat_R$
    \State Remove all marked examples to obtain a new training set $S' \subseteq S$
    \State \textbf{Output}: $S'$ such that $\hhat = \argmin{h \in \cH} \robustloss(h, S')$ satisfies $\robustloss(\hhat, \cD) \lesssim \epsclean$ with probability $1-\delta$
\end{algorithmic}
\end{algorithm}

We have two cases to consider based on the number of outliers in our training set. Let $m$ be the total number of examples in our training set.

\paragraph{Case 1 -- $\alpha m \le \max\inbraces{\nfrac{2}{3\epsclean} \cdot \logv{\nfrac{1}{\delta}}, 24\logv{\nfrac{2}{\delta}}}$}

It is easy to see that $\cL(\hstar, S) \le \alpha$. Using this, we have
\begin{align*}
    \cL(\hstar, S) \le \alpha \frac{2}{3\epsclean \cdot m} \cdot \logv{\frac{1}{\delta}} \lesssim \frac{\epsclean}{\vcdim{\cH} + \logv{\nfrac{1}{\delta}}} \cdot \logv{\frac{1}{\delta}} < \epsclean,
\end{align*}
which implies that we exit the routine via the early-return. From uniform convergence, this implies that with probability $1-\delta$ over the draws of $S$, we have 
\begin{align*}
    \robustloss\inparen{\argmin{h \in \cH} \robustloss\inparen{h, S'}, \cD} \lesssim \epsclean.
\end{align*}
In the other case, we write
\begin{align*}
    \cL(\hstar, S) \le \alpha \le \frac{24\logv{\nfrac{2}{\delta}}}{m} \lesssim \frac{\epsclean^2\logv{\nfrac{1}{\delta}}}{\vcdim{\cH} + \logv{\nfrac{1}{\delta}}} \lesssim \epsclean^2 \le \epsclean,
\end{align*}
and the rest follows from a similar argument.

\paragraph{Case 2 -- $\alpha m \ge \max\inbraces{\nfrac{2}{3\epsclean} \cdot \logv{\nfrac{1}{\delta}}, 24\logv{\nfrac{2}{\delta}}}$}

Let $\tau = \delta$; we make this rewrite to help simplify the various failure events.

\paragraph{Part 1 -- Partitioning does not affect outlier balance.} Define indicator random variables $X_i$ such that $X_i$ is $1$ if and only if example $i$ ends up in $S_R$. We want to show that
$$\prv{\sum_{i \in \Sbackdoor} X_i \notin \insquare{0.5, 1.5}\alpha\cdot\nfrac{m}{2}} \le \tau.$$
Although the $X_i$ are not independent, they are negatively associated, so we can still use the Chernoff Bound and the fact that the number of outliers $\alpha m \ge 24\logv{\nfrac{2}{\tau}}$:
\begin{align*}
    \prv{\sum_{i \in \Sbackdoor} X_i \notin \insquare{0.5, 1.5}\alpha\cdot\nfrac{m}{2}} &\le 2\expv{-\frac{\nfrac{\alpha}{2} \cdot m \cdot \nfrac{1}{4}}{3}} \le 2\expv{-\frac{\alpha m}{24}} \le \tau
\end{align*}
Moreover, if $S_L$ has a $\insquare{\nfrac{\alpha}{2}, \nfrac{3\alpha}{2}}$ fraction of outliers, then it also follows that $S_R$ has a $\insquare{\nfrac{\alpha}{2}, \nfrac{3\alpha}{2}}$ fraction of outliers. Thus, this step succeeds with probability $1-\tau$.

\paragraph{Part 2 -- Approximately correctly marking points.} We now move onto showing that $\hhat_L$ deletes most outliers from $S_R$ while deleting few clean points. Recall that $\hhat_L$ satisfies $\robustloss(\hhat_L, \cD) \le \epsclean$ with probability $1-\delta$. Thus, we have that $\hhat_L$ labels the outliers as opposite the target label with probability at least $1-\epsclean$. Since we have that the number of outliers $\alpha m \ge \nfrac{2}{3\epsclean} \cdot \logv{\nfrac{1}{\tau}}$, we have from Chernoff Bound that (let $X_i$ be the indicator random variable that is $1$ when $\hhat_L$ classifies a backdoored example as the target label)
\begin{align*}
    \prv{\sum_{i \in \Sbackdoor \cap S_R} X_i \ge 2 \cdot \inparen{\epsclean \cdot \frac{3}{2}\alpha m}} \le \expv{-\epsclean \cdot \frac{3}{2}\alpha m} \le \tau.
\end{align*}
Thus, with probability $1-2\tau$, we mark all but at most $\epsclean \cdot 6\alpha m$ outliers across both $S_R$ and $S_L$; since we impose that $\alpha \lesssim 1$, we have that we delete all but a $c\epsclean$ fraction of outliers for some universal constant $c$.

It remains to show that we do not delete too many good points. Since $\hhat_L$ has true error at most $\epsclean$ and using the fact that $m(1-\nfrac{\alpha}{2}) \ge m(1-\alpha) \ge m\alpha \ge \frac{2\logv{\nfrac{1}{\tau}}}{\epsclean}$, from the Chernoff Bound, we have (let $X_i$ be the indicator random variable that is $1$ when $\hhat_L$ misclassifies a clean example)
\begin{align*}
    \prv{\sum_{i \in \Sclean \cap S_R} X_i \ge 2\cdot\inparen{\epsclean \cdot \inparen{1-\nfrac{\alpha}{2}} \cdot \frac{m}{2}}} \le \expv{-\epsclean \cdot \inparen{1-\nfrac{\alpha}{2}} \cdot \frac{m}{2}} \le \tau.
\end{align*}
From a union bound over the runs of $\hhat_L$ and $\hhat_R$, we have that with probability $1-2\tau$, we mark at most $2m\epsclean\cdot\inparen{1-\nfrac{\alpha}{2}} \le 2m\epsclean$ clean points for deletion. From a union bound, we have that this whole step succeeds with probability $1-4\tau-2\delta$.

\paragraph{Part 3 -- There exists a low-error classifier.} At this stage, we have a training set $S'$ that has at least $m(1-2\epsclean)$ clean points and at most $\epsclean \cdot 6\alpha m$ outliers. Recall that $\hstar$ incurs robust loss on none of the clean points and incurs robust loss on every outlier. This implies that $\hstar$ has empirical robust loss at most
\begin{align*}
    \frac{\epsclean \cdot 6\alpha m}{m(1-2\epsclean)} = \frac{6\alpha\epsclean}{1-2\epsclean} \le 2\epsclean,
\end{align*}
where we use the fact that we pick $\epsclean \le \nfrac{1}{10} < \nfrac{1}{4}$ and $\alpha \le \nfrac{1}{6}$. From this, it follows that $\hhat = \argmin{h \in \cH} \robustloss(h, S')$ satisfies $\robustloss(\hhat, S') \le 2\epsclean$.

\paragraph{Part 4 -- Generalizing from $S'$ to $\cD$.} We now have to argue that $\robustloss(\hhat, S') \le 2\epsclean$ implies $\robustloss(\hhat, \cD) \lesssim \epsclean$. Recall that we have drawn enough samples to achieve uniform convergence (see \cite{Cullina2018-os} and \cite{Montasser2020-ir}); in particular, assuming that our previous steps succeeded in removing very few points from $\Sclean$, then for all $h \in \cH$, we have with probability $1-\delta$ that
$$\abs{\robustloss(h, \cD) - \robustloss(h, \Sclean)} \le \epsclean.$$
Observe that we have deleted at most $m\cdot 2\epsclean$ points from $\Sclean$. Let $\Sclean' \coloneqq S' \cap \Sclean$ (i.e., the surviving members of $\Sclean$ from our filtering procedure). We now use \Cref{claim:generalize_to_filter_interm} and triangle inequality to write
\begin{align*}
    \abs{\robustloss(h, \Sclean') - \robustloss(h, \cD)} \le &\abs{\robustloss(h, \Sclean) - \robustloss(h, \Sclean')} + \\
    &\abs{\robustloss(h, \cD) - \robustloss(h, \Sclean)}\\
    < &\quad4\epsclean.
\end{align*}
Next, consider some $\hhat$ satisfying $\robustloss(\hhat, S') \le 2\epsclean$ (which must exist, as per our argument in Part 3), and observe that
\begin{align*}
    \robustloss(\hhat, S') &\ge (1-2\epsclean)\robustloss(\hhat, S' \cap \Sclean) + 2\epsclean\robustloss(\hhat, S' \cap \Sbackdoor) \\
    &\ge (1-2\epsclean)\robustloss(\hhat, \Sclean').
\end{align*}
This implies
\begin{align*}
    \robustloss(\hhat, \Sclean') &\le \frac{2\epsclean}{1-2\epsclean} = 2\epsclean\inparen{\frac{1}{1-2\epsclean}} \le \frac{5\epsclean}{2}.
\end{align*}
We now use the fact that $\abs{\robustloss(h, \Sclean') - \robustloss(h, \cD)} < 4\epsclean$ to arrive at the conclusion that $\robustloss(h,\cD) < \nfrac{13}{2}\cdot\epsclean$, which is the statement of \Cref{thm:generalize_to_filter}.

The constants in the statement of \Cref{thm:generalize_to_filter} follow from setting $\tau = \delta$.
\end{proof}

\newpage

\section{Numerical trials}
\label{sec:backdoor_numerical_trials}

In this section, we present a practical use case for \Cref{thm:certify_backdoor}.

Recall that, at a high level, \Cref{thm:certify_backdoor} states that under certain assumptions, minimizing robust loss on the corrupted training set will either:
\begin{enumerate}
    \item Result in a low robust loss, which will imply from uniform convergence that the resulting classifier is robust to adversarial (and therefore backdoor) perturbations.
    \item Result in a high robust loss, which will be noticeable at training time.
\end{enumerate}
This suggests that practitioners can use adversarial training on a training set which may be backdoored and use the resulting robust loss value to make a decision about whether to deploy the classifier. To empirically validate this approach, we run this procedure (i.e., some variant of \Cref{alg:certify_backdoor}) on the MNIST handwritten digit classification task\footnote{We select MNIST because one can achieve a reasonably robust classifier on the clean version of the dataset. This helps us underscore the difference between the robust loss at train time with and without backdoors in the training set. Moreover, this allows us to explore a setting where our assumptions in \Cref{thm:certify_backdoor} might not hold -- in particular, it's not clear that we have enough data to attain uniform convergence for the binary loss and $\robustloss$, and it's not clear how to efficiently minimize $\robustloss$.}(see \cite{lecun-mnisthandwrittendigit-2010}). Here, the learner wishes to recover a neural network robust to small $\ell_\infty$ perturbations and where the adversary is allowed to make a small $\ell_\infty$-norm watermark.

\paragraph{Disclaimers.} As far as we are aware, the MNIST dataset does not contain personally identifiable information or objectionable content. The MNIST dataset is made available under the terms of the Creative Commons Attribution-Share Alike 3.0 License.

\paragraph{Reproducibility.} We have included all the code to generate these results in the supplementary material. Our code can be found at \url{https://github.com/narenmanoj/mnist-adv-training}.\footnote{Some of our code is derived from the GitHub repositories \url{https://github.com/MadryLab/backdoor_data_poisoning} and \url{https://github.com/skmda37/Adversarial_Machine_Learning_Tensorflow}.}. Our code is tested and working with TensorFlow 2.4.1, CUDA 11.0, NVIDIA RTX 2080Ti, and the Google Colab GPU runtime.

\subsection{MNIST using neural networks}
\label{subs:exp_mnist}

\subsubsection{Scenario}

Recall that the MNIST dataset consists of $10$ classes, where each corresponds to a handwritten digit in $\inbraces{0, \dots, 9}$. The classification task here is to recover a classifier that, upon receiving an image of a handwritten digit, correctly identifies which digit is present in the image.

In our example use case, an adversary picks a target label $t \in \inbraces{0, \dots, 9}$ and a small additive watermark. If the true classifier is $\hstar(\vx)$, then the adversary wants the learner to find a classifier $\hhat$ maximizing $\prvv{\vx\sim\cD | \hstar(\vx) \neq t}{\hhat(\vx) = t}$. In other words, this can be seen as a ``many-to-one'' attack, where the adversary is corrupting examples whose labels are not $t$ in order to induce a classification of $t$. The adversary is allowed to inject some number of examples into the training set such that the resulting fraction of corrupted examples in the training set is at most $\alpha$.

We will experimentally demonstrate that the learner can use the intuition behind \Cref{thm:certify_backdoor} to either recover a reasonably robust classifier or detect the presence of significant corruptions in the training set. Specifically, the learner can optimize a proxy for the robust loss via adversarial training using $\ell_\infty$ bounded adversarial examples, as done by \cite{Madry2017-ep}.

\paragraph{Instantiation of relevant problem parameters.} Let $\cH$ be the set of neural networks with architecture as shown in Table \ref{table:architecture}. Let $\cX$ be the set of images of handwritten digits; we represent these as vectors in $\insquare{0,1}^{784}$. We define $\fadv$ as
$$\inbraces{\patch{\vx} \suchthat \norm{\vx - \patch{\vx}}_\infty \le 0.3 \text{ and } \patch{\vx} - \vx = \mathsf{pattern}},$$
where $\mathsf{pattern}$ is the shape of the backdoor (we use an ``X'' shape in the top left corner of the image, inspired by \cite{Tran2018-bf}). We let the maximum $\ell_\infty$ perturbation be at most $0.3$ since this parameter has been historically used in training and evaluating robust networks on MNIST (see \cite{Madry2017-ep}). In our setup, we demonstrate that these parameters suffice to yield a successful backdoor attack on a vanilla training procedure (described in greater detail in a subsequent paragraph).

Although it is not clear how to efficiently exactly calculate and minimize $\robustloss$, we will approximate $\robustloss$ by calculating $\ell_\infty$-perturbed adversarial examples using a Projected Gradient Descent (PGD) attack. To minimize $\robustloss$, we use adversarial training as described in \cite{Madry2017-ep}. Generating \Cref{table:mnist_full} takes roughly 155 minutes using our implementation of this procedure with TensorFlow 2.4.1 running on the GPU runtime freely available via Google Colab. We list all our relevant optimization and other experimental parameters in \Cref{table:exp_params}.

\begin{table}[h]
\caption{Neural network architecture used in experiments. We implemented this architecture using the Keras API of TensorFlow 2.4.1.\label{table:architecture}}
\centering
\begin{tabular}{|l|l|}
\hline
\multicolumn{1}{|c|}{\textbf{Layer}}   & \multicolumn{1}{c|}{\textbf{Parameters}}                                                                                       \\ \hline
\texttt{Conv2D}       & \texttt{filters=32}, \texttt{kernel\_size=(3,3)},\texttt{activation='relu'} \\
\texttt{MaxPooling2D} & \texttt{pool\_size=(2,2)}                                                                                     \\
\texttt{Conv2D}       & \texttt{filters=64},\texttt{kernel\_size=(3,3)},\texttt{activation='relu'}  \\
\texttt{Flatten}      &                                                                                                                                \\
\texttt{Dense}        & \texttt{units=1024},\texttt{activation='relu'}                                               \\
\texttt{Dense}        & \texttt{units=10},\texttt{activation='softmax'}                                              \\ \hline
\end{tabular}
\end{table}

\begin{table}[h]
\caption{Experimental hyperparameters. We made no effort to optimize these hyperparameters; indeed, many of these are simply the default arguments for the respective TensorFlow functions.\label{table:exp_params}}
\centering
\begin{tabular}{|l|l|}
\hline
\multicolumn{1}{|c|}{\textbf{Property}} & \multicolumn{1}{c|}{\textbf{Details}}                                       \\ \hline
Epochs                                  & 2                                                                           \\
Validation Split                        & None                                                                        \\
Batch Size                              & 32                                                                          \\
Loss                                    & Sparse Categorical Cross Entropy                                            \\
Optimizer                               & RMSProp (step size = $0.001$, $\rho$ = 0.9, momentum = 0, $\eps = 10^{-7}$) \\
NumPy Random Seed                       & 4321                                                                        \\
TensorFlow Random Seed                  & 1234                                                                        \\
PGD Attack                              & $\eps = 0.3$, step size = $0.01$, iterations = $40$, restarts = $10$        \\ \hline
\end{tabular}
\end{table}

\paragraph{Optimization details.} See \Cref{table:exp_params} for all relevant hyperparameters and see \Cref{table:architecture} for the architecture we use.

For the ``Vanilla Training'' procedure, we use no adversarial training and simply use our optimizer to minimize our loss directly. For the ``PGD-Adversarial Training'' procedure, we use adversarial training with a PGD adversary.

In our implementation of adversarial training, we compute adversarial examples for each image in each batch using the PGD attack and we minimize our surrogate loss on this new batch. This is sufficient to attain a classifier with estimated robust loss of around $0.08$ on an uncorrupted training set. 

\subsubsection{Goals and evaluation methods}

We want to observe the impact of adding backdoor examples and the impact of running adversarial training on varied values of $\alpha$ (the fraction of the training set that is corrupted). 

To do so, we fix a value for $\alpha$ and a target label $t$ and inject enough backdoor examples such that exactly an $\alpha$ fraction of the resulting training set contains corrupted examples. Then, we evaluate the train and test robust losses on the training set with and without adversarial training to highlight the difference in robust loss observable to the learner. As sanity checks, we also include binary losses and test set metrics. For the full set of metrics we collect, see \Cref{table:mnist_full}.

To avoid out-of-memory issues when computing the robust loss on the full training set (roughly $60000$ training examples and their adversarial examples), we sample $5000$ training set examples uniformly at random from the full training set and compute the robust loss on these examples. By Hoeffding's Inequality \cite{vershynin_2018}, this means that with probability $0.99$ over the choice of the subsampled training set, the difference between our reported statistic and its population value is at most $\sim 0.02$.

\subsubsection{Results and discussion}

\begin{table}[H]
\caption{Results with MNIST with a target label $t = 0$ and backdoor pattern ``X.'' In each cell, the top number represents the respective value when the network was trained without any kind of robust training, and the bottom number represents the respective value when the network was trained using adversarial training as per \cite{Madry2017-ep}. For example, at $\alpha = 0.05$, for Vanilla Training, the training $0-1$ loss is only $0.01$, but the training robust loss is $1.00$, whereas for PGD-Adversarial Training, the training $0-1$ loss is $0.07$ and the training robust loss is $0.13$. The Backdoor Success Rate is our estimate of $\prvv{x \sim \cD || y \neq t}{\patch{x} = t}$, which may be less than the value of the robust loss.\label{table:mnist_full}}
\centering
\begin{tabular}{|c|c|l|l|l|l|l|}
\hline
\multicolumn{2}{|c|}{$\alpha$}                                    & 0.00 & 0.05 & 0.15 & 0.20 & 0.30 \\ \hline
\multirow{2}{*}{Training $0-1$ Loss}   & Vanilla Training         & 0.01 & 0.01 & 0.01 & 0.01 & 0.01 \\
                                       & PGD-Adversarial Training & 0.02 & 0.07 & 0.17 & 0.22 & 0.33 \\ \hline
\multirow{2}{*}{Training Robust Loss}  & Vanilla Training         & 1.00 & 1.00 & 1.00 & 1.00 & 1.00 \\
                                       & PGD-Adversarial Training & 0.09 & 0.13 & 0.24 & 0.27 & 0.41 \\ \hline
\multirow{2}{*}{Testing $0-1$ Loss}    & Vanilla Training         & 0.01 & 0.01 & 0.01 & 0.02 & 0.01 \\
                                       & PGD-Adversarial Training & 0.02 & 0.03 & 0.03 & 0.03 & 0.06 \\ \hline
\multirow{2}{*}{Testing Robust Loss}   & Vanilla Training         & 1.00 & 1.00 & 1.00 & 1.00 & 1.00 \\
                                       & PGD-Adversarial Training & 0.09 & 0.09 & 0.11 & 0.10 & 0.19 \\ \hline
\multirow{2}{*}{Backdoor Success Rate} & Vanilla Training         & 0.00 & 1.00 & 1.00 & 1.00 & 1.00 \\
                                       & PGD-Adversarial Training & 0.00 & 0.00 & 0.01 & 0.00 & 0.05 \\ \hline
\end{tabular}
\end{table}

See \Cref{table:mnist_full} for sample results from our trials. Over runs of the same experiment with varied target labels $t$, we attain similar results; see \Cref{subs:full_tables} for the full results. We now discuss the key takeaways from this numerical trial.

\paragraph{Training robust loss increases with $\alpha$.} Observe that our proxy for $\robustloss(\hhat, S)$ increases as $\alpha$ increases. This is consistent with the intuition from \Cref{thm:certify_backdoor} in that a highly corrupted training set is unlikely to have low robust loss. Hence, if the learner expects a reasonably low robust loss and fails to observe this during training, then the learner can reject the training set, particularly at high $\alpha$.

\paragraph{Smaller $\alpha$ and adversarial training defeats backdoor.} On the other hand, notice that at smaller values of $\alpha$ (particularly $\alpha \le 0.20$), the learner can still recover a classifier with minimal decrease in robust accuracy. Furthermore, there is not an appreciable decrease in natural accuracy either when using adversarial training on a minimally corrupted training set. Interestingly, even at large $\alpha$, the test-time robust loss and binary losses are not too high when adversarial training was used. Furthermore, the test-time robust loss attained at $\alpha > 0$ is certainly better than that obtained when adversarial training is not used, even at $\alpha = 0$. Hence, although the practitioner cannot certify that the learned model is robust without a clean validation set, the learned model does tend to be fairly robust.

\paragraph{Backdoor is successful with vanilla training.} Finally, as a sanity check, notice that when we use vanilla training, the backdoor trigger induces a targeted misclassification very reliably, even at $\alpha = 0.05$. Furthermore, the training and testing error on clean data is very low, which indicates that the learner would have failed to detect the fact that the model had been corrupted had they checked only the training and testing errors before deployment.

\paragraph{Prior empirical work.} The work of \cite{Borgnia2021-vk} empirically shows the power of data augmentation in defending against backdoored training sets. Although their implementation of data augmentation is different from ours\footnote{Observe that our implementation of adversarial training can be seen as a form of adaptive data augmentation.}, their work still demonstrates that attempting to minimize some proxy for the robust loss can lead to a classifier robust to backdoors at test time. However, our evaluation also demonstrates that classifiers trained using adversarial training can be robust against test-time adversarial attacks, in addition to being robust to train-time backdoor attacks. Furthermore, our empirical results indicate that the train-time robust loss can serve as a good indicator for whether a significant number of backdoors are in the training set.

\newpage

\subsubsection{Results for all target labels}
\label{subs:full_tables}

Here, we present tables of the form of \Cref{table:mnist_full} for all choices of target label $t \in \inbraces{0, \dots, 9}$. Notice that the key takeaways remain the same across all target labels.

\begin{table}[H]
\caption{Results with MNIST with a target label $t = 0$ and backdoor pattern ``X.''}
\centering
\begin{tabular}{|c|c|l|l|l|l|l|}
\hline
\multicolumn{2}{|c|}{$\alpha$}                                    & 0.00 & 0.05 & 0.15 & 0.20 & 0.30 \\ \hline
\multirow{2}{*}{Training $0-1$ Loss}   & Vanilla Training         & 0.01 & 0.01 & 0.01 & 0.01 & 0.01 \\
                                       & PGD-Adversarial Training & 0.02 & 0.07 & 0.17 & 0.22 & 0.33 \\ \hline
\multirow{2}{*}{Training Robust Loss}  & Vanilla Training         & 1.00 & 1.00 & 1.00 & 1.00 & 1.00 \\
                                       & PGD-Adversarial Training & 0.09 & 0.13 & 0.24 & 0.27 & 0.41 \\ \hline
\multirow{2}{*}{Testing $0-1$ Loss}    & Vanilla Training         & 0.01 & 0.01 & 0.01 & 0.02 & 0.01 \\
                                       & PGD-Adversarial Training & 0.02 & 0.03 & 0.03 & 0.03 & 0.06 \\ \hline
\multirow{2}{*}{Testing Robust Loss}   & Vanilla Training         & 1.00 & 1.00 & 1.00 & 1.00 & 1.00 \\
                                       & PGD-Adversarial Training & 0.09 & 0.09 & 0.11 & 0.10 & 0.19 \\ \hline
\multirow{2}{*}{Backdoor Success Rate} & Vanilla Training         & 0.00 & 1.00 & 1.00 & 1.00 & 1.00 \\
                                       & PGD-Adversarial Training & 0.00 & 0.00 & 0.01 & 0.00 & 0.05 \\ \hline
\end{tabular}
\end{table}

\begin{table}[H]
\caption{Results with MNIST with a target label $t = 1$ and backdoor pattern ``X.''}
\centering
\begin{tabular}{|c|c|l|l|l|l|l|}
\hline
\multicolumn{2}{|c|}{$\alpha$}                                    & 0.00 & 0.05 & 0.15 & 0.20 & 0.30 \\ \hline
\multirow{2}{*}{Training $0-1$ Loss}   & Vanilla Training         & 0.01 & 0.01 & 0.01 & 0.01 & 0.01 \\
                                       & PGD-Adversarial Training & 0.02 & 0.07 & 0.17 & 0.23 & 0.32 \\ \hline
\multirow{2}{*}{Training Robust Loss}  & Vanilla Training         & 1.00 & 1.00 & 1.00 & 1.00 & 1.00 \\
                                       & PGD-Adversarial Training & 0.08 & 0.12 & 0.23 & 0.32 & 0.38 \\ \hline
\multirow{2}{*}{Testing $0-1$ Loss}    & Vanilla Training         & 0.01 & 0.01 & 0.01 & 0.01 & 0.01 \\
                                       & PGD-Adversarial Training & 0.02 & 0.02 & 0.03 & 0.04 & 0.05 \\ \hline
\multirow{2}{*}{Testing Robust Loss}   & Vanilla Training         & 1.00 & 1.00 & 1.00 & 1.00 & 1.00 \\
                                       & PGD-Adversarial Training & 0.09 & 0.08 & 0.11 & 0.13 & 0.14 \\ \hline
\multirow{2}{*}{Backdoor Success Rate} & Vanilla Training         & 0.00 & 1.00 & 1.00 & 1.00 & 1.00 \\
                                       & PGD-Adversarial Training & 0.00 & 0.00 & 0.00 & 0.02 & 0.03 \\ \hline
\end{tabular}
\end{table}

\begin{table}[H]
\caption{Results with MNIST with a target label $t = 2$ and backdoor pattern ``X.''}
\centering
\begin{tabular}{|c|c|l|l|l|l|l|}
\hline
\multicolumn{2}{|c|}{$\alpha$}                                    & 0.00 & 0.05 & 0.15 & 0.20 & 0.30 \\ \hline
\multirow{2}{*}{Training $0-1$ Loss}   & Vanilla Training         & 0.01 & 0.01 & 0.01 & 0.01 & 0.00 \\
                                       & PGD-Adversarial Training & 0.02 & 0.07 & 0.17 & 0.22 & 0.32 \\ \hline
\multirow{2}{*}{Training Robust Loss}  & Vanilla Training         & 1.00 & 1.00 & 1.00 & 1.00 & 1.00 \\
                                       & PGD-Adversarial Training & 0.08 & 0.13 & 0.23 & 0.28 & 0.38 \\ \hline
\multirow{2}{*}{Testing $0-1$ Loss}    & Vanilla Training         & 0.01 & 0.02 & 0.01 & 0.02 & 0.01 \\
                                       & PGD-Adversarial Training & 0.02 & 0.03 & 0.03 & 0.03 & 0.05 \\ \hline
\multirow{2}{*}{Testing Robust Loss}   & Vanilla Training         & 1.00 & 1.00 & 1.00 & 1.00 & 1.00 \\
                                       & PGD-Adversarial Training & 0.09 & 0.09 & 0.10 & 0.10 & 0.14 \\ \hline
\multirow{2}{*}{Backdoor Success Rate} & Vanilla Training         & 0.00 & 1.00 & 1.00 & 1.00 & 1.00 \\
                                       & PGD-Adversarial Training & 0.00 & 0.00 & 0.00 & 0.01 & 0.04 \\ \hline
\end{tabular}
\end{table}

\begin{table}[H]
\caption{Results with MNIST with a target label $t = 3$ and backdoor pattern ``X.''}
\centering
\begin{tabular}{|c|c|l|l|l|l|l|}
\hline
\multicolumn{2}{|c|}{$\alpha$}                                    & 0.00 & 0.05 & 0.15 & 0.20 & 0.30 \\ \hline
\multirow{2}{*}{Training $0-1$ Loss}   & Vanilla Training         & 0.01 & 0.01 & 0.01 & 0.01 & 0.01 \\
                                       & PGD-Adversarial Training & 0.02 & 0.07 & 0.18 & 0.23 & 0.32 \\ \hline
\multirow{2}{*}{Training Robust Loss}  & Vanilla Training         & 1.00 & 1.00 & 1.00 & 1.00 & 1.00 \\
                                       & PGD-Adversarial Training & 0.08 & 0.13 & 0.23 & 0.28 & 0.38 \\ \hline
\multirow{2}{*}{Testing $0-1$ Loss}    & Vanilla Training         & 0.01 & 0.01 & 0.01 & 0.02 & 0.02 \\
                                       & PGD-Adversarial Training & 0.02 & 0.02 & 0.03 & 0.04 & 0.05 \\ \hline
\multirow{2}{*}{Testing Robust Loss}   & Vanilla Training         & 1.00 & 1.00 & 1.00 & 1.00 & 1.00 \\
                                       & PGD-Adversarial Training & 0.09 & 0.09 & 0.11 & 0.11 & 0.13 \\ \hline
\multirow{2}{*}{Backdoor Success Rate} & Vanilla Training         & 0.00 & 1.00 & 1.00 & 1.00 & 1.00 \\
                                       & PGD-Adversarial Training & 0.00 & 0.01 & 0.00 & 0.01 & 0.03 \\ \hline
\end{tabular}
\end{table}

\begin{table}[H]
\caption{Results with MNIST with a target label $t = 4$ and backdoor pattern ``X.''}
\centering
\begin{tabular}{|c|c|l|l|l|l|l|}
\hline
\multicolumn{2}{|c|}{$\alpha$}                                    & 0.00 & 0.05 & 0.15 & 0.20 & 0.30 \\ \hline
\multirow{2}{*}{Training $0-1$ Loss}   & Vanilla Training         & 0.01 & 0.01 & 0.01 & 0.01 & 0.01 \\
                                       & PGD-Adversarial Training & 0.02 & 0.07 & 0.17 & 0.22 & 0.32 \\ \hline
\multirow{2}{*}{Training Robust Loss}  & Vanilla Training         & 1.00 & 1.00 & 1.00 & 1.00 & 1.00 \\
                                       & PGD-Adversarial Training & 0.08 & 0.13 & 0.24 & 0.27 & 0.42 \\ \hline
\multirow{2}{*}{Testing $0-1$ Loss}    & Vanilla Training         & 0.01 & 0.01 & 0.01 & 0.01 & 0.01 \\
                                       & PGD-Adversarial Training & 0.02 & 0.02 & 0.03 & 0.03 & 0.05 \\ \hline
\multirow{2}{*}{Testing Robust Loss}   & Vanilla Training         & 1.00 & 1.00 & 1.00 & 1.00 & 1.00 \\
                                       & PGD-Adversarial Training & 0.08 & 0.09 & 0.11 & 0.10 & 0.15 \\ \hline
\multirow{2}{*}{Backdoor Success Rate} & Vanilla Training         & 0.00 & 1.00 & 1.00 & 1.00 & 1.00 \\
                                       & PGD-Adversarial Training & 0.00 & 0.00 & 0.01 & 0.01 & 0.04 \\ \hline
\end{tabular}
\end{table}

\begin{table}[H]
\caption{Results with MNIST with a target label $t = 5$ and backdoor pattern ``X.''}
\centering
\begin{tabular}{|c|c|l|l|l|l|l|}
\hline
\multicolumn{2}{|c|}{$\alpha$}                                    & 0.00 & 0.05 & 0.15 & 0.20 & 0.30 \\ \hline
\multirow{2}{*}{Training $0-1$ Loss}   & Vanilla Training         & 0.01 & 0.01 & 0.01 & 0.01 & 0.01 \\
                                       & PGD-Adversarial Training & 0.02 & 0.07 & 0.17 & 0.22 & 0.33 \\ \hline
\multirow{2}{*}{Training Robust Loss}  & Vanilla Training         & 1.00 & 1.00 & 1.00 & 1.00 & 1.00 \\
                                       & PGD-Adversarial Training & 0.07 & 0.13 & 0.23 & 0.28 & 0.41 \\ \hline
\multirow{2}{*}{Testing $0-1$ Loss}    & Vanilla Training         & 0.01 & 0.01 & 0.01 & 0.02 & 0.02 \\
                                       & PGD-Adversarial Training & 0.02 & 0.03 & 0.03 & 0.03 & 0.06 \\ \hline
\multirow{2}{*}{Testing Robust Loss}   & Vanilla Training         & 1.00 & 1.00 & 1.00 & 1.00 & 1.00 \\
                                       & PGD-Adversarial Training & 0.08 & 0.09 & 0.11 & 0.10 & 0.16 \\ \hline
\multirow{2}{*}{Backdoor Success Rate} & Vanilla Training         & 0.00 & 1.00 & 1.00 & 1.00 & 1.00 \\
                                       & PGD-Adversarial Training & 0.00 & 0.00 & 0.01 & 0.01 & 0.05 \\ \hline
\end{tabular}
\end{table}

\begin{table}[H]
\caption{Results with MNIST with a target label $t = 6$ and backdoor pattern ``X.''}
\centering
\begin{tabular}{|c|c|l|l|l|l|l|}
\hline
\multicolumn{2}{|c|}{$\alpha$}                                    & 0.00 & 0.05 & 0.15 & 0.20 & 0.30 \\ \hline
\multirow{2}{*}{Training $0-1$ Loss}   & Vanilla Training         & 0.01 & 0.01 & 0.01 & 0.01 & 0.01 \\
                                       & PGD-Adversarial Training & 0.02 & 0.07 & 0.17 & 0.22 & 0.33 \\ \hline
\multirow{2}{*}{Training Robust Loss}  & Vanilla Training         & 1.00 & 1.00 & 1.00 & 1.00 & 1.00 \\
                                       & PGD-Adversarial Training & 0.08 & 0.12 & 0.24 & 0.27 & 0.40 \\ \hline
\multirow{2}{*}{Testing $0-1$ Loss}    & Vanilla Training         & 0.01 & 0.02 & 0.01 & 0.01 & 0.01 \\
                                       & PGD-Adversarial Training & 0.02 & 0.03 & 0.03 & 0.03 & 0.06 \\ \hline
\multirow{2}{*}{Testing Robust Loss}   & Vanilla Training         & 1.00 & 1.00 & 1.00 & 1.00 & 1.00 \\
                                       & PGD-Adversarial Training & 0.09 & 0.09 & 0.12 & 0.10 & 0.16 \\ \hline
\multirow{2}{*}{Backdoor Success Rate} & Vanilla Training         & 0.00 & 1.00 & 1.00 & 1.00 & 1.00 \\
                                       & PGD-Adversarial Training & 0.00 & 0.00 & 0.01 & 0.01 & 0.04 \\ \hline
\end{tabular}
\end{table}

\begin{table}[H]
\caption{Results with MNIST with a target label $t = 7$ and backdoor pattern ``X.''}
\centering
\begin{tabular}{|c|c|l|l|l|l|l|}
\hline
\multicolumn{2}{|c|}{$\alpha$}                                    & 0.00 & 0.05 & 0.15 & 0.20 & 0.30 \\ \hline
\multirow{2}{*}{Training $0-1$ Loss}   & Vanilla Training         & 0.01 & 0.01 & 0.01 & 0.01 & 0.01 \\
                                       & PGD-Adversarial Training & 0.02 & 0.07 & 0.18 & 0.22 & 0.32 \\ \hline
\multirow{2}{*}{Training Robust Loss}  & Vanilla Training         & 1.00 & 1.00 & 1.00 & 1.00 & 1.00 \\
                                       & PGD-Adversarial Training & 0.07 & 0.12 & 0.25 & 0.29 & 0.39 \\ \hline
\multirow{2}{*}{Testing $0-1$ Loss}    & Vanilla Training         & 0.01 & 0.01 & 0.01 & 0.02 & 0.01 \\
                                       & PGD-Adversarial Training & 0.02 & 0.03 & 0.03 & 0.03 & 0.04 \\ \hline
\multirow{2}{*}{Testing Robust Loss}   & Vanilla Training         & 1.00 & 1.00 & 1.00 & 1.00 & 1.00 \\
                                       & PGD-Adversarial Training & 0.08 & 0.08 & 0.11 & 0.10 & 0.13 \\ \hline
\multirow{2}{*}{Backdoor Success Rate} & Vanilla Training         & 0.00 & 1.00 & 1.00 & 1.00 & 1.00 \\
                                       & PGD-Adversarial Training & 0.00 & 0.00 & 0.00 & 0.00 & 0.03 \\ \hline
\end{tabular}
\end{table}

\begin{table}[H]
\caption{Results with MNIST with a target label $t = 8$ and backdoor pattern ``X.''}
\centering
\begin{tabular}{|c|c|l|l|l|l|l|}
\hline
\multicolumn{2}{|c|}{$\alpha$}                                    & 0.00 & 0.05 & 0.15 & 0.20 & 0.30 \\ \hline
\multirow{2}{*}{Training $0-1$ Loss}   & Vanilla Training         & 0.01 & 0.01 & 0.01 & 0.01 & 0.01 \\
                                       & PGD-Adversarial Training & 0.02 & 0.07 & 0.17 & 0.22 & 0.32 \\ \hline
\multirow{2}{*}{Training Robust Loss}  & Vanilla Training         & 1.00 & 1.00 & 1.00 & 1.00 & 1.00 \\
                                       & PGD-Adversarial Training & 0.08 & 0.14 & 0.23 & 0.28 & 0.41 \\ \hline
\multirow{2}{*}{Testing $0-1$ Loss}    & Vanilla Training         & 0.01 & 0.01 & 0.01 & 0.01 & 0.01 \\
                                       & PGD-Adversarial Training & 0.02 & 0.03 & 0.03 & 0.03 & 0.05 \\ \hline
\multirow{2}{*}{Testing Robust Loss}   & Vanilla Training         & 1.00 & 1.00 & 1.00 & 1.00 & 1.00 \\
                                       & PGD-Adversarial Training & 0.08 & 0.09 & 0.11 & 0.10 & 0.17 \\ \hline
\multirow{2}{*}{Backdoor Success Rate} & Vanilla Training         & 0.00 & 1.00 & 1.00 & 1.00 & 1.00 \\
                                       & PGD-Adversarial Training & 0.00 & 0.00 & 0.01 & 0.01 & 0.05 \\ \hline
\end{tabular}
\end{table}

\begin{table}[H]
\caption{Results with MNIST with a target label $t = 9$ and backdoor pattern ``X.''}
\centering
\begin{tabular}{|c|c|l|l|l|l|l|}
\hline
\multicolumn{2}{|c|}{$\alpha$}                                    & 0.00 & 0.05 & 0.15 & 0.20 & 0.30 \\ \hline
\multirow{2}{*}{Training $0-1$ Loss}   & Vanilla Training         & 0.01 & 0.01 & 0.01 & 0.01 & 0.01 \\
                                       & PGD-Adversarial Training & 0.02 & 0.07 & 0.17 & 0.22 & 0.33 \\ \hline
\multirow{2}{*}{Training Robust Loss}  & Vanilla Training         & 1.00 & 1.00 & 1.00 & 1.00 & 1.00 \\
                                       & PGD-Adversarial Training & 0.08 & 0.13 & 0.23 & 0.29 & 0.43 \\ \hline
\multirow{2}{*}{Testing $0-1$ Loss}    & Vanilla Training         & 0.01 & 0.01 & 0.01 & 0.01 & 0.01 \\
                                       & PGD-Adversarial Training & 0.02 & 0.03 & 0.03 & 0.04 & 0.06 \\ \hline
\multirow{2}{*}{Testing Robust Loss}   & Vanilla Training         & 1.00 & 1.00 & 1.00 & 1.00 & 1.00 \\
                                       & PGD-Adversarial Training & 0.09 & 0.10 & 0.11 & 0.11 & 0.20 \\ \hline
\multirow{2}{*}{Backdoor Success Rate} & Vanilla Training         & 0.00 & 1.00 & 1.00 & 1.00 & 1.00 \\
                                       & PGD-Adversarial Training & 0.01 & 0.01 & 0.01 & 0.01 & 0.06 \\ \hline
\end{tabular}
\end{table}
\chapter{Spectral clustering in semirandom stochastic block models\label{chapter:cluster}}

In this chapter, we study the robustness of spectral clustering algorithms under helpful model misspecification. This chapter is based on joint work with Aditya Bhaskara, Agastya Vibhuti Jha, Michael Kapralov, Davide Mazzali, and Weronika Wrzos-Kaminska \cite{bjkmmw24}.

\section{Introduction}
\label{sec:cluster_sbm_intro}

Graph partitioning or clustering is a fundamental unsupervised learning primitive. In a graph partitioning problem, one seeks to identify clusters of vertices that are highly internally connected and sparsely connected to the outside. This task is of particular significance when the given graph presents a latent community structure. In this setting, the goal is to recover the communities as accurately as possible. Various statistical models that attempt to capture this situation have been proposed and studied in the literature. Perhaps the most popular of these is the Symmetric Stochastic Block Model (SSBM) \cite{hll83}.

Following the notation of previous works \cite{afwz17,dls20}, in this chapter we describe an SSBM with specifications $n, P_1, P_2, p, q$, where $n$ is an even positive integer, $P_1$ and $P_2$ are a partitioning of the vertex set $V=\inbraces{1,\dots,n}$ into subsets of equal size, and $p$ and $q$ are probabilities. Without loss of generality, we may assume that the partitions $P_1$ and $P_2$ consist of vertices $1,\dots,n/2$ and $n/2+1,\dots,n$, respectively. Hence, with a mild abuse of notation, we write an SSBM with parameters $n, p, q$ only and write it as $\mathsf{SSBM}(n,p,q)$. Now, let $\mathsf{SSBM}(n,p,q)$ be a distribution over random undirected graphs $G=(V,E)$ where each edge $(v,w) \in P_1 \times P_1$ and $(v,w) \in P_2 \times P_2$ (which we refer to as ``internal edges'') appears independently with probability $p$, and each edge $(v,w) \in P_1 \times P_2$ (which we refer to as ``crossing edges'') appears independently with probability $q$. When $p \gg q$, there should be many more internal edges than crossing edges. Hence, we expect the community structure to become more evident as $p$ tends away from $q$. 

In such scenarios, our general algorithmic goal is to efficiently identify $P_1$ and $P_2$ when given $G$ without any community labels. This task is hereafter referred to as the \textit{graph bisection problem}. In this work, we will be interested in \textit{exact recovery}, also known as \textit{strong consistency}, in which we want an algorithm that, with probability at least $1-1/n$ over the randomness of the instance, exactly returns the partition $\{P_1,P_2\}$ for all $n$ sufficiently large. Other approximate notions of recovery (such as almost exact, partial, and weak recovery) are also well-studied but are beyond the scope of this work.

Although the $\mathsf{SSBM}(n,p,q)$ distribution over graphs is a useful starting point for algorithm design and has led to a deep theory about when recovery is possible and of what nature \cite{abbe_survey}, it may not be representative of all scenarios in which we should expect our algorithms to succeed. To remedy this, researchers have proposed several different random graph models that may be more reflective of properties satisfied by real-world networks. These include the geometric block model~\cite{gnw24}, the Gaussian mixture block model \cite{ls24}, and others.

In this chapter, we take a different perspective to graph generation by considering various \textit{semirandom models}. At a high level, a semirandom model for a statistical problem interpolates between an average-case input (for example produced by a model such as the SSBM) and a worst-case input, in a way that still allows for a meaningful notion of ground-truth solution. In our context of graph bisection, this can be achieved by an adversary adding internal edges or by the distribution of internal edges itself being nonhomogeneous (i.e., every internal edge $(v,w)$ appears independently with probability $p_{vw} \ge p$, where the $p_{vw}$ may be chosen adversarially for each internal edge). Researchers have studied similar semirandom models for graph bisection  \cite{fk01,mmv12,mpw16,moitra_sbm_2021, cdm24} and other statistical problems such as classification under Massart noise \cite{mn07}, detecting a planted clique in a random graph \cite{fk01, csv17, mmt20, BKS23}, sparse recovery \cite{kllst23}, and top-$K$ ranking \cite{ycom24}.

These modeling modifications are not necessarily meant to capture a real-world data generation process. Rather, they are a useful testbed with which we can determine whether commonly used algorithms have overfit to statistical assumptions present in the model.  In particular, observe that these changes in model specification are ostensibly helpful, in that increasing the number of internal edges should only enhance the community structure. Perhaps surprisingly, it is known that a number of natural algorithms that succeed in the SSBM setting no longer work under such helpful modifications \cite{moitra_sbm_2021}. Therefore, it is natural to ask which algorithms for graph bisection are robust in semirandom models.

At this point, the performance of approaches based on convex programming is well-understood in various semirandom models \cite{fk01,mmv12,mpw16,moitra_sbm_2021,cdm24}. However, in practice, it is impractical to run such an algorithm due to computational costs. Another class of algorithms, that we call \textit{spectral algorithms}, is more widely used in practice. Loosely speaking, a spectral algorithm constructs a matrix $\mM$ that is a function of the graph $G$ and outputs a clustering arising from the embedding of the vertices determined by the eigenvectors of $\mM$. Popular choices of matrices include the unnormalized Laplacian $\mL_G$ and the normalized Laplacian $\cL_G$ (we will formally define and intuit these notions in the sequel) \cite{vl07}. This is because structural properties of both $\mL_G$ and $\cL_G$ imply that the second smallest eigenvalue of each, denoted as $\lambda_2(\mL_G)$ and $\lambda_2(\cL_G)$, serves as a continuous proxy for connectivity, and the corresponding eigenvector, $\vu_2(\mL_G)$ and $\vu_2(\cL_G)$, has entries whose signs reveal a lot of information about the underlying community structure. This motivates \Cref{alg:spectral_general}. It can be run, for example, with $\mathsf{Matrix}(G) \coloneqq \mL_G$ or $\mathsf{Matrix}(G) \coloneqq \cL_G$. Following this discussion, we arrive at the question we study in this chapter.

\begin{algorithm}[t]
\caption{\textsf{SpectralBisection}: given $G=(V,E)$, outputs a bipartition of $V$}
\label{alg:spectral_general}
\begin{algorithmic}[1]
\Procedure{$\mathsf{SpectralBisection}$}{$G$} \Comment{$G=(V,E)$ is the input graph}
\State $\mM \gets \mathsf{Matrix}(G)$ \Comment{$\mM \in \mathbb{R}^{V \times V}$ is a matrix with real eigenvalues}
\State $((\lambda_i, \vu_i))_{i = 1}^n \gets $ eigenvalue-eigenvector pairs of $\mM$ with $\lambda_1 \le \dots \le \lambda_n$ \Comment{$n = |V|$}
\State $S \gets \{v \in V: \, \vu_2[v] < 0\}$
\State \Return $\{S,V\setminus S\}$
\EndProcedure
\end{algorithmic}
\end{algorithm}

\begin{question}
    \textit{Under which semirandom models do the Laplacian-based spectral algorithms, using the second eigenvector of $\mL_G$ or $\cL_G$, exactly recover the ground-truth communities $P_1$ and $P_2$?}
\end{question}

\smallpar{Main contributions.} Our results show a surprising difference in the robustness of spectral bisection when considering the normalized versus the unnormalized Laplacian. We summarize our results below:
\begin{itemize}
    \item Consider a nonhomogeneous symmetric stochastic block model with parameters $q < p < \pbar$, where every internal edge appears independently with probability $p_{uv} \in [p,\pbar]$ and every crossing edge appears independently with probability $q$. We show that under an appropriate spectral gap condition, the spectral algorithm with the unnormalized Laplacian exactly recovers the communities $P_1$ and $P_2$. Moreover, this holds even if an adversary plants $\ll np$ internal edges per vertex prior to the edge sampling phase.
    \item Consider a stronger semirandom model where the subgraphs on the two communities $P_1$ and $P_2$ are adversarially chosen and the crossing edges are sampled independently with probability $q$. We show that if the graph is sufficiently dense and satisfies a spectral gap condition, then the spectral algorithm with the unnormalized Laplacian exactly recovers the communities~$P_1$~and~$P_2$.  
    \item We show that there is a family of instances from a nonhomogeneous symmetric stochastic block model in which the spectral algorithm achieves exact recovery with the unnormalized Laplacian, but incurs a constant error rate with the normalized Laplacian. This is surprising because it contradicts conventional wisdom that normalized spectral clustering should be favored over unnormalized spectral clustering \cite{vl07}.
\end{itemize}

We also numerically complement our findings via experiments on various parameter settings.

\smallpar{Outline.} The rest of this chapter is organized as follows. In \Cref{sec:cluster_results}, we more formally define our semirandom models, the Laplacians $\mL$ and $\cL$, and formally state our results. In \Cref{sec:cluster_sbm_chapter_overview}, we give sketches of the proofs of our results. In \Cref{sec:cluster_experiments}, we show results from numerical trials suggested by our theory. In \Cref{sec:cluster_concentration,sec:cluster_leave_one_out} we prove important auxiliary lemmas we need for our results. In \Cref{sec:cluster_strong_consistency}, we prove our robustness results for the unnormalized Laplacian. In \Cref{sec:cluster_inconsistency}, we prove our inconsistency result for the normalized Laplacian.
In \Cref{sec:cluster_experiments_more}, we give additional numerical trials and discussion.

\section{Models and main results}
\label{sec:cluster_results}

In this chapter, we study unnormalized and normalized spectral clustering in several semirandom SSBMs. These models permit a richer family of graphs than the SSBM alone.

\smallpar{Matrices related to graphs.} Throughout this chapter, all graphs are to be interpreted as being undirected, and we assume that the vertices of an $n$-vertex graph coincide with  the set $\{1,\dots, n\}$. With this in mind, we begin with defining various matrices associated with graphs, building up to the unnormalized and normalized Laplacians, which are central to the family of algorithms we analyze (\Cref{alg:spectral_general}).

\begin{definition}[Adjacency matrix]
\label{def:adjacency}
Let $G=(V,E)$ be a graph. The adjacency matrix $\mA_G \in \mathbb{R}^{V \times V}$ of $G$ is the matrix with entries defined as $\mA_G[v,w]=\indicator{(v,w) \in E}$.
\end{definition}

\begin{definition}[Degree matrix]
\label{def:degree_matrix}
Let $G=(V,E)$ be a graph.  The degree matrix $\mD_G \in \mathbb{R}^{V \times V}$ of $G$ is the diagonal matrix with entries defined as $\mD_G[v,v]=\vd_G[v]$, where $\vd_G[v]$ is the degree of $v$. 
\end{definition}

\begin{definition}[Unnormalized Laplacian]
\label{def:unnormalized_laplacian}
Let $G=(V,E)$ be a graph. The unnormalized Laplacian $\mL_G \in \mathbb{R}^{V \times V}$  of $G$ is the matrix defined as \smash{$\mL_G \coloneqq \mD_G - \mA_G = \sum_{(v,w) \in E} (\ve_v-\ve_w)(\ve_v-\ve_w)^\top$}, where $\ve_i$ denotes the $i$-th standard basis vector.
\end{definition}

\begin{definition}[Normalized Laplacians]
\label{def:normalized_laplacian}
Let $G=(V,E)$ be a graph. The symmetric normalized Laplacian $\cL_{G,\mathsf{sym}} \in \mathbb{R}^{V \times V}$ and the random walk Laplacian $\cL_{G,\mathsf{rw}} \in \mathbb{R}^{V \times V}$ of $G$ are defined as
\begin{align*}
    \cL_{G, \mathsf{sym}} &\coloneqq \mI-\mD_G^{-1/2}\mA_G\mD_G^{-1/2} \, , & \cL_{G, \mathsf{rw}} &\coloneqq \mI - \mD_G^{-1}\mA_G.
\end{align*}
\end{definition}

For all notions above, when the graph $G$ is clear from context, we omit the subscript $G$. Furthermore, when we discuss normalized Laplacians, we intend its symmetric version $\cL_{\mathsf{sym}}$ unless otherwise stated. So, we omit this subscript as well and simply write $\cL$.

Next, we define the spectral bisection algorithms. We will discuss some intuition for why these algorithms are reasonable heuristics in \Cref{sec:cluster_sbm_chapter_overview}.

\begin{definition}[Unnormalized and normalized spectral bisection]
\label{def:spectral_bisection}
Let $G=(V,E)$ be a graph, and let its unnormalized and normalized Laplacians be $\mL$ and $\cL$, respectively. We refer to the algorithm resulting from running \Cref{alg:spectral_general} on $G$ with $\mathsf{Matrix}(G) \coloneqq \mL_G$ as \textup{unnormalized spectral bisection}. We refer to the algorithm resulting from running \Cref{alg:spectral_general} on $G$ with $\mathsf{Matrix}(G) = \cL_G$ as \textup{normalized spectral bisection}.
\end{definition}

Our goal is to understand when the above algorithms, applied to a graph with a latent community structure, achieve \textit{exact recovery} or \textit{strong consistency}, defined as follows.

\begin{definition}
\label{def:exact_recovery}
Let $\{P_1,P_2\}$ be a partitioning of $V=\{1,\dots,n\}$, and let $\mathcal{D} \coloneqq \mathcal{D}(\{P_1,P_2\})$ be a distribution over $n$-vertex graphs $G=(V,E)$. We say that an algorithm is \textup{strongly consistent} or achieves \textup{exact recovery} on $\cD$ if given a graph $G \sim \cD$ it outputs the correct partitioning $\inbraces{P_1,P_2}$ with probability at least~$1-1/n$ over the randomness of $G$.
\end{definition}

\subsection{Nonhomogeneous symmetric stochastic block model} Our first model is a family of nonhomogeneous symmetric stochastic block models, defined below.

\begin{model}[Nonhomogeneous symmetric stochastic block model]
\label{def:nonhomogeneous}
Let $n$ be an even positive integer, $V = \inbraces{1,\dots,n}$, $\{P_1,P_2\}$ be a partitioning of $V$ into two equally-sized subsets, and $q < p \le \pbar$ be probabilities. Let $\cD$ be any probability distribution over graphs $G=(V,E)$ such that for every $(v,w) \in P_1 \times P_1$ and $(v,w) \in P_2 \times P_2$, the edge $(v,w)$ appears in $E$ independently with some probability $p_{vw} \in [p,\pbar]$, and for every $(v,w) \in P_1 \times P_2$, the edge $(v,w)$ appears in $E$ independently with probability $q$. We call such $\cD$ a nonhomogeneous symmetric stochastic block model (which we will abbreviate as NSSBM). We call the set of all such $\cD$ the family of nonhomogeneous stochastic block models with parameters $p, \pbar, q$, written as $\mathsf{NSSBM}(n,p,\pbar,q)$.
\end{model}

To visualize \Cref{def:nonhomogeneous}, consider the expected adjacency matrix of some NSSBM distribution. We then have the relations 
\begin{align*}
    \left[\begin{array}{ c | c }
    p \cdot \mJ_{n/2} & q \cdot \mJ_{n/2} \\
    \hline
    q \cdot \mJ_{n/2} & p \cdot \mJ_{n/2}
  \end{array}\right] \le \left[\begin{array}{ c | c }
    \mP_{P_1} & q \cdot \mJ_{n/2} \\
    \hline
    q \cdot \mJ_{n/2} & \mP_{P_2}
  \end{array}\right] \le \left[\begin{array}{ c | c }
    \pbar \cdot \mJ_{n/2} & q \cdot \mJ_{n/2} \\
    \hline
    q \cdot \mJ_{n/2} & \pbar \cdot \mJ_{n/2}
  \end{array}\right] \, ,
\end{align*}
where the leftmost matrix denotes the expected adjacency matrix of $\mathsf{SSBM}(n,p,q)$, the rightmost matrix denotes the expected adjacency matrix of $\mathsf{SSBM}(n,\pbar,q)$, $\mJ_k$ denotes the $k \times k$ all-ones matrix, and $\mP_{P_1}$ and $\mP_{P_2}$ denote the edge probability matrices for edges internal to $P_1$ and $P_2$, respectively.

The above also shows that the rank of the expected adjacency matrix for $\mathsf{SSBM}(n,p,q)$ is $2$. However, the rank for the expected adjacency matrix for some NSSBM distribution may be as large as $\Omega(n)$. Perhaps surprisingly, this will turn out to be unimportant for our entrywise eigenvector perturbation analysis. In particular, the tools we use were originally designed for low-rank signal matrices or spiked low-rank signal matrices \cite{afwz17,dls20,bv24}, but we will see that they can be adapted to the signal matrices we consider.

The NSSBM family generalizes the symmetric stochastic block model described in the previous section -- this is attained by setting $p_{vw} = p$ for all internal edges $(v,w)$. However, it can also encode biases for certain graph properties. For instance, a distribution from the NSSBM family may encode the idea that certain subsets of $P_1$ are expected to be denser than $P_1$ as a whole.

With this definition in hand, we are ready to formally state our first technical result in \Cref{mainthm:nonhomogeneous}.

\begin{restatable}{mainthm}{nonhomogeneous}
\label{mainthm:nonhomogeneous}
Let $p,\pbar,q$ be probabilities such that $q < p \le \pbar$ and such that $\alpha \coloneqq \pbar/(p-q)$ is an arbitrary constant. Let $\cD \in \mathsf{NSSBM}(n,p,\pbar,q)$. Let $n \ge N(\alpha)$ where the function $N(\alpha)$ only depends on $\alpha$. There exists a universal constant $C>0$ such that if
\begin{align*}
    n(p-q) \ge C\inparen{\sqrt{n\pbar\log n}+\log n}, \tag{gap condition}
\end{align*}
then unnormalized spectral bisection is strongly consistent on $\cD$.
\end{restatable}

We prove \Cref{mainthm:nonhomogeneous} in \Cref{sec:cluster_proof_nonhomogeneous}. In fact, we show a somewhat stronger statement -- in addition to the process described above, we also allow the adversary to, before sampling the graph, set a small number of the $p_{vw}$ to $1$ (at most $n\pbar/\log \log n$ edges per vertex). We detail this further in \Cref{sec:cluster_proof_nonhomogeneous}.

We now remark on the tightness of our gap condition in \Cref{mainthm:nonhomogeneous}. A work of \citet{abh16} identifies an exact information-theoretic threshold above which exact recovery with high probability is possible and below which no algorithm can be strongly consistent. In particular, the threshold states that for any $p$ and $q$ satisfying \smash{$\sqrt{p}-\sqrt{q} > \sqrt{2\log n / n}$}, exact recovery is possible, and when $p$ and $q$ do not satisfy this, exact recovery is information-theoretically impossible. Furthermore, \citet{fk01} prove that the information-theoretic threshold does not change in a somewhat stronger semirandom model that includes the NSSBM family. Additionally, \citet{dls20} show that unnormalized spectral bisection is strongly consistent all the way to this threshold in the special case where the graph is drawn from $\mathsf{SSBM}(n,p,q)$. By contrast, our gap condition holds in the same critical degree regime as in the information-theoretic threshold (namely, $p =\Theta
(\log n /n$)) but our constant is not optimal. We incur this constant loss because for the sake of presentation, we opt for a cleaner argument that can handle the nonhomogeneity and generalizes more readily across degree regimes. To our knowledge, none of these features are present in prior work analyzing spectral methods in an SSBM setting \cite{afwz17,dls20}.

\subsection{Deterministic clusters model}

Given \Cref{mainthm:nonhomogeneous}, it is natural to ask what happens if we allow the adversary full control over the structure of the graphs in $P_1$ and $P_2$ instead of simply allowing the adversary to perturb the edge probabilities.
In this section, we answer this question. We first describe a more adversarial semirandom model than the NSSBM family. We call this model the \textit{deterministic clusters} model, defined as follows.

\begin{model}[Deterministic clusters model]
\label{def:det_clusters}
Let $n$ be an even positive integer, $V = \inbraces{1,\dots,n}$, $\{P_1,P_2\}$ be a partitioning of $V$ into two equally-sized subsets, $q$ be a probability, and $\din$ be an integer degree lower bound.  Consider a graph $G=(V,E)$ generated according to the following process.
\begin{enumerate}
    \item The adversary chooses arbitrarily graphs $G[P_1]$ and $G[P_2]$ with minimum degree $\din$;
    \item Nature samples every edge $(v,w) \in P_1 \times P_2$ to be in $E$ independently with probability $q$.
    \item The adversary arbitrarily adds edges $(v,w) \in P_1 \times P_1$ and $(v,w) \in P_2 \times P_2$ to $E$ after observing the edges sampled by nature.
\end{enumerate}
We call a distribution $\mathcal{D}$ of graphs generated according to the above process a deterministic clusters model (DCM). We call the set of all such $\mathcal{D}$ the family of deterministic clusters models with parameters $\din$ and $q$, written as $\mathsf{DCM}(n, \din, q)$.
\end{model}

The DCM graph generation process is heavily motivated by the one studied by \citet{mmv12}. This model is much more flexible than the SSBM and NSSBM settings in that the graphs the adversary draws on $P_1$ and $P_2$ are allowed to look very far from random graphs. This means the DCM is a particularly good benchmark for algorithms to ensure they are not implicitly using properties of random graphs that might not hold in the worst case. 

Within the DCM setting, we have \Cref{mainthm:sqrtngap}.

\begin{restatable}{mainthm}{sqrtngap}
\label{mainthm:sqrtngap}
Let $q$ be a probability and $\din$ be an integer, and let $\cD \in \mathsf{DCM}(n, \din, q)$. For $G \sim \cD$, let $\widehat{\mL}$ denote the expectation of $\mL$ after step (2) but before step (3) in \cref{def:det_clusters}. There exists constants $C_1,C_2,C_3>0$ such that for all $n$ sufficiently large, if
\begin{equation*}
    \din \ge C_1 \cdot \inparen{\frac{nq}{2}+\sqrt{n}} \quad \text{and} \quad \lambda_3(\widehat{\mL})-\lambda_2(\widehat{\mL}) \ge \sqrt{n}+C_2nq+C_3\inparen{\sqrt{nq\log n}+\log n} \, ,
\end{equation*}
then unnormalized spectral bisection is strongly consistent on $\cD$.
\end{restatable}

We prove \Cref{mainthm:sqrtngap} in \Cref{sec:cluster_proof_sqrtngap}. We remark that, as in \Cref{mainthm:nonhomogeneous}, the constants that appear in \Cref{mainthm:sqrtngap} are somewhat arbitrary. They are chosen to make our proofs cleaner and can likely be optimized.

As a basic application of \Cref{mainthm:sqrtngap}, note that in the SSBM, if $p = \omega(1/\sqrt{n})$ and $q = 1/\sqrt{n}$, then for $n$ sufficiently large, with high probability, the resulting graph satisfies the conditions needed to apply \Cref{mainthm:sqrtngap}. For a more interesting example, let $P_1$ and $P_2$ be two $d$-regular spectral expanders with $d = \omega(\sqrt{n})$ and let $q \le 1/\sqrt{n}$. On top of both of these two graph classes, one can further allow arbitrary edge insertions inside $P_1$ and $P_2$ while still being guaranteed exact recovery from unnormalized spectral bisection.

\subsection{Inconsistency of normalized spectral clustering}

Notice that in \Cref{mainthm:nonhomogeneous} and \Cref{mainthm:sqrtngap}, we only address the strong consistency of the unnormalized Laplacian in our nonhomogeneous and semirandom models. But what happens when we run spectral bisection with the \textit{normalized} Laplacian?

In \Cref{mainthm:inconsistency}, we prove that there is a subfamily of instances belonging to $\mathsf{NSSBM}(n,p,\pbar,q)$ with $\pbar=6p,q=p/2$ on which unnormalized spectral bisection is strongly consistent (following from \Cref{mainthm:nonhomogeneous}) but normalized spectral clustering is inconsistent in a rather strong sense. Thus, one cannot obtain results similar to \Cref{mainthm:nonhomogeneous} and \Cref{mainthm:sqrtngap} for normalized spectral bisection.

\begin{restatable}{mainthm}{inconsistency}
\label{mainthm:inconsistency}
For all $n$ sufficiently large, there exists a nonhomogeneous stochastic block model such that unnormalized spectral bisection is strongly consistent whereas normalized spectral bisection (both symmetric and random-walk) incurs a misclassification rate of at least $24\%$ with probability $ 1-1/n$.
\end{restatable}

We prove \Cref{mainthm:inconsistency} in \Cref{sec:cluster_inconsistency}. Furthermore, we expect that it is straightforward to adapt the example in \Cref{mainthm:inconsistency} to prove an analogous result for our DCM setting.

The result of \Cref{mainthm:inconsistency} may run counter to conventional wisdom, which suggests that normalized spectral clustering should be favored over the unnormalized variant \cite{vl07}. Perhaps a more nuanced view in light of \Cref{mainthm:nonhomogeneous} and \Cref{mainthm:sqrtngap} is that that the normalized Laplacian and its eigenvectors enjoy stronger concentration guarantees \cite{sb15,dls20}, but the unnormalized Laplacian's second eigenvector is more robust to monotone adversarial changes.

\subsection{Open problems}
Perhaps the most natural follow-up question inspired by our results is to determine whether the restriction that every internal edge probability $p_{vw} \le \pbar$ can be lifted entirely while still maintaining strong consistency of the unnormalized Laplacian (\Cref{mainthm:sqrtngap}). Another exciting direction for future work is to lower the degree and/or spectral gap requirement present in our results in the DCM setting (\Cref{mainthm:sqrtngap}). Finally, we only study insertion-only monotone adversaries, as crossing edge deletions change the second eigenvector of the expected Laplacian. It would be illuminating to understand the robustness of Laplacian-based spectral algorithms against a monotone adversary that is also allowed to delete crossing edges. We are optimistic that the answers to one or more of these questions will further improve our understanding of the robustness of spectral clustering to ``helpful'' model misspecification.

\section{Analysis sketch}
\label{sec:cluster_sbm_chapter_overview}
First, let us give some intuition as to why one may expect that unnormalized spectral bisection is robust against our monotone adversaries. Here and in the sequel, let $\utwostar = [\onev_{n/2} \oplus -\onev_{n/2}]/\sqrt{n}$, where $\onev_{k}$ denotes the all-$1$s vector in $k$ dimensions and $\oplus$ denotes vector concatenation. Let $\mL$ be the unnormalized Laplacian of the graph we want to partition, $\Lstar \coloneqq \exv{\mL}$, $\mE \coloneqq \mL-\Lstar$, and $\lambda_i^{\star} \coloneqq \lambda_i(\Lstar)$ for $1 \le i \le n$. For an edge $(v,w)$, let $\ve_{vw} \coloneqq \ve_{v}-\ve_{w}$, so that $\ve_{vw}$ is an edge incidence vector corresponding to the edge $(v,w)$. Let $p_{vw}$ be the probability that the edge $(v,w)$ appears in $G$ and observe that $\Lstar$ can be written as
\begin{align*}
    \Lstar = \sum_{(v,w) \in E_\text{internal}} p_{vw} \cdot \ve_{vw}\ve_{vw}^T + \sum_{(v,w) \in E_\text{crossing}} q \cdot \ve_{vw}\ve_{vw}^T \, ,
\end{align*}
where $E_\text{internal}=(P_1\times P_1) \cup (P_2\times P_2)$ and $E_\text{crossing}=P_1\times P_2$.
We can verify that $\utwostar$ is an eigenvector of $\Lstar$ -- indeed, we do so in \Cref{lemma:u2expect}. And, for now, assume that $\utwostar$ does correspond to the second smallest eigenvalue of $\Lstar$ (in our NSSBM family, this is easily ensured by enforcing $p > q$). Moreover, for every internal edge $(v,w) \in E_\text{internal}$, we have $\ip{\ve_{vw},\utwostar}=0$. Hence, any changes in internal edges do not change the fact that $\utwostar$ is an eigenvector of the perturbed matrix. Thus, if the sampled $\mL$ is close enough to $\Lstar$, then it is plausible that the second eigenvector of $\mL$, denoted as $\vu_2$, is pretty close to $\utwostar$. In fact, the following conceptually stronger statement holds. If the subgraph formed by selecting just the crossing edges of $G$ is regular, then $\utwostar$ is an eigenvector of $\mL$. This follows from the fact that $\utwostar$ is an eigenvector of the unnormalized Laplacian of any regular bipartite graph where both sides have size $n/2$ and the previous observation that every internal edge is orthogonal to $\utwostar$.

To make this perturbation idea more formal, we recall the Davis-Kahan Theorem. Loosely, it states that $\norm{\vu_2-\utwostar}_2 \lesssim \norm{(\mL-\Lstar)\utwostar}_2/(\lambda_3^{\star}-\lambda_2^{\star})$ (we give a more formal statement in \Cref{lemma:dk_easy}). Expanding the entrywise absolute value $\abs{(\mL-\Lstar)\utwostar}$ reveals that its entries can be expressed as $2\abs{\vdout[v]-\exv{\vdout[v]}}/\sqrt{n}$, where $\vdout[v]$ denotes the number of edges incident to $v$ crossing to the opposite community as $v$. This is unaffected by any increase in the number of edges incident to $v$ that stay within the same community as $v$, denoted as $\vdin[v]$. Hence, regardless of how many internal edges we add before sampling or what substructures they encourage/create, if we have $\lambda_2^{\star} \ll \lambda_3^{\star}$, then we get $\norm{\vu_2-\utwostar}_2 \le o(1)$. This immediately implies that $\vu_2$ is a correct classifier on all but an $o(1)$ fraction of the vertices.

\smallpar{Entrywise analysis of $\vu_2$ and NSSBM strong consistency.} In order to achieve strong consistency, we need that for all $n$ sufficiently large, $\vu_2$ is a perfect classifier. Unfortunately, the above argument does not immediately give that. In particular, in the density and spectral gap regimes we consider, the bound of $o(1)$ yielded by the Davis-Kahan theorem is not sufficiently small to directly yield $\norm{\utwostar-\vu_2}_2 \ll 1/\sqrt{n}$. Instead, we carry out an entrywise analysis of $\vu_2$. A general framework for doing so is given by \citet{afwz17} and is adapted to the unnormalized and normalized Laplacians by \citet{dls20}.

At a high level, we adapt the analysis of \citet{dls20} to our setting. We consider the intermediate estimator vector $\inparen{\mD-\lambda_2\mI}^{-1}\mA\utwostar$. This is a natural choice because we can verify $(\mD-\lambda_2\mI)^{-1}\mA\vu_2=\vu_2$. We will see that it is enough to show that this intermediate estimator correctly classifies all the vertices while satisfying $|(\mD-\lambda_2\mI)^{-1}\mA(\utwostar-\vu_2)| \le |\inparen{\mD-\lambda_2\mI}^{-1}\mA\utwostar|$ (again, the absolute value is taken entrywise). With this in mind, taking some entry indexed by $v \in V$ and multiplying both sides by $\vd[v]-\lambda_2$ (which we will show is positive with high probability), we see that it is enough to show
\begin{align}
     \abs{\ip{\va_v,\utwostar-\vu_2}} \le \abs{\ip{\va_v,\utwostar}} = \frac{\abs{\vdin[v]-\vdout[v]}}{\sqrt{n}},\label{eq:overview_master_lemma}
\end{align}
where $\va_v$ denotes the $v$-th row of $\mA$. The advantage of this rewrite is that the right hand side can be uniformly bounded, so it is enough to control the left hand side. 

To argue about the left hand side of \eqref{eq:overview_master_lemma}, it may be tempting to use the fact that $\va_v$ is a Bernoulli random vector and use Bernstein's inequality to argue about the sum of rescalings of these Bernoulli random variables. Unfortunately, we cannot do this since $\vu_2$ and $\va_v$ are dependent. To resolve this, we use a leave-one-out trick~\cite{afwz17,bv24}. We can think of this as leaving out the vertex $v$ corresponding to the entry we want to analyze and sampling the edges incident to the rest of the vertices. The second eigenvector of the resulting \smash{$\mL^{(v)}$}, denoted as \smash{$\vu_2^{(v)}$}, is a very good proxy for $\vu_2$ and is independent from $\va_v$. Hence, we may complete the proof of \Cref{mainthm:nonhomogeneous}. 

One of our main observations is that although this style of analysis was originally built for low-rank signal matrices \cite{afwz17,bv24}, it can be adapted to handle the nonhomogeneity inside $P_1$ and $P_2$. In particular, the nonhomogeneity we permit in the NSSBM family may make $\Lstar$ look very far from a spiked low-rank signal matrix. Furthermore, our entrywise analysis of eigenvectors under perturbations is one of the first that we are aware of that moves beyond analyzing low-rank signal matrices or spiked low-rank signal matrices.

\smallpar{Extension to deterministic clusters.} To prove \Cref{mainthm:sqrtngap}, we start again at \eqref{eq:overview_master_lemma}. An alternate way to upper bound the left hand side is to use the Cauchy-Schwarz inequality. A variant of the Davis-Kahan theorem gives us control over $\norm{\vu_2-\utwostar}_2$ while $\norm{\va_v}_2 = \sqrt{\vd[v]}$. The advantage of this is that we get a worst-case upper bound on the left hand side of \eqref{eq:overview_master_lemma} -- it holds no matter what edges orthogonal to $\utwostar$ are inserted before or after nature samples the crossing edges (which are precisely the internal edges). Combining these and using the fact that the right hand side of \eqref{eq:overview_master_lemma} is increasing in $\vdin[v]$ (and increases faster than $\norm{\va_v}_2 = \sqrt{\vd[v]}$) allows us to complete the proof of \Cref{mainthm:sqrtngap}.

\smallpar{Inconsistency of normalized spectral bisection.} Finally, we describe the family of hard instances we use to prove \Cref{mainthm:inconsistency}. To motivate this family of instances, recall that by the graph version of Cheeger's inequality, the second eigenvalue of $\cL$ and the corresponding eigenvector can be used to find a sparse cut in $G$. Thus, if we create sparse cuts inside $P_1$ that are sparser than the cut formed by separating $P_1$ and $P_2$, then conceivably the normalized Laplacian's second eigenvector may return the new sparser cut.

To make this formal, consider the following graph structure. Let $n$ be a multiple of $4$. Let $L_1$ consist of indices $1,\dots,n/4$, $L_2$ consist of indices $n/4+1,\dots,n/2$, and $R$ consist of indices $n/2+1,\dots,n$. Consider the block structure induced by the matrix $\Astar = \exv{\mA}$ shown in \Cref{table:hard_a_overview}.

\begin{table}
\centering
\begin{tabular}{l|cccl}
                     & \multicolumn{1}{l}{$L_1$}               & \multicolumn{1}{l}{$L_2$}                                    & \multicolumn{2}{l}{$R$}                                                     \\ \hline
$L_1$                & $Kp \cdot \mathbbm{1}_{n/4 \times n/4}$ & \multicolumn{1}{c|}{$p \cdot \mathbbm{1}_{n/4 \times n/4}$}  & \multicolumn{2}{c}{\multirow{2}{*}{$q \cdot \mathbbm{1}_{n/2 \times n/2}$}} \\
$L_2$                & $p \cdot \mathbbm{1}_{n/4 \times n/4}$  & \multicolumn{1}{c|}{$Kp \cdot \mathbbm{1}_{n/4 \times n/4}$} & \multicolumn{2}{c}{}                                                        \\ \cline{2-5} 
\multirow{2}{*}{$R$} & \multicolumn{2}{c|}{\multirow{2}{*}{$q \cdot \mathbbm{1}_{n/2 \times n/2}$}}                           & \multicolumn{2}{c}{\multirow{2}{*}{$p \cdot \mathbbm{1}_{n/2 \times n/2}$}} \\
                     & \multicolumn{2}{c|}{}                                                                                  & \multicolumn{2}{c}{}                                                       
\end{tabular}
\caption{$\Astar$ for \Cref{mainthm:inconsistency} is defined to have the above block structure.\label{table:hard_a_overview}}
\end{table}

Intuitively, as $K$ gets larger, the cut separating $L_1$ from $V \setminus L_1$ becomes sparser. From Cheeger's inequality, this witnesses a small $\lambda_2(\cL)$ and therefore the corresponding $\vu_2(\cL)$ may return the cut $L_1, V \setminus L_1$. We formally prove that this is indeed what happens when $K$ is a sufficiently large constant and then \Cref{mainthm:inconsistency} follows.

\section{Numerical trials}
\label{sec:cluster_experiments}
We programmatically generate synthetic graphs that help illustrate our theoretical findings using the libraries NetworkX 3.3 (BSD 3-Clause license), SciPy 1.13.0 (BSD 3-Clause License), and NumPy 1.26.4 (modified BSD license) \cite{networkx,scipy,numpy}. We ran all our experiments on a free Google Colab instance with the CPU runtime, and each experiment takes under one hour to run. In this section we focus on a setting that allows relating \cref{mainthm:nonhomogeneous} and \cref{mainthm:inconsistency}, and defer more experiments that investigate both NSSBM and DCM graphs to \cref{sec:cluster_experiments_more}.

To put \cref{mainthm:nonhomogeneous} and \cref{mainthm:inconsistency} in perspective, we consider graphs generated following the process outlined in the proof of \cref{mainthm:inconsistency}, which gives rise to the following benchmark distribution.

\smallpar{Benchmark distribution.} Let $n$ be divisible by $4$ and let $\{P_1,P-2\}$ be a partitioning of $V=[n]$ into two equally-sized subsets. Let $\{L_1,L_2\}$ be a bipartition of $P_1$ such that $|L_1|=|L_2|=n/4$ and call $L=P_1, \, R=P_2$ for convenience as in the proof of \cref{mainthm:inconsistency}. Then, for some $p,\pbar,q \in [0,1]$ such that $q \le p \le \pbar$, consider the distribution $\mathcal{D}_{p,\pbar,q}$ over graphs $G=(V,E)$ obtained by sampling every edge $(u,v) \in (L_1 \times L_1) \cup (L_2 \times L_2)$ independently with probability $\pbar$, every edge $(u,v) \in (L_1 \times L_2) \cup (R \times R)$ independently with probability $p$, and every edge $(u,v) \in L \times R$ independently with probability $q$. One can see that $\mathcal{D}_{p,\pbar,q}$ is in fact in the set $\mathsf{NSSBM}(n,p,\pbar,q)$.

\begin{figure}
    \centering
    \includegraphics[width=0.31\columnwidth]{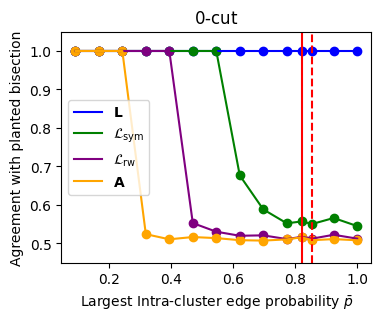}
    \includegraphics[width=0.33\columnwidth]{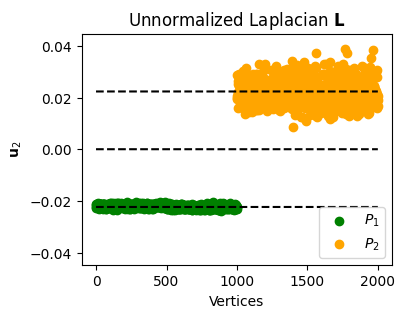}
    \includegraphics[width=0.33\columnwidth]{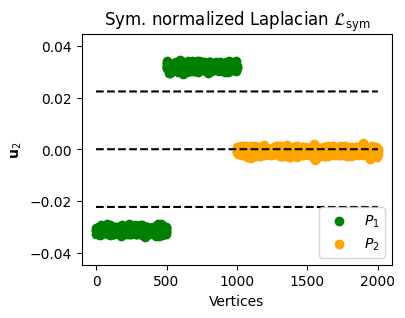}
    \includegraphics[width=0.31\columnwidth]{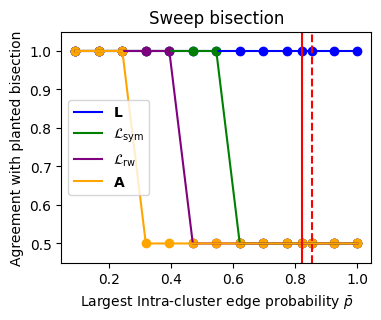}
    \includegraphics[width=0.33\columnwidth]{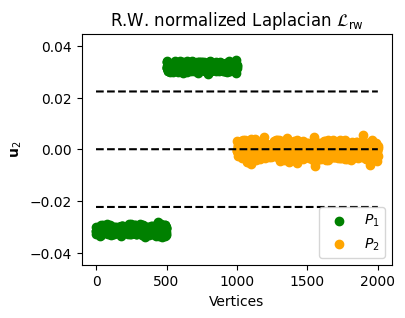}
    \includegraphics[width=0.31\columnwidth]{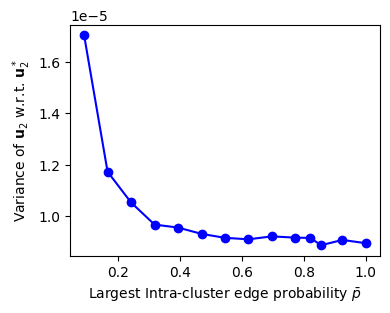}
    \caption{\textbf{Top left, bottom left}: Agreement with the planted bisection of the bipartition obtained from several matrices associated with an input graph generated from a distribution in $\mathsf{NSSBM}(n,p,\pbar,q)$ for fixed values of $n,p,q$ and varying values of $\pbar$. In the top left plot, the bipartition is the $0$-cut of the second eigenvector, as in \cref{alg:spectral_general}. In the bottom left plot, the bipartition is the sweep cut of the first $n/2$ vertices in the second eigenvector. The dashed vertical line corresponds to $\pbar_{\textsf{max}}=\pbar_{\textsf{max}}(n,p,q)$ (see~\eqref{eq:defpbarmax}), and the solid vertical line corresponds to $\pbar_{\textsf{thr}}=\pbar_{\textsf{thr}}(n,p,q)$ (see~\eqref{eq:defpbarthr}). \textbf{Top middle, top right, bottom middle}: Embedding of the vertices given by the second eigenvector $\vu_2$ of several matrices associated with a graph sampled from $\mathcal{D}_{p,\pbar,q}$ with $\pbar = \pbar_{\textsf{thr}}$. Horizontal dashed lines, from top to bottom, correspond to $1/\sqrt{n}, 0, -1/\sqrt{n}$ respectively.\\ \textbf{Bottom right}: Variance of the embedding in the second eigenvector $\vu_2$ of the unnormalized Laplacian with respect to the ideal eigenvector $\utwostar$ (see~\eqref{eq:variance}), for input graphs generated from a distribution in $\mathsf{NSSBM}(n,p,\pbar,q)$ with fixed values of $n,p,q$ and varying values of $\pbar$.}
    \label{fig:embedding_nssbm}
\end{figure}

\smallpar{Setup.} Let us fix $n=2000$, $p=24 \log n /n$, $q = 8 \log n /n$. For varying values of $\pbar$ in the range $[p,1]$, we sample $t=10$ independent draws $G$ from $\mathcal{D}_{p,\pbar,q}$. For each of them, we run spectral bisection (i.e. \cref{alg:spectral_general}) with matrices $\mL,\cL_{\mathsf{sym}},\cL_{\mathsf{rw}},\mA$. Then, we compute the \textit{agreement} of the bipartition hence obtained (with respect to the planted bisection), that is the fraction of correctly classified  vertices. We average the agreement across the $t$ independent draws. The results are shown in the top left plot of \cref{fig:embedding_nssbm}. Another natural way to get a bipartition of $V$ from the eigenvector is a \textit{sweep cut}. In a sweep cut, we sort the entries of $\vu_2$ and take the vertices corresponding to the smallest $n/2$ entries to be on one side of the bisection and put the remaining on the other side. The average agreement obtained in this other fashion is shown in the bottom left plot of \cref{fig:embedding_nssbm}.
\smallpar{Theoretical framing.} As per \cref{mainthm:nonhomogeneous},  we expect unnormalized spectral bisection to achieve exact recovery (i.e. agreement equal to $1$) whenever $\pbar \le \pbar_{\textsf{max}}$, where
\begin{equation}
\label{eq:defpbarmax}
    \pbar_{\textsf{max}} = \frac{\inparen{n(p-q)-\log n}^2}{n \log n}
\end{equation}
is obtained by rearranging the precondition of \cref{mainthm:nonhomogeneous}, ignoring the constants and disregarding the fact that $\alpha$ should be $O(1)$. On the contrary, the proof of \cref{mainthm:inconsistency} shows that normalized spectral bisection misclassifies a constant fraction of vertices provided that $p/q \ge 2$ (which our choice of parameters satisfies) and $\pbar \ge \pbar_{\textsf{thr}}$, where
\begin{equation}
\label{eq:defpbarthr}
    \pbar_{\textsf{thr}} = 3\cdot p^2/q \, .
\end{equation}
In \cref{fig:embedding_nssbm}, the solid vertical line corresponds to the value of $\pbar_{\textsf{thr}}$ on the $x$-axis, and the dashed vertical line corresponds to the value of $\pbar_{\textsf{max}}$ on the $x$-axis. In particular, observe that in our setting $\pbar_{\textsf{thr}}<\pbar_{\textsf{max}}$, so there is an interval of values for $\pbar$ where we expect \cref{mainthm:nonhomogeneous} and \cref{mainthm:inconsistency} to apply simultaneously.

\smallpar{Empirical evidence: consistency.} One can see from the top left plot in \cref{fig:embedding_nssbm} that the agreement of unnormalized spectral bisection is $100\%$ for all values of $\pbar$, even beyond $\pbar_{\textsf{thr}}$ and $\pbar_{\textsf{max}}$. On the other hand, the agreement of the bipartition obtained from all other matrices (hence including normalized spectral bisection) drops below $70\%$ well before the threshold $\pbar_{\textsf{thr}}$ predicted by \cref{mainthm:inconsistency}. From the right plot in \cref{fig:embedding_nssbm}, we see that computing the bipartition by taking a sweep cut of $n/2$ vertices does not change the results -- $\vu_2$ of the unnormalized Laplacian continues to achieve $100\%$ agreement, while for all other matrices the corresponding $\vu_2$ remains inconsistent.

\smallpar{Empirical evidence: embedding variance.} From the setting of the experiment we just illustrated, observe that as we increase $\pbar$, we expect the subgraph $G[L]$ to have increasing volume. As illustrated in \cref{fig:embedding_nssbm}, this seems to correlate with a decrease in the ``variance'' of the second eigenvector $\vu_2$ of the unnormalized Laplacian with respect to the ideal second eigenvector $\utwostar$. More precisely, we compute the average distance squared of the embedding of a vertex in $\vu_2$ from its ideal embedding in $\utwostar$, i.e. the quantity
\begin{equation}
\label{eq:variance}
    \min_{s \in \{\pm 1\}} \frac{1}{n}\norm{\vu_2- s \cdot \utwostar}^2_2 \, .
\end{equation}
This suggests that not only does the second eigenvector of the unnormalized Laplacian remain robust to monotone adversaries, but it actually concentrates more strongly around the ideal embedding~$\utwostar$.

\smallpar{Empirical evidence: example embedding.} Let us fix the value $\pbar = \pbar_{\textsf{thr}}$, for which we see in \cref{fig:varying_clique_cuts} that all matrices except the unnormalized Laplacian fail to recover the planted bisection. We generate a graph from $\mathcal{D}_{p,\pbar,q}$, and plot how the vertices are embedded in the real line by the second eigenvector of all the matrices we consider. The result is shown in \cref{fig:embedding_nssbm}, where the three horizontal dashed lines, from top to bottom, respectively correspond to the value of $1/\sqrt{n}, 0, -1/\sqrt{n}$ on the $y$-axis.

\subsection{Related work}
\label{sec:sbm_related_work}

\smallpar{Community detection.} Community detection has garnered significant attention in theoretical computer science, statistics, and data science. For a general overview of recent progress and related literature, see the survey by \citet{abbe_survey}. In what follows, we discuss the works we believe are most related to what we study in this chapter.

As mentioned in the introduction, perhaps the most fundamental and well-studied model is the symmetric stochastic block model (SSBM), due to \cite{hll83}. The celebrated work of \citet{abh16} gives sharp bounds on the threshold for exact recovery for the SSBM setting. They complement their result by showing that SDP based methods can achieve the information theoretic lower bound for the planted bisection problem, even with a monotone adversary \cite{moitra_sbm_2021}. A line of work \cite{afwz17, dls20} demonstrates that natural spectral algorithms achieve exact recovery for the SSBM all the way to the information-theoretic threshold.

\smallpar{Generalizations of the symmetric stochastic block model.} Since the introduction of SBMs \cite{hll83}, numerous variants have been proposed that are designed to better reflect real-world graph properties. For instance, real-life social networks are likely to contain triangles. To address this, \citet{BS17} introduced a spatial stochastic block model, sometimes known as the geometric stochastic block model (GSBM). Other variations were introduced in the works of \cite{GMPS18, GMPS20}. Subsequent work studies the performance of spectral algorithms on certain Gaussian or Geometric Mixture block models \cite{ABARS20, ABD21, ls24,gnw24}.

Studying community detection with a semirandom model approaches this modeling question differently. Rather than implicitly encouraging a particular structure within the clusters like the models just mentioned, a semirandom adversary (including the ones we study in this chapter) can more directly test the robustness of the algorithm to specially designed substructures.

\smallpar{Semirandom and monotone adversaries.} As far as we are aware, \citet{BS95} were the first to introduce a semirandom model. Within this model, they studied graph coloring problems. \citet{fk01} demonstrated that semidefinite programming methods can accurately recover communities up to a certain threshold, even in the semi-random setting. Other problems, such as detecting a planted clique \cite{Jerrum92, Kuc95, BHK16}, have also been studied in the semi-random model of \cite{fk01}. In the setting of planted clique, a natural spectral algorithm fails against monotone adversaries \cite{mmt20, BKS23}.
Monotone adversaries and semirandom models have also been extensively studied for other statistical and algorithmic problems \cite{AV18, kllst23, GC23, bglmsy23}. Finally, \cite{sl17} shows that a spectral heuristic due to \citet{boppana87} is robust under a monotone adversary that is allowed to both insert internal edges and delete crossing edges. However, as far as we are aware, this algorithm does not fit in the framework of \Cref{alg:spectral_general}.

We remark that the models we study in this chapter are most closely related to models studied by \cite{mn07} and \cite{mmv12}. In particular, allowing increased internal edge probabilities is analogous to Massart noise in classification problems, and our model with adversarially chosen internal edges can be seen as the same model as that studied in \cite{mmv12} (although without allowing crossing edge deletions). Finally, note that \citet{cdm24} give a near-linear time algorithm for graph clustering in the model of \cite{mmv12}, though they do not explicitly show their algorithm is strongly consistent on instances that are information-theoretically exactly recoverable.

\section{Deferred proofs}
\label{sec:cluster_sbm_deferred_proofs}

In this section, we build the tools we need to prove \Cref{mainthm:nonhomogeneous}, \Cref{mainthm:sqrtngap}, and\Cref{mainthm:inconsistency}. Throughout, it will be helpful to refer to the overview (\Cref{sec:cluster_sbm_chapter_overview}) for a proof roadmap.

\smallpar{Notation in the proofs.} In all proofs, we adopt the notation used in the technical overview (\Cref{sec:cluster_sbm_chapter_overview}). Additionally, for a vertex $v \in V$, let $P(v)$ denote the community that $v$ belongs to.

\subsection{Concentration inequalities}
\label{sec:cluster_concentration}
Our proof strategy for \Cref{mainthm:nonhomogeneous} and \Cref{mainthm:sqrtngap} is to appeal to \Cref{lemma:strong_consistency}, which guarantees strong consistency provided that $\vd[v] - \lambda_2 > 0$, $\vdin[v] > \vdout[v]$, and $\abs{\ip{\va_v,\utwostar-\vu_2}} \le (\vdin[v]-\vdout[v])/\sqrt{n}$ for all vertices $v$. Proving that the first two conditions hold is relatively easy. In the setting of \Cref{mainthm:nonhomogeneous}, it essentially follows from concentration of the degrees, which is proved in \Cref{sec:cluster_concentration_degrees}. In the setting of \Cref{mainthm:sqrtngap}, it follows from the assumptions of the Theorem. Proving that the third condition holds is the main technical challenge. 

For all three parts, our proofs rely on several auxiliary concentration results. We prove these in \Cref{sec:cluster_concentration_laplacian} and \Cref{sec:cluster_evect_perturbations}.

We extensively use the following variants of Bernstein's Inequality, which can be derived from \cite[Theorem 2.8.4]{vershynin_2018}.

\begin{lemma}
\label{lemma:chernoff}
Let $X = \sum_{i=1}^m X_i$, where $X_i = 1$ with probability $p_i$ and $X_i = 0$ with probability $1-p_i$ and all the $X_i$ are independent. Let $\mu = \exv{X}$. Then, for all $t > 0$ we have
\begin{align*}
    \prv{\abs{X - \mu} \ge t} &\le 2\expv{-\min\inbraces{\frac{t^2}{4\sum_{i=1}^m p_i(1-p_i)}, \frac{3t}{4}}}.
\end{align*}
\end{lemma}

\noindent
From this, we get the following very useful corollary. 
\begin{lemma}
\label{lemma:chernoff_useful}
Let $X = \sum_{i=1}^m X_i$, where $X_i = 1$ with probability $p_i$ and $X_i = 0$ with probability $1-p_i$ and all the $X_i$ are independent. Let $\mu = \exv{X}$. Then, for all $t > 0$, with probability at least $1-\delta$ we have
\begin{align*}
    \abs{X-\mu} \le \sqrt{4\sum_{i=1}^m p_i(1-p_i)\logv{\nfrac{2}{\delta}}}+\nfrac{4}{3}\logv{\nfrac{2}{\delta}}.
\end{align*}
\end{lemma}

\subsection{Concentration of degrees}\label{sec:cluster_concentration_degrees}

In this Section, we give concentration statements regarding the number of internal vertices incident to each vertex and the number of crossing edges incident to each vertex. We then compare these against $\lambda_2$.

\begin{lemma}
\label{lemma:dout}
Suppose the crossing edges are sampled identically and independently with probability $q$. Then, for some universal constant $C>0$, with probability at least $ 1-\delta$ we have that
\begin{align*}
    \forall v \in V, \quad \abs{\vdout[v]-\exv{\vdout[v]}} \le C\inparen{\sqrt{nq\logv{\nfrac{n}{\delta}}}+\logv{\nfrac{n}{\delta}}}.
\end{align*}
\end{lemma}
\begin{proof}[Proof of \Cref{lemma:dout}]
Choose some $v \in V$. Consider the random variable $\vdout[v]$. Using \Cref{lemma:chernoff_useful}, we have that there is a constant $C>0$ such that with probability at least $ 1-\delta/n$ one has
\begin{align*}
    \abs{\vdout[v]-\exv{\vdout[v]}} \le C\inparen{\sqrt{4nq/2\logv{\nfrac{2n}{\delta}}}+\logv{\nfrac{2n}{\delta}}}.
\end{align*}
Taking a union bound over all $n$ vertices completes the proof of \Cref{lemma:dout}.
\end{proof}
\noindent
Note that \Cref{lemma:dout} above applies in both the settings of \Cref{mainthm:nonhomogeneous} and \Cref{mainthm:sqrtngap}.

\begin{lemma}
\label{lemma:din}
Suppose the internal edges are sampled independently with probabilities $p_{vw}$ such that $p \le p_{vw} \le \pbar$. Then, for some universal constant $C>0$, with probability $\ge 1-\delta$ we have that
\begin{align*}
     \forall v \in V, \quad \abs{\vdin[v]-\exv{\vdin[v]}} \le C\inparen{\sqrt{\sum_{w \in P(v)\setminus \{v\}} p_{vw}(1-p_{vw})\logv{\nfrac{n}{\delta}}}+\logv{\nfrac{n}{\delta}}}.
\end{align*}
\end{lemma}
\begin{proof}[Proof of \Cref{lemma:din}]
As before, choose some $v \in V$ and consider the random variable $\vdin[v]$. By \Cref{lemma:chernoff_useful}, we have that there is a constant $C>0$ such that with probability at least $ 1-\delta/n$ one has
\begin{align*}
    \abs{\vdin[v]-\exv{\vdin[v]}} \le C\inparen{\sqrt{4\sum_{w \in P(v) \setminus \{v\}} p_{vw}(1-p_{vw})\logv{\nfrac{2n}{\delta}}}+\logv{\nfrac{2n}{\delta}}}.
\end{align*}
Taking a union bound over all $n$ vertices completes the proof of \Cref{lemma:din}.
\end{proof}
Combining the above two lemmas, we obtain a lower-bound on $ \vdin[v]-\vdout[v]$. In particular, the following lemma implies that in the setting of \Cref{mainthm:nonhomogeneous}, we have $\vdin[v]>\vdout[v]$. This will be required for applying \Cref{lemma:strong_consistency}. 
\begin{lemma}
\label{lemma:degree_diff_easy}
There exists a universal constant $C>0$ such that with probability $\ge 1-\delta$, in the same settings as \Cref{lemma:dout} and \Cref{lemma:din} and assuming the gap condition in \Cref{mainthm:nonhomogeneous}, if $p \ge q$, then for all $v \in V$ we have
\begin{align*}
    \vdin[v]-\vdout[v] \ge \frac{n(p-q)}{2} - C\inparen{\sqrt{np\logv{\nfrac{n}{\delta}}}+\logv{\nfrac{n}{\delta}}}.
\end{align*}
\end{lemma}
\begin{proof}[Proof of \Cref{lemma:degree_diff_easy}]
Let $v \in V$. First, we call \Cref{lemma:dout} with a failure probability of $\delta/(2n)$ to conclude that
\begin{align*}
    \vdout[v] \le \frac{nq}{2} + C_{\ref{lemma:dout}}\inparen{\sqrt{\frac{nq}{2}\logv{\nfrac{2n^2}{\delta}}}+\logv{\nfrac{2n^2}{\delta}}}.
\end{align*}
Next, we call \Cref{lemma:din} with a failure probability of $\delta/(2n)$ to conclude that
\begin{align*}
    \vdin[v] &\ge \sum_{w \in P(v)\setminus \{v\}} p_{vw} - C_{\ref{lemma:din}}\inparen{\sqrt{\sum_{w \in P(v)\setminus \{v\}} p_{vw}(1-p_{vw})\logv{\nfrac{2n^2}{\delta}}}+\logv{\nfrac{2n^2}{\delta}}} \\
    &\ge \sum_{w \in P(v)\setminus \{v\}} p_{vw} - C_{\ref{lemma:din}}\inparen{\sqrt{\sum_{w \in P(v)\setminus \{v\}} p_{vw}\logv{\nfrac{2n^2}{\delta}}}+\logv{\nfrac{2n^2}{\delta}}} \\
    &\ge \frac{np}{2} - 2C_{\ref{lemma:din}}\inparen{\sqrt{\frac{np}{2}\logv{\nfrac{n^2}{\delta}}}+\logv{\nfrac{2n^2}{\delta}}}.
\end{align*}
where the last line uses the fact that $x - c\sqrt{x}$ is increasing in $x$ whenever $x \ge c^2/4$ and $c>0$. We subtract and conclude the proof of \Cref{lemma:degree_diff_easy} by a union bound.
\end{proof}
The following lemma will be useful for lower-bounding $\vd[v] - \lambda_2$ in \Cref{mainthm:nonhomogeneous}. 
\begin{lemma}
\label{lemma:degree_vs_lambda}
Suppose every crossing edge appears independently with probability $q$. Then, with probability $\ge 1-\delta$, for all $v \in V$ we have
\begin{align*}
    \lambda_2 \le 2\vdout[v] + C\inparen{\sqrt{nq\logv{\nfrac{n}{\delta}}} + \logv{\nfrac{n}{\delta}}}.
\end{align*}
\end{lemma}
\begin{proof}[Proof of \Cref{lemma:degree_vs_lambda}]
Observe that with probability at least $ 1-\delta$, $\vdout[w] - \exv{\vdout[v]} \le \sqrt{2nq \log (\nfrac{2n}{\delta})} + 2\log (2n/\delta)$ for all $w \in V$ by \cref{lemma:chernoff_useful}. Then, for every $v \in V$ we have
\begin{align*}
    \frac{2}{n}\sum_{w \in P(v)} \vdout[w] - \vdout[v] &= \inparen{\frac{2}{n}\sum_{w \in P(v)} \vdout[w] - \exv{\vdout[v]}} + \inparen{\exv{\vdout[v]} - \vdout[v]} \\
    &\le \abs{\frac{2}{n}\sum_{w \in P(v)} \vdout[w] - \exv{\vdout[v]}} + \abs{\exv{\vdout[v]} - \vdout[v]} \\
    &\le \sqrt{2nq\logv{\nfrac{2n}{\delta}}} + \sqrt{2nq\logv{\nfrac{2n}{\delta}}} + 4\logv{\nfrac{2n}{\delta}} \\
    &\le 3\sqrt{nq\logv{\nfrac{n}{\delta}}} + 10\logv{\nfrac{n}{\delta}}.
\end{align*}
Next, by the min-max principle, we have
\begin{align*}
    \lambda_2 \le \sum_{(w, w') \in E} \inparen{\utwostar[w]-\utwostar[w']}^2 = \frac{4}{n}\sum_{w \in P(v)} \vdout[w].
\end{align*}
Combining everything, we get
\begin{align*}
    \lambda_2 \le 2\inparen{\frac{2}{n}\sum_{w \in P(v)} \vdout[w]} \le 2\inparen{\vdout[v] + 3\sqrt{nq\logv{\nfrac{n}{\delta}}}+10\logv{\nfrac{n}{\delta}}},
\end{align*}
completing the proof of \Cref{lemma:degree_vs_lambda}.
\end{proof}
We can now lower-bound $\vd[v] - \lambda_2$. Note that the following lower bound implies that $\vd[v] > \lambda_2$, as required by \Cref{lemma:strong_consistency}. 
\begin{lemma}
\label{lemma:dv_lambda2_positive}
In the setting of \Cref{mainthm:nonhomogeneous}, with probability $\ge 1-\delta$, for all $v \in V$, we have $\vd[v] - \lambda_2 > n(p-q)/4$.
\end{lemma}
\begin{proof}[Proof of \Cref{lemma:dv_lambda2_positive}]
Recall that the gap condition in \Cref{mainthm:nonhomogeneous} tells us that $p$ and $q$ are such that for a universal constant $C$,
\begin{align*}
    n(p-q) \ge C\inparen{\sqrt{np\logv{\nfrac{n}{\delta}}}+\logv{\nfrac{n}{\delta}}}.
\end{align*}
We have for all $n$ sufficiently large (specifically, $n \ge N(\alpha,\delta)$ for some $N$ that is a function only of the constant $\alpha$, and we take $\delta \ge 1/n^{O(1)}$) that with probability at least $1-\delta$,
\begin{align*}
    \vd[v] - \lambda_2 &= \vdin[v]-\vdout[v] + (2\vdout[v]-\lambda_2) \\
    & \ge \vdin[v]-\vdout[v] - C_{\ref{lemma:degree_vs_lambda}}\inparen{\sqrt{nq\logv{\nfrac{n}{\delta}}}+ \logv{\nfrac{n}{\delta}}} \\
    &\ge \frac{n(p-q)}{2} - \inparen{C_{\ref{lemma:degree_diff_easy}}+C_{\ref{lemma:degree_vs_lambda}}}\inparen{\sqrt{np\logv{\nfrac{n}{\delta}}}+\logv{\nfrac{n}{\delta}}},
\end{align*}
so insisting
\begin{align*}
    \frac{n(p-q)}{4} \ge \inparen{C_{\ref{lemma:degree_diff_easy}}+C_{\ref{lemma:degree_vs_lambda}}}\inparen{\sqrt{np\logv{\nfrac{n}{\delta}}}+\logv{\nfrac{n}{\delta}}}+1
\end{align*}
gives the condition required to complete the proof of \Cref{lemma:dv_lambda2_positive}.
\end{proof}
The following technical lemma will be useful for upper-bounding $\maxnorm{\vu_2}$ in \Cref{lemma:utwo_inf_norm_bootstrap}. 
\begin{lemma}
\label{lemma:npbar_vs_denom}
In the setting of \Cref{mainthm:nonhomogeneous}, there exists a universal constant $C$ such that with probability $\ge 1-\delta$, for all $v \in V$ we have
\begin{align*}
    \frac{n\pbar+\logv{\nfrac{n}{\delta}}}{\vd[v]-\lambda_2} \le 4\alpha+C.
\end{align*}
\end{lemma}
\begin{proof}[Proof of \Cref{lemma:npbar_vs_denom}]

By \Cref{lemma:dv_lambda2_positive}, we have with probability $\ge 1-\delta$ that for all $v \in V$,
\begin{align*}
    \vd[v]-\lambda_2 \ge \frac{n(p-q)}{4}.
\end{align*}
 This gives
\begin{align*}
    \frac{n\pbar+\logv{\nfrac{n}{\delta}}}{\vd[v]-\lambda_2} \le \frac{4(n\pbar+\logv{\nfrac{n}{\delta}})}{n(p-q)} = \frac{4\pbar}{p-q} + \frac{4\logv{\nfrac{n}{\delta}}}{n(p-q)} \le 4\alpha + C.
\end{align*}
This completes the proof of \Cref{lemma:npbar_vs_denom}.
\end{proof}

\subsection{Concentration of Laplacian and eigenvalue perturbations}\label{sec:cluster_concentration_laplacian}

For the matrix concentration lemmas, we need a result due to \citet{llv16}. We reproduce it below.

\begin{lemma}[{\cite[Theorem 2.1]{llv16}}]
\label{lemma:llv_concentration}
Consider a random graph from the model $G(n,\inbraces{p_{ij}})$. Let $d = \max_{ij} np_{ij}$. For any $r \ge 1$, the following holds with probability at least $1-n^{-r}$ for a universal constant $C$. Consider any subset consisting of $10n/d$ vertices, and reduce the weights of the edges incident to those vertices in an arbitrary way. Let $d'$ be the maximal degree of the resulting graph. Then, the adjacency matrix $\mA'$ of the new weighted graph satisfies
\begin{align*}
    \opnorm{\mA'-\exv{\mA}} \le Cr^{3/2}\inparen{\sqrt{d} + \sqrt{d'}}.
\end{align*}
Moreover, the same holds for $d'$ being the maximal $\ell_2$ norm of the rows of $\mA'$.
\end{lemma}

\begin{lemma}
\label{lemma:laplacian_concentration}
Let $\mL$ be a Laplacian sampled from the nonhomogeneous Erd\H{o}s-R\'enyi model where each edge $(i,j)$ is present independently with probability $p_{ij}$. Then, there exists a universal constant $C$ such that for all $n$ sufficiently large, with probability $\ge 1-\delta$ for any $\delta \ge n^{-10}$,
\begin{align*}
    \opnorm{\mL-\exv{\mL}} \le C\inparen{\sqrt{n\max_{(i,j)\colon p_{ij}\neq 1} p_{ij}\logv{\nfrac{n}{\delta}}}+\logv{\nfrac{n}{\delta}}}.
\end{align*}
\end{lemma}
\begin{proof}[Proof of \Cref{lemma:laplacian_concentration}]
Without loss of generality, for all $p_{ij}$ that are $1$, reset their probabilities to $0$. To see that this is valid, let $\mL'$ be a Laplacian sampled from this modified distribution and notice that $\mL'-\exv{\mL'}=\mL-\exv{\mL}$.

By \Cref{lemma:llv_concentration} and \cref{lemma:chernoff_useful}, we have with probability $\ge 1-\delta/2$ that
\begin{align*}
    \opnorm{\mA-\exv{\mA}} & \le 200C_{\ref{lemma:llv_concentration}} \sqrt{2n\max_{ij} p_{ij} + C_{\ref{lemma:chernoff_useful}}\inparen{\sqrt{n\max_{ij} p_{ij}\log (\nfrac{8n}{\delta})}+\log (\nfrac{8n}{\delta})}} \\
    & \le 400 C_{\ref{lemma:llv_concentration}} C_{\ref{lemma:chernoff_useful}} \sqrt{n\max_{ij} p_{ij} +\log (\nfrac{8n}{\delta})} \\
    & \le 400 C_{\ref{lemma:llv_concentration}} C_{\ref{lemma:chernoff_useful}} \inparen{\sqrt{n\max_{ij} p_{ij}\log (\nfrac{8n}{\delta})} +\log (\nfrac{8n}{\delta})}\, 
\end{align*}
and by \Cref{lemma:dout} and \Cref{lemma:din}, we have with probability $1-\delta/2$ that
\begin{align*}
    \opnorm{\mD-\exv{\mD}} &\le \max_{v \in V} \abs{\vdout[v]-\exv{\vdout[v]}} + \max_{v \in V} \abs{\vdin[v]-\exv{\vdin[v]}} \\
    &\le 2\max\inbraces{C_{\ref{lemma:dout}},C_{\ref{lemma:din}}}\inparen{\sqrt{n\max_{ij} p_{ij}\logv{\nfrac{2n}{\delta}}}+\logv{\nfrac{2n}{\delta}}}
\end{align*}
Now, observe that with probability $\ge 1-\delta$ (following from a union bound),
\begin{align*}
     \opnorm{\mL-\exv{\mL}} &= \opnorm{\mD-\exv{\mD} - (\mA - \exv{\mA})} \le \opnorm{\mD-\exv{\mD}} + \opnorm{\mA-\exv{\mA}} \\
     &\le 800 C_{\ref{lemma:llv_concentration}} C_{\ref{lemma:chernoff_useful}}\max\inbraces{C_{\ref{lemma:dout}},C_{\ref{lemma:din}}}\inparen{\sqrt{n\max_{ij} p_{ij}\logv{\nfrac{8n}{\delta}}}+\logv{\nfrac{8n}{\delta}}} ,
\end{align*}
completing the proof of \Cref{lemma:laplacian_concentration}.
\end{proof}
By applying the above lemma, we can show that there is a gap between $\lambda_3$ and $\lambda_2^{\star}$, which will allow us to apply Davis-Kahan style bounds. More concretely, \Cref{lemma:lambda3_lambda2star} and \Cref{lemma:eutwostar_norm}, together with \Cref{lemma:utwo_to_utwostar}, show that 
$\norm{u_2 - u_2^{\star}}_2$ is small. This will be useful for proving that in the context for \Cref{mainthm:nonhomogeneous}, the condition $\abs{\ip{\va_v,\utwostar-\vu_2}} \le (\vdin[v]-\vdout[v])/\sqrt{n}$ in \Cref{lemma:strong_consistency} is satisfied. 
\begin{lemma}
\label{lemma:lambda3_lambda2star}
In the setting of \Cref{mainthm:nonhomogeneous}, there exists a universal constant $C$ such that the following holds.

Let $p$ and $q$ be such that we have
\begin{align*}
    n(p-q) \ge C\inparen{\sqrt{n\pbar\logv{\nfrac{n}{\delta}}}+\log(\nfrac{n}{\delta})}.
\end{align*}
Then, for any $\delta \ge n^{-10}$, with probability $\ge 1-\delta$, we have $\lambda_3 - \lambda_2^{\star} \ge n(p-q)/4$.
\end{lemma}
\begin{proof}[Proof of \Cref{lemma:lambda3_lambda2star}]
By Weyl's inequality and \Cref{lemma:laplacian_concentration}, we have with probability $\ge1-\delta$ that
\begin{align*}
    \lambda_3-\lambda_2^{\star} \ge \lambda_3^{\star}-\lambda_2^{\star} - \opnorm{\mL-\Lstar} \ge \frac{n(p-q)}{2} - C_{\ref{lemma:laplacian_concentration}}\inparen{\sqrt{n\pbar\logv{\nfrac{n}{\delta}}}+\log(\nfrac{n}{\delta})}.
\end{align*}
Let $C \ge 4C_{\ref{lemma:laplacian_concentration}}$. Then,
\begin{align*}
    \frac{n(p-q)}{4} \ge C_{\ref{lemma:laplacian_concentration}}\inparen{\sqrt{n\pbar\logv{\nfrac{n}{\delta}}}+\log(\nfrac{n}{\delta})}.
\end{align*}
Subtracting completes the proof of \Cref{lemma:lambda3_lambda2star}.
\end{proof}
Next, we bound $\norm{\mE\utwostar}_2$, which we will need in order to apply our Davis-Kahan style bound in \Cref{lemma:utwo_to_utwostar}. We remark that \Cref{lemma:eutwostar_norm} below holds both in the setting of \Cref{mainthm:nonhomogeneous} and of \Cref{mainthm:sqrtngap}. 
\begin{lemma}
\label{lemma:eutwostar_norm}
Suppose each crossing edge in our graph appears independently with probability $q$. There exists a universal constant $C$ such that for all $n$ sufficiently large, with probability $\ge 1-\delta$, we have
\begin{align*}
    \norm{\mE\utwostar}_2 \le C\inparen{\frac{\logv{\nfrac{1}{\delta}}}{\log n}}^{3/2}\inparen{\sqrt{nq} + (nq\logv{\nfrac{n}{\delta}})^{1/4} + \sqrt{\logv{\nfrac{n}{\delta}}}}.
\end{align*}
\end{lemma}
\begin{proof}[Proof of \Cref{lemma:eutwostar_norm}]
Observe that $\abs{\mE\utwostar} = 2\abs{\vdout - \exv{\vdout}}/\sqrt{n}$. By \Cref{lemma:dout}, for all $v \in V$, with probability $\ge 1-\delta/2$, we have $\vdout[v] \le nq/2 + C_{\ref{lemma:dout}}\inparen{\sqrt{nq\cdot \logv{2n/\delta}} + \logv{\nfrac{2n}{\delta}}}$.

So, if we let $\mA_{\mathsf{out}}$ and $\mA_{\mathsf{out}}^{\star}$ denote the adjacency matrices consisting only of the crossing edges and the expected value of that, respectively, then invoking \Cref{lemma:llv_concentration}, with probability $\ge 1-\delta$, we have
\begin{align*}
    \norm{\mE\utwostar}_2 &= \frac{2\norm{\vdout - \exv{\vdout}}_2}{\sqrt{n}} = \frac{2\norm{\inparen{\mA_{\mathsf{out}}-\mA_{\mathsf{out}}^{\star}}\onev}_2}{\sqrt{n}} \le 2\opnorm{\mA_{\mathsf{out}}-\mA_{\mathsf{out}}^{\star}} \\
    &\le 2C_{\ref{lemma:llv_concentration}}\inparen{\frac{\logv{\nfrac{2}{\delta}}}{\log n}}^{3/2}\inparen{\sqrt{\frac{nq}{2}} + \sqrt{C_{\ref{lemma:dout}}}\sqrt{nq + \sqrt{nq\logv{\nfrac{2n}{\delta}}}+\logv{\nfrac{2n}{\delta}}}},
\end{align*}
completing the proof of \Cref{lemma:eutwostar_norm}.
\end{proof}
Finally, we apply \Cref{lemma:llv_concentration} in order to bound bound $\norm{\va_v-\astar_v}_2 $. 
\begin{lemma}
\label{lemma:arow_norm}
In the setting of \Cref{mainthm:nonhomogeneous}, with probability $\ge 1-\delta$, we have
\begin{align*}
    \norm{\va_v-\astar_v}_2 \le C\inparen{\frac{\logv{\nfrac{1}{\delta}}}{\log n}}^{3/2}\inparen{\sqrt{n\pbar}+(n\pbar\logv{\nfrac{n}{\delta}})^{1/4}+\sqrt{\logv{\nfrac{n}{\delta}}}}.
\end{align*}
\end{lemma}
\begin{proof}[Proof of \Cref{lemma:arow_norm}]
We use a similar proof to that of \Cref{lemma:eutwostar_norm}. Indeed, invoke \Cref{lemma:llv_concentration} (observe that we can set $p_{ij}$ for the deterministic internal edges to $0$ as they do not affect $\mA-\exv{\mA}$) and notice that
\begin{align*}
    \norm{\va_v-\astar_v}_2 \le \opnorm{\mA-\mA^{\star}} \le C_{\ref{lemma:llv_concentration}}\inparen{\frac{\logv{\nfrac{2}{\delta}}}{\log n}}^{3/2}\inparen{\sqrt{n\pbar}+(n\pbar\logv{\nfrac{2n}{\delta}})^{1/4}+\sqrt{\logv{\nfrac{2n}{\delta}}}},
\end{align*}
where we used $d' \le n(\pbar+q)/2 + 2\max\inbraces{C_{\ref{lemma:dout}},C_{\ref{lemma:din}}}\inparen{\sqrt{n\pbar\logv{\nfrac{2n}{\delta}}}+\logv{\nfrac{2n}{\delta}}}$ from combining \Cref{lemma:dout} and \Cref{lemma:din}. This completes the proof of \Cref{lemma:arow_norm}.
\end{proof}

\subsection{Eigenvector perturbations}\label{sec:cluster_evect_perturbations}

In this Appendix, we give our Euclidean norm eigenvector perturbation bounds.

First, we verify that $\utwostar$ is indeed the second eigenvector of $\Lstar$.

\begin{lemma}
\label{lemma:u2expect}
In the setting of \Cref{mainthm:nonhomogeneous}, we have $\Lstar\utwostar = \lambda_2(\Lstar)\utwostar = nq\utwostar$, where $\Lstar = \exv{\mL}$. 

In the setting of \Cref{mainthm:sqrtngap}, we have $\Lstar\utwostar = \lambda_2(\Lstar)\utwostar = nq\utwostar$, where $\Lstar$ denotes the Laplacian matrix that agrees with $\mL$ on all internal edges and agrees with $\exv{\mL}$ on all crossing edges.
\end{lemma}
\begin{proof}[Proof of \Cref{lemma:u2expect}]
In both cases, one can check that $\utwostar$ is an eigenvector of $\Lstar$ with eigenvalue $nq$: for any $v \in P_2$ (i.e. $\utwostar[v] = -1/\sqrt{n}$ without loss of generality), one has
\begin{equation*}
    \left(\Lstar \utwostar\right)_v =\frac{1}{\sqrt{n}}\left( -(\vdin[v]+nq/2)-\sum_{w \in P_1: \{v,w\} \in E}(-1) + \sum_{w \in P_2} (-q)\right) = -\frac{nq}{\sqrt{n}} = nq \cdot \utwostar[v]\, .
\end{equation*}
By virtue of the above observations, it suffices to argue that $nq < \lambda_3(\Lstar) \le \dots \le \lambda_n(\Lstar)$. 

In the setting of \Cref{mainthm:nonhomogeneous}, we claim $\lambda_3^{\star} \geq \frac{n(p+q)}{2} > nq$. This is because because $p_{vw} \ge p$, which implies that if we consider $\Lstar_1$ to be the expected Laplacian for $\mathsf{SSBM}(n,p,q)$ and $\Lstar_2$ to be the expected Laplacian for $\mathsf{NSSBM}(n,p,\pbar,q)$, then $\Lstar_2 \succeq \Lstar_1$.. 

In the setting of \Cref{mainthm:sqrtngap}, we have $\lambda_3(\Lhat) - \lambda_2(\Lhat) >  nq$, by the theorem assumption. Since $\Lstar$ is obtained from $\Lhat$ by adding the adversarial edges, we have $\lambda_i(\Lstar) \ge \lambda_i(\Lhat)$ for all $i$. In particular, we have $\lambda_3(\Lstar) \ge \lambda_3(\Lhat) = \lambda_2(\Lhat)  + (\lambda_3(\Lhat) - \lambda_2(\Lhat)) > nq$, where the last inequality is using the fact $\lambda_2(\Lhat) \ge 0$. Therefore, $nq$ must be the second eigenvalue of $\Lstar$, completing the proof of \Cref{lemma:u2expect}.
\end{proof}

Next, we prove a general Davis-Kahan style bound. 
\begin{lemma}
\label{lemma:dk_easy}
Let $\mL$ and $\Lhat$ be two weighted Laplacian matrices. Let $\vu_2$ and $\utwohat$ be the second eigenvectors of $\mL$ and $\Lhat$, respectively. Then,
\begin{align*}
    \norm{\vu_2-\utwohat}_2 \le \sqrt{2} \cdot \min\inbraces{\frac{\norm{(\Lhat-\mL)\vu_2}_2}{\abs{\lambda_3(\Lhat)-\lambda_2(\mL)}},\frac{\norm{(\Lhat-\mL)\utwohat}_2}{\abs{\lambda_3(\mL)-\lambda_2(\Lhat)}}}
\end{align*}
\end{lemma}
\begin{proof}[Proof of \Cref{lemma:dk_easy}]
One can get this sort of guarantee from variants of the Davis-Kahan theorem, but it is more illuminating to write an eigenvalue decomposition and observe it from there. Without loss of generality, assume that $\ip{\utwohat,\vu_2}\ge 0$ (indeed, otherwise we can always negate $\utwohat$ if this is not the case). 
Notice that
\begin{align*}
    \norm{(\Lhat-\mL)\vu_2}_2^2 &= \norm{\inparen{\Lhat-\lambda_2(\mL)\mI}\vu_2}_2^2 \\
    &= (\lambda_2(\Lhat)-\lambda_2(\mL))^2\ip{\utwohat, \vu_2}^2 + \sum_{i=3}^n \inparen{\lambda_i(\Lhat)-\lambda_2(\mL)}^2\ip{\widehat{\vu_i},\vu_2}^2 \\
    &\ge \sum_{i=3}^n \inparen{\lambda_3(\Lhat)-\lambda_2(\mL)}^2 \ip{\widehat{\vu_i},\vu_2}^2 = \inparen{\lambda_3(\Lhat)-\lambda_2(\mL)}^2\inparen{1-\ip{\utwohat,\vu_2}^2},
\end{align*}
which rearranges to 
\begin{align*}
    \ip{\utwohat,\vu_2}^2 \ge 1 - \inparen{\frac{\norm{(\Lhat-\mL)\vu_2}_2}{\lambda_3(\Lhat)-\lambda_2(\mL)}}^2.
\end{align*}
Now, if  $\norm{(\Lhat-\mL)\vu_2}_2 \geq |\lambda_3(\Lhat)-\lambda_2(\mL)|$, then the condition  $\norm{\vu_2-\utwohat}_2 \leq \sqrt{2}\cdot  \frac{\norm{(\Lhat-\mL)\vu_2}_2}{\abs{\lambda_3(\Lhat)-\lambda_2(\mL)}}$ is trivially satisfied, since $\norm{\vu_2-\utwohat}_2 \leq \sqrt{2 - 2 \ip{\utwohat,\vu_2}} \leq \sqrt{2}$.
Otherwise, taking the square roots of both sides, we obtain 

\begin{align*}
    \ip{\utwohat,\vu_2} \ge \sqrt{1 - \inparen{\frac{\norm{(\Lhat-\mL)\vu_2}_2}{\lambda_3(\Lhat)-\lambda_2(\mL)}}^2}, 
\end{align*}
which gives
\begin{align*}
    \norm{\utwohat-\vu_2}_2^2 = 2-2\ip{\utwohat,\vu_2} \le 2-2\sqrt{1 - \inparen{\frac{\norm{(\Lhat-\mL)\vu_2}_2}{\lambda_3(\Lhat)-\lambda_2(\mL)}}^2} \le 2 \cdot \inparen{\frac{\norm{(\Lhat-\mL)\vu_2}_2}{\lambda_3(\Lhat)-\lambda_2(\mL)}}^2.
\end{align*}
Taking the square root of both sides and repeating this argument by exchanging the roles of $\mL$ and $\Lhat$ yields the statement of \Cref{lemma:dk_easy}.
\end{proof}
This immediately implies the following upper-bound on  $\norm{\vu_2-\utwostar}_2$. We will use it repeatedly, both in  \cref{mainthm:nonhomogeneous} and \Cref{mainthm:sqrtngap}. 
\begin{lemma}
\label{lemma:utwo_to_utwostar}
We have
\begin{align*}
    \norm{\vu_2-\utwostar}_2 \le \sqrt{2} \cdot \frac{\norm{\mE\utwostar}_2}{\abs{\lambda_3-\lambda_2^{\star}}}.
\end{align*}
\end{lemma}
\begin{proof}
\Cref{lemma:utwo_to_utwostar} immediately follows from \Cref{lemma:dk_easy} by letting $\Lhat = \Lstar$.
\end{proof}
Combining with \Cref{lemma:lambda3_lambda2star} and \Cref{lemma:eutwostar_norm}, we can now upper-bound $\norm{\vu_2-\utwostar}_2$ in the setting of \Cref{mainthm:nonhomogeneous}. 
\begin{lemma}
\label{lemma:utwo_to_utwostar_nssbm}
In the setting of \Cref{mainthm:nonhomogeneous}, there exists a universal constant $C$ such that, for $\delta \ge 3n^{-10}$, with probability $\ge 1-\delta$, we have
\begin{align*}
    \norm{\vu_2-\utwostar}_2 \le \frac{C}{\sqrt{\logv{\nfrac{n}{\delta}}}}.
\end{align*}
\end{lemma}
\begin{proof}[Proof of \Cref{lemma:utwo_to_utwostar_nssbm}]
Using \Cref{lemma:utwo_to_utwostar}, \Cref{lemma:lambda3_lambda2star}  and \Cref{lemma:eutwostar_norm}, we have
\begin{align*}
    \norm{\vu_2-\utwostar}_2 \le \frac{400\sqrt{2}C_{\ref{lemma:eutwostar_norm}}\inparen{\sqrt{nq}+\inparen{nq\logv{\nfrac{3n}{\delta}}}^{1/4}+\sqrt{\logv{\nfrac{3n}{\delta}}}}}{n(p-q)}.
\end{align*}
At this point, it is enough to show that there exists a universal constant $C$ such that
\begin{align*}
    Cn(p-q) \ge 400\sqrt{2}C_{\ref{lemma:eutwostar_norm}}\inparen{\sqrt{nq\logv{\nfrac{n}{\delta}}}+\inparen{nq}^{1/4}\inparen{\logv{\nfrac{n}{\delta}}}^{3/4}+\logv{\nfrac{n}{\delta}}}.
\end{align*}
To see this, note that for any two nonnegative real numbers we have $2a^{1/4}b^{1/4} \le \sqrt{b} + \sqrt{a}$, which implies $2a^{1/4}b^{3/4} \le b + \sqrt{ab}$. Let $a = nq$ and $b = \logv{\nfrac{3n}{\delta}}$, and we get
\begin{align*}
    &\quad 400\sqrt{2}C_{\ref{lemma:eutwostar_norm}}\inparen{\sqrt{nq\logv{\nfrac{3n}{\delta}}}+\inparen{nq}^{1/4}\inparen{\logv{\nfrac{n}{\delta}}}^{3/4}+\logv{\nfrac{3n}{\delta}}} \\
    &\le 800\sqrt{2}C_{\ref{lemma:eutwostar_norm}}\inparen{\sqrt{nq\logv{\nfrac{3n}{\delta}}}+\logv{\nfrac{3n}{\delta}}} \\
    &\le 800\sqrt{2}C_{\ref{lemma:eutwostar_norm}}\inparen{\sqrt{n\pbar\logv{\nfrac{3n}{\delta}}}+\logv{\nfrac{3n}{\delta}}} \le C n(p-q),
\end{align*}
where the last inequality follows from the assumption we gave in \Cref{mainthm:nonhomogeneous}. We therefore conclude the proof of \Cref{lemma:utwo_to_utwostar_nssbm}.
\end{proof}
Next, we prove $\ell_1$ norm concentration for the rows of $\mA$ and for the rows of $\mL$ in the setting of \Cref{mainthm:nonhomogeneous}. We will use this in \Cref{lemma:loo_close}, where we will bound $ \norm{\vu_2^{(v)}-\vu_2}_{2}$. Here $\vu_2^{(v)}$ denotes the second eigenvector of the leave-one-out Laplacian $\mL^{(v)}$. 
\begin{lemma}
\label{lemma:l_one_norm}
In the setting of \Cref{mainthm:nonhomogeneous}, there exists a universal constant $C$ such that with probability $\ge 1-\delta$, for all $v \in V$, we have
\begin{align*}
    \norm{\va_v-\astar_v}_1 &\le C\inparen{n\pbar+\sqrt{n\pbar\logv{\nfrac{n}{\delta}}}+\logv{\nfrac{n}{\delta}}} \\
    \norm{\vl_{v}-\exv{\vl_{v}}}_1 &\le C\inparen{n\pbar+\sqrt{n\pbar\logv{\nfrac{n}{\delta}}}+\logv{\nfrac{n}{\delta}}}.
\end{align*}
\end{lemma}
\begin{proof}[Proof of \Cref{lemma:l_one_norm}]
It is easy to see that
\begin{align*}
    \norm{\vl_v-\exv{\vl_v}}_1 = \abs{\vd[v] - \exv{\vd[v]}} + \norm{\va_v-\astar_v}_1.
\end{align*}
Let us consider the second term above. By \Cref{lemma:din} and \Cref{lemma:dout}, we have with probability $\ge 1-\delta/2$ that for all $v \in V$
\begin{align*}
    \norm{\va_v-\astar_v}_1 &\le \norm{\va_v}_1 + \norm{\astar_v}_1 \\
    &\le 2\inparen{\frac{n\pbar}{2} + \max\inbraces{C_{\ref{lemma:dout}},C_{\ref{lemma:din}}}\inparen{\sqrt{n\pbar\logv{\nfrac{4n}{\delta}}}+\logv{\nfrac{4n}{\delta}}}} + n\pbar \\
    &= 2n\pbar+2\max\inbraces{C_{\ref{lemma:dout}},C_{\ref{lemma:din}}}\inparen{\sqrt{n\pbar\logv{\nfrac{4n}{\delta}}}+\logv{\nfrac{4n}{\delta}}}.
\end{align*}
Finally, by \Cref{lemma:dout} and \Cref{lemma:din}, we have with probability $1-\delta/2$ that for all $v \in V$,
\begin{align*}
    \abs{\vd[v]-\exv{\vd[v]}} &\le \max_{v \in V} \abs{\vdout[v]-\exv{\vdout[v]}} + \max_{v \in V} \abs{\vdin[v]-\exv{\vdin[v]}} \\
    &\le 2\max\inbraces{C_{\ref{lemma:dout}},C_{\ref{lemma:din}}}\inparen{\sqrt{n\pbar\logv{\nfrac{4n}{\delta}}}+\logv{\nfrac{4n}{\delta}}}
\end{align*}
Adding everything up means that with probability $\ge 1-\delta$, for all $v \in V$, we have
\begin{align*}
    \norm{\vl_v - \exv{\vl_v}}_1 \le 2n\pbar + 4\max\inbraces{C_{\ref{lemma:din}},C_{\ref{lemma:dout}}}\inparen{\sqrt{n\pbar\logv{\nfrac{4n}{\delta}}}+\logv{\nfrac{4n}{\delta}}},
\end{align*}
which completes the proof of \Cref{lemma:l_one_norm}.
\end{proof}
Having established \Cref{lemma:l_one_norm}, we can now upper-bound $\norm{\vu_2^{(v)}-\vu_2}_{2}$.
\begin{lemma}
\label{lemma:loo_close}
In the setting of \Cref{mainthm:nonhomogeneous}, for $\delta \ge 2n^{-9}$ with probability $\ge 1-\delta$, for all $v \in V$, we have
\begin{align*}
    \norm{\vu_2^{(v)}-\vu_2}_{2} \le \maxnorm{\vu_2} \cdot \frac{C\inparen{\pbar+\sqrt{\pbar\logv{\nfrac{n}{\delta}}/n}+\logv{\nfrac{n}{\delta}}/n}}{p-q}
\end{align*}
\end{lemma}
\begin{proof}[Proof of \Cref{lemma:loo_close}]
Recall that the gap condition in \Cref{mainthm:nonhomogeneous} means that $p$ and $q$ are such that for a universal constant $C$,
\begin{align*}
    n(p-q) \ge C\inparen{\sqrt{n\pbar\logv{\nfrac{n}{\delta}}}+\logv{\nfrac{n}{\delta}}}.
\end{align*}
To appeal to \Cref{lemma:dk_easy}, we need to understand the entries of the matrix $\mL-\mL^{(v)}$. It is easy to see that this matrix only has nonzero entries on the diagonal and in the $v$th row and column. There, the $v$th row and column of $\mL-\mL^{(v)}$ are exactly equal to those of $\mL-\Lstar$. Moreover, the $w\neq v$th diagonal entry of $\mL-\mL^{(v)}$ is exactly $\indicator{(v,w) \in E} - p_{vw}$.

Hence, we have
\begin{align*}
    &\quad \norm{\inparen{\mL-\mL^{(v)}}\vu_2}_2\\
    &= \inparen{\sum_{w=1}^n \ip{\inparen{\mL-\mL^{(v)}}_w,\vu_2}^2}^{1/2} \\
    &= \inparen{\ip{\inparen{\mL-\Lstar}_v,\vu_2}^2 + \sum_{w \neq v} \inparen{\inparen{\va_v[w] - p_{vw}}\vu_2[w]-\inparen{\va_v[w] - p_{vw}}\vu_2[v]}^2}^{1/2} \\
    &\le \abs{\ip{\inparen{\mL-\Lstar}_v,\vu_2}} + \inparen{\sum_{w \neq v} \inparen{\inparen{\va_v[w] - p_{vw}}\vu_2[w]-\inparen{\va_v[w] - p_{vw}}\vu_2[v]}^2}^{1/2}\\
    &\le \inparen{\norm{\vl_v - \exv{\vl_v}}_1 + 2\norm{\va_v-\astar_v}_2} \cdot \maxnorm{\vu_2} \\
    &\le \inparen{\norm{\vl_v - \exv{\vl_v}}_1 + 2\norm{\va_v-\astar_v}_1} \cdot \maxnorm{\vu_2} \\
    &\le \maxnorm{\vu_2} \cdot 3C_{\ref{lemma:l_one_norm}}\inparen{n\pbar+\sqrt{n\pbar\logv{\nfrac{2n^2}{\delta}}}+\logv{\nfrac{2n^2}{\delta}}}.
\end{align*}
Now, let $C \ge 8C_{\ref{lemma:laplacian_concentration}}$. Using \Cref{lemma:laplacian_concentration} to understand the concentration of sampling the graph except edges incident to $v$, along with Weyl's inequality, we have with probability $\ge 1-\delta$ that for all $v \in V$ and for all $n$ sufficiently large,
\begin{align*}
    \abs{\lambda_3^{(v)} - \lambda_2} &\ge \inparen{\lambda_3^{(v)} - \lambda_3^{\star}} - \inparen{\lambda_2 - \lambda_2^{\star}} + \inparen{\lambda_3^{\star} - \lambda_2^{\star}} \\
    &\ge - 2\inparen{C_{\ref{lemma:laplacian_concentration}}\sqrt{n\pbar\logv{\nfrac{2n^2}{\delta}}}+\logv{\nfrac{2n^2}{\delta}}}+\frac{n(p-q)}{2} \ge \frac{n(p-q)}{4}.
\end{align*}

Now, using \Cref{lemma:dk_easy}, we get
\begin{align*}
    \norm{\vu_2^{(v)}-\vu_2}_2 &\le \frac{\norm{\inparen{\mL-\mL^{(v)}}\vu_2}_2}{|\lambda_3^{(v)} - \lambda_2|} \le \maxnorm{\vu_2} \cdot \frac{12C_{\ref{lemma:l_one_norm}}\inparen{n\pbar+\sqrt{n\pbar\logv{\nfrac{n}{\delta}}}+\logv{\nfrac{n}{\delta}}}}{n(p-q)} \\
    &\le \maxnorm{\vu_2} \cdot \frac{12C_{\ref{lemma:l_one_norm}}\inparen{\pbar+\sqrt{\pbar\logv{\nfrac{2n^2}{\delta}}/n}+\logv{\nfrac{2n^2}{\delta}}/n}}{p-q},
\end{align*}
completing the proof of \Cref{lemma:loo_close}.
\end{proof}

\subsection{Leave-one-out and bootstrap}\label{sec:cluster_leave_one_out}

The main goal of this section is to establish an upper-bound on $\abs{\ip{\va_v-\astar_v, \vu_2-\utwostar}}$ in the setting of \Cref{mainthm:nonhomogeneous}. To this end, we will need the following concentration inequality from \cite{afwz17}. 
\begin{lemma}[Lemma 7 from \cite{afwz17}]
\label{lemma:row_concentration}
Let $\vw \in \R^n$ and $X_i \sim \mathsf{Ber}(p_i)$. Let $p \ge p_i$ for all $i \in [n]$. Let $X \in \R^n$ be the vector formed by stacking the $X_i$. Then,
\begin{align*}
    \prv{\abs{\ip{\vw, X - \exv{X}}} \ge \frac{(2+a)pn}{\max\inparen{1, \logv{\frac{\sqrt{n}\norm{\vw}_{\infty}}{\norm{\vw}_2}}}} \cdot \norm{\vw}_{\infty}} \le 2\expv{-anp}.
\end{align*}
\end{lemma}

\begin{lemma}
\label{lemma:random_fluctuations}
In the setting of \Cref{mainthm:nonhomogeneous}, suppose $\va_v$ is such that $\va_v[w] \sim \mathsf{Bernoulli}(p_{vw})$ and let $\pbar \geq \max_{w \colon p_{vw }\neq 1} p_{vw}$. With probability $\ge 1-\delta$ for $\delta \ge 1/n^2$, for all $v \in V$, we have
\begin{align*}
    \abs{\ip{\va_v-\astar_v, \vu_2-\utwostar}} \le C\inparen{n\pbar + \logv{\nfrac{n}{\delta}}}\inparen{\frac{\norm{\vu_2}_{\infty}}{\log\log n} + \frac{1}{\sqrt{n}\log\log n}}.
\end{align*}
\end{lemma}
\begin{proof}[Proof of \Cref{lemma:random_fluctuations}]
Ideally, one would treat $\vu_2-\utwostar$ as fixed and then apply Bernstein's inequality to argue that the sum of centered Bernoulli random variables as written above concentrates well. Unfortunately, since $\vu_2$ depends on $\va_v-\astar_v$, we cannot express this inner product as the sum of independent random variables.

To resolve this, we use the leave-one-out method. Let $\vu_2^{(v)}$ be the second eigenvector of the leave-one-out Laplacian $\mL^{(v)}$ of $\mA^{(v)}$, where $\mA^{(v)}$ is chosen to agree with $\mA$ everywhere except for the $v$th row and $v$th column. The $v$th row and $v$th column of $\mA^{(v)}$ are replaced with those of $\Astar$. Now, $\va_v$ does not depend on $\mL^{(v)}$ and therefore $\vu_2^{(v)}$.

We therefore write
\begin{align*}
    \abs{\ip{\va_v-\astar_v, \vu_2-\utwostar}} &\le \abs{\ip{\va_v-\astar_v, \vu_2-\vu_2^{(v)}}} + \abs{\ip{\va_v-\astar_v, \vu_2^{(v)}-\utwostar}} \\
    &\le \norm{\va_v-\astar_v}_2 \cdot \norm{\vu_2^{(v)}-\vu_2}_2 + \abs{\ip{\va_v-\astar_v, \vu_2^{(v)}-\utwostar}} \\
    &\le \norm{\va_v-\astar_v}_2 \cdot \frac{C_{\ref{lemma:loo_close}}\pbar}{p-q}\maxnorm{\vu_2} + \abs{\ip{\va_v-\astar_v, \vu_2^{(v)}-\utwostar}}.
\end{align*}
To bound the rightmost term of the RHS, we use Lemma 7 of \cite{afwz17}, reproduced in \Cref{lemma:row_concentration}. In that, let $\vw \coloneqq \vu_2^{(v)} - \utwostar$. Let $a = \frac{1}{n\pbar}\logv{\nfrac{20n}{\delta}}$ so that $2 \exp(-2an\pbar) \leq \delta/(10n)$. Note that for the deterministic entries, we have $\va_v-\astar_v = 1-1 = 0$, so in \Cref{lemma:row_concentration}, we can set $X_w \sim \mathsf{Ber}(0)$ for these entries. Now, by \Cref{lemma:row_concentration}, with probability $\ge 1-\delta/n$, we have
\begin{align}
    \abs{\ip{\vu_2^{(v)}-\utwostar, \va_v-\va_v^{\star}}} \le \frac{2n\pbar + \logv{\frac{20n}{\delta}}}{\max\inparen{1, \logv{\frac{\sqrt{n}\norm{\vw}_{\infty}}{\norm{\vw}_2}}}} \cdot \norm{\vw}_{\infty}.\label{eq:row_concentrate_extraction}
\end{align}
Let us first bound $\maxnorm{\vw} = \maxnorm{\vu_2^{(v)}-\utwostar}$. We write
\begin{align}
    \maxnorm{\vu_2^{(v)}-\utwostar} &\le \maxnorm{\vu_2^{(v)}-\vu_2} + \maxnorm{\vu_2 - \utwostar} \\
    &\le \norm{\vu_2^{(v)}-\vu_2}_2 + \maxnorm{\vu_2}+\maxnorm{\utwostar} \\
    &\le 2\max\inbraces{C_{\ref{lemma:loo_close}}(\alpha,\delta)\maxnorm{\vu_2},\frac{1}{\sqrt{n}}}.\label{eq:row_concentrate_w_maxnorm}
\end{align}
In what follows, we omit the arguments $\alpha$ and $\delta$ in mentions of $C_{\ref{lemma:loo_close}}$. Next, using \Cref{lemma:utwo_to_utwostar_nssbm}, the triangle inequality, and $\delta \ge 1/n^3$, we have
\begin{align*}
    \norm{\vw}_2 = \norm{\vu_2^{(v)}-\utwostar}_2 \le C_{\ref{lemma:loo_close}}\norm{\vu_2}_{\infty} + \frac{4C_{\ref{lemma:utwo_to_utwostar_nssbm}}}{\sqrt{\log n}}.
\end{align*}
We now have two cases based on the value of $\sqrt{n} \cdot \frac{\maxnorm{\vu_2^{(v)}-\utwostar}}{\norm{\vu_2^{(v)}-\utwostar}_2}$.

\smallpar{Case 1 -- $\vw$ is not too ``flat.''} Let us first handle the case where
\begin{align*}
    \frac{\sqrt{n} \cdot \maxnorm{\vu_2^{(v)}-\utwostar}}{\norm{\vu_2^{(v)}-\utwostar}_2} \ge \sqrt{\log n}.
\end{align*}
We plug this into \eqref{eq:row_concentrate_extraction} and get
\begin{align*}
    \abs{\ip{\vu_2^{(v)}-\utwostar, \va_v-\va_v^{\star}}} &\le \frac{2n\pbar + \logv{\frac{20n}{\delta}}}{\max\inparen{1, \logv{\frac{\sqrt{n}\norm{\vw}_{\infty}}{\norm{\vw}_2}}}} \cdot \norm{\vw}_{\infty} \\
    &\le 4 \cdot \frac{n\pbar + \logv{\nfrac{20n}{\delta}}}{\log\log n}\inparen{C_{\ref{lemma:loo_close}}\maxnorm{\vu_2} + \frac{1}{\sqrt{n}}},
\end{align*}
where the last inequality follows from \eqref{eq:row_concentrate_w_maxnorm}.

\smallpar{Case 2 -- $\vw$ is ``flat.''} We now assume
\begin{align*}
    \frac{\sqrt{n} \cdot \maxnorm{\vu_2^{(v)}-\utwostar}}{\norm{\vu_2^{(v)}-\utwostar}_2} \le \sqrt{\log n}.
\end{align*}
We can easily check that the function
\begin{align*}
    \frac{x}{\max\inparen{1, \log x}}
\end{align*}
is increasing, so its maximum will be attained at the largest value of $x$ in the domain. Let $x = \sqrt{n}\norm{\vw}_{\infty}/\norm{\vw}_2$ and write
\begin{align*}
    &\quad \frac{2n\pbar + \logv{\frac{20n}{\delta}}}{\max\inparen{1, \logv{\frac{\sqrt{n}\norm{\vw}_{\infty}}{\norm{\vw}_2}}}} \cdot \norm{\vw}_{\infty} \\
    &= \frac{2n\pbar + \logv{\frac{20n}{\delta}}}{\max\inparen{1, \logv{\frac{\sqrt{n}\norm{\vw}_{\infty}}{\norm{\vw}_2}}}} \cdot \frac{\sqrt{n}\norm{\vw}_{\infty}}{\norm{\vw}_2} \cdot \frac{\norm{\vw}_2}{\sqrt{n}} \\
    &\le \frac{2n\pbar + \logv{\frac{20n}{\delta}}}{\log\log n} \cdot \sqrt{\frac{\log n}{n}} \cdot \norm{\vw}_2 \\
    &\le \frac{2n\pbar + \logv{\frac{20n}{\delta}}}{\log\log n} \cdot \sqrt{\frac{\log n}{n}} \cdot C_{\ref{lemma:loo_close}}\inparen{\maxnorm{\vu_2} + \frac{1}{\sqrt{\logv{\nfrac{20n^2}{\delta}}}}} \\
    &= C_{\ref{lemma:loo_close}}\inparen{\frac{2n\pbar + \logv{\frac{20n}{\delta}}}{\log\log n} \cdot \sqrt{\frac{\log n}{n}}\maxnorm{\vu_2} + \frac{2n\pbar + \logv{\frac{20n}{\delta}}}{\sqrt{n} \cdot \log\log n}}.
\end{align*}
All of this tells us that
\begin{align*}
    \abs{\ip{\va_v-\astar_v, \vu_2^{(v)}-\utwostar}} \le 4C_{\ref{lemma:loo_close}} \cdot \inparen{n\pbar + \logv{\nfrac{20n}{\delta}}}\inparen{\frac{\norm{\vu_2}_{\infty}}{\log\log n} + \frac{1}{\sqrt{n}\log\log n}}.
\end{align*}
It remains to handle the term
\begin{align*}
    \norm{\va_v-\astar_v}_2 \cdot \maxnorm{\vu_2}.
\end{align*}
Indeed, using \Cref{lemma:arow_norm}, we have with probability $\ge 1-\delta$ that
\begin{align*}
    \norm{\va_v-\astar_v}_2 \cdot \maxnorm{\vu_2} \le C_{\ref{lemma:arow_norm}}\inparen{\frac{\logv{\nfrac{20n}{\delta}}}{\log n}}^{3/2}\sqrt{n\pbar} \cdot \maxnorm{\vu_2}.
\end{align*}
Combining everything tells us that
\begin{align*}
    \abs{\ip{\va_v-\astar_v,\vu_2-\utwostar}} &\le 30C_{\ref{lemma:arow_norm}}\inparen{\frac{\logv{\nfrac{20n}{\delta}}}{\log n}}^{3/2}\sqrt{n\pbar} \cdot \maxnorm{\vu_2} \\
    &\quad + 4C_{\ref{lemma:loo_close}}\cdot\inparen{n\pbar + \logv{\nfrac{20n}{\delta}}}\inparen{\frac{\norm{\vu_2}_{\infty}}{\log\log n} + \frac{1}{\sqrt{n}\log\log n}} \\
    &\le C\inparen{n\pbar + \logv{\nfrac{20n}{\delta}}}\inparen{\frac{\norm{\vu_2}_{\infty}}{\log\log n} + \frac{1}{\sqrt{n}\log\log n}}.
\end{align*}
Taking a union bound over all $v \in V$ concludes the proof of \Cref{lemma:random_fluctuations}.
\end{proof}
Finally, we establish an upper-bound on $\maxnorm{\vu_2}$. This will be used repeatedly in the proof of \Cref{mainthm:nonhomogeneous}. 
\begin{lemma}
\label{lemma:utwo_inf_norm_bootstrap}
In the same setting as \Cref{mainthm:nonhomogeneous}, with probability $\ge1-\delta$ for $\delta \ge 10n^2$, we have for some constant $C(\alpha,\delta)$ that
\begin{align*}
    \maxnorm{\vu_2} \le \frac{C(\alpha,\delta)}{\sqrt{n}}.
\end{align*}
\end{lemma}
\begin{proof}[Proof of \Cref{lemma:utwo_inf_norm_bootstrap}]
First, observe that
\begin{align*}
    (\mD-\mA)\vu_2 = \lambda_2\vu_2,
\end{align*}
which means that
\begin{align*}
    \inparen{\mD-\lambda_2\mI}^{-1}\mA\vu_2 = \vu_2.
\end{align*}
By \Cref{lemma:dv_lambda2_positive}, with probability $\ge 1-\delta$, for all $v \in V$ we have
\begin{align*}
    \vd[v] - \lambda_2 \ge \frac{n(p-q)}{4}.
\end{align*}
Combining with \Cref{lemma:degree_vs_lambda}, we have
\begin{align*}
    \frac{\vdin[v]-\vdout[v]}{\vdin[v]-\vdout[v] + \inparen{2\vdout[v]-\lambda_2}} &= 1-\frac{2\vdout[v]-\lambda_2}{\vdin[v]-\vdout[v] + \inparen{2\vdout[v]-\lambda_2}} \\
    &\le 1+\frac{C_{\ref{lemma:degree_vs_lambda}}\inparen{\sqrt{nq\logv{\nfrac{10n}{\delta}}}+\logv{\nfrac{10n}{\delta}}}}{\vdin[v]-\vdout[v] + \inparen{2\vdout[v]-\lambda_2}} \\
    &\le 1+\frac{4C_{\ref{lemma:degree_vs_lambda}}\inparen{\sqrt{nq\logv{\nfrac{10n}{\delta}}}+\logv{\nfrac{10n}{\delta}}}}{n(p-q)} \le C',
\end{align*}
for some constant $C'>0$, where the penultimate line follows from \Cref{lemma:dv_lambda2_positive} and the last line follows from the gap assumption in \Cref{mainthm:nonhomogeneous}. Furthermore, by \Cref{lemma:npbar_vs_denom} and \Cref{lemma:utwo_to_utwostar_nssbm}, we have with probability $\ge1-\delta$ that for all $v \in V$,
\begin{align*}
    \frac{\abs{\ip{\astar_v,\utwostar-\vu_2}}}{\vd[v]-\lambda_2} \le \frac{\pbar\sqrt{n}}{\vd[v]-\lambda_2} \cdot \frac{C_{\ref{lemma:utwo_to_utwostar_nssbm}}}{\sqrt{\logv{\nfrac{10n}{\delta}}}} \le \frac{C_{\ref{lemma:npbar_vs_denom}}(\alpha) \cdot C_{\ref{lemma:utwo_to_utwostar_nssbm}}}{\sqrt{n\logv{\nfrac{10n}{\delta}}}}.
\end{align*}

Now, using \Cref{lemma:npbar_vs_denom} (and using \Cref{lemma:dv_lambda2_positive} to ensure that $\vd[v]-\lambda_2 > 0$ for all $v \in V$), we have
\begin{align*}
    \maxnorm{\vu_2} &= \maxnorm{\inparen{\mD-\lambda_2\mI}^{-1}\mA\vu_2} \\
    &= \maxnorm{\inparen{\mD-\lambda_2\mI}^{-1}\mA\vu_2 - \inparen{\mD-\lambda_2\mI}^{-1}\mA\utwostar + \inparen{\mD-\lambda_2\mI}^{-1}\mA\utwostar} \\
    &\le \maxnorm{\inparen{\mD-\lambda_2\mI}^{-1}\mA\utwostar} + \maxnorm{\inparen{\mD-\lambda_2\mI}^{-1}\mA(\utwostar-\vu_2)} \\
    &= \max_{1 \le v \le n} \frac{\abs{\ip{\va_v,\utwostar}}}{\vd[v]-\lambda_2} + \max_{1 \le v \le n} \frac{\abs{\ip{\va_v,\utwostar-\vu_2}}}{\vd[v]-\lambda_2} \\
    &= \frac{1}{\sqrt{n}}\inparen{\max_{1 \le v \le n} \frac{\abs{\vdin[v]-\vdout[v]}}{\vd[v]-\lambda_2}} + \max_{1 \le v \le n} \frac{\abs{\ip{\va_v,\utwostar-\vu_2}}}{\vd[v]-\lambda_2} \\
    &\le \frac{C}{\sqrt{n}} + \max_{1 \le v \le n} \frac{\abs{\ip{\va_v-\astar_v,\utwostar-\vu_2}}}{\vd[v]-\lambda_2} + \max_{1 \le v \le n} \frac{\abs{\ip{\astar_v,\utwostar-\vu_2}}}{\vd[v]-\lambda_2} \\
    &\le \frac{C}{\sqrt{n}} + \frac{C_{\ref{lemma:random_fluctuations}}\inparen{n\pbar+\logv{\nfrac{10n}{\delta}}}}{\vd[v]-\lambda_2} \cdot \inparen{\frac{1}{\sqrt{n}\log\log n} + \frac{\maxnorm{\vu_2}}{\log\log n}}+\frac{C_{\ref{lemma:npbar_vs_denom}}(\alpha) \cdot C_{\ref{lemma:utwo_to_utwostar_nssbm}}}{\sqrt{n\logv{\nfrac{10n}{\delta}}}} \\
    &\le \frac{C}{\sqrt{n}} +C_{\ref{lemma:random_fluctuations}} \cdot C_{\ref{lemma:npbar_vs_denom}}(\alpha) \cdot \inparen{\frac{1}{\sqrt{n}\log\log n} + \frac{\maxnorm{\vu_2}}{\log\log n}} + \frac{C_{\ref{lemma:npbar_vs_denom}}(\alpha) \cdot C_{\ref{lemma:utwo_to_utwostar_nssbm}}}{\sqrt{n\logv{\nfrac{10n}{\delta}}}}. 
\end{align*}
Note that any $n$ large enough
\begin{align*}
     \frac{C_{\ref{lemma:random_fluctuations}} \cdot C_{\ref{lemma:npbar_vs_denom}}(\alpha) \cdot \maxnorm{\vu_2}}{\log\log n} \le \frac{\maxnorm{\vu_2}}{2}.
\end{align*}
Thus, rearranging and solving for $\maxnorm{\vu_2}$ yields 
\begin{align*}
    \maxnorm{\vu_2} &\le 2\inparen{\frac{C}{\sqrt{n}} +C_{\ref{lemma:random_fluctuations}} \cdot C_{\ref{lemma:npbar_vs_denom}}(\alpha) \cdot \inparen{\frac{1}{\sqrt{n}\log\log n}} + \frac{C_{\ref{lemma:npbar_vs_denom}}(\alpha) \cdot C_{\ref{lemma:utwo_to_utwostar_nssbm}}}{\sqrt{n\logv{\nfrac{10n}{\delta}}}}},
\end{align*}
completing the proof of \Cref{lemma:utwo_inf_norm_bootstrap}.
\end{proof}

\subsection{Strong consistency of unnormalized spectral bisection}
\label{sec:cluster_strong_consistency}

In this section, we prove our main positive results \Cref{mainthm:nonhomogeneous} and \Cref{mainthm:sqrtngap}. It will be helpful to recall the proof sketches given in \Cref{sec:cluster_sbm_chapter_overview} while reading this section.

At a high level, the proof plan is as follows.
\begin{enumerate}
    \item We first establish a sufficient condition for a particular vertex to be classified correctly. We can think of this as simultaneously showing that the intermediate estimator $(\mD-\lambda_2\mI)^{-1}\mA\utwostar$ is strongly consistent and that the corresponding ``noise'' term $(\mD-\lambda_2\mI)^{-1}\mA(\utwostar-\vu_2)$ is a lower-order term in comparison to this. For a more formal way to see this, see \Cref{lemma:strong_consistency}.
    \item For the proof of \Cref{mainthm:nonhomogeneous}, the main technical challenge in showing that the noise term above is small amounts to analyzing the random quantity $\abs{\ip{\va_v,\utwostar-\vu_2}}$. This is where we will have to use the leave-one-out method to decouple the dependence between $\va_v$ and $\vu_2$. The relevant lemmas for the leave-one-out analysis are \Cref{lemma:random_fluctuations} and \Cref{lemma:utwo_inf_norm_bootstrap}.
    \item Finally, for the proof of \Cref{mainthm:sqrtngap}, we again appeal to \Cref{lemma:strong_consistency} but use a different approach to show that the noise term is small.
\end{enumerate}

\subsubsection{A sufficient condition for exact recovery and proof}

The main result of this subsection is \Cref{lemma:strong_consistency}, which gives a general condition under which a particular vertex will be classified correctly. The proofs of \Cref{mainthm:nonhomogeneous} and \Cref{mainthm:sqrtngap} will follow by invoking \Cref{lemma:strong_consistency}. We remark that the point of this lemma is mostly conceptual; the crux of the analysis lies in establishing that these conditions are satisfied our models.

\begin{lemma}
\label{lemma:strong_consistency}
Let $v \in V$ be some vertex. If $\vd[w] - \lambda_2 > 0$ for all $w \in V$, $\vdin[v] > \vdout[v]$, and $\abs{\ip{\va_v,\utwostar-\vu_2}} \le (\vdin[v]-\vdout[v])/\sqrt{n}$, then $\signv{\vu_2[v]} = \signv{\utwostar[v]}$, i.e., $\vu_2$ correctly classifies vertex $v$.
\end{lemma}
The goal of the rest of this section is to prove \Cref{lemma:strong_consistency}.

Our approach is to study the intermediate estimator
\begin{align*}
    \inparen{\mD-\lambda_2\mI}^{-1}\mA\utwostar.
\end{align*}
At a high level, our goal is to show that this correctly classifies all the vertices with high probability and also is very close to $\vu_2$ in $\ell_{\infty}$ norm with high probability. \citet{dls20} used this intermediate estimator to prove the strong consistency of unnormalized spectral bisection for $\mathsf{SBM}(n,p,q)$ instances.

Next, we show that this estimator is consistent and prove \Cref{lemma:strong_consistency}.

\begin{proof}[Proof of \Cref{lemma:strong_consistency}]
Observe that
\begin{align*}
    \vu_2 = \inparen{\mD-\lambda_2\mI}^{-1}\mA\utwostar - \inparen{\mD-\lambda_2\mI}^{-1}\mA\inparen{\utwostar-\vu_2}.
\end{align*}
Without loss of generality, suppose $v \in P_1$. In particular, this means that $\utwostar[v] = 1/\sqrt{n}$. Our goal is to show that $\vu_2[v] > 0$. And, as per the above, this means that it is enough to show that
\begin{align*}
    \inparen{\inparen{\mD-\lambda_2\mI}^{-1}\mA\utwostar}[v] \ge \inparen{\inparen{\mD-\lambda_2\mI}^{-1}\mA\inparen{\utwostar-\vu_2}}[v],
\end{align*}
or equivalently, using the fact that $\vd[v]-\lambda_2 > 0$, 
\begin{align*}
    \ip{\va_v,\utwostar} \ge \ip{\va_v,\utwostar-\vu_2},
\end{align*}
where $\va_v$ denotes the $v$-th row of $A$. 
To see that the above holds, use the fact that we know that $\vdin[v]-\vdout[v] > 0$, which gives
\begin{align*}
    \ip{\va_v,\utwostar} = \frac{\vdin[v]-\vdout[v]}{\sqrt{n}} \ge \abs{\ip{\va_v,\utwostar-\vu_2}} \ge \ip{\va_v,\utwostar-\vu_2}.
\end{align*}
This is exactly what we needed, and we conclude the proof of \Cref{lemma:strong_consistency}.
\end{proof}

\subsection{Proofs of main results}

At this point, we are ready to prove our main results.

\subsubsection{Nonhomogeneous symmetric stochastic block model (Proof of \texorpdfstring{\Cref{mainthm:nonhomogeneous}}{Theorem 1})}
\label{sec:cluster_proof_nonhomogeneous}

We are finally ready to prove \Cref{mainthm:nonhomogeneous}. For convenience, we reproduce its statement here.

\nonhomogeneous*

\begin{proof}[Proof of \Cref{mainthm:nonhomogeneous}]
As mentioned in \Cref{sec:cluster_results}, we actually prove a slightly stronger statement -- we will allow the adversary to set at most $n\pbar/\log \log n$ of the $p_{vw}$ to $1$ per vertex $v$ (in other words, the adversary can commit to at most $n\pbar/\log\log n$ edges per vertex that are guaranteed to appear in the final graph).

Our plan is to apply \Cref{lemma:strong_consistency}. In order to do so, we start with showing that for all $v$, we have $\vdin[v] > \vdout[v]$. By \Cref{lemma:degree_diff_easy}, with probability $\ge 1-\delta$, we have for all $v \in V$ that

\begin{align*}
    \vdin[v]-\vdout[v] \ge \frac{n(p-q)}{2} - C_{\ref{lemma:degree_diff_easy}}\inparen{\sqrt{np\logv{\nfrac{n}{\delta}}}+\logv{\nfrac{n}{\delta}}} > 0.
\end{align*}
Additionally, by \Cref{lemma:dv_lambda2_positive}, we have for all $v$ that $\vd[v] > \lambda_2$. 

The final item we need is to show that for all $v \in V$, we have $\abs{\ip{\va_v,\utwostar-\vu_2}} \le \abs{\ip{\va_v,\utwostar}}$. Observe that
\begin{align*}
    \abs{\ip{\va_v,\utwostar-\vu_2}} &\le \abs{\ip{\astar_v,\utwostar-\vu_2}} + \abs{\ip{\va_v-\astar_v,\utwostar-\vu_2}},
\end{align*}
where $\astar_v$ denotes the $v$-th row of $\exv{\mA}$. 
We handle the terms one at a time. First, note that by \Cref{lemma:lambda3_lambda2star}, with probability $\ge 1-\delta$, we have
\begin{align*}
    \lambda_3-\lambda_2^{\star} \ge \frac{n(p-q)}{4}.
\end{align*}
Now, let $\mE := \mL - \exv{\mL}$, and let $\astar_v[\rand] \in \R^{V}$ correspond to the vector that entrywise agrees with $\astar_v$ wherever $\astar_v$ is not $1$ and is zero elsewhere. This corresponds to the edges incident to $v$ that will be sampled randomly from the distribution over graphs. This means that for all $n \ge N(\delta)$ and choosing $\delta \ge 1/(10n)$, we have

\begin{align*}
    \abs{\ip{\astar_v[\rand],\utwostar-\vu_2}} &\le \norm{\astar_v[\rand]}_2 \cdot \frac{\sqrt{2}\norm{\mE\utwostar}_2}{\abs{\lambda_3-\lambda_2^{\star}}} & \text{( \Cref{lemma:utwo_to_utwostar})}\\
    &\le \pbar\sqrt{n} \cdot \frac{40\sqrt{2}C_{\ref{lemma:eutwostar_norm}}\inparen{\sqrt{nq}+(nq\log n)^{1/4}+\sqrt{\log n}}}{n(p-q)} & \text{(\Cref{lemma:eutwostar_norm,lemma:lambda3_lambda2star})} \\
    &\le \frac{1000 C_{\ref{lemma:eutwostar_norm}}n\pbar}{\sqrt{n}\log\log n} & \text{(gap in \Cref{mainthm:nonhomogeneous})}
\end{align*}
To handle the oblivious insertions, let $\ddet \in \R^V$ denote the degree vector that counts the number of deterministic edges inserted incident to $v$, for all $v \in V$. Under this notation, we have
\begin{align*}
    \abs{\ip{\astar_v-\astar_v[\rand],\utwostar-\vu_2}} \le \ddet[v] \cdot \maxnorm{\vu_2-\utwostar} \le \frac{n\pbar}{\sqrt{n}\log\log n} + \frac{n\pbar\maxnorm{\vu_2}}{\log\log n}.
\end{align*}
where the last inequality follows from using $\maxnorm{\vu_2-\utwostar} \le \maxnorm{\vu_2}+\maxnorm{\utwostar}$.
Combining yields
\begin{align*}
    \abs{\ip{\astar_v,\utwostar-\vu_2}} &\le C' \frac{ n\pbar}{\sqrt{n}\log\log n} + \frac{n\pbar\maxnorm{\vu_2}}{\log\log n},
\end{align*}
for some constant $C' >0$. Now, notice that for all $n$ sufficiently large,
\begin{align*}
    &\quad \abs{\ip{\va_v-\astar_v,\utwostar-\vu_2}}\\
    &\le C_{\ref{lemma:random_fluctuations}}\inparen{n\pbar + \logv{\nfrac{n}{\delta}}}\inparen{\frac{\norm{\vu_2}_{\infty}}{\log\log n} + \frac{1}{\sqrt{n}\log\log n}} & \text{(\Cref{lemma:random_fluctuations})} \\
    &\le C_{\ref{lemma:random_fluctuations}}\inparen{n\pbar + \logv{\nfrac{n}{\delta}}}\inparen{\frac{\frac{C_{\ref{lemma:utwo_inf_norm_bootstrap}}(\alpha,\delta)}{\sqrt{n}}}{\log\log n} + \frac{1}{\sqrt{n}\log\log n}} & \text{(\Cref{lemma:utwo_inf_norm_bootstrap})} \\
    &\le \frac{C_1(\alpha,\delta)\cdot(n\pbar + \logv{\nfrac{n}{\delta}})}{\sqrt{n}\log\log n}.
\end{align*}
Adding yields for $n \ge N(\alpha,\delta)$,
\begin{align*}
    \abs{\ip{\va_v,\utwostar-\vu_2}} &\le \abs{\ip{\astar_v,\utwostar-\vu_2}} + \abs{\ip{\va_v-\astar_v,\utwostar-\vu_2}} \\
    &\le \frac{C_2(\alpha,\delta)\cdot(n\pbar + \logv{\nfrac{n}{\delta}})}{\sqrt{n}\log\log n} \\
    &\le \frac{1}{\sqrt{n}} \cdot \inparen{\frac{n(p-q)}{2} - C_{\ref{lemma:degree_diff_easy}}\inparen{\sqrt{np\logv{\nfrac{n}{\delta}}}+\logv{\nfrac{n}{\delta}}}} & \text{(gap condition)} \\
    &\le  \frac{\vdin[v]-\vdout[v]}{\sqrt{n}} = \abs{\ip{\va_v,\utwostar}},
\end{align*}
which means we satisfy the conditions required by \Cref{lemma:strong_consistency}. Taking a union bound over all our (constantly many) probabilistic statements, setting $\delta = \Theta(1/n)$, and rescaling completes the proof of \Cref{mainthm:nonhomogeneous}.
\end{proof}

\subsubsection{Deterministic clusters model}
\label{sec:cluster_proof_sqrtngap}

For convenience, we reproduce the statement of \Cref{mainthm:sqrtngap} here.
\sqrtngap*
\begin{proof}[Proof of \Cref{mainthm:sqrtngap}]
In this proof, let $\Lstar$ be the Laplacian matrix that agrees with $\mL$ on all internal edges and agrees with $\exv{\mL}$ on all crossing edges. Let $\mL^{(\text{cross})}$ denote the Laplacian matrix corresponding to the cross edges, so we can write $ \Lstar = \mL - \mL^{(\text{cross})} + \exv{\mL^{(\text{cross})}}$. Although $\Lstar \neq \exv{\mL}$ due to the adaptive adversary, by \Cref{lemma:u2expect}, we still have $\Lstar\utwostar = \lambda_2^{\star}\utwostar = nq\utwostar$. Moreover, $(\mL-\Lstar)\utwostar$ is the vector whose entries are of the form $2(\vdout[v]-\exv{\vdout[v]})/\sqrt{n}$. Thus, we will be able to apply \Cref{lemma:utwo_to_utwostar} and \Cref{lemma:eutwostar_norm} later on. Finally, observe that $\lambda_i(\Lstar) \ge \lambda_i(\widehat{\mL})$ for all $i \ge 3$ and $\lambda_2(\widehat{\mL})=\lambda_2(\Lstar)=nq$. Thus, one can use the spectral gap $\lambda_3(\widehat{\mL})-\lambda_2(\widehat{\mL})$ to reason about $\lambda_3^{\star}-\lambda_2^{\star}$.

Let $\delta \ge 1/(10n)$. We will apply \Cref{lemma:strong_consistency} to get strong consistency. First, let us verify that $\vd[v] > \lambda_2$ for all $v$. Applying \Cref{lemma:laplacian_concentration} to the matrix $\mL^{(\text{cross})}$ gives
\begin{align*}
    \opnorm{ \mL - \Lstar} = \opnorm{\mL^{(\text{cross})} - \exv{\mL^{(\text{cross})}}} \le C_{\ref{lemma:laplacian_concentration}}\inparen{\sqrt{nq\logv{\nfrac{n}{\delta}}} + \log(n/\delta)}.
\end{align*}
Thus, using Weyl's inequality, for $n > N(\delta)$, we have
\begin{align*}
    \vd[v] - \lambda_2 &\ge \vdin[v]-\lambda^{\star}_2 -  \opnorm{ \mL - \Lstar} \\
    &\ge C_1 \frac{nq}{2}+C_1\sqrt{n}-nq-C_{\ref{lemma:laplacian_concentration}}\inparen{\sqrt{nq\logv{\nfrac{n}{\delta}}} + \log(n/\delta)} >0.
\end{align*}
Next, we verify that $\vdin[v] > \vdout[v]$ for all $v$. By \Cref{lemma:dout}, with probability $\ge 1-\delta$, for all $v \in V$, we have
\begin{align*}
    \abs{\vdout[v] - \frac{nq}{2}} \le C_{\ref{lemma:dout}}\inparen{\sqrt{nq\logv{\nfrac{n}{\delta}}}+\logv{\nfrac{n}{\delta}}}.
\end{align*}
So for $n > N(\delta)$, we obtain
\begin{align*}
    \vdin[v]-  \vdout[v] \ge C_1 \frac{nq}{2}+C_1\sqrt{n} - \frac{nq}{2} - C_{\ref{lemma:dout}}\inparen{\sqrt{nq\logv{\nfrac{n}{\delta}}}+\logv{\nfrac{n}{\delta}}} >0. 
\end{align*}
Here, in the last inequality we used the fact that $\sqrt{nq\logv{\nfrac{n}{\delta}}} \leq \max\{nq, \log(n/\delta)\}$ . 

Finally, we need to show that for all $v \in V$,
\begin{align*}
    \abs{\ip{\va_v,\utwostar-\vu_2}} \le \abs{\ip{\va_v,\utwostar}} = \frac{\vdin[v]-\vdout[v]}{\sqrt{n}}.
\end{align*}
By Cauchy-Schwarz, we have
\begin{align*}
    \abs{\ip{\va_v,\utwostar-\vu_2}} \le \norm{\va_v}_2 \cdot \norm{\utwostar-\vu_2}_2 = \sqrt{\vdin[v]+\vdout[v]} \cdot \norm{\utwostar-\vu_2}_2.
\end{align*}
Thus, it is enough to show that for all $v \in V$ we get
\begin{align*}
    \sqrt{n}\norm{\utwostar-\vu_2}_2 \le \frac{\vdin[v]-\vdout[v]}{\sqrt{\vdin[v]+\vdout[v]}}.
\end{align*}

Observe that the RHS above is a decreasing function in $\vdout[v]$ and an increasing function in $\vdin[v]$. 

Now, by \Cref{lemma:utwo_to_utwostar} and \Cref{lemma:eutwostar_norm}, we have
\begin{align}
    \sqrt{n}\norm{\utwostar-\vu_2}_2 \le \frac{\sqrt{n}\norm{\mE\utwostar}_2}{\abs{\lambda_3-\lambda_2^{\star}}} \le \frac{6C_{\ref{lemma:eutwostar_norm}}\sqrt{n}\inparen{\sqrt{nq}+(nq\logv{\nfrac{n}{\delta}})^{1/4}+\sqrt{\logv{\nfrac{n}{\delta}}}}}{\abs{\lambda_3-\lambda_2^{\star}}}.\label{eq:thm2_dk}
\end{align}
We now do casework on the value of $q$.

\smallpar{Case 1: $q \le \logv{\nfrac{n}{\delta}} / n$.} Carrying on from \eqref{eq:thm2_dk} and applying \Cref{lemma:laplacian_concentration} (we can set $p_{ij}$ for the deterministic internal edges to $0$ as they do not affect $\mL-\exv{\mL}$) along with Weyl's inequality, for all $n \ge N(\delta)$ we have
\begin{align*}
    \sqrt{n}\norm{\utwostar-\vu_2}_2 &\le \frac{18C_{\ref{lemma:eutwostar_norm}}\sqrt{n\logv{\nfrac{n}{\delta}}}}{\abs{\lambda_3-\lambda_2^{\star}}} \le \frac{18C_{\ref{lemma:eutwostar_norm}}\sqrt{n\logv{\nfrac{n}{\delta}}}}{\sqrt{n}-3C_{\ref{lemma:laplacian_concentration}}\logv{\nfrac{n}{\delta}}} \\
    &\le C\sqrt{\logv{\nfrac{n}{\delta}}} \ll \frac{\vdin[v]-\vdout[v]}{\sqrt{\vdin[v]+\vdout[v]}}, 
\end{align*}
as required. Here the last inequality follows using the fact that $\vdin[v] \ge C_1 \inparen{ \frac{nq}{2} + \sqrt{n}}$ and $\vdout[v] \le \frac{nq}{2} + 2C_{\ref{lemma:dout}}\log(n/\delta)$. 

\smallpar{Case 2: $\logv{\nfrac{n}{\delta}}/n \le q$.} Similar to the previous case, we get
\begin{align}
     \sqrt{n}\norm{\utwostar-\vu_2}_2 &\le \frac{18C_{\ref{lemma:eutwostar_norm}}\sqrt{n}\cdot\sqrt{nq}}{\abs{\lambda_3-\lambda_2^{\star}}} \le \frac{18 C_{\ref{lemma:eutwostar_norm}}\sqrt{n}\cdot\sqrt{nq}}{\sqrt{n}+(C_2 - 2C_{\ref{lemma:laplacian_concentration}})nq}\\
     &\le 18 C_{\ref{lemma:eutwostar_norm}}\cdot\max\inbraces{\sqrt{nq},\frac{1}{(C_2 - 2C_{\ref{lemma:laplacian_concentration}})\sqrt{q}}}.\label{eq:thm2_lhs}
\end{align}
Additionally, we can use the conclusion of \Cref{lemma:dout} to write with probability $\ge 1-\delta$ for all $v \in V$ and $n \ge N(\delta)$ that
\begin{align}
     \frac{\vdin[v]-\vdout[v]}{\sqrt{\vdin[v]+\vdout[v]}}
&\ge \frac{(C_1/2-2C_{\ref{lemma:dout}} -1/2)nq+C_1\sqrt{n}}{\sqrt{(C_1/2 + 2C_{\ref{lemma:dout}}+1/2)nq}} \\
&\ge \frac{C_1/2-2C_{\ref{lemma:dout}} -1/2}{\sqrt{C_1/2 + 2C_{\ref{lemma:dout}}+1/2}}\max\inbraces{\sqrt{nq},\sqrt{\frac{1}{q}}}.\label{eq:thm2_degreefn}
\end{align}
From this, it is clear that one can choose constants $C_1$ and $C_2$ such that \eqref{eq:thm2_lhs} is at most \eqref{eq:thm2_degreefn}. Taking a union bound over all our (constantly many) probabilistic statements, setting $\delta = \Theta(1/n)$, and rescaling completes the proof of \Cref{mainthm:sqrtngap}.
\end{proof}

\subsection{Inconsistency of normalized spectral bisection}
\label{sec:cluster_inconsistency}

In this section, we design a family of problem instances on which unnormalized spectral bisection is strongly consistent whereas normalized spectral bisection is inconsistent. Specifically, our goal is to prove \Cref{mainthm:inconsistency}.

\inconsistency*

\subsubsection{The nested block example}

We first state the family of instances on which we will prove our inconsistency results. Let $n$ be a multiple of $4$. Let $L_1$ consist of indices $1,\dots,n/4$, $L_2$ consist of indices $n/4+1,\dots,n/2$, and $R$ consist of indices $n/2+1,\dots,n$.

As mentioned in \Cref{sec:cluster_sbm_chapter_overview}, consider the following block structure determined by the $\Astar$ written below, where $q < p$ and $K \ge 3p/q$.
\begin{table}[H]
\centering
\begin{tabular}{l|cccl}
                     & \multicolumn{1}{l}{$L_1$}               & \multicolumn{1}{l}{$L_2$}                                    & \multicolumn{2}{l}{$R$}                                                     \\ \hline
$L_1$                & $Kp \cdot \mathbbm{1}_{n/4 \times n/4}$ & \multicolumn{1}{c|}{$p \cdot \mathbbm{1}_{n/4 \times n/4}$}  & \multicolumn{2}{c}{\multirow{2}{*}{$q \cdot \mathbbm{1}_{n/2 \times n/2}$}} \\
$L_2$                & $p \cdot \mathbbm{1}_{n/4 \times n/4}$  & \multicolumn{1}{c|}{$Kp \cdot \mathbbm{1}_{n/4 \times n/4}$} & \multicolumn{2}{c}{}                                                        \\ \cline{2-5} 
\multirow{2}{*}{$R$} & \multicolumn{2}{c|}{\multirow{2}{*}{$q \cdot \mathbbm{1}_{n/2 \times n/2}$}}                           & \multicolumn{2}{c}{\multirow{2}{*}{$p \cdot \mathbbm{1}_{n/2 \times n/2}$}} \\
                     & \multicolumn{2}{c|}{}                                                                                  & \multicolumn{2}{c}{}                                                       
\end{tabular}
\caption{$\Astar$ is defined to have the above block structure.}
\end{table}
We will draw our instances from the nonhomogeneous stochastic block model according to the probabilities prescribed above. Note that within the two clusters $L \coloneqq L_1 \cup L_2$ and $R$, each edge appears with probability at least $p$. Moreover, each edge in $L \times R$ appears with probability exactly $q$. However, there are also two subcommunities $L_1$ and $L_2$ that appear within $L$. Furthermore, observe that unnormalized spectral bisection is consistent on this family of examples with probability $\ge 1-1/n$ by \Cref{mainthm:nonhomogeneous}.

\subsubsection{Technical lemmas}

We next show some technical statements that we will need later in the proof of \Cref{mainthm:inconsistency}.

\begin{lemma}
\label{lemma:split_matrix_diagonal_opnorm}
Let $\mM \in \R^{k \times k}$. Then,
\begin{align*}
    \opnorm{\mM} \le \max_{i \le k} \abs{\mM[i][i]} + k\max_{i \neq j}\abs{\mM[i][j]}.
\end{align*}
\end{lemma}
\begin{proof}[Proof of \Cref{lemma:split_matrix_diagonal_opnorm}]
For a matrix $\mN \in \R^{k \times k}$, it is easy to check that
\begin{align*}
    \opnorm{\mN} \le \fnorm{\mN} \le k\max_{i, j \le k} \abs{\mN[i][j]}.
\end{align*}
Next, let $\diag{\mM}$ denote the matrix that agrees with $\mM$ on the diagonal and is $0$ elsewhere. Notice that
\begin{align*}
    \opnorm{\mM} &\le \opnorm{\diag{\mM}} + \opnorm{\mM - \diag{\mM}} \le \max_{i \le k} \abs{\mM[i][i]} + k\max_{i \neq j}\abs{\mM[i][j]},
\end{align*}
completing the proof of \Cref{lemma:split_matrix_diagonal_opnorm}.
\end{proof}

\begin{lemma}
\label{lemma:laplacian_diff_decreasing}
Let $\eps_x$ be a constant where $0 \le \eps_x < x$. Let $\eps_y$ be defined similarly. The function $f(x,y)$ defined as
\begin{align*}
    f(x,y) &\coloneqq \frac{1}{\sqrt{x-\eps_x}\sqrt{y-\eps_y}} - \frac{1}{\sqrt{x}\sqrt{y}}
\end{align*}
is decreasing in $x$ and $y$.
\end{lemma}
\begin{proof}[Proof of \Cref{lemma:laplacian_diff_decreasing}]
It is enough to just check the inequality for $x$. We take the derivative of $f(x,y)$ with respect to $x$ and get
\begin{align*}
    \frac{1}{2}\inparen{-\frac{1}{(x-\eps_x)^{3/2}(y-\eps_y)^{1/2}}+\frac{1}{x^{3/2}y^{1/2}}} < 0,
\end{align*}
where the inequality follows from observing $0 < x-\eps_x \le x$ and similarly for $y$. This completes the proof of \Cref{lemma:laplacian_diff_decreasing}.
\end{proof}

\subsubsection{Proof of \texorpdfstring{\Cref{mainthm:inconsistency}}{Theorem 3}}

First, we construct $\cLstar$.

\begin{lemma}
\label{lemma:hard_case_lstar}
Let $\cLstar \coloneqq \mI-\inparen{\Dstar}^{-1/2}\Astar\inparen{\Dstar}^{-1/2}$. Then, $\mI-\cLstar$ has the following block structure.
\begin{table}[H]
\centering
\resizebox{0.98\columnwidth}{!}{
\begin{tabular}{l|cccl}
                     & \multicolumn{1}{l}{$L_1$}                                                                              & \multicolumn{1}{l}{$L_2$}                                                                                                  & \multicolumn{2}{l}{$R$}                                                                                                                                                            \\ \hline
$L_1$                & $\frac{Kp}{\frac{n}{2} \cdot \inparen{p \cdot \frac{K+1}{2} + q}} \cdot \mathbbm{1}_{n/4 \times n/4}$  & \multicolumn{1}{c|}{$\frac{p}{\frac{n}{2} \cdot \inparen{p \cdot \frac{K+1}{2} + q}} \cdot \mathbbm{1}_{n/4 \times n/4}$}  & \multicolumn{2}{c}{\multirow{2}{*}{\smash{$\frac{q}{\sqrt{\frac{n}{2} \cdot \inparen{p \cdot \frac{K+1}{2} + q}\cdot\frac{n}{2}\cdot\inparen{p+q}}} \cdot \mathbbm{1}_{n/2 \times n/2}$}}} \\
$L_2$                & $\frac{p}{{\frac{n}{2} \cdot \inparen{p \cdot \frac{K+1}{2} + q}}} \cdot \mathbbm{1}_{n/4 \times n/4}$ & \multicolumn{1}{c|}{$\frac{Kp}{\frac{n}{2} \cdot \inparen{p \cdot \frac{K+1}{2} + q}} \cdot \mathbbm{1}_{n/4 \times n/4}$} & \multicolumn{2}{c}{}                                                                                                                                                               \\ \cline{2-5} 
\multirow{2}{*}{$R$} & \multicolumn{2}{c|}{\multirow{2}{*}{$\frac{q}{\sqrt{\frac{n}{2} \cdot \inparen{p \cdot \frac{K+1}{2} + q}\cdot\frac{n}{2}\cdot\inparen{p+q}}} \cdot \mathbbm{1}_{n/2 \times n/2}$}}                                                 & \multicolumn{2}{c}{\multirow{2}{*}{\smash{$\frac{p}{\frac{n}{2} \cdot \inparen{p+q}} \cdot \mathbbm{1}_{n/2 \times n/2}$}}}                                                                \\
                     & \multicolumn{2}{c|}{}                                                                                                                                                                                                               & \multicolumn{2}{c}{}                                                                                                                                                              
\end{tabular}}
\end{table}
\end{lemma}
\begin{proof}[Proof of \Cref{lemma:hard_case_lstar}]
It is easy to see that for any $v \in L$, we have $\vd^{\star}[v] = \frac{n}{2} \cdot \inparen{p \cdot \frac{K+1}{2} + q}$ and for any $v \in R$, we have $\vd^{\star}[v] = \frac{n}{2} \cdot \inparen{p + q}$. \Cref{lemma:hard_case_lstar} follows by noting that every element of $\mI-\cLstar$ is of the form $\astar_i[j]/\sqrt{\vd^{\star}[i]\vd^{\star}[j]}$.
\end{proof}

Next, we analyze the eigenvalues and eigenvectors of $\cLstar$.

\begin{lemma}
\label{lemma:hard_case_lstar_evects}
Up to normalization and sign, the eigenvector-eigenvalue pairs of $\mI-\cLstar$ corresponding to the nonzero eigenvalues of $\mI-\cLstar$ are
\begin{align*}
    (\lambda_1^{\star},\vu_1^{\star}) &= \inparen{1, \insquare{\onev_{n/4} \oplus \onev_{n/4} \oplus y_{+} \cdot \onev_{n/4} \oplus y_{+} \cdot \onev_{n/4}}} \\
    (\lambda_2^{\star}, \utwostar) &= \inparen{\frac{(K-1)p}{2\inparen{p \cdot \frac{K+1}{2}+q}} ,\insquare{\onev_{n/4} \oplus -\onev_{n/4} \oplus 0_{n/4} \oplus 0_{n/4}}} \\
    (\lambda_3^{\star}, \vu_3^{\star}) &= \inparen{-1+p\inparen{\frac{1}{p+q}+\frac{K+1}{p(K+1)+2q}}, \insquare{\onev_{n/4} \oplus \onev_{n/4} \oplus y_{-} \cdot \onev_{n/4} \oplus y_{-} \cdot \onev_{n/4}}}
\end{align*}
where $y_{+}$ and $y_{-}$ are chosen according to the formulas
\begin{align*}
    y_{+} &= \sqrt{\frac{2(p+q)}{p(K+1)+2q}} & y_{-} &= -\sqrt{\frac{p(K+1)+2q}{2(p+q)}}.
\end{align*}
Moreover, we have $\lambda_1^{\star} > \lambda_2^{\star} > \lambda_3^{\star} > 0$ and
\begin{align*}
    \lambda_2^{\star}-\lambda_3^{\star} \ge 1 - \frac{p^2(K+3)+4pq}{p^2(K+3)+4pq+2q^2}.
\end{align*}
\end{lemma}
\begin{proof}[Proof of \Cref{lemma:hard_case_lstar_evects}]
As we can see from \Cref{lemma:hard_case_lstar}, $\mI-\cLstar$ is a matrix whose rank is at most $3$, since it can be constructed by carefully repeating $3$ distinct column vectors. Thus, it can have at most $3$ nonzero eigenvalues. In what follows, we consider the case where $K > 1$ so that there are exactly $3$ nonzero eigenvalues. 

The next step is to confirm that the stated eigenvalue-eigenvector pairs are in fact valid. We begin with $\vu_1^{\star}$. Every entry in the first $n/2$ entries of $(\mI-\cLstar)\vu_1^{\star}$ can be expressed as
\begin{align*}
    &\quad\frac{n}{4} \cdot \frac{Kp}{\frac{n}{2} \cdot \inparen{p \cdot \frac{K+1}{2}+q}} + \frac{n}{4} \cdot \frac{p}{\frac{n}{2} \cdot \inparen{p \cdot \frac{K+1}{2}+q}} + \frac{n}{2} \cdot \inparen{\frac{q\cdot\sqrt{\frac{2(p+q)}{p(K+1)+2q}}}{\sqrt{\frac{n}{2} \cdot \inparen{p \cdot \frac{K+1}{2} + q}\cdot\frac{n}{2}\cdot\inparen{p+q}}}} \\
    &= \frac{(K+1)p}{(K+1)p+2q} + \frac{q\cdot\sqrt{\frac{2(p+q)}{p(K+1)+2q}}}{\sqrt{\inparen{p \cdot \frac{K+1}{2} + q}\inparen{p+q}}} = \frac{(K+1)p}{(K+1)p+2q} + \frac{q\cdot\sqrt{\frac{2}{p(K+1)+2q}}}{\sqrt{\inparen{p \cdot \frac{K+1}{2} + q}}} \\
    &= \frac{(K+1)p}{(K+1)p+2q} + \frac{2q}{(K+1)p+2q} = 1,
\end{align*}
and every entry in the second $n/2$ entries of $(\mI-\cLstar)\vu_1^{\star}$ can be expressed as
\begin{align*}
    &\quad \frac{n}{2} \cdot \frac{q}{\sqrt{\frac{n}{2} \cdot \inparen{p \cdot \frac{K+1}{2} + q}\cdot\frac{n}{2}\cdot\inparen{p+q}}} + \frac{n}{2} \cdot \frac{p}{\frac{n}{2} \cdot (p+q)} \cdot \sqrt{\frac{2(p+q)}{p(K+1)+2q}} \\
    &= \frac{q}{\sqrt{\inparen{p \cdot \frac{K+1}{2} + q}\inparen{p+q}}} +\frac{p}{(p+q)} \cdot \sqrt{\frac{2(p+q)}{p(K+1)+2q}} \\
    &= \frac{q}{\sqrt{\inparen{p \cdot \frac{K+1}{2} + q}\inparen{p+q}}} + p\cdot \sqrt{\frac{1}{\inparen{p \cdot \frac{K+1}{2} + q}\inparen{p+q}}} \\
    &= \frac{\sqrt{p+q}}{\sqrt{p \cdot \frac{K+1}{2} + q}} = \sqrt{\frac{2(p+q)}{p(K+1)+2q}} = y_{+}.
\end{align*}
For $\vu_2^{\star}$, we can use the block structure and easily verify
\begin{align*}
    \inparen{\mI-\cLstar}\utwostar = \frac{n}{4} \cdot \frac{(K-1)p}{\frac{n}{2} \cdot \inparen{p \cdot \frac{K+1}{2}+q}} \insquare{\onev_{n/4} \oplus -\onev_{n/4} \oplus 0_{n/4} \oplus 0_{n/4}} = \lambda_2^{\star}\utwostar.
\end{align*}
We now address $\vu_3^{\star}$. The first $n/2$ entries of $(\mI-\cLstar)\vu_3^{\star}$ are
\begin{align*}
    &\quad\frac{n}{4} \cdot \frac{Kp}{\frac{n}{2} \cdot \inparen{p \cdot \frac{K+1}{2}+q}} + \frac{n}{4} \cdot \frac{p}{\frac{n}{2} \cdot \inparen{p \cdot \frac{K+1}{2}+q}} + \frac{n}{2} \cdot \inparen{\frac{q\cdot-\sqrt{\frac{p(K+1)+2q}{2(p+q)}}}{\sqrt{\frac{n}{2} \cdot \inparen{p \cdot \frac{K+1}{2} + q}\cdot\frac{n}{2}\cdot\inparen{p+q}}}} \\
    &= \frac{(K+1)p}{(K+1)p+2q} + \inparen{\frac{q\cdot-\sqrt{\frac{1}{p+q}}}{\sqrt{p+q}}} = \frac{(K+1)p}{(K+1)p+2q} - \frac{q}{p+q} = \lambda_3^{\star},
\end{align*}
and the second $n/2$ entries of $(\mI-\cLstar)\vu_3^{\star}$ are
\begin{align*}
    &\quad \frac{n}{2} \cdot \frac{q}{\sqrt{\frac{n}{2} \cdot \inparen{p \cdot \frac{K+1}{2} + q}\cdot\frac{n}{2}\cdot\inparen{p+q}}} + \frac{n}{2} \cdot \frac{p}{\frac{n}{2}\cdot\inparen{p+q}} \cdot -\sqrt{\frac{p(K+1)+2q}{2(p+q)}} \\
    &= \frac{q}{\sqrt{\inparen{p \cdot \frac{K+1}{2} + q}\inparen{p+q}}} - \frac{p}{\inparen{p+q}} \cdot \sqrt{\frac{p(K+1)+2q}{2(p+q)}} \\
    &= -\sqrt{\frac{p(K+1)+2q}{2(p+q)}}\inparen{\frac{-2q}{p(K+1)+2q}+\frac{p}{p+q}} = y_{-} \cdot \lambda_3^{\star}.
\end{align*}
Finally, it remains to check that $1 > \lambda_2^{\star} > \lambda_3^{\star} > 0$. The fact that $\lambda_2^{\star} < 1$ easily follows from using $p + q > 0$. To prepare to bound $\lambda_2^{\star}-\lambda_3^{\star}$, we first use $p \ge q$ to establish
\begin{align*}
    p^2-pq+2q^2 = p(p-q)+2q^2 \ge 2q^2.
\end{align*}
This implies
\begin{align*}
   pq(K-1)+2q^2 \ge 3p^2-pq+2q^2 = 2p^2 + (p^2-pq+2q^2) \ge 2p^2+2q^2,
\end{align*}
which rearranges to
\begin{align*}
    p^2(K+1)+pq(K+3)+2q^2 \ge p^2(K+3)+4pq + 2q^2.
\end{align*}
Next, we write
\begin{align*}
    \lambda_2^{\star} - \lambda_3^{\star} &= \inparen{\frac{(K-1)p}{2\inparen{p \cdot \frac{K+1}{2}+q}}} - \inparen{-1+p\inparen{\frac{1}{p+q}+\frac{K+1}{p(K+1)+2q}}} \\
    &= 1 - \frac{p}{p+q}-\frac{2p}{p(K+1)+2q} = 1 - \inparen{\frac{p^2(K+1)+2pq+2p^2+2pq}{(p+q)(p(K+1)+2q)}} \\
    &= 1 - \frac{p^2(K+3)+4pq}{p^2(K+1)+pq(K+3)+2q^2} \ge 1 - \frac{p^2(K+3)+4pq}{p^2(K+3)+4pq+2q^2} > 0.
\end{align*}
Finally, to show $\lambda_3^{\star} > 0$, we write
\begin{align*}
    \lambda_3^{\star}+1 = \frac{p}{p+q} + \frac{p(K+1)}{p(K+1)+2q} > \frac{2p}{p+q} > 1,
\end{align*}
which allows us to complete the proof of \Cref{lemma:hard_case_lstar_evects}.
\end{proof}

Next, we argue that studying $\cL^{\star}$, which is formed by taking into account the weighted self-loops, gives us an understanding that is not too far from that of $\cL^{\star}_{\mathsf{nl}}$, which is formed by setting $p_{vv} = 0$ for all $v \in V$.

\begin{lemma}
\label{lemma:normalized_delete_sl}
Let $\mP$ be the diagonal matrix where $\mP[v,v]=p_{vv}$. Let $\cL^{\star}_{\mathsf{nl}}$ be the normalized Laplacian of the graph formed by $\Astar - \mP$. Then, we have
\begin{align*}
    \opnorm{\cL^{\star}-\cL^{\star}_{\mathsf{nl}}} \le \frac{6K}{n-2}. 
\end{align*}
\end{lemma}
\begin{proof}[Proof of \Cref{lemma:normalized_delete_sl}]
Recall $\Lstar \coloneqq \mD^{\star}-\Astar$. Let $\mD^{\star}_{\mathsf{nl}}$ be defined analogously to $\cL^{\star}_{\mathsf{nl}}$. Observe that we have
\begin{align*}
    \cL^{\star} &= \inparen{\mD^{\star}}^{-1/2}\Lstar\inparen{\mD^{\star}}^{-1/2} \\
    \cL^{\star}_{\mathsf{nl}} &= \inparen{\mD^{\star}_{\mathsf{nl}}}^{-1/2}\Lstar\inparen{\mD^{\star}_{\mathsf{nl}}}^{-1/2}.
\end{align*}
From this, we see that writing down the $v,w$th entry of the difference gives
\begin{align*}
    \inparen{\cL^{\star}_{\mathsf{nl}}-\cL^{\star}}[v,w] &= \Lstar[v,w]\inparen{\frac{1}{\sqrt{(\vd^{\star}[v]-p_{vv})(\vd^{\star}[w]-p_{ww})}}-\frac{1}{\sqrt{\vd^{\star}[v]\vd^{\star}[w]}}}.
\end{align*}
This resolves to different forms based on whether $v=w$. When $v=w$, evaluating the formula gives
\begin{align*}
    \inparen{\cL^{\star}_{\mathsf{nl}}-\cL^{\star}}[v,v] = \frac{\vd^{\star}[v]-p_{vv}}{\vd^{\star}[v]-p_{vv}} - \frac{\vd^{\star}[v]-p_{vv}}{\vd^{\star}[v]} = \frac{p_{vv}}{\vd^{\star}[v]}.
\end{align*}
When $v\neq w$, we apply \Cref{lemma:laplacian_diff_decreasing} and get
\begin{align*}
    \abs{\inparen{\cL^{\star}_{\mathsf{nl}}-\cL^{\star}}[v,w]} &= p_{vw}\inparen{\frac{1}{\sqrt{(\vd^{\star}[v]-p_{vv})(\vd^{\star}[w]-p_{ww})}}-\frac{1}{\sqrt{\vd^{\star}[v]\vd^{\star}[w]}}} \\
    &\le Kp\inparen{\frac{1}{np/2-p}-\frac{1}{np/2}} = \frac{4K}{n^2-2n}.
\end{align*}
Using this analysis and applying \Cref{lemma:split_matrix_diagonal_opnorm} gives
\begin{align*}
    \opnorm{\cL^{\star}_{\mathsf{nl}}-\cL^{\star}} &\le \max_{v \in V} \frac{p_{vv}}{\vd^{\star}[v]} + n\max_{v\neq w} p_{vw}\inparen{\frac{1}{\sqrt{(\vd^{\star}[v]-p_{vv})(\vd^{\star}[w]-p_{ww})}}-\frac{1}{\sqrt{\vd^{\star}[v]\vd^{\star}[w]}}} \\
    &\le \frac{2K}{n} + n\max_{v\neq w} p_{vw}\inparen{\frac{1}{\sqrt{(\vd^{\star}[v]-p_{vv})(\vd^{\star}[w]-p_{ww})}}-\frac{1}{\sqrt{\vd^{\star}[v]\vd^{\star}[w]}}} \\
    &\le \frac{2K}{n} + \frac{4K}{n-2} \le \frac{6K}{n-2},
\end{align*}
completing the proof of \Cref{lemma:normalized_delete_sl}.
\end{proof}

This gives \Cref{lemma:normalized_delete_sl_u2}, which means we can use $\utwostar$ as a suitable proxy for $\signv{\vu_2(\cL_{\mathsf{nl}}^{\star})}$.

\begin{lemma}
\label{lemma:normalized_delete_sl_u2}
There exists a constant $C(\alpha,K)$ depending on $\alpha$ and $K$ such that we have
\begin{align*}
    \norm{\vu_2(\cL_{\mathsf{nl}}^{\star}) - \utwostar}_{\infty} \le \frac{C(\alpha,K)}{n}.
\end{align*}
This implies that for all $n$ sufficiently large, we have $\signv{\vu_2(\cL_{\mathsf{nl}}^{\star})} = \signv{\utwostar}$.
\end{lemma}
\begin{proof}[Proof of \Cref{lemma:normalized_delete_sl_u2}]
By \Cref{lemma:hard_case_lstar_evects}, Weyl's inequality, and \Cref{lemma:normalized_delete_sl}, we know that for all  $n$ sufficiently large,
\begin{align*}
    \lambda_2^{\star} - \lambda_3(\cL_{\mathsf{nl}}^{\star}) &= \inparen{\lambda_2^{\star}-\lambda_3^{\star}} + (\lambda_3^{\star}-\lambda_3(\cL_{\mathsf{nl}}^{\star})) \\
    &\ge \inparen{1 - \frac{p^2(K+3)+4pq}{p^2(K+3)+4pq+2q^2}} - \frac{C_{\ref{lemma:normalized_delete_sl}}}{n} \ge C_1(\alpha,K).
\end{align*}
Combining this with \Cref{lemma:normalized_delete_sl} again, the Davis-Kahan inequality tells us that
\begin{align*}
    \norm{\vu_2(\cL_{\mathsf{nl}}^{\star}) - \utwostar}_{\infty} \le \norm{\vu_2(\cL_{\mathsf{nl}}^{\star}) - \utwostar}_{2} \le \frac{C_2(\alpha,K)}{n},
\end{align*}
and then using the fact that $\maxnorm{\utwostar} = 1/\sqrt{n}$ (arising from \Cref{lemma:hard_case_lstar_evects}) completes the proof of \Cref{lemma:normalized_delete_sl_u2}.
\end{proof}

We are now ready to prove the inconsistency of normalized spectral bisection on the nested block examples. 

\begin{proof}[Proof of \Cref{mainthm:inconsistency}]

Let $G$ be a graph drawn from the nested block example. We choose $p$ and $q$ such that $p \gtrsim \log n / n$ and $p/q = \alpha \ge 2$ where $\alpha$ is some constant and such that $p$ and $q$ both satisfy the conditions of \Cref{mainthm:nonhomogeneous}. Let $K \ge 3\alpha$. Observe that the true communities are $L$ and $R$. We will show that bisection based on $\vu_2$ of $\mI-\cL$ (corresponding to the eigenvector associated with the second smallest eigenvalue of $\cL$) will attain a large misclassification rate. In particular, based on our calculation in \Cref{lemma:hard_case_lstar_evects}, we expect that $\vu_2$ will output a bisection that places $L_1$ and $L_2$ into separate clusters. On the other hand, by \Cref{mainthm:nonhomogeneous}, for all $n$ large enough, the unnormalized spectral bisection algorithm will be strongly consistent.

First, observe that it is enough to prove the inconsistency result just for the symmetric normalized Laplacian. Indeed, observe that if $\vu_2$ is an eigenvector of $\mI-\cL = \mD^{-1/2}\mA\mD^{-1/2}$, then we have
\begin{align*}
    \lambda_2\mD^{-1/2}\vu_2 = \mD^{-1}\mA\mD^{-1/2}\vu_2 = \mD^{-1}\mA(\mD^{-1/2}\vu_2),
\end{align*}
which shows that $\mD^{-1/2}\vu_2$ must be the eigenvector of the random-walk normalized Laplacian $\mI-\mD^{-1}\mA$ corresponding to eigenvalue $\lambda_2$. Since $\mD$ is a positive diagonal matrix, it does not change the signs of $\vu_2$ and therefore the output of the normalized spectral bisection algorithm is the same.

Our general approach to prove the inconsistency is to use the Davis-Kahan Theorem, a bound on $\opnorm{\cL-\cLstar_{\mathsf{nl}}}$, and a bound on the gap $\lambda_2^{\star}-\lambda_3$. Let $\vd_{\min}$ be the minimum degree of the graph given by adjacency matrix $\mA$ and let $\vd_{\min}^{\star}$ be the minimum weighted degree of the graph given by the adjacency matrix $\mA^{\star}$. First, using \cite[Theorem 3.1]{dls20}, we have with probability $1-n^{-r}$ for some constant $r \ge 1$ and constants $C(r)$ and $C$ (the latter of which does not depend on $r$), for all $n$ sufficiently large,
\begin{align*}
    \opnorm{\cL-\cLstar_{\mathsf{nl}}} &\le \frac{C(r)\inparen{n\max_{(i,j)} p_{ij}}^{5/2}}{\min\inbraces{\vd_{\min},\vd_{\min}^{\star}}^3} \\
    &\le \frac{C(r)\inparen{n \cdot Kp}^{5/2}}{\min\inbraces{n(p+q)/3, n(p+q)/3-C\sqrt{n(p+q)\log n}}^3} \\
    &\le \frac{C_1(r,\alpha)K^{5/2}(np)^{5/2}}{\inparen{np}^3} = \frac{C_1(r,\alpha)K^{5/2}}{\sqrt{np}}.
\end{align*}
Next, we invoke \Cref{lemma:normalized_delete_sl} to write
\begin{align*}
    \lambda_2(\cL_{\mathsf{nl}}^{\star}) - \lambda_3 &= \inparen{\lambda_2^{\star}-\lambda_3^{\star}} + (\lambda_3^{\star}-\lambda_3) + \inparen{\lambda_2(\cL_{\mathsf{nl}}^{\star}) - \lambda_2^{\star}} \\
    &\ge \inparen{1 - \frac{p^2(K+3)+4pq}{p^2(K+3)+4pq+2q^2}} - \frac{C_2(r,\alpha)K^{5/2}}{\sqrt{np}} - \frac{C_{\ref{lemma:normalized_delete_sl}}}{n} \ge C_g(\alpha,K),
\end{align*}
where the last line denotes a positive constant depending on $q$ and $K$ (this constant will always be positive for sufficiently large $n$, as we showed that $\lambda_2^{\star}-\lambda_3^{\star} > 0$ in \Cref{lemma:hard_case_lstar_evects}).

Putting everything together, we get by the Davis-Kahan theorem that some signing of $\vu_2$ satisfies
\begin{align*}
    \norm{\vu_2-\vu_2(\cL_{\mathsf{nl}}^{\star})}_2 \le \frac{\opnorm{\cL-\cLstar_{\mathsf{nl}}}}{\min\inbraces{\abs{\lambda_2(\cLstar_{\mathsf{nl}})-\lambda_3},1-\lambda_2(\cLstar_{\mathsf{nl}})}} \le \frac{C_3(r)K^{5/2}}{C_g'(\alpha,K)\sqrt{np}} \le \frac{C_4(r,\alpha,K)}{\sqrt{np}}.
\end{align*}
Now, consider the subset of coordinates of $\vu_2$ belonging to $L_1$. Suppose $m$ of these coordinates do not agree in sign with $\utwostar$. To maximize $m$, each of these coordinates in $\vu_2$ should be $0$, so using this reasoning and applying \Cref{lemma:normalized_delete_sl_u2} means the total $\ell_2$ error can be bounded (using \Cref{lemma:normalized_delete_sl}) as
\begin{align*}
    m\inparen{\frac{1}{\sqrt{\nfrac{n}{2}}}-\frac{C_{\ref{lemma:normalized_delete_sl}}}{n}}^2 \le \norm{\vu_2-\vu_2(\cL_{\mathsf{nl}}^{\star})}_2^2 \le \frac{C_4(r,\alpha,K)^2}{np}.
\end{align*}
This means the number of coordinates $m$ on which $\vu_2$ and $\utwostar$ disagree on is at most
\begin{align*}
    \frac{n \cdot C_5(r,\alpha,K)^2}{2np},
\end{align*}
and therefore the misclassification rate of $\vu_2$ with respect to the true labeling induced by $L$ and $R$ must be at least
\begin{align*}
    \frac{\frac{n}{4}-\frac{n \cdot C_5(r,\alpha,K)^2}{2np}}{n} = \frac{1}{4} - \frac{C_5(r,\alpha,K)^2}{2np}.
\end{align*}
Since $p \gtrsim \log n / n$, this completes the proof of \Cref{mainthm:inconsistency}.
\end{proof}

\section{Additional experiments}
\label{sec:cluster_experiments_more}
In this section, we show more numerical trials that complement those discussed in \cref{sec:cluster_experiments}.

\subsection{Varying edge probabilities in an NSSBM}
\label{subsec:cluster_experimentnssbm}
 In \Cref{sec:cluster_experiments}, we investigated the behavior of an NSSBM model by fixing the values of $p,q$ and varying the largest edge probability $\pbar$. Here, we take an alternative approach, and instead fix $\pbar$ and vary the values of $p$ and $q$.

\begin{figure}
    \centering
    \includegraphics[width=0.45\columnwidth]{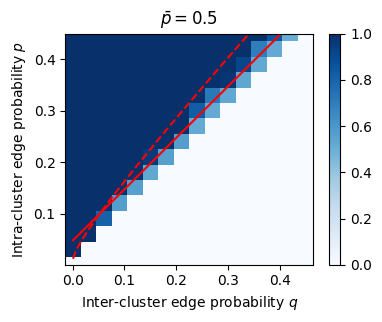}
    \includegraphics[width=0.45\columnwidth]{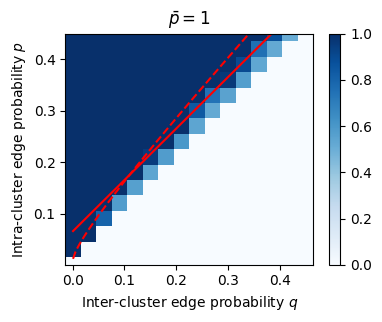}
    \caption{Agreement with the planted bisection of the bipartition obtained from unnormalized spectral bisection, for graphs generated from a distribution in $\mathsf{NSSBM}(n,p,\pbar,q)$ for fixed values of $n,\pbar$ and varying values of $p>q$. The left plot uses $\pbar=1/2$, the right plot uses $\pbar=1$. The solid red curves plot the function $p_{\mathsf{thr}}(q)$ (see~\eqref{eq:thr}), and the dashed red curves plot the function $p_{\mathsf{info}}(q)$ (see~\eqref{eq:info}).}
    \label{fig:varying_pairs}
\end{figure}

\smallpar{Setup.} Let us fix $n=2000$, $\pbar \in  \{1/2, 1\}$. For varying $p,q$ in the range $[1/n,9/20]$ such that $p > q$, we sample $t=3$ independent draws $G$ from the same benchmark distribution $\mathcal{D}_{p,\pbar,q}$ used in \cref{sec:cluster_experiments}. For each of them, we compute the agreement of the bipartition obtained by unnormalized spectral bisection with respect to the planted bisection. For each $(p,q)$, we plot the average agreement across the $t$ independent draws. The results are shown in \cref{fig:varying_pairs}, where in the left and right plot we ran the experiments with $\pbar=1/2$ and $\pbar=1$ respectively. The lower diagonal of these plots, where $p \le q$, is artificially set to $0$.

\smallpar{Theoretical framing.} According to \cref{mainthm:nonhomogeneous}, fixing the value of $\pbar \in \{1/2,1\}$, we obtain that unnormalized spectral bisection achieves exact recovery provided that for $q \in [1/n,9/20]$ one has $p \ge p_{\mathsf{thr}}(q)$ where
\begin{equation}
\label{eq:thr}
    p_{\mathsf{thr}}(q) = \frac{\sqrt{\pbar \log n}}{\sqrt{n}}+q
\end{equation}
is obtained by rearranging the precondition of \cref{mainthm:nonhomogeneous}, ignoring the constants, and disregarding the fact that $\alpha$ should be $O(1)$. The solid red curve in \cref{fig:varying_pairs} plots $p_{\mathsf{thr}}(q)$ as a function of $q$. For comparison, the information-theoretic threshold for SSBM \cite{abh16} demands that $p \ge p_{\mathsf{info}}(q)$ where
\begin{equation}
\label{eq:info}
    p_{\mathsf{info}}(q) = \inparen{\sqrt{2} \sqrt{\frac{\log n}{n}}+\sqrt{q}}^2 \, .
\end{equation}
The dashed red curve in \cref{fig:varying_pairs} plots $p_{\mathsf{info}}(q)$ as a function of $q$.

\smallpar{Empirical evidence.} From \cref{fig:varying_pairs}, one can see that our experiments reflect the behavior predicted by \cref{mainthm:nonhomogeneous} quite closely, although empirically we achieve $100\%$ agreement slightly above $p_{\mathsf{thr}}(q)$ (i.e. the solid red curve). However, this is likely due to the constant factors from \cref{mainthm:nonhomogeneous} that we ignored, and also $n=2000$ is plausibly too small to show asymptotic behaviors. Nevertheless, we do achieve $100\%$ agreement consistently as soon as we surpass the information-theoretic threshold $p_{\mathsf{info}}(q)$: in the regime of our experiment, it appears that the unnormalized Laplacian is robust all the way to the optimal threshold for exact recovery in the SSBM.

\subsection{Varying the size of  a planted clique in a DCM}
In some sense, the experiments from \Cref{sec:cluster_experiments} and \cref{subsec:cluster_experimentnssbm} can be thought of as experiments for the deterministic clusters model too. This is because each realization of the internal edges gives rise to a different DCM distribution (see \cref{sec:cluster_results}). We complement our previous discussion by illustrating the behavior of certain families of DCM distributions that are conceptually different than those considered in \Cref{sec:cluster_experiments}.

\smallpar{Benchmark distribution.} Let $n$ be divisible by $4$ and let $\{P_1,P-2\}$ be a partitioning of $V=[n]$ into two equally-sized subsets. Fix $ p \in [0,1]$. For some set $S \subseteq P_1$ such that $S=\{1,\dots,|S|\}$ (for simplicity), let $G_2 = (P_2,E_2)\sim \mathsf{ER}(n/2,p)$ be a graph drawn from the Erd\H{o}s-R\'enyi distribution with sampling rate~$p$, and let $G_1 = (P_1,E_1)\sim \mathsf{ERPC}(n/2,p,S)$ be also a graph drawn from the Erd\H{o}s-R\'enyi distribution with sampling rate $p$ where we additionally plant a clique on the vertices $S$. Fixing $G_1,G_2$, for $q \in [0,1]$ we consider the distribution \smash{$\mathcal{D}_{q}^{G_1,G_2}$} over graphs $G=(V,E)$ where $G[P_1]=G_1$, $G[P_2]=G_2$, and every edge $(u,v) \in P_1 \times P_2$ is sampled independently with probability $q$. One can see that $\mathcal{D}_{q}^{G_1,G_2}$ is in fact in the set $\mathsf{DCM}(n,d_{\mathsf{in}},q)$ for some $d_{\mathsf{in}}$.

\begin{figure}
    \centering
    \includegraphics[width=0.45\columnwidth]{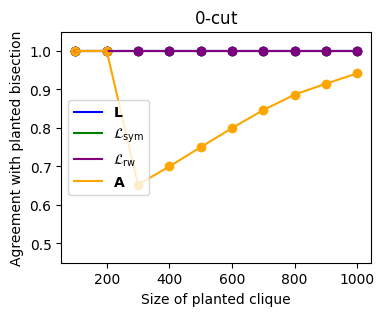}
    \includegraphics[width=0.45\columnwidth]{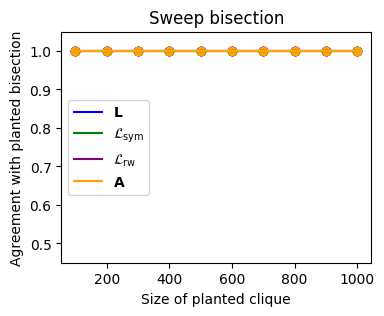}
    \caption{Agreement with the planted bisection of the bipartition obtained from several matrices associated with an input graph generated from a distribution $\mathcal{D}_{q}^{G_1,G_2} \in \mathsf{DCM}(n,d_{\mathsf{in}},q)$ for fixed values of $n,q$ and varying the size of the planted clique $S$. In the left plot, the bipartition is the $0$-cut of the second eigenvector, as in \cref{alg:spectral_general}. In the right plot, the bipartition is the sweep cut of the first $n/2$ vertices in the second eigenvector.}
    \label{fig:varying_clique_cuts}
\end{figure}

\smallpar{Setup.} Let us fix $n=2000$, $p=9/\sqrt{n}$, $q = 1/\sqrt{n}$. For varying values of $|S|$ in the range $[|P_1|/10,|P_1|]$, we sample $G_1 = (P_1,E_1)\sim \mathsf{ERPC}(n/2,p,S)$ and $G_2 = (P_2,E_2)\sim \mathsf{ER}(n/2,p)$, and then draw $t=10$ independent samples $G$ from $\mathcal{D}_{q}^{G_1,G_2}$. For each sample $G$, we run spectral bisection (i.e. \cref{alg:spectral_general}) with matrices $\mL,\cL_{\mathsf{sym}},\cL_{\mathsf{rw}},\mA$. Then, we compute the {agreement} of the bipartition hence obtained with respect to the planted bisection, and average it out across the $t$ independent draws. The results are shown in the left plot of \cref{fig:varying_clique_cuts}. Again, another natural way to get a bipartition of $V$ from the eigenvector is a {sweep cut}, and the average agreements that this results in are shown in the right plot of \cref{fig:varying_clique_cuts}.

\begin{figure}
    \centering
    \includegraphics[width=0.45\columnwidth]{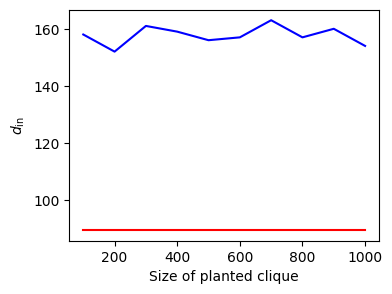}
    \includegraphics[width=0.45\columnwidth]{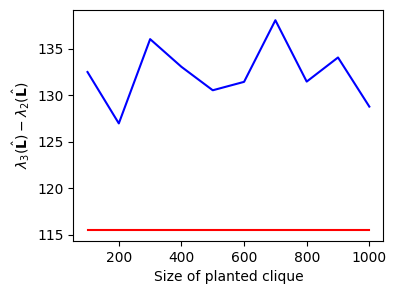}
    \caption{The minimum in-cluster degree $d_{\mathsf{in}}$ and the spectral gap $\lambda_3(\widehat{\mL})-\lambda_2(\widehat{\mL})$ of distributions \smash{$\mathcal{D}_{q}^{G_1,G_2} \in \mathsf{DCM}(n,d_{\mathsf{in}},q)$} with fixed values of $n,q$ and varying the size of the planted clique $S$. The red horizontal line on the left corresponds to the value $nq+\sqrt{n}$, the red horizontal line on the right corresponds to the value $\sqrt{n} + nq+\sqrt{nq \log n} +\log n$.}
    \label{fig:varying_clique_params}
\end{figure}

\smallpar{Theoretical framing.} Ignoring the constants, \cref{mainthm:sqrtngap} guarantees that exact recovery is achieved by unnormalized spectral bisection as long as \smash{$d_{\mathsf{in}} \ge nq+\sqrt{n}$} and \smash{$\lambda_3(\widehat{\mL})-\lambda_2(\widehat{\mL}) \ge \sqrt{n} + nq+\sqrt{nq \log n} +\log n$}, where \smash{$\widehat{\mL}$} is the expected Laplacian of $\mathcal{D}_{q}^{G_1,G_2}$. For each clique size that we consider, \cref{fig:varying_clique_params} shows the minimum in-cluster degree of the graphs $G_1,G_2$ that we draw (in the left plot), and the spectral gap \smash{$\lambda_3(\widehat{\mL})-\lambda_2(\widehat{\mL})$}. The red horizontal lines in the left and right plot respectively correspond to the value of $nq+\sqrt{n}$ and $\sqrt{n} + nq+\sqrt{nq \log n} +\log n$ on the $y$-axis, indicating the lower bound on $d_{\mathsf{in}}$ and \smash{$\lambda_3(\widehat{\mL})-\lambda_2(\widehat{\mL})$} demanded by \cref{mainthm:sqrtngap}.

\smallpar{Empirical evidence: consistency.} From \cref{fig:varying_clique_params}, one can see that all the distributions $\mathcal{D}_{q}^{G_1,G_2}$ that we use roughly meet the requirement of \cref{mainthm:sqrtngap}. Indeed, in the left plot of \cref{fig:varying_clique_cuts} one sees that unnormalized spectral bisection consistently achieves exact recovery for all clique sizes. On the contrary, the bipartition obtained by running spectral bisection with the adjacency matrix $\mA$ misclassifies a fraction of the vertices for certain sizes of the planted clique. Nevertheless, the sweep cut obtained from all the matrices recovers the planted bisection exactly.

\begin{figure}[H]
    \centering
    \includegraphics[width=0.45\columnwidth]{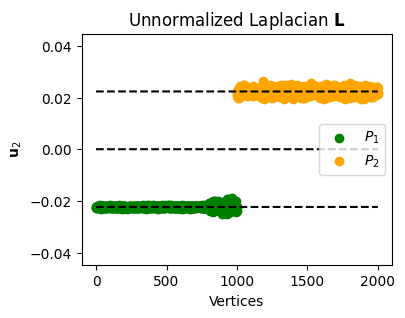}
    \includegraphics[width=0.45\columnwidth]{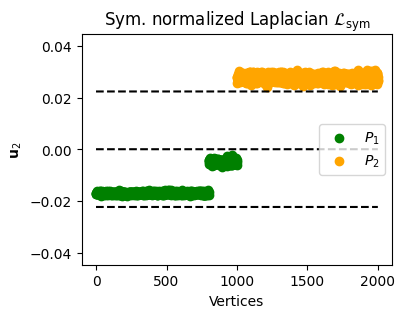}
    \includegraphics[width=0.45\columnwidth]{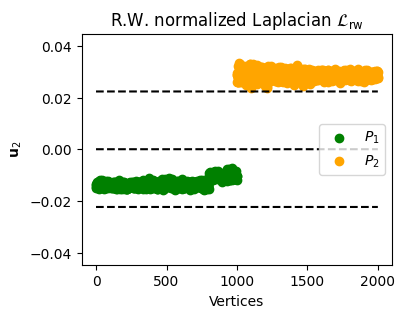}
    \includegraphics[width=0.45\columnwidth]{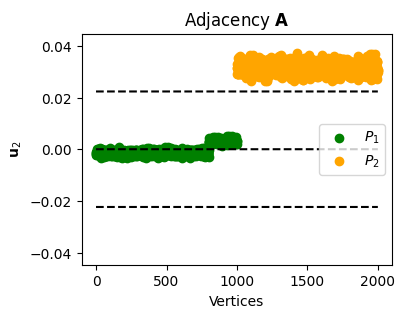}
    \caption{Embedding of the vertices given by the second eigenvector $\vu_2$ of several matrices associated with a graph sampled from a distribution \smash{$\mathcal{D}_{q}^{G_1,G_2} \in \mathsf{DCM}(n,d_{\mathsf{in}},q)$}, with the size of the planted clique set to $|S|=2/5 \cdot n$. Horizontal dashed lines, from top to bottom, correspond to $1/\sqrt{n}, 0, -1/\sqrt{n}$ respectively.}
    \label{fig:embedding}
\end{figure}

\smallpar{Empirical evidence: example embedding.} Let us fix the value $|S|=800$ for the size of the planted clique, for which we see in \cref{fig:varying_clique_cuts} that the adjacency matrix fails to recover the planted bisection. We generate a graph from a distribution $\mathcal{D}_{q}^{G_1,G_2}$ with clique size $|S|=800$, and plot how the vertices are embedded in the real line by the second eigenvector of all the matrices we consider. The result is shown in \cref{fig:embedding}, where the three horizontal dashed lines, from top to bottom, respectively correspond to the value of $1/\sqrt{n}, 0, -1/\sqrt{n}$ on the $y$-axis. Graphically, one can see that the embedding in the unnormalized Laplacian is indeed the one that moves the least away from the values $\pm 1/\sqrt{n}$, and in fact the vertices $\{1,\dots, 800\} \subseteq P_1$ where we plant the clique concentrate even more around $1/\sqrt{n}$. This is a phenomenon related to the one illustrated by \cref{fig:embedding_nssbm}. Finally, one can see from the embedding that splitting vertices around $0$ does result in misclassifying a fraction of the vertices for the adjacency matrix. However, taking a sweep cut that splits the vertices into two equally sized parts recovers the planted bisection for all matrices. This reflects the results shown in \cref{fig:varying_clique_cuts}.

\printbibliography[heading=bibintoc]

@book{dk_book,
  title={Algorithmic high-dimensional robust statistics},
  author={Diakonikolas, Ilias and Kane, Daniel M},
  year={2023},
  publisher={Cambridge university press}
}

@article{bn21,
  title={Noise in Classification},
  author={Balcan, Maria-Florina and Haghtalab, Nika},
  journal={Beyond the Worst-Case Analysis of Algorithms},
  pages={361},
  year={2021},
  publisher={Cambridge University Press}
}

@INPROCEEDINGS{Borgnia2021-vk,
  title     = "Strong Data Augmentation Sanitizes Poisoning and Backdoor
               Attacks Without an Accuracy Tradeoff",
  booktitle = "{ICASSP} 2021 - 2021 {IEEE} International Conference on
               Acoustics, Speech and Signal Processing ({ICASSP})",
  author    = "Borgnia, Eitan and Cherepanova, Valeriia and Fowl, Liam and
               Ghiasi, Amin and Geiping, Jonas and Goldblum, Micah and
               Goldstein, Tom and Gupta, Arjun",
  pages     = "3855--3859",
  month     =  jun,
  year      =  2021,
eprint={2011.09527},
      archivePrefix={arXiv},
      primaryClass={cs.CR},
}

@article{Feldman2020-ex,
  title={What neural networks memorize and why: Discovering the long tail via influence estimation},
  author={Feldman, Vitaly and Zhang, Chiyuan},
  journal={Advances in Neural Information Processing Systems},
  volume={33},
  pages={2881--2891},
  year={2020},
archivePrefix = "arXiv",
  primaryClass  = "cs.LG",
  eprint        = "2008.03703"
}

@INPROCEEDINGS{Arpit2017-yz,
  title     = "A Closer Look at Memorization in Deep Networks",
  booktitle = "Proceedings of the 34th International Conference on Machine
               Learning",
  author    = "Arpit, Devansh and Jastrz{\k e}bski, Stanis{\l}aw and Ballas,
               Nicolas and Krueger, David and Bengio, Emmanuel and Kanwal,
               Maxinder S and Maharaj, Tegan and Fischer, Asja and Courville,
               Aaron and Bengio, Yoshua and Lacoste-Julien, Simon",
  editor    = "Precup, Doina and Teh, Yee Whye",
  publisher = "PMLR",
  volume    =  70,
  pages     = "233--242",
  series    = "Proceedings of Machine Learning Research",
  year      =  2017,
eprint={1706.05394},
      archivePrefix={arXiv},
      primaryClass={stat.ML},
}

@ARTICLE{Gu2017-rj,
  title         = "{BadNets}: Identifying Vulnerabilities in the Machine
                   Learning Model Supply Chain",
  author        = "Gu, Tianyu and Dolan-Gavitt, Brendan and Garg, Siddharth",
  month         =  aug,
  year          =  2017,
  archivePrefix = "arXiv",
  primaryClass  = "cs.CR",
  eprint        = "1708.06733"
}

@article{lecun-mnisthandwrittendigit-2010,
  author = {LeCun, Yann and Cortes, Corinna},
  howpublished = {http://yann.lecun.com/exdb/mnist/},
  lastchecked = {2016-01-14 14:24:11},
  title = {{MNIST} handwritten digit database},
  url = {http://yann.lecun.com/exdb/mnist/},
  username = {mhwombat},
  year = 2010
}

@INPROCEEDINGS{Ribeiro2016-ig,
  title     = "``Why Should {I} Trust You?'': Explaining the Predictions of Any
               Classifier",
  booktitle = "Proceedings of the 22nd {ACM} {SIGKDD} International Conference
               on Knowledge Discovery and Data Mining",
  author    = "Ribeiro, Marco Tulio and Singh, Sameer and Guestrin, Carlos",
  publisher = "Association for Computing Machinery",
  pages     = "1135--1144",
  series    = "KDD '16",
  month     =  aug,
  year      =  2016,
  address   = "New York, NY, USA",
  location  = "San Francisco, California, USA",
eprint={1602.04938},
      archivePrefix={arXiv},
      primaryClass={cs.LG},
}

@inproceedings{
Xiao2020-ju,
title={Noise or Signal: The Role of Image Backgrounds in Object Recognition},
author={Kai Yuanqing Xiao and Logan Engstrom and Andrew Ilyas and Aleksander Madry},
booktitle={International Conference on Learning Representations},
year={2021},
url={https://openreview.net/forum?id=gl3D-xY7wLq},
eprint={2006.09994},
      archivePrefix={arXiv},
      primaryClass={cs.CV},
}

@INPROCEEDINGS{Ilyas2019-ot,
  title     = "Adversarial Examples Are Not Bugs, They Are Features",
  booktitle = "Advances in Neural Information Processing Systems",
  author    = "Ilyas, Andrew and Santurkar, Shibani and Tsipras, Dimitris and
               Engstrom, Logan and Tran, Brandon and Madry, Aleksander",
  editor    = "Wallach, H and Larochelle, H and Beygelzimer, A and
               d' Alch{\'e}-Buc, F and Fox, E and
               Garnett, R",
  publisher = "Curran Associates, Inc.",
  volume    =  32,
  year      =  2019,
eprint={1905.02175},
      archivePrefix={arXiv},
      primaryClass={stat.ML},
}

@book{boyd2004convex,
  title={Convex optimization},
  author={Boyd, Stephen P and Vandenberghe, Lieven},
  year={2004},
  publisher={Cambridge university press}
}

@article{Cullina2018-os,
  title={Pac-learning in the presence of adversaries},
  author={Cullina, Daniel and Bhagoji, Arjun Nitin and Mittal, Prateek},
  journal={Advances in Neural Information Processing Systems},
  volume={31},
  year={2018},
eprint={1806.01471},
      archivePrefix={arXiv},
      primaryClass={stat.ML},
}

@INPROCEEDINGS{Montasser2020-ir,
  title     = "Efficiently Learning Adversarially Robust Halfspaces with Noise",
  booktitle = "Proceedings of the 37th International Conference on Machine
               Learning",
  author    = "Montasser, Omar and Goel, Surbhi and Diakonikolas, Ilias and
               Srebro, Nathan",
  editor    = "Iii, Hal Daum{\'e} and Singh, Aarti",
  publisher = "PMLR",
  volume    =  119,
  pages     = "7010--7021",
  series    = "Proceedings of Machine Learning Research",
  year      =  2020,
eprint={2005.07652},
      archivePrefix={arXiv},
      primaryClass={cs.LG},
}

@misc{Turner2019-jc,
  title         = "{Label-Consistent} Backdoor Attacks",
  author        = "Turner, Alexander and Tsipras, Dimitris and Madry,
                   Aleksander",
  month         =  dec,
  year          =  2019,
  archivePrefix = "arXiv",
  primaryClass  = "stat.ML",
  eprint        = "1912.02771"
}

@ARTICLE{Gu2019-ip,
  title    = "{BadNets}: Evaluating Backdooring Attacks on Deep Neural Networks",
  author   = "Gu, Tianyu and Liu, Kang and Dolan-Gavitt, Brendan and Garg,
              Siddharth",
  journal  = "IEEE Access",
  volume   =  7,
  pages    = "47230--47244",
  year     =  2019
}

@misc{Li2020-my,
  title         = "Backdoor Learning: A Survey",
  author        = "Li, Yiming and Wu, Baoyuan and Jiang, Yong and Li, Zhifeng
                   and Xia, Shu-Tao",
  month         =  jul,
  year          =  2020,
  archivePrefix = "arXiv",
  primaryClass  = "cs.CR",
  eprint        = "2007.08745"
}

@inproceedings{Truong2020-dk,
  title={Systematic evaluation of backdoor data poisoning attacks on image classifiers},
  author={Truong, Loc and Jones, Chace and Hutchinson, Brian and August, Andrew and Praggastis, Brenda and Jasper, Robert and Nichols, Nicole and Tuor, Aaron},
  booktitle={Proceedings of the IEEE/CVF conference on computer vision and pattern recognition workshops},
  pages={788--789},
  year={2020},
archivePrefix = "arXiv",
  primaryClass  = "cs.CV",
  eprint        = "2004.11514"
}

@INPROCEEDINGS{Adi2018-fz,
  title     = "Turning your weakness into a strength: watermarking deep neural
               networks by backdooring",
  booktitle = "Proceedings of the 27th {USENIX} Conference on Security
               Symposium",
  author    = "Adi, Yossi and Baum, Carsten and Cisse, Moustapha and Pinkas,
               Benny and Keshet, Joseph",
  publisher = "USENIX Association",
  pages     = "1615--1631",
  series    = "SEC'18",
  month     =  aug,
  year      =  2018,
  address   = "USA",
  location  = "Baltimore, MD, USA",
eprint={1802.04633},
      archivePrefix={arXiv},
      primaryClass={cs.LG},
}

@misc{Chen2017-kq,
      title={Targeted Backdoor Attacks on Deep Learning Systems Using Data Poisoning}, 
      author={Xinyun Chen and Chang Liu and Bo Li and Kimberly Lu and Dawn Song},
      year={2017},
      eprint={1712.05526},
      archivePrefix={arXiv},
      primaryClass={cs.CR}
}

@inproceedings{Wang2020-yt,
 author = {Wang, Hongyi and Sreenivasan, Kartik and Rajput, Shashank and Vishwakarma, Harit and Agarwal, Saurabh and Sohn, Jy-yong and Lee, Kangwook and Papailiopoulos, Dimitris},
 booktitle = {Advances in Neural Information Processing Systems},
 editor = {H. Larochelle and M. Ranzato and R. Hadsell and M. F. Balcan and H. Lin},
 pages = {16070--16084},
 publisher = {Curran Associates, Inc.},
 title = {Attack of the Tails: Yes, You Really Can Backdoor Federated Learning},
 volume = {33},
 year = {2020},
eprint={2007.05084},
      archivePrefix={arXiv},
      primaryClass={cs.LG},
}

@misc{Madry2017-ep,
  title         = "Towards Deep Learning Models Resistant to Adversarial
                   Attacks",
  author        = "Madry, Aleksander and Makelov, Aleksandar and Schmidt,
                   Ludwig and Tsipras, Dimitris and Vladu, Adrian",
  month         =  jun,
  year          =  2017,
  archivePrefix = "arXiv",
  primaryClass  = "stat.ML",
  eprint        = "1706.06083"
}

@INPROCEEDINGS{Montasser2019-ro,
  title     = "{VC} Classes are Adversarially Robustly Learnable, but Only
               Improperly",
  booktitle = "Proceedings of the {Thirty-Second} Conference on Learning Theory",
  author    = "Montasser, Omar and Hanneke, Steve and Srebro, Nathan",
  editor    = "Beygelzimer, Alina and Hsu, Daniel",
  publisher = "PMLR",
  volume    =  99,
  pages     = "2512--2530",
  series    = "Proceedings of Machine Learning Research",
  year      =  2019,
  address   = "Phoenix, USA",
eprint={1902.04217},
      archivePrefix={arXiv},
      primaryClass={cs.LG},
}

@ARTICLE{Saha2019-ce,
  title    = "Hidden Trigger Backdoor Attacks",
  author   = "Saha, Aniruddha and Subramanya, Akshayvarun and Pirsiavash, Hamed",
  journal  = "AAAI",
  volume   =  34,
  number   =  07,
  pages    = "11957--11965",
  month    =  apr,
  year     =  2020,
  language = "en",
eprint={1910.00033},
      archivePrefix={arXiv},
      primaryClass={cs.CV},
}

@INPROCEEDINGS{Shen2018-jx,
  title     = "Learning with Bad Training Data via Iterative Trimmed Loss
               Minimization",
  booktitle = "Proceedings of the 36th International Conference on Machine
               Learning",
  author    = "Shen, Yanyao and Sanghavi, Sujay",
  editor    = "Chaudhuri, Kamalika and Salakhutdinov, Ruslan",
  publisher = "PMLR",
  volume    =  97,
  pages     = "5739--5748",
  series    = "Proceedings of Machine Learning Research",
  year      =  2019,
eprint={1810.11874},
      archivePrefix={arXiv},
      primaryClass={cs.LG},
}

@BOOK{Shalev-Shwartz2014-oj,
  title     = "Understanding Machine Learning: From Theory to Algorithms",
  author    = "Shalev-Shwartz, Shai and Ben-David, Shai",
  publisher = "Cambridge University Press",
  month     =  may,
  year      =  2014,
  language  = "en"
}

@INPROCEEDINGS{Tran2018-bf,
  title     = "Spectral Signatures in Backdoor Attacks",
  booktitle = "Advances in Neural Information Processing Systems",
  author    = "Tran, Brandon and Li, Jerry and Madry, Aleksander",
  editor    = "Bengio, S and Wallach, H and Larochelle, H and Grauman, K and
               Cesa-Bianchi, N and Garnett, R",
  publisher = "Curran Associates, Inc.",
  volume    =  31,
  year      =  2018,
eprint={1811.00636},
      archivePrefix={arXiv},
      primaryClass={cs.LG},
}

@InProceedings{sfkm24,
  title = 	 {Faster Convergence with MultiWay Preferences},
  author =       {Saha, Aadirupa and Feldman, Vitaly and Mansour, Yishay and Koren, Tomer},
  booktitle = 	 {Proceedings of The 27th International Conference on Artificial Intelligence and Statistics},
  pages = 	 {433--441},
  year = 	 {2024},
  editor = 	 {Dasgupta, Sanjoy and Mandt, Stephan and Li, Yingzhen},
  volume = 	 {238},
  series = 	 {Proceedings of Machine Learning Research},
  month = 	 {05},
  publisher =    {PMLR},
  url = 	 {https://proceedings.mlr.press/v238/saha24a.html},
eprint={2312.11788},
      archivePrefix={arXiv},
      primaryClass={cs.LG},
}

@article{dasgupta,
author = {Dasgupta, Sanjoy and Gupta, Anupam},
title = {An elementary proof of a theorem of Johnson and Lindenstrauss},
journal = {Random Structures \& Algorithms},
volume = {22},
number = {1},
pages = {60-65},
doi = {https://doi.org/10.1002/rsa.10073},
year = {2003}
}

@inbook{moitra_2021, place={Cambridge}, title={Semirandom Stochastic Block Models}, DOI={10.1017/9781108637435.014}, booktitle={Beyond the Worst-Case Analysis of Algorithms}, publisher={Cambridge University Press}, author={Moitra, Ankur}, editor={Roughgarden, TimEditor}, year={2021}, pages={212–233}}

@inbook{feige_2021, place={Cambridge}, title={Introduction to Semirandom Models}, DOI={10.1017/9781108637435.013}, booktitle={Beyond the Worst-Case Analysis of Algorithms}, publisher={Cambridge University Press}, author={Feige, Uriel}, editor={Roughgarden, TimEditor}, year={2021}, pages={189–211}}

@article{shivaswamy2011online,
  title={Online learning with preference feedback},
  author={Shivaswamy, Pannagadatta K and Joachims, Thorsten},
  journal={arXiv preprint arXiv:1111.0712},
  year={2011}
}

@inproceedings{saha2022efficient,
  title={Efficient and optimal algorithms for contextual dueling bandits under realizability},
  author={Saha, Aadirupa and Krishnamurthy, Akshay},
  booktitle={International Conference on Algorithmic Learning Theory},
  pages={968--994},
  year={2022},
  organization={PMLR}
}

@inproceedings{komiyama2015regret,
  title={Regret lower bound and optimal algorithm in dueling bandit problem},
  author={Komiyama, Junpei and Honda, Junya and Kashima, Hisashi and Nakagawa, Hiroshi},
  booktitle={Conference on learning theory},
  pages={1141--1154},
  year={2015},
  organization={PMLR}
}

@article{ouyang2022training,
  title={Training language models to follow instructions with human feedback},
  author={Ouyang, Long and Wu, Jeffrey and Jiang, Xu and Almeida, Diogo and Wainwright, Carroll and Mishkin, Pamela and Zhang, Chong and Agarwal, Sandhini and Slama, Katarina and Ray, Alex and others},
  journal={Advances in Neural Information Processing Systems},
  volume={35},
  pages={27730--27744},
  year={2022}
}

@inproceedings{chu2011contextual,
  title={Contextual bandits with linear payoff functions},
  author={Chu, Wei and Li, Lihong and Reyzin, Lev and Schapire, Robert},
  booktitle={Proceedings of the Fourteenth International Conference on Artificial Intelligence and Statistics},
  pages={208--214},
  year={2011},
  organization={JMLR Workshop and Conference Proceedings}
}

@article{bobadilla2013recommender,
  title={Recommender systems survey},
  author={Bobadilla, Jes{\'u}s and Ortega, Fernando and Hernando, Antonio and Guti{\'e}rrez, Abraham},
  journal={Knowledge-based systems},
  volume={46},
  pages={109--132},
  year={2013},
  publisher={Elsevier}
}

@inproceedings{ailon2014reducing,
  title={Reducing dueling bandits to cardinal bandits},
  author={Ailon, Nir and Karnin, Zohar and Joachims, Thorsten},
  booktitle={International Conference on Machine Learning},
  pages={856--864},
  year={2014},
  organization={PMLR}
}

@misc{bubeck2015convex,
      title={Convex Optimization: Algorithms and Complexity}, 
      author={Sébastien Bubeck},
      year={2015},
      eprint={1405.4980},
      archivePrefix={arXiv},
      primaryClass={math.OC}
}

@book{vershynin_2018, place={Cambridge}, series={Cambridge Series in Statistical and Probabilistic Mathematics}, title={High-Dimensional Probability: An Introduction with Applications in Data Science}, DOI={10.1017/9781108231596}, publisher={Cambridge University Press}, author={Vershynin, Roman}, year={2018}, collection={Cambridge Series in Statistical and Probabilistic Mathematics}}

@article{lobel2018multidimensional,
  title={Multidimensional binary search for contextual decision-making},
  author={Lobel, Ilan and Leme, Renato Paes and Vladu, Adrian},
  journal={Operations Research},
  volume={66},
  number={5},
  pages={1346--1361},
  year={2018},
  publisher={INFORMS}
}

@misc{lls20,
  doi = {10.48550/ARXIV.2003.01703},
  
  url = {https://arxiv.org/abs/2003.01703},
  
  author = {Liu, Allen and Leme, Renato Paes and Schneider, Jon},
  
  title = {Optimal Contextual Pricing and Extensions},
  
  publisher = {arXiv},
  
  year = {2020},
  
  copyright = {arXiv.org perpetual, non-exclusive license}
}

@misc{ls18,
  doi = {10.48550/ARXIV.1804.03195},
  
  url = {https://arxiv.org/abs/1804.03195},
  
  author = {Leme, Renato Paes and Schneider, Jon},
  
  title = {Contextual Search via Intrinsic Volumes},
  
  publisher = {arXiv},
  
  year = {2018},
  
  copyright = {arXiv.org perpetual, non-exclusive license}
}

@misc{li2019stochastic,
      title={Stochastic Linear Optimization with Adversarial Corruption}, 
      author={Yingkai Li and Edmund Y. Lou and Liren Shan},
      year={2019},
      eprint={1909.02109},
      archivePrefix={arXiv},
      primaryClass={cs.LG}
}

@misc{gollapudi2021contextual,
  doi = {10.48550/ARXIV.2106.04819},
  
  url = {https://arxiv.org/abs/2106.04819},
  
  author = {Gollapudi, Sreenivas and Guruganesh, Guru and Kollias, Kostas and Manurangsi, Pasin and Leme, Renato Paes and Schneider, Jon},
  
  title = {Contextual Recommendations and Low-Regret Cutting-Plane Algorithms},
  
  publisher = {arXiv},
  
  year = {2021},
  
  copyright = {Creative Commons Attribution 4.0 International}
}

@misc{bfl21,
  doi = {10.48550/ARXIV.2106.14015},
  
  url = {https://arxiv.org/abs/2106.14015},
  
  author = {Besbes, Omar and Fonseca, Yuri and Lobel, Ilan},
  
  title = {Contextual Inverse Optimization: Offline and Online Learning},
  
  publisher = {arXiv},
  
  year = {2021},
  
  copyright = {Creative Commons Attribution 4.0 International}
}

@inproceedings{bv06,
  title={Learning from revealed preference},
  author={Beigman, Eyal and Vohra, Rakesh},
  booktitle={Proceedings of the 7th ACM Conference on Electronic Commerce},
  pages={36--42},
  year={2006}
}

@inproceedings{skm21,
  title={Dueling convex optimization},
  author={Saha, Aadirupa and Koren, Tomer and Mansour, Yishay},
  booktitle={International Conference on Machine Learning},
  pages={9245--9254},
  year={2021},
  organization={PMLR}
}

@article{jamieson2012query,
  title={Query complexity of derivative-free optimization},
  author={Jamieson, Kevin G and Nowak, Robert and Recht, Ben},
  journal={Advances in Neural Information Processing Systems},
  volume={25},
  year={2012}
}

@article{saha2022dueling,
  title={Dueling Convex Optimization with General Preferences},
  author={Saha, Aadirupa and Koren, Tomer and Mansour, Yishay},
  journal={arXiv preprint arXiv:2210.02562},
  year={2022}
}

@article{yue2012k,
  title={The k-armed dueling bandits problem},
  author={Yue, Yisong and Broder, Josef and Kleinberg, Robert and Joachims, Thorsten},
  journal={Journal of Computer and System Sciences},
  volume={78},
  number={5},
  pages={1538--1556},
  year={2012},
  publisher={Elsevier}
}

@inproceedings{dudik2015contextual,
  title={Contextual dueling bandits},
  author={Dud{\'i}k, Miroslav and Hofmann, Katja and Schapire, Robert E and Slivkins, Aleksandrs and Zoghi, Masrour},
  booktitle={Conference on Learning Theory},
  pages={563--587},
  year={2015},
  organization={PMLR}
}

@article{davenport20141,
  title={1-bit matrix completion},
  author={Davenport, Mark A and Plan, Yaniv and Van Den Berg, Ewout and Wootters, Mary},
  journal={Information and Inference: A Journal of the IMA},
  volume={3},
  number={3},
  pages={189--223},
  year={2014},
  publisher={OUP}
}

@article{blum1995coloring,
  title={Coloring random and semi-random k-colorable graphs},
  author={Blum, Avrim and Spencer, Joel},
  journal={Journal of Algorithms},
  volume={19},
  number={2},
  pages={204--234},
  year={1995},
  publisher={Elsevier}
}

@article{kelner2022semi,
  title={Semi-random sparse recovery in nearly-linear time},
  author={Kelner, Jonathan A and Li, Jerry and Liu, Allen and Sidford, Aaron and Tian, Kevin},
  journal={arXiv preprint arXiv:2203.04002},
  year={2022}
}

@inproceedings{moitra2016robust,
  title={How robust are reconstruction thresholds for community detection?},
  author={Moitra, Ankur and Perry, William and Wein, Alexander S},
  booktitle={Proceedings of the forty-eighth annual ACM symposium on Theory of Computing},
  pages={828--841},
  year={2016}
}

@inproceedings{cheng2018non,
  title={Non-convex matrix completion against a semi-random adversary},
  author={Cheng, Yu and Ge, Rong},
  booktitle={Conference On Learning Theory},
  pages={1362--1394},
  year={2018},
  organization={PMLR}
}

@article{massart2006risk,
  title={Risk bounds for statistical learning},
  author={Massart, Pascal and N{\'e}d{\'e}lec, {\'E}lodie},
  year={2006}
}

@article{diakonikolas2019distribution,
  title={Distribution-independent pac learning of halfspaces with massart noise},
  author={Diakonikolas, Ilias and Gouleakis, Themis and Tzamos, Christos},
  journal={Advances in Neural Information Processing Systems},
  volume={32},
  year={2019}
}

@inproceedings{yao77,
  title={Probabilistic computations: Toward a unified measure of complexity},
  author={Yao, Andrew Chi-Chin},
  booktitle={18th Annual Symposium on Foundations of Computer Science (sfcs 1977)},
  pages={222--227},
  year={1977},
  organization={IEEE Computer Society}
}

@article{j16,
  title={Recommender systems—beyond matrix completion},
  author={Jannach, Dietmar and Resnick, Paul and Tuzhilin, Alexander and Zanker, Markus},
  journal={Communications of the ACM},
  volume={59},
  number={11},
  pages={94--102},
  year={2016},
  publisher={ACM New York, NY, USA}
}

@misc{ksj18,
      title={Global linear convergence of Newton's method without strong-convexity or Lipschitz gradients}, 
      author={Sai Praneeth Karimireddy and Sebastian U. Stich and Martin Jaggi},
      year={2018},
      eprint={1806.00413},
      archivePrefix={arXiv},
      primaryClass={cs.LG},
      url={https://arxiv.org/abs/1806.00413}, 
}

@InProceedings{syls16,
  title = 	 {AdaDelay: Delay Adaptive Distributed Stochastic Optimization},
  author = 	 {Sra, Suvrit and Yu, Adams Wei and Li, Mu and Smola, Alex},
  booktitle = 	 {Proceedings of the 19th International Conference on Artificial Intelligence and Statistics},
  pages = 	 {957--965},
  year = 	 {2016},
  editor = 	 {Gretton, Arthur and Robert, Christian C.},
  volume = 	 {51},
  series = 	 {Proceedings of Machine Learning Research},
  address = 	 {Cadiz, Spain},
  month = 	 {05},
  publisher =    {PMLR},
eprint={1508.05003},
      archivePrefix={arXiv},
      primaryClass={stat.ML},
}

@misc{ajk24,
      title={Acceleration Meets Inverse Maintenance: Faster $\ell_{\infty}$-Regression}, 
      author={Deeksha Adil and Shunhua Jiang and Rasmus Kyng},
      year={2024},
      eprint={2409.20030},
      archivePrefix={arXiv},
      primaryClass={cs.DS},
      url={https://arxiv.org/abs/2409.20030}, 
}

@inproceedings{bcll18,
author = {Bubeck, S\'{e}bastien and Cohen, Michael B. and Lee, Yin Tat and Li, Yuanzhi},
title = {An homotopy method for lp regression provably beyond self-concordance and in input-sparsity time},
year = {2018},
isbn = {9781450355599},
publisher = {Association for Computing Machinery},
address = {New York, NY, USA},
booktitle = {Proceedings of the 50th Annual ACM SIGACT Symposium on Theory of Computing},
pages = {1130–1137},
numpages = {8},
location = {Los Angeles, CA, USA},
series = {STOC 2018},
eprint={1711.01328},
      archivePrefix={arXiv},
      primaryClass={math.OC},
}

@incollection{john1948,
	title        = {Extremum problems with inequalities as subsidiary conditions},
	author       = {John, Fritz},
	year         = 1948,
	booktitle    = {Studies and Essays Presented to R. Courant on his 60th Birthday},
	publisher    = {Interscience Publishers, Inc},
	pages        = {187--204}
}

@book{nn94,
  title={Interior Point Polynomial Algorithms in Convex Programming},
  author={Nesterov, Y. and Nemirovskii, A.},
  isbn={9780898715156},
  lccn={93038912},
  series={SIAM studies in applied and numerical mathematics: Society for Industrial and Applied Mathematics},
  year={1994},
  publisher={Society for Industrial and Applied Mathematics}
}

@inproceedings{chjjs22,
author = {Carmon, Yair and Hausler, Danielle and Jambulapati, Arun and Jin, Yujia and Sidford, Aaron},
title = {Optimal and adaptive monteiro-svaiter acceleration},
year = {2022},
isbn = {9781713871088},
publisher = {Curran Associates Inc.},
address = {Red Hook, NY, USA},
booktitle = {Proceedings of the 36th International Conference on Neural Information Processing Systems},
articleno = {1479},
numpages = {13},
location = {New Orleans, LA, USA},
series = {NIPS '22},
eprint={2205.15371},
archivePrefix={arXiv},
primaryClass={math.OC},
}

@inbook{bjlls19,
author = {Bubeck, S\'{e}bastien and Jiang, Qijia and Lee, Yin Tat and Li, Yuanzhi and Sidford, Aaron},
title = {Complexity of highly parallel non-smooth convex optimization},
year = {2019},
publisher = {Curran Associates Inc.},
address = {Red Hook, NY, USA},
booktitle = {Proceedings of the 33rd International Conference on Neural Information Processing Systems},
articleno = {1247},
numpages = {10},
eprint={1906.10655},
archivePrefix={arXiv},
primaryClass={math.OC},
}

@inbook{akps19,
author = {Deeksha Adil and Rasmus Kyng and Richard Peng and Sushant Sachdeva},
title = {Iterative Refinement for {$\ell_p$}-norm Regression},
booktitle = {Proceedings of the 2019 Annual ACM-SIAM Symposium on Discrete Algorithms (SODA)},
chapter = {},
pages = {1405-1424},
year = {2019},
doi = {10.1137/1.9781611975482.86},
URL = {https://epubs.siam.org/doi/abs/10.1137/1.9781611975482.86},
eprint = {https://epubs.siam.org/doi/pdf/10.1137/1.9781611975482.86},
    eprint={1901.06764},
    archivePrefix={arXiv},
    primaryClass={cs.DS},
}

@article{lfn18,
  title={Relatively smooth convex optimization by first-order methods, and applications},
  author={Lu, Haihao and Freund, Robert M and Nesterov, Yurii},
  journal={SIAM Journal on Optimization},
  volume={28},
  number={1},
  pages={333--354},
  year={2018},
  publisher={SIAM},
eprint={1610.05708},
      archivePrefix={arXiv},
      primaryClass={math.OC},
}

@article{qsz02,
title = {A primal–dual algorithm for minimizing a sum of Euclidean norms},
journal = {Journal of Computational and Applied Mathematics},
volume = {138},
number = {1},
pages = {127-150},
year = {2002},
issn = {0377-0427},
doi = {https://doi.org/10.1016/S0377-0427(01)00357-0},
url = {https://www.sciencedirect.com/science/article/pii/S0377042701003570},
author = {Liqun Qi and Defeng Sun and Guanglu Zhou}
}

@inproceedings{bllssw21,
author = {van den Brand, Jan and Lee, Yin Tat and Liu, Yang P. and Saranurak, Thatchaphol and Sidford, Aaron and Song, Zhao and Wang, Di},
title = {Minimum cost flows, MDPs, and ${\ell_1}$-regression in nearly linear time for dense instances},
year = {2021},
isbn = {9781450380539},
publisher = {Association for Computing Machinery},
address = {New York, NY, USA},
url = {https://doi.org/10.1145/3406325.3451108},
doi = {10.1145/3406325.3451108},
booktitle = {Proceedings of the 53rd Annual ACM SIGACT Symposium on Theory of Computing},
pages = {859–869},
numpages = {11},
location = {Virtual, Italy},
series = {STOC 2021},
eprint={2101.05719},
      archivePrefix={arXiv},
      primaryClass={cs.DS},
}

@article{and96,
author = {Andersen, Knud D.},
title = {An Efficient Newton Barrier Method for Minimizing a Sum of Euclidean Norms},
journal = {SIAM Journal on Optimization},
volume = {6},
number = {1},
pages = {74-95},
year = {1996},
doi = {10.1137/0806006},
URL = {https://doi.org/10.1137/0806006},
eprint = {https://doi.org/10.1137/0806006}
}

@article{berk2021fairness,
  title={Fairness in criminal justice risk assessments: The state of the art},
  author={Berk, Richard and Heidari, Hoda and Jabbari, Shahin and Kearns, Michael and Roth, Aaron},
  journal={Sociological Methods \& Research},
  volume={50},
  number={1},
  pages={3--44},
  year={2021},
  publisher={Sage Publications Sage CA: Los Angeles, CA}
}

@inproceedings{selbst2019fairness,
  title={Fairness and abstraction in sociotechnical systems},
  author={Selbst, Andrew D and Boyd, Danah and Friedler, Sorelle A and Venkatasubramanian, Suresh and Vertesi, Janet},
  booktitle={Proceedings of the conference on fairness, accountability, and transparency},
  pages={59--68},
  year={2019}
}

@book{conn2000trust,
  title={Trust region methods},
  author={Conn, Andrew R and Gould, Nicholas IM and Toint, Philippe L},
  year={2000},
  publisher={SIAM}
}

@inproceedings{veale2018fairness,
  title={Fairness and accountability design needs for algorithmic support in high-stakes public sector decision-making},
  author={Veale, Michael and Van Kleek, Max and Binns, Reuben},
  booktitle={Proceedings of the 2018 chi conference on human factors in computing systems},
  pages={1--14},
  year={2018}
}

@article{data2014seizing,
  title={Seizing Opportunities, Preserving Values},
  author={Data, Big},
  journal={The White House Report Washington},
  year={2014}
}

@misc{chouldechova2016fair,
  title={Fair prediction with disparate impact: A study of bias in recidivism prediction instruments. Big Data, 5 (2), 153-163},
  author={Chouldechova, Alexandra},
  year={2016}
}

@inproceedings{agarwal2019fair,
  title={Fair regression: Quantitative definitions and reduction-based algorithms},
  author={Agarwal, Alekh and Dud{\'\i}k, Miroslav and Wu, Zhiwei Steven},
  booktitle={International Conference on Machine Learning},
  pages={120--129},
  year={2019},
  organization={PMLR}
}

@article{barocas2016big,
  title={Big data's disparate impact},
  author={Barocas, Solon and Selbst, Andrew D},
  journal={Calif. L. Rev.},
  volume={104},
  pages={671},
  year={2016},
  publisher={HeinOnline}
}

@article{corbett2023measure,
  title={The measure and mismeasure of fairness},
  author={Corbett-Davies, Sam and Gaebler, Johann D and Nilforoshan, Hamed and Shroff, Ravi and Goel, Sharad},
  journal={The Journal of Machine Learning Research},
  volume={24},
  number={1},
  pages={14730--14846},
  year={2023},
  publisher={JMLRORG}
}

@article{cohen2024fairness,
  title={Fairness of linear regression in decision making},
  author={Cohen-Addad, Vincent and Gavva, Surya Teja and Karthik, CS and Mathieu, Claire and Namrata},
  journal={International journal of data science and analytics},
  volume={18},
  number={3},
  pages={337--347},
  year={2024},
  publisher={Springer}
}

@article{duchi2016statistics,
  title={Statistics of robust optimization: A generalized empirical likelihood approach. arXiv},
  author={Duchi, J and Glynn, P and Namkoong, Hongseok},
  journal={Machine Learning},
  year={2016}
}

@article{chen2022fair,
  title={Fair assortment planning},
  author={Chen, Qinyi and Golrezaei, Negin and Susan, Fransisca and Baskoro, Edy},
  journal={arXiv preprint arXiv:2208.07341},
  year={2022}
}

@article{zhang2024stochastic,
  title={Stochastic approximation approaches to group distributionally robust optimization},
  author={Zhang, Lijun and Zhao, Peng and Zhuang, Zhen-Hua and Yang, Tianbao and Zhou, Zhi-Hua},
  journal={Advances in Neural Information Processing Systems},
  volume={36},
  year={2024}
}

@inproceedings{kleinberg2018algorithmic,
  title={Algorithmic fairness},
  author={Kleinberg, Jon and Ludwig, Jens and Mullainathan, Sendhil and Rambachan, Ashesh},
  booktitle={Aea papers and proceedings},
  volume={108},
  pages={22--27},
  year={2018}
}

@article{soma2022optimal,
  title={Optimal algorithms for group distributionally robust optimization and beyond},
  author={Soma, Tasuku and Gatmiry, Khashayar and Jegelka, Stefanie},
  journal={arXiv preprint arXiv:2212.13669},
  year={2022}
}

@article{golrezaei2024online,
  title={Online Combinatorial Optimization with Group Fairness Constraints},
  author={Golrezaei, Negin and Niazadeh, Rad and Patel, Kumar Kshitij and Susan, Fransisca},
  journal={Available at SSRN 4824251},
  year={2024}
}

@InProceedings{aakmrz22,
  title = 	 {Active Sampling for Min-Max Fairness},
  author =       {Abernethy, Jacob D and Awasthi, Pranjal and Kleindessner, Matth{\"a}us and Morgenstern, Jamie and Russell, Chris and Zhang, Jie},
  booktitle = 	 {Proceedings of the 39th International Conference on Machine Learning},
  pages = 	 {53--65},
  year = 	 {2022},
  editor = 	 {Chaudhuri, Kamalika and Jegelka, Stefanie and Song, Le and Szepesvari, Csaba and Niu, Gang and Sabato, Sivan},
  volume = 	 {162},
  series = 	 {Proceedings of Machine Learning Research},
  month = 	 {07},
  publisher =    {PMLR},
  url = 	 {https://proceedings.mlr.press/v162/abernethy22a.html}
}

@MISC {var_stackexchange,
    TITLE = {Bound the variance of the product of two random varables.},
    AUTHOR = {Davide Giraudo},
    HOWPUBLISHED = {Mathematics Stack Exchange},
    URL = {https://math.stackexchange.com/q/1044864},
    YEAR={2014},
    MONTH={11},
}

@article{ms13,
author = {Monteiro, Renato D. C. and Svaiter, B. F.},
title = {An Accelerated Hybrid Proximal Extragradient Method for Convex Optimization and Its Implications to Second-Order Methods},
journal = {SIAM Journal on Optimization},
volume = {23},
number = {2},
pages = {1092-1125},
year = {2013},
doi = {10.1137/110833786},
URL = { 
        https://doi.org/10.1137/110833786
},
eprint = { 
        https://doi.org/10.1137/110833786
}
}

@misc{
svwz24,
title={On Socially Fair Regression and Low-Rank Approximation},
author={Zhao Song and Ali Vakilian and David Woodruff and Samson Zhou},
year={2024},
url={https://openreview.net/forum?id=KJHUYWviZ6}
}

@misc{ob18,
      title={Finite-sample analysis of M-estimators using self-concordance}, 
      author={Dmitrii Ostrovskii and Francis Bach},
      year={2020},
      eprint={1810.06838},
      archivePrefix={arXiv},
      primaryClass={math.ST},
      url={https://arxiv.org/abs/1810.06838}, 
}

@inproceedings{msbacon,
author = {Carmon, Yair and Jambulapati, Arun and Jiang, Qijia and Jin, Yujia and Lee, Yin Tat and Sidford, Aaron and Tian, Kevin},
title = {Acceleration with a ball optimization oracle},
year = {2020},
isbn = {9781713829546},
publisher = {Curran Associates Inc.},
address = {Red Hook, NY, USA},
booktitle = {Proceedings of the 34th International Conference on Neural Information Processing Systems},
articleno = {1599},
numpages = {12},
series = {NIPS '20},
eprint={2003.08078},
      archivePrefix={arXiv},
        primaryClass={math.OC}
}

@inbook{mo23,
      title={The Change-of-Measure Method, Block Lewis Weights, and Approximating Matrix Block Norms}, 
      author={Naren Sarayu Manoj and Max Ovsiankin},
      year={2025},
      eprint={2311.10013},
      archivePrefix={arXiv},
      primaryClass={math.FA},
      booktitle = {Proceedings of the 2025 Annual ACM-SIAM Symposium on Discrete Algorithms (SODA)}
}

@inproceedings{cd21,
  title={Query complexity of least absolute deviation regression via robust uniform convergence},
  author={Chen, Xue and Dereziński, Michał},
  booktitle={Conference on Learning Theory},
  pages={1144--1179},
  year={2021},
  organization={PMLR},
  eprint={2102.02322},
  archivePrefix={arXiv},
  primaryClass={cs.LG}
}

@inproceedings{clmps16,
  title={Geometric median in nearly linear time},
  author={Cohen, Michael B and Lee, Yin Tat and Miller, Gary and Pachocki, Jakub and Sidford, Aaron},
  booktitle={Proceedings of the forty-eighth annual ACM symposium on Theory of Computing},
  pages={9--21},
  year={2016},
  eprint={1606.05225},
  archivePrefix={arXiv},
  primaryClass={cs.DS}
}

@article{acco00,
  title={An efficient primal-dual interior-point method for minimizing a sum of Euclidean norms},
  author={Andersen, Knud D and Christiansen, Edmund and Conn, Andrew R and Overton, Michael L},
  journal={SIAM Journal on Scientific Computing},
  volume={22},
  number={1},
  pages={243--262},
  year={2000},
  publisher={SIAM}
}

@article{xy97,
  title={An efficient algorithm for minimizing a sum of Euclidean norms with applications},
  author={Xue, Guoliang and Ye, Yinyu},
  journal={SIAM Journal on Optimization},
  volume={7},
  number={4},
  pages={1017--1036},
  year={1997},
  publisher={SIAM}
}

@article{rudelson1996,
  title={Random vectors in the isotropic position},
  author={Rudelson, Mark},
  journal={Journal of Functional Analysis},
  volume={164},
  number={1},
  pages={60--72},
  year={1999},
  publisher={Elsevier},
  eprint={math/9608208},
  archivePrefix={arXiv},
  primaryClass={math.MG}
}

@book{Wainwright_2019, place={Cambridge}, series={Cambridge Series in Statistical and Probabilistic Mathematics}, title={High-Dimensional Statistics: A Non-Asymptotic Viewpoint}, publisher={Cambridge University Press}, author={Wainwright, Martin J.}, year={2019}, collection={Cambridge Series in Statistical and Probabilistic Mathematics}}

@article{yy06,
  title={Model selection and estimation in regression with grouped variables},
  author={Yuan, Ming and Lin, Yi},
  journal={Journal of the Royal Statistical Society Series B: Statistical Methodology},
  volume={68},
  number={1},
  pages={49--67},
  year={2006},
  publisher={Oxford University Press}
}

@article{bach07,
  title={Consistency of the group lasso and multiple kernel learning.},
  author={Bach, Francis R},
  journal={Journal of Machine Learning Research},
  volume={9},
  number={6},
  year={2008},
  eprint={0707.3390},
  archivePrefix={arXiv},
  primaryClass={cs.LG}
}

@article{sra12,
  title={Fast projections onto mixed-norm balls with applications},
  author={Sra, Suvrit},
  journal={Data Mining and Knowledge Discovery},
  volume={25},
  pages={358--377},
  year={2012},
  publisher={Springer},
  eprint={1204.1437},
  archivePrefix={arXiv},
  primaryClass={stat.ML}
}

@article{sfht13,
 ISSN = {10618600},
 URL = {http://www.jstor.org/stable/43304828},
 author = {Noah Simon and Jerome Friedman and Trevor Hastie and Robert Tibshirani},
 journal = {Journal of Computational and Graphical Statistics},
 number = {2},
 pages = {231--245},
 publisher = {[American Statistical Association, Taylor & Francis, Ltd., Institute of Mathematical Statistics, Interface Foundation of America]},
 title = {A Sparse-Group Lasso},
 volume = {22},
 year = {2013}
}

@inproceedings{hrr22,
  title={The vector balancing constant for zonotopes},
  author={Bozzai, Rainie and Reis, Victor and Rothvoss, Thomas},
  booktitle={2023 IEEE 64th Annual Symposium on Foundations of Computer Science (FOCS)},
  pages={1292--1300},
  year={2023},
  organization={IEEE},
  eprint={2210.16460},
  archivePrefix={arXiv},
  primaryClass={math.MG}
}

@article{mu19,
  title={Entropy numbers of finite dimensional mixed-norm balls and function space embeddings with small mixed smoothness},
  author={Mayer, Sebastian and Ullrich, Tino},
  journal={Constructive Approximation},
  volume={53},
  pages={249--279},
  year={2021},
  publisher={Springer},
  eprint={1904.04619},
  archivePrefix={arXiv},
  primaryClass={math.FA}
}

@article{kv15,
  title={Volumes of unit balls of mixed sequence spaces},
  author={Kempka, Henning and Vyb\'{i}ral, Jan},
  journal={Mathematische Nachrichten},
  volume={290},
  number={8-9},
  pages={1317--1327},
  year={2017},
  publisher={Wiley Online Library},
  eprint={1505.05867},
  archivePrefix={arXiv},
  primaryClass={math.FA}
}

@misc{ps12,
      title={Combinatorial Inequalities and Subspaces of L1}, 
      author={Joscha Prochno and Carsten Schuett},
      year={2012},
      eprint={1204.6025},
      archivePrefix={arXiv},
      primaryClass={math.FA}
}

@misc{jkb22,
      title={Limit theorems for mixed-norm sequence spaces with applications to volume distribution}, 
      author={Michael Juhos and Zakhar Kabluchko and Joscha Prochno},
      year={2022},
      eprint={2209.08937},
      archivePrefix={arXiv},
      primaryClass={math.PR}
}

@article{ptj87,
 author = {Alain Pajor and Nicole Tomczak-Jaegermann},
 journal = {Proceedings of the American Mathematical Society},
 number = {4},
 pages = {637--642},
 publisher = {American Mathematical Society},
 title = {Subspaces of Small Codimension of Finite-Dimensional Banach Spaces},
 volume = {97},
 year = {1986}
}

@misc{kla23,
      title={Logarithmic bounds for isoperimetry and slices of convex sets}, 
      author={Bo'az Klartag},
      year={2023},
      eprint={2303.14938},
      archivePrefix={arXiv},
      primaryClass={math.FA}
}

@article{js01,
author = {Johnson, William and Schechtman, Gideon},
year = {2000},
month = {09},
pages = {},
title = {Finite dimensional subspaces of {$L_p$}},
volume = {1},
journal = {Handbook of the Geometry of Banach Spaces}
}

@article{Lewis1978,
author = {Lewis, D.},
journal = {Studia Mathematica},
language = {eng},
number = {2},
pages = {207-212},
title = {Finite dimensional subspaces of {$L_{p}$}},
url = {http://eudml.org/doc/218208},
volume = {63},
year = {1978},
}

@inproceedings{wy23_sens,
  title={Sharper Bounds for {$\ell_p$} Sensitivity Sampling},
  author={Woodruff, David and Yasuda, Taisuke},
  booktitle={International Conference on Machine Learning},
  pages={37238--37272},
  year={2023},
  organization={PMLR},
  eprint={2306.00732},
  archivePrefix={arXiv},
  primaryClass={cs.DS}
}

@inproceedings{jlls23,
  title={Sparsifying sums of norms},
  author={Jambulapati, Arun and Lee, James R and Liu, Yang P and Sidford, Aaron},
  booktitle={2023 IEEE 64th Annual Symposium on Foundations of Computer Science (FOCS)},
  pages={1953--1962},
  year={2023},
  organization={IEEE},
  eprint={2305.09049},
  archivePrefix={arXiv},
  primaryClass={cs.DS}
}

@inproceedings{ccly19,
  title={A near-optimal algorithm for approximating the john ellipsoid},
  author={Cohen, Michael B and Cousins, Ben and Lee, Yin Tat and Yang, Xin},
  booktitle={Conference on Learning Theory},
  pages={849--873},
  year={2019},
  organization={PMLR},
  eprint={1905.11580},
  archivePrefix={arXiv},
  primaryClass={cs.DS}
}

@article{blm89,
  title={Approximation of zonoids by zonotopes},
  author={Bourgain, Jean and Lindenstrauss, Joram and Milman, Vitali},
  year={1989}
}

@article{sz01,
  title={Embedding Subspaces of {$L_p$} into {$\ell_n^p$}, {$0< p< 1$}},
  author={Schechtman, Gideon and Zvavitch, Artem},
  journal={Mathematische Nachrichten},
  volume={227},
  number={1},
  pages={133--142},
  year={2001},
  publisher={Wiley Online Library}
}

@article{vh15,
  title={Chaining, interpolation and convexity II: The contraction principle},
  author={van Handel, Ramon},
  journal={The Annals of Probability},
  volume={46},
  number={3},
  pages={1764--1805},
  year={2018},
  publisher={JSTOR},
  eprint={1610.05199},
  archivePrefix={arXiv},
  primaryClass={math.PR}
}

@article{lww19,
  title={Tight bounds for the subspace sketch problem with applications},
  author={Li, Yi and Wang, Ruosong and Woodruff, David P},
  journal={SIAM Journal on Computing},
  volume={50},
  number={4},
  pages={1287--1335},
  year={2021},
  publisher={SIAM},
  eprint={1904.05543},
  archivePrefix={arXiv},
  primaryClass={cs.DS}
}

@inproceedings{ss08,
  title={Graph sparsification by effective resistances},
  author={Spielman, Daniel A and Srivastava, Nikhil},
  booktitle={Proceedings of the fortieth annual ACM symposium on Theory of computing},
  pages={563--568},
  year={2008},
  eprint={0803.0929},
  archivePrefix={arXiv},
  primaryClass={cs.DS}
}

@inproceedings{fl11,
  title={A unified framework for approximating and clustering data},
  author={Feldman, Dan and Langberg, Michael},
  booktitle={Proceedings of the forty-third annual ACM symposium on Theory of computing},
  pages={569--578},
  year={2011},
  eprint={1106.1379},
  archivePrefix={arXiv},
  primaryClass={cs.LG}
}

@book{rothvoss,
  title={Asymptotic Convex Geometry},
  author={Thomas Rothvoss},
  url={https://sites.math.washington.edu/~rothvoss/lecturenotes/AsymptoticConvexGeometry-24-FEB-2023.pdf},
  year={2023},
}

@book{tal21,
   title =     {Upper and Lower Bounds for Stochastic Processes: Decomposition Theorems},
   author =    {Michel Talagrand},
   publisher = {Springer},
   isbn =      {9783030825942},
   year =      {2021},
   series =    {Ergebnisse der Mathematik und ihrer Grenzgebiete. 3. Folge / A Series of Modern Surveys in Mathematics, 60; 60},
   edition =   {2nd ed. 2021},
   volume =    {},
   url =       {https://link.springer.com/book/10.1007/978-3-030-82595-9}
}

@inproceedings{cp15,
  title={{$\ell_p$} row sampling by lewis weights},
  author={Cohen, Michael B and Peng, Richard},
  booktitle={Proceedings of the forty-seventh annual ACM symposium on Theory of computing},
  pages={183--192},
  year={2015},
  eprint={1412.0588},
  archivePrefix={arXiv},
  primaryClass={cs.DS}
}

@inproceedings{mmwy21,
  title={Active Linear Regression for {$\ell_p$} Norms and Beyond}, 
  author={Musco, Cameron and Musco, Christopher and Woodruff, David P and Yasuda, Taisuke},
  booktitle={2022 IEEE 63rd Annual Symposium on Foundations of Computer Science (FOCS)},
  pages={744--753},
  year={2022},
  organization={IEEE},
  eprint={2111.04888},
  archivePrefix={arXiv},
  primaryClass={cs.LG}
}

@inproceedings{lee22,
  title={Spectral hypergraph sparsification via chaining},
  author={Lee, James R},
  booktitle={Proceedings of the 55th Annual ACM Symposium on Theory of Computing},
  pages={207--218},
  year={2023},
  eprint={2209.04539},
  archivePrefix={arXiv},
  primaryClass={math.PR}
}

@inproceedings{kkty21a,
  title={Spectral hypergraph sparsifiers of nearly linear size},
  author={Kapralov, Michael and Krauthgamer, Robert and Tardos, Jakab and Yoshida, Yuichi},
  booktitle={2021 IEEE 62nd Annual Symposium on Foundations of Computer Science (FOCS)},
  pages={1159--1170},
  year={2022},
  organization={IEEE},
  eprint={2106.02353},
  archivePrefix={arXiv},
  primaryClass={cs.DS}
}

@misc{ls19,
      title={Solving Linear Programs with Sqrt(rank) Linear System Solves}, 
      author={Yin Tat Lee and Aaron Sidford},
      year={2019},
      eprint={1910.08033},
      archivePrefix={arXiv},
      primaryClass={cs.DS}
}

@inproceedings{wy22_oneshot,
  title={Online lewis weight sampling},
  author={Woodruff, David P and Yasuda, Taisuke},
  booktitle={Proceedings of the 2023 Annual ACM-SIAM Symposium on Discrete Algorithms (SODA)},
  pages={4622--4666},
  year={2023},
  organization={SIAM},
  eprint={2207.08268},
  archivePrefix={arXiv},
  primaryClass={cs.DS}
}

@inproceedings{wy22,
  title={High-dimensional geometric streaming in polynomial space},
  author={Woodruff, David P and Yasuda, Taisuke},
  booktitle={2022 IEEE 63rd Annual Symposium on Foundations of Computer Science (FOCS)},
  pages={732--743},
  year={2022},
  organization={IEEE},
 eprint={2204.03790},
  archivePrefix={arXiv},
  primaryClass={cs.DS}
}

@inproceedings{jls21,
  title={Improved iteration complexities for overconstrained p-norm regression},
  author={Jambulapati, Arun and Liu, Yang P and Sidford, Aaron},
  booktitle={Proceedings of the 54th Annual ACM SIGACT Symposium on Theory of Computing},
  pages={529--542},
  year={2022},
  eprint={2111.01848},
  archivePrefix={arXiv},
  primaryClass={cs.DS}
}

@inproceedings{jls22,
  title={Chaining, group leverage score overestimates, and fast spectral hypergraph sparsification},
  author={Jambulapati, Arun and Liu, Yang P and Sidford, Aaron},
  booktitle={Proceedings of the 55th Annual ACM Symposium on Theory of Computing},
  pages={196--206},
  year={2023},
  eprint={2209.10539},
  archivePrefix={arXiv},
  primaryClass={cs.DS}
}

@article{nie2010efficient,
  title={Efficient and robust feature selection via joint {$\ell_{2, 1}$}-norms minimization},
  author={Nie, Feiping and Huang, Heng and Cai, Xiao and Ding, Chris},
  journal={Advances in neural information processing systems},
  volume={23},
  year={2010}
}

@article{chen2018tight,
  title={Tight bounds for collaborative pac learning via multiplicative weights},
  author={Chen, Jiecao and Zhang, Qin and Zhou, Yuan},
  journal={Advances in neural information processing systems},
  volume={31},
  year={2018},
eprint={1805.09217},
      archivePrefix={arXiv},
      primaryClass={cs.LG},
}

@article{nguyen2018improved,
  title={Improved algorithms for collaborative PAC learning},
  author={Nguyen, Huy and Zakynthinou, Lydia},
  journal={Advances in Neural Information Processing Systems},
  volume={31},
  year={2018},
eprint={1805.08356},
      archivePrefix={arXiv},
      primaryClass={cs.LG},
}

@inproceedings{zhang2024optimal,
  title={Optimal multi-distribution learning},
  author={Zhang, Zihan and Zhan, Wenhao and Chen, Yuxin and Du, Simon S and Lee, Jason D},
  booktitle={The Thirty Seventh Annual Conference on Learning Theory},
  pages={5220--5223},
  year={2024},
  organization={PMLR},
eprint={2312.05134},
      archivePrefix={arXiv},
      primaryClass={cs.LG},
}

@article{haghtalab2022demand,
  title={On-demand sampling: Learning optimally from multiple distributions},
  author={Haghtalab, Nika and Jordan, Michael and Zhao, Eric},
  journal={Advances in Neural Information Processing Systems},
  volume={35},
  pages={406--419},
  year={2022},
eprint={2210.12529},
      archivePrefix={arXiv},
      primaryClass={cs.LG},
}

@inproceedings{mohri2019agnostic,
  title={Agnostic federated learning},
  author={Mohri, Mehryar and Sivek, Gary and Suresh, Ananda Theertha},
  booktitle={International conference on machine learning},
  pages={4615--4625},
  year={2019},
  organization={PMLR},
eprint={1902.00146},
      archivePrefix={arXiv},
      primaryClass={cs.LG},
}

@article{hanneke2019value,
  title={On the value of target data in transfer learning},
  author={Hanneke, Steve and Kpotufe, Samory},
  journal={Advances in Neural Information Processing Systems},
  volume={32},
  year={2019},
eprint={2002.04747},
      archivePrefix={arXiv},
      primaryClass={cs.LG},
}

@article{bullins2021stochastic,
  title={A stochastic newton algorithm for distributed convex optimization},
  author={Bullins, Brian and Patel, Kshitij and Shamir, Ohad and Srebro, Nathan and Woodworth, Blake E},
  journal={Advances in Neural Information Processing Systems},
  volume={34},
  pages={26818--26830},
  year={2021},
eprint={2110.02954},
      archivePrefix={arXiv},
      primaryClass={math.OC},
}

@inproceedings{mmrha17,
  title={Communication-efficient learning of deep networks from decentralized data},
  author={McMahan, Brendan and Moore, Eider and Ramage, Daniel and Hampson, Seth and y Arcas, Blaise Aguera},
  booktitle={Artificial intelligence and statistics},
  pages={1273--1282},
  year={2017},
  organization={PMLR},
eprint={1602.05629},
      archivePrefix={arXiv},
      primaryClass={cs.LG},
}

@inproceedings{woodworth2020local,
  title={Is local SGD better than minibatch SGD?},
  author={Woodworth, Blake and Patel, Kumar Kshitij and Stich, Sebastian and Dai, Zhen and Bullins, Brian and Mcmahan, Brendan and Shamir, Ohad and Srebro, Nathan},
  booktitle={International Conference on Machine Learning},
  pages={10334--10343},
  year={2020},
  organization={PMLR},
eprint={2002.07839},
      archivePrefix={arXiv},
      primaryClass={cs.LG},
}

@InProceedings{patel2024limits,
  title = 	 {The Limits and Potentials of Local SGD for Distributed Heterogeneous Learning with Intermittent Communication},
  author =       {Patel, Kumar Kshitij and Glasgow, Margalit and Zindari, Ali and Wang, Lingxiao and Stich, Sebastian U and Cheng, Ziheng and Joshi, Nirmit and Srebro, Nathan},
  booktitle = 	 {Proceedings of Thirty Seventh Conference on Learning Theory},
  pages = 	 {4115--4157},
  year = 	 {2024},
  editor = 	 {Agrawal, Shipra and Roth, Aaron},
  volume = 	 {247},
  series = 	 {Proceedings of Machine Learning Research},
  month = 	 {07},
  publisher =    {PMLR},
  url = 	 {https://proceedings.mlr.press/v247/patel24a.html},
eprint={2405.11667},
      archivePrefix={arXiv},
      primaryClass={cs.LG},
}

@book{vapnik2013nature,
  title={The nature of statistical learning theory},
  author={Vapnik, Vladimir},
  year={2013},
  publisher={Springer science \& business media}
}

@article{bertsimas2011price,
  title={The price of fairness},
  author={Bertsimas, Dimitris and Farias, Vivek F and Trichakis, Nikolaos},
  journal={Operations research},
  volume={59},
  number={1},
  pages={17--31},
  year={2011},
  publisher={INFORMS}
}

@article{asadpour2022sequential,
  title={Sequential Submodular Maximization and Applications to Ranking an Assortment of Products},
  author={Asadpour, Arash and Niazadeh, Rad and Saberi, Amin and Shameli, Ali},
  journal={Operations Research},
  year={2022},
  publisher={INFORMS},
eprint={2002.09458},
      archivePrefix={arXiv},
      primaryClass={cs.GT},
}

@misc{chouldechova2018frontiers,
      title={The Frontiers of Fairness in Machine Learning}, 
      author={Alexandra Chouldechova and Aaron Roth},
      year={2018},
      eprint={1810.08810},
      archivePrefix={arXiv},
      primaryClass={cs.LG},
      url={https://arxiv.org/abs/1810.08810}, 
}

@article{ben2013robust,
  title={Robust solutions of optimization problems affected by uncertain probabilities},
  author={Ben-Tal, Aharon and Den Hertog, Dick and De Waegenaere, Anja and Melenberg, Bertrand and Rennen, Gijs},
  journal={Management Science},
  volume={59},
  number={2},
  pages={341--357},
  year={2013},
  publisher={INFORMS}
}

@article{blum2017collaborative,
  title={Collaborative PAC learning},
  author={Blum, Avrim and Haghtalab, Nika and Procaccia, Ariel D and Qiao, Mingda},
  journal={Advances in Neural Information Processing Systems},
  volume={30},
  year={2017}
}

@article{levy2020large,
  title={Large-scale methods for distributionally robust optimization},
  author={Levy, Daniel and Carmon, Yair and Duchi, John C and Sidford, Aaron},
  journal={Advances in Neural Information Processing Systems},
  volume={33},
  pages={8847--8860},
  year={2020},
eprint={2010.05893},
      archivePrefix={arXiv},
      primaryClass={math.OC},
}

@inproceedings{
sagawa2019distributionally,
title={Distributionally Robust Neural Networks},
author={Shiori Sagawa and Pang Wei Koh and Tatsunori B. Hashimoto and Percy Liang},
booktitle={International Conference on Learning Representations},
year={2020},
url={https://openreview.net/forum?id=ryxGuJrFvS},
eprint={1911.08731},
      archivePrefix={arXiv},
      primaryClass={cs.LG},
}

@article{rahmattalabi2019exploring,
  title={Exploring algorithmic fairness in robust graph covering problems},
  author={Rahmattalabi, Aida and Vayanos, Phebe and Fulginiti, Anthony and Rice, Eric and Wilder, Bryan and Yadav, Amulya and Tambe, Milind},
  journal={Advances in neural information processing systems},
  volume={32},
  year={2019},
eprint={2006.06865},
      archivePrefix={arXiv},
      primaryClass={math.OC},
}

@inproceedings{chierichetti2019matroids,
  title={Matroids, matchings, and fairness},
  author={Chierichetti, Flavio and Kumar, Ravi and Lattanzi, Silvio and Vassilvtiskii, Sergei},
  booktitle={The 22nd International Conference on Artificial Intelligence and Statistics},
  pages={2212--2220},
  year={2019},
  organization={PMLR}
}

@inproceedings{dwork2012fairness,
  title={Fairness through awareness},
  author={Dwork, Cynthia and Hardt, Moritz and Pitassi, Toniann and Reingold, Omer and Zemel, Richard},
  booktitle={Proceedings of the 3rd innovations in theoretical computer science conference},
  pages={214--226},
  year={2012}
}

@article{loi2019philosophical,
  title={A Philosophical Theory of Fairness for Prediction-Based Decisions.},
  author={Loi, Michele and Herlitz, Anders and Heidari, Hoda},
  journal={Available at SSRN 3450300},
  year={2019}
}

@article{hardt2016equality,
  title={Equality of opportunity in supervised learning},
  author={Hardt, Moritz and Price, Eric and Srebro, Nati},
  journal={Advances in neural information processing systems},
  volume={29},
  year={2016},
eprint={1610.02413},
      archivePrefix={arXiv},
      primaryClass={cs.LG},
}

@inproceedings{kearns2018preventing,
  title={Preventing fairness gerrymandering: Auditing and learning for subgroup fairness},
  author={Kearns, Michael and Neel, Seth and Roth, Aaron and Wu, Zhiwei Steven},
  booktitle={International conference on machine learning},
  pages={2564--2572},
  year={2018},
  organization={PMLR}
}

@inproceedings{kearns2019empirical,
  title={An empirical study of rich subgroup fairness for machine learning},
  author={Kearns, Michael and Neel, Seth and Roth, Aaron and Wu, Zhiwei Steven},
  booktitle={Proceedings of the conference on fairness, accountability, and transparency},
  pages={100--109},
  year={2019}
}

@inproceedings{singh2018fairness,
  title={Fairness of exposure in rankings},
  author={Singh, Ashudeep and Joachims, Thorsten},
  booktitle={Proceedings of the 24th ACM SIGKDD International Conference on Knowledge Discovery \& Data Mining},
  pages={2219--2228},
  year={2018}
}

@article{singh2019policy,
  title={Policy learning for fairness in ranking},
  author={Singh, Ashudeep and Joachims, Thorsten},
  journal={Advances in neural information processing systems},
  volume={32},
  year={2019}
}

@inproceedings{biega2018equity,
  title={Equity of attention: Amortizing individual fairness in rankings},
  author={Biega, Asia J and Gummadi, Krishna P and Weikum, Gerhard},
  booktitle={The 41st international acm sigir conference on research \& development in information retrieval},
  pages={405--414},
  year={2018}
}

@article{bertsimas2012efficiency,
  title={On the efficiency-fairness trade-off},
  author={Bertsimas, Dimitris and Farias, Vivek F and Trichakis, Nikolaos},
  journal={Management Science},
  volume={58},
  number={12},
  pages={2234--2250},
  year={2012},
  publisher={INFORMS}
}

@inproceedings{manshadi2021fair,
  title={Fair dynamic rationing},
  author={Manshadi, Vahideh and Niazadeh, Rad and Rodilitz, Scott},
  booktitle={Proceedings of the 22nd ACM Conference on Economics and Computation},
  pages={694--695},
  year={2021}
}

@inproceedings{balseiro2021regularized,
  title={Regularized online allocation problems: Fairness and beyond},
  author={Balseiro, Santiago and Lu, Haihao and Mirrokni, Vahab},
  booktitle={International Conference on Machine Learning},
  pages={630--639},
  year={2021},
  organization={PMLR}
}

@book{miettinen1999nonlinear,
  title={Nonlinear multiobjective optimization},
  author={Miettinen, Kaisa},
  volume={12},
  year={1999},
  publisher={Springer Science \& Business Media}
}

@article{hooker2012combining,
  title={Combining equity and utilitarianism in a mathematical programming model},
  author={Hooker, John N and Williams, H Paul},
  journal={Management Science},
  volume={58},
  number={9},
  pages={1682--1693},
  year={2012},
  publisher={INFORMS}
}

@article{mulvany2021fair,
  title={Fair scheduling of heterogeneous customer populations},
  author={Mulvany, Justin and Randhawa, Ramandeep S},
  journal={Available at SSRN 3803016},
  year={2021}
}

@article{kusner2017counterfactual,
  title={Counterfactual fairness},
  author={Kusner, Matt J and Loftus, Joshua and Russell, Chris and Silva, Ricardo},
  journal={Advances in neural information processing systems},
  volume={30},
  year={2017}
}

@inproceedings{ustun2019fairness,
  title={Fairness without harm: Decoupled classifiers with preference guarantees},
  author={Ustun, Berk and Liu, Yang and Parkes, David},
  booktitle={International Conference on Machine Learning},
  pages={6373--6382},
  year={2019},
  organization={PMLR}
}

@inproceedings{donahue2020fairness,
  title={Fairness and utilization in allocating resources with uncertain demand},
  author={Donahue, Kate and Kleinberg, Jon},
  booktitle={Proceedings of the 2020 conference on fairness, accountability, and transparency},
  pages={658--668},
  year={2020}
}

@book{ehrgott2005multicriteria,
  title={Multicriteria optimization},
  author={Ehrgott, Matthias},
  volume={491},
  year={2005},
  publisher={Springer Science \& Business Media}
}

@inproceedings{goel2018non,
  title={Non-discriminatory machine learning through convex fairness criteria},
  author={Goel, Naman and Yaghini, Mohammad and Faltings, Boi},
  booktitle={Proceedings of the 2018 AAAI/ACM Conference on AI, Ethics, and Society},
  pages={116--116},
  year={2018}
}

@article{gupta2022socially,
  title={Socially fair and hierarchical facility location problems},
  author={Gupta, Swati and Moondra, Jai and Singh, Mohit},
  journal={arXiv preprint arXiv:2211.14873},
  year={2022}
}

@article{ma2023fairness,
  title={Fairness maximization among offline agents in online-matching markets},
  author={Ma, Will and Xu, Pan and Xu, Yifan},
  journal={ACM Transactions on Economics and Computation},
  volume={10},
  number={4},
  pages={1--27},
  year={2023},
  publisher={ACM New York, NY}
}

@inproceedings{calders2009building,
  title={Building classifiers with independency constraints},
  author={Calders, Toon and Kamiran, Faisal and Pechenizkiy, Mykola},
  booktitle={2009 IEEE international conference on data mining workshops},
  pages={13--18},
  year={2009},
  organization={IEEE}
}

@inproceedings{kasy2021fairness,
  title={Fairness, equality, and power in algorithmic decision-making},
  author={Kasy, Maximilian and Abebe, Rediet},
  booktitle={Proceedings of the 2021 ACM Conference on Fairness, Accountability, and Transparency},
  pages={576--586},
  year={2021}
}

@inproceedings{abebe2020roles,
  title={Roles for computing in social change},
  author={Abebe, Rediet and Barocas, Solon and Kleinberg, Jon and Levy, Karen and Raghavan, Manish and Robinson, David G},
  booktitle={Proceedings of the 2020 conference on fairness, accountability, and transparency},
  pages={252--260},
  year={2020}
}

@inproceedings{
cdm24,
title={A Near-Linear Time Approximation Algorithm for Beyond-Worst-Case Graph Clustering},
author={Vincent Cohen-Addad and Tommaso d'Orsi and Aida Mousavifar},
booktitle={Forty-first International Conference on Machine Learning},
year={2024},
url={https://openreview.net/forum?id=MSFxOMM0gK},
eprint={2406.04857},
      archivePrefix={arXiv},
      primaryClass={cs.DS},
}

@article{hll83,
title = {Stochastic blockmodels: First steps},
journal = {Social Networks},
volume = {5},
number = {2},
pages = {109-137},
year = {1983},
issn = {0378-8733},
doi = {https://doi.org/10.1016/0378-8733(83)90021-7},
url = {https://www.sciencedirect.com/science/article/pii/0378873383900217},
author = {Paul W. Holland and Kathryn Blackmond Laskey and Samuel Leinhardt}
}

@Article{         numpy,
 title         = {Array programming with {NumPy}},
 author        = {Charles R. Harris and K. Jarrod Millman and St{\'{e}}fan J.
                 van der Walt and Ralf Gommers and Pauli Virtanen and David
                 Cournapeau and Eric Wieser and Julian Taylor and Sebastian
                 Berg and Nathaniel J. Smith and Robert Kern and Matti Picus
                 and Stephan Hoyer and Marten H. van Kerkwijk and Matthew
                 Brett and Allan Haldane and Jaime Fern{\'{a}}ndez del
                 R{\'{i}}o and Mark Wiebe and Pearu Peterson and Pierre
                 G{\'{e}}rard-Marchant and Kevin Sheppard and Tyler Reddy and
                 Warren Weckesser and Hameer Abbasi and Christoph Gohlke and
                 Travis E. Oliphant},
 year          = {2020},
 month         = sep,
 journal       = {Nature},
 volume        = {585},
 number        = {7825},
 pages         = {357--362},
 doi           = {10.1038/s41586-020-2649-2},
 publisher     = {Springer Science and Business Media {LLC}},
 url           = {https://doi.org/10.1038/s41586-020-2649-2}
}

@techreport{networkx,
  title={Exploring network structure, dynamics, and function using NetworkX},
  author={Hagberg, Aric and Swart, Pieter and S Chult, Daniel},
  year={2008},
  institution={Los Alamos National Lab.(LANL), Los Alamos, NM (United States)}
}

@ARTICLE{scipy,
  author  = {Virtanen, Pauli and Gommers, Ralf and Oliphant, Travis E. and
            Haberland, Matt and Reddy, Tyler and Cournapeau, David and
            Burovski, Evgeni and Peterson, Pearu and Weckesser, Warren and
            Bright, Jonathan and {van der Walt}, St{\'e}fan J. and
            Brett, Matthew and Wilson, Joshua and Millman, K. Jarrod and
            Mayorov, Nikolay and Nelson, Andrew R. J. and Jones, Eric and
            Kern, Robert and Larson, Eric and Carey, C J and
            Polat, {\.I}lhan and Feng, Yu and Moore, Eric W. and
            {VanderPlas}, Jake and Laxalde, Denis and Perktold, Josef and
            Cimrman, Robert and Henriksen, Ian and Quintero, E. A. and
            Harris, Charles R. and Archibald, Anne M. and
            Ribeiro, Ant{\^o}nio H. and Pedregosa, Fabian and
            {van Mulbregt}, Paul and {SciPy 1.0 Contributors}},
  title   = {{{SciPy} 1.0: Fundamental Algorithms for Scientific
            Computing in Python}},
  journal = {Nature Methods},
  year    = {2020},
  volume  = {17},
  pages   = {261--272},
  adsurl  = {https://rdcu.be/b08Wh},
  doi     = {10.1038/s41592-019-0686-2},
}

@INPROCEEDINGS{boppana87,
  author={Boppana, Ravi B.},
  booktitle={28th Annual Symposium on Foundations of Computer Science (sfcs 1987)}, 
  title={Eigenvalues and graph bisection: An average-case analysis}, 
  year={1987},
  volume={},
  number={},
  pages={280-285},
  doi={10.1109/SFCS.1987.22}}

@InProceedings{sl17,
  author =	{Schuster, Martin R. and Liskiewicz, Maciej},
  title =	{{New Abilities and Limitations of Spectral Graph Bisection}},
  booktitle =	{25th Annual European Symposium on Algorithms (ESA 2017)},
  pages =	{66:1--66:15},
  series =	{Leibniz International Proceedings in Informatics (LIPIcs)},
  ISBN =	{978-3-95977-049-1},
  ISSN =	{1868-8969},
  year =	{2017},
  volume =	{87},
  editor =	{Pruhs, Kirk and Sohler, Christian},
  publisher =	{Schloss Dagstuhl -- Leibniz-Zentrum f{\"u}r Informatik},
  address =	{Dagstuhl, Germany},
  URL =		{https://drops.dagstuhl.de/entities/document/10.4230/LIPIcs.ESA.2017.66},
  URN =		{urn:nbn:de:0030-drops-78658},
  doi =		{10.4230/LIPIcs.ESA.2017.66}
}

@inbook{bv24,
author = {Abhinav Bhardwaj and Van Vu},
title = {Matrix Perturbation: Davis-Kahan in the Infinity Norm},
booktitle = {Proceedings of the 2024 Annual ACM-SIAM Symposium on Discrete Algorithms (SODA)},
chapter = {},
pages = {880-934},
year = {2024},
month = {01},
eprint={2304.00328},
      archivePrefix={arXiv},
      primaryClass={math.PR},
}

@ARTICLE{abh16,
  author={Abbe, Emmanuel and Bandeira, Afonso S. and Hall, Georgina},
  journal={IEEE Transactions on Information Theory}, 
  title={Exact Recovery in the Stochastic Block Model}, 
  year={2016},
  volume={62},
  number={1},
  pages={471-487},
  doi={10.1109/TIT.2015.2490670},
eprint={1405.3267},
      archivePrefix={arXiv},
      primaryClass={cs.SI},
}

@InProceedings{kllst23,
  title = 	 {Semi-Random Sparse Recovery in Nearly-Linear Time},
  author =       {Kelner, Jonathan and Li, Jerry and Liu, Allen X. and Sidford, Aaron and Tian, Kevin},
  booktitle = 	 {Proceedings of Thirty Sixth Conference on Learning Theory},
  pages = 	 {2352--2398},
  year = 	 {2023},
  editor = 	 {Neu, Gergely and Rosasco, Lorenzo},
  volume = 	 {195},
  series = 	 {Proceedings of Machine Learning Research},
  month = 	 {07},
  publisher =    {PMLR},
  url = 	 {https://proceedings.mlr.press/v195/kelner23a.html},
eprint={2203.04002},
      archivePrefix={arXiv},
      primaryClass={cs.DS},
}

@article{mn07,
   title={Risk bounds for statistical learning},
   volume={34},
   ISSN={0090-5364},
   url={http://dx.doi.org/10.1214/009053606000000786},
   DOI={10.1214/009053606000000786},
   number={5},
   journal={The Annals of Statistics},
   publisher={Institute of Mathematical Statistics},
   author={Massart, Pascal and Nédélec, Élodie},
   year={2006},
   month={10},
}

@InProceedings{ycom24,
  title = 	 {Top-$K$ ranking with a monotone adversary},
  author =       {Yang, Yuepeng and Chen, Antares and Orecchia, Lorenzo and Ma, Cong},
  booktitle = 	 {Proceedings of Thirty Seventh Conference on Learning Theory},
  pages = 	 {5123--5162},
  year = 	 {2024},
  editor = 	 {Agrawal, Shipra and Roth, Aaron},
  volume = 	 {247},
  series = 	 {Proceedings of Machine Learning Research},
  month = 	 {07},
  publisher =    {PMLR},
  url = 	 {https://proceedings.mlr.press/v247/yang24b.html},
    eprint={2402.07445},
      archivePrefix={arXiv},
      primaryClass={stat.ML}
}

@misc{ls24,
      title={Spectral clustering in the Gaussian mixture block model}, 
      author={Shuangping Li and Tselil Schramm},
      year={2024},
      eprint={2305.00979},
      archivePrefix={arXiv},
      primaryClass={stat.ML}
}

@inbook{gnw24,
author = {Julia Gaudio and Xiaochun Niu and Ermin Wei},
title = {Exact Community Recovery in the Geometric SBM},
booktitle = {Proceedings of the 2024 Annual ACM-SIAM Symposium on Discrete Algorithms (SODA)},
chapter = {},
pages = {2158-2184},
doi = {10.1137/1.9781611977912.78},
URL = {https://epubs.siam.org/doi/abs/10.1137/1.9781611977912.78},
eprint = {https://epubs.siam.org/doi/pdf/10.1137/1.9781611977912.78},
year={2024},
eprint={2307.11196},
      archivePrefix={arXiv},
      primaryClass={cs.SI},
}

@article{abbe_survey,
  author  = {Emmanuel Abbe},
  title   = {Community Detection and Stochastic Block Models: Recent Developments},
  journal = {Journal of Machine Learning Research},
  year    = {2018},
  volume  = {18},
  number  = {177},
  pages   = {1--86},
  url     = {http://jmlr.org/papers/v18/16-480.html},
eprint={1703.10146},
      archivePrefix={arXiv},
      primaryClass={math.PR},
}

@inbook{moitra_sbm_2021, place={Cambridge}, title={Semirandom Stochastic Block Models}, DOI={10.1017/9781108637435.014}, booktitle={Beyond the Worst-Case Analysis of Algorithms}, publisher={Cambridge University Press}, author={Moitra, Ankur}, editor={Roughgarden, Tim}, year={2021}, pages={212–233}}

@article{gv15,
  title={Community detection in sparse networks via Grothendieck’s inequality},
  author={Gu{\'e}don, Olivier and Vershynin, Roman},
  journal={Probability Theory and Related Fields},
  volume={165},
  number={3},
  pages={1025--1049},
  year={2016},
  publisher={Springer},
eprint={1411.4686},
      archivePrefix={arXiv},
      primaryClass={math.ST}
}

@inproceedings{mpw16,
author = {Moitra, Ankur and Perry, William and Wein, Alexander S.},
title = {How robust are reconstruction thresholds for community detection?},
year = {2016},
isbn = {9781450341325},
publisher = {Association for Computing Machinery},
address = {New York, NY, USA},
url = {https://doi.org/10.1145/2897518.2897573},
doi = {10.1145/2897518.2897573},
booktitle = {Proceedings of the Forty-Eighth Annual ACM Symposium on Theory of Computing},
pages = {828–841},
numpages = {14},
location = {Cambridge, MA, USA},
series = {STOC '16},
eprint={1511.01473},
      archivePrefix={arXiv},
      primaryClass={cs.DS},
}

@article{vl07,
  title={A tutorial on spectral clustering},
  author={Von Luxburg, Ulrike},
  journal={Statistics and computing},
  volume={17},
  pages={395--416},
  year={2007},
  publisher={Springer},
eprint={0711.0189},
      archivePrefix={arXiv},
      primaryClass={cs.DS},
}

@article{sb15,
   title={Role of normalization in spectral clustering for stochastic blockmodels},
   volume={43},
   ISSN={0090-5364},
   url={http://dx.doi.org/10.1214/14-AOS1285},
   DOI={10.1214/14-aos1285},
   number={3},
   journal={The Annals of Statistics},
   publisher={Institute of Mathematical Statistics},
   author={Sarkar, Purnamrita and Bickel, Peter J.},
   year={2015},
   month={06},
eprint={1310.1495},
      archivePrefix={arXiv},
      primaryClass={stat.ML},
}

@article{afwz17,
  title={Entrywise eigenvector analysis of random matrices with low expected rank},
  author={Abbe, Emmanuel and Fan, Jianqing and Wang, Kaizheng and Zhong, Yiqiao},
  journal={Annals of statistics},
  volume={48},
  number={3},
  pages={1452},
  year={2020},
  publisher={NIH Public Access},
eprint={1709.09565},
      archivePrefix={arXiv},
      primaryClass={math.ST}
}

@article{fk01,
author = {Feige, Uriel and Kilian, Joe},
title = {Heuristics for Semirandom Graph Problems},
year = {2001},
issue_date = {December 2001},
publisher = {Academic Press, Inc.},
address = {USA},
volume = {63},
number = {4},
issn = {0022-0000},
url = {https://doi.org/10.1006/jcss.2001.1773},
doi = {10.1006/jcss.2001.1773},
journal = {J. Comput. Syst. Sci.},
month = {12},
pages = {639–671},
numpages = {33}
}

@inproceedings{mmv12,
author = {Makarychev, Konstantin and Makarychev, Yury and Vijayaraghavan, Aravindan},
title = {Approximation algorithms for semi-random partitioning problems},
year = {2012},
isbn = {9781450312455},
publisher = {Association for Computing Machinery},
address = {New York, NY, USA},
booktitle = {Proceedings of the Forty-Fourth Annual ACM Symposium on Theory of Computing},
pages = {367–384},
numpages = {18},
location = {New York, New York, USA},
series = {STOC '12},
eprint={1205.2234},
archivePrefix={arXiv},
primaryClass={cs.DS}
}

@article{llv16,
  title={Concentration and regularization of random graphs},
  author={Le, Can M and Levina, Elizaveta and Vershynin, Roman},
  journal={Random Structures \& Algorithms},
  volume={51},
  number={3},
  pages={538--561},
  year={2017},
  publisher={Wiley Online Library},
eprint={1506.00669},
      archivePrefix={arXiv},
      primaryClass={math.PR}
}

@article{dls20,
  title={Strong consistency, graph laplacians, and the stochastic block model},
  author={Deng, Shaofeng and Ling, Shuyang and Strohmer, Thomas},
  journal={Journal of Machine Learning Research},
  volume={22},
  number={117},
  pages={1--44},
  year={2021},
eprint={2004.09780},
      archivePrefix={arXiv},
      primaryClass={stat.ML}
}

@inproceedings{mmt20,
author = {Theo McKenzie and Hermish Mehta and Luca Trevisan},
title = {A new algorithm for the robust semi-random independent set problem},
year = {2020},
publisher = {Society for Industrial and Applied Mathematics},
address = {USA},
booktitle = {Proceedings of the Thirty-First Annual ACM-SIAM Symposium on Discrete Algorithms},
pages = {738–746},
numpages = {9},
location = {Salt Lake City, Utah},
series = {SODA '20},
eprint={1808.03633},
      archivePrefix={arXiv},
      primaryClass={cs.DS},
}

@inproceedings{csv17,
author = {Charikar, Moses and Steinhardt, Jacob and Valiant, Gregory},
title = {Learning from untrusted data},
year = {2017},
isbn = {9781450345286},
publisher = {Association for Computing Machinery},
address = {New York, NY, USA},
booktitle = {Proceedings of the 49th Annual ACM SIGACT Symposium on Theory of Computing},
pages = {47–60},
numpages = {14},
location = {Montreal, Canada},
series = {STOC 2017},
eprint={1611.02315},
      archivePrefix={arXiv},
      primaryClass={cs.LG},
}

@article{BS95,
title = {Coloring Random and Semi-Random k-Colorable Graphs},
journal = {Journal of Algorithms},
volume = {19},
number = {2},
pages = {204-234},
year = {1995},
issn = {0196-6774},
doi = {https://doi.org/10.1006/jagm.1995.1034},
url = {https://www.sciencedirect.com/science/article/pii/S0196677485710346},
author = {Avrim Blum and Joel Spencer}
}

@INPROCEEDINGS{BS17,
  author={Sankararaman, Abishek and Baccelli, François},
  booktitle={2017 55th Annual Allerton Conference on Communication, Control, and Computing (Allerton)}, 
  title={Community detection on euclidean random graphs}, 
  year={2017},
  volume={},
  number={},
  pages={510-517},
  doi={10.1109/ALLERTON.2017.8262780}
}

@article{ABARS20,
author = {Abbe, Emmanuel and Boix-Adser\`{a}, Enric and Ralli, Peter and Sandon, Colin},
title = {Graph Powering and Spectral Robustness},
journal = {SIAM Journal on Mathematics of Data Science},
volume = {2},
number = {1},
pages = {132-157},
year = {2020},
doi = {10.1137/19M1257135},

URL = {
        https://doi.org/10.1137/19M1257135
},
eprint = { 
        https://doi.org/10.1137/19M1257135
},
    archivePrefix={arXiv},
    eprint={1809.04818},
    primaryClass={cs.DS}
}

@INPROCEEDINGS{GMPS18,
  author={Galhotra, Sainyam and Pal, Soumyabrata and Mazumdar, Arya and Saha, Barna},
  booktitle={2018 56th Annual Allerton Conference on Communication, Control, and Computing (Allerton)}, 
  title={The Geometric Block Model and Applications}, 
  year={2018},
  volume={},
  number={},
  pages={1147-1150},
  doi={10.1109/ALLERTON.2018.8635938}}

@InProceedings{GMPS20,
  author =	{Galhotra, Sainyam and Mazumdar, Arya and Pal, Soumyabrata and Saha, Barna},
  title =	{{Connectivity of Random Annulus Graphs and the Geometric Block Model}},
  booktitle =	{Approximation, Randomization, and Combinatorial Optimization. Algorithms and Techniques (APPROX/RANDOM 2019)},
  pages =	{53:1--53:23},
  series =	{Leibniz International Proceedings in Informatics (LIPIcs)},
  ISBN =	{978-3-95977-125-2},
  ISSN =	{1868-8969},
  year =	{2019},
  volume =	{145},
  editor =	{Achlioptas, Dimitris and V\'{e}gh, L\'{a}szl\'{o} A.},
  publisher =	{Schloss Dagstuhl -- Leibniz-Zentrum f{\"u}r Informatik},
  address =	{Dagstuhl, Germany},
eprint={1804.05013},
      archivePrefix={arXiv},
      primaryClass={cs.DM}
}

@article{ABD21,
	author = {Avrachenkov, Konstantin and Bobu, Andrei and Dreveton, Maximilien},
        year = {2021},
        month = {03},
	doi = {10.1007/s00041-021-09825-2},
	id = {Avrachenkov2021},
	journal = {Journal of Fourier Analysis and Applications},
	number = {2},
	pages = {22},
	title = {Higher-Order Spectral Clustering for Geometric Graphs},
	url = {https://doi.org/10.1007/s00041-021-09825-2},
	volume = {27},
	year = {2021},
archivePrefix={arXiv},
      primaryClass={cs.LG},
      eprint={2009.11353},
}

@article{BHK16,
  title={A nearly tight sum-of-squares lower bound for the planted clique problem},
  author={Barak, Boaz and Hopkins, Samuel and Kelner, Jonathan and Kothari, Pravesh K and Moitra, Ankur and Potechin, Aaron},
  journal={SIAM Journal on Computing},
  volume={48},
  number={2},
  pages={687--735},
  year={2019},
  publisher={SIAM},
eprint={1604.03084},
      archivePrefix={arXiv},
      primaryClass={cs.CC},
}

@article{Jerrum92,
  author       = {Mark Jerrum},
  title        = {Large Cliques Elude the Metropolis Process},
  journal      = {Random Struct. Algorithms},
  volume       = {3},
  number       = {4},
  pages        = {347--360},
  year         = {1992},
  url          = {https://doi.org/10.1002/rsa.3240030402},
  doi          = {10.1002/RSA.3240030402},
  bibsource    = {dblp computer science bibliography, https://dblp.org}
}

@article{Kuc95,
	author = {Lud{\v e}k Ku{\v c}era},
	doi = {https://doi.org/10.1016/0166-218X(94)00103-K},
	issn = {0166-218X},
	journal = {Discrete Applied Mathematics},
	note = {Combinatorial optimization 1992},
	number = {2},
	pages = {193-212},
	title = {Expected complexity of graph partitioning problems},
	url = {https://www.sciencedirect.com/science/article/pii/0166218X9400103K},
	volume = {57},
	year = {1995}
}

@inproceedings{BKS23,
author = {Buhai, Rares-Darius and Kothari, Pravesh K. and Steurer, David},
title = {Algorithms Approaching the Threshold for Semi-random Planted Clique},
year = {2023},
isbn = {9781450399135},
publisher = {Association for Computing Machinery},
address = {New York, NY, USA},
booktitle = {Proceedings of the 55th Annual ACM Symposium on Theory of Computing},
pages = {1918–1926},
numpages = {9},
location = {Orlando, FL, USA},
series = {STOC 2023},
eprint={2212.05619},
      archivePrefix={arXiv},
      primaryClass={cs.DS}
}

@InProceedings{bglmsy23,
      title={Dueling Optimization with a Monotone Adversary}, 
      author={Avrim Blum and Meghal Gupta and Gene Li and Naren Sarayu Manoj and Aadirupa Saha and Yuanyuan Yang},
      year={2024},
      month={02},
      booktitle = {Proceedings of Thirty Fifth Conference on Algorithmic Learning Theory (ALT)},
      archivePrefix={arXiv},
      primaryClass={cs.DS},
      eprint={2311.11185},
}

@InProceedings{AV18,
  title = 	 {Clustering Semi-Random Mixtures of {G}aussians},
  author =       {Vijayaraghavan, Aravindan and Awasthi, Pranjal},
  booktitle = 	 {Proceedings of the 35th International Conference on Machine Learning},
  pages = 	 {5055--5064},
  year = 	 {2018},
  editor = 	 {Dy, Jennifer and Krause, Andreas},
  volume = 	 {80},
  series = 	 {Proceedings of Machine Learning Research},
  month = 	 {07},
  publisher =    {PMLR},
  url = 	 {https://proceedings.mlr.press/v80/vijayaraghavan18a.html},
eprint={1711.08841},
      archivePrefix={arXiv},
      primaryClass={cs.DS},
}

@inproceedings{
GC23,
title={Robust Matrix Sensing in the Semi-Random Model},
author={Xing Gao and Yu Cheng},
booktitle={Thirty-seventh Conference on Neural Information Processing Systems},
year={2023},
url={https://openreview.net/forum?id=nSr2epejn2}
}

@InProceedings{ls20,
  title = 	 {Costly Zero Order Oracles},
  author =       {Paes Leme, Renato and Schneider, Jon},
  booktitle = 	 {Proceedings of the Conference on Learning Theory},
  pages = 	 {3120--3132},
  year = 	 {2020},
  volume = 	 {125},
  month = 	 {07},
}

@misc{lu20,
      title={A Note on {J}ohn Simplex with Positive Dilation}, 
      author={Zhou Lu},
      year={2020},
      eprint={2012.03427},
      archivePrefix={arXiv},
      primaryClass={math.MG}
}

@ARTICLE{barvinok2012thrifty,
  author={Barvinok, Alexander},
  journal={International Mathematics Research Notices}, 
  title={Thrifty Approximations of Convex Bodies by Polytopes}, 
  year={2014},
  volume={2014},
  number={16},
  pages={4341-4356},
  doi={10.1093/imrn/rnt078},
eprint={1206.3993},
      archivePrefix={arXiv},
      primaryClass={math.MG},
}

@article{gk10,
title = {Determinants and the volumes of parallelotopes and zonotopes},
journal = {Linear Algebra and its Applications},
volume = {433},
number = {1},
pages = {28-40},
year = {2010},
author = {Eugene Gover and Nishan Krikorian},
}

@article{clarkson10,
author = {Clarkson, Kenneth L.},
title = {Coresets, Sparse Greedy Approximation, and the {Frank--Wolfe} Algorithm},
year = {2010},
issue_date = {August 2010},
volume = {6},
number = {4},
journal = {ACM Trans. Algorithms},
month = {09},
articleno = {63},
numpages = {30}
}

@inproceedings{
bmv23,
title={Tight Bounds for Volumetric Spanners and Applications},
author={Aditya Bhaskara and Sepideh Mahabadi and Ali Vakilian},
booktitle={Conference on Neural Information Processing Systems},
year={2023},
eprint={2310.00175},
      archivePrefix={arXiv},
      primaryClass={cs.DS},
}

@article{ty07,
author = {Todd, Michael J. and Yildirim, E. Alper},
title = {On {K}hachiyan's Algorithm for the Computation of Minimum-Volume Enclosing Ellipsoids},
year = {2007},
issue_date = {August, 2007},
volume = {155},
number = {13},
journal = {Discrete Appl. Math.},
month = {08},
pages = {1731–1744},
numpages = {14},
}

@article{ky05,
author = {Kumar, P. and Yildirim, E. A.},
title = {Minimum-Volume Enclosing Ellipsoids and Core Sets},
year = {2005},
issue_date = {July      2005},
publisher = {Plenum Press},
address = {USA},
volume = {126},
number = {1},
journal = {J. Optim. Theory Appl.},
month = {07},
pages = {1–21},
numpages = {21}
}

@InProceedings{blum2017approximate,
  author =	{Avrim Blum and Vladimir Braverman and Ananya Kumar and Harry Lang and Lin F. Yang},
  title =	{Approximate Convex Hull of Data Streams},
  booktitle =	{Proceedings of the International Colloquium on Automata, Languages, and  Programming (ICALP)},
  pages =	{21:1--21:13},
  year =	{2018},
  volume =	{107},
eprint={1712.04564},
      archivePrefix={arXiv},
      primaryClass={cs.CG},
}

@article{agarwal2005geometric,
  title={Geometric approximation via coresets},
  author={Agarwal, Pankaj K and Har-Peled, Sariel and Varadarajan, Kasturi R},
  journal={Combinatorial and computational geometry},
  volume={52},
  number={1},
  pages={1--30},
  year={2005}
}

@inproceedings{mmo22,
	title        = {Streaming Algorithms for Ellipsoidal Approximation of Convex Polytopes},
	author       = {Makarychev, Yury and Manoj, Naren Sarayu and Ovsiankin, Max},
	year         = 2022,
	booktitle    = {Proceedings of the Conference on Learning Theory},
        archivePrefix={arXiv},
          primaryClass={cs.DS},
          eprint={2206.07250},
}

@inproceedings{woodruff2022high,
  title={High-dimensional geometric streaming in polynomial space},
  author={Woodruff, David P and Yasuda, Taisuke},
  booktitle={Proceedings of the Symposium on Foundations of Computer Science},
  pages={732--743},
  year={2022},
eprint={2204.03790},
      archivePrefix={arXiv},
      primaryClass={cs.DS},
}

@book{todd16,
	title        = {Minimum-Volume Ellipsoids: Theory and Algorithms},
	author       = {Todd, Michael J.},
	year         = 2016
}

@inbook{horn1991,
	title        = {Topics in Matrix Analysis},
	author       = {Horn, Roger A. and Johnson, Charles R.},
	year         = 1991,
	publisher    = {Cambridge University Press},
	place        = {Cambridge}
}

@article{stange08,
	title        = {On the efficient update of the Singular Value Decomposition},
	author       = {Stange, Peter},
	year         = 2008,
	journal      = {PAMM},
	volume       = 8,
	number       = 1,
	pages        = {10827--10828}
}

@inproceedings{cv15,
	title        = {Bypassing {KLS}: {G}aussian Cooling and an {$O^{*}(N^3)$} Volume Algorithm},
	author       = {Cousins, Benjamin and Vempala, Santosh},
	year         = 2015,
	booktitle    = {Proceedings of the Symposium on Theory of Computing},
	pages        = {539–548},
eprint={1409.6011},
      archivePrefix={arXiv},
      primaryClass={cs.DS},
}

@article{boyd97,
	title        = {Obstacle Collision Detection Using Best Ellipsoid Fit},
	author       = {Rimon, Elon and Boyd, Stephen P.},
	year         = 1997,
	journal      = {J. Intell. Robotics Syst.},
	publisher    = {Kluwer Academic Publishers},
	address      = {USA},
	volume       = 18,
	number       = 2,
	pages        = {105–126},
	issue_date   = {February 1997},
	numpages     = 22
}

@inproceedings{he21,
	title        = {Reducing Isotropy and Volume to {KLS}: An {$O^{*}(n^3\psi^2)$} Volume Algorithm},
	author       = {Jia, He and Laddha, Aditi and Lee, Yin Tat and Vempala, Santosh},
	year         = 2021,
	booktitle    = {Proceedings of the Symposium on Theory of Computing},
	pages        = {961–974},
	numpages     = 14,
        eprint={2008.02146},
      archivePrefix={arXiv},
      primaryClass={cs.DS},
}

@inproceedings{mukhopadhyay2010approximate,
  title={Approximate ellipsoid in the streaming model},
  author={Mukhopadhyay, Asish and Sarker, Animesh and Switzer, Tom},
  booktitle={International Conference on Combinatorial Optimization and Applications},
  pages={401--413},
  year={2010}
}

@inproceedings{mukhopadhyayapproximate,
  title={Approximate minimum spanning ellipse in the streaming model},
  author={Mukhopadhyay, Asish and Greene, Eugene and Sarker, Animesh and Switzer, Tom},
  booktitle={The 7th Japan Conference on Computational Geometry and Graphs},
  year={2009}
}

@article{nesterov2008rounding,
  title={Rounding of convex sets and efficient gradient methods for linear programming problems},
  author={Nesterov, Yurii},
  journal={Optimisation Methods and Software},
  volume={23},
  number={1},
  pages={109--128},
  year={2008},
  publisher={Taylor \& Francis}
}

@article{howard1997john,
  title={The {John} ellipsoid theorem},
  author={Howard, Ralph},
  journal={University of South Carolina},
  year={1997}
}

@inproceedings{agarwal2010streaming,
  title={Streaming algorithms for extent problems in high dimensions},
  author={Agarwal, Pankaj K and Sharathkumar, R},
  booktitle={Proceedings of the  Symposium on Discrete Algorithms},
  pages={1481--1489},
  year={2010},
}

@book{artstein2015asymptotic,
  title={Asymptotic Geometric Analysis, Part I},
  author={Artstein-Avidan, Shiri and Giannopoulos, Apostolos and Milman, Vitali D},
  volume={202},
  year={2015},
  publisher={American Mathematical Soc.}
}

@inproceedings{bourgain2006estimates,
  title={Estimates related to Steiner symmetrizations},
  author={Bourgain, J and Lindenstrauss, J and Milman, V},
  booktitle={Geometric Aspects of Functional Analysis: Israel Seminar (GAFA) 1987--88},
  pages={264--273},
  year={1987},
  organization={Springer}
}

@article{bjkmmw24,
  title={On the Robustness of Spectral Algorithms for Semirandom Stochastic Block Models}, 
  author={Aditya Bhaskara and Agastya Vibhuti Jha and Michael Kapralov and {Naren} {Sarayu} {Manoj} and Davide Mazzali and Weronika Wrzos-Kaminska},
  year={2024},
  journal={Advances in Neural Information Processing Systems (NeurIPS)},
  volume={37},
    eprint={2412.14315},
      archivePrefix={arXiv},
      primaryClass={stat.ML},
}

@inproceedings{
mp24,
title={A Second-Order Algorithm for Empirical Group Distributionally Robust Regression},
author={{Naren} {Sarayu} {Manoj} and Kumar Kshitij Patel},
booktitle={OPT 2024: Optimization for Machine Learning},
year={2024},
url={https://openreview.net/forum?id=CPFpHb3vBm}
}

@InProceedings{mmo23,
      title={Near-Optimal Streaming Ellipsoidal Rounding for General Convex Polytopes}, 
      author={Yury Makarychev and {Naren} {Sarayu} {Manoj} and Max Ovsiankin},
      year={2024},
      month={06},
      booktitle = {Proceedings of Fifty Sixth Annual ACM Symposium on Theory of Computing (STOC)},
      archivePrefix={arXiv},
      primaryClass={cs.DS},
      eprint={2311.09460},
}

@article{bm21,
  title={Excess Capacity and Backdoor Poisoning},
  author={{Manoj}, {Naren} {Sarayu} and Blum, Avrim},
  journal={Advances in Neural Information Processing Systems (NeurIPS)},
  volume={34},
  year={2021},
  archivePrefix={arXiv},
  primaryClass={cs.LG},
  eprint={2109.00685},
}

@article{len83,
  title={Integer programming with a fixed number of variables},
  author={Lenstra Jr, Hendrik W},
  journal={Mathematics of operations research},
  volume={8},
  number={4},
  pages={538--548},
  year={1983},
  publisher={Informs}
}

@article{kha80,
  title={Polynomial algorithms in linear programming},
  author={Khachiyan, Leonid G},
  journal={USSR Computational Mathematics and Mathematical Physics},
  volume={20},
  number={1},
  pages={53--72},
  year={1980},
  publisher={Elsevier}
}

@article{valiant84,
author = {Valiant, L. G.},
title = {A Theory of the Learnable},
year = {1984},
issue_date = {Nov. 1984},
publisher = {Association for Computing Machinery},
address = {New York, NY, USA},
volume = {27},
number = {11},
issn = {0001-0782},
url = {https://doi.org/10.1145/1968.1972},
doi = {10.1145/1968.1972},
journal = {Commun. ACM},
month = {11},
pages = {1134–1142},
numpages = {9}
}

\end{document}